\documentclass[11pt,a4paper,reqno]{amsart}
\usepackage[a4paper,left=2.0cm,right=2.0cm,top=2.3cm,bottom=2.3cm]{geometry}
\usepackage{amssymb,amsmath,amsthm}
\usepackage{enumitem}
\usepackage{placeins}
\usepackage{tikz}
\usetikzlibrary{arrows.meta}
\definecolor{cbpInk}{RGB}{37,53,66}
\definecolor{cbpFluid}{RGB}{231,240,242}
\definecolor{cbpCore}{RGB}{151,61,45}
\definecolor{cbpFlow}{RGB}{35,102,117}
\tikzset{cbpwall/.style={draw=cbpInk,line width=.75pt},
cbpinner/.style={draw=cbpCore,line width=.9pt},
cbpleader/.style={draw=cbpInk!75,line width=.5pt},
cbpflow/.style={-{Stealth[length=2.1mm,width=1.4mm]},draw=cbpFlow,line width=.85pt},
cbpradius/.style={-{Stealth[length=1.15mm,width=.9mm]},draw=cbpInk,line width=.6pt},
cbpfigure/.style={x=1cm,y=1cm,font=\small,line cap=round,
line join=round,line width=.55pt,>={Stealth[length=2.1mm,width=1.4mm]}}}
\tikzset{cbpmesh/.style={fill=cbpFlow!#1!white,draw=cbpFlow!#1!white,line width=.12pt},
cbpaxis/.style={-{Stealth[length=1.7mm,width=1.2mm]},
draw=cbpInk!70,line width=.5pt}}
\usepackage[linkcolor=blue,citecolor=cyan,urlcolor=blue,colorlinks,hypertexnames=false]{hyperref}
\usepackage[T1]{fontenc}
\usepackage{lmodern}

\newtheorem{theorem}{Theorem}[section]
\newtheorem{lemma}[theorem]{Lemma}
\newtheorem{proposition}[theorem]{Proposition}
\newtheorem{corollary}[theorem]{Corollary}
\newtheorem{remark}[theorem]{Remark}

\newcommand{\R}{\mathbb R}
\newcommand{\T}{\mathbb T}

\numberwithin{equation}{section}
\allowdisplaybreaks
\def\cbpOuterSection{(-2.05,-.22)
.. controls (-2.25,.75) and (-1.45,1.25) .. (-.52,.96)
.. controls (.38,.67) and (.68,1.00) .. (1.25,.53)
.. controls (1.94,-.29) and (.96,-.95) .. (.16,-.78)
.. controls (-.56,-.60) and (-1.74,-1.19) .. (-2.05,-.22)--cycle}

\newcommand{\cbpRotatedOuter}{
\path[cbpmesh=24.9] (0.3118,2.6260) -- (0.0000,2.6318) -- (0.0000,2.6353) -- (0.3170,2.6294) -- cycle;
\draw[draw=cbpInk!70,line width=.45pt] (0.0000,2.6353) -- (0.3170,2.6294);
\path[cbpmesh=24.7] (0.0000,2.6318) -- (-0.3118,2.6260) -- (-0.3170,2.6294) -- (0.0000,2.6353) -- cycle;
\draw[draw=cbpInk!70,line width=.45pt] (-0.3170,2.6294) -- (0.0000,2.6353);
\path[cbpmesh=25.2] (0.6216,2.6088) -- (0.3118,2.6260) -- (0.3170,2.6294) -- (0.6320,2.6119) -- cycle;
\draw[draw=cbpInk!70,line width=.45pt] (0.3170,2.6294) -- (0.6320,2.6119);
\path[cbpmesh=24.4] (-0.3118,2.6260) -- (-0.6216,2.6088) -- (-0.6320,2.6119) -- (-0.3170,2.6294) -- cycle;
\draw[draw=cbpInk!70,line width=.45pt] (-0.6320,2.6119) -- (-0.3170,2.6294);
\path[cbpmesh=25.4] (0.9276,2.5802) -- (0.6216,2.6088) -- (0.6320,2.6119) -- (0.9431,2.5829) -- cycle;
\draw[draw=cbpInk!70,line width=.45pt] (0.6320,2.6119) -- (0.9431,2.5829);
\path[cbpmesh=24.0] (-0.6216,2.6088) -- (-0.9276,2.5802) -- (-0.9431,2.5829) -- (-0.6320,2.6119) -- cycle;
\draw[draw=cbpInk!70,line width=.45pt] (-0.9431,2.5829) -- (-0.6320,2.6119);
\path[cbpmesh=23.0] (0.0000,2.6219) -- (-0.3066,2.6162) -- (-0.3118,2.6260) -- (0.0000,2.6318) -- cycle;
\path[cbpmesh=23.2] (0.3066,2.6162) -- (0.0000,2.6219) -- (0.0000,2.6318) -- (0.3118,2.6260) -- cycle;
\path[cbpmesh=23.5] (0.6113,2.5993) -- (0.3066,2.6162) -- (0.3118,2.6260) -- (0.6216,2.6088) -- cycle;
\path[cbpmesh=22.8] (-0.3066,2.6162) -- (-0.6113,2.5993) -- (-0.6216,2.6088) -- (-0.3118,2.6260) -- cycle;
\path[cbpmesh=25.6] (1.2279,2.5405) -- (0.9276,2.5802) -- (0.9431,2.5829) -- (1.2485,2.5424) -- cycle;
\draw[draw=cbpInk!70,line width=.45pt] (0.9431,2.5829) -- (1.2485,2.5424);
\path[cbpmesh=23.7] (-0.9276,2.5802) -- (-1.2279,2.5405) -- (-1.2485,2.5424) -- (-0.9431,2.5829) -- cycle;
\draw[draw=cbpInk!70,line width=.45pt] (-1.2485,2.5424) -- (-0.9431,2.5829);
\path[cbpmesh=23.6] (0.9123,2.5712) -- (0.6113,2.5993) -- (0.6216,2.6088) -- (0.9276,2.5802) -- cycle;
\path[cbpmesh=22.5] (-0.6113,2.5993) -- (-0.9123,2.5712) -- (-0.9276,2.5802) -- (-0.6216,2.6088) -- cycle;
\path[cbpmesh=21.8] (0.3014,2.6007) -- (0.0000,2.6063) -- (0.0000,2.6219) -- (0.3066,2.6162) -- cycle;
\path[cbpmesh=21.6] (0.0000,2.6063) -- (-0.3014,2.6007) -- (-0.3066,2.6162) -- (0.0000,2.6219) -- cycle;
\path[cbpmesh=21.9] (0.6010,2.5841) -- (0.3014,2.6007) -- (0.3066,2.6162) -- (0.6113,2.5993) -- cycle;
\path[cbpmesh=21.4] (-0.3014,2.6007) -- (-0.6010,2.5841) -- (-0.6113,2.5993) -- (-0.3066,2.6162) -- cycle;
\path[cbpmesh=22.2] (-0.9123,2.5712) -- (-1.2076,2.5321) -- (-1.2279,2.5405) -- (-0.9276,2.5802) -- cycle;
\path[cbpmesh=23.8] (1.2076,2.5321) -- (0.9123,2.5712) -- (0.9276,2.5802) -- (1.2279,2.5405) -- cycle;
\path[cbpmesh=21.2] (-0.6010,2.5841) -- (-0.8969,2.5565) -- (-0.9123,2.5712) -- (-0.6113,2.5993) -- cycle;
\path[cbpmesh=22.1] (0.8969,2.5565) -- (0.6010,2.5841) -- (0.6113,2.5993) -- (0.9123,2.5712) -- cycle;
\path[cbpmesh=20.6] (0.2963,2.5804) -- (0.0000,2.5859) -- (0.0000,2.6063) -- (0.3014,2.6007) -- cycle;
\path[cbpmesh=20.5] (0.0000,2.5859) -- (-0.2963,2.5804) -- (-0.3014,2.6007) -- (0.0000,2.6063) -- cycle;
\path[cbpmesh=23.3] (-1.2279,2.5405) -- (-1.5206,2.4898) -- (-1.5461,2.4909) -- (-1.2485,2.5424) -- cycle;
\draw[draw=cbpInk!70,line width=.45pt] (-1.5461,2.4909) -- (-1.2485,2.5424);
\path[cbpmesh=25.8] (1.5206,2.4898) -- (1.2279,2.5405) -- (1.2485,2.5424) -- (1.5461,2.4909) -- cycle;
\draw[draw=cbpInk!70,line width=.45pt] (1.2485,2.5424) -- (1.5461,2.4909);
\path[cbpmesh=20.4] (-0.2963,2.5804) -- (-0.5907,2.5641) -- (-0.6010,2.5841) -- (-0.3014,2.6007) -- cycle;
\path[cbpmesh=20.7] (0.5907,2.5641) -- (0.2963,2.5804) -- (0.3014,2.6007) -- (0.6010,2.5841) -- cycle;
\path[cbpmesh=21.0] (-0.8969,2.5565) -- (-1.1873,2.5180) -- (-1.2076,2.5321) -- (-0.9123,2.5712) -- cycle;
\path[cbpmesh=22.2] (1.1873,2.5180) -- (0.8969,2.5565) -- (0.9123,2.5712) -- (1.2076,2.5321) -- cycle;
\path[cbpmesh=20.8] (0.8815,2.5369) -- (0.5907,2.5641) -- (0.6010,2.5841) -- (0.8969,2.5565) -- cycle;
\path[cbpmesh=20.2] (-0.5907,2.5641) -- (-0.8815,2.5369) -- (-0.8969,2.5565) -- (-0.6010,2.5841) -- cycle;
\path[cbpmesh=19.8] (0.2910,2.5561) -- (0.0000,2.5615) -- (0.0000,2.5859) -- (0.2963,2.5804) -- cycle;
\path[cbpmesh=19.7] (0.0000,2.5615) -- (-0.2910,2.5561) -- (-0.2963,2.5804) -- (0.0000,2.5859) -- cycle;
\path[cbpmesh=23.9] (1.4954,2.4822) -- (1.2076,2.5321) -- (1.2279,2.5405) -- (1.5206,2.4898) -- cycle;
\path[cbpmesh=21.9] (-1.2076,2.5321) -- (-1.4954,2.4822) -- (-1.5206,2.4898) -- (-1.2279,2.5405) -- cycle;
\path[cbpmesh=19.8] (0.5802,2.5400) -- (0.2910,2.5561) -- (0.2963,2.5804) -- (0.5907,2.5641) -- cycle;
\path[cbpmesh=19.6] (-0.2910,2.5561) -- (-0.5802,2.5400) -- (-0.5907,2.5641) -- (-0.2963,2.5804) -- cycle;
\path[cbpmesh=25.9] (1.8040,2.4285) -- (1.5206,2.4898) -- (1.5461,2.4909) -- (1.8342,2.4285) -- cycle;
\draw[draw=cbpInk!70,line width=.45pt] (1.5461,2.4909) -- (1.8342,2.4285);
\path[cbpmesh=22.9] (-1.5206,2.4898) -- (-1.8040,2.4285) -- (-1.8342,2.4285) -- (-1.5461,2.4909) -- cycle;
\draw[draw=cbpInk!70,line width=.45pt] (-1.8342,2.4285) -- (-1.5461,2.4909);
\path[cbpmesh=20.1] (-0.8815,2.5369) -- (-1.1668,2.4991) -- (-1.1873,2.5180) -- (-0.8969,2.5565) -- cycle;
\path[cbpmesh=20.9] (1.1668,2.4991) -- (0.8815,2.5369) -- (0.8969,2.5565) -- (1.1873,2.5180) -- cycle;
\path[cbpmesh=19.9] (0.8658,2.5134) -- (0.5802,2.5400) -- (0.5907,2.5641) -- (0.8815,2.5369) -- cycle;
\path[cbpmesh=19.5] (-0.5802,2.5400) -- (-0.8658,2.5134) -- (-0.8815,2.5369) -- (-0.5907,2.5641) -- cycle;
\path[cbpmesh=19.2] (0.2856,2.5286) -- (0.0000,2.5338) -- (0.0000,2.5615) -- (0.2910,2.5561) -- cycle;
\path[cbpmesh=19.2] (0.0000,2.5338) -- (-0.2856,2.5286) -- (-0.2910,2.5561) -- (0.0000,2.5615) -- cycle;
\path[cbpmesh=22.3] (1.4703,2.4690) -- (1.1873,2.5180) -- (1.2076,2.5321) -- (1.4954,2.4822) -- cycle;
\path[cbpmesh=20.8] (-1.1873,2.5180) -- (-1.4703,2.4690) -- (-1.4954,2.4822) -- (-1.2076,2.5321) -- cycle;
\path[cbpmesh=19.3] (0.5694,2.5128) -- (0.2856,2.5286) -- (0.2910,2.5561) -- (0.5802,2.5400) -- cycle;
\path[cbpmesh=19.1] (-0.2856,2.5286) -- (-0.5694,2.5128) -- (-0.5802,2.5400) -- (-0.2910,2.5561) -- cycle;
\path[cbpmesh=28.0] (2.1109,2.3558) -- (1.8342,2.4285) -- (1.8650,2.4213) -- (2.1464,2.3473) -- cycle;
\draw[draw=cbpInk!70,line width=.45pt] (1.8342,2.4285) -- (1.8650,2.4213);
\draw[draw=cbpInk!70,line width=.45pt] (1.8650,2.4213) -- (2.1464,2.3473);
\path[cbpmesh=23.8] (-1.8342,2.4285) -- (-2.1109,2.3558) -- (-2.1464,2.3473) -- (-1.8650,2.4213) -- cycle;
\draw[draw=cbpInk!70,line width=.45pt] (-2.1464,2.3473) -- (-1.8650,2.4213);
\draw[draw=cbpInk!70,line width=.45pt] (-1.8650,2.4213) -- (-1.8342,2.4285);
\path[cbpmesh=21.6] (-1.4954,2.4822) -- (-1.7741,2.4219) -- (-1.8040,2.4285) -- (-1.5206,2.4898) -- cycle;
\path[cbpmesh=24.0] (1.7741,2.4219) -- (1.4954,2.4822) -- (1.5206,2.4898) -- (1.8040,2.4285) -- cycle;
\path[cbpmesh=19.4] (-0.8658,2.5134) -- (-1.1461,2.4763) -- (-1.1668,2.4991) -- (-0.8815,2.5369) -- cycle;
\path[cbpmesh=19.9] (1.1461,2.4763) -- (0.8658,2.5134) -- (0.8815,2.5369) -- (1.1668,2.4991) -- cycle;
\path[cbpmesh=19.3] (0.8497,2.4866) -- (0.5694,2.5128) -- (0.5802,2.5400) -- (0.8658,2.5134) -- cycle;
\path[cbpmesh=19.0] (-0.5694,2.5128) -- (-0.8497,2.4866) -- (-0.8658,2.5134) -- (-0.5802,2.5400) -- cycle;
\path[cbpmesh=18.9] (0.2800,2.4986) -- (0.0000,2.5038) -- (0.0000,2.5338) -- (0.2856,2.5286) -- cycle;
\path[cbpmesh=18.9] (0.0000,2.5038) -- (-0.2800,2.4986) -- (-0.2856,2.5286) -- (0.0000,2.5338) -- cycle;
\path[cbpmesh=19.9] (-1.1668,2.4991) -- (-1.4450,2.4510) -- (-1.4703,2.4690) -- (-1.1873,2.5180) -- cycle;
\path[cbpmesh=21.0] (1.4450,2.4510) -- (1.1668,2.4991) -- (1.1873,2.5180) -- (1.4703,2.4690) -- cycle;
\path[cbpmesh=18.9] (0.5583,2.4831) -- (0.2800,2.4986) -- (0.2856,2.5286) -- (0.5694,2.5128) -- cycle;
\path[cbpmesh=18.8] (-0.2800,2.4986) -- (-0.5583,2.4831) -- (-0.5694,2.5128) -- (-0.2856,2.5286) -- cycle;
\path[cbpmesh=22.5] (-1.8040,2.4285) -- (-2.0762,2.3569) -- (-2.1109,2.3558) -- (-1.8342,2.4285) -- cycle;
\path[cbpmesh=26.0] (2.0762,2.3569) -- (1.8040,2.4285) -- (1.8342,2.4285) -- (2.1109,2.3558) -- cycle;
\path[cbpmesh=19.0] (-0.8497,2.4866) -- (-1.1248,2.4502) -- (-1.1461,2.4763) -- (-0.8658,2.5134) -- cycle;
\path[cbpmesh=19.3] (1.1248,2.4502) -- (0.8497,2.4866) -- (0.8658,2.5134) -- (1.1461,2.4763) -- cycle;
\path[cbpmesh=20.6] (-1.4703,2.4690) -- (-1.7443,2.4097) -- (-1.7741,2.4219) -- (-1.4954,2.4822) -- cycle;
\path[cbpmesh=22.3] (1.7443,2.4097) -- (1.4703,2.4690) -- (1.4954,2.4822) -- (1.7741,2.4219) -- cycle;
\path[cbpmesh=19.0] (0.8331,2.4575) -- (0.5583,2.4831) -- (0.5694,2.5128) -- (0.8497,2.4866) -- cycle;
\path[cbpmesh=18.8] (-0.5583,2.4831) -- (-0.8331,2.4575) -- (-0.8497,2.4866) -- (-0.5694,2.5128) -- cycle;
\path[cbpmesh=18.8] (0.2742,2.4670) -- (0.0000,2.4721) -- (0.0000,2.5038) -- (0.2800,2.4986) -- cycle;
\path[cbpmesh=18.8] (0.0000,2.4721) -- (-0.2742,2.4670) -- (-0.2800,2.4986) -- (0.0000,2.5038) -- cycle;
\path[cbpmesh=20.0] (1.4193,2.4289) -- (1.1461,2.4763) -- (1.1668,2.4991) -- (1.4450,2.4510) -- cycle;
\path[cbpmesh=19.3] (-1.1461,2.4763) -- (-1.4193,2.4289) -- (-1.4450,2.4510) -- (-1.1668,2.4991) -- cycle;
\path[cbpmesh=18.9] (0.5467,2.4519) -- (0.2742,2.4670) -- (0.2800,2.4986) -- (0.5583,2.4831) -- cycle;
\path[cbpmesh=18.8] (-0.2742,2.4670) -- (-0.5467,2.4519) -- (-0.5583,2.4831) -- (-0.2800,2.4986) -- cycle;
\path[cbpmesh=21.0] (1.7142,2.3927) -- (1.4450,2.4510) -- (1.4703,2.4690) -- (1.7443,2.4097) -- cycle;
\path[cbpmesh=19.7] (-1.4450,2.4510) -- (-1.7142,2.3927) -- (-1.7443,2.4097) -- (-1.4703,2.4690) -- cycle;
\path[cbpmesh=18.7] (-0.8331,2.4575) -- (-1.1028,2.4218) -- (-1.1248,2.4502) -- (-0.8497,2.4866) -- cycle;
\path[cbpmesh=19.0] (1.1028,2.4218) -- (0.8331,2.4575) -- (0.8497,2.4866) -- (1.1248,2.4502) -- cycle;
\path[cbpmesh=21.3] (-1.7741,2.4219) -- (-2.0418,2.3515) -- (-2.0762,2.3569) -- (-1.8040,2.4285) -- cycle;
\path[cbpmesh=24.1] (2.0418,2.3515) -- (1.7741,2.4219) -- (1.8040,2.4285) -- (2.0762,2.3569) -- cycle;
\path[cbpmesh=18.9] (0.8159,2.4267) -- (0.5467,2.4519) -- (0.5583,2.4831) -- (0.8331,2.4575) -- cycle;
\path[cbpmesh=18.7] (-0.5467,2.4519) -- (-0.8159,2.4267) -- (-0.8331,2.4575) -- (-0.5583,2.4831) -- cycle;
\path[cbpmesh=28.1] (2.3747,2.2730) -- (2.1109,2.3558) -- (2.1464,2.3473) -- (2.4146,2.2632) -- cycle;
\draw[draw=cbpInk!70,line width=.45pt] (2.1464,2.3473) -- (2.4146,2.2632);
\path[cbpmesh=23.3] (-2.1109,2.3558) -- (-2.3747,2.2730) -- (-2.4146,2.2632) -- (-2.1464,2.3473) -- cycle;
\draw[draw=cbpInk!70,line width=.45pt] (-2.4146,2.2632) -- (-2.1464,2.3473);
\path[cbpmesh=18.9] (0.2681,2.4347) -- (0.0000,2.4396) -- (0.0000,2.4721) -- (0.2742,2.4670) -- cycle;
\path[cbpmesh=18.9] (0.0000,2.4396) -- (-0.2681,2.4347) -- (-0.2742,2.4670) -- (0.0000,2.4721) -- cycle;
\path[cbpmesh=19.4] (1.3929,2.4037) -- (1.1248,2.4502) -- (1.1461,2.4763) -- (1.4193,2.4289) -- cycle;
\path[cbpmesh=18.9] (-1.1248,2.4502) -- (-1.3929,2.4037) -- (-1.4193,2.4289) -- (-1.1461,2.4763) -- cycle;
\path[cbpmesh=19.0] (0.5346,2.4198) -- (0.2681,2.4347) -- (0.2742,2.4670) -- (0.5467,2.4519) -- cycle;
\path[cbpmesh=18.9] (-0.2681,2.4347) -- (-0.5346,2.4198) -- (-0.5467,2.4519) -- (-0.2742,2.4670) -- cycle;
\path[cbpmesh=20.0] (1.6837,2.3717) -- (1.4193,2.4289) -- (1.4450,2.4510) -- (1.7142,2.3927) -- cycle;
\path[cbpmesh=19.2] (-1.4193,2.4289) -- (-1.6837,2.3717) -- (-1.7142,2.3927) -- (-1.4450,2.4510) -- cycle;
\path[cbpmesh=20.3] (-1.7443,2.4097) -- (-2.0075,2.3405) -- (-2.0418,2.3515) -- (-1.7741,2.4219) -- cycle;
\path[cbpmesh=22.4] (2.0075,2.3405) -- (1.7443,2.4097) -- (1.7741,2.4219) -- (2.0418,2.3515) -- cycle;
\path[cbpmesh=18.7] (-0.8159,2.4267) -- (-1.0800,2.3918) -- (-1.1028,2.4218) -- (-0.8331,2.4575) -- cycle;
\path[cbpmesh=18.9] (1.0800,2.3918) -- (0.8159,2.4267) -- (0.8331,2.4575) -- (1.1028,2.4218) -- cycle;
\path[cbpmesh=19.0] (0.7978,2.3953) -- (0.5346,2.4198) -- (0.5467,2.4519) -- (0.8159,2.4267) -- cycle;
\path[cbpmesh=18.8] (-0.5346,2.4198) -- (-0.7978,2.3953) -- (-0.8159,2.4267) -- (-0.5467,2.4519) -- cycle;
\path[cbpmesh=18.7] (-1.1028,2.4218) -- (-1.3657,2.3762) -- (-1.3929,2.4037) -- (-1.1248,2.4502) -- cycle;
\path[cbpmesh=19.0] (1.3657,2.3762) -- (1.1028,2.4218) -- (1.1248,2.4502) -- (1.3929,2.4037) -- cycle;
\path[cbpmesh=26.0] (2.3356,2.2755) -- (2.0762,2.3569) -- (2.1109,2.3558) -- (2.3747,2.2730) -- cycle;
\path[cbpmesh=22.0] (-2.0762,2.3569) -- (-2.3356,2.2755) -- (-2.3747,2.2730) -- (-2.1109,2.3558) -- cycle;
\path[cbpmesh=19.2] (0.2618,2.4023) -- (0.0000,2.4071) -- (0.0000,2.4396) -- (0.2681,2.4347) -- cycle;
\path[cbpmesh=19.2] (0.0000,2.4071) -- (-0.2618,2.4023) -- (-0.2681,2.4347) -- (0.0000,2.4396) -- cycle;
\path[cbpmesh=19.1] (-0.2618,2.4023) -- (-0.5219,2.3878) -- (-0.5346,2.4198) -- (-0.2681,2.4347) -- cycle;
\path[cbpmesh=19.3] (0.5219,2.3878) -- (0.2618,2.4023) -- (0.2681,2.4347) -- (0.5346,2.4198) -- cycle;
\path[cbpmesh=18.8] (-1.3929,2.4037) -- (-1.6525,2.3476) -- (-1.6837,2.3717) -- (-1.4193,2.4289) -- cycle;
\path[cbpmesh=19.4] (1.6525,2.3476) -- (1.3929,2.4037) -- (1.4193,2.4289) -- (1.6837,2.3717) -- cycle;
\path[cbpmesh=21.0] (1.9729,2.3247) -- (1.7142,2.3927) -- (1.7443,2.4097) -- (2.0075,2.3405) -- cycle;
\path[cbpmesh=19.6] (-1.7142,2.3927) -- (-1.9729,2.3247) -- (-2.0075,2.3405) -- (-1.7443,2.4097) -- cycle;
\path[cbpmesh=18.8] (-0.7978,2.3953) -- (-1.0561,2.3611) -- (-1.0800,2.3918) -- (-0.8159,2.4267) -- cycle;
\path[cbpmesh=19.0] (1.0561,2.3611) -- (0.7978,2.3953) -- (0.8159,2.4267) -- (1.0800,2.3918) -- cycle;
\path[cbpmesh=18.9] (1.3374,2.3472) -- (1.0800,2.3918) -- (1.1028,2.4218) -- (1.3657,2.3762) -- cycle;
\path[cbpmesh=18.6] (-1.0800,2.3918) -- (-1.3374,2.3472) -- (-1.3657,2.3762) -- (-1.1028,2.4218) -- cycle;
\path[cbpmesh=20.9] (-2.0418,2.3515) -- (-2.2969,2.2715) -- (-2.3356,2.2755) -- (-2.0762,2.3569) -- cycle;
\path[cbpmesh=24.1] (2.2969,2.2715) -- (2.0418,2.3515) -- (2.0762,2.3569) -- (2.3356,2.2755) -- cycle;
\path[cbpmesh=19.3] (0.7789,2.3638) -- (0.5219,2.3878) -- (0.5346,2.4198) -- (0.7978,2.3953) -- cycle;
\path[cbpmesh=19.1] (-0.5219,2.3878) -- (-0.7789,2.3638) -- (-0.7978,2.3953) -- (-0.5346,2.4198) -- cycle;
\path[cbpmesh=19.7] (0.2550,2.3707) -- (0.0000,2.3754) -- (0.0000,2.4071) -- (0.2618,2.4023) -- cycle;
\path[cbpmesh=19.6] (0.0000,2.3754) -- (-0.2550,2.3707) -- (-0.2618,2.4023) -- (0.0000,2.4071) -- cycle;
\path[cbpmesh=18.7] (-1.3657,2.3762) -- (-1.6202,2.3211) -- (-1.6525,2.3476) -- (-1.3929,2.4037) -- cycle;
\path[cbpmesh=19.0] (1.6202,2.3211) -- (1.3657,2.3762) -- (1.3929,2.4037) -- (1.6525,2.3476) -- cycle;
\path[cbpmesh=20.1] (1.9378,2.3049) -- (1.6837,2.3717) -- (1.7142,2.3927) -- (1.9729,2.3247) -- cycle;
\path[cbpmesh=19.1] (-1.6837,2.3717) -- (-1.9378,2.3049) -- (-1.9729,2.3247) -- (-1.7142,2.3927) -- cycle;
\path[cbpmesh=28.1] (2.6238,2.1808) -- (2.3747,2.2730) -- (2.4146,2.2632) -- (2.6679,2.1694) -- cycle;
\draw[draw=cbpInk!70,line width=.45pt] (2.4146,2.2632) -- (2.6679,2.1694);
\path[cbpmesh=22.7] (-2.3747,2.2730) -- (-2.6238,2.1808) -- (-2.6679,2.1694) -- (-2.4146,2.2632) -- cycle;
\draw[draw=cbpInk!70,line width=.45pt] (-2.6679,2.1694) -- (-2.4146,2.2632);
\path[cbpmesh=19.7] (0.5085,2.3567) -- (0.2550,2.3707) -- (0.2618,2.4023) -- (0.5219,2.3878) -- cycle;
\path[cbpmesh=19.5] (-0.2550,2.3707) -- (-0.5085,2.3567) -- (-0.5219,2.3878) -- (-0.2618,2.4023) -- cycle;
\path[cbpmesh=19.0] (-0.7789,2.3638) -- (-1.0310,2.3305) -- (-1.0561,2.3611) -- (-0.7978,2.3953) -- cycle;
\path[cbpmesh=19.4] (1.0310,2.3305) -- (0.7789,2.3638) -- (0.7978,2.3953) -- (1.0561,2.3611) -- cycle;
\path[cbpmesh=20.0] (-2.0075,2.3405) -- (-2.2583,2.2618) -- (-2.2969,2.2715) -- (-2.0418,2.3515) -- cycle;
\path[cbpmesh=22.4] (2.2583,2.2618) -- (2.0075,2.3405) -- (2.0418,2.3515) -- (2.2969,2.2715) -- cycle;
\path[cbpmesh=18.7] (-1.0561,2.3611) -- (-1.3079,2.3175) -- (-1.3374,2.3472) -- (-1.0800,2.3918) -- cycle;
\path[cbpmesh=19.1] (1.3079,2.3175) -- (1.0561,2.3611) -- (1.0800,2.3918) -- (1.3374,2.3472) -- cycle;
\path[cbpmesh=19.8] (0.7588,2.3333) -- (0.5085,2.3567) -- (0.5219,2.3878) -- (0.7789,2.3638) -- cycle;
\path[cbpmesh=19.4] (-0.5085,2.3567) -- (-0.7588,2.3333) -- (-0.7789,2.3638) -- (-0.5219,2.3878) -- cycle;
\path[cbpmesh=19.4] (1.9018,2.2820) -- (1.6525,2.3476) -- (1.6837,2.3717) -- (1.9378,2.3049) -- cycle;
\path[cbpmesh=18.8] (-1.6525,2.3476) -- (-1.9018,2.2820) -- (-1.9378,2.3049) -- (-1.6837,2.3717) -- cycle;
\path[cbpmesh=18.9] (1.5867,2.2932) -- (1.3374,2.3472) -- (1.3657,2.3762) -- (1.6202,2.3211) -- cycle;
\path[cbpmesh=18.6] (-1.3374,2.3472) -- (-1.5867,2.2932) -- (-1.6202,2.3211) -- (-1.3657,2.3762) -- cycle;
\path[cbpmesh=26.0] (2.5806,2.1848) -- (2.3356,2.2755) -- (2.3747,2.2730) -- (2.6238,2.1808) -- cycle;
\path[cbpmesh=21.6] (-2.3356,2.2755) -- (-2.5806,2.1848) -- (-2.6238,2.1808) -- (-2.3747,2.2730) -- cycle;
\path[cbpmesh=20.3] (0.2479,2.3408) -- (0.0000,2.3454) -- (0.0000,2.3754) -- (0.2550,2.3707) -- cycle;
\path[cbpmesh=20.2] (0.0000,2.3454) -- (-0.2479,2.3408) -- (-0.2550,2.3707) -- (0.0000,2.3754) -- cycle;
\path[cbpmesh=21.1] (2.2194,2.2473) -- (1.9729,2.3247) -- (2.0075,2.3405) -- (2.2583,2.2618) -- cycle;
\path[cbpmesh=19.4] (-1.9729,2.3247) -- (-2.2194,2.2473) -- (-2.2583,2.2618) -- (-2.0075,2.3405) -- cycle;
\path[cbpmesh=20.4] (0.4942,2.3271) -- (0.2479,2.3408) -- (0.2550,2.3707) -- (0.5085,2.3567) -- cycle;
\path[cbpmesh=20.1] (-0.2479,2.3408) -- (-0.4942,2.3271) -- (-0.5085,2.3567) -- (-0.2550,2.3707) -- cycle;
\path[cbpmesh=19.4] (1.2768,2.2879) -- (1.0310,2.3305) -- (1.0561,2.3611) -- (1.3079,2.3175) -- cycle;
\path[cbpmesh=18.9] (-1.0310,2.3305) -- (-1.2768,2.2879) -- (-1.3079,2.3175) -- (-1.0561,2.3611) -- cycle;
\path[cbpmesh=19.9] (1.0045,2.3008) -- (0.7588,2.3333) -- (0.7789,2.3638) -- (1.0310,2.3305) -- cycle;
\path[cbpmesh=19.3] (-0.7588,2.3333) -- (-1.0045,2.3008) -- (-1.0310,2.3305) -- (-0.7789,2.3638) -- cycle;
\path[cbpmesh=18.6] (-1.6202,2.3211) -- (-1.8647,2.2569) -- (-1.9018,2.2820) -- (-1.6525,2.3476) -- cycle;
\path[cbpmesh=19.0] (1.8647,2.2569) -- (1.6202,2.3211) -- (1.6525,2.3476) -- (1.9018,2.2820) -- cycle;
\path[cbpmesh=19.1] (1.5516,2.2647) -- (1.3079,2.3175) -- (1.3374,2.3472) -- (1.5867,2.2932) -- cycle;
\path[cbpmesh=18.7] (-1.3079,2.3175) -- (-1.5516,2.2647) -- (-1.5867,2.2932) -- (-1.3374,2.3472) -- cycle;
\path[cbpmesh=20.4] (0.7376,2.3044) -- (0.4942,2.3271) -- (0.5085,2.3567) -- (0.7588,2.3333) -- cycle;
\path[cbpmesh=19.9] (-0.4942,2.3271) -- (-0.7376,2.3044) -- (-0.7588,2.3333) -- (-0.5085,2.3567) -- cycle;
\path[cbpmesh=24.1] (2.5379,2.1823) -- (2.2969,2.2715) -- (2.3356,2.2755) -- (2.5806,2.1848) -- cycle;
\path[cbpmesh=20.6] (-2.2969,2.2715) -- (-2.5379,2.1823) -- (-2.5806,2.1848) -- (-2.3356,2.2755) -- cycle;
\path[cbpmesh=20.1] (2.1799,2.2289) -- (1.9378,2.3049) -- (1.9729,2.3247) -- (2.2194,2.2473) -- cycle;
\path[cbpmesh=19.0] (-1.9378,2.3049) -- (-2.1799,2.2289) -- (-2.2194,2.2473) -- (-1.9729,2.3247) -- cycle;
\path[cbpmesh=21.0] (0.2403,2.3133) -- (0.0000,2.3177) -- (0.0000,2.3454) -- (0.2479,2.3408) -- cycle;
\path[cbpmesh=20.8] (0.0000,2.3177) -- (-0.2403,2.3133) -- (-0.2479,2.3408) -- (0.0000,2.3454) -- cycle;
\path[cbpmesh=18.6] (-1.5867,2.2932) -- (-1.8261,2.2303) -- (-1.8647,2.2569) -- (-1.6202,2.3211) -- cycle;
\path[cbpmesh=18.9] (1.8261,2.2303) -- (1.5867,2.2932) -- (1.6202,2.3211) -- (1.8647,2.2569) -- cycle;
\path[cbpmesh=19.9] (1.2439,2.2593) -- (1.0045,2.3008) -- (1.0310,2.3305) -- (1.2768,2.2879) -- cycle;
\path[cbpmesh=19.2] (-1.0045,2.3008) -- (-1.2439,2.2593) -- (-1.2768,2.2879) -- (-1.0310,2.3305) -- cycle;
\path[cbpmesh=20.5] (0.9763,2.2728) -- (0.7376,2.3044) -- (0.7588,2.3333) -- (1.0045,2.3008) -- cycle;
\path[cbpmesh=19.8] (-0.7376,2.3044) -- (-0.9763,2.2728) -- (-1.0045,2.3008) -- (-0.7588,2.3333) -- cycle;
\path[cbpmesh=20.7] (-0.2403,2.3133) -- (-0.4791,2.3000) -- (-0.4942,2.3271) -- (-0.2479,2.3408) -- cycle;
\path[cbpmesh=21.1] (0.4791,2.3000) -- (0.2403,2.3133) -- (0.2479,2.3408) -- (0.4942,2.3271) -- cycle;
\path[cbpmesh=22.4] (2.4952,2.1741) -- (2.2583,2.2618) -- (2.2969,2.2715) -- (2.5379,2.1823) -- cycle;
\path[cbpmesh=19.8] (-2.2583,2.2618) -- (-2.4952,2.1741) -- (-2.5379,2.1823) -- (-2.2969,2.2715) -- cycle;
\path[cbpmesh=18.9] (-1.2768,2.2879) -- (-1.5147,2.2364) -- (-1.5516,2.2647) -- (-1.3079,2.3175) -- cycle;
\path[cbpmesh=19.4] (1.5147,2.2364) -- (1.2768,2.2879) -- (1.3079,2.3175) -- (1.5516,2.2647) -- cycle;
\path[cbpmesh=19.4] (2.1395,2.2075) -- (1.9018,2.2820) -- (1.9378,2.3049) -- (2.1799,2.2289) -- cycle;
\path[cbpmesh=18.7] (-1.9018,2.2820) -- (-2.1395,2.2075) -- (-2.1799,2.2289) -- (-1.9378,2.3049) -- cycle;
\path[cbpmesh=28.0] (2.8568,2.0797) -- (2.6238,2.1808) -- (2.6679,2.1694) -- (2.9047,2.0666) -- cycle;
\draw[draw=cbpInk!70,line width=.45pt] (2.6679,2.1694) -- (2.9047,2.0666);
\path[cbpmesh=22.2] (-2.6238,2.1808) -- (-2.8568,2.0797) -- (-2.9047,2.0666) -- (-2.6679,2.1694) -- cycle;
\draw[draw=cbpInk!70,line width=.45pt] (-2.9047,2.0666) -- (-2.6679,2.1694);
\path[cbpmesh=21.2] (0.7150,2.2780) -- (0.4791,2.3000) -- (0.4942,2.3271) -- (0.7376,2.3044) -- cycle;
\path[cbpmesh=20.5] (-0.4791,2.3000) -- (-0.7150,2.2780) -- (-0.7376,2.3044) -- (-0.4942,2.3271) -- cycle;
\path[cbpmesh=18.6] (-1.5516,2.2647) -- (-1.7857,2.2032) -- (-1.8261,2.2303) -- (-1.5867,2.2932) -- cycle;
\path[cbpmesh=19.1] (1.7857,2.2032) -- (1.5516,2.2647) -- (1.5867,2.2932) -- (1.8261,2.2303) -- cycle;
\path[cbpmesh=21.1] (2.4523,2.1612) -- (2.2194,2.2473) -- (2.2583,2.2618) -- (2.4952,2.1741) -- cycle;
\path[cbpmesh=19.2] (-2.2194,2.2473) -- (-2.4523,2.1612) -- (-2.4952,2.1741) -- (-2.2583,2.2618) -- cycle;
\path[cbpmesh=18.5] (-1.8647,2.2569) -- (-2.0977,2.1838) -- (-2.1395,2.2075) -- (-1.9018,2.2820) -- cycle;
\path[cbpmesh=19.0] (2.0977,2.1838) -- (1.8647,2.2569) -- (1.9018,2.2820) -- (2.1395,2.2075) -- cycle;
\path[cbpmesh=20.6] (1.2091,2.2325) -- (0.9763,2.2728) -- (1.0045,2.3008) -- (1.2439,2.2593) -- cycle;
\path[cbpmesh=19.7] (-0.9763,2.2728) -- (-1.2091,2.2325) -- (-1.2439,2.2593) -- (-1.0045,2.3008) -- cycle;
\path[cbpmesh=26.0] (2.8097,2.0854) -- (2.5806,2.1848) -- (2.6238,2.1808) -- (2.8568,2.0797) -- cycle;
\path[cbpmesh=21.1] (-2.5806,2.1848) -- (-2.8097,2.0854) -- (-2.8568,2.0797) -- (-2.6238,2.1808) -- cycle;
\path[cbpmesh=19.9] (1.4757,2.2091) -- (1.2439,2.2593) -- (1.2768,2.2879) -- (1.5147,2.2364) -- cycle;
\path[cbpmesh=19.1] (-1.2439,2.2593) -- (-1.4757,2.2091) -- (-1.5147,2.2364) -- (-1.2768,2.2879) -- cycle;
\path[cbpmesh=21.8] (0.2322,2.2890) -- (0.0000,2.2933) -- (0.0000,2.3177) -- (0.2403,2.3133) -- cycle;
\path[cbpmesh=21.6] (0.0000,2.2933) -- (-0.2322,2.2890) -- (-0.2403,2.3133) -- (0.0000,2.3177) -- cycle;
\path[cbpmesh=21.3] (0.9464,2.2474) -- (0.7150,2.2780) -- (0.7376,2.3044) -- (0.9763,2.2728) -- cycle;
\path[cbpmesh=20.4] (-0.7150,2.2780) -- (-0.9464,2.2474) -- (-0.9763,2.2728) -- (-0.7376,2.3044) -- cycle;
\path[cbpmesh=21.4] (-0.2322,2.2890) -- (-0.4630,2.2762) -- (-0.4791,2.3000) -- (-0.2403,2.3133) -- cycle;
\path[cbpmesh=21.9] (0.4630,2.2762) -- (0.2322,2.2890) -- (0.2403,2.3133) -- (0.4791,2.3000) -- cycle;
\path[cbpmesh=18.8] (-2.1799,2.2289) -- (-2.4086,2.1443) -- (-2.4523,2.1612) -- (-2.2194,2.2473) -- cycle;
\path[cbpmesh=20.1] (2.4086,2.1443) -- (2.1799,2.2289) -- (2.2194,2.2473) -- (2.4523,2.1612) -- cycle;
\path[cbpmesh=18.8] (-1.5147,2.2364) -- (-1.7432,2.1763) -- (-1.7857,2.2032) -- (-1.5516,2.2647) -- cycle;
\path[cbpmesh=19.4] (1.7432,2.1763) -- (1.5147,2.2364) -- (1.5516,2.2647) -- (1.7857,2.2032) -- cycle;
\path[cbpmesh=18.9] (2.0543,2.1587) -- (1.8261,2.2303) -- (1.8647,2.2569) -- (2.0977,2.1838) -- cycle;
\path[cbpmesh=18.5] (-1.8261,2.2303) -- (-2.0543,2.1587) -- (-2.0977,2.1838) -- (-1.8647,2.2569) -- cycle;
\path[cbpmesh=20.2] (-2.5379,2.1823) -- (-2.7632,2.0845) -- (-2.8097,2.0854) -- (-2.5806,2.1848) -- cycle;
\path[cbpmesh=24.1] (2.7632,2.0845) -- (2.5379,2.1823) -- (2.5806,2.1848) -- (2.8097,2.0854) -- cycle;
\path[cbpmesh=22.1] (0.6909,2.2549) -- (0.4630,2.2762) -- (0.4791,2.3000) -- (0.7150,2.2780) -- cycle;
\path[cbpmesh=21.2] (-0.4630,2.2762) -- (-0.6909,2.2549) -- (-0.7150,2.2780) -- (-0.4791,2.3000) -- cycle;
\path[cbpmesh=22.7] (-2.9047,2.0666) -- (-3.1237,1.9554) -- (-3.1767,1.9321) -- (-2.9540,2.0452) -- cycle;
\draw[draw=cbpInk!70,line width=.45pt] (-3.1767,1.9321) -- (-2.9540,2.0452);
\draw[draw=cbpInk!70,line width=.45pt] (-2.9540,2.0452) -- (-2.9047,2.0666);
\path[cbpmesh=30.0] (3.1237,1.9554) -- (2.9047,2.0666) -- (2.9540,2.0452) -- (3.1767,1.9321) -- cycle;
\draw[draw=cbpInk!70,line width=.45pt] (2.9047,2.0666) -- (2.9540,2.0452);
\draw[draw=cbpInk!70,line width=.45pt] (2.9540,2.0452) -- (3.1767,1.9321);
\path[cbpmesh=20.6] (1.4344,2.1837) -- (1.2091,2.2325) -- (1.2439,2.2593) -- (1.4757,2.2091) -- cycle;
\path[cbpmesh=19.5] (-1.2091,2.2325) -- (-1.4344,2.1837) -- (-1.4757,2.2091) -- (-1.2439,2.2593) -- cycle;
\path[cbpmesh=20.2] (-0.9464,2.2474) -- (-1.1720,2.2083) -- (-1.2091,2.2325) -- (-0.9763,2.2728) -- cycle;
\path[cbpmesh=21.4] (1.1720,2.2083) -- (0.9464,2.2474) -- (0.9763,2.2728) -- (1.2091,2.2325) -- cycle;
\path[cbpmesh=18.6] (-2.1395,2.2075) -- (-2.3639,2.1244) -- (-2.4086,2.1443) -- (-2.1799,2.2289) -- cycle;
\path[cbpmesh=19.4] (2.3639,2.1244) -- (2.1395,2.2075) -- (2.1799,2.2289) -- (2.4086,2.1443) -- cycle;
\path[cbpmesh=19.1] (2.0088,2.1332) -- (1.7857,2.2032) -- (1.8261,2.2303) -- (2.0543,2.1587) -- cycle;
\path[cbpmesh=18.6] (-1.7857,2.2032) -- (-2.0088,2.1332) -- (-2.0543,2.1587) -- (-1.8261,2.2303) -- cycle;
\path[cbpmesh=19.5] (-2.4952,2.1741) -- (-2.7168,2.0780) -- (-2.7632,2.0845) -- (-2.5379,2.1823) -- cycle;
\path[cbpmesh=22.4] (2.7168,2.0780) -- (2.4952,2.1741) -- (2.5379,2.1823) -- (2.7632,2.0845) -- cycle;
\path[cbpmesh=19.0] (-1.4757,2.2091) -- (-1.6984,2.1506) -- (-1.7432,2.1763) -- (-1.5147,2.2364) -- cycle;
\path[cbpmesh=20.0] (1.6984,2.1506) -- (1.4757,2.2091) -- (1.5147,2.2364) -- (1.7432,2.1763) -- cycle;
\path[cbpmesh=21.0] (-0.6909,2.2549) -- (-0.9146,2.2253) -- (-0.9464,2.2474) -- (-0.7150,2.2780) -- cycle;
\path[cbpmesh=22.2] (0.9146,2.2253) -- (0.6909,2.2549) -- (0.7150,2.2780) -- (0.9464,2.2474) -- cycle;
\path[cbpmesh=22.7] (0.2236,2.2688) -- (0.0000,2.2729) -- (0.0000,2.2933) -- (0.2322,2.2890) -- cycle;
\path[cbpmesh=22.5] (0.0000,2.2729) -- (-0.2236,2.2688) -- (-0.2322,2.2890) -- (0.0000,2.2933) -- cycle;
\path[cbpmesh=19.0] (2.3178,2.1023) -- (2.0977,2.1838) -- (2.1395,2.2075) -- (2.3639,2.1244) -- cycle;
\path[cbpmesh=18.5] (-2.0977,2.1838) -- (-2.3178,2.1023) -- (-2.3639,2.1244) -- (-2.1395,2.2075) -- cycle;
\path[cbpmesh=21.7] (-2.8568,2.0797) -- (-3.0721,1.9704) -- (-3.1237,1.9554) -- (-2.9047,2.0666) -- cycle;
\path[cbpmesh=28.0] (3.0721,1.9704) -- (2.8568,2.0797) -- (2.9047,2.0666) -- (3.1237,1.9554) -- cycle;
\path[cbpmesh=22.8] (0.4458,2.2565) -- (0.2236,2.2688) -- (0.2322,2.2890) -- (0.4630,2.2762) -- cycle;
\path[cbpmesh=22.2] (-0.2236,2.2688) -- (-0.4458,2.2565) -- (-0.4630,2.2762) -- (-0.2322,2.2890) -- cycle;
\path[cbpmesh=19.0] (-2.4523,2.1612) -- (-2.6700,2.0667) -- (-2.7168,2.0780) -- (-2.4952,2.1741) -- cycle;
\path[cbpmesh=21.1] (2.6700,2.0667) -- (2.4523,2.1612) -- (2.4952,2.1741) -- (2.7168,2.0780) -- cycle;
\path[cbpmesh=19.4] (1.9611,2.1080) -- (1.7432,2.1763) -- (1.7857,2.2032) -- (2.0088,2.1332) -- cycle;
\path[cbpmesh=18.7] (-1.7432,2.1763) -- (-1.9611,2.1080) -- (-2.0088,2.1332) -- (-1.7857,2.2032) -- cycle;
\path[cbpmesh=21.4] (1.3904,2.1610) -- (1.1720,2.2083) -- (1.2091,2.2325) -- (1.4344,2.1837) -- cycle;
\path[cbpmesh=20.0] (-1.1720,2.2083) -- (-1.3904,2.1610) -- (-1.4344,2.1837) -- (-1.2091,2.2325) -- cycle;
\path[cbpmesh=23.0] (0.6652,2.2360) -- (0.4458,2.2565) -- (0.4630,2.2762) -- (0.6909,2.2549) -- cycle;
\path[cbpmesh=22.0] (-0.4458,2.2565) -- (-0.6652,2.2360) -- (-0.6909,2.2549) -- (-0.4630,2.2762) -- cycle;
\path[cbpmesh=22.3] (1.1326,2.1876) -- (0.9146,2.2253) -- (0.9464,2.2474) -- (1.1720,2.2083) -- cycle;
\path[cbpmesh=20.8] (-0.9146,2.2253) -- (-1.1326,2.1876) -- (-1.1720,2.2083) -- (-0.9464,2.2474) -- cycle;
\path[cbpmesh=19.4] (-1.4344,2.1837) -- (-1.6508,2.1268) -- (-1.6984,2.1506) -- (-1.4757,2.2091) -- cycle;
\path[cbpmesh=20.6] (1.6508,2.1268) -- (1.4344,2.1837) -- (1.4757,2.2091) -- (1.6984,2.1506) -- cycle;
\path[cbpmesh=18.9] (2.2698,2.0790) -- (2.0543,2.1587) -- (2.0977,2.1838) -- (2.3178,2.1023) -- cycle;
\path[cbpmesh=18.5] (-2.0543,2.1587) -- (-2.2698,2.0790) -- (-2.3178,2.1023) -- (-2.0977,2.1838) -- cycle;
\path[cbpmesh=20.7] (-2.8097,2.0854) -- (-3.0215,1.9778) -- (-3.0721,1.9704) -- (-2.8568,2.0797) -- cycle;
\path[cbpmesh=26.0] (3.0215,1.9778) -- (2.8097,2.0854) -- (2.8568,2.0797) -- (3.0721,1.9704) -- cycle;
\path[cbpmesh=20.1] (2.6225,2.0515) -- (2.4086,2.1443) -- (2.4523,2.1612) -- (2.6700,2.0667) -- cycle;
\path[cbpmesh=18.7] (-2.4086,2.1443) -- (-2.6225,2.0515) -- (-2.6700,2.0667) -- (-2.4523,2.1612) -- cycle;
\path[cbpmesh=21.7] (-0.6652,2.2360) -- (-0.8806,2.2075) -- (-0.9146,2.2253) -- (-0.6909,2.2549) -- cycle;
\path[cbpmesh=23.1] (0.8806,2.2075) -- (0.6652,2.2360) -- (0.6909,2.2549) -- (0.9146,2.2253) -- cycle;
\path[cbpmesh=18.9] (-1.6984,2.1506) -- (-1.9106,2.0840) -- (-1.9611,2.1080) -- (-1.7432,2.1763) -- cycle;
\path[cbpmesh=20.0] (1.9106,2.0840) -- (1.6984,2.1506) -- (1.7432,2.1763) -- (1.9611,2.1080) -- cycle;
\path[cbpmesh=23.1] (-3.1767,1.9321) -- (-3.3798,1.8112) -- (-3.4381,1.7765) -- (-3.2316,1.8996) -- cycle;
\draw[draw=cbpInk!70,line width=.45pt] (-3.4381,1.7765) -- (-3.2316,1.8996);
\draw[draw=cbpInk!70,line width=.45pt] (-3.2316,1.8996) -- (-3.1767,1.9321);
\path[cbpmesh=31.7] (3.3798,1.8112) -- (3.1767,1.9321) -- (3.2316,1.8996) -- (3.4381,1.7765) -- cycle;
\draw[draw=cbpInk!70,line width=.45pt] (3.1767,1.9321) -- (3.2316,1.8996);
\draw[draw=cbpInk!70,line width=.45pt] (3.2316,1.8996) -- (3.4381,1.7765);
\path[cbpmesh=23.6] (0.2144,2.2534) -- (0.0000,2.2574) -- (0.0000,2.2729) -- (0.2236,2.2688) -- cycle;
\path[cbpmesh=23.3] (0.0000,2.2574) -- (-0.2144,2.2534) -- (-0.2236,2.2688) -- (0.0000,2.2729) -- cycle;
\path[cbpmesh=19.1] (2.2196,2.0552) -- (2.0088,2.1332) -- (2.0543,2.1587) -- (2.2698,2.0790) -- cycle;
\path[cbpmesh=18.5] (-2.0088,2.1332) -- (-2.2196,2.0552) -- (-2.2698,2.0790) -- (-2.0543,2.1587) -- cycle;
\path[cbpmesh=19.9] (-2.7632,2.0845) -- (-2.9715,1.9787) -- (-3.0215,1.9778) -- (-2.8097,2.0854) -- cycle;
\path[cbpmesh=24.1] (2.9715,1.9787) -- (2.7632,2.0845) -- (2.8097,2.0854) -- (3.0215,1.9778) -- cycle;
\path[cbpmesh=19.4] (2.5738,2.0333) -- (2.3639,2.1244) -- (2.4086,2.1443) -- (2.6225,2.0515) -- cycle;
\path[cbpmesh=18.5] (-2.3639,2.1244) -- (-2.5738,2.0333) -- (-2.6225,2.0515) -- (-2.4086,2.1443) -- cycle;
\path[cbpmesh=20.6] (-1.1326,2.1876) -- (-1.3436,2.1419) -- (-1.3904,2.1610) -- (-1.1720,2.2083) -- cycle;
\path[cbpmesh=22.3] (1.3436,2.1419) -- (1.1326,2.1876) -- (1.1720,2.2083) -- (1.3904,2.1610) -- cycle;
\path[cbpmesh=23.8] (0.4274,2.2416) -- (0.2144,2.2534) -- (0.2236,2.2688) -- (0.4458,2.2565) -- cycle;
\path[cbpmesh=23.1] (-0.2144,2.2534) -- (-0.4274,2.2416) -- (-0.4458,2.2565) -- (-0.2236,2.2688) -- cycle;
\path[cbpmesh=21.5] (1.6003,2.1059) -- (1.3904,2.1610) -- (1.4344,2.1837) -- (1.6508,2.1268) -- cycle;
\path[cbpmesh=19.8] (-1.3904,2.1610) -- (-1.6003,2.1059) -- (-1.6508,2.1268) -- (-1.4344,2.1837) -- cycle;
\path[cbpmesh=38] (0.1401,0.0895) -- (0.0000,0.0921) -- (0.0000,0.0925) -- (0.1328,0.0900) -- cycle;
\draw[draw=cbpInk!70,line width=.45pt] (0.1401,0.0895) -- (0.0000,0.0921);
\path[cbpmesh=38] (0.0000,0.0921) -- (-0.1401,0.0895) -- (-0.1328,0.0900) -- (0.0000,0.0925) -- cycle;
\draw[draw=cbpInk!70,line width=.45pt] (0.0000,0.0921) -- (-0.1401,0.0895);
\path[cbpmesh=23.2] (1.0905,2.1711) -- (0.8806,2.2075) -- (0.9146,2.2253) -- (1.1326,2.1876) -- cycle;
\path[cbpmesh=21.5] (-0.8806,2.2075) -- (-1.0905,2.1711) -- (-1.1326,2.1876) -- (-0.9146,2.2253) -- cycle;
\path[cbpmesh=38] (-0.1401,0.0895) -- (-0.2793,0.0818) -- (-0.2647,0.0827) -- (-0.1328,0.0900) -- cycle;
\draw[draw=cbpInk!70,line width=.45pt] (-0.1401,0.0895) -- (-0.2793,0.0818);
\draw[draw=cbpInk!70,line width=.45pt] (-0.2793,0.0818) -- (-0.2647,0.0827);
\path[cbpmesh=38] (0.2793,0.0818) -- (0.1401,0.0895) -- (0.1328,0.0900) -- (0.2647,0.0827) -- cycle;
\draw[draw=cbpInk!70,line width=.45pt] (0.2793,0.0818) -- (0.1401,0.0895);
\draw[draw=cbpInk!70,line width=.45pt] (0.2647,0.0827) -- (0.2793,0.0818);
\path[cbpmesh=24.0] (0.6378,2.2219) -- (0.4274,2.2416) -- (0.4458,2.2565) -- (0.6652,2.2360) -- cycle;
\path[cbpmesh=22.8] (-0.4274,2.2416) -- (-0.6378,2.2219) -- (-0.6652,2.2360) -- (-0.4458,2.2565) -- cycle;
\path[cbpmesh=22.1] (-3.1237,1.9554) -- (-3.3234,1.8365) -- (-3.3798,1.8112) -- (-3.1767,1.9321) -- cycle;
\path[cbpmesh=29.8] (3.3234,1.8365) -- (3.1237,1.9554) -- (3.1767,1.9321) -- (3.3798,1.8112) -- cycle;
\path[cbpmesh=19.2] (-2.7168,2.0780) -- (-2.9215,1.9740) -- (-2.9715,1.9787) -- (-2.7632,2.0845) -- cycle;
\path[cbpmesh=22.4] (2.9215,1.9740) -- (2.7168,2.0780) -- (2.7632,2.0845) -- (2.9715,1.9787) -- cycle;
\path[cbpmesh=20.7] (1.8571,2.0621) -- (1.6508,2.1268) -- (1.6984,2.1506) -- (1.9106,2.0840) -- cycle;
\path[cbpmesh=19.2] (-1.6508,2.1268) -- (-1.8571,2.0621) -- (-1.9106,2.0840) -- (-1.6984,2.1506) -- cycle;
\path[cbpmesh=18.4] (-2.3178,2.1023) -- (-2.5235,2.0130) -- (-2.5738,2.0333) -- (-2.3639,2.1244) -- cycle;
\path[cbpmesh=19.0] (2.5235,2.0130) -- (2.3178,2.1023) -- (2.3639,2.1244) -- (2.5738,2.0333) -- cycle;
\path[cbpmesh=18.6] (-1.9611,2.1080) -- (-2.1668,2.0318) -- (-2.2196,2.0552) -- (-2.0088,2.1332) -- cycle;
\path[cbpmesh=19.4] (2.1668,2.0318) -- (1.9611,2.1080) -- (2.0088,2.1332) -- (2.2196,2.0552) -- cycle;
\path[cbpmesh=24.1] (0.8442,2.1946) -- (0.6378,2.2219) -- (0.6652,2.2360) -- (0.8806,2.2075) -- cycle;
\path[cbpmesh=22.5] (-0.6378,2.2219) -- (-0.8442,2.1946) -- (-0.8806,2.2075) -- (-0.6652,2.2360) -- cycle;
\path[cbpmesh=21.0] (2.8712,1.9645) -- (2.6700,2.0667) -- (2.7168,2.0780) -- (2.9215,1.9740) -- cycle;
\path[cbpmesh=18.8] (-2.6700,2.0667) -- (-2.8712,1.9645) -- (-2.9215,1.9740) -- (-2.7168,2.0780) -- cycle;
\path[cbpmesh=18.4] (-2.2698,2.0790) -- (-2.4713,1.9915) -- (-2.5235,2.0130) -- (-2.3178,2.1023) -- cycle;
\path[cbpmesh=18.9] (2.4713,1.9915) -- (2.2698,2.0790) -- (2.3178,2.1023) -- (2.5235,2.0130) -- cycle;
\path[cbpmesh=20.3] (-1.3436,2.1419) -- (-1.5464,2.0886) -- (-1.6003,2.1059) -- (-1.3904,2.1610) -- cycle;
\path[cbpmesh=22.4] (1.5464,2.0886) -- (1.3436,2.1419) -- (1.3904,2.1610) -- (1.6003,2.1059) -- cycle;
\path[cbpmesh=21.1] (-3.0721,1.9704) -- (-3.2685,1.8534) -- (-3.3234,1.8365) -- (-3.1237,1.9554) -- cycle;
\path[cbpmesh=27.9] (3.2685,1.8534) -- (3.0721,1.9704) -- (3.1237,1.9554) -- (3.3234,1.8365) -- cycle;
\path[cbpmesh=23.3] (1.2937,2.1271) -- (1.0905,2.1711) -- (1.1326,2.1876) -- (1.3436,2.1419) -- cycle;
\path[cbpmesh=21.2] (-1.0905,2.1711) -- (-1.2937,2.1271) -- (-1.3436,2.1419) -- (-1.1326,2.1876) -- cycle;
\path[cbpmesh=20.0] (2.1110,2.0098) -- (1.9106,2.0840) -- (1.9611,2.1080) -- (2.1668,2.0318) -- cycle;
\path[cbpmesh=18.8] (-1.9106,2.0840) -- (-2.1110,2.0098) -- (-2.1668,2.0318) -- (-1.9611,2.1080) -- cycle;
\path[cbpmesh=22.8] (0.0000,2.2371) -- (-0.2046,2.2333) -- (-0.2144,2.2534) -- (0.0000,2.2574) -- cycle;
\path[cbpmesh=23.0] (0.2046,2.2333) -- (0.0000,2.2371) -- (0.0000,2.2574) -- (0.2144,2.2534) -- cycle;
\path[cbpmesh=21.5] (1.8002,2.0431) -- (1.6003,2.1059) -- (1.6508,2.1268) -- (1.8571,2.0621) -- cycle;
\path[cbpmesh=19.6] (-1.6003,2.1059) -- (-1.8002,2.0431) -- (-1.8571,2.0621) -- (-1.6508,2.1268) -- cycle;
\path[cbpmesh=38] (0.1328,0.0900) -- (0.0000,0.0925) -- (0.0000,0.1112) -- (0.1261,0.1088) -- cycle;
\path[cbpmesh=38] (0.0000,0.0925) -- (-0.1328,0.0900) -- (-0.1261,0.1088) -- (0.0000,0.1112) -- cycle;
\path[cbpmesh=22.6] (-0.2046,2.2333) -- (-0.4079,2.2220) -- (-0.4274,2.2416) -- (-0.2144,2.2534) -- cycle;
\path[cbpmesh=23.2] (0.4079,2.2220) -- (0.2046,2.2333) -- (0.2144,2.2534) -- (0.4274,2.2416) -- cycle;
\path[cbpmesh=20.0] (2.8201,1.9511) -- (2.6225,2.0515) -- (2.6700,2.0667) -- (2.8712,1.9645) -- cycle;
\path[cbpmesh=18.6] (-2.6225,2.0515) -- (-2.8201,1.9511) -- (-2.8712,1.9645) -- (-2.6700,2.0667) -- cycle;
\path[cbpmesh=38] (0.2647,0.0827) -- (0.1328,0.0900) -- (0.1261,0.1088) -- (0.2514,0.1019) -- cycle;
\path[cbpmesh=38] (-0.1328,0.0900) -- (-0.2647,0.0827) -- (-0.2514,0.1019) -- (-0.1261,0.1088) -- cycle;
\path[cbpmesh=24.3] (1.0455,2.1597) -- (0.8442,2.1946) -- (0.8806,2.2075) -- (1.0905,2.1711) -- cycle;
\path[cbpmesh=22.2] (-0.8442,2.1946) -- (-1.0455,2.1597) -- (-1.0905,2.1711) -- (-0.8806,2.2075) -- cycle;
\path[cbpmesh=18.4] (-2.2196,2.0552) -- (-2.4166,1.9697) -- (-2.4713,1.9915) -- (-2.2698,2.0790) -- cycle;
\path[cbpmesh=19.1] (2.4166,1.9697) -- (2.2196,2.0552) -- (2.2698,2.0790) -- (2.4713,1.9915) -- cycle;
\path[cbpmesh=20.2] (-3.0215,1.9778) -- (-3.2147,1.8628) -- (-3.2685,1.8534) -- (-3.0721,1.9704) -- cycle;
\path[cbpmesh=25.9] (3.2147,1.8628) -- (3.0215,1.9778) -- (3.0721,1.9704) -- (3.2685,1.8534) -- cycle;
\path[cbpmesh=38] (-0.2647,0.0827) -- (-0.3951,0.0705) -- (-0.3751,0.0903) -- (-0.2514,0.1019) -- cycle;
\draw[draw=cbpInk!70,line width=.45pt] (-0.2647,0.0827) -- (-0.3951,0.0705);
\path[cbpmesh=38] (0.3951,0.0705) -- (0.2647,0.0827) -- (0.2514,0.1019) -- (0.3751,0.0903) -- cycle;
\draw[draw=cbpInk!70,line width=.45pt] (0.3951,0.0705) -- (0.2647,0.0827);
\path[cbpmesh=22.3] (-0.4079,2.2220) -- (-0.6087,2.2033) -- (-0.6378,2.2219) -- (-0.4274,2.2416) -- cycle;
\path[cbpmesh=23.4] (0.6087,2.2033) -- (0.4079,2.2220) -- (0.4274,2.2416) -- (0.6378,2.2219) -- cycle;
\path[cbpmesh=20.7] (2.0519,1.9900) -- (1.8571,2.0621) -- (1.9106,2.0840) -- (2.1110,2.0098) -- cycle;
\path[cbpmesh=19.0] (-1.8571,2.0621) -- (-2.0519,1.9900) -- (-2.1110,2.0098) -- (-1.9106,2.0840) -- cycle;
\path[cbpmesh=18.4] (-2.5738,2.0333) -- (-2.7678,1.9348) -- (-2.8201,1.9511) -- (-2.6225,2.0515) -- cycle;
\path[cbpmesh=19.4] (2.7678,1.9348) -- (2.5738,2.0333) -- (2.6225,2.0515) -- (2.8201,1.9511) -- cycle;
\path[cbpmesh=38] (0.5230,0.0536) -- (0.3951,0.0705) -- (0.3751,0.0903) -- (0.4965,0.0742) -- cycle;
\draw[draw=cbpInk!70,line width=.45pt] (0.5230,0.0536) -- (0.3951,0.0705);
\path[cbpmesh=38] (-0.3951,0.0705) -- (-0.5230,0.0536) -- (-0.4965,0.0742) -- (-0.3751,0.0903) -- cycle;
\draw[draw=cbpInk!70,line width=.45pt] (-0.3951,0.0705) -- (-0.5230,0.0536);
\path[cbpmesh=20.9] (-1.2937,2.1271) -- (-1.4889,2.0758) -- (-1.5464,2.0886) -- (-1.3436,2.1419) -- cycle;
\path[cbpmesh=23.4] (1.4889,2.0758) -- (1.2937,2.1271) -- (1.3436,2.1419) -- (1.5464,2.0886) -- cycle;
\path[cbpmesh=18.5] (-2.1668,2.0318) -- (-2.3591,1.9484) -- (-2.4166,1.9697) -- (-2.2196,2.0552) -- cycle;
\path[cbpmesh=19.4] (2.3591,1.9484) -- (2.1668,2.0318) -- (2.2196,2.0552) -- (2.4166,1.9697) -- cycle;
\path[cbpmesh=19.5] (-2.9715,1.9787) -- (-3.1615,1.8656) -- (-3.2147,1.8628) -- (-3.0215,1.9778) -- cycle;
\path[cbpmesh=24.0] (3.1615,1.8656) -- (2.9715,1.9787) -- (3.0215,1.9778) -- (3.2147,1.8628) -- cycle;
\path[cbpmesh=23.6] (0.8057,2.1772) -- (0.6087,2.2033) -- (0.6378,2.2219) -- (0.8442,2.1946) -- cycle;
\path[cbpmesh=22.1] (-0.6087,2.2033) -- (-0.8057,2.1772) -- (-0.8442,2.1946) -- (-0.6378,2.2219) -- cycle;
\path[cbpmesh=22.4] (1.7396,2.0280) -- (1.5464,2.0886) -- (1.6003,2.1059) -- (1.8002,2.0431) -- cycle;
\path[cbpmesh=20.0] (-1.5464,2.0886) -- (-1.7396,2.0280) -- (-1.8002,2.0431) -- (-1.6003,2.1059) -- cycle;
\path[cbpmesh=31.5] (3.5620,1.6832) -- (3.3798,1.8112) -- (3.4381,1.7765) -- (3.6235,1.6463) -- cycle;
\draw[draw=cbpInk!70,line width=.45pt] (3.4381,1.7765) -- (3.6235,1.6463);
\path[cbpmesh=22.4] (-3.3798,1.8112) -- (-3.5620,1.6832) -- (-3.6235,1.6463) -- (-3.4381,1.7765) -- cycle;
\draw[draw=cbpInk!70,line width=.45pt] (-3.6235,1.6463) -- (-3.4381,1.7765);
\path[cbpmesh=24.4] (1.2403,2.1176) -- (1.0455,2.1597) -- (1.0905,2.1711) -- (1.2937,2.1271) -- cycle;
\path[cbpmesh=21.8] (-1.0455,2.1597) -- (-1.2403,2.1176) -- (-1.2937,2.1271) -- (-1.0905,2.1711) -- cycle;
\path[cbpmesh=38] (-0.5230,0.0536) -- (-0.6476,0.0320) -- (-0.6149,0.0537) -- (-0.4965,0.0742) -- cycle;
\draw[draw=cbpInk!70,line width=.45pt] (-0.5230,0.0536) -- (-0.6476,0.0320);
\path[cbpmesh=38] (0.6476,0.0320) -- (0.5230,0.0536) -- (0.4965,0.0742) -- (0.6149,0.0537) -- cycle;
\draw[draw=cbpInk!70,line width=.45pt] (0.6476,0.0320) -- (0.5230,0.0536);
\path[cbpmesh=18.4] (-2.5235,2.0130) -- (-2.7137,1.9164) -- (-2.7678,1.9348) -- (-2.5738,2.0333) -- cycle;
\path[cbpmesh=19.0] (2.7137,1.9164) -- (2.5235,2.0130) -- (2.5738,2.0333) -- (2.7678,1.9348) -- cycle;
\path[cbpmesh=38] (0.1261,0.1088) -- (0.0000,0.1112) -- (0.0000,0.1503) -- (0.1200,0.1481) -- cycle;
\path[cbpmesh=38] (0.0000,0.1112) -- (-0.1261,0.1088) -- (-0.1200,0.1481) -- (0.0000,0.1503) -- cycle;
\path[cbpmesh=21.5] (0.1950,2.2025) -- (0.0000,2.2061) -- (0.0000,2.2371) -- (0.2046,2.2333) -- cycle;
\path[cbpmesh=21.4] (0.0000,2.2061) -- (-0.1950,2.2025) -- (-0.2046,2.2333) -- (0.0000,2.2371) -- cycle;
\path[cbpmesh=38] (0.2514,0.1019) -- (0.1261,0.1088) -- (0.1200,0.1481) -- (0.2392,0.1415) -- cycle;
\path[cbpmesh=38] (-0.1261,0.1088) -- (-0.2514,0.1019) -- (-0.2392,0.1415) -- (-0.1200,0.1481) -- cycle;
\path[cbpmesh=19.0] (-2.9215,1.9740) -- (-3.1083,1.8628) -- (-3.1615,1.8656) -- (-2.9715,1.9787) -- cycle;
\path[cbpmesh=22.3] (3.1083,1.8628) -- (2.9215,1.9740) -- (2.9715,1.9787) -- (3.1615,1.8656) -- cycle;
\path[cbpmesh=19.4] (-1.8002,2.0431) -- (-1.9891,1.9732) -- (-2.0519,1.9900) -- (-1.8571,2.0621) -- cycle;
\path[cbpmesh=21.5] (1.9891,1.9732) -- (1.8002,2.0431) -- (1.8571,2.0621) -- (2.0519,1.9900) -- cycle;
\path[cbpmesh=21.6] (0.3888,2.1917) -- (0.1950,2.2025) -- (0.2046,2.2333) -- (0.4079,2.2220) -- cycle;
\path[cbpmesh=21.2] (-0.1950,2.2025) -- (-0.3888,2.1917) -- (-0.4079,2.2220) -- (-0.2046,2.2333) -- cycle;
\path[cbpmesh=18.7] (-2.1110,2.0098) -- (-2.2984,1.9285) -- (-2.3591,1.9484) -- (-2.1668,2.0318) -- cycle;
\path[cbpmesh=20.0] (2.2984,1.9285) -- (2.1110,2.0098) -- (2.1668,2.0318) -- (2.3591,1.9484) -- cycle;
\path[cbpmesh=23.7] (0.9978,2.1439) -- (0.8057,2.1772) -- (0.8442,2.1946) -- (1.0455,2.1597) -- cycle;
\path[cbpmesh=21.8] (-0.8057,2.1772) -- (-0.9978,2.1439) -- (-1.0455,2.1597) -- (-0.8442,2.1946) -- cycle;
\path[cbpmesh=38] (0.7683,0.0059) -- (0.6476,0.0320) -- (0.6149,0.0537) -- (0.7295,0.0289) -- cycle;
\draw[draw=cbpInk!70,line width=.45pt] (0.7683,0.0059) -- (0.6476,0.0320);
\path[cbpmesh=38] (-0.6476,0.0320) -- (-0.7683,0.0059) -- (-0.7295,0.0289) -- (-0.6149,0.0537) -- cycle;
\draw[draw=cbpInk!70,line width=.45pt] (-0.6476,0.0320) -- (-0.7683,0.0059);
\path[cbpmesh=38] (-0.2514,0.1019) -- (-0.3751,0.0903) -- (-0.3570,0.1305) -- (-0.2392,0.1415) -- cycle;
\path[cbpmesh=38] (0.3751,0.0903) -- (0.2514,0.1019) -- (0.2392,0.1415) -- (0.3570,0.1305) -- cycle;
\path[cbpmesh=18.4] (-2.4713,1.9915) -- (-2.6576,1.8969) -- (-2.7137,1.9164) -- (-2.5235,2.0130) -- cycle;
\path[cbpmesh=18.9] (2.6576,1.8969) -- (2.4713,1.9915) -- (2.5235,2.0130) -- (2.7137,1.9164) -- cycle;
\path[cbpmesh=21.5] (-3.3234,1.8365) -- (-3.5026,1.7106) -- (-3.5620,1.6832) -- (-3.3798,1.8112) -- cycle;
\path[cbpmesh=29.7] (3.5026,1.7106) -- (3.3234,1.8365) -- (3.3798,1.8112) -- (3.5620,1.6832) -- cycle;
\path[cbpmesh=21.8] (0.5802,2.1739) -- (0.3888,2.1917) -- (0.4079,2.2220) -- (0.6087,2.2033) -- cycle;
\path[cbpmesh=21.0] (-0.3888,2.1917) -- (-0.5802,2.1739) -- (-0.6087,2.2033) -- (-0.4079,2.2220) -- cycle;
\path[cbpmesh=38] (-0.3751,0.0903) -- (-0.4965,0.0742) -- (-0.4726,0.1152) -- (-0.3570,0.1305) -- cycle;
\path[cbpmesh=38] (0.4965,0.0742) -- (0.3751,0.0903) -- (0.3570,0.1305) -- (0.4726,0.1152) -- cycle;
\path[cbpmesh=18.6] (-2.8712,1.9645) -- (-3.0548,1.8552) -- (-3.1083,1.8628) -- (-2.9215,1.9740) -- cycle;
\path[cbpmesh=21.0] (3.0548,1.8552) -- (2.8712,1.9645) -- (2.9215,1.9740) -- (3.1083,1.8628) -- cycle;
\path[cbpmesh=23.4] (1.6749,2.0174) -- (1.4889,2.0758) -- (1.5464,2.0886) -- (1.7396,2.0280) -- cycle;
\path[cbpmesh=20.6] (-1.4889,2.0758) -- (-1.6749,2.0174) -- (-1.7396,2.0280) -- (-1.5464,2.0886) -- cycle;
\path[cbpmesh=21.5] (-1.2403,2.1176) -- (-1.4275,2.0684) -- (-1.4889,2.0758) -- (-1.2937,2.1271) -- cycle;
\path[cbpmesh=24.4] (1.4275,2.0684) -- (1.2403,2.1176) -- (1.2937,2.1271) -- (1.4889,2.0758) -- cycle;
\path[cbpmesh=38] (0.8843,-0.0246) -- (0.7683,0.0059) -- (0.7295,0.0289) -- (0.8395,-0.0000) -- cycle;
\draw[draw=cbpInk!70,line width=.45pt] (0.8843,-0.0246) -- (0.7683,0.0059);
\path[cbpmesh=38] (-0.7683,0.0059) -- (-0.8843,-0.0246) -- (-0.8395,-0.0000) -- (-0.7295,0.0289) -- cycle;
\draw[draw=cbpInk!70,line width=.45pt] (-0.7683,0.0059) -- (-0.8843,-0.0246);
\path[cbpmesh=18.4] (-2.4166,1.9697) -- (-2.5988,1.8771) -- (-2.6576,1.8969) -- (-2.4713,1.9915) -- cycle;
\path[cbpmesh=19.1] (2.5988,1.8771) -- (2.4166,1.9697) -- (2.4713,1.9915) -- (2.6576,1.8969) -- cycle;
\path[cbpmesh=20.7] (2.2341,1.9109) -- (2.0519,1.9900) -- (2.1110,2.0098) -- (2.2984,1.9285) -- cycle;
\path[cbpmesh=18.9] (-2.0519,1.9900) -- (-2.2341,1.9109) -- (-2.2984,1.9285) -- (-2.1110,2.0098) -- cycle;
\path[cbpmesh=38] (-0.4965,0.0742) -- (-0.6149,0.0537) -- (-0.5853,0.0957) -- (-0.4726,0.1152) -- cycle;
\path[cbpmesh=38] (0.6149,0.0537) -- (0.4965,0.0742) -- (0.4726,0.1152) -- (0.5853,0.0957) -- cycle;
\path[cbpmesh=20.8] (-0.5802,2.1739) -- (-0.7680,2.1490) -- (-0.8057,2.1772) -- (-0.6087,2.2033) -- cycle;
\path[cbpmesh=21.9] (0.7680,2.1490) -- (0.5802,2.1739) -- (0.6087,2.2033) -- (0.8057,2.1772) -- cycle;
\path[cbpmesh=20.6] (-3.2685,1.8534) -- (-3.4447,1.7296) -- (-3.5026,1.7106) -- (-3.3234,1.8365) -- cycle;
\path[cbpmesh=27.7] (3.4447,1.7296) -- (3.2685,1.8534) -- (3.3234,1.8365) -- (3.5026,1.7106) -- cycle;
\path[cbpmesh=19.8] (-1.7396,2.0280) -- (-1.9221,1.9604) -- (-1.9891,1.9732) -- (-1.8002,2.0431) -- cycle;
\path[cbpmesh=22.4] (1.9221,1.9604) -- (1.7396,2.0280) -- (1.8002,2.0431) -- (1.9891,1.9732) -- cycle;
\path[cbpmesh=23.8] (1.1837,2.1037) -- (0.9978,2.1439) -- (1.0455,2.1597) -- (1.2403,2.1176) -- cycle;
\path[cbpmesh=21.5] (-0.9978,2.1439) -- (-1.1837,2.1037) -- (-1.2403,2.1176) -- (-1.0455,2.1597) -- cycle;
\path[cbpmesh=18.4] (-2.8201,1.9511) -- (-3.0004,1.8437) -- (-3.0548,1.8552) -- (-2.8712,1.9645) -- cycle;
\path[cbpmesh=20.0] (3.0004,1.8437) -- (2.8201,1.9511) -- (2.8712,1.9645) -- (3.0548,1.8552) -- cycle;
\path[cbpmesh=38] (0.1200,0.1481) -- (0.0000,0.1503) -- (0.0000,0.2122) -- (0.1146,0.2101) -- cycle;
\path[cbpmesh=38] (0.0000,0.1503) -- (-0.1200,0.1481) -- (-0.1146,0.2101) -- (0.0000,0.2122) -- cycle;
\path[cbpmesh=38] (0.2392,0.1415) -- (0.1200,0.1481) -- (0.1146,0.2101) -- (0.2285,0.2037) -- cycle;
\path[cbpmesh=38] (-0.1200,0.1481) -- (-0.2392,0.1415) -- (-0.2285,0.2037) -- (-0.1146,0.2101) -- cycle;
\path[cbpmesh=20.0] (0.1856,2.1615) -- (0.0000,2.1650) -- (0.0000,2.2061) -- (0.1950,2.2025) -- cycle;
\path[cbpmesh=20.0] (0.0000,2.1650) -- (-0.1856,2.1615) -- (-0.1950,2.2025) -- (0.0000,2.2061) -- cycle;
\path[cbpmesh=38] (0.7295,0.0289) -- (0.6149,0.0537) -- (0.5853,0.0957) -- (0.6943,0.0721) -- cycle;
\path[cbpmesh=38] (-0.6149,0.0537) -- (-0.7295,0.0289) -- (-0.6943,0.0721) -- (-0.5853,0.0957) -- cycle;
\path[cbpmesh=18.4] (-2.3591,1.9484) -- (-2.5370,1.8580) -- (-2.5988,1.8771) -- (-2.4166,1.9697) -- cycle;
\path[cbpmesh=19.4] (2.5370,1.8580) -- (2.3591,1.9484) -- (2.4166,1.9697) -- (2.5988,1.8771) -- cycle;
\path[cbpmesh=38] (0.3570,0.1305) -- (0.2392,0.1415) -- (0.2285,0.2037) -- (0.3411,0.1932) -- cycle;
\path[cbpmesh=38] (-0.2392,0.1415) -- (-0.3570,0.1305) -- (-0.3411,0.1932) -- (-0.2285,0.2037) -- cycle;
\path[cbpmesh=19.4] (2.9447,1.8294) -- (2.7678,1.9348) -- (2.8201,1.9511) -- (3.0004,1.8437) -- cycle;
\path[cbpmesh=18.4] (-2.7678,1.9348) -- (-2.9447,1.8294) -- (-3.0004,1.8437) -- (-2.8201,1.9511) -- cycle;
\path[cbpmesh=20.1] (0.3701,2.1513) -- (0.1856,2.1615) -- (0.1950,2.2025) -- (0.3888,2.1917) -- cycle;
\path[cbpmesh=19.9] (-0.1856,2.1615) -- (-0.3701,2.1513) -- (-0.3888,2.1917) -- (-0.1950,2.2025) -- cycle;
\path[cbpmesh=38] (0.9947,-0.0593) -- (0.8843,-0.0246) -- (0.8395,-0.0000) -- (0.9444,-0.0329) -- cycle;
\draw[draw=cbpInk!70,line width=.45pt] (0.9947,-0.0593) -- (0.8843,-0.0246);
\path[cbpmesh=38] (-0.8843,-0.0246) -- (-0.9947,-0.0593) -- (-0.9444,-0.0329) -- (-0.8395,-0.0000) -- cycle;
\draw[draw=cbpInk!70,line width=.45pt] (-0.8843,-0.0246) -- (-0.9947,-0.0593);
\path[cbpmesh=19.8] (-3.2147,1.8628) -- (-3.3880,1.7410) -- (-3.4447,1.7296) -- (-3.2685,1.8534) -- cycle;
\path[cbpmesh=25.7] (3.3880,1.7410) -- (3.2147,1.8628) -- (3.2685,1.8534) -- (3.4447,1.7296) -- cycle;
\path[cbpmesh=20.6] (-0.7680,2.1490) -- (-0.9510,2.1173) -- (-0.9978,2.1439) -- (-0.8057,2.1772) -- cycle;
\path[cbpmesh=22.0] (0.9510,2.1173) -- (0.7680,2.1490) -- (0.8057,2.1772) -- (0.9978,2.1439) -- cycle;
\path[cbpmesh=19.2] (-1.9891,1.9732) -- (-2.1657,1.8966) -- (-2.2341,1.9109) -- (-2.0519,1.9900) -- cycle;
\path[cbpmesh=21.5] (2.1657,1.8966) -- (1.9891,1.9732) -- (2.0519,1.9900) -- (2.2341,1.9109) -- cycle;
\path[cbpmesh=38.0] (0.4726,0.1152) -- (0.3570,0.1305) -- (0.3411,0.1932) -- (0.4515,0.1786) -- cycle;
\path[cbpmesh=38] (-0.3570,0.1305) -- (-0.4726,0.1152) -- (-0.4515,0.1786) -- (-0.3411,0.1932) -- cycle;
\path[cbpmesh=22.5] (-3.6235,1.6463) -- (-3.7866,1.5096) -- (-3.8536,1.4600) -- (-3.6877,1.5991) -- cycle;
\draw[draw=cbpInk!70,line width=.45pt] (-3.8536,1.4600) -- (-3.6877,1.5991);
\draw[draw=cbpInk!70,line width=.45pt] (-3.6877,1.5991) -- (-3.6235,1.6463);
\path[cbpmesh=32.9] (3.7866,1.5096) -- (3.6235,1.6463) -- (3.6877,1.5991) -- (3.8536,1.4600) -- cycle;
\draw[draw=cbpInk!70,line width=.45pt] (3.6235,1.6463) -- (3.6877,1.5991);
\draw[draw=cbpInk!70,line width=.45pt] (3.6877,1.5991) -- (3.8536,1.4600);
\path[cbpmesh=20.2] (0.5523,2.1343) -- (0.3701,2.1513) -- (0.3888,2.1917) -- (0.5802,2.1739) -- cycle;
\path[cbpmesh=19.7] (-0.3701,2.1513) -- (-0.5523,2.1343) -- (-0.5802,2.1739) -- (-0.3888,2.1917) -- cycle;
\path[cbpmesh=24.5] (1.6058,2.0124) -- (1.4275,2.0684) -- (1.4889,2.0758) -- (1.6749,2.0174) -- cycle;
\path[cbpmesh=21.1] (-1.4275,2.0684) -- (-1.6058,2.0124) -- (-1.6749,2.0174) -- (-1.4889,2.0758) -- cycle;
\path[cbpmesh=18.3] (-2.7137,1.9164) -- (-2.8872,1.8131) -- (-2.9447,1.8294) -- (-2.7678,1.9348) -- cycle;
\path[cbpmesh=19.0] (2.8872,1.8131) -- (2.7137,1.9164) -- (2.7678,1.9348) -- (2.9447,1.8294) -- cycle;
\path[cbpmesh=38] (-0.7295,0.0289) -- (-0.8395,-0.0000) -- (-0.7991,0.0445) -- (-0.6943,0.0721) -- cycle;
\path[cbpmesh=38] (0.8395,-0.0000) -- (0.7295,0.0289) -- (0.6943,0.0721) -- (0.7991,0.0445) -- cycle;
\path[cbpmesh=21.1] (-1.1837,2.1037) -- (-1.3624,2.0567) -- (-1.4275,2.0684) -- (-1.2403,2.1176) -- cycle;
\path[cbpmesh=23.8] (1.3624,2.0567) -- (1.1837,2.1037) -- (1.2403,2.1176) -- (1.4275,2.0684) -- cycle;
\path[cbpmesh=23.4] (1.8506,1.9524) -- (1.6749,2.0174) -- (1.7396,2.0280) -- (1.9221,1.9604) -- cycle;
\path[cbpmesh=20.3] (-1.6749,2.0174) -- (-1.8506,1.9524) -- (-1.9221,1.9604) -- (-1.7396,2.0280) -- cycle;
\path[cbpmesh=18.6] (-2.2984,1.9285) -- (-2.4717,1.8405) -- (-2.5370,1.8580) -- (-2.3591,1.9484) -- cycle;
\path[cbpmesh=20.0] (2.4717,1.8405) -- (2.2984,1.9285) -- (2.3591,1.9484) -- (2.5370,1.8580) -- cycle;
\path[cbpmesh=19.1] (-3.1615,1.8656) -- (-3.3319,1.7458) -- (-3.3880,1.7410) -- (-3.2147,1.8628) -- cycle;
\path[cbpmesh=23.9] (3.3319,1.7458) -- (3.1615,1.8656) -- (3.2147,1.8628) -- (3.3880,1.7410) -- cycle;
\path[cbpmesh=37.7] (0.5853,0.0957) -- (0.4726,0.1152) -- (0.4515,0.1786) -- (0.5591,0.1600) -- cycle;
\path[cbpmesh=38] (-0.4726,0.1152) -- (-0.5853,0.0957) -- (-0.5591,0.1600) -- (-0.4515,0.1786) -- cycle;
\path[cbpmesh=20.3] (0.7311,2.1106) -- (0.5523,2.1343) -- (0.5802,2.1739) -- (0.7680,2.1490) -- cycle;
\path[cbpmesh=19.6] (-0.5523,2.1343) -- (-0.7311,2.1106) -- (-0.7680,2.1490) -- (-0.5802,2.1739) -- cycle;
\path[cbpmesh=38] (-0.9947,-0.0593) -- (-1.0991,-0.0979) -- (-1.0435,-0.0696) -- (-0.9444,-0.0329) -- cycle;
\draw[draw=cbpInk!70,line width=.45pt] (-0.9947,-0.0593) -- (-1.0991,-0.0979);
\path[cbpmesh=38] (1.0991,-0.0979) -- (0.9947,-0.0593) -- (0.9444,-0.0329) -- (1.0435,-0.0696) -- cycle;
\draw[draw=cbpInk!70,line width=.45pt] (1.0991,-0.0979) -- (0.9947,-0.0593);
\path[cbpmesh=20.4] (-0.9510,2.1173) -- (-1.1283,2.0789) -- (-1.1837,2.1037) -- (-0.9978,2.1439) -- cycle;
\path[cbpmesh=22.0] (1.1283,2.0789) -- (0.9510,2.1173) -- (0.9978,2.1439) -- (1.1837,2.1037) -- cycle;
\path[cbpmesh=18.3] (-2.6576,1.8969) -- (-2.8274,1.7957) -- (-2.8872,1.8131) -- (-2.7137,1.9164) -- cycle;
\path[cbpmesh=18.9] (2.8274,1.7957) -- (2.6576,1.8969) -- (2.7137,1.9164) -- (2.8872,1.8131) -- cycle;
\path[cbpmesh=36.0] (0.1146,0.2101) -- (0.0000,0.2122) -- (0.0000,0.2989) -- (0.1100,0.2969) -- cycle;
\path[cbpmesh=36.6] (0.0000,0.2122) -- (-0.1146,0.2101) -- (-0.1100,0.2969) -- (0.0000,0.2989) -- cycle;
\path[cbpmesh=21.7] (-3.5620,1.6832) -- (-3.7223,1.5488) -- (-3.7866,1.5096) -- (-3.6235,1.6463) -- cycle;
\path[cbpmesh=31.3] (3.7223,1.5488) -- (3.5620,1.6832) -- (3.6235,1.6463) -- (3.7866,1.5096) -- cycle;
\path[cbpmesh=18.7] (-3.1083,1.8628) -- (-3.2759,1.7450) -- (-3.3319,1.7458) -- (-3.1615,1.8656) -- cycle;
\path[cbpmesh=22.3] (3.2759,1.7450) -- (3.1083,1.8628) -- (3.1615,1.8656) -- (3.3319,1.7458) -- cycle;
\path[cbpmesh=22.4] (2.0927,1.8864) -- (1.9221,1.9604) -- (1.9891,1.9732) -- (2.1657,1.8966) -- cycle;
\path[cbpmesh=19.5] (-1.9221,1.9604) -- (-2.0927,1.8864) -- (-2.1657,1.8966) -- (-1.9891,1.9732) -- cycle;
\path[cbpmesh=35.4] (0.2285,0.2037) -- (0.1146,0.2101) -- (0.1100,0.2969) -- (0.2194,0.2908) -- cycle;
\path[cbpmesh=37.2] (-0.1146,0.2101) -- (-0.2285,0.2037) -- (-0.2194,0.2908) -- (-0.1100,0.2969) -- cycle;
\path[cbpmesh=37.4] (0.6943,0.0721) -- (0.5853,0.0957) -- (0.5591,0.1600) -- (0.6633,0.1374) -- cycle;
\path[cbpmesh=38] (-0.5853,0.0957) -- (-0.6943,0.0721) -- (-0.6633,0.1374) -- (-0.5591,0.1600) -- cycle;
\path[cbpmesh=18.7] (0.0000,2.1142) -- (-0.1765,2.1109) -- (-0.1856,2.1615) -- (0.0000,2.1650) -- cycle;
\path[cbpmesh=18.7] (0.1765,2.1109) -- (0.0000,2.1142) -- (0.0000,2.1650) -- (0.1856,2.1615) -- cycle;
\path[cbpmesh=38] (-0.8395,-0.0000) -- (-0.9444,-0.0329) -- (-0.8989,0.0132) -- (-0.7991,0.0445) -- cycle;
\path[cbpmesh=38] (0.9444,-0.0329) -- (0.8395,-0.0000) -- (0.7991,0.0445) -- (0.8989,0.0132) -- cycle;
\path[cbpmesh=20.6] (2.4025,1.8254) -- (2.2341,1.9109) -- (2.2984,1.9285) -- (2.4717,1.8405) -- cycle;
\path[cbpmesh=18.7] (-2.2341,1.9109) -- (-2.4025,1.8254) -- (-2.4717,1.8405) -- (-2.2984,1.9285) -- cycle;
\path[cbpmesh=18.7] (0.3520,2.1011) -- (0.1765,2.1109) -- (0.1856,2.1615) -- (0.3701,2.1513) -- cycle;
\path[cbpmesh=18.6] (-0.1765,2.1109) -- (-0.3520,2.1011) -- (-0.3701,2.1513) -- (-0.1856,2.1615) -- cycle;
\path[cbpmesh=34.9] (0.3411,0.1932) -- (0.2285,0.2037) -- (0.2194,0.2908) -- (0.3274,0.2807) -- cycle;
\path[cbpmesh=37.9] (-0.2285,0.2037) -- (-0.3411,0.1932) -- (-0.3274,0.2807) -- (-0.2194,0.2908) -- cycle;
\path[cbpmesh=19.5] (-0.7311,2.1106) -- (-0.9054,2.0804) -- (-0.9510,2.1173) -- (-0.7680,2.1490) -- cycle;
\path[cbpmesh=20.3] (0.9054,2.0804) -- (0.7311,2.1106) -- (0.7680,2.1490) -- (0.9510,2.1173) -- cycle;
\path[cbpmesh=18.3] (-2.5988,1.8771) -- (-2.7649,1.7782) -- (-2.8274,1.7957) -- (-2.6576,1.8969) -- cycle;
\path[cbpmesh=19.1] (2.7649,1.7782) -- (2.5988,1.8771) -- (2.6576,1.8969) -- (2.8274,1.7957) -- cycle;
\path[cbpmesh=18.4] (-3.0548,1.8552) -- (-3.2195,1.7394) -- (-3.2759,1.7450) -- (-3.1083,1.8628) -- cycle;
\path[cbpmesh=21.0] (3.2195,1.7394) -- (3.0548,1.8552) -- (3.1083,1.8628) -- (3.2759,1.7450) -- cycle;
\path[cbpmesh=34.5] (0.4515,0.1786) -- (0.3411,0.1932) -- (0.3274,0.2807) -- (0.4334,0.2667) -- cycle;
\path[cbpmesh=38] (-0.3411,0.1932) -- (-0.4515,0.1786) -- (-0.4334,0.2667) -- (-0.3274,0.2807) -- cycle;
\path[cbpmesh=18.7] (0.5253,2.0850) -- (0.3520,2.1011) -- (0.3701,2.1513) -- (0.5523,2.1343) -- cycle;
\path[cbpmesh=18.6] (-0.3520,2.1011) -- (-0.5253,2.0850) -- (-0.5523,2.1343) -- (-0.3701,2.1513) -- cycle;
\path[cbpmesh=23.9] (1.5326,2.0033) -- (1.3624,2.0567) -- (1.4275,2.0684) -- (1.6058,2.0124) -- cycle;
\path[cbpmesh=20.8] (-1.3624,2.0567) -- (-1.5326,2.0033) -- (-1.6058,2.0124) -- (-1.4275,2.0684) -- cycle;
\path[cbpmesh=38] (-0.6943,0.0721) -- (-0.7991,0.0445) -- (-0.7634,0.1111) -- (-0.6633,0.1374) -- cycle;
\path[cbpmesh=37.2] (0.7991,0.0445) -- (0.6943,0.0721) -- (0.6633,0.1374) -- (0.7634,0.1111) -- cycle;
\path[cbpmesh=38] (1.1967,-0.1402) -- (1.0991,-0.0979) -- (1.0435,-0.0696) -- (1.1362,-0.1098) -- cycle;
\draw[draw=cbpInk!70,line width=.45pt] (1.1967,-0.1402) -- (1.0991,-0.0979);
\draw[draw=cbpInk!70,line width=.45pt] (1.1362,-0.1098) -- (1.1967,-0.1402);
\path[cbpmesh=38] (-1.0991,-0.0979) -- (-1.1967,-0.1402) -- (-1.1362,-0.1098) -- (-1.0435,-0.0696) -- cycle;
\draw[draw=cbpInk!70,line width=.45pt] (-1.0991,-0.0979) -- (-1.1967,-0.1402);
\draw[draw=cbpInk!70,line width=.45pt] (-1.1967,-0.1402) -- (-1.1362,-0.1098);
\path[cbpmesh=24.5] (1.7743,1.9501) -- (1.6058,2.0124) -- (1.6749,2.0174) -- (1.8506,1.9524) -- cycle;
\path[cbpmesh=20.8] (-1.6058,2.0124) -- (-1.7743,1.9501) -- (-1.8506,1.9524) -- (-1.6749,2.0174) -- cycle;
\path[cbpmesh=20.8] (-3.5026,1.7106) -- (-3.6602,1.5785) -- (-3.7223,1.5488) -- (-3.5620,1.6832) -- cycle;
\path[cbpmesh=29.4] (3.6602,1.5785) -- (3.5026,1.7106) -- (3.5620,1.6832) -- (3.7223,1.5488) -- cycle;
\path[cbpmesh=22.1] (1.2985,2.0342) -- (1.1283,2.0789) -- (1.1837,2.1037) -- (1.3624,2.0567) -- cycle;
\path[cbpmesh=20.1] (-1.1283,2.0789) -- (-1.2985,2.0342) -- (-1.3624,2.0567) -- (-1.1837,2.1037) -- cycle;
\path[cbpmesh=38] (1.0435,-0.0696) -- (0.9444,-0.0329) -- (0.8989,0.0132) -- (0.9932,-0.0217) -- cycle;
\path[cbpmesh=38] (-0.9444,-0.0329) -- (-1.0435,-0.0696) -- (-0.9932,-0.0217) -- (-0.8989,0.0132) -- cycle;
\path[cbpmesh=34.2] (0.5591,0.1600) -- (0.4515,0.1786) -- (0.4334,0.2667) -- (0.5367,0.2488) -- cycle;
\path[cbpmesh=38] (-0.4515,0.1786) -- (-0.5591,0.1600) -- (-0.5367,0.2488) -- (-0.4334,0.2667) -- cycle;
\path[cbpmesh=18.3] (-3.0004,1.8437) -- (-3.1622,1.7301) -- (-3.2195,1.7394) -- (-3.0548,1.8552) -- cycle;
\path[cbpmesh=20.0] (3.1622,1.7301) -- (3.0004,1.8437) -- (3.0548,1.8552) -- (3.2195,1.7394) -- cycle;
\path[cbpmesh=18.9] (-2.1657,1.8966) -- (-2.3289,1.8137) -- (-2.4025,1.8254) -- (-2.2341,1.9109) -- cycle;
\path[cbpmesh=21.5] (2.3289,1.8137) -- (2.1657,1.8966) -- (2.2341,1.9109) -- (2.4025,1.8254) -- cycle;
\path[cbpmesh=18.4] (-2.5370,1.8580) -- (-2.6992,1.7615) -- (-2.7649,1.7782) -- (-2.5988,1.8771) -- cycle;
\path[cbpmesh=19.4] (2.6992,1.7615) -- (2.5370,1.8580) -- (2.5988,1.8771) -- (2.7649,1.7782) -- cycle;
\path[cbpmesh=18.6] (-0.5253,2.0850) -- (-0.6953,2.0625) -- (-0.7311,2.1106) -- (-0.5523,2.1343) -- cycle;
\path[cbpmesh=18.7] (0.6953,2.0625) -- (0.5253,2.0850) -- (0.5523,2.1343) -- (0.7311,2.1106) -- cycle;
\path[cbpmesh=19.9] (-1.8506,1.9524) -- (-2.0149,1.8811) -- (-2.0927,1.8864) -- (-1.9221,1.9604) -- cycle;
\path[cbpmesh=23.4] (2.0149,1.8811) -- (1.8506,1.9524) -- (1.9221,1.9604) -- (2.0927,1.8864) -- cycle;
\path[cbpmesh=20.4] (1.0741,2.0439) -- (0.9054,2.0804) -- (0.9510,2.1173) -- (1.1283,2.0789) -- cycle;
\path[cbpmesh=19.4] (-0.9054,2.0804) -- (-1.0741,2.0439) -- (-1.1283,2.0789) -- (-0.9510,2.1173) -- cycle;
\path[cbpmesh=17.5] (0.1678,2.0511) -- (0.0000,2.0542) -- (0.0000,2.1142) -- (0.1765,2.1109) -- cycle;
\path[cbpmesh=17.6] (0.0000,2.0542) -- (-0.1678,2.0511) -- (-0.1765,2.1109) -- (0.0000,2.1142) -- cycle;
\path[cbpmesh=38] (-0.7991,0.0445) -- (-0.8989,0.0132) -- (-0.8588,0.0812) -- (-0.7634,0.1111) -- cycle;
\path[cbpmesh=37.1] (0.8989,0.0132) -- (0.7991,0.0445) -- (0.7634,0.1111) -- (0.8588,0.0812) -- cycle;
\path[cbpmesh=20.0] (-3.4447,1.7296) -- (-3.5998,1.5996) -- (-3.6602,1.5785) -- (-3.5026,1.7106) -- cycle;
\path[cbpmesh=27.5] (3.5998,1.5996) -- (3.4447,1.7296) -- (3.5026,1.7106) -- (3.6602,1.5785) -- cycle;
\path[cbpmesh=33.2] (0.1100,0.2969) -- (0.0000,0.2989) -- (0.0000,0.4127) -- (0.1063,0.4108) -- cycle;
\path[cbpmesh=33.8] (0.0000,0.2989) -- (-0.1100,0.2969) -- (-0.1063,0.4108) -- (0.0000,0.4127) -- cycle;
\path[cbpmesh=38] (-0.5591,0.1600) -- (-0.6633,0.1374) -- (-0.6368,0.2272) -- (-0.5367,0.2488) -- cycle;
\path[cbpmesh=33.9] (0.6633,0.1374) -- (0.5591,0.1600) -- (0.5367,0.2488) -- (0.6368,0.2272) -- cycle;
\path[cbpmesh=19.3] (3.1035,1.7179) -- (2.9447,1.8294) -- (3.0004,1.8437) -- (3.1622,1.7301) -- cycle;
\path[cbpmesh=18.3] (-2.9447,1.8294) -- (-3.1035,1.7179) -- (-3.1622,1.7301) -- (-3.0004,1.8437) -- cycle;
\path[cbpmesh=32.6] (0.2194,0.2908) -- (0.1100,0.2969) -- (0.1063,0.4108) -- (0.2120,0.4049) -- cycle;
\path[cbpmesh=34.5] (-0.1100,0.2969) -- (-0.2194,0.2908) -- (-0.2120,0.4049) -- (-0.1063,0.4108) -- cycle;
\path[cbpmesh=17.5] (0.3345,2.0419) -- (0.1678,2.0511) -- (0.1765,2.1109) -- (0.3520,2.1011) -- cycle;
\path[cbpmesh=17.6] (-0.1678,2.0511) -- (-0.3345,2.0419) -- (-0.3520,2.1011) -- (-0.1765,2.1109) -- cycle;
\path[cbpmesh=18.4] (-2.4717,1.8405) -- (-2.6297,1.7464) -- (-2.6992,1.7615) -- (-2.5370,1.8580) -- cycle;
\path[cbpmesh=19.9] (2.6297,1.7464) -- (2.4717,1.8405) -- (2.5370,1.8580) -- (2.6992,1.7615) -- cycle;
\path[cbpmesh=18.8] (0.8611,2.0337) -- (0.6953,2.0625) -- (0.7311,2.1106) -- (0.9054,2.0804) -- cycle;
\path[cbpmesh=18.6] (-0.6953,2.0625) -- (-0.8611,2.0337) -- (-0.9054,2.0804) -- (-0.7311,2.1106) -- cycle;
\path[cbpmesh=35.2] (-0.2194,0.2908) -- (-0.3274,0.2807) -- (-0.3163,0.3952) -- (-0.2120,0.4049) -- cycle;
\path[cbpmesh=32.1] (0.3274,0.2807) -- (0.2194,0.2908) -- (0.2120,0.4049) -- (0.3163,0.3952) -- cycle;
\path[cbpmesh=17.4] (0.4992,2.0265) -- (0.3345,2.0419) -- (0.3520,2.1011) -- (0.5253,2.0850) -- cycle;
\path[cbpmesh=17.7] (-0.3345,2.0419) -- (-0.4992,2.0265) -- (-0.5253,2.0850) -- (-0.3520,2.1011) -- cycle;
\path[cbpmesh=38] (-1.0435,-0.0696) -- (-1.1362,-0.1098) -- (-1.0814,-0.0600) -- (-0.9932,-0.0217) -- cycle;
\path[cbpmesh=38] (1.1362,-0.1098) -- (1.0435,-0.0696) -- (0.9932,-0.0217) -- (1.0814,-0.0600) -- cycle;
\path[cbpmesh=19.2] (-2.0927,1.8864) -- (-2.2505,1.8062) -- (-2.3289,1.8137) -- (-2.1657,1.8966) -- cycle;
\path[cbpmesh=22.4] (2.2505,1.8062) -- (2.0927,1.8864) -- (2.1657,1.8966) -- (2.3289,1.8137) -- cycle;
\path[cbpmesh=19.3] (-3.3880,1.7410) -- (-3.5405,1.6132) -- (-3.5998,1.5996) -- (-3.4447,1.7296) -- cycle;
\path[cbpmesh=25.6] (3.5405,1.6132) -- (3.3880,1.7410) -- (3.4447,1.7296) -- (3.5998,1.5996) -- cycle;
\path[cbpmesh=19.0] (3.0429,1.7037) -- (2.8872,1.8131) -- (2.9447,1.8294) -- (3.1035,1.7179) -- cycle;
\path[cbpmesh=18.3] (-2.8872,1.8131) -- (-3.0429,1.7037) -- (-3.1035,1.7179) -- (-2.9447,1.8294) -- cycle;
\path[cbpmesh=36.0] (-0.3274,0.2807) -- (-0.4334,0.2667) -- (-0.4187,0.3816) -- (-0.3163,0.3952) -- cycle;
\path[cbpmesh=31.7] (0.4334,0.2667) -- (0.3274,0.2807) -- (0.3163,0.3952) -- (0.4187,0.3816) -- cycle;
\path[cbpmesh=38] (-0.6633,0.1374) -- (-0.7634,0.1111) -- (-0.7328,0.2019) -- (-0.6368,0.2272) -- cycle;
\path[cbpmesh=33.8] (0.7634,0.1111) -- (0.6633,0.1374) -- (0.6368,0.2272) -- (0.7328,0.2019) -- cycle;
\path[cbpmesh=22.1] (1.4608,1.9833) -- (1.2985,2.0342) -- (1.3624,2.0567) -- (1.5326,2.0033) -- cycle;
\path[cbpmesh=19.9] (-1.2985,2.0342) -- (-1.4608,1.9833) -- (-1.5326,2.0033) -- (-1.3624,2.0567) -- cycle;
\path[cbpmesh=20.5] (-1.5326,2.0033) -- (-1.6934,1.9438) -- (-1.7743,1.9501) -- (-1.6058,2.0124) -- cycle;
\path[cbpmesh=23.9] (1.6934,1.9438) -- (1.5326,2.0033) -- (1.6058,2.0124) -- (1.7743,1.9501) -- cycle;
\path[cbpmesh=38] (-0.8989,0.0132) -- (-0.9932,-0.0217) -- (-0.9488,0.0478) -- (-0.8588,0.0812) -- cycle;
\path[cbpmesh=37.1] (0.9932,-0.0217) -- (0.8989,0.0132) -- (0.8588,0.0812) -- (0.9488,0.0478) -- cycle;
\path[cbpmesh=20.4] (1.2362,2.0013) -- (1.0741,2.0439) -- (1.1283,2.0789) -- (1.2985,2.0342) -- cycle;
\path[cbpmesh=19.2] (-1.0741,2.0439) -- (-1.2362,2.0013) -- (-1.2985,2.0342) -- (-1.1283,2.0789) -- cycle;
\path[cbpmesh=17.4] (0.6608,2.0051) -- (0.4992,2.0265) -- (0.5253,2.0850) -- (0.6953,2.0625) -- cycle;
\path[cbpmesh=17.7] (-0.4992,2.0265) -- (-0.6608,2.0051) -- (-0.6953,2.0625) -- (-0.5253,2.0850) -- cycle;
\path[cbpmesh=20.4] (-1.7743,1.9501) -- (-1.9318,1.8817) -- (-2.0149,1.8811) -- (-1.8506,1.9524) -- cycle;
\path[cbpmesh=24.5] (1.9318,1.8817) -- (1.7743,1.9501) -- (1.8506,1.9524) -- (2.0149,1.8811) -- cycle;
\path[cbpmesh=36.9] (-0.4334,0.2667) -- (-0.5367,0.2488) -- (-0.5185,0.3643) -- (-0.4187,0.3816) -- cycle;
\path[cbpmesh=31.3] (0.5367,0.2488) -- (0.4334,0.2667) -- (0.4187,0.3816) -- (0.5185,0.3643) -- cycle;
\path[cbpmesh=18.5] (-2.4025,1.8254) -- (-2.5561,1.7340) -- (-2.6297,1.7464) -- (-2.4717,1.8405) -- cycle;
\path[cbpmesh=20.6] (2.5561,1.7340) -- (2.4025,1.8254) -- (2.4717,1.8405) -- (2.6297,1.7464) -- cycle;
\path[cbpmesh=18.3] (-2.8274,1.7957) -- (-2.9799,1.6886) -- (-3.0429,1.7037) -- (-2.8872,1.8131) -- cycle;
\path[cbpmesh=18.9] (2.9799,1.6886) -- (2.8274,1.7957) -- (2.8872,1.8131) -- (3.0429,1.7037) -- cycle;
\path[cbpmesh=18.8] (1.0215,1.9990) -- (0.8611,2.0337) -- (0.9054,2.0804) -- (1.0741,2.0439) -- cycle;
\path[cbpmesh=18.5] (-0.8611,2.0337) -- (-1.0215,1.9990) -- (-1.0741,2.0439) -- (-0.9054,2.0804) -- cycle;
\path[cbpmesh=18.8] (-3.3319,1.7458) -- (-3.4818,1.6202) -- (-3.5405,1.6132) -- (-3.3880,1.7410) -- cycle;
\path[cbpmesh=23.8] (3.4818,1.6202) -- (3.3319,1.7458) -- (3.3880,1.7410) -- (3.5405,1.6132) -- cycle;
\path[cbpmesh=16.5] (0.1594,1.9827) -- (0.0000,1.9857) -- (0.0000,2.0542) -- (0.1678,2.0511) -- cycle;
\path[cbpmesh=16.7] (0.0000,1.9857) -- (-0.1594,1.9827) -- (-0.1678,2.0511) -- (0.0000,2.0542) -- cycle;
\path[cbpmesh=21.8] (-3.7866,1.5096) -- (-3.9263,1.3673) -- (-3.9958,1.3152) -- (-3.8536,1.4600) -- cycle;
\draw[draw=cbpInk!70,line width=.45pt] (-3.9958,1.3152) -- (-3.8536,1.4600);
\path[cbpmesh=32.5] (3.9263,1.3673) -- (3.7866,1.5096) -- (3.8536,1.4600) -- (3.9958,1.3152) -- cycle;
\draw[draw=cbpInk!70,line width=.45pt] (3.8536,1.4600) -- (3.9958,1.3152);
\path[cbpmesh=16.4] (0.3178,1.9739) -- (0.1594,1.9827) -- (0.1678,2.0511) -- (0.3345,2.0419) -- cycle;
\path[cbpmesh=16.8] (-0.1594,1.9827) -- (-0.3178,1.9739) -- (-0.3345,2.0419) -- (-0.1678,2.0511) -- cycle;
\path[cbpmesh=33.7] (0.8588,0.0812) -- (0.7634,0.1111) -- (0.7328,0.2019) -- (0.8244,0.1732) -- cycle;
\path[cbpmesh=38] (-0.7634,0.1111) -- (-0.8588,0.0812) -- (-0.8244,0.1732) -- (-0.7328,0.2019) -- cycle;
\path[cbpmesh=29.4] (0.1063,0.4108) -- (0.0000,0.4127) -- (0.0000,0.5473) -- (0.1047,0.5454) -- cycle;
\path[cbpmesh=30.0] (0.0000,0.4127) -- (-0.1063,0.4108) -- (-0.1047,0.5454) -- (0.0000,0.5473) -- cycle;
\path[cbpmesh=38] (1.2218,-0.1533) -- (1.1362,-0.1098) -- (1.0814,-0.0600) -- (1.1629,-0.1014) -- cycle;
\draw[draw=cbpInk!70,line width=.45pt] (1.2218,-0.1533) -- (1.1362,-0.1098);
\path[cbpmesh=38] (-1.1362,-0.1098) -- (-1.2218,-0.1533) -- (-1.1629,-0.1014) -- (-1.0814,-0.0600) -- cycle;
\draw[draw=cbpInk!70,line width=.45pt] (-1.1362,-0.1098) -- (-1.2218,-0.1533);
\path[cbpmesh=31.1] (0.6368,0.2272) -- (0.5367,0.2488) -- (0.5185,0.3643) -- (0.6152,0.3434) -- cycle;
\path[cbpmesh=37.8] (-0.5367,0.2488) -- (-0.6368,0.2272) -- (-0.6152,0.3434) -- (-0.5185,0.3643) -- cycle;
\path[cbpmesh=28.8] (0.2120,0.4049) -- (0.1063,0.4108) -- (0.1047,0.5454) -- (0.2088,0.5396) -- cycle;
\path[cbpmesh=30.8] (-0.1063,0.4108) -- (-0.2120,0.4049) -- (-0.2088,0.5396) -- (-0.1047,0.5454) -- cycle;
\path[cbpmesh=17.4] (0.8183,1.9778) -- (0.6608,2.0051) -- (0.6953,2.0625) -- (0.8611,2.0337) -- cycle;
\path[cbpmesh=17.8] (-0.6608,2.0051) -- (-0.8183,1.9778) -- (-0.8611,2.0337) -- (-0.6953,2.0625) -- cycle;
\path[cbpmesh=18.3] (-2.7649,1.7782) -- (-2.9140,1.6735) -- (-2.9799,1.6886) -- (-2.8274,1.7957) -- cycle;
\path[cbpmesh=19.0] (2.9140,1.6735) -- (2.7649,1.7782) -- (2.8274,1.7957) -- (2.9799,1.6886) -- cycle;
\path[cbpmesh=19.6] (-2.0149,1.8811) -- (-2.1668,1.8040) -- (-2.2505,1.8062) -- (-2.0927,1.8864) -- cycle;
\path[cbpmesh=23.4] (2.1668,1.8040) -- (2.0149,1.8811) -- (2.0927,1.8864) -- (2.2505,1.8062) -- cycle;
\path[cbpmesh=16.3] (0.4742,1.9593) -- (0.3178,1.9739) -- (0.3345,2.0419) -- (0.4992,2.0265) -- cycle;
\path[cbpmesh=16.9] (-0.3178,1.9739) -- (-0.4742,1.9593) -- (-0.4992,2.0265) -- (-0.3345,2.0419) -- cycle;
\path[cbpmesh=18.4] (-3.2759,1.7450) -- (-3.4233,1.6214) -- (-3.4818,1.6202) -- (-3.3319,1.7458) -- cycle;
\path[cbpmesh=22.2] (3.4233,1.6214) -- (3.2759,1.7450) -- (3.3319,1.7458) -- (3.4818,1.6202) -- cycle;
\path[cbpmesh=38] (-0.9932,-0.0217) -- (-1.0814,-0.0600) -- (-1.0331,0.0113) -- (-0.9488,0.0478) -- cycle;
\path[cbpmesh=37.2] (1.0814,-0.0600) -- (0.9932,-0.0217) -- (0.9488,0.0478) -- (1.0331,0.0113) -- cycle;
\path[cbpmesh=28.3] (0.3163,0.3952) -- (0.2120,0.4049) -- (0.2088,0.5396) -- (0.3116,0.5300) -- cycle;
\path[cbpmesh=31.5] (-0.2120,0.4049) -- (-0.3163,0.3952) -- (-0.3116,0.5300) -- (-0.2088,0.5396) -- cycle;
\path[cbpmesh=18.7] (-2.3289,1.8137) -- (-2.4778,1.7250) -- (-2.5561,1.7340) -- (-2.4025,1.8254) -- cycle;
\path[cbpmesh=21.4] (2.4778,1.7250) -- (2.3289,1.8137) -- (2.4025,1.8254) -- (2.5561,1.7340) -- cycle;
\path[cbpmesh=20.4] (1.3906,1.9528) -- (1.2362,2.0013) -- (1.2985,2.0342) -- (1.4608,1.9833) -- cycle;
\path[cbpmesh=19.1] (-1.2362,2.0013) -- (-1.3906,1.9528) -- (-1.4608,1.9833) -- (-1.2985,2.0342) -- cycle;
\path[cbpmesh=27.8] (0.4187,0.3816) -- (0.3163,0.3952) -- (0.3116,0.5300) -- (0.4124,0.5166) -- cycle;
\path[cbpmesh=32.4] (-0.3163,0.3952) -- (-0.4187,0.3816) -- (-0.4124,0.5166) -- (-0.3116,0.5300) -- cycle;
\path[cbpmesh=21.0] (-3.7223,1.5488) -- (-3.8596,1.4089) -- (-3.9263,1.3673) -- (-3.7866,1.5096) -- cycle;
\path[cbpmesh=31.0] (3.8596,1.4089) -- (3.7223,1.5488) -- (3.7866,1.5096) -- (3.9263,1.3673) -- cycle;
\path[cbpmesh=19.6] (-1.4608,1.9833) -- (-1.6140,1.9266) -- (-1.6934,1.9438) -- (-1.5326,2.0033) -- cycle;
\path[cbpmesh=22.1] (1.6140,1.9266) -- (1.4608,1.9833) -- (1.5326,2.0033) -- (1.6934,1.9438) -- cycle;
\path[cbpmesh=16.2] (0.6277,1.9390) -- (0.4742,1.9593) -- (0.4992,2.0265) -- (0.6608,2.0051) -- cycle;
\path[cbpmesh=17.1] (-0.4742,1.9593) -- (-0.6277,1.9390) -- (-0.6608,2.0051) -- (-0.4992,2.0265) -- cycle;
\path[cbpmesh=30.9] (0.7328,0.2019) -- (0.6368,0.2272) -- (0.6152,0.3434) -- (0.7080,0.3190) -- cycle;
\path[cbpmesh=38] (-0.6368,0.2272) -- (-0.7328,0.2019) -- (-0.7080,0.3190) -- (-0.6152,0.3434) -- cycle;
\path[cbpmesh=18.8] (1.1757,1.9585) -- (1.0215,1.9990) -- (1.0741,2.0439) -- (1.2362,2.0013) -- cycle;
\path[cbpmesh=18.5] (-1.0215,1.9990) -- (-1.1757,1.9585) -- (-1.2362,2.0013) -- (-1.0741,2.0439) -- cycle;
\path[cbpmesh=38] (-0.8588,0.0812) -- (-0.9488,0.0478) -- (-0.9109,0.1412) -- (-0.8244,0.1732) -- cycle;
\path[cbpmesh=33.6] (0.9488,0.0478) -- (0.8588,0.0812) -- (0.8244,0.1732) -- (0.9109,0.1412) -- cycle;
\path[cbpmesh=18.2] (-3.2195,1.7394) -- (-3.3644,1.6180) -- (-3.4233,1.6214) -- (-3.2759,1.7450) -- cycle;
\path[cbpmesh=20.9] (3.3644,1.6180) -- (3.2195,1.7394) -- (3.2759,1.7450) -- (3.4233,1.6214) -- cycle;
\path[cbpmesh=18.3] (-2.6992,1.7615) -- (-2.8447,1.6592) -- (-2.9140,1.6735) -- (-2.7649,1.7782) -- cycle;
\path[cbpmesh=19.4] (2.8447,1.6592) -- (2.6992,1.7615) -- (2.7649,1.7782) -- (2.9140,1.6735) -- cycle;
\path[cbpmesh=16.0] (0.0000,1.9090) -- (-0.1514,1.9063) -- (-0.1594,1.9827) -- (0.0000,1.9857) -- cycle;
\path[cbpmesh=15.8] (0.1514,1.9063) -- (0.0000,1.9090) -- (0.0000,1.9857) -- (0.1594,1.9827) -- cycle;
\path[cbpmesh=20.1] (-1.6934,1.9438) -- (-1.8437,1.8785) -- (-1.9318,1.8817) -- (-1.7743,1.9501) -- cycle;
\path[cbpmesh=23.9] (1.8437,1.8785) -- (1.6934,1.9438) -- (1.7743,1.9501) -- (1.9318,1.8817) -- cycle;
\path[cbpmesh=17.9] (-0.8183,1.9778) -- (-0.9708,1.9448) -- (-1.0215,1.9990) -- (-0.8611,2.0337) -- cycle;
\path[cbpmesh=17.3] (0.9708,1.9448) -- (0.8183,1.9778) -- (0.8611,2.0337) -- (1.0215,1.9990) -- cycle;
\path[cbpmesh=15.6] (0.3019,1.8979) -- (0.1514,1.9063) -- (0.1594,1.9827) -- (0.3178,1.9739) -- cycle;
\path[cbpmesh=16.3] (-0.1514,1.9063) -- (-0.3019,1.8979) -- (-0.3178,1.9739) -- (-0.1594,1.9827) -- cycle;
\path[cbpmesh=33.3] (-0.4187,0.3816) -- (-0.5185,0.3643) -- (-0.5107,0.4996) -- (-0.4124,0.5166) -- cycle;
\path[cbpmesh=27.5] (0.5185,0.3643) -- (0.4187,0.3816) -- (0.4124,0.5166) -- (0.5107,0.4996) -- cycle;
\path[cbpmesh=28.5] (0.0000,0.5473) -- (-0.1047,0.5454) -- (-0.1041,0.6786) -- (0.0000,0.6805) -- cycle;
\path[cbpmesh=27.8] (0.1047,0.5454) -- (0.0000,0.5473) -- (0.0000,0.6805) -- (0.1041,0.6786) -- cycle;
\path[cbpmesh=18.2] (-3.1622,1.7301) -- (-3.3045,1.6108) -- (-3.3644,1.6180) -- (-3.2195,1.7394) -- cycle;
\path[cbpmesh=19.9] (3.3045,1.6108) -- (3.1622,1.7301) -- (3.2195,1.7394) -- (3.3644,1.6180) -- cycle;
\path[cbpmesh=38] (-1.2218,-0.1533) -- (-1.2999,-0.1998) -- (-1.2373,-0.1456) -- (-1.1629,-0.1014) -- cycle;
\draw[draw=cbpInk!70,line width=.45pt] (-1.2218,-0.1533) -- (-1.2999,-0.1998);
\path[cbpmesh=38] (1.2999,-0.1998) -- (1.2218,-0.1533) -- (1.1629,-0.1014) -- (1.2373,-0.1456) -- cycle;
\draw[draw=cbpInk!70,line width=.45pt] (1.2999,-0.1998) -- (1.2218,-0.1533);
\path[cbpmesh=15.5] (0.4505,1.8840) -- (0.3019,1.8979) -- (0.3178,1.9739) -- (0.4742,1.9593) -- cycle;
\path[cbpmesh=16.5] (-0.3019,1.8979) -- (-0.4505,1.8840) -- (-0.4742,1.9593) -- (-0.3178,1.9739) -- cycle;
\path[cbpmesh=37.3] (1.1629,-0.1014) -- (1.0814,-0.0600) -- (1.0331,0.0113) -- (1.1109,-0.0283) -- cycle;
\path[cbpmesh=38] (-1.0814,-0.0600) -- (-1.1629,-0.1014) -- (-1.1109,-0.0283) -- (-1.0331,0.0113) -- cycle;
\path[cbpmesh=16.2] (0.7774,1.9131) -- (0.6277,1.9390) -- (0.6608,2.0051) -- (0.8183,1.9778) -- cycle;
\path[cbpmesh=17.3] (-0.6277,1.9390) -- (-0.7774,1.9131) -- (-0.8183,1.9778) -- (-0.6608,2.0051) -- cycle;
\path[cbpmesh=27.2] (0.2088,0.5396) -- (0.1047,0.5454) -- (0.1041,0.6786) -- (0.2076,0.6729) -- cycle;
\path[cbpmesh=29.2] (-0.1047,0.5454) -- (-0.2088,0.5396) -- (-0.2076,0.6729) -- (-0.1041,0.6786) -- cycle;
\path[cbpmesh=19.0] (-2.2505,1.8062) -- (-2.3944,1.7206) -- (-2.4778,1.7250) -- (-2.3289,1.8137) -- cycle;
\path[cbpmesh=22.3] (2.3944,1.7206) -- (2.2505,1.8062) -- (2.3289,1.8137) -- (2.4778,1.7250) -- cycle;
\path[cbpmesh=20.2] (-3.6602,1.5785) -- (-3.7952,1.4410) -- (-3.8596,1.4089) -- (-3.7223,1.5488) -- cycle;
\path[cbpmesh=29.2] (3.7952,1.4410) -- (3.6602,1.5785) -- (3.7223,1.5488) -- (3.8596,1.4089) -- cycle;
\path[cbpmesh=30.7] (0.8244,0.1732) -- (0.7328,0.2019) -- (0.7080,0.3190) -- (0.7964,0.2913) -- cycle;
\path[cbpmesh=38] (-0.7328,0.2019) -- (-0.8244,0.1732) -- (-0.7964,0.2913) -- (-0.7080,0.3190) -- cycle;
\path[cbpmesh=18.3] (-2.6297,1.7464) -- (-2.7715,1.6468) -- (-2.8447,1.6592) -- (-2.6992,1.7615) -- cycle;
\path[cbpmesh=19.9] (2.7715,1.6468) -- (2.6297,1.7464) -- (2.6992,1.7615) -- (2.8447,1.6592) -- cycle;
\path[cbpmesh=20.0] (-1.9318,1.8817) -- (-2.0774,1.8078) -- (-2.1668,1.8040) -- (-2.0149,1.8811) -- cycle;
\path[cbpmesh=24.4] (2.0774,1.8078) -- (1.9318,1.8817) -- (2.0149,1.8811) -- (2.1668,1.8040) -- cycle;
\path[cbpmesh=27.2] (0.6152,0.3434) -- (0.5185,0.3643) -- (0.5107,0.4996) -- (0.6059,0.4790) -- cycle;
\path[cbpmesh=34.2] (-0.5185,0.3643) -- (-0.6152,0.3434) -- (-0.6059,0.4790) -- (-0.5107,0.4996) -- cycle;
\path[cbpmesh=30.0] (-0.2088,0.5396) -- (-0.3116,0.5300) -- (-0.3098,0.6633) -- (-0.2076,0.6729) -- cycle;
\path[cbpmesh=26.7] (0.3116,0.5300) -- (0.2088,0.5396) -- (0.2076,0.6729) -- (0.3098,0.6633) -- cycle;
\path[cbpmesh=38] (-0.9488,0.0478) -- (-1.0331,0.0113) -- (-0.9918,0.1061) -- (-0.9109,0.1412) -- cycle;
\path[cbpmesh=33.7] (1.0331,0.0113) -- (0.9488,0.0478) -- (0.9109,0.1412) -- (0.9918,0.1061) -- cycle;
\path[cbpmesh=18.2] (-3.1035,1.7179) -- (-3.2432,1.6008) -- (-3.3045,1.6108) -- (-3.1622,1.7301) -- cycle;
\path[cbpmesh=19.3] (3.2432,1.6008) -- (3.1035,1.7179) -- (3.1622,1.7301) -- (3.3045,1.6108) -- cycle;
\path[cbpmesh=16.8] (-0.4505,1.8840) -- (-0.5963,1.8647) -- (-0.6277,1.9390) -- (-0.4742,1.9593) -- cycle;
\path[cbpmesh=15.4] (0.5963,1.8647) -- (0.4505,1.8840) -- (0.4742,1.9593) -- (0.6277,1.9390) -- cycle;
\path[cbpmesh=15.4] (0.1439,1.8222) -- (0.0000,1.8248) -- (0.0000,1.9090) -- (0.1514,1.9063) -- cycle;
\path[cbpmesh=15.7] (0.0000,1.8248) -- (-0.1439,1.8222) -- (-0.1514,1.9063) -- (0.0000,1.9090) -- cycle;
\path[cbpmesh=30.8] (-0.3116,0.5300) -- (-0.4124,0.5166) -- (-0.4101,0.6500) -- (-0.3098,0.6633) -- cycle;
\path[cbpmesh=26.3] (0.4124,0.5166) -- (0.3116,0.5300) -- (0.3098,0.6633) -- (0.4101,0.6500) -- cycle;
\path[cbpmesh=18.8] (1.3226,1.9124) -- (1.1757,1.9585) -- (1.2362,2.0013) -- (1.3906,1.9528) -- cycle;
\path[cbpmesh=18.5] (-1.1757,1.9585) -- (-1.3226,1.9124) -- (-1.3906,1.9528) -- (-1.2362,2.0013) -- cycle;
\path[cbpmesh=17.9] (-0.9708,1.9448) -- (-1.1173,1.9063) -- (-1.1757,1.9585) -- (-1.0215,1.9990) -- cycle;
\path[cbpmesh=17.3] (1.1173,1.9063) -- (0.9708,1.9448) -- (1.0215,1.9990) -- (1.1757,1.9585) -- cycle;
\path[cbpmesh=15.2] (0.2869,1.8142) -- (0.1439,1.8222) -- (0.1514,1.9063) -- (0.3019,1.8979) -- cycle;
\path[cbpmesh=16.0] (-0.1439,1.8222) -- (-0.2869,1.8142) -- (-0.3019,1.8979) -- (-0.1514,1.9063) -- cycle;
\path[cbpmesh=20.4] (1.5365,1.8988) -- (1.3906,1.9528) -- (1.4608,1.9833) -- (1.6140,1.9266) -- cycle;
\path[cbpmesh=18.9] (-1.3906,1.9528) -- (-1.5365,1.8988) -- (-1.6140,1.9266) -- (-1.4608,1.9833) -- cycle;
\path[cbpmesh=19.5] (-3.5998,1.5996) -- (-3.7326,1.4644) -- (-3.7952,1.4410) -- (-3.6602,1.5785) -- cycle;
\path[cbpmesh=27.3] (3.7326,1.4644) -- (3.5998,1.5996) -- (3.6602,1.5785) -- (3.7952,1.4410) -- cycle;
\path[cbpmesh=17.5] (-0.7774,1.9131) -- (-0.9223,1.8817) -- (-0.9708,1.9448) -- (-0.8183,1.9778) -- cycle;
\path[cbpmesh=16.1] (0.9223,1.8817) -- (0.7774,1.9131) -- (0.8183,1.9778) -- (0.9708,1.9448) -- cycle;
\path[cbpmesh=18.4] (-2.5561,1.7340) -- (-2.6939,1.6371) -- (-2.7715,1.6468) -- (-2.6297,1.7464) -- cycle;
\path[cbpmesh=20.6] (2.6939,1.6371) -- (2.5561,1.7340) -- (2.6297,1.7464) -- (2.7715,1.6468) -- cycle;
\path[cbpmesh=27.0] (0.7080,0.3190) -- (0.6152,0.3434) -- (0.6059,0.4790) -- (0.6974,0.4550) -- cycle;
\path[cbpmesh=35.2] (-0.6152,0.3434) -- (-0.7080,0.3190) -- (-0.6974,0.4550) -- (-0.6059,0.4790) -- cycle;
\path[cbpmesh=18.2] (-3.0429,1.7037) -- (-3.1798,1.5889) -- (-3.2432,1.6008) -- (-3.1035,1.7179) -- cycle;
\path[cbpmesh=19.0] (3.1798,1.5889) -- (3.0429,1.7037) -- (3.1035,1.7179) -- (3.2432,1.6008) -- cycle;
\path[cbpmesh=38] (-0.8244,0.1732) -- (-0.9109,0.1412) -- (-0.8800,0.2603) -- (-0.7964,0.2913) -- cycle;
\path[cbpmesh=30.7] (0.9109,0.1412) -- (0.8244,0.1732) -- (0.7964,0.2913) -- (0.8800,0.2603) -- cycle;
\path[cbpmesh=25.9] (0.5107,0.4996) -- (0.4124,0.5166) -- (0.4101,0.6500) -- (0.5079,0.6331) -- cycle;
\path[cbpmesh=31.7] (-0.4124,0.5166) -- (-0.5107,0.4996) -- (-0.5079,0.6331) -- (-0.4101,0.6500) -- cycle;
\path[cbpmesh=26.2] (0.1041,0.6786) -- (0.0000,0.6805) -- (0.0000,0.8119) -- (0.1045,0.8100) -- cycle;
\path[cbpmesh=26.9] (0.0000,0.6805) -- (-0.1041,0.6786) -- (-0.1045,0.8100) -- (0.0000,0.8119) -- cycle;
\path[cbpmesh=15.0] (0.4281,1.8011) -- (0.2869,1.8142) -- (0.3019,1.8979) -- (0.4505,1.8840) -- cycle;
\path[cbpmesh=16.4] (-0.2869,1.8142) -- (-0.4281,1.8011) -- (-0.4505,1.8840) -- (-0.3019,1.8979) -- cycle;
\path[cbpmesh=19.4] (-1.6140,1.9266) -- (-1.7573,1.8644) -- (-1.8437,1.8785) -- (-1.6934,1.9438) -- cycle;
\path[cbpmesh=22.1] (1.7573,1.8644) -- (1.6140,1.9266) -- (1.6934,1.9438) -- (1.8437,1.8785) -- cycle;
\path[cbpmesh=27.6] (-0.1041,0.6786) -- (-0.2076,0.6729) -- (-0.2084,0.8042) -- (-0.1045,0.8100) -- cycle;
\path[cbpmesh=25.6] (0.2076,0.6729) -- (0.1041,0.6786) -- (0.1045,0.8100) -- (0.2084,0.8042) -- cycle;
\path[cbpmesh=15.2] (0.7385,1.8401) -- (0.5963,1.8647) -- (0.6277,1.9390) -- (0.7774,1.9131) -- cycle;
\path[cbpmesh=17.0] (-0.5963,1.8647) -- (-0.7385,1.8401) -- (-0.7774,1.9131) -- (-0.6277,1.9390) -- cycle;
\path[cbpmesh=19.3] (-2.1668,1.8040) -- (-2.3053,1.7215) -- (-2.3944,1.7206) -- (-2.2505,1.8062) -- cycle;
\path[cbpmesh=23.3] (2.3053,1.7215) -- (2.1668,1.8040) -- (2.2505,1.8062) -- (2.3944,1.7206) -- cycle;
\path[cbpmesh=37.5] (1.2373,-0.1456) -- (1.1629,-0.1014) -- (1.1109,-0.0283) -- (1.1820,-0.0706) -- cycle;
\path[cbpmesh=38] (-1.1629,-0.1014) -- (-1.2373,-0.1456) -- (-1.1820,-0.0706) -- (-1.1109,-0.0283) -- cycle;
\path[cbpmesh=38] (-1.2999,-0.1998) -- (-1.3700,-0.2490) -- (-1.3040,-0.1925) -- (-1.2373,-0.1456) -- cycle;
\draw[draw=cbpInk!70,line width=.45pt] (-1.2999,-0.1998) -- (-1.3700,-0.2490);
\path[cbpmesh=38] (1.3700,-0.2490) -- (1.2999,-0.1998) -- (1.2373,-0.1456) -- (1.3040,-0.1925) -- cycle;
\draw[draw=cbpInk!70,line width=.45pt] (1.3700,-0.2490) -- (1.2999,-0.1998);
\path[cbpmesh=25.1] (0.3098,0.6633) -- (0.2076,0.6729) -- (0.2084,0.8042) -- (0.3110,0.7946) -- cycle;
\path[cbpmesh=28.4] (-0.2076,0.6729) -- (-0.3098,0.6633) -- (-0.3110,0.7946) -- (-0.2084,0.8042) -- cycle;
\path[cbpmesh=18.2] (-2.9799,1.6886) -- (-3.1140,1.5762) -- (-3.1798,1.5889) -- (-3.0429,1.7037) -- cycle;
\path[cbpmesh=18.9] (3.1140,1.5762) -- (2.9799,1.6886) -- (3.0429,1.7037) -- (3.1798,1.5889) -- cycle;
\path[cbpmesh=18.9] (-3.5405,1.6132) -- (-3.6711,1.4802) -- (-3.7326,1.4644) -- (-3.5998,1.5996) -- cycle;
\path[cbpmesh=25.4] (3.6711,1.4802) -- (3.5405,1.6132) -- (3.5998,1.5996) -- (3.7326,1.4644) -- cycle;
\path[cbpmesh=33.8] (1.1109,-0.0283) -- (1.0331,0.0113) -- (0.9918,0.1061) -- (1.0665,0.0681) -- cycle;
\path[cbpmesh=38] (-1.0331,0.0113) -- (-1.1109,-0.0283) -- (-1.0665,0.0681) -- (-0.9918,0.1061) -- cycle;
\path[cbpmesh=15.4] (0.1369,1.7311) -- (0.0000,1.7336) -- (0.0000,1.8248) -- (0.1439,1.8222) -- cycle;
\path[cbpmesh=15.7] (0.0000,1.7336) -- (-0.1369,1.7311) -- (-0.1439,1.8222) -- (0.0000,1.8248) -- cycle;
\path[cbpmesh=16.7] (-0.4281,1.8011) -- (-0.5667,1.7827) -- (-0.5963,1.8647) -- (-0.4505,1.8840) -- cycle;
\path[cbpmesh=14.8] (0.5667,1.7827) -- (0.4281,1.8011) -- (0.4505,1.8840) -- (0.5963,1.8647) -- cycle;
\path[cbpmesh=25.6] (0.6059,0.4790) -- (0.5107,0.4996) -- (0.5079,0.6331) -- (0.6025,0.6126) -- cycle;
\path[cbpmesh=32.7] (-0.5107,0.4996) -- (-0.6059,0.4790) -- (-0.6025,0.6126) -- (-0.5079,0.6331) -- cycle;
\path[cbpmesh=19.8] (-1.8437,1.8785) -- (-1.9827,1.8080) -- (-2.0774,1.8078) -- (-1.9318,1.8817) -- cycle;
\path[cbpmesh=23.8] (1.9827,1.8080) -- (1.8437,1.8785) -- (1.9318,1.8817) -- (2.0774,1.8078) -- cycle;
\path[cbpmesh=15.1] (0.2729,1.7235) -- (0.1369,1.7311) -- (0.1439,1.8222) -- (0.2869,1.8142) -- cycle;
\path[cbpmesh=16.1] (-0.1369,1.7311) -- (-0.2729,1.7235) -- (-0.2869,1.8142) -- (-0.1439,1.8222) -- cycle;
\path[cbpmesh=26.9] (0.7964,0.2913) -- (0.7080,0.3190) -- (0.6974,0.4550) -- (0.7845,0.4276) -- cycle;
\path[cbpmesh=36.2] (-0.7080,0.3190) -- (-0.7964,0.2913) -- (-0.7845,0.4276) -- (-0.6974,0.4550) -- cycle;
\path[cbpmesh=24.7] (0.4101,0.6500) -- (0.3098,0.6633) -- (0.3110,0.7946) -- (0.4116,0.7813) -- cycle;
\path[cbpmesh=29.2] (-0.3098,0.6633) -- (-0.4101,0.6500) -- (-0.4116,0.7813) -- (-0.3110,0.7946) -- cycle;
\path[cbpmesh=18.5] (-2.4778,1.7250) -- (-2.6114,1.6312) -- (-2.6939,1.6371) -- (-2.5561,1.7340) -- cycle;
\path[cbpmesh=21.4] (2.6114,1.6312) -- (2.4778,1.7250) -- (2.5561,1.7340) -- (2.6939,1.6371) -- cycle;
\path[cbpmesh=17.6] (-0.9223,1.8817) -- (-1.0614,1.8451) -- (-1.1173,1.9063) -- (-0.9708,1.9448) -- cycle;
\path[cbpmesh=16.1] (1.0614,1.8451) -- (0.9223,1.8817) -- (0.9708,1.9448) -- (1.1173,1.9063) -- cycle;
\path[cbpmesh=24.6] (0.1045,0.8100) -- (0.0000,0.8119) -- (0.0000,0.9409) -- (0.1058,0.9389) -- cycle;
\path[cbpmesh=25.2] (0.0000,0.8119) -- (-0.1045,0.8100) -- (-0.1058,0.9389) -- (0.0000,0.9409) -- cycle;
\path[cbpmesh=18.0] (-1.1173,1.9063) -- (-1.2569,1.8625) -- (-1.3226,1.9124) -- (-1.1757,1.9585) -- cycle;
\path[cbpmesh=17.3] (1.2569,1.8625) -- (1.1173,1.9063) -- (1.1757,1.9585) -- (1.3226,1.9124) -- cycle;
\path[cbpmesh=14.8] (0.4072,1.7110) -- (0.2729,1.7235) -- (0.2869,1.8142) -- (0.4281,1.8011) -- cycle;
\path[cbpmesh=16.5] (-0.2729,1.7235) -- (-0.4072,1.7110) -- (-0.4281,1.8011) -- (-0.2869,1.8142) -- cycle;
\path[cbpmesh=30.8] (0.9918,0.1061) -- (0.9109,0.1412) -- (0.8800,0.2603) -- (0.9581,0.2264) -- cycle;
\path[cbpmesh=38] (-0.9109,0.1412) -- (-0.9918,0.1061) -- (-0.9581,0.2264) -- (-0.8800,0.2603) -- cycle;
\path[cbpmesh=18.2] (-2.9140,1.6735) -- (-3.0451,1.5635) -- (-3.1140,1.5762) -- (-2.9799,1.6886) -- cycle;
\path[cbpmesh=19.0] (3.0451,1.5635) -- (2.9140,1.6735) -- (2.9799,1.6886) -- (3.1140,1.5762) -- cycle;
\path[cbpmesh=17.3] (-0.7385,1.8401) -- (-0.8761,1.8103) -- (-0.9223,1.8817) -- (-0.7774,1.9131) -- cycle;
\path[cbpmesh=15.2] (0.8761,1.8103) -- (0.7385,1.8401) -- (0.7774,1.9131) -- (0.9223,1.8817) -- cycle;
\path[cbpmesh=24.0] (0.2084,0.8042) -- (0.1045,0.8100) -- (0.1058,0.9389) -- (0.2110,0.9331) -- cycle;
\path[cbpmesh=26.0] (-0.1045,0.8100) -- (-0.2084,0.8042) -- (-0.2110,0.9331) -- (-0.1058,0.9389) -- cycle;
\path[cbpmesh=18.4] (-3.4818,1.6202) -- (-3.6103,1.4893) -- (-3.6711,1.4802) -- (-3.5405,1.6132) -- cycle;
\path[cbpmesh=23.6] (3.6103,1.4893) -- (3.4818,1.6202) -- (3.5405,1.6132) -- (3.6711,1.4802) -- cycle;
\path[cbpmesh=18.8] (1.4613,1.8611) -- (1.3226,1.9124) -- (1.3906,1.9528) -- (1.5365,1.8988) -- cycle;
\path[cbpmesh=18.4] (-1.3226,1.9124) -- (-1.4613,1.8611) -- (-1.5365,1.8988) -- (-1.3906,1.9528) -- cycle;
\path[cbpmesh=24.3] (0.5079,0.6331) -- (0.4101,0.6500) -- (0.4116,0.7813) -- (0.5098,0.7643) -- cycle;
\path[cbpmesh=30.1] (-0.4101,0.6500) -- (-0.5079,0.6331) -- (-0.5098,0.7643) -- (-0.4116,0.7813) -- cycle;
\path[cbpmesh=14.6] (0.7018,1.7593) -- (0.5667,1.7827) -- (0.5963,1.8647) -- (0.7385,1.8401) -- cycle;
\path[cbpmesh=17.1] (-0.5667,1.7827) -- (-0.7018,1.7593) -- (-0.7385,1.8401) -- (-0.5963,1.8647) -- cycle;
\path[cbpmesh=15.7] (0.1304,1.6334) -- (0.0000,1.6358) -- (0.0000,1.7336) -- (0.1369,1.7311) -- cycle;
\path[cbpmesh=16.1] (0.0000,1.6358) -- (-0.1304,1.6334) -- (-0.1369,1.7311) -- (0.0000,1.7336) -- cycle;
\path[cbpmesh=25.4] (0.6974,0.4550) -- (0.6059,0.4790) -- (0.6025,0.6126) -- (0.6935,0.5887) -- cycle;
\path[cbpmesh=33.6] (-0.6059,0.4790) -- (-0.6974,0.4550) -- (-0.6935,0.5887) -- (-0.6025,0.6126) -- cycle;
\path[cbpmesh=26.7] (-0.2084,0.8042) -- (-0.3110,0.7946) -- (-0.3149,0.9234) -- (-0.2110,0.9331) -- cycle;
\path[cbpmesh=23.5] (0.3110,0.7946) -- (0.2084,0.8042) -- (0.2110,0.9331) -- (0.3149,0.9234) -- cycle;
\path[cbpmesh=15.3] (0.2600,1.6262) -- (0.1304,1.6334) -- (0.1369,1.7311) -- (0.2729,1.7235) -- cycle;
\path[cbpmesh=16.6] (-0.1304,1.6334) -- (-0.2600,1.6262) -- (-0.2729,1.7235) -- (-0.1369,1.7311) -- cycle;
\path[cbpmesh=14.6] (0.5391,1.6935) -- (0.4072,1.7110) -- (0.4281,1.8011) -- (0.5667,1.7827) -- cycle;
\path[cbpmesh=17.0] (-0.4072,1.7110) -- (-0.5391,1.6935) -- (-0.5667,1.7827) -- (-0.4281,1.8011) -- cycle;
\path[cbpmesh=19.6] (-2.0774,1.8078) -- (-2.2102,1.7287) -- (-2.3053,1.7215) -- (-2.1668,1.8040) -- cycle;
\path[cbpmesh=24.3] (2.2102,1.7287) -- (2.0774,1.8078) -- (2.1668,1.8040) -- (2.3053,1.7215) -- cycle;
\path[cbpmesh=20.4] (1.6729,1.8396) -- (1.5365,1.8988) -- (1.6140,1.9266) -- (1.7573,1.8644) -- cycle;
\path[cbpmesh=18.8] (-1.5365,1.8988) -- (-1.6729,1.8396) -- (-1.7573,1.8644) -- (-1.6140,1.9266) -- cycle;
\path[cbpmesh=37.2] (-0.7964,0.2913) -- (-0.8800,0.2603) -- (-0.8668,0.3972) -- (-0.7845,0.4276) -- cycle;
\path[cbpmesh=26.8] (0.8800,0.2603) -- (0.7964,0.2913) -- (0.7845,0.4276) -- (0.8668,0.3972) -- cycle;
\path[cbpmesh=18.2] (-2.8447,1.6592) -- (-2.9727,1.5519) -- (-3.0451,1.5635) -- (-2.9140,1.6735) -- cycle;
\path[cbpmesh=19.4] (2.9727,1.5519) -- (2.8447,1.6592) -- (2.9140,1.6735) -- (3.0451,1.5635) -- cycle;
\path[cbpmesh=38] (-1.2373,-0.1456) -- (-1.3040,-0.1925) -- (-1.2457,-0.1153) -- (-1.1820,-0.0706) -- cycle;
\path[cbpmesh=37.7] (1.3040,-0.1925) -- (1.2373,-0.1456) -- (1.1820,-0.0706) -- (1.2457,-0.1153) -- cycle;
\path[cbpmesh=18.2] (-3.4233,1.6214) -- (-3.5496,1.4928) -- (-3.6103,1.4893) -- (-3.4818,1.6202) -- cycle;
\path[cbpmesh=22.0] (3.5496,1.4928) -- (3.4233,1.6214) -- (3.4818,1.6202) -- (3.6103,1.4893) -- cycle;
\path[cbpmesh=22.9] (0.1058,0.9389) -- (0.0000,0.9409) -- (0.0000,1.0670) -- (0.1080,1.0650) -- cycle;
\path[cbpmesh=23.6] (0.0000,0.9409) -- (-0.1058,0.9389) -- (-0.1080,1.0650) -- (0.0000,1.0670) -- cycle;
\path[cbpmesh=23.0] (0.4116,0.7813) -- (0.3110,0.7946) -- (0.3149,0.9234) -- (0.4168,0.9099) -- cycle;
\path[cbpmesh=27.6] (-0.3110,0.7946) -- (-0.4116,0.7813) -- (-0.4168,0.9099) -- (-0.3149,0.9234) -- cycle;
\path[cbpmesh=21.1] (-3.9263,1.3673) -- (-4.0418,1.2204) -- (-4.1133,1.1656) -- (-3.9958,1.3152) -- cycle;
\draw[draw=cbpInk!70,line width=.45pt] (-4.1133,1.1656) -- (-3.9958,1.3152);
\path[cbpmesh=32.2] (4.0418,1.2204) -- (3.9263,1.3673) -- (3.9958,1.3152) -- (4.1133,1.1656) -- cycle;
\draw[draw=cbpInk!70,line width=.45pt] (3.9958,1.3152) -- (4.1133,1.1656);
\path[cbpmesh=38] (-1.1109,-0.0283) -- (-1.1820,-0.0706) -- (-1.1347,0.0275) -- (-1.0665,0.0681) -- cycle;
\path[cbpmesh=34.0] (1.1820,-0.0706) -- (1.1109,-0.0283) -- (1.0665,0.0681) -- (1.1347,0.0275) -- cycle;
\path[cbpmesh=15.0] (0.3880,1.6142) -- (0.2600,1.6262) -- (0.2729,1.7235) -- (0.4072,1.7110) -- cycle;
\path[cbpmesh=17.1] (-0.2600,1.6262) -- (-0.3880,1.6142) -- (-0.4072,1.7110) -- (-0.2729,1.7235) -- cycle;
\path[cbpmesh=18.7] (-2.3944,1.7206) -- (-2.5235,1.6299) -- (-2.6114,1.6312) -- (-2.4778,1.7250) -- cycle;
\path[cbpmesh=22.3] (2.5235,1.6299) -- (2.3944,1.7206) -- (2.4778,1.7250) -- (2.6114,1.6312) -- cycle;
\path[cbpmesh=24.0] (0.6025,0.6126) -- (0.5079,0.6331) -- (0.5098,0.7643) -- (0.6048,0.7437) -- cycle;
\path[cbpmesh=31.0] (-0.5079,0.6331) -- (-0.6025,0.6126) -- (-0.6048,0.7437) -- (-0.5098,0.7643) -- cycle;
\path[cbpmesh=22.4] (0.2110,0.9331) -- (0.1058,0.9389) -- (0.1080,1.0650) -- (0.2153,1.0590) -- cycle;
\path[cbpmesh=24.3] (-0.1058,0.9389) -- (-0.2110,0.9331) -- (-0.2153,1.0590) -- (-0.1080,1.0650) -- cycle;
\path[cbpmesh=38] (-1.3700,-0.2490) -- (-1.4317,-0.3007) -- (-1.3627,-0.2417) -- (-1.3040,-0.1925) -- cycle;
\draw[draw=cbpInk!70,line width=.45pt] (-1.3700,-0.2490) -- (-1.4317,-0.3007);
\draw[draw=cbpInk!70,line width=.45pt] (-1.4317,-0.3007) -- (-1.3627,-0.2417);
\path[cbpmesh=38] (1.4317,-0.3007) -- (1.3700,-0.2490) -- (1.3040,-0.1925) -- (1.3627,-0.2417) -- cycle;
\draw[draw=cbpInk!70,line width=.45pt] (1.4317,-0.3007) -- (1.3700,-0.2490);
\draw[draw=cbpInk!70,line width=.45pt] (1.3627,-0.2417) -- (1.4317,-0.3007);
\path[cbpmesh=16.3] (0.1245,1.5297) -- (0.0000,1.5320) -- (0.0000,1.6358) -- (0.1304,1.6334) -- cycle;
\path[cbpmesh=16.8] (0.0000,1.5320) -- (-0.1245,1.5297) -- (-0.1304,1.6334) -- (0.0000,1.6358) -- cycle;
\path[cbpmesh=15.9] (0.2483,1.5228) -- (0.1245,1.5297) -- (0.1304,1.6334) -- (0.2600,1.6262) -- cycle;
\path[cbpmesh=17.3] (-0.1245,1.5297) -- (-0.2483,1.5228) -- (-0.2600,1.6262) -- (-0.1304,1.6334) -- cycle;
\path[cbpmesh=25.1] (-0.2110,0.9331) -- (-0.3149,0.9234) -- (-0.3213,1.0491) -- (-0.2153,1.0590) -- cycle;
\path[cbpmesh=21.9] (0.3149,0.9234) -- (0.2110,0.9331) -- (0.2153,1.0590) -- (0.3213,1.0491) -- cycle;
\path[cbpmesh=14.5] (0.8326,1.7310) -- (0.7018,1.7593) -- (0.7385,1.8401) -- (0.8761,1.8103) -- cycle;
\path[cbpmesh=17.5] (-0.7018,1.7593) -- (-0.8326,1.7310) -- (-0.8761,1.8103) -- (-0.7385,1.8401) -- cycle;
\path[cbpmesh=17.6] (-0.8761,1.8103) -- (-1.0083,1.7755) -- (-1.0614,1.8451) -- (-0.9223,1.8817) -- cycle;
\path[cbpmesh=15.1] (1.0083,1.7755) -- (0.8761,1.8103) -- (0.9223,1.8817) -- (1.0614,1.8451) -- cycle;
\path[cbpmesh=25.3] (0.7845,0.4276) -- (0.6974,0.4550) -- (0.6935,0.5887) -- (0.7801,0.5615) -- cycle;
\path[cbpmesh=34.6] (-0.6974,0.4550) -- (-0.7845,0.4276) -- (-0.7801,0.5615) -- (-0.6935,0.5887) -- cycle;
\path[cbpmesh=14.4] (0.6676,1.6712) -- (0.5391,1.6935) -- (0.5667,1.7827) -- (0.7018,1.7593) -- cycle;
\path[cbpmesh=17.5] (-0.5391,1.6935) -- (-0.6676,1.6712) -- (-0.7018,1.7593) -- (-0.5667,1.7827) -- cycle;
\path[cbpmesh=38] (-0.9918,0.1061) -- (-1.0665,0.0681) -- (-1.0304,0.1898) -- (-0.9581,0.2264) -- cycle;
\path[cbpmesh=30.9] (1.0665,0.0681) -- (0.9918,0.1061) -- (0.9581,0.2264) -- (1.0304,0.1898) -- cycle;
\path[cbpmesh=19.1] (-1.7573,1.8644) -- (-1.8898,1.7971) -- (-1.9827,1.8080) -- (-1.8437,1.8785) -- cycle;
\path[cbpmesh=22.1] (1.8898,1.7971) -- (1.7573,1.8644) -- (1.8437,1.8785) -- (1.9827,1.8080) -- cycle;
\path[cbpmesh=22.0] (0.0000,1.0670) -- (-0.1080,1.0650) -- (-0.1110,1.1876) -- (0.0000,1.1897) -- cycle;
\path[cbpmesh=21.3] (0.1080,1.0650) -- (0.0000,1.0670) -- (0.0000,1.1897) -- (0.1110,1.1876) -- cycle;
\path[cbpmesh=22.7] (0.5098,0.7643) -- (0.4116,0.7813) -- (0.4168,0.9099) -- (0.5161,0.8927) -- cycle;
\path[cbpmesh=28.4] (-0.4116,0.7813) -- (-0.5098,0.7643) -- (-0.5161,0.8927) -- (-0.4168,0.9099) -- cycle;
\path[cbpmesh=18.1] (-3.3644,1.6180) -- (-3.4885,1.4916) -- (-3.5496,1.4928) -- (-3.4233,1.6214) -- cycle;
\path[cbpmesh=20.8] (3.4885,1.4916) -- (3.3644,1.6180) -- (3.4233,1.6214) -- (3.5496,1.4928) -- cycle;
\path[cbpmesh=16.1] (1.1940,1.8035) -- (1.0614,1.8451) -- (1.1173,1.9063) -- (1.2569,1.8625) -- cycle;
\path[cbpmesh=17.8] (-1.0614,1.8451) -- (-1.1940,1.8035) -- (-1.2569,1.8625) -- (-1.1173,1.9063) -- cycle;
\path[cbpmesh=17.6] (-0.3880,1.6142) -- (-0.5136,1.5976) -- (-0.5391,1.6935) -- (-0.4072,1.7110) -- cycle;
\path[cbpmesh=14.7] (0.5136,1.5976) -- (0.3880,1.6142) -- (0.4072,1.7110) -- (0.5391,1.6935) -- cycle;
\path[cbpmesh=17.3] (0.1193,1.4205) -- (0.0000,1.4227) -- (0.0000,1.5320) -- (0.1245,1.5297) -- cycle;
\path[cbpmesh=17.8] (0.0000,1.4227) -- (-0.1193,1.4205) -- (-0.1245,1.5297) -- (0.0000,1.5320) -- cycle;
\path[cbpmesh=20.8] (0.2153,1.0590) -- (0.1080,1.0650) -- (0.1110,1.1876) -- (0.2213,1.1815) -- cycle;
\path[cbpmesh=22.7] (-0.1080,1.0650) -- (-0.2153,1.0590) -- (-0.2213,1.1815) -- (-0.1110,1.1876) -- cycle;
\path[cbpmesh=18.2] (-2.7715,1.6468) -- (-2.8962,1.5422) -- (-2.9727,1.5519) -- (-2.8447,1.6592) -- cycle;
\path[cbpmesh=19.9] (2.8962,1.5422) -- (2.7715,1.6468) -- (2.8447,1.6592) -- (2.9727,1.5519) -- cycle;
\path[cbpmesh=15.5] (0.3706,1.5114) -- (0.2483,1.5228) -- (0.2600,1.6262) -- (0.3880,1.6142) -- cycle;
\path[cbpmesh=17.9] (-0.2483,1.5228) -- (-0.3706,1.5114) -- (-0.3880,1.6142) -- (-0.2600,1.6262) -- cycle;
\path[cbpmesh=16.8] (0.2379,1.4139) -- (0.1193,1.4205) -- (0.1245,1.5297) -- (0.2483,1.5228) -- cycle;
\path[cbpmesh=18.4] (-0.1193,1.4205) -- (-0.2379,1.4139) -- (-0.2483,1.5228) -- (-0.1245,1.5297) -- cycle;
\path[cbpmesh=19.8] (0.1110,1.1876) -- (0.0000,1.1897) -- (0.0000,1.3084) -- (0.1148,1.3063) -- cycle;
\path[cbpmesh=20.4] (0.0000,1.1897) -- (-0.1110,1.1876) -- (-0.1148,1.3063) -- (0.0000,1.3084) -- cycle;
\path[cbpmesh=21.4] (0.4168,0.9099) -- (0.3149,0.9234) -- (0.3213,1.0491) -- (0.4254,1.0354) -- cycle;
\path[cbpmesh=25.9] (-0.3149,0.9234) -- (-0.4168,0.9099) -- (-0.4254,1.0354) -- (-0.3213,1.0491) -- cycle;
\path[cbpmesh=18.4] (0.1148,1.3063) -- (0.0000,1.3084) -- (0.0000,1.4227) -- (0.1193,1.4205) -- cycle;
\path[cbpmesh=19.0] (0.0000,1.3084) -- (-0.1148,1.3063) -- (-0.1193,1.4205) -- (0.0000,1.4227) -- cycle;
\path[cbpmesh=18.1] (-1.2569,1.8625) -- (-1.3888,1.8137) -- (-1.4613,1.8611) -- (-1.3226,1.9124) -- cycle;
\path[cbpmesh=17.3] (1.3888,1.8137) -- (1.2569,1.8625) -- (1.3226,1.9124) -- (1.4613,1.8611) -- cycle;
\path[cbpmesh=23.8] (0.6935,0.5887) -- (0.6025,0.6126) -- (0.6048,0.7437) -- (0.6960,0.7198) -- cycle;
\path[cbpmesh=32.0] (-0.6025,0.6126) -- (-0.6935,0.5887) -- (-0.6960,0.7198) -- (-0.6048,0.7437) -- cycle;
\path[cbpmesh=20.3] (0.3213,1.0491) -- (0.2153,1.0590) -- (0.2213,1.1815) -- (0.3303,1.1713) -- cycle;
\path[cbpmesh=23.4] (-0.2153,1.0590) -- (-0.3213,1.0491) -- (-0.3303,1.1713) -- (-0.2213,1.1815) -- cycle;
\path[cbpmesh=19.3] (0.2213,1.1815) -- (0.1110,1.1876) -- (0.1148,1.3063) -- (0.2289,1.2999) -- cycle;
\path[cbpmesh=21.1] (-0.1110,1.1876) -- (-0.2213,1.1815) -- (-0.2289,1.2999) -- (-0.1148,1.3063) -- cycle;
\path[cbpmesh=17.9] (0.2289,1.2999) -- (0.1148,1.3063) -- (0.1193,1.4205) -- (0.2379,1.4139) -- cycle;
\path[cbpmesh=19.6] (-0.1148,1.3063) -- (-0.2289,1.2999) -- (-0.2379,1.4139) -- (-0.1193,1.4205) -- cycle;
\path[cbpmesh=38] (-0.8800,0.2603) -- (-0.9581,0.2264) -- (-0.9437,0.3638) -- (-0.8668,0.3972) -- cycle;
\path[cbpmesh=26.9] (0.9581,0.2264) -- (0.8800,0.2603) -- (0.8668,0.3972) -- (0.9437,0.3638) -- cycle;
\path[cbpmesh=20.3] (-3.8596,1.4089) -- (-3.9731,1.2645) -- (-4.0418,1.2204) -- (-3.9263,1.3673) -- cycle;
\path[cbpmesh=30.6] (3.9731,1.2645) -- (3.8596,1.4089) -- (3.9263,1.3673) -- (4.0418,1.2204) -- cycle;
\path[cbpmesh=18.1] (-3.3045,1.6108) -- (-3.4264,1.4866) -- (-3.4885,1.4916) -- (-3.3644,1.6180) -- cycle;
\path[cbpmesh=19.9] (3.4264,1.4866) -- (3.3045,1.6108) -- (3.3644,1.6180) -- (3.4885,1.4916) -- cycle;
\path[cbpmesh=16.4] (0.3550,1.4029) -- (0.2379,1.4139) -- (0.2483,1.5228) -- (0.3706,1.5114) -- cycle;
\path[cbpmesh=19.0] (-0.2379,1.4139) -- (-0.3550,1.4029) -- (-0.3706,1.5114) -- (-0.2483,1.5228) -- cycle;
\path[cbpmesh=18.5] (-0.3706,1.5114) -- (-0.4905,1.4955) -- (-0.5136,1.5976) -- (-0.3880,1.6142) -- cycle;
\path[cbpmesh=15.2] (0.4905,1.4955) -- (0.3706,1.5114) -- (0.3880,1.6142) -- (0.5136,1.5976) -- cycle;
\path[cbpmesh=22.4] (0.6048,0.7437) -- (0.5098,0.7643) -- (0.5161,0.8927) -- (0.6123,0.8719) -- cycle;
\path[cbpmesh=29.4] (-0.5098,0.7643) -- (-0.6048,0.7437) -- (-0.6123,0.8719) -- (-0.5161,0.8927) -- cycle;
\path[cbpmesh=18.8] (0.3303,1.1713) -- (0.2213,1.1815) -- (0.2289,1.2999) -- (0.3416,1.2894) -- cycle;
\path[cbpmesh=21.8] (-0.2213,1.1815) -- (-0.3303,1.1713) -- (-0.3416,1.2894) -- (-0.2289,1.2999) -- cycle;
\path[cbpmesh=18.1] (-0.5136,1.5976) -- (-0.6361,1.5764) -- (-0.6676,1.6712) -- (-0.5391,1.6935) -- cycle;
\path[cbpmesh=14.5] (0.6361,1.5764) -- (0.5136,1.5976) -- (0.5391,1.6935) -- (0.6676,1.6712) -- cycle;
\path[cbpmesh=20.3] (-0.2289,1.2999) -- (-0.3416,1.2894) -- (-0.3550,1.4029) -- (-0.2379,1.4139) -- cycle;
\path[cbpmesh=17.5] (0.3416,1.2894) -- (0.2289,1.2999) -- (0.2379,1.4139) -- (0.3550,1.4029) -- cycle;
\path[cbpmesh=18.9] (-2.3053,1.7215) -- (-2.4296,1.6341) -- (-2.5235,1.6299) -- (-2.3944,1.7206) -- cycle;
\path[cbpmesh=23.2] (2.4296,1.6341) -- (2.3053,1.7215) -- (2.3944,1.7206) -- (2.5235,1.6299) -- cycle;
\path[cbpmesh=26.8] (-0.4168,0.9099) -- (-0.5161,0.8927) -- (-0.5268,1.0178) -- (-0.4254,1.0354) -- cycle;
\path[cbpmesh=21.1] (0.5161,0.8927) -- (0.4168,0.9099) -- (0.4254,1.0354) -- (0.5268,1.0178) -- cycle;
\path[cbpmesh=18.0] (-0.6676,1.6712) -- (-0.7920,1.6443) -- (-0.8326,1.7310) -- (-0.7018,1.7593) -- cycle;
\path[cbpmesh=14.2] (0.7920,1.6443) -- (0.6676,1.6712) -- (0.7018,1.7593) -- (0.8326,1.7310) -- cycle;
\path[cbpmesh=18.4] (-1.4613,1.8611) -- (-1.5911,1.8048) -- (-1.6729,1.8396) -- (-1.5365,1.8988) -- cycle;
\path[cbpmesh=18.8] (1.5911,1.8048) -- (1.4613,1.8611) -- (1.5365,1.8988) -- (1.6729,1.8396) -- cycle;
\path[cbpmesh=24.2] (-0.3213,1.0491) -- (-0.4254,1.0354) -- (-0.4372,1.1571) -- (-0.3303,1.1713) -- cycle;
\path[cbpmesh=19.9] (0.4254,1.0354) -- (0.3213,1.0491) -- (0.3303,1.1713) -- (0.4372,1.1571) -- cycle;
\path[cbpmesh=19.4] (-1.9827,1.8080) -- (-2.1095,1.7325) -- (-2.2102,1.7287) -- (-2.0774,1.8078) -- cycle;
\path[cbpmesh=23.8] (2.1095,1.7325) -- (1.9827,1.8080) -- (2.0774,1.8078) -- (2.2102,1.7287) -- cycle;
\path[cbpmesh=25.3] (0.8668,0.3972) -- (0.7845,0.4276) -- (0.7801,0.5615) -- (0.8620,0.5313) -- cycle;
\path[cbpmesh=35.7] (-0.7845,0.4276) -- (-0.8668,0.3972) -- (-0.8620,0.5313) -- (-0.7801,0.5615) -- cycle;
\path[cbpmesh=18.2] (-2.6939,1.6371) -- (-2.8151,1.5355) -- (-2.8962,1.5422) -- (-2.7715,1.6468) -- cycle;
\path[cbpmesh=20.5] (2.8151,1.5355) -- (2.6939,1.6371) -- (2.7715,1.6468) -- (2.8962,1.5422) -- cycle;
\path[cbpmesh=17.9] (-0.8326,1.7310) -- (-0.9582,1.6980) -- (-1.0083,1.7755) -- (-0.8761,1.8103) -- cycle;
\path[cbpmesh=14.4] (0.9582,1.6980) -- (0.8326,1.7310) -- (0.8761,1.8103) -- (1.0083,1.7755) -- cycle;
\path[cbpmesh=38] (-1.1820,-0.0706) -- (-1.2457,-0.1153) -- (-1.1959,-0.0155) -- (-1.1347,0.0275) -- cycle;
\path[cbpmesh=34.3] (1.2457,-0.1153) -- (1.1820,-0.0706) -- (1.1347,0.0275) -- (1.1959,-0.0155) -- cycle;
\path[cbpmesh=19.7] (-0.3550,1.4029) -- (-0.4700,1.3877) -- (-0.4905,1.4955) -- (-0.3706,1.5114) -- cycle;
\path[cbpmesh=16.0] (0.4700,1.3877) -- (0.3550,1.4029) -- (0.3706,1.5114) -- (0.4905,1.4955) -- cycle;
\path[cbpmesh=38] (-1.3040,-0.1925) -- (-1.3627,-0.2417) -- (-1.3018,-0.1623) -- (-1.2457,-0.1153) -- cycle;
\path[cbpmesh=38] (1.3627,-0.2417) -- (1.3040,-0.1925) -- (1.2457,-0.1153) -- (1.3018,-0.1623) -- cycle;
\path[cbpmesh=18.1] (-3.2432,1.6008) -- (-3.3628,1.4789) -- (-3.4264,1.4866) -- (-3.3045,1.6108) -- cycle;
\path[cbpmesh=19.3] (3.3628,1.4789) -- (3.2432,1.6008) -- (3.3045,1.6108) -- (3.4264,1.4866) -- cycle;
\path[cbpmesh=22.6] (-0.3303,1.1713) -- (-0.4372,1.1571) -- (-0.4521,1.2748) -- (-0.3416,1.2894) -- cycle;
\path[cbpmesh=18.4] (0.4372,1.1571) -- (0.3303,1.1713) -- (0.3416,1.2894) -- (0.4521,1.2748) -- cycle;
\path[cbpmesh=21.0] (-0.3416,1.2894) -- (-0.4521,1.2748) -- (-0.4700,1.3877) -- (-0.3550,1.4029) -- cycle;
\path[cbpmesh=17.1] (0.4521,1.2748) -- (0.3416,1.2894) -- (0.3550,1.4029) -- (0.4700,1.3877) -- cycle;
\path[cbpmesh=15.1] (1.1343,1.7360) -- (1.0083,1.7755) -- (1.0614,1.8451) -- (1.1940,1.8035) -- cycle;
\path[cbpmesh=17.9] (-1.0083,1.7755) -- (-1.1343,1.7360) -- (-1.1940,1.8035) -- (-1.0614,1.8451) -- cycle;
\path[cbpmesh=33.0] (-0.6935,0.5887) -- (-0.7801,0.5615) -- (-0.7830,0.6925) -- (-0.6960,0.7198) -- cycle;
\path[cbpmesh=23.7] (0.7801,0.5615) -- (0.6935,0.5887) -- (0.6960,0.7198) -- (0.7830,0.6925) -- cycle;
\path[cbpmesh=14.9] (0.6075,1.4753) -- (0.4905,1.4955) -- (0.5136,1.5976) -- (0.6361,1.5764) -- cycle;
\path[cbpmesh=19.1] (-0.4905,1.4955) -- (-0.6075,1.4753) -- (-0.6361,1.5764) -- (-0.5136,1.5976) -- cycle;
\path[cbpmesh=31.1] (1.1347,0.0275) -- (1.0665,0.0681) -- (1.0304,0.1898) -- (1.0962,0.1505) -- cycle;
\path[cbpmesh=38] (-1.0665,0.0681) -- (-1.1347,0.0275) -- (-1.0962,0.1505) -- (-1.0304,0.1898) -- cycle;
\path[cbpmesh=19.5] (0.5268,1.0178) -- (0.4254,1.0354) -- (0.4372,1.1571) -- (0.5414,1.1391) -- cycle;
\path[cbpmesh=25.1] (-0.4254,1.0354) -- (-0.5268,1.0178) -- (-0.5414,1.1391) -- (-0.4372,1.1571) -- cycle;
\path[cbpmesh=19.6] (-3.7952,1.4410) -- (-3.9069,1.2989) -- (-3.9731,1.2645) -- (-3.8596,1.4089) -- cycle;
\path[cbpmesh=28.9] (3.9069,1.2989) -- (3.7952,1.4410) -- (3.8596,1.4089) -- (3.9731,1.2645) -- cycle;
\path[cbpmesh=30.3] (-0.6048,0.7437) -- (-0.6960,0.7198) -- (-0.7047,0.8476) -- (-0.6123,0.8719) -- cycle;
\path[cbpmesh=22.2] (0.6960,0.7198) -- (0.6048,0.7437) -- (0.6123,0.8719) -- (0.7047,0.8476) -- cycle;
\path[cbpmesh=27.7] (-0.5161,0.8927) -- (-0.6123,0.8719) -- (-0.6249,0.9966) -- (-0.5268,1.0178) -- cycle;
\path[cbpmesh=20.8] (0.6123,0.8719) -- (0.5161,0.8927) -- (0.5268,1.0178) -- (0.6249,0.9966) -- cycle;
\path[cbpmesh=18.7] (-1.6729,1.8396) -- (-1.7990,1.7756) -- (-1.8898,1.7971) -- (-1.7573,1.8644) -- cycle;
\path[cbpmesh=20.4] (1.7990,1.7756) -- (1.6729,1.8396) -- (1.7573,1.8644) -- (1.8898,1.7971) -- cycle;
\path[cbpmesh=18.7] (-0.6361,1.5764) -- (-0.7546,1.5507) -- (-0.7920,1.6443) -- (-0.6676,1.6712) -- cycle;
\path[cbpmesh=14.3] (0.7546,1.5507) -- (0.6361,1.5764) -- (0.6676,1.6712) -- (0.7920,1.6443) -- cycle;
\path[cbpmesh=18.2] (-3.1798,1.5889) -- (-3.2971,1.4695) -- (-3.3628,1.4789) -- (-3.2432,1.6008) -- cycle;
\path[cbpmesh=19.0] (3.2971,1.4695) -- (3.1798,1.5889) -- (3.2432,1.6008) -- (3.3628,1.4789) -- cycle;
\path[cbpmesh=18.0] (-1.1940,1.8035) -- (-1.3193,1.7572) -- (-1.3888,1.8137) -- (-1.2569,1.8625) -- cycle;
\path[cbpmesh=16.0] (1.3193,1.7572) -- (1.1940,1.8035) -- (1.2569,1.8625) -- (1.3888,1.8137) -- cycle;
\path[cbpmesh=20.4] (-0.4700,1.3877) -- (-0.5820,1.3683) -- (-0.6075,1.4753) -- (-0.4905,1.4955) -- cycle;
\path[cbpmesh=15.7] (0.5820,1.3683) -- (0.4700,1.3877) -- (0.4905,1.4955) -- (0.6075,1.4753) -- cycle;
\path[cbpmesh=23.4] (-0.4372,1.1571) -- (-0.5414,1.1391) -- (-0.5599,1.2561) -- (-0.4521,1.2748) -- cycle;
\path[cbpmesh=18.1] (0.5414,1.1391) -- (0.4372,1.1571) -- (0.4521,1.2748) -- (0.5599,1.2561) -- cycle;
\path[cbpmesh=27.0] (1.0304,0.1898) -- (0.9581,0.2264) -- (0.9437,0.3638) -- (1.0149,0.3277) -- cycle;
\path[cbpmesh=38] (-0.9581,0.2264) -- (-1.0304,0.1898) -- (-1.0149,0.3277) -- (-0.9437,0.3638) -- cycle;
\path[cbpmesh=21.8] (-0.4521,1.2748) -- (-0.5599,1.2561) -- (-0.5820,1.3683) -- (-0.4700,1.3877) -- cycle;
\path[cbpmesh=16.8] (0.5599,1.2561) -- (0.4521,1.2748) -- (0.4700,1.3877) -- (0.5820,1.3683) -- cycle;
\path[cbpmesh=18.3] (-2.6114,1.6312) -- (-2.7289,1.5327) -- (-2.8151,1.5355) -- (-2.6939,1.6371) -- cycle;
\path[cbpmesh=21.3] (2.7289,1.5327) -- (2.6114,1.6312) -- (2.6939,1.6371) -- (2.8151,1.5355) -- cycle;
\path[cbpmesh=18.5] (-0.7920,1.6443) -- (-0.9115,1.6129) -- (-0.9582,1.6980) -- (-0.8326,1.7310) -- cycle;
\path[cbpmesh=14.1] (0.9115,1.6129) -- (0.7920,1.6443) -- (0.8326,1.7310) -- (0.9582,1.6980) -- cycle;
\path[cbpmesh=36.7] (-0.8668,0.3972) -- (-0.9437,0.3638) -- (-0.9385,0.4980) -- (-0.8620,0.5313) -- cycle;
\path[cbpmesh=25.3] (0.9437,0.3638) -- (0.8668,0.3972) -- (0.8620,0.5313) -- (0.9385,0.4980) -- cycle;
\path[cbpmesh=19.2] (0.6249,0.9966) -- (0.5268,1.0178) -- (0.5414,1.1391) -- (0.6423,1.1173) -- cycle;
\path[cbpmesh=25.9] (-0.5268,1.0178) -- (-0.6249,0.9966) -- (-0.6423,1.1173) -- (-0.5414,1.1391) -- cycle;
\path[cbpmesh=14.7] (0.7207,1.4508) -- (0.6075,1.4753) -- (0.6361,1.5764) -- (0.7546,1.5507) -- cycle;
\path[cbpmesh=19.8] (-0.6075,1.4753) -- (-0.7207,1.4508) -- (-0.7546,1.5507) -- (-0.6361,1.5764) -- cycle;
\path[cbpmesh=18.9] (3.2289,1.4592) -- (3.1140,1.5762) -- (3.1798,1.5889) -- (3.2971,1.4695) -- cycle;
\path[cbpmesh=18.2] (-3.1140,1.5762) -- (-3.2289,1.4592) -- (-3.2971,1.4695) -- (-3.1798,1.5889) -- cycle;
\path[cbpmesh=18.2] (-1.3888,1.8137) -- (-1.5121,1.7602) -- (-1.5911,1.8048) -- (-1.4613,1.8611) -- cycle;
\path[cbpmesh=17.3] (1.5121,1.7602) -- (1.3888,1.8137) -- (1.4613,1.8611) -- (1.5911,1.8048) -- cycle;
\path[cbpmesh=18.9] (-3.7326,1.4644) -- (-3.8424,1.3247) -- (-3.9069,1.2989) -- (-3.7952,1.4410) -- cycle;
\path[cbpmesh=27.0] (3.8424,1.3247) -- (3.7326,1.4644) -- (3.7952,1.4410) -- (3.9069,1.2989) -- cycle;
\path[cbpmesh=19.2] (-2.2102,1.7287) -- (-2.3294,1.6449) -- (-2.4296,1.6341) -- (-2.3053,1.7215) -- cycle;
\path[cbpmesh=24.2] (2.3294,1.6449) -- (2.2102,1.7287) -- (2.3053,1.7215) -- (2.4296,1.6341) -- cycle;
\path[cbpmesh=18.3] (-0.9582,1.6980) -- (-1.0779,1.6604) -- (-1.1343,1.7360) -- (-1.0083,1.7755) -- cycle;
\path[cbpmesh=14.4] (1.0779,1.6604) -- (0.9582,1.6980) -- (1.0083,1.7755) -- (1.1343,1.7360) -- cycle;
\path[cbpmesh=34.1] (-0.7801,0.5615) -- (-0.8620,0.5313) -- (-0.8651,0.6621) -- (-0.7830,0.6925) -- cycle;
\path[cbpmesh=23.6] (0.8620,0.5313) -- (0.7801,0.5615) -- (0.7830,0.6925) -- (0.8651,0.6621) -- cycle;
\path[cbpmesh=28.6] (-0.6123,0.8719) -- (-0.7047,0.8476) -- (-0.7192,0.9718) -- (-0.6249,0.9966) -- cycle;
\path[cbpmesh=20.6] (0.7047,0.8476) -- (0.6123,0.8719) -- (0.6249,0.9966) -- (0.7192,0.9718) -- cycle;
\path[cbpmesh=31.3] (-0.6960,0.7198) -- (-0.7830,0.6925) -- (-0.7928,0.8200) -- (-0.7047,0.8476) -- cycle;
\path[cbpmesh=22.1] (0.7830,0.6925) -- (0.6960,0.7198) -- (0.7047,0.8476) -- (0.7928,0.8200) -- cycle;
\path[cbpmesh=24.2] (-0.5414,1.1391) -- (-0.6423,1.1173) -- (-0.6643,1.2335) -- (-0.5599,1.2561) -- cycle;
\path[cbpmesh=17.8] (0.6423,1.1173) -- (0.5414,1.1391) -- (0.5599,1.2561) -- (0.6643,1.2335) -- cycle;
\path[cbpmesh=21.1] (-0.5820,1.3683) -- (-0.6905,1.3449) -- (-0.7207,1.4508) -- (-0.6075,1.4753) -- cycle;
\path[cbpmesh=15.5] (0.6905,1.3449) -- (0.5820,1.3683) -- (0.6075,1.4753) -- (0.7207,1.4508) -- cycle;
\path[cbpmesh=18.9] (-1.8898,1.7971) -- (-2.0106,1.7252) -- (-2.1095,1.7325) -- (-1.9827,1.8080) -- cycle;
\path[cbpmesh=22.0] (2.0106,1.7252) -- (1.8898,1.7971) -- (1.9827,1.8080) -- (2.1095,1.7325) -- cycle;
\path[cbpmesh=22.6] (-0.5599,1.2561) -- (-0.6643,1.2335) -- (-0.6905,1.3449) -- (-0.5820,1.3683) -- cycle;
\path[cbpmesh=16.5] (0.6643,1.2335) -- (0.5599,1.2561) -- (0.5820,1.3683) -- (0.6905,1.3449) -- cycle;
\path[cbpmesh=19.4] (-0.7546,1.5507) -- (-0.8685,1.5208) -- (-0.9115,1.6129) -- (-0.7920,1.6443) -- cycle;
\path[cbpmesh=14.2] (0.8685,1.5208) -- (0.7546,1.5507) -- (0.7920,1.6443) -- (0.9115,1.6129) -- cycle;
\path[cbpmesh=38] (-1.2457,-0.1153) -- (-1.3018,-0.1623) -- (-1.2497,-0.0606) -- (-1.1959,-0.0155) -- cycle;
\path[cbpmesh=34.6] (1.3018,-0.1623) -- (1.2457,-0.1153) -- (1.1959,-0.0155) -- (1.2497,-0.0606) -- cycle;
\path[cbpmesh=18.2] (-3.0451,1.5635) -- (-3.1575,1.4491) -- (-3.2289,1.4592) -- (-3.1140,1.5762) -- cycle;
\path[cbpmesh=19.0] (3.1575,1.4491) -- (3.0451,1.5635) -- (3.1140,1.5762) -- (3.2289,1.4592) -- cycle;
\path[cbpmesh=15.0] (1.2533,1.6920) -- (1.1343,1.7360) -- (1.1940,1.8035) -- (1.3193,1.7572) -- cycle;
\path[cbpmesh=18.2] (-1.1343,1.7360) -- (-1.2533,1.6920) -- (-1.3193,1.7572) -- (-1.1940,1.8035) -- cycle;
\path[cbpmesh=38] (-1.1347,0.0275) -- (-1.1959,-0.0155) -- (-1.1553,0.1090) -- (-1.0962,0.1505) -- cycle;
\path[cbpmesh=31.4] (1.1959,-0.0155) -- (1.1347,0.0275) -- (1.0962,0.1505) -- (1.1553,0.1090) -- cycle;
\path[cbpmesh=38] (-1.3627,-0.2417) -- (-1.4129,-0.2929) -- (-1.3498,-0.2112) -- (-1.3018,-0.1623) -- cycle;
\draw[draw=cbpInk!70,line width=.45pt] (-1.3627,-0.2417) -- (-1.4129,-0.2929);
\path[cbpmesh=38] (1.4129,-0.2929) -- (1.3627,-0.2417) -- (1.3018,-0.1623) -- (1.3498,-0.2112) -- cycle;
\draw[draw=cbpInk!70,line width=.45pt] (1.4129,-0.2929) -- (1.3627,-0.2417);
\path[cbpmesh=22.2] (2.6370,1.5347) -- (2.5235,1.6299) -- (2.6114,1.6312) -- (2.7289,1.5327) -- cycle;
\path[cbpmesh=18.4] (-2.5235,1.6299) -- (-2.6370,1.5347) -- (-2.7289,1.5327) -- (-2.6114,1.6312) -- cycle;
\path[cbpmesh=18.4] (-3.6711,1.4802) -- (-3.7791,1.3428) -- (-3.8424,1.3247) -- (-3.7326,1.4644) -- cycle;
\path[cbpmesh=25.2] (3.7791,1.3428) -- (3.6711,1.4802) -- (3.7326,1.4644) -- (3.8424,1.3247) -- cycle;
\path[cbpmesh=19.1] (0.7192,0.9718) -- (0.6249,0.9966) -- (0.6423,1.1173) -- (0.7393,1.0918) -- cycle;
\path[cbpmesh=26.9] (-0.6249,0.9966) -- (-0.7192,0.9718) -- (-0.7393,1.0918) -- (-0.6423,1.1173) -- cycle;
\path[cbpmesh=18.3] (-1.5911,1.8048) -- (-1.7110,1.7438) -- (-1.7990,1.7756) -- (-1.6729,1.8396) -- cycle;
\path[cbpmesh=18.8] (1.7110,1.7438) -- (1.5911,1.8048) -- (1.6729,1.8396) -- (1.7990,1.7756) -- cycle;
\path[cbpmesh=38] (-1.0304,0.1898) -- (-1.0962,0.1505) -- (-1.0798,0.2890) -- (-1.0149,0.3277) -- cycle;
\path[cbpmesh=27.2] (1.0962,0.1505) -- (1.0304,0.1898) -- (1.0149,0.3277) -- (1.0798,0.2890) -- cycle;
\path[cbpmesh=14.6] (0.8294,1.4222) -- (0.7207,1.4508) -- (0.7546,1.5507) -- (0.8685,1.5208) -- cycle;
\path[cbpmesh=20.5] (-0.7207,1.4508) -- (-0.8294,1.4222) -- (-0.8685,1.5208) -- (-0.7546,1.5507) -- cycle;
\path[cbpmesh=14.0] (1.0254,1.5772) -- (0.9115,1.6129) -- (0.9582,1.6980) -- (1.0779,1.6604) -- cycle;
\path[cbpmesh=19.0] (-0.9115,1.6129) -- (-1.0254,1.5772) -- (-1.0779,1.6604) -- (-0.9582,1.6980) -- cycle;
\path[cbpmesh=29.6] (-0.7047,0.8476) -- (-0.7928,0.8200) -- (-0.8091,0.9436) -- (-0.7192,0.9718) -- cycle;
\path[cbpmesh=20.5] (0.7928,0.8200) -- (0.7047,0.8476) -- (0.7192,0.9718) -- (0.8091,0.9436) -- cycle;
\path[cbpmesh=37.8] (-0.9437,0.3638) -- (-1.0149,0.3277) -- (-1.0092,0.4621) -- (-0.9385,0.4980) -- cycle;
\path[cbpmesh=25.4] (1.0149,0.3277) -- (0.9437,0.3638) -- (0.9385,0.4980) -- (1.0092,0.4621) -- cycle;
\path[cbpmesh=25.1] (-0.6423,1.1173) -- (-0.7393,1.0918) -- (-0.7645,1.2072) -- (-0.6643,1.2335) -- cycle;
\path[cbpmesh=17.6] (0.7393,1.0918) -- (0.6423,1.1173) -- (0.6643,1.2335) -- (0.7645,1.2072) -- cycle;
\path[cbpmesh=18.1] (-2.9727,1.5519) -- (-3.0824,1.4402) -- (-3.1575,1.4491) -- (-3.0451,1.5635) -- cycle;
\path[cbpmesh=19.3] (3.0824,1.4402) -- (2.9727,1.5519) -- (3.0451,1.5635) -- (3.1575,1.4491) -- cycle;
\path[cbpmesh=21.9] (-0.6905,1.3449) -- (-0.7947,1.3175) -- (-0.8294,1.4222) -- (-0.7207,1.4508) -- cycle;
\path[cbpmesh=15.3] (0.7947,1.3175) -- (0.6905,1.3449) -- (0.7207,1.4508) -- (0.8294,1.4222) -- cycle;
\path[cbpmesh=32.4] (-0.7830,0.6925) -- (-0.8651,0.6621) -- (-0.8759,0.7892) -- (-0.7928,0.8200) -- cycle;
\path[cbpmesh=22.0] (0.8651,0.6621) -- (0.7830,0.6925) -- (0.7928,0.8200) -- (0.8759,0.7892) -- cycle;
\path[cbpmesh=16.1] (1.4364,1.7063) -- (1.3193,1.7572) -- (1.3888,1.8137) -- (1.5121,1.7602) -- cycle;
\path[cbpmesh=18.2] (-1.3193,1.7572) -- (-1.4364,1.7063) -- (-1.5121,1.7602) -- (-1.3888,1.8137) -- cycle;
\path[cbpmesh=35.1] (-0.8620,0.5313) -- (-0.9385,0.4980) -- (-0.9419,0.6287) -- (-0.8651,0.6621) -- cycle;
\path[cbpmesh=23.7] (0.9385,0.4980) -- (0.8620,0.5313) -- (0.8651,0.6621) -- (0.9419,0.6287) -- cycle;
\path[cbpmesh=23.5] (-0.6643,1.2335) -- (-0.7645,1.2072) -- (-0.7947,1.3175) -- (-0.6905,1.3449) -- cycle;
\path[cbpmesh=16.4] (0.7645,1.2072) -- (0.6643,1.2335) -- (0.6905,1.3449) -- (0.7947,1.3175) -- cycle;
\path[cbpmesh=23.6] (-4.1133,1.1656) -- (-4.2055,1.0122) -- (-4.2853,0.8976) -- (-4.1914,1.0539) -- cycle;
\draw[draw=cbpInk!70,line width=.45pt] (-4.2853,0.8976) -- (-4.1914,1.0539);
\draw[draw=cbpInk!70,line width=.45pt] (-4.1914,1.0539) -- (-4.1133,1.1656);
\path[cbpmesh=37.6] (4.2055,1.0122) -- (4.1133,1.1656) -- (4.1914,1.0539) -- (4.2853,0.8976) -- cycle;
\draw[draw=cbpInk!70,line width=.45pt] (4.1133,1.1656) -- (4.1914,1.0539);
\draw[draw=cbpInk!70,line width=.45pt] (4.1914,1.0539) -- (4.2853,0.8976);
\path[cbpmesh=18.1] (-3.6103,1.4893) -- (-3.7165,1.3542) -- (-3.7791,1.3428) -- (-3.6711,1.4802) -- cycle;
\path[cbpmesh=23.4] (3.7165,1.3542) -- (3.6103,1.4893) -- (3.6711,1.4802) -- (3.7791,1.3428) -- cycle;
\path[cbpmesh=18.7] (-1.0779,1.6604) -- (-1.1910,1.6186) -- (-1.2533,1.6920) -- (-1.1343,1.7360) -- cycle;
\path[cbpmesh=14.4] (1.1910,1.6186) -- (1.0779,1.6604) -- (1.1343,1.7360) -- (1.2533,1.6920) -- cycle;
\path[cbpmesh=19.1] (-2.1095,1.7325) -- (-2.2232,1.6526) -- (-2.3294,1.6449) -- (-2.2102,1.7287) -- cycle;
\path[cbpmesh=23.7] (2.2232,1.6526) -- (2.1095,1.7325) -- (2.2102,1.7287) -- (2.3294,1.6449) -- cycle;
\path[cbpmesh=20.0] (-0.8685,1.5208) -- (-0.9770,1.4868) -- (-1.0254,1.5772) -- (-0.9115,1.6129) -- cycle;
\path[cbpmesh=14.1] (0.9770,1.4868) -- (0.8685,1.5208) -- (0.9115,1.6129) -- (1.0254,1.5772) -- cycle;
\path[cbpmesh=18.9] (0.8091,0.9436) -- (0.7192,0.9718) -- (0.7393,1.0918) -- (0.8316,1.0628) -- cycle;
\path[cbpmesh=27.8] (-0.7192,0.9718) -- (-0.8091,0.9436) -- (-0.8316,1.0628) -- (-0.7393,1.0918) -- cycle;
\path[cbpmesh=18.5] (-1.7990,1.7756) -- (-1.9140,1.7071) -- (-2.0106,1.7252) -- (-1.8898,1.7971) -- cycle;
\path[cbpmesh=20.3] (1.9140,1.7071) -- (1.7990,1.7756) -- (1.8898,1.7971) -- (2.0106,1.7252) -- cycle;
\path[cbpmesh=18.1] (-2.8962,1.5422) -- (-3.0031,1.4334) -- (-3.0824,1.4402) -- (-2.9727,1.5519) -- cycle;
\path[cbpmesh=19.8] (3.0031,1.4334) -- (2.8962,1.5422) -- (2.9727,1.5519) -- (3.0824,1.4402) -- cycle;
\path[cbpmesh=18.6] (-2.4296,1.6341) -- (-2.5390,1.5425) -- (-2.6370,1.5347) -- (-2.5235,1.6299) -- cycle;
\path[cbpmesh=23.1] (2.5390,1.5425) -- (2.4296,1.6341) -- (2.5235,1.6299) -- (2.6370,1.5347) -- cycle;
\path[cbpmesh=17.9] (-3.5496,1.4928) -- (-3.6540,1.3599) -- (-3.7165,1.3542) -- (-3.6103,1.4893) -- cycle;
\path[cbpmesh=21.9] (3.6540,1.3599) -- (3.5496,1.4928) -- (3.6103,1.4893) -- (3.7165,1.3542) -- cycle;
\path[cbpmesh=14.5] (0.9330,1.3897) -- (0.8294,1.4222) -- (0.8685,1.5208) -- (0.9770,1.4868) -- cycle;
\path[cbpmesh=21.2] (-0.8294,1.4222) -- (-0.9330,1.3897) -- (-0.9770,1.4868) -- (-0.8685,1.5208) -- cycle;
\path[cbpmesh=17.3] (1.6261,1.7023) -- (1.5121,1.7602) -- (1.5911,1.8048) -- (1.7110,1.7438) -- cycle;
\path[cbpmesh=18.3] (-1.5121,1.7602) -- (-1.6261,1.7023) -- (-1.7110,1.7438) -- (-1.5911,1.8048) -- cycle;
\path[cbpmesh=17.5] (0.8316,1.0628) -- (0.7393,1.0918) -- (0.7645,1.2072) -- (0.8600,1.1772) -- cycle;
\path[cbpmesh=26.1] (-0.7393,1.0918) -- (-0.8316,1.0628) -- (-0.8600,1.1772) -- (-0.7645,1.2072) -- cycle;
\path[cbpmesh=30.6] (-0.7928,0.8200) -- (-0.8759,0.7892) -- (-0.8940,0.9122) -- (-0.8091,0.9436) -- cycle;
\path[cbpmesh=20.4] (0.8759,0.7892) -- (0.7928,0.8200) -- (0.8091,0.9436) -- (0.8940,0.9122) -- cycle;
\path[cbpmesh=38] (-1.3018,-0.1623) -- (-1.3498,-0.2112) -- (-1.2958,-0.1076) -- (-1.2497,-0.0606) -- cycle;
\path[cbpmesh=35.0] (1.3498,-0.2112) -- (1.3018,-0.1623) -- (1.2497,-0.0606) -- (1.2958,-0.1076) -- cycle;
\path[cbpmesh=31.7] (1.2497,-0.0606) -- (1.1959,-0.0155) -- (1.1553,0.1090) -- (1.2073,0.0654) -- cycle;
\path[cbpmesh=38] (-1.1959,-0.0155) -- (-1.2497,-0.0606) -- (-1.2073,0.0654) -- (-1.1553,0.1090) -- cycle;
\path[cbpmesh=15.2] (0.8940,1.2863) -- (0.7947,1.3175) -- (0.8294,1.4222) -- (0.9330,1.3897) -- cycle;
\path[cbpmesh=22.7] (-0.7947,1.3175) -- (-0.8940,1.2863) -- (-0.9330,1.3897) -- (-0.8294,1.4222) -- cycle;
\path[cbpmesh=18.6] (-1.2533,1.6920) -- (-1.3645,1.6437) -- (-1.4364,1.7063) -- (-1.3193,1.7572) -- cycle;
\path[cbpmesh=15.1] (1.3645,1.6437) -- (1.2533,1.6920) -- (1.3193,1.7572) -- (1.4364,1.7063) -- cycle;
\path[cbpmesh=16.3] (0.8600,1.1772) -- (0.7645,1.2072) -- (0.7947,1.3175) -- (0.8940,1.2863) -- cycle;
\path[cbpmesh=24.3] (-0.7645,1.2072) -- (-0.8600,1.1772) -- (-0.8940,1.2863) -- (-0.7947,1.3175) -- cycle;
\path[cbpmesh=14.0] (1.1329,1.5374) -- (1.0254,1.5772) -- (1.0779,1.6604) -- (1.1910,1.6186) -- cycle;
\path[cbpmesh=19.5] (-1.0254,1.5772) -- (-1.1329,1.5374) -- (-1.1910,1.6186) -- (-1.0779,1.6604) -- cycle;
\path[cbpmesh=38] (-1.0962,0.1505) -- (-1.1553,0.1090) -- (-1.1380,0.2481) -- (-1.0798,0.2890) -- cycle;
\path[cbpmesh=27.5] (1.1553,0.1090) -- (1.0962,0.1505) -- (1.0798,0.2890) -- (1.1380,0.2481) -- cycle;
\path[cbpmesh=33.4] (-0.8651,0.6621) -- (-0.9419,0.6287) -- (-0.9537,0.7554) -- (-0.8759,0.7892) -- cycle;
\path[cbpmesh=22.1] (0.9419,0.6287) -- (0.8651,0.6621) -- (0.8759,0.7892) -- (0.9537,0.7554) -- cycle;
\path[cbpmesh=17.9] (-3.4885,1.4916) -- (-3.5911,1.3610) -- (-3.6540,1.3599) -- (-3.5496,1.4928) -- cycle;
\path[cbpmesh=20.7] (3.5911,1.3610) -- (3.4885,1.4916) -- (3.5496,1.4928) -- (3.6540,1.3599) -- cycle;
\path[cbpmesh=38] (-1.4129,-0.2929) -- (-1.4545,-0.3458) -- (-1.3895,-0.2618) -- (-1.3498,-0.2112) -- cycle;
\draw[draw=cbpInk!70,line width=.45pt] (-1.4129,-0.2929) -- (-1.4545,-0.3458);
\draw[draw=cbpInk!70,line width=.45pt] (-1.4545,-0.3458) -- (-1.3895,-0.2618);
\path[cbpmesh=38] (1.4545,-0.3458) -- (1.4129,-0.2929) -- (1.3498,-0.2112) -- (1.3895,-0.2618) -- cycle;
\draw[draw=cbpInk!70,line width=.45pt] (1.4545,-0.3458) -- (1.4129,-0.2929);
\draw[draw=cbpInk!70,line width=.45pt] (1.3895,-0.2618) -- (1.4545,-0.3458);
\path[cbpmesh=38] (-1.0149,0.3277) -- (-1.0798,0.2890) -- (-1.0737,0.4237) -- (-1.0092,0.4621) -- cycle;
\path[cbpmesh=25.7] (1.0798,0.2890) -- (1.0149,0.3277) -- (1.0092,0.4621) -- (1.0737,0.4237) -- cycle;
\path[cbpmesh=36.2] (-0.9385,0.4980) -- (-1.0092,0.4621) -- (-1.0129,0.5927) -- (-0.9419,0.6287) -- cycle;
\path[cbpmesh=23.8] (1.0092,0.4621) -- (0.9385,0.4980) -- (0.9419,0.6287) -- (1.0129,0.5927) -- cycle;
\path[cbpmesh=18.1] (-2.8151,1.5355) -- (-2.9190,1.4297) -- (-3.0031,1.4334) -- (-2.8962,1.5422) -- cycle;
\path[cbpmesh=20.4] (2.9190,1.4297) -- (2.8151,1.5355) -- (2.8962,1.5422) -- (3.0031,1.4334) -- cycle;
\path[cbpmesh=20.4] (-4.0418,1.2204) -- (-4.1324,1.0696) -- (-4.2055,1.0122) -- (-4.1133,1.1656) -- cycle;
\path[cbpmesh=31.8] (4.1324,1.0696) -- (4.0418,1.2204) -- (4.1133,1.1656) -- (4.2055,1.0122) -- cycle;
\path[cbpmesh=18.9] (0.8940,0.9122) -- (0.8091,0.9436) -- (0.8316,1.0628) -- (0.9189,1.0305) -- cycle;
\path[cbpmesh=28.8] (-0.8091,0.9436) -- (-0.8940,0.9122) -- (-0.9189,1.0305) -- (-0.8316,1.0628) -- cycle;
\path[cbpmesh=14.1] (1.0795,1.4488) -- (0.9770,1.4868) -- (1.0254,1.5772) -- (1.1329,1.5374) -- cycle;
\path[cbpmesh=20.6] (-0.9770,1.4868) -- (-1.0795,1.4488) -- (-1.1329,1.5374) -- (-1.0254,1.5772) -- cycle;
\path[cbpmesh=17.9] (-3.4264,1.4866) -- (-3.5272,1.3584) -- (-3.5911,1.3610) -- (-3.4885,1.4916) -- cycle;
\path[cbpmesh=19.8] (3.5272,1.3584) -- (3.4264,1.4866) -- (3.4885,1.4916) -- (3.5911,1.3610) -- cycle;
\path[cbpmesh=18.6] (-2.0106,1.7252) -- (-2.1190,1.6490) -- (-2.2232,1.6526) -- (-2.1095,1.7325) -- cycle;
\path[cbpmesh=22.0] (2.1190,1.6490) -- (2.0106,1.7252) -- (2.1095,1.7325) -- (2.2232,1.6526) -- cycle;
\path[cbpmesh=18.3] (-1.7110,1.7438) -- (-1.8204,1.6787) -- (-1.9140,1.7071) -- (-1.7990,1.7756) -- cycle;
\path[cbpmesh=18.8] (1.8204,1.6787) -- (1.7110,1.7438) -- (1.7990,1.7756) -- (1.9140,1.7071) -- cycle;
\path[cbpmesh=14.4] (1.2967,1.5727) -- (1.1910,1.6186) -- (1.2533,1.6920) -- (1.3645,1.6437) -- cycle;
\path[cbpmesh=19.2] (-1.1910,1.6186) -- (-1.2967,1.5727) -- (-1.3645,1.6437) -- (-1.2533,1.6920) -- cycle;
\path[cbpmesh=24.1] (2.4342,1.5571) -- (2.3294,1.6449) -- (2.4296,1.6341) -- (2.5390,1.5425) -- cycle;
\path[cbpmesh=18.9] (-2.3294,1.6449) -- (-2.4342,1.5571) -- (-2.5390,1.5425) -- (-2.4296,1.6341) -- cycle;
\path[cbpmesh=18.4] (-1.4364,1.7063) -- (-1.5447,1.6514) -- (-1.6261,1.7023) -- (-1.5121,1.7602) -- cycle;
\path[cbpmesh=16.1] (1.5447,1.6514) -- (1.4364,1.7063) -- (1.5121,1.7602) -- (1.6261,1.7023) -- cycle;
\path[cbpmesh=17.5] (0.9189,1.0305) -- (0.8316,1.0628) -- (0.8600,1.1772) -- (0.9503,1.1438) -- cycle;
\path[cbpmesh=27.0] (-0.8316,1.0628) -- (-0.9189,1.0305) -- (-0.9503,1.1438) -- (-0.8600,1.1772) -- cycle;
\path[cbpmesh=14.5] (1.0309,1.3534) -- (0.9330,1.3897) -- (0.9770,1.4868) -- (1.0795,1.4488) -- cycle;
\path[cbpmesh=22.0] (-0.9330,1.3897) -- (-1.0309,1.3534) -- (-1.0795,1.4488) -- (-0.9770,1.4868) -- cycle;
\path[cbpmesh=31.6] (-0.8759,0.7892) -- (-0.9537,0.7554) -- (-0.9733,0.8777) -- (-0.8940,0.9122) -- cycle;
\path[cbpmesh=20.5] (0.9537,0.7554) -- (0.8759,0.7892) -- (0.8940,0.9122) -- (0.9733,0.8777) -- cycle;
\path[cbpmesh=18.1] (-2.7289,1.5327) -- (-2.8296,1.4301) -- (-2.9190,1.4297) -- (-2.8151,1.5355) -- cycle;
\path[cbpmesh=21.2] (2.8296,1.4301) -- (2.7289,1.5327) -- (2.8151,1.5355) -- (2.9190,1.4297) -- cycle;
\path[cbpmesh=18.0] (-3.3628,1.4789) -- (-3.4617,1.3531) -- (-3.5272,1.3584) -- (-3.4264,1.4866) -- cycle;
\path[cbpmesh=19.2] (3.4617,1.3531) -- (3.3628,1.4789) -- (3.4264,1.4866) -- (3.5272,1.3584) -- cycle;
\path[cbpmesh=25.2] (-0.8600,1.1772) -- (-0.9503,1.1438) -- (-0.9877,1.2516) -- (-0.8940,1.2863) -- cycle;
\path[cbpmesh=16.2] (0.9503,1.1438) -- (0.8600,1.1772) -- (0.8940,1.2863) -- (0.9877,1.2516) -- cycle;
\path[cbpmesh=15.2] (0.9877,1.2516) -- (0.8940,1.2863) -- (0.9330,1.3897) -- (1.0309,1.3534) -- cycle;
\path[cbpmesh=23.5] (-0.8940,1.2863) -- (-0.9877,1.2516) -- (-1.0309,1.3534) -- (-0.9330,1.3897) -- cycle;
\path[cbpmesh=19.7] (-3.9731,1.2645) -- (-4.0622,1.1162) -- (-4.1324,1.0696) -- (-4.0418,1.2204) -- cycle;
\path[cbpmesh=30.2] (4.0622,1.1162) -- (3.9731,1.2645) -- (4.0418,1.2204) -- (4.1324,1.0696) -- cycle;
\path[cbpmesh=34.5] (-0.9419,0.6287) -- (-1.0129,0.5927) -- (-1.0256,0.7189) -- (-0.9537,0.7554) -- cycle;
\path[cbpmesh=22.2] (1.0129,0.5927) -- (0.9419,0.6287) -- (0.9537,0.7554) -- (1.0256,0.7189) -- cycle;
\path[cbpmesh=38] (-1.2497,-0.0606) -- (-1.2958,-0.1076) -- (-1.2519,0.0201) -- (-1.2073,0.0654) -- cycle;
\path[cbpmesh=32.2] (1.2958,-0.1076) -- (1.2497,-0.0606) -- (1.2073,0.0654) -- (1.2519,0.0201) -- cycle;
\path[cbpmesh=24.0] (1.0737,0.4237) -- (1.0092,0.4621) -- (1.0129,0.5927) -- (1.0777,0.5541) -- cycle;
\path[cbpmesh=37.3] (-1.0092,0.4621) -- (-1.0737,0.4237) -- (-1.0777,0.5541) -- (-1.0129,0.5927) -- cycle;
\path[cbpmesh=18.1] (-3.2971,1.4695) -- (-3.3941,1.3460) -- (-3.4617,1.3531) -- (-3.3628,1.4789) -- cycle;
\path[cbpmesh=18.9] (3.3941,1.3460) -- (3.2971,1.4695) -- (3.3628,1.4789) -- (3.4617,1.3531) -- cycle;
\path[cbpmesh=27.9] (1.2073,0.0654) -- (1.1553,0.1090) -- (1.1380,0.2481) -- (1.1892,0.2052) -- cycle;
\path[cbpmesh=38] (-1.1553,0.1090) -- (-1.2073,0.0654) -- (-1.1892,0.2052) -- (-1.1380,0.2481) -- cycle;
\path[cbpmesh=14.1] (1.2335,1.4937) -- (1.1329,1.5374) -- (1.1910,1.6186) -- (1.2967,1.5727) -- cycle;
\path[cbpmesh=20.1] (-1.1329,1.5374) -- (-1.2335,1.4937) -- (-1.2967,1.5727) -- (-1.1910,1.6186) -- cycle;
\path[cbpmesh=38] (-1.3498,-0.2112) -- (-1.3895,-0.2618) -- (-1.3339,-0.1561) -- (-1.2958,-0.1076) -- cycle;
\path[cbpmesh=35.5] (1.3895,-0.2618) -- (1.3498,-0.2112) -- (1.2958,-0.1076) -- (1.3339,-0.1561) -- cycle;
\path[cbpmesh=38] (-1.0798,0.2890) -- (-1.1380,0.2481) -- (-1.1316,0.3830) -- (-1.0737,0.4237) -- cycle;
\path[cbpmesh=26.0] (1.1380,0.2481) -- (1.0798,0.2890) -- (1.0737,0.4237) -- (1.1316,0.3830) -- cycle;
\path[cbpmesh=29.8] (-0.8940,0.9122) -- (-0.9733,0.8777) -- (-1.0005,0.9951) -- (-0.9189,1.0305) -- cycle;
\path[cbpmesh=19.0] (0.9733,0.8777) -- (0.8940,0.9122) -- (0.9189,1.0305) -- (1.0005,0.9951) -- cycle;
\path[cbpmesh=18.9] (-1.3645,1.6437) -- (-1.4674,1.5915) -- (-1.5447,1.6514) -- (-1.4364,1.7063) -- cycle;
\path[cbpmesh=15.1] (1.4674,1.5915) -- (1.3645,1.6437) -- (1.4364,1.7063) -- (1.5447,1.6514) -- cycle;
\path[cbpmesh=17.3] (1.7300,1.6404) -- (1.6261,1.7023) -- (1.7110,1.7438) -- (1.8204,1.6787) -- cycle;
\path[cbpmesh=18.3] (-1.6261,1.7023) -- (-1.7300,1.6404) -- (-1.8204,1.6787) -- (-1.7110,1.7438) -- cycle;
\path[cbpmesh=18.2] (-2.6370,1.5347) -- (-2.7343,1.4356) -- (-2.8296,1.4301) -- (-2.7289,1.5327) -- cycle;
\path[cbpmesh=22.0] (2.7343,1.4356) -- (2.6370,1.5347) -- (2.7289,1.5327) -- (2.8296,1.4301) -- cycle;
\path[cbpmesh=18.1] (-3.2289,1.4592) -- (-3.3239,1.3383) -- (-3.3941,1.3460) -- (-3.2971,1.4695) -- cycle;
\path[cbpmesh=18.9] (3.3239,1.3383) -- (3.2289,1.4592) -- (3.2971,1.4695) -- (3.3941,1.3460) -- cycle;
\path[cbpmesh=21.3] (-1.0795,1.4488) -- (-1.1753,1.4072) -- (-1.2335,1.4937) -- (-1.1329,1.5374) -- cycle;
\path[cbpmesh=14.1] (1.1753,1.4072) -- (1.0795,1.4488) -- (1.1329,1.5374) -- (1.2335,1.4937) -- cycle;
\path[cbpmesh=18.4] (-1.9140,1.7071) -- (-2.0173,1.6346) -- (-2.1190,1.6490) -- (-2.0106,1.7252) -- cycle;
\path[cbpmesh=20.3] (2.0173,1.6346) -- (1.9140,1.7071) -- (2.0106,1.7252) -- (2.1190,1.6490) -- cycle;
\path[cbpmesh=19.0] (-3.9069,1.2989) -- (-3.9944,1.1532) -- (-4.0622,1.1162) -- (-3.9731,1.2645) -- cycle;
\path[cbpmesh=28.5] (3.9944,1.1532) -- (3.9069,1.2989) -- (3.9731,1.2645) -- (4.0622,1.1162) -- cycle;
\path[cbpmesh=28.0] (-0.9189,1.0305) -- (-1.0005,0.9951) -- (-1.0346,1.1072) -- (-0.9503,1.1438) -- cycle;
\path[cbpmesh=17.5] (1.0005,0.9951) -- (0.9189,1.0305) -- (0.9503,1.1438) -- (1.0346,1.1072) -- cycle;
\path[cbpmesh=18.7] (-2.2232,1.6526) -- (-2.3232,1.5687) -- (-2.4342,1.5571) -- (-2.3294,1.6449) -- cycle;
\path[cbpmesh=23.5] (2.3232,1.5687) -- (2.2232,1.6526) -- (2.3294,1.6449) -- (2.4342,1.5571) -- cycle;
\path[cbpmesh=20.6] (1.0256,0.7189) -- (0.9537,0.7554) -- (0.9733,0.8777) -- (1.0467,0.8404) -- cycle;
\path[cbpmesh=32.7] (-0.9537,0.7554) -- (-1.0256,0.7189) -- (-1.0467,0.8404) -- (-0.9733,0.8777) -- cycle;
\path[cbpmesh=14.5] (1.1225,1.3137) -- (1.0309,1.3534) -- (1.0795,1.4488) -- (1.1753,1.4072) -- cycle;
\path[cbpmesh=22.7] (-1.0309,1.3534) -- (-1.1225,1.3137) -- (-1.1753,1.4072) -- (-1.0795,1.4488) -- cycle;
\path[cbpmesh=26.1] (-0.9503,1.1438) -- (-1.0346,1.1072) -- (-1.0754,1.2135) -- (-0.9877,1.2516) -- cycle;
\path[cbpmesh=16.3] (1.0346,1.1072) -- (0.9503,1.1438) -- (0.9877,1.2516) -- (1.0754,1.2135) -- cycle;
\path[cbpmesh=18.1] (-3.1575,1.4491) -- (-3.2504,1.3309) -- (-3.3239,1.3383) -- (-3.2289,1.4592) -- cycle;
\path[cbpmesh=19.0] (3.2504,1.3309) -- (3.1575,1.4491) -- (3.2289,1.4592) -- (3.3239,1.3383) -- cycle;
\path[cbpmesh=24.4] (-0.9877,1.2516) -- (-1.0754,1.2135) -- (-1.1225,1.3137) -- (-1.0309,1.3534) -- cycle;
\path[cbpmesh=15.2] (1.0754,1.2135) -- (0.9877,1.2516) -- (1.0309,1.3534) -- (1.1225,1.3137) -- cycle;
\path[cbpmesh=14.4] (1.3945,1.5230) -- (1.2967,1.5727) -- (1.3645,1.6437) -- (1.4674,1.5915) -- cycle;
\path[cbpmesh=19.6] (-1.2967,1.5727) -- (-1.3945,1.5230) -- (-1.4674,1.5915) -- (-1.3645,1.6437) -- cycle;
\path[cbpmesh=35.5] (-1.0129,0.5927) -- (-1.0777,0.5541) -- (-1.0912,0.6799) -- (-1.0256,0.7189) -- cycle;
\path[cbpmesh=22.4] (1.0777,0.5541) -- (1.0129,0.5927) -- (1.0256,0.7189) -- (1.0912,0.6799) -- cycle;
\path[cbpmesh=18.4] (-3.8424,1.3247) -- (-3.9285,1.1813) -- (-3.9944,1.1532) -- (-3.9069,1.2989) -- cycle;
\path[cbpmesh=26.7] (3.9285,1.1813) -- (3.8424,1.3247) -- (3.9069,1.2989) -- (3.9944,1.1532) -- cycle;
\path[cbpmesh=24.3] (1.1316,0.3830) -- (1.0737,0.4237) -- (1.0777,0.5541) -- (1.1358,0.5133) -- cycle;
\path[cbpmesh=38] (-1.0737,0.4237) -- (-1.1316,0.3830) -- (-1.1358,0.5133) -- (-1.0777,0.5541) -- cycle;
\path[cbpmesh=18.0] (-3.0824,1.4402) -- (-3.1730,1.3248) -- (-3.2504,1.3309) -- (-3.1575,1.4491) -- cycle;
\path[cbpmesh=19.3] (3.1730,1.3248) -- (3.0824,1.4402) -- (3.1575,1.4491) -- (3.2504,1.3309) -- cycle;
\path[cbpmesh=18.6] (-1.5447,1.6514) -- (-1.6435,1.5926) -- (-1.7300,1.6404) -- (-1.6261,1.7023) -- cycle;
\path[cbpmesh=16.1] (1.6435,1.5926) -- (1.5447,1.6514) -- (1.6261,1.7023) -- (1.7300,1.6404) -- cycle;
\path[cbpmesh=30.8] (-0.9733,0.8777) -- (-1.0467,0.8404) -- (-1.0759,0.9568) -- (-1.0005,0.9951) -- cycle;
\path[cbpmesh=19.1] (1.0467,0.8404) -- (0.9733,0.8777) -- (1.0005,0.9951) -- (1.0759,0.9568) -- cycle;
\path[cbpmesh=18.3] (-2.5390,1.5425) -- (-2.6326,1.4471) -- (-2.7343,1.4356) -- (-2.6370,1.5347) -- cycle;
\path[cbpmesh=23.0] (2.6326,1.4471) -- (2.5390,1.5425) -- (2.6370,1.5347) -- (2.7343,1.4356) -- cycle;
\path[cbpmesh=38] (-1.1380,0.2481) -- (-1.1892,0.2052) -- (-1.1826,0.3403) -- (-1.1316,0.3830) -- cycle;
\path[cbpmesh=26.3] (1.1892,0.2052) -- (1.1380,0.2481) -- (1.1316,0.3830) -- (1.1826,0.3403) -- cycle;
\path[cbpmesh=38] (-1.2073,0.0654) -- (-1.2519,0.0201) -- (-1.2331,0.1605) -- (-1.1892,0.2052) -- cycle;
\path[cbpmesh=28.4] (1.2519,0.0201) -- (1.2073,0.0654) -- (1.1892,0.2052) -- (1.2331,0.1605) -- cycle;
\path[cbpmesh=38] (-1.2958,-0.1076) -- (-1.3339,-0.1561) -- (-1.2887,-0.0268) -- (-1.2519,0.0201) -- cycle;
\path[cbpmesh=32.7] (1.3339,-0.1561) -- (1.2958,-0.1076) -- (1.2519,0.0201) -- (1.2887,-0.0268) -- cycle;
\path[cbpmesh=24.4] (-4.2853,0.8976) -- (-4.3528,0.7382) -- (-4.4129,0.6221) -- (-4.3444,0.7836) -- cycle;
\draw[draw=cbpInk!70,line width=.45pt] (-4.4129,0.6221) -- (-4.3444,0.7836);
\draw[draw=cbpInk!70,line width=.45pt] (-4.3444,0.7836) -- (-4.2853,0.8976);
\path[cbpmesh=38] (4.3528,0.7382) -- (4.2853,0.8976) -- (4.3444,0.7836) -- (4.4129,0.6221) -- cycle;
\draw[draw=cbpInk!70,line width=.45pt] (4.2853,0.8976) -- (4.3444,0.7836);
\draw[draw=cbpInk!70,line width=.45pt] (4.3444,0.7836) -- (4.4129,0.6221);
\path[cbpmesh=20.7] (-1.2335,1.4937) -- (-1.3265,1.4465) -- (-1.3945,1.5230) -- (-1.2967,1.5727) -- cycle;
\path[cbpmesh=14.1] (1.3265,1.4465) -- (1.2335,1.4937) -- (1.2967,1.5727) -- (1.3945,1.5230) -- cycle;
\path[cbpmesh=18.3] (-1.8204,1.6787) -- (-1.9185,1.6097) -- (-2.0173,1.6346) -- (-1.9140,1.7071) -- cycle;
\path[cbpmesh=18.8] (1.9185,1.6097) -- (1.8204,1.6787) -- (1.9140,1.7071) -- (2.0173,1.6346) -- cycle;
\path[cbpmesh=38] (-1.3895,-0.2618) -- (-1.4206,-0.3136) -- (-1.3638,-0.2058) -- (-1.3339,-0.1561) -- cycle;
\draw[draw=cbpInk!70,line width=.45pt] (-1.3895,-0.2618) -- (-1.4206,-0.3136);
\draw[draw=cbpInk!70,line width=.45pt] (-1.4206,-0.3136) -- (-1.3638,-0.2058);
\path[cbpmesh=36.1] (1.4206,-0.3136) -- (1.3895,-0.2618) -- (1.3339,-0.1561) -- (1.3638,-0.2058) -- cycle;
\draw[draw=cbpInk!70,line width=.45pt] (1.4206,-0.3136) -- (1.3895,-0.2618);
\draw[draw=cbpInk!70,line width=.45pt] (1.3638,-0.2058) -- (1.4206,-0.3136);
\path[cbpmesh=18.0] (-3.7791,1.3428) -- (-3.8638,1.2018) -- (-3.9285,1.1813) -- (-3.8424,1.3247) -- cycle;
\path[cbpmesh=24.9] (3.8638,1.2018) -- (3.7791,1.3428) -- (3.8424,1.3247) -- (3.9285,1.1813) -- cycle;
\path[cbpmesh=28.9] (-1.0005,0.9951) -- (-1.0759,0.9568) -- (-1.1126,1.0676) -- (-1.0346,1.1072) -- cycle;
\path[cbpmesh=17.7] (1.0759,0.9568) -- (1.0005,0.9951) -- (1.0346,1.1072) -- (1.1126,1.0676) -- cycle;
\path[cbpmesh=18.0] (-3.0031,1.4334) -- (-3.0914,1.3210) -- (-3.1730,1.3248) -- (-3.0824,1.4402) -- cycle;
\path[cbpmesh=19.7] (3.0914,1.3210) -- (3.0031,1.4334) -- (3.0824,1.4402) -- (3.1730,1.3248) -- cycle;
\path[cbpmesh=21.9] (2.2143,1.5691) -- (2.1190,1.6490) -- (2.2232,1.6526) -- (2.3232,1.5687) -- cycle;
\path[cbpmesh=18.4] (-2.1190,1.6490) -- (-2.2143,1.5691) -- (-2.3232,1.5687) -- (-2.2232,1.6526) -- cycle;
\path[cbpmesh=33.7] (-1.0256,0.7189) -- (-1.0912,0.6799) -- (-1.1136,0.8006) -- (-1.0467,0.8404) -- cycle;
\path[cbpmesh=20.8] (1.0912,0.6799) -- (1.0256,0.7189) -- (1.0467,0.8404) -- (1.1136,0.8006) -- cycle;
\path[cbpmesh=14.2] (1.2639,1.3622) -- (1.1753,1.4072) -- (1.2335,1.4937) -- (1.3265,1.4465) -- cycle;
\path[cbpmesh=22.0] (-1.1753,1.4072) -- (-1.2639,1.3622) -- (-1.3265,1.4465) -- (-1.2335,1.4937) -- cycle;
\path[cbpmesh=19.2] (-1.4674,1.5915) -- (-1.5612,1.5356) -- (-1.6435,1.5926) -- (-1.5447,1.6514) -- cycle;
\path[cbpmesh=15.2] (1.5612,1.5356) -- (1.4674,1.5915) -- (1.5447,1.6514) -- (1.6435,1.5926) -- cycle;
\path[cbpmesh=27.1] (-1.0346,1.1072) -- (-1.1126,1.0676) -- (-1.1565,1.1724) -- (-1.0754,1.2135) -- cycle;
\path[cbpmesh=16.4] (1.1126,1.0676) -- (1.0346,1.1072) -- (1.0754,1.2135) -- (1.1565,1.1724) -- cycle;
\path[cbpmesh=14.6] (1.2071,1.2707) -- (1.1225,1.3137) -- (1.1753,1.4072) -- (1.2639,1.3622) -- cycle;
\path[cbpmesh=23.5] (-1.1225,1.3137) -- (-1.2071,1.2707) -- (-1.2639,1.3622) -- (-1.1753,1.4072) -- cycle;
\path[cbpmesh=25.2] (-1.0754,1.2135) -- (-1.1565,1.1724) -- (-1.2071,1.2707) -- (-1.1225,1.3137) -- cycle;
\path[cbpmesh=15.4] (1.1565,1.1724) -- (1.0754,1.2135) -- (1.1225,1.3137) -- (1.2071,1.2707) -- cycle;
\path[cbpmesh=17.7] (-3.7165,1.3542) -- (-3.7998,1.2155) -- (-3.8638,1.2018) -- (-3.7791,1.3428) -- cycle;
\path[cbpmesh=23.2] (3.7998,1.2155) -- (3.7165,1.3542) -- (3.7791,1.3428) -- (3.8638,1.2018) -- cycle;
\path[cbpmesh=36.6] (-1.0777,0.5541) -- (-1.1358,0.5133) -- (-1.1500,0.6385) -- (-1.0912,0.6799) -- cycle;
\path[cbpmesh=22.7] (1.1358,0.5133) -- (1.0777,0.5541) -- (1.0912,0.6799) -- (1.1500,0.6385) -- cycle;
\path[cbpmesh=17.9] (-2.9190,1.4297) -- (-3.0048,1.3205) -- (-3.0914,1.3210) -- (-3.0031,1.4334) -- cycle;
\path[cbpmesh=20.3] (3.0048,1.3205) -- (2.9190,1.4297) -- (3.0031,1.4334) -- (3.0914,1.3210) -- cycle;
\path[cbpmesh=18.5] (-2.4342,1.5571) -- (-2.5240,1.4656) -- (-2.6326,1.4471) -- (-2.5390,1.5425) -- cycle;
\path[cbpmesh=23.9] (2.5240,1.4656) -- (2.4342,1.5571) -- (2.5390,1.5425) -- (2.6326,1.4471) -- cycle;
\path[cbpmesh=18.4] (-1.7300,1.6404) -- (-1.8233,1.5748) -- (-1.9185,1.6097) -- (-1.8204,1.6787) -- cycle;
\path[cbpmesh=17.4] (1.8233,1.5748) -- (1.7300,1.6404) -- (1.8204,1.6787) -- (1.9185,1.6097) -- cycle;
\path[cbpmesh=31.9] (-1.0467,0.8404) -- (-1.1136,0.8006) -- (-1.1446,0.9159) -- (-1.0759,0.9568) -- cycle;
\path[cbpmesh=19.3] (1.1136,0.8006) -- (1.0467,0.8404) -- (1.0759,0.9568) -- (1.1446,0.9159) -- cycle;
\path[cbpmesh=38] (-1.1316,0.3830) -- (-1.1826,0.3403) -- (-1.1869,0.4704) -- (-1.1358,0.5133) -- cycle;
\path[cbpmesh=24.7] (1.1826,0.3403) -- (1.1316,0.3830) -- (1.1358,0.5133) -- (1.1869,0.4704) -- cycle;
\path[cbpmesh=17.6] (-3.6540,1.3599) -- (-3.7359,1.2236) -- (-3.7998,1.2155) -- (-3.7165,1.3542) -- cycle;
\path[cbpmesh=21.8] (3.7359,1.2236) -- (3.6540,1.3599) -- (3.7165,1.3542) -- (3.7998,1.2155) -- cycle;
\path[cbpmesh=20.1] (-1.3945,1.5230) -- (-1.4836,1.4699) -- (-1.5612,1.5356) -- (-1.4674,1.5915) -- cycle;
\path[cbpmesh=14.5] (1.4836,1.4699) -- (1.3945,1.5230) -- (1.4674,1.5915) -- (1.5612,1.5356) -- cycle;
\path[cbpmesh=22.7] (-4.2055,1.0122) -- (-4.2718,0.8558) -- (-4.3528,0.7382) -- (-4.2853,0.8976) -- cycle;
\path[cbpmesh=37.0] (4.2718,0.8558) -- (4.2055,1.0122) -- (4.2853,0.8976) -- (4.3528,0.7382) -- cycle;
\path[cbpmesh=26.8] (1.2331,0.1605) -- (1.1892,0.2052) -- (1.1826,0.3403) -- (1.2262,0.2959) -- cycle;
\path[cbpmesh=38] (-1.1892,0.2052) -- (-1.2331,0.1605) -- (-1.2262,0.2959) -- (-1.1826,0.3403) -- cycle;
\path[cbpmesh=38] (-1.2519,0.0201) -- (-1.2887,-0.0268) -- (-1.2693,0.1144) -- (-1.2331,0.1605) -- cycle;
\path[cbpmesh=28.9] (1.2887,-0.0268) -- (1.2519,0.0201) -- (1.2331,0.1605) -- (1.2693,0.1144) -- cycle;
\path[cbpmesh=18.2] (-2.0173,1.6346) -- (-2.1080,1.5585) -- (-2.2143,1.5691) -- (-2.1190,1.6490) -- cycle;
\path[cbpmesh=20.3] (2.1080,1.5585) -- (2.0173,1.6346) -- (2.1190,1.6490) -- (2.2143,1.5691) -- cycle;
\path[cbpmesh=17.9] (-2.8296,1.4301) -- (-2.9128,1.3242) -- (-3.0048,1.3205) -- (-2.9190,1.4297) -- cycle;
\path[cbpmesh=21.1] (2.9128,1.3242) -- (2.8296,1.4301) -- (2.9190,1.4297) -- (3.0048,1.3205) -- cycle;
\path[cbpmesh=38] (-1.3339,-0.1561) -- (-1.3638,-0.2058) -- (-1.3176,-0.0749) -- (-1.2887,-0.0268) -- cycle;
\path[cbpmesh=33.3] (1.3638,-0.2058) -- (1.3339,-0.1561) -- (1.2887,-0.0268) -- (1.3176,-0.0749) -- cycle;
\path[cbpmesh=17.7] (-3.5911,1.3610) -- (-3.6716,1.2270) -- (-3.7359,1.2236) -- (-3.6540,1.3599) -- cycle;
\path[cbpmesh=20.6] (3.6716,1.2270) -- (3.5911,1.3610) -- (3.6540,1.3599) -- (3.7359,1.2236) -- cycle;
\path[cbpmesh=29.9] (-1.0759,0.9568) -- (-1.1446,0.9159) -- (-1.1837,1.0252) -- (-1.1126,1.0676) -- cycle;
\path[cbpmesh=17.8] (1.1446,0.9159) -- (1.0759,0.9568) -- (1.1126,1.0676) -- (1.1837,1.0252) -- cycle;
\path[cbpmesh=21.2] (-1.3265,1.4465) -- (-1.4113,1.3960) -- (-1.4836,1.4699) -- (-1.3945,1.5230) -- cycle;
\path[cbpmesh=14.2] (1.4113,1.3960) -- (1.3265,1.4465) -- (1.3945,1.5230) -- (1.4836,1.4699) -- cycle;
\path[cbpmesh=21.1] (1.1500,0.6385) -- (1.0912,0.6799) -- (1.1136,0.8006) -- (1.1737,0.7584) -- cycle;
\path[cbpmesh=34.8] (-1.0912,0.6799) -- (-1.1500,0.6385) -- (-1.1737,0.7584) -- (-1.1136,0.8006) -- cycle;
\path[cbpmesh=18.8] (-1.6435,1.5926) -- (-1.7321,1.5303) -- (-1.8233,1.5748) -- (-1.7300,1.6404) -- cycle;
\path[cbpmesh=16.2] (1.7321,1.5303) -- (1.6435,1.5926) -- (1.7300,1.6404) -- (1.8233,1.5748) -- cycle;
\path[cbpmesh=17.8] (-3.5272,1.3584) -- (-3.6062,1.2268) -- (-3.6716,1.2270) -- (-3.5911,1.3610) -- cycle;
\path[cbpmesh=19.8] (3.6062,1.2268) -- (3.5272,1.3584) -- (3.5911,1.3610) -- (3.6716,1.2270) -- cycle;
\path[cbpmesh=28.0] (-1.1126,1.0676) -- (-1.1837,1.0252) -- (-1.2304,1.1283) -- (-1.1565,1.1724) -- cycle;
\path[cbpmesh=16.6] (1.1837,1.0252) -- (1.1126,1.0676) -- (1.1565,1.1724) -- (1.2304,1.1283) -- cycle;
\path[cbpmesh=22.7] (-1.2639,1.3622) -- (-1.3447,1.3141) -- (-1.4113,1.3960) -- (-1.3265,1.4465) -- cycle;
\path[cbpmesh=14.3] (1.3447,1.3141) -- (1.2639,1.3622) -- (1.3265,1.4465) -- (1.4113,1.3960) -- cycle;
\path[cbpmesh=18.4] (-2.3232,1.5687) -- (-2.4089,1.4814) -- (-2.5240,1.4656) -- (-2.4342,1.5571) -- cycle;
\path[cbpmesh=23.4] (2.4089,1.4814) -- (2.3232,1.5687) -- (2.4342,1.5571) -- (2.5240,1.4656) -- cycle;
\path[cbpmesh=15.5] (1.2304,1.1283) -- (1.1565,1.1724) -- (1.2071,1.2707) -- (1.2842,1.2248) -- cycle;
\path[cbpmesh=26.1] (-1.1565,1.1724) -- (-1.2304,1.1283) -- (-1.2842,1.2248) -- (-1.2071,1.2707) -- cycle;
\path[cbpmesh=24.3] (-1.2071,1.2707) -- (-1.2842,1.2248) -- (-1.3447,1.3141) -- (-1.2639,1.3622) -- cycle;
\path[cbpmesh=14.8] (1.2842,1.2248) -- (1.2071,1.2707) -- (1.2639,1.3622) -- (1.3447,1.3141) -- cycle;
\path[cbpmesh=23.1] (1.1869,0.4704) -- (1.1358,0.5133) -- (1.1500,0.6385) -- (1.2017,0.5952) -- cycle;
\path[cbpmesh=37.7] (-1.1358,0.5133) -- (-1.1869,0.4704) -- (-1.2017,0.5952) -- (-1.1500,0.6385) -- cycle;
\path[cbpmesh=17.9] (-2.7343,1.4356) -- (-2.8147,1.3332) -- (-2.9128,1.3242) -- (-2.8296,1.4301) -- cycle;
\path[cbpmesh=21.9] (2.8147,1.3332) -- (2.7343,1.4356) -- (2.8296,1.4301) -- (2.9128,1.3242) -- cycle;
\path[cbpmesh=18.0] (-3.4617,1.3531) -- (-3.5393,1.2239) -- (-3.6062,1.2268) -- (-3.5272,1.3584) -- cycle;
\path[cbpmesh=19.2] (3.5393,1.2239) -- (3.4617,1.3531) -- (3.5272,1.3584) -- (3.6062,1.2268) -- cycle;
\path[cbpmesh=19.7] (-4.1324,1.0696) -- (-4.1975,0.9159) -- (-4.2718,0.8558) -- (-4.2055,1.0122) -- cycle;
\path[cbpmesh=31.3] (4.1975,0.9159) -- (4.1324,1.0696) -- (4.2055,1.0122) -- (4.2718,0.8558) -- cycle;
\path[cbpmesh=18.2] (-1.9185,1.6097) -- (-2.0049,1.5374) -- (-2.1080,1.5585) -- (-2.0173,1.6346) -- cycle;
\path[cbpmesh=18.7] (2.0049,1.5374) -- (1.9185,1.6097) -- (2.0173,1.6346) -- (2.1080,1.5585) -- cycle;
\path[cbpmesh=32.9] (-1.1136,0.8006) -- (-1.1737,0.7584) -- (-1.2064,0.8725) -- (-1.1446,0.9159) -- cycle;
\path[cbpmesh=19.6] (1.1737,0.7584) -- (1.1136,0.8006) -- (1.1446,0.9159) -- (1.2064,0.8725) -- cycle;
\path[cbpmesh=19.5] (-1.5612,1.5356) -- (-1.6454,1.4765) -- (-1.7321,1.5303) -- (-1.6435,1.5926) -- cycle;
\path[cbpmesh=15.3] (1.6454,1.4765) -- (1.5612,1.5356) -- (1.6435,1.5926) -- (1.7321,1.5303) -- cycle;
\path[cbpmesh=38] (-1.1826,0.3403) -- (-1.2262,0.2959) -- (-1.2307,0.4258) -- (-1.1869,0.4704) -- cycle;
\path[cbpmesh=25.2] (1.2262,0.2959) -- (1.1826,0.3403) -- (1.1869,0.4704) -- (1.2307,0.4258) -- cycle;
\path[cbpmesh=18.1] (-3.3941,1.3460) -- (-3.4702,1.2194) -- (-3.5393,1.2239) -- (-3.4617,1.3531) -- cycle;
\path[cbpmesh=18.9] (3.4702,1.2194) -- (3.3941,1.3460) -- (3.4617,1.3531) -- (3.5393,1.2239) -- cycle;
\path[cbpmesh=38] (-1.2331,0.1605) -- (-1.2693,0.1144) -- (-1.2623,0.2500) -- (-1.2262,0.2959) -- cycle;
\path[cbpmesh=27.3] (1.2693,0.1144) -- (1.2331,0.1605) -- (1.2262,0.2959) -- (1.2623,0.2500) -- cycle;
\path[cbpmesh=18.1] (-3.3239,1.3383) -- (-3.3983,1.2143) -- (-3.4702,1.2194) -- (-3.3941,1.3460) -- cycle;
\path[cbpmesh=18.8] (3.3983,1.2143) -- (3.3239,1.3383) -- (3.3941,1.3460) -- (3.4702,1.2194) -- cycle;
\path[cbpmesh=38] (-1.2887,-0.0268) -- (-1.3176,-0.0749) -- (-1.2978,0.0670) -- (-1.2693,0.1144) -- cycle;
\path[cbpmesh=29.5] (1.3176,-0.0749) -- (1.2887,-0.0268) -- (1.2693,0.1144) -- (1.2978,0.0670) -- cycle;
\path[cbpmesh=30.9] (-1.1446,0.9159) -- (-1.2064,0.8725) -- (-1.2476,0.9804) -- (-1.1837,1.0252) -- cycle;
\path[cbpmesh=18.1] (1.2064,0.8725) -- (1.1446,0.9159) -- (1.1837,1.0252) -- (1.2476,0.9804) -- cycle;
\path[cbpmesh=14.7] (1.5636,1.4137) -- (1.4836,1.4699) -- (1.5612,1.5356) -- (1.6454,1.4765) -- cycle;
\path[cbpmesh=20.5] (-1.4836,1.4699) -- (-1.5636,1.4137) -- (-1.6454,1.4765) -- (-1.5612,1.5356) -- cycle;
\path[cbpmesh=18.0] (-2.6326,1.4471) -- (-2.7101,1.3485) -- (-2.8147,1.3332) -- (-2.7343,1.4356) -- cycle;
\path[cbpmesh=22.8] (2.7101,1.3485) -- (2.6326,1.4471) -- (2.7343,1.4356) -- (2.8147,1.3332) -- cycle;
\path[cbpmesh=18.1] (-2.2143,1.5691) -- (-2.2960,1.4859) -- (-2.4089,1.4814) -- (-2.3232,1.5687) -- cycle;
\path[cbpmesh=21.8] (2.2960,1.4859) -- (2.2143,1.5691) -- (2.3232,1.5687) -- (2.4089,1.4814) -- cycle;
\path[cbpmesh=35.8] (-1.1500,0.6385) -- (-1.2017,0.5952) -- (-1.2265,0.7141) -- (-1.1737,0.7584) -- cycle;
\path[cbpmesh=21.5] (1.2017,0.5952) -- (1.1500,0.6385) -- (1.1737,0.7584) -- (1.2265,0.7141) -- cycle;
\path[cbpmesh=38] (-1.3638,-0.2058) -- (-1.3853,-0.2565) -- (-1.3383,-0.1239) -- (-1.3176,-0.0749) -- cycle;
\draw[draw=cbpInk!70,line width=.45pt] (-1.3638,-0.2058) -- (-1.3853,-0.2565);
\draw[draw=cbpInk!70,line width=.45pt] (-1.3853,-0.2565) -- (-1.3383,-0.1239);
\path[cbpmesh=33.9] (1.3853,-0.2565) -- (1.3638,-0.2058) -- (1.3176,-0.0749) -- (1.3383,-0.1239) -- cycle;
\draw[draw=cbpInk!70,line width=.45pt] (1.3853,-0.2565) -- (1.3638,-0.2058);
\draw[draw=cbpInk!70,line width=.45pt] (1.3383,-0.1239) -- (1.3853,-0.2565);
\path[cbpmesh=19.0] (-4.0622,1.1162) -- (-4.1262,0.9652) -- (-4.1975,0.9159) -- (-4.1324,1.0696) -- cycle;
\path[cbpmesh=29.8] (4.1262,0.9652) -- (4.0622,1.1162) -- (4.1324,1.0696) -- (4.1975,0.9159) -- cycle;
\path[cbpmesh=17.4] (1.9053,1.5061) -- (1.8233,1.5748) -- (1.9185,1.6097) -- (2.0049,1.5374) -- cycle;
\path[cbpmesh=18.5] (-1.8233,1.5748) -- (-1.9053,1.5061) -- (-2.0049,1.5374) -- (-1.9185,1.6097) -- cycle;
\path[cbpmesh=18.0] (-3.2504,1.3309) -- (-3.3232,1.2097) -- (-3.3983,1.2143) -- (-3.3239,1.3383) -- cycle;
\path[cbpmesh=18.9] (3.3232,1.2097) -- (3.2504,1.3309) -- (3.3239,1.3383) -- (3.3983,1.2143) -- cycle;
\path[cbpmesh=16.8] (1.2476,0.9804) -- (1.1837,1.0252) -- (1.2304,1.1283) -- (1.2968,1.0817) -- cycle;
\path[cbpmesh=28.9] (-1.1837,1.0252) -- (-1.2476,0.9804) -- (-1.2968,1.0817) -- (-1.2304,1.1283) -- cycle;
\path[cbpmesh=21.8] (-1.4113,1.3960) -- (-1.4874,1.3425) -- (-1.5636,1.4137) -- (-1.4836,1.4699) -- cycle;
\path[cbpmesh=14.4] (1.4874,1.3425) -- (1.4113,1.3960) -- (1.4836,1.4699) -- (1.5636,1.4137) -- cycle;
\path[cbpmesh=23.3] (-1.3447,1.3141) -- (-1.4172,1.2632) -- (-1.4874,1.3425) -- (-1.4113,1.3960) -- cycle;
\path[cbpmesh=14.5] (1.4172,1.2632) -- (1.3447,1.3141) -- (1.4113,1.3960) -- (1.4874,1.3425) -- cycle;
\path[cbpmesh=26.9] (-1.2304,1.1283) -- (-1.2968,1.0817) -- (-1.3535,1.1761) -- (-1.2842,1.2248) -- cycle;
\path[cbpmesh=15.8] (1.2968,1.0817) -- (1.2304,1.1283) -- (1.2842,1.2248) -- (1.3535,1.1761) -- cycle;
\path[cbpmesh=38] (-1.1869,0.4704) -- (-1.2307,0.4258) -- (-1.2461,0.5500) -- (-1.2017,0.5952) -- cycle;
\path[cbpmesh=23.6] (1.2307,0.4258) -- (1.1869,0.4704) -- (1.2017,0.5952) -- (1.2461,0.5500) -- cycle;
\path[cbpmesh=17.9] (-3.1730,1.3248) -- (-3.2442,1.2065) -- (-3.3232,1.2097) -- (-3.2504,1.3309) -- cycle;
\path[cbpmesh=19.2] (3.2442,1.2065) -- (3.1730,1.3248) -- (3.2504,1.3309) -- (3.3232,1.2097) -- cycle;
\path[cbpmesh=25.1] (-1.2842,1.2248) -- (-1.3535,1.1761) -- (-1.4172,1.2632) -- (-1.3447,1.3141) -- cycle;
\path[cbpmesh=15.0] (1.3535,1.1761) -- (1.2842,1.2248) -- (1.3447,1.3141) -- (1.4172,1.2632) -- cycle;
\path[cbpmesh=18.4] (-3.9944,1.1532) -- (-4.0573,1.0046) -- (-4.1262,0.9652) -- (-4.0622,1.1162) -- cycle;
\path[cbpmesh=28.1] (4.0573,1.0046) -- (3.9944,1.1532) -- (4.0622,1.1162) -- (4.1262,0.9652) -- cycle;
\path[cbpmesh=33.9] (-1.1737,0.7584) -- (-1.2265,0.7141) -- (-1.2606,0.8270) -- (-1.2064,0.8725) -- cycle;
\path[cbpmesh=19.9] (1.2265,0.7141) -- (1.1737,0.7584) -- (1.2064,0.8725) -- (1.2606,0.8270) -- cycle;
\path[cbpmesh=19.0] (-1.7321,1.5303) -- (-1.8100,1.4650) -- (-1.9053,1.5061) -- (-1.8233,1.5748) -- cycle;
\path[cbpmesh=16.3] (1.8100,1.4650) -- (1.7321,1.5303) -- (1.8233,1.5748) -- (1.9053,1.5061) -- cycle;
\path[cbpmesh=18.1] (-2.5240,1.4656) -- (-2.5983,1.3711) -- (-2.7101,1.3485) -- (-2.6326,1.4471) -- cycle;
\path[cbpmesh=23.7] (2.5983,1.3711) -- (2.5240,1.4656) -- (2.6326,1.4471) -- (2.7101,1.3485) -- cycle;
\path[cbpmesh=38] (-1.2262,0.2959) -- (-1.2623,0.2500) -- (-1.2669,0.3798) -- (-1.2307,0.4258) -- cycle;
\path[cbpmesh=25.7] (1.2623,0.2500) -- (1.2262,0.2959) -- (1.2307,0.4258) -- (1.2669,0.3798) -- cycle;
\path[cbpmesh=25.5] (-4.4129,0.6221) -- (-4.4541,0.4585) -- (-4.4947,0.3428) -- (-4.4531,0.5079) -- cycle;
\draw[draw=cbpInk!70,line width=.45pt] (-4.4947,0.3428) -- (-4.4531,0.5079);
\draw[draw=cbpInk!70,line width=.45pt] (-4.4531,0.5079) -- (-4.4129,0.6221);
\path[cbpmesh=38] (4.4541,0.4585) -- (4.4129,0.6221) -- (4.4531,0.5079) -- (4.4947,0.3428) -- cycle;
\draw[draw=cbpInk!70,line width=.45pt] (4.4129,0.6221) -- (4.4531,0.5079);
\draw[draw=cbpInk!70,line width=.45pt] (4.4531,0.5079) -- (4.4947,0.3428);
\path[cbpmesh=17.8] (-3.0914,1.3210) -- (-3.1607,1.2057) -- (-3.2442,1.2065) -- (-3.1730,1.3248) -- cycle;
\path[cbpmesh=19.7] (3.1607,1.2057) -- (3.0914,1.3210) -- (3.1730,1.3248) -- (3.2442,1.2065) -- cycle;
\path[cbpmesh=18.1] (-2.1080,1.5585) -- (-2.1858,1.4793) -- (-2.2960,1.4859) -- (-2.2143,1.5691) -- cycle;
\path[cbpmesh=20.2] (2.1858,1.4793) -- (2.1080,1.5585) -- (2.2143,1.5691) -- (2.2960,1.4859) -- cycle;
\path[cbpmesh=38] (-1.2693,0.1144) -- (-1.2978,0.0670) -- (-1.2905,0.2029) -- (-1.2623,0.2500) -- cycle;
\path[cbpmesh=27.9] (1.2978,0.0670) -- (1.2693,0.1144) -- (1.2623,0.2500) -- (1.2905,0.2029) -- cycle;
\path[cbpmesh=17.9] (-3.9285,1.1813) -- (-3.9904,1.0353) -- (-4.0573,1.0046) -- (-3.9944,1.1532) -- cycle;
\path[cbpmesh=26.4] (3.9904,1.0353) -- (3.9285,1.1813) -- (3.9944,1.1532) -- (4.0573,1.0046) -- cycle;
\path[cbpmesh=18.5] (1.2606,0.8270) -- (1.2064,0.8725) -- (1.2476,0.9804) -- (1.3037,0.9333) -- cycle;
\path[cbpmesh=31.9] (-1.2064,0.8725) -- (-1.2606,0.8270) -- (-1.3037,0.9333) -- (-1.2476,0.9804) -- cycle;
\path[cbpmesh=19.9] (-1.6454,1.4765) -- (-1.7194,1.4144) -- (-1.8100,1.4650) -- (-1.7321,1.5303) -- cycle;
\path[cbpmesh=15.4] (1.7194,1.4144) -- (1.6454,1.4765) -- (1.7321,1.5303) -- (1.8100,1.4650) -- cycle;
\path[cbpmesh=17.7] (-3.0048,1.3205) -- (-3.0722,1.2084) -- (-3.1607,1.2057) -- (-3.0914,1.3210) -- cycle;
\path[cbpmesh=20.3] (3.0722,1.2084) -- (3.0048,1.3205) -- (3.0914,1.3210) -- (3.1607,1.2057) -- cycle;
\path[cbpmesh=38] (-1.3176,-0.0749) -- (-1.3383,-0.1239) -- (-1.3182,0.0188) -- (-1.2978,0.0670) -- cycle;
\path[cbpmesh=30.2] (1.3383,-0.1239) -- (1.3176,-0.0749) -- (1.2978,0.0670) -- (1.3182,0.0188) -- cycle;
\path[cbpmesh=21.9] (1.2461,0.5500) -- (1.2017,0.5952) -- (1.2265,0.7141) -- (1.2717,0.6681) -- cycle;
\path[cbpmesh=36.9] (-1.2017,0.5952) -- (-1.2461,0.5500) -- (-1.2717,0.6681) -- (-1.2265,0.7141) -- cycle;
\path[cbpmesh=18.2] (-2.0049,1.5374) -- (-2.0788,1.4620) -- (-2.1858,1.4793) -- (-2.1080,1.5585) -- cycle;
\path[cbpmesh=18.7] (2.0788,1.4620) -- (2.0049,1.5374) -- (2.1080,1.5585) -- (2.1858,1.4793) -- cycle;
\path[cbpmesh=17.6] (-3.8638,1.2018) -- (-3.9247,1.0581) -- (-3.9904,1.0353) -- (-3.9285,1.1813) -- cycle;
\path[cbpmesh=24.6] (3.9247,1.0581) -- (3.8638,1.2018) -- (3.9285,1.1813) -- (3.9904,1.0353) -- cycle;
\path[cbpmesh=29.8] (-1.2476,0.9804) -- (-1.3037,0.9333) -- (-1.3551,1.0328) -- (-1.2968,1.0817) -- cycle;
\path[cbpmesh=17.2] (1.3037,0.9333) -- (1.2476,0.9804) -- (1.2968,1.0817) -- (1.3551,1.0328) -- cycle;
\path[cbpmesh=14.8] (1.6340,1.3548) -- (1.5636,1.4137) -- (1.6454,1.4765) -- (1.7194,1.4144) -- cycle;
\path[cbpmesh=21.0] (-1.5636,1.4137) -- (-1.6340,1.3548) -- (-1.7194,1.4144) -- (-1.6454,1.4765) -- cycle;
\path[cbpmesh=18.0] (-2.4089,1.4814) -- (-2.4798,1.3912) -- (-2.5983,1.3711) -- (-2.5240,1.4656) -- cycle;
\path[cbpmesh=23.2] (2.4798,1.3912) -- (2.4089,1.4814) -- (2.5240,1.4656) -- (2.5983,1.3711) -- cycle;
\path[cbpmesh=17.7] (-2.9128,1.3242) -- (-2.9781,1.2156) -- (-3.0722,1.2084) -- (-3.0048,1.3205) -- cycle;
\path[cbpmesh=21.0] (2.9781,1.2156) -- (2.9128,1.3242) -- (3.0048,1.3205) -- (3.0722,1.2084) -- cycle;
\path[cbpmesh=38] (-1.2307,0.4258) -- (-1.2669,0.3798) -- (-1.2827,0.5034) -- (-1.2461,0.5500) -- cycle;
\path[cbpmesh=24.1] (1.2669,0.3798) -- (1.2307,0.4258) -- (1.2461,0.5500) -- (1.2827,0.5034) -- cycle;
\path[cbpmesh=27.8] (-1.2968,1.0817) -- (-1.3551,1.0328) -- (-1.4144,1.1251) -- (-1.3535,1.1761) -- cycle;
\path[cbpmesh=16.1] (1.3551,1.0328) -- (1.2968,1.0817) -- (1.3535,1.1761) -- (1.4144,1.1251) -- cycle;
\path[cbpmesh=22.4] (-1.4874,1.3425) -- (-1.5543,1.2864) -- (-1.6340,1.3548) -- (-1.5636,1.4137) -- cycle;
\path[cbpmesh=14.6] (1.5543,1.2864) -- (1.4874,1.3425) -- (1.5636,1.4137) -- (1.6340,1.3548) -- cycle;
\path[cbpmesh=23.5] (-4.3528,0.7382) -- (-4.3935,0.5769) -- (-4.4541,0.4585) -- (-4.4129,0.6221) -- cycle;
\path[cbpmesh=38] (4.3935,0.5769) -- (4.3528,0.7382) -- (4.4129,0.6221) -- (4.4541,0.4585) -- cycle;
\path[cbpmesh=17.4] (-3.7998,1.2155) -- (-3.8597,1.0743) -- (-3.9247,1.0581) -- (-3.8638,1.2018) -- cycle;
\path[cbpmesh=23.0] (3.8597,1.0743) -- (3.7998,1.2155) -- (3.8638,1.2018) -- (3.9247,1.0581) -- cycle;
\path[cbpmesh=34.9] (-1.2265,0.7141) -- (-1.2717,0.6681) -- (-1.3072,0.7796) -- (-1.2606,0.8270) -- cycle;
\path[cbpmesh=20.4] (1.2717,0.6681) -- (1.2265,0.7141) -- (1.2606,0.8270) -- (1.3072,0.7796) -- cycle;
\path[cbpmesh=25.8] (-1.3535,1.1761) -- (-1.4144,1.1251) -- (-1.4810,1.2097) -- (-1.4172,1.2632) -- cycle;
\path[cbpmesh=15.2] (1.4144,1.1251) -- (1.3535,1.1761) -- (1.4172,1.2632) -- (1.4810,1.2097) -- cycle;
\path[cbpmesh=24.0] (-1.4172,1.2632) -- (-1.4810,1.2097) -- (-1.5543,1.2864) -- (-1.4874,1.3425) -- cycle;
\path[cbpmesh=14.7] (1.4810,1.2097) -- (1.4172,1.2632) -- (1.4874,1.3425) -- (1.5543,1.2864) -- cycle;
\path[cbpmesh=38] (-1.2623,0.2500) -- (-1.2905,0.2029) -- (-1.2953,0.3325) -- (-1.2669,0.3798) -- cycle;
\path[cbpmesh=26.3] (1.2905,0.2029) -- (1.2623,0.2500) -- (1.2669,0.3798) -- (1.2953,0.3325) -- cycle;
\path[cbpmesh=17.4] (1.9756,1.4345) -- (1.9053,1.5061) -- (2.0049,1.5374) -- (2.0788,1.4620) -- cycle;
\path[cbpmesh=18.6] (-1.9053,1.5061) -- (-1.9756,1.4345) -- (-2.0788,1.4620) -- (-2.0049,1.5374) -- cycle;
\path[cbpmesh=17.6] (-2.8147,1.3332) -- (-2.8778,1.2282) -- (-2.9781,1.2156) -- (-2.9128,1.3242) -- cycle;
\path[cbpmesh=21.8] (2.8778,1.2282) -- (2.8147,1.3332) -- (2.9128,1.3242) -- (2.9781,1.2156) -- cycle;
\path[cbpmesh=17.4] (-3.7359,1.2236) -- (-3.7948,1.0847) -- (-3.8597,1.0743) -- (-3.7998,1.2155) -- cycle;
\path[cbpmesh=21.6] (3.7948,1.0847) -- (3.7359,1.2236) -- (3.7998,1.2155) -- (3.8597,1.0743) -- cycle;
\path[cbpmesh=38] (-1.2978,0.0670) -- (-1.3182,0.0188) -- (-1.3109,0.1549) -- (-1.2905,0.2029) -- cycle;
\path[cbpmesh=28.6] (1.3182,0.0188) -- (1.2978,0.0670) -- (1.2905,0.2029) -- (1.3109,0.1549) -- cycle;
\path[cbpmesh=17.9] (-2.2960,1.4859) -- (-2.3636,1.3999) -- (-2.4798,1.3912) -- (-2.4089,1.4814) -- cycle;
\path[cbpmesh=21.6] (2.3636,1.3999) -- (2.2960,1.4859) -- (2.4089,1.4814) -- (2.4798,1.3912) -- cycle;
\path[cbpmesh=32.8] (-1.2606,0.8270) -- (-1.3072,0.7796) -- (-1.3518,0.8843) -- (-1.3037,0.9333) -- cycle;
\path[cbpmesh=18.9] (1.3072,0.7796) -- (1.2606,0.8270) -- (1.3037,0.9333) -- (1.3518,0.8843) -- cycle;
\path[cbpmesh=17.5] (-3.6716,1.2270) -- (-3.7294,1.0905) -- (-3.7948,1.0847) -- (-3.7359,1.2236) -- cycle;
\path[cbpmesh=20.5] (3.7294,1.0905) -- (3.6716,1.2270) -- (3.7359,1.2236) -- (3.7948,1.0847) -- cycle;
\path[cbpmesh=37.9] (-1.2461,0.5500) -- (-1.2827,0.5034) -- (-1.3091,0.6204) -- (-1.2717,0.6681) -- cycle;
\path[cbpmesh=22.5] (1.2827,0.5034) -- (1.2461,0.5500) -- (1.2717,0.6681) -- (1.3091,0.6204) -- cycle;
\path[cbpmesh=38] (-1.3383,-0.1239) -- (-1.3508,-0.1735) -- (-1.3305,-0.0301) -- (-1.3182,0.0188) -- cycle;
\draw[draw=cbpInk!70,line width=.45pt] (-1.3383,-0.1239) -- (-1.3508,-0.1735);
\draw[draw=cbpInk!70,line width=.45pt] (-1.3508,-0.1735) -- (-1.3305,-0.0301);
\path[cbpmesh=30.9] (1.3508,-0.1735) -- (1.3383,-0.1239) -- (1.3182,0.0188) -- (1.3305,-0.0301) -- cycle;
\draw[draw=cbpInk!70,line width=.45pt] (1.3508,-0.1735) -- (1.3383,-0.1239);
\draw[draw=cbpInk!70,line width=.45pt] (1.3305,-0.0301) -- (1.3508,-0.1735);
\path[cbpmesh=19.2] (-1.8100,1.4650) -- (-1.8768,1.3969) -- (-1.9756,1.4345) -- (-1.9053,1.5061) -- cycle;
\path[cbpmesh=16.3] (1.8768,1.3969) -- (1.8100,1.4650) -- (1.9053,1.5061) -- (1.9756,1.4345) -- cycle;
\path[cbpmesh=17.7] (-3.6062,1.2268) -- (-3.6631,1.0927) -- (-3.7294,1.0905) -- (-3.6716,1.2270) -- cycle;
\path[cbpmesh=19.7] (3.6631,1.0927) -- (3.6062,1.2268) -- (3.6716,1.2270) -- (3.7294,1.0905) -- cycle;
\path[cbpmesh=17.7] (-2.7101,1.3485) -- (-2.7708,1.2474) -- (-2.8778,1.2282) -- (-2.8147,1.3332) -- cycle;
\path[cbpmesh=22.6] (2.7708,1.2474) -- (2.7101,1.3485) -- (2.8147,1.3332) -- (2.8778,1.2282) -- cycle;
\path[cbpmesh=21.9] (-4.2718,0.8558) -- (-4.3117,0.6975) -- (-4.3935,0.5769) -- (-4.3528,0.7382) -- cycle;
\path[cbpmesh=36.3] (4.3117,0.6975) -- (4.2718,0.8558) -- (4.3528,0.7382) -- (4.3935,0.5769) -- cycle;
\path[cbpmesh=30.7] (-1.3037,0.9333) -- (-1.3518,0.8843) -- (-1.4051,0.9819) -- (-1.3551,1.0328) -- cycle;
\path[cbpmesh=17.6] (1.3518,0.8843) -- (1.3037,0.9333) -- (1.3551,1.0328) -- (1.4051,0.9819) -- cycle;
\path[cbpmesh=17.9] (-3.5393,1.2239) -- (-3.5951,1.0923) -- (-3.6631,1.0927) -- (-3.6062,1.2268) -- cycle;
\path[cbpmesh=19.1] (3.5951,1.0923) -- (3.5393,1.2239) -- (3.6062,1.2268) -- (3.6631,1.0927) -- cycle;
\path[cbpmesh=20.2] (-1.7194,1.4144) -- (-1.7828,1.3498) -- (-1.8768,1.3969) -- (-1.8100,1.4650) -- cycle;
\path[cbpmesh=15.5] (1.7828,1.3498) -- (1.7194,1.4144) -- (1.8100,1.4650) -- (1.8768,1.3969) -- cycle;
\path[cbpmesh=38] (-1.2669,0.3798) -- (-1.2953,0.3325) -- (-1.3115,0.4555) -- (-1.2827,0.5034) -- cycle;
\path[cbpmesh=24.7] (1.2953,0.3325) -- (1.2669,0.3798) -- (1.2827,0.5034) -- (1.3115,0.4555) -- cycle;
\path[cbpmesh=17.9] (-2.1858,1.4793) -- (-2.2501,1.3975) -- (-2.3636,1.3999) -- (-2.2960,1.4859) -- cycle;
\path[cbpmesh=20.1] (2.2501,1.3975) -- (2.1858,1.4793) -- (2.2960,1.4859) -- (2.3636,1.3999) -- cycle;
\path[cbpmesh=18.0] (-3.4702,1.2194) -- (-3.5249,1.0904) -- (-3.5951,1.0923) -- (-3.5393,1.2239) -- cycle;
\path[cbpmesh=18.9] (3.5249,1.0904) -- (3.4702,1.2194) -- (3.5393,1.2239) -- (3.5951,1.0923) -- cycle;
\path[cbpmesh=28.6] (-1.3551,1.0328) -- (-1.4051,0.9819) -- (-1.4666,1.0719) -- (-1.4144,1.1251) -- cycle;
\path[cbpmesh=16.4] (1.4051,0.9819) -- (1.3551,1.0328) -- (1.4144,1.1251) -- (1.4666,1.0719) -- cycle;
\path[cbpmesh=35.9] (-1.2717,0.6681) -- (-1.3091,0.6204) -- (-1.3456,0.7307) -- (-1.3072,0.7796) -- cycle;
\path[cbpmesh=20.9] (1.3091,0.6204) -- (1.2717,0.6681) -- (1.3072,0.7796) -- (1.3456,0.7307) -- cycle;
\path[cbpmesh=15.0] (1.6943,1.2934) -- (1.6340,1.3548) -- (1.7194,1.4144) -- (1.7828,1.3498) -- cycle;
\path[cbpmesh=21.4] (-1.6340,1.3548) -- (-1.6943,1.2934) -- (-1.7828,1.3498) -- (-1.7194,1.4144) -- cycle;
\path[cbpmesh=18.0] (-3.3983,1.2143) -- (-3.4519,1.0880) -- (-3.5249,1.0904) -- (-3.4702,1.2194) -- cycle;
\path[cbpmesh=18.8] (3.4519,1.0880) -- (3.3983,1.2143) -- (3.4702,1.2194) -- (3.5249,1.0904) -- cycle;
\path[cbpmesh=38] (-1.2905,0.2029) -- (-1.3109,0.1549) -- (-1.3157,0.2844) -- (-1.2953,0.3325) -- cycle;
\path[cbpmesh=27.0] (1.3109,0.1549) -- (1.2905,0.2029) -- (1.2953,0.3325) -- (1.3157,0.2844) -- cycle;
\path[cbpmesh=26.6] (-1.4144,1.1251) -- (-1.4666,1.0719) -- (-1.5356,1.1541) -- (-1.4810,1.2097) -- cycle;
\path[cbpmesh=15.6] (1.4666,1.0719) -- (1.4144,1.1251) -- (1.4810,1.2097) -- (1.5356,1.1541) -- cycle;
\path[cbpmesh=22.9] (-1.5543,1.2864) -- (-1.6117,1.2280) -- (-1.6943,1.2934) -- (-1.6340,1.3548) -- cycle;
\path[cbpmesh=14.8] (1.6117,1.2280) -- (1.5543,1.2864) -- (1.6340,1.3548) -- (1.6943,1.2934) -- cycle;
\path[cbpmesh=17.8] (-2.5983,1.3711) -- (-2.6565,1.2742) -- (-2.7708,1.2474) -- (-2.7101,1.3485) -- cycle;
\path[cbpmesh=23.5] (2.6565,1.2742) -- (2.5983,1.3711) -- (2.7101,1.3485) -- (2.7708,1.2474) -- cycle;
\path[cbpmesh=24.7] (-1.4810,1.2097) -- (-1.5356,1.1541) -- (-1.6117,1.2280) -- (-1.5543,1.2864) -- cycle;
\path[cbpmesh=15.0] (1.5356,1.1541) -- (1.4810,1.2097) -- (1.5543,1.2864) -- (1.6117,1.2280) -- cycle;
\path[cbpmesh=18.0] (-3.3232,1.2097) -- (-3.3755,1.0861) -- (-3.4519,1.0880) -- (-3.3983,1.2143) -- cycle;
\path[cbpmesh=18.9] (3.3755,1.0861) -- (3.3232,1.2097) -- (3.3983,1.2143) -- (3.4519,1.0880) -- cycle;
\path[cbpmesh=38] (-1.3182,0.0188) -- (-1.3305,-0.0301) -- (-1.3231,0.1063) -- (-1.3109,0.1549) -- cycle;
\path[cbpmesh=29.3] (1.3305,-0.0301) -- (1.3182,0.0188) -- (1.3109,0.1549) -- (1.3231,0.1063) -- cycle;
\path[cbpmesh=18.2] (-2.0788,1.4620) -- (-2.1400,1.3842) -- (-2.2501,1.3975) -- (-2.1858,1.4793) -- cycle;
\path[cbpmesh=18.7] (2.1400,1.3842) -- (2.0788,1.4620) -- (2.1858,1.4793) -- (2.2501,1.3975) -- cycle;
\path[cbpmesh=19.0] (-4.1975,0.9159) -- (-4.2367,0.7604) -- (-4.3117,0.6975) -- (-4.2718,0.8558) -- cycle;
\path[cbpmesh=30.8] (4.2367,0.7604) -- (4.1975,0.9159) -- (4.2718,0.8558) -- (4.3117,0.6975) -- cycle;
\path[cbpmesh=17.9] (-3.2442,1.2065) -- (-3.2953,1.0858) -- (-3.3755,1.0861) -- (-3.3232,1.2097) -- cycle;
\path[cbpmesh=19.2] (3.2953,1.0858) -- (3.2442,1.2065) -- (3.3232,1.2097) -- (3.3755,1.0861) -- cycle;
\path[cbpmesh=33.8] (-1.3072,0.7796) -- (-1.3456,0.7307) -- (-1.3916,0.8337) -- (-1.3518,0.8843) -- cycle;
\path[cbpmesh=19.4] (1.3456,0.7307) -- (1.3072,0.7796) -- (1.3518,0.8843) -- (1.3916,0.8337) -- cycle;
\path[cbpmesh=38] (-1.2827,0.5034) -- (-1.3115,0.4555) -- (-1.3385,0.5716) -- (-1.3091,0.6204) -- cycle;
\path[cbpmesh=23.0] (1.3115,0.4555) -- (1.2827,0.5034) -- (1.3091,0.6204) -- (1.3385,0.5716) -- cycle;
\path[cbpmesh=17.7] (-3.1607,1.2057) -- (-3.2105,1.0882) -- (-3.2953,1.0858) -- (-3.2442,1.2065) -- cycle;
\path[cbpmesh=19.6] (3.2105,1.0882) -- (3.1607,1.2057) -- (3.2442,1.2065) -- (3.2953,1.0858) -- cycle;
\path[cbpmesh=17.7] (-2.4798,1.3912) -- (-2.5354,1.2987) -- (-2.6565,1.2742) -- (-2.5983,1.3711) -- cycle;
\path[cbpmesh=23.0] (2.5354,1.2987) -- (2.4798,1.3912) -- (2.5983,1.3711) -- (2.6565,1.2742) -- cycle;
\path[cbpmesh=18.7] (-1.9756,1.4345) -- (-2.0338,1.3605) -- (-2.1400,1.3842) -- (-2.0788,1.4620) -- cycle;
\path[cbpmesh=17.5] (2.0338,1.3605) -- (1.9756,1.4345) -- (2.0788,1.4620) -- (2.1400,1.3842) -- cycle;
\path[cbpmesh=31.6] (-1.3518,0.8843) -- (-1.3916,0.8337) -- (-1.4465,0.9293) -- (-1.4051,0.9819) -- cycle;
\path[cbpmesh=18.0] (1.3916,0.8337) -- (1.3518,0.8843) -- (1.4051,0.9819) -- (1.4465,0.9293) -- cycle;
\path[cbpmesh=18.4] (-4.1262,0.9652) -- (-4.1647,0.8123) -- (-4.2367,0.7604) -- (-4.1975,0.9159) -- cycle;
\path[cbpmesh=29.3] (4.1647,0.8123) -- (4.1262,0.9652) -- (4.1975,0.9159) -- (4.2367,0.7604) -- cycle;
\path[cbpmesh=38] (-1.2953,0.3325) -- (-1.3157,0.2844) -- (-1.3321,0.4068) -- (-1.3115,0.4555) -- cycle;
\path[cbpmesh=25.3] (1.3157,0.2844) -- (1.2953,0.3325) -- (1.3115,0.4555) -- (1.3321,0.4068) -- cycle;
\path[cbpmesh=17.6] (-3.0722,1.2084) -- (-3.1206,1.0941) -- (-3.2105,1.0882) -- (-3.1607,1.2057) -- cycle;
\path[cbpmesh=20.1] (3.1206,1.0941) -- (3.0722,1.2084) -- (3.1607,1.2057) -- (3.2105,1.0882) -- cycle;
\path[cbpmesh=26.8] (-4.4947,0.3428) -- (-4.5086,0.1767) -- (-4.5306,0.0635) -- (-4.5166,0.2304) -- cycle;
\draw[draw=cbpInk!70,line width=.45pt] (-4.5306,0.0635) -- (-4.5166,0.2304);
\draw[draw=cbpInk!70,line width=.45pt] (-4.5166,0.2304) -- (-4.4947,0.3428);
\path[cbpmesh=38] (4.5086,0.1767) -- (4.4947,0.3428) -- (4.5166,0.2304) -- (4.5306,0.0635) -- cycle;
\draw[draw=cbpInk!70,line width=.45pt] (4.4947,0.3428) -- (4.5166,0.2304);
\draw[draw=cbpInk!70,line width=.45pt] (4.5166,0.2304) -- (4.5306,0.0635);
\path[cbpmesh=19.4] (-1.8768,1.3969) -- (-1.9320,1.3267) -- (-2.0338,1.3605) -- (-1.9756,1.4345) -- cycle;
\path[cbpmesh=16.4] (1.9320,1.3267) -- (1.8768,1.3969) -- (1.9756,1.4345) -- (2.0338,1.3605) -- cycle;
\path[cbpmesh=36.8] (-1.3091,0.6204) -- (-1.3385,0.5716) -- (-1.3758,0.6805) -- (-1.3456,0.7307) -- cycle;
\path[cbpmesh=21.4] (1.3385,0.5716) -- (1.3091,0.6204) -- (1.3456,0.7307) -- (1.3758,0.6805) -- cycle;
\path[cbpmesh=29.5] (-1.4051,0.9819) -- (-1.4465,0.9293) -- (-1.5097,1.0170) -- (-1.4666,1.0719) -- cycle;
\path[cbpmesh=16.9] (1.4465,0.9293) -- (1.4051,0.9819) -- (1.4666,1.0719) -- (1.5097,1.0170) -- cycle;
\path[cbpmesh=38] (-1.3109,0.1549) -- (-1.3231,0.1063) -- (-1.3280,0.2356) -- (-1.3157,0.2844) -- cycle;
\path[cbpmesh=27.7] (1.3231,0.1063) -- (1.3109,0.1549) -- (1.3157,0.2844) -- (1.3280,0.2356) -- cycle;
\path[cbpmesh=17.5] (-2.9781,1.2156) -- (-3.0250,1.1048) -- (-3.1206,1.0941) -- (-3.0722,1.2084) -- cycle;
\path[cbpmesh=20.8] (3.0250,1.1048) -- (2.9781,1.2156) -- (3.0722,1.2084) -- (3.1206,1.0941) -- cycle;
\path[cbpmesh=17.9] (-4.0573,1.0046) -- (-4.0952,0.8543) -- (-4.1647,0.8123) -- (-4.1262,0.9652) -- cycle;
\path[cbpmesh=27.7] (4.0952,0.8543) -- (4.0573,1.0046) -- (4.1262,0.9652) -- (4.1647,0.8123) -- cycle;
\path[cbpmesh=17.6] (-2.3636,1.3999) -- (-2.4165,1.3117) -- (-2.5354,1.2987) -- (-2.4798,1.3912) -- cycle;
\path[cbpmesh=21.5] (2.4165,1.3117) -- (2.3636,1.3999) -- (2.4798,1.3912) -- (2.5354,1.2987) -- cycle;
\path[cbpmesh=20.5] (-1.7828,1.3498) -- (-1.8353,1.2831) -- (-1.9320,1.3267) -- (-1.8768,1.3969) -- cycle;
\path[cbpmesh=15.7] (1.8353,1.2831) -- (1.7828,1.3498) -- (1.8768,1.3969) -- (1.9320,1.3267) -- cycle;
\path[cbpmesh=27.4] (-1.4666,1.0719) -- (-1.5097,1.0170) -- (-1.5808,1.0966) -- (-1.5356,1.1541) -- cycle;
\path[cbpmesh=16.0] (1.5097,1.0170) -- (1.4666,1.0719) -- (1.5356,1.1541) -- (1.5808,1.0966) -- cycle;
\path[cbpmesh=38] (-1.3305,-0.0301) -- (-1.3347,-0.0793) -- (-1.3272,0.0575) -- (-1.3231,0.1063) -- cycle;
\draw[draw=cbpInk!70,line width=.45pt] (-1.3305,-0.0301) -- (-1.3347,-0.0793);
\draw[draw=cbpInk!70,line width=.45pt] (-1.3347,-0.0793) -- (-1.3272,0.0575);
\path[cbpmesh=30.1] (1.3347,-0.0793) -- (1.3305,-0.0301) -- (1.3231,0.1063) -- (1.3272,0.0575) -- cycle;
\draw[draw=cbpInk!70,line width=.45pt] (1.3347,-0.0793) -- (1.3305,-0.0301);
\draw[draw=cbpInk!70,line width=.45pt] (1.3272,0.0575) -- (1.3347,-0.0793);
\path[cbpmesh=21.9] (-1.6943,1.2934) -- (-1.7441,1.2299) -- (-1.8353,1.2831) -- (-1.7828,1.3498) -- cycle;
\path[cbpmesh=15.2] (1.7441,1.2299) -- (1.6943,1.2934) -- (1.7828,1.3498) -- (1.8353,1.2831) -- cycle;
\path[cbpmesh=25.3] (-1.5356,1.1541) -- (-1.5808,1.0966) -- (-1.6591,1.1677) -- (-1.6117,1.2280) -- cycle;
\path[cbpmesh=15.4] (1.5808,1.0966) -- (1.5356,1.1541) -- (1.6117,1.2280) -- (1.6591,1.1677) -- cycle;
\path[cbpmesh=23.5] (-1.6117,1.2280) -- (-1.6591,1.1677) -- (-1.7441,1.2299) -- (-1.6943,1.2934) -- cycle;
\path[cbpmesh=15.1] (1.6591,1.1677) -- (1.6117,1.2280) -- (1.6943,1.2934) -- (1.7441,1.2299) -- cycle;
\path[cbpmesh=17.4] (-2.8778,1.2282) -- (-2.9231,1.1212) -- (-3.0250,1.1048) -- (-2.9781,1.2156) -- cycle;
\path[cbpmesh=21.6] (2.9231,1.1212) -- (2.8778,1.2282) -- (2.9781,1.2156) -- (3.0250,1.1048) -- cycle;
\path[cbpmesh=38] (-1.3115,0.4555) -- (-1.3321,0.4068) -- (-1.3596,0.5218) -- (-1.3385,0.5716) -- cycle;
\path[cbpmesh=23.7] (1.3321,0.4068) -- (1.3115,0.4555) -- (1.3385,0.5716) -- (1.3596,0.5218) -- cycle;
\path[cbpmesh=34.7] (-1.3456,0.7307) -- (-1.3758,0.6805) -- (-1.4228,0.7818) -- (-1.3916,0.8337) -- cycle;
\path[cbpmesh=19.9] (1.3758,0.6805) -- (1.3456,0.7307) -- (1.3916,0.8337) -- (1.4228,0.7818) -- cycle;
\path[cbpmesh=17.4] (-3.9904,1.0353) -- (-4.0276,0.8874) -- (-4.0952,0.8543) -- (-4.0573,1.0046) -- cycle;
\path[cbpmesh=26.0] (4.0276,0.8874) -- (3.9904,1.0353) -- (4.0573,1.0046) -- (4.0952,0.8543) -- cycle;
\path[cbpmesh=24.6] (-4.4541,0.4585) -- (-4.4679,0.2939) -- (-4.5086,0.1767) -- (-4.4947,0.3428) -- cycle;
\path[cbpmesh=38] (4.4679,0.2939) -- (4.4541,0.4585) -- (4.4947,0.3428) -- (4.5086,0.1767) -- cycle;
\path[cbpmesh=17.8] (-2.2501,1.3975) -- (-2.3005,1.3135) -- (-2.4165,1.3117) -- (-2.3636,1.3999) -- cycle;
\path[cbpmesh=20.0] (2.3005,1.3135) -- (2.2501,1.3975) -- (2.3636,1.3999) -- (2.4165,1.3117) -- cycle;
\path[cbpmesh=17.2] (-3.9247,1.0581) -- (-3.9613,0.9127) -- (-4.0276,0.8874) -- (-3.9904,1.0353) -- cycle;
\path[cbpmesh=24.3] (3.9613,0.9127) -- (3.9247,1.0581) -- (3.9904,1.0353) -- (4.0276,0.8874) -- cycle;
\path[cbpmesh=17.4] (-2.7708,1.2474) -- (-2.8145,1.1444) -- (-2.9231,1.1212) -- (-2.8778,1.2282) -- cycle;
\path[cbpmesh=22.4] (2.8145,1.1444) -- (2.7708,1.2474) -- (2.8778,1.2282) -- (2.9231,1.1212) -- cycle;
\path[cbpmesh=38] (-1.3157,0.2844) -- (-1.3280,0.2356) -- (-1.3446,0.3574) -- (-1.3321,0.4068) -- cycle;
\path[cbpmesh=26.1] (1.3280,0.2356) -- (1.3157,0.2844) -- (1.3321,0.4068) -- (1.3446,0.3574) -- cycle;
\path[cbpmesh=32.5] (-1.3916,0.8337) -- (-1.4228,0.7818) -- (-1.4789,0.8753) -- (-1.4465,0.9293) -- cycle;
\path[cbpmesh=18.5] (1.4228,0.7818) -- (1.3916,0.8337) -- (1.4465,0.9293) -- (1.4789,0.8753) -- cycle;
\path[cbpmesh=17.1] (-3.8597,1.0743) -- (-3.8957,0.9312) -- (-3.9613,0.9127) -- (-3.9247,1.0581) -- cycle;
\path[cbpmesh=22.7] (3.8957,0.9312) -- (3.8597,1.0743) -- (3.9247,1.0581) -- (3.9613,0.9127) -- cycle;
\path[cbpmesh=18.1] (-2.1400,1.3842) -- (-2.1879,1.3044) -- (-2.3005,1.3135) -- (-2.2501,1.3975) -- cycle;
\path[cbpmesh=18.7] (2.1879,1.3044) -- (2.1400,1.3842) -- (2.2501,1.3975) -- (2.3005,1.3135) -- cycle;
\path[cbpmesh=37.8] (-1.3385,0.5716) -- (-1.3596,0.5218) -- (-1.3974,0.6293) -- (-1.3758,0.6805) -- cycle;
\path[cbpmesh=22.1] (1.3596,0.5218) -- (1.3385,0.5716) -- (1.3758,0.6805) -- (1.3974,0.6293) -- cycle;
\path[cbpmesh=38] (-1.3231,0.1063) -- (-1.3272,0.0575) -- (-1.3321,0.1865) -- (-1.3280,0.2356) -- cycle;
\draw[draw=cbpInk!70,line width=.45pt] (-1.3272,0.0575) -- (-1.3321,0.1865);
\path[cbpmesh=28.5] (1.3272,0.0575) -- (1.3231,0.1063) -- (1.3280,0.2356) -- (1.3321,0.1865) -- cycle;
\draw[draw=cbpInk!70,line width=.45pt] (1.3321,0.1865) -- (1.3272,0.0575);
\path[cbpmesh=17.1] (-3.7948,1.0847) -- (-3.8302,0.9441) -- (-3.8957,0.9312) -- (-3.8597,1.0743) -- cycle;
\path[cbpmesh=21.4] (3.8302,0.9441) -- (3.7948,1.0847) -- (3.8597,1.0743) -- (3.8957,0.9312) -- cycle;
\path[cbpmesh=30.3] (-1.4465,0.9293) -- (-1.4789,0.8753) -- (-1.5435,0.9607) -- (-1.5097,1.0170) -- cycle;
\path[cbpmesh=17.3] (1.4789,0.8753) -- (1.4465,0.9293) -- (1.5097,1.0170) -- (1.5435,0.9607) -- cycle;
\path[cbpmesh=17.4] (-2.6565,1.2742) -- (-2.6984,1.1754) -- (-2.8145,1.1444) -- (-2.7708,1.2474) -- cycle;
\path[cbpmesh=23.3] (2.6984,1.1754) -- (2.6565,1.2742) -- (2.7708,1.2474) -- (2.8145,1.1444) -- cycle;
\path[cbpmesh=17.3] (-3.7294,1.0905) -- (-3.7643,0.9523) -- (-3.8302,0.9441) -- (-3.7948,1.0847) -- cycle;
\path[cbpmesh=20.3] (3.7643,0.9523) -- (3.7294,1.0905) -- (3.7948,1.0847) -- (3.8302,0.9441) -- cycle;
\path[cbpmesh=22.7] (-4.3935,0.5769) -- (-4.4071,0.4146) -- (-4.4679,0.2939) -- (-4.4541,0.4585) -- cycle;
\path[cbpmesh=38] (4.4071,0.4146) -- (4.3935,0.5769) -- (4.4541,0.4585) -- (4.4679,0.2939) -- cycle;
\path[cbpmesh=18.7] (-2.0338,1.3605) -- (-2.0793,1.2847) -- (-2.1879,1.3044) -- (-2.1400,1.3842) -- cycle;
\path[cbpmesh=17.5] (2.0793,1.2847) -- (2.0338,1.3605) -- (2.1400,1.3842) -- (2.1879,1.3044) -- cycle;
\path[cbpmesh=17.6] (-3.6631,1.0927) -- (-3.6973,0.9569) -- (-3.7643,0.9523) -- (-3.7294,1.0905) -- cycle;
\path[cbpmesh=19.6] (3.6973,0.9569) -- (3.6631,1.0927) -- (3.7294,1.0905) -- (3.7643,0.9523) -- cycle;
\path[cbpmesh=28.1] (-1.5097,1.0170) -- (-1.5435,0.9607) -- (-1.6162,1.0376) -- (-1.5808,1.0966) -- cycle;
\path[cbpmesh=16.4] (1.5435,0.9607) -- (1.5097,1.0170) -- (1.5808,1.0966) -- (1.6162,1.0376) -- cycle;
\path[cbpmesh=38] (-1.3321,0.4068) -- (-1.3446,0.3574) -- (-1.3723,0.4715) -- (-1.3596,0.5218) -- cycle;
\path[cbpmesh=24.4] (1.3446,0.3574) -- (1.3321,0.4068) -- (1.3596,0.5218) -- (1.3723,0.4715) -- cycle;
\path[cbpmesh=19.6] (-1.9320,1.3267) -- (-1.9753,1.2546) -- (-2.0793,1.2847) -- (-2.0338,1.3605) -- cycle;
\path[cbpmesh=16.6] (1.9753,1.2546) -- (1.9320,1.3267) -- (2.0338,1.3605) -- (2.0793,1.2847) -- cycle;
\path[cbpmesh=17.8] (-3.5951,1.0923) -- (-3.6287,0.9591) -- (-3.6973,0.9569) -- (-3.6631,1.0927) -- cycle;
\path[cbpmesh=19.1] (3.6287,0.9591) -- (3.5951,1.0923) -- (3.6631,1.0927) -- (3.6973,0.9569) -- cycle;
\path[cbpmesh=35.6] (-1.3758,0.6805) -- (-1.3974,0.6293) -- (-1.4452,0.7289) -- (-1.4228,0.7818) -- cycle;
\path[cbpmesh=20.5] (1.3974,0.6293) -- (1.3758,0.6805) -- (1.4228,0.7818) -- (1.4452,0.7289) -- cycle;
\path[cbpmesh=26.0] (-1.5808,1.0966) -- (-1.6162,1.0376) -- (-1.6963,1.1058) -- (-1.6591,1.1677) -- cycle;
\path[cbpmesh=15.8] (1.6162,1.0376) -- (1.5808,1.0966) -- (1.6591,1.1677) -- (1.6963,1.1058) -- cycle;
\path[cbpmesh=17.9] (-3.5249,1.0904) -- (-3.5578,0.9598) -- (-3.6287,0.9591) -- (-3.5951,1.0923) -- cycle;
\path[cbpmesh=18.8] (3.5578,0.9598) -- (3.5249,1.0904) -- (3.5951,1.0923) -- (3.6287,0.9591) -- cycle;
\path[cbpmesh=20.8] (-1.8353,1.2831) -- (-1.8764,1.2146) -- (-1.9753,1.2546) -- (-1.9320,1.3267) -- cycle;
\path[cbpmesh=15.9] (1.8764,1.2146) -- (1.8353,1.2831) -- (1.9320,1.3267) -- (1.9753,1.2546) -- cycle;
\path[cbpmesh=17.4] (-2.5354,1.2987) -- (-2.5753,1.2045) -- (-2.6984,1.1754) -- (-2.6565,1.2742) -- cycle;
\path[cbpmesh=22.8] (2.5753,1.2045) -- (2.5354,1.2987) -- (2.6565,1.2742) -- (2.6984,1.1754) -- cycle;
\path[cbpmesh=24.0] (-1.6591,1.1677) -- (-1.6963,1.1058) -- (-1.7832,1.1649) -- (-1.7441,1.2299) -- cycle;
\path[cbpmesh=15.4] (1.6963,1.1058) -- (1.6591,1.1677) -- (1.7441,1.2299) -- (1.7832,1.1649) -- cycle;
\path[cbpmesh=22.3] (-1.7441,1.2299) -- (-1.7832,1.1649) -- (-1.8764,1.2146) -- (-1.8353,1.2831) -- cycle;
\path[cbpmesh=15.5] (1.7832,1.1649) -- (1.7441,1.2299) -- (1.8353,1.2831) -- (1.8764,1.2146) -- cycle;
\path[cbpmesh=18.0] (-3.4519,1.0880) -- (-3.4841,0.9601) -- (-3.5578,0.9598) -- (-3.5249,1.0904) -- cycle;
\path[cbpmesh=18.8] (3.4841,0.9601) -- (3.4519,1.0880) -- (3.5249,1.0904) -- (3.5578,0.9598) -- cycle;
\path[cbpmesh=17.9] (-3.3755,1.0861) -- (-3.4071,0.9610) -- (-3.4841,0.9601) -- (-3.4519,1.0880) -- cycle;
\path[cbpmesh=18.9] (3.4071,0.9610) -- (3.3755,1.0861) -- (3.4519,1.0880) -- (3.4841,0.9601) -- cycle;
\path[cbpmesh=38] (-1.3280,0.2356) -- (-1.3321,0.1865) -- (-1.3487,0.3077) -- (-1.3446,0.3574) -- cycle;
\path[cbpmesh=26.9] (1.3321,0.1865) -- (1.3280,0.2356) -- (1.3446,0.3574) -- (1.3487,0.3077) -- cycle;
\path[cbpmesh=17.8] (-3.2953,1.0858) -- (-3.3261,0.9637) -- (-3.4071,0.9610) -- (-3.3755,1.0861) -- cycle;
\path[cbpmesh=19.1] (3.3261,0.9637) -- (3.2953,1.0858) -- (3.3755,1.0861) -- (3.4071,0.9610) -- cycle;
\path[cbpmesh=33.4] (-1.4228,0.7818) -- (-1.4452,0.7289) -- (-1.5022,0.8204) -- (-1.4789,0.8753) -- cycle;
\path[cbpmesh=19.1] (1.4452,0.7289) -- (1.4228,0.7818) -- (1.4789,0.8753) -- (1.5022,0.8204) -- cycle;
\path[cbpmesh=21.2] (-4.3117,0.6975) -- (-4.3250,0.5382) -- (-4.4071,0.4146) -- (-4.3935,0.5769) -- cycle;
\path[cbpmesh=35.7] (4.3250,0.5382) -- (4.3117,0.6975) -- (4.3935,0.5769) -- (4.4071,0.4146) -- cycle;
\path[cbpmesh=38] (-1.3596,0.5218) -- (-1.3723,0.4715) -- (-1.4105,0.5775) -- (-1.3974,0.6293) -- cycle;
\path[cbpmesh=22.8] (1.3723,0.4715) -- (1.3596,0.5218) -- (1.3974,0.6293) -- (1.4105,0.5775) -- cycle;
\path[cbpmesh=17.4] (-2.4165,1.3117) -- (-2.4546,1.2219) -- (-2.5753,1.2045) -- (-2.5354,1.2987) -- cycle;
\path[cbpmesh=21.3] (2.4546,1.2219) -- (2.4165,1.3117) -- (2.5354,1.2987) -- (2.5753,1.2045) -- cycle;
\path[cbpmesh=17.6] (-3.2105,1.0882) -- (-3.2405,0.9692) -- (-3.3261,0.9637) -- (-3.2953,1.0858) -- cycle;
\path[cbpmesh=19.5] (3.2405,0.9692) -- (3.2105,1.0882) -- (3.2953,1.0858) -- (3.3261,0.9637) -- cycle;
\path[cbpmesh=17.4] (-3.1206,1.0941) -- (-3.1497,0.9785) -- (-3.2405,0.9692) -- (-3.2105,1.0882) -- cycle;
\path[cbpmesh=20.0] (3.1497,0.9785) -- (3.1206,1.0941) -- (3.2105,1.0882) -- (3.2405,0.9692) -- cycle;
\path[cbpmesh=31.1] (-1.4789,0.8753) -- (-1.5022,0.8204) -- (-1.5679,0.9033) -- (-1.5435,0.9607) -- cycle;
\path[cbpmesh=17.9] (1.5022,0.8204) -- (1.4789,0.8753) -- (1.5435,0.9607) -- (1.5679,0.9033) -- cycle;
\path[cbpmesh=17.6] (-2.3005,1.3135) -- (-2.3368,1.2280) -- (-2.4546,1.2219) -- (-2.4165,1.3117) -- cycle;
\path[cbpmesh=19.9] (2.3368,1.2280) -- (2.3005,1.3135) -- (2.4165,1.3117) -- (2.4546,1.2219) -- cycle;
\path[cbpmesh=17.3] (-3.0250,1.1048) -- (-3.0533,0.9927) -- (-3.1497,0.9785) -- (-3.1206,1.0941) -- cycle;
\path[cbpmesh=20.7] (3.0533,0.9927) -- (3.0250,1.1048) -- (3.1206,1.0941) -- (3.1497,0.9785) -- cycle;
\path[cbpmesh=38] (-1.3446,0.3574) -- (-1.3487,0.3077) -- (-1.3765,0.4208) -- (-1.3723,0.4715) -- cycle;
\path[cbpmesh=25.2] (1.3487,0.3077) -- (1.3446,0.3574) -- (1.3723,0.4715) -- (1.3765,0.4208) -- cycle;
\path[cbpmesh=18.4] (-4.2367,0.7604) -- (-4.2498,0.6038) -- (-4.3250,0.5382) -- (-4.3117,0.6975) -- cycle;
\path[cbpmesh=30.2] (4.2498,0.6038) -- (4.2367,0.7604) -- (4.3117,0.6975) -- (4.3250,0.5382) -- cycle;
\path[cbpmesh=36.5] (-1.3974,0.6293) -- (-1.4105,0.5775) -- (-1.4587,0.6754) -- (-1.4452,0.7289) -- cycle;
\path[cbpmesh=21.2] (1.4105,0.5775) -- (1.3974,0.6293) -- (1.4452,0.7289) -- (1.4587,0.6754) -- cycle;
\path[cbpmesh=28.8] (-1.5435,0.9607) -- (-1.5679,0.9033) -- (-1.6417,0.9775) -- (-1.6162,1.0376) -- cycle;
\path[cbpmesh=16.9] (1.5679,0.9033) -- (1.5435,0.9607) -- (1.6162,1.0376) -- (1.6417,0.9775) -- cycle;
\path[cbpmesh=28.5] (-4.5306,0.0635) -- (-4.5166,-0.1034) -- (-4.5210,-0.2123) -- (-4.5350,-0.0452) -- cycle;
\draw[draw=cbpInk!70,line width=.45pt] (-4.5210,-0.2123) -- (-4.5350,-0.0452);
\draw[draw=cbpInk!70,line width=.45pt] (-4.5350,-0.0452) -- (-4.5306,0.0635);
\path[cbpmesh=38] (4.5166,-0.1034) -- (4.5306,0.0635) -- (4.5350,-0.0452) -- (4.5210,-0.2123) -- cycle;
\draw[draw=cbpInk!70,line width=.45pt] (4.5306,0.0635) -- (4.5350,-0.0452);
\draw[draw=cbpInk!70,line width=.45pt] (4.5350,-0.0452) -- (4.5210,-0.2123);
\path[cbpmesh=17.1] (-2.9231,1.1212) -- (-2.9505,1.0129) -- (-3.0533,0.9927) -- (-3.0250,1.1048) -- cycle;
\path[cbpmesh=21.4] (2.9505,1.0129) -- (2.9231,1.1212) -- (3.0250,1.1048) -- (3.0533,0.9927) -- cycle;
\path[cbpmesh=18.1] (-2.1879,1.3044) -- (-2.2224,1.2230) -- (-2.3368,1.2280) -- (-2.3005,1.3135) -- cycle;
\path[cbpmesh=18.7] (2.2224,1.2230) -- (2.1879,1.3044) -- (2.3005,1.3135) -- (2.3368,1.2280) -- cycle;
\path[cbpmesh=26.6] (-1.6162,1.0376) -- (-1.6417,0.9775) -- (-1.7230,1.0427) -- (-1.6963,1.1058) -- cycle;
\path[cbpmesh=16.2] (1.6417,0.9775) -- (1.6162,1.0376) -- (1.6963,1.1058) -- (1.7230,1.0427) -- cycle;
\path[cbpmesh=38] (-1.3321,0.1865) -- (-1.3280,0.1375) -- (-1.3446,0.2580) -- (-1.3487,0.3077) -- cycle;
\draw[draw=cbpInk!70,line width=.45pt] (-1.3321,0.1865) -- (-1.3280,0.1375);
\draw[draw=cbpInk!70,line width=.45pt] (-1.3280,0.1375) -- (-1.3446,0.2580);
\path[cbpmesh=27.7] (1.3280,0.1375) -- (1.3321,0.1865) -- (1.3487,0.3077) -- (1.3446,0.2580) -- cycle;
\draw[draw=cbpInk!70,line width=.45pt] (1.3280,0.1375) -- (1.3321,0.1865);
\draw[draw=cbpInk!70,line width=.45pt] (1.3446,0.2580) -- (1.3280,0.1375);
\path[cbpmesh=18.8] (-2.0793,1.2847) -- (-2.1121,1.2073) -- (-2.2224,1.2230) -- (-2.1879,1.3044) -- cycle;
\path[cbpmesh=17.6] (2.1121,1.2073) -- (2.0793,1.2847) -- (2.1879,1.3044) -- (2.2224,1.2230) -- cycle;
\path[cbpmesh=17.1] (-2.8145,1.1444) -- (-2.8408,1.0401) -- (-2.9505,1.0129) -- (-2.9231,1.1212) -- cycle;
\path[cbpmesh=22.2] (2.8408,1.0401) -- (2.8145,1.1444) -- (2.9231,1.1212) -- (2.9505,1.0129) -- cycle;
\path[cbpmesh=17.8] (-4.1647,0.8123) -- (-4.1776,0.6584) -- (-4.2498,0.6038) -- (-4.2367,0.7604) -- cycle;
\path[cbpmesh=28.8] (4.1776,0.6584) -- (4.1647,0.8123) -- (4.2367,0.7604) -- (4.2498,0.6038) -- cycle;
\path[cbpmesh=24.5] (-1.6963,1.1058) -- (-1.7230,1.0427) -- (-1.8113,1.0986) -- (-1.7832,1.1649) -- cycle;
\path[cbpmesh=15.8] (1.7230,1.0427) -- (1.6963,1.1058) -- (1.7832,1.1649) -- (1.8113,1.0986) -- cycle;
\path[cbpmesh=19.8] (-1.9753,1.2546) -- (-2.0064,1.1812) -- (-2.1121,1.2073) -- (-2.0793,1.2847) -- cycle;
\path[cbpmesh=16.7] (2.0064,1.1812) -- (1.9753,1.2546) -- (2.0793,1.2847) -- (2.1121,1.2073) -- cycle;
\path[cbpmesh=34.2] (-1.4452,0.7289) -- (-1.4587,0.6754) -- (-1.5162,0.7647) -- (-1.5022,0.8204) -- cycle;
\path[cbpmesh=19.8] (1.4587,0.6754) -- (1.4452,0.7289) -- (1.5022,0.8204) -- (1.5162,0.7647) -- cycle;
\path[cbpmesh=22.7] (-1.7832,1.1649) -- (-1.8113,1.0986) -- (-1.9060,1.1448) -- (-1.8764,1.2146) -- cycle;
\path[cbpmesh=15.8] (1.8113,1.0986) -- (1.7832,1.1649) -- (1.8764,1.2146) -- (1.9060,1.1448) -- cycle;
\path[cbpmesh=21.1] (-1.8764,1.2146) -- (-1.9060,1.1448) -- (-2.0064,1.1812) -- (-1.9753,1.2546) -- cycle;
\path[cbpmesh=16.1] (1.9060,1.1448) -- (1.8764,1.2146) -- (1.9753,1.2546) -- (2.0064,1.1812) -- cycle;
\path[cbpmesh=38] (-1.3723,0.4715) -- (-1.3765,0.4208) -- (-1.4149,0.5254) -- (-1.4105,0.5775) -- cycle;
\path[cbpmesh=23.5] (1.3765,0.4208) -- (1.3723,0.4715) -- (1.4105,0.5775) -- (1.4149,0.5254) -- cycle;
\path[cbpmesh=17.1] (-2.6984,1.1754) -- (-2.7236,1.0754) -- (-2.8408,1.0401) -- (-2.8145,1.1444) -- cycle;
\path[cbpmesh=23.0] (2.7236,1.0754) -- (2.6984,1.1754) -- (2.8145,1.1444) -- (2.8408,1.0401) -- cycle;
\path[cbpmesh=17.4] (-4.0952,0.8543) -- (-4.1079,0.7029) -- (-4.1776,0.6584) -- (-4.1647,0.8123) -- cycle;
\path[cbpmesh=27.2] (4.1079,0.7029) -- (4.0952,0.8543) -- (4.1647,0.8123) -- (4.1776,0.6584) -- cycle;
\path[cbpmesh=26.0] (-4.5086,0.1767) -- (-4.4947,0.0107) -- (-4.5166,-0.1034) -- (-4.5306,0.0635) -- cycle;
\path[cbpmesh=38] (4.4947,0.0107) -- (4.5086,0.1767) -- (4.5306,0.0635) -- (4.5166,-0.1034) -- cycle;
\path[cbpmesh=31.8] (-1.5022,0.8204) -- (-1.5162,0.7647) -- (-1.5825,0.8452) -- (-1.5679,0.9033) -- cycle;
\path[cbpmesh=18.5] (1.5162,0.7647) -- (1.5022,0.8204) -- (1.5679,0.9033) -- (1.5825,0.8452) -- cycle;
\path[cbpmesh=38] (-1.3487,0.3077) -- (-1.3446,0.2580) -- (-1.3723,0.3700) -- (-1.3765,0.4208) -- cycle;
\draw[draw=cbpInk!70,line width=.45pt] (-1.3446,0.2580) -- (-1.3723,0.3700);
\path[cbpmesh=26.0] (1.3446,0.2580) -- (1.3487,0.3077) -- (1.3765,0.4208) -- (1.3723,0.3700) -- cycle;
\draw[draw=cbpInk!70,line width=.45pt] (1.3723,0.3700) -- (1.3446,0.2580);
\path[cbpmesh=17.0] (-4.0276,0.8874) -- (-4.0401,0.7386) -- (-4.1079,0.7029) -- (-4.0952,0.8543) -- cycle;
\path[cbpmesh=25.6] (4.0401,0.7386) -- (4.0276,0.8874) -- (4.0952,0.8543) -- (4.1079,0.7029) -- cycle;
\path[cbpmesh=17.1] (-2.5753,1.2045) -- (-2.5994,1.1090) -- (-2.7236,1.0754) -- (-2.6984,1.1754) -- cycle;
\path[cbpmesh=22.5] (2.5994,1.1090) -- (2.5753,1.2045) -- (2.6984,1.1754) -- (2.7236,1.0754) -- cycle;
\path[cbpmesh=37.3] (-1.4105,0.5775) -- (-1.4149,0.5254) -- (-1.4632,0.6215) -- (-1.4587,0.6754) -- cycle;
\path[cbpmesh=21.9] (1.4149,0.5254) -- (1.4105,0.5775) -- (1.4587,0.6754) -- (1.4632,0.6215) -- cycle;
\path[cbpmesh=16.8] (-3.9613,0.9127) -- (-3.9736,0.7663) -- (-4.0401,0.7386) -- (-4.0276,0.8874) -- cycle;
\path[cbpmesh=24.0] (3.9736,0.7663) -- (3.9613,0.9127) -- (4.0276,0.8874) -- (4.0401,0.7386) -- cycle;
\path[cbpmesh=29.5] (-1.5679,0.9033) -- (-1.5825,0.8452) -- (-1.6570,0.9167) -- (-1.6417,0.9775) -- cycle;
\path[cbpmesh=17.4] (1.5825,0.8452) -- (1.5679,0.9033) -- (1.6417,0.9775) -- (1.6570,0.9167) -- cycle;
\path[cbpmesh=16.8] (-3.8957,0.9312) -- (-3.9078,0.7873) -- (-3.9736,0.7663) -- (-3.9613,0.9127) -- cycle;
\path[cbpmesh=22.5] (3.9078,0.7873) -- (3.8957,0.9312) -- (3.9613,0.9127) -- (3.9736,0.7663) -- cycle;
\path[cbpmesh=23.8] (-4.4679,0.2939) -- (-4.4541,0.1294) -- (-4.4947,0.0107) -- (-4.5086,0.1767) -- cycle;
\path[cbpmesh=38] (4.4541,0.1294) -- (4.4679,0.2939) -- (4.5086,0.1767) -- (4.4947,0.0107) -- cycle;
\path[cbpmesh=17.2] (-2.4546,1.2219) -- (-2.4775,1.1309) -- (-2.5994,1.1090) -- (-2.5753,1.2045) -- cycle;
\path[cbpmesh=21.1] (2.4775,1.1309) -- (2.4546,1.2219) -- (2.5753,1.2045) -- (2.5994,1.1090) -- cycle;
\path[cbpmesh=16.9] (-3.8302,0.9441) -- (-3.8421,0.8025) -- (-3.9078,0.7873) -- (-3.8957,0.9312) -- cycle;
\path[cbpmesh=21.2] (3.8421,0.8025) -- (3.8302,0.9441) -- (3.8957,0.9312) -- (3.9078,0.7873) -- cycle;
\path[cbpmesh=27.2] (-1.6417,0.9775) -- (-1.6570,0.9167) -- (-1.7391,0.9789) -- (-1.7230,1.0427) -- cycle;
\path[cbpmesh=16.6] (1.6570,0.9167) -- (1.6417,0.9775) -- (1.7230,1.0427) -- (1.7391,0.9789) -- cycle;
\path[cbpmesh=17.2] (-3.7643,0.9523) -- (-3.7759,0.8132) -- (-3.8421,0.8025) -- (-3.8302,0.9441) -- cycle;
\path[cbpmesh=20.2] (3.7759,0.8132) -- (3.7643,0.9523) -- (3.8302,0.9441) -- (3.8421,0.8025) -- cycle;
\path[cbpmesh=38] (-1.3765,0.4208) -- (-1.3723,0.3700) -- (-1.4105,0.4733) -- (-1.4149,0.5254) -- cycle;
\path[cbpmesh=24.3] (1.3723,0.3700) -- (1.3765,0.4208) -- (1.4149,0.5254) -- (1.4105,0.4733) -- cycle;
\path[cbpmesh=35.0] (-1.4587,0.6754) -- (-1.4632,0.6215) -- (-1.5209,0.7087) -- (-1.5162,0.7647) -- cycle;
\path[cbpmesh=20.4] (1.4632,0.6215) -- (1.4587,0.6754) -- (1.5162,0.7647) -- (1.5209,0.7087) -- cycle;
\path[cbpmesh=17.5] (-3.6973,0.9569) -- (-3.7087,0.8203) -- (-3.7759,0.8132) -- (-3.7643,0.9523) -- cycle;
\path[cbpmesh=19.5] (3.7087,0.8203) -- (3.6973,0.9569) -- (3.7643,0.9523) -- (3.7759,0.8132) -- cycle;
\path[cbpmesh=17.5] (-2.3368,1.2280) -- (-2.3586,1.1414) -- (-2.4775,1.1309) -- (-2.4546,1.2219) -- cycle;
\path[cbpmesh=19.8] (2.3586,1.1414) -- (2.3368,1.2280) -- (2.4546,1.2219) -- (2.4775,1.1309) -- cycle;
\path[cbpmesh=17.7] (-3.6287,0.9591) -- (-3.6399,0.8250) -- (-3.7087,0.8203) -- (-3.6973,0.9569) -- cycle;
\path[cbpmesh=19.0] (3.6399,0.8250) -- (3.6287,0.9591) -- (3.6973,0.9569) -- (3.7087,0.8203) -- cycle;
\path[cbpmesh=25.0] (-1.7230,1.0427) -- (-1.7391,0.9789) -- (-1.8282,1.0314) -- (-1.8113,1.0986) -- cycle;
\path[cbpmesh=16.2] (1.7391,0.9789) -- (1.7230,1.0427) -- (1.8113,1.0986) -- (1.8282,1.0314) -- cycle;
\path[cbpmesh=17.9] (-3.5578,0.9598) -- (-3.5688,0.8283) -- (-3.6399,0.8250) -- (-3.6287,0.9591) -- cycle;
\path[cbpmesh=18.8] (3.5688,0.8283) -- (3.5578,0.9598) -- (3.6287,0.9591) -- (3.6399,0.8250) -- cycle;
\path[cbpmesh=17.9] (-3.4841,0.9601) -- (-3.4949,0.8313) -- (-3.5688,0.8283) -- (-3.5578,0.9598) -- cycle;
\path[cbpmesh=18.7] (3.4949,0.8313) -- (3.4841,0.9601) -- (3.5578,0.9598) -- (3.5688,0.8283) -- cycle;
\path[cbpmesh=18.1] (-2.2224,1.2230) -- (-2.2432,1.1407) -- (-2.3586,1.1414) -- (-2.3368,1.2280) -- cycle;
\path[cbpmesh=18.6] (2.2432,1.1407) -- (2.2224,1.2230) -- (2.3368,1.2280) -- (2.3586,1.1414) -- cycle;
\path[cbpmesh=17.9] (-3.4071,0.9610) -- (-3.4176,0.8351) -- (-3.4949,0.8313) -- (-3.4841,0.9601) -- cycle;
\path[cbpmesh=18.8] (3.4176,0.8351) -- (3.4071,0.9610) -- (3.4841,0.9601) -- (3.4949,0.8313) -- cycle;
\path[cbpmesh=23.1] (-1.8113,1.0986) -- (-1.8282,1.0314) -- (-1.9238,1.0742) -- (-1.9060,1.1448) -- cycle;
\path[cbpmesh=16.1] (1.8282,1.0314) -- (1.8113,1.0986) -- (1.9060,1.1448) -- (1.9238,1.0742) -- cycle;
\path[cbpmesh=17.7] (-3.3261,0.9637) -- (-3.3363,0.8408) -- (-3.4176,0.8351) -- (-3.4071,0.9610) -- cycle;
\path[cbpmesh=19.1] (3.3363,0.8408) -- (3.3261,0.9637) -- (3.4071,0.9610) -- (3.4176,0.8351) -- cycle;
\path[cbpmesh=18.9] (-2.1121,1.2073) -- (-2.1318,1.1291) -- (-2.2432,1.1407) -- (-2.2224,1.2230) -- cycle;
\path[cbpmesh=17.6] (2.1318,1.1291) -- (2.1121,1.2073) -- (2.2224,1.2230) -- (2.2432,1.1407) -- cycle;
\path[cbpmesh=17.5] (-3.2405,0.9692) -- (-3.2505,0.8494) -- (-3.3363,0.8408) -- (-3.3261,0.9637) -- cycle;
\path[cbpmesh=19.4] (3.2505,0.8494) -- (3.2405,0.9692) -- (3.3261,0.9637) -- (3.3363,0.8408) -- cycle;
\path[cbpmesh=21.4] (-1.9060,1.1448) -- (-1.9238,1.0742) -- (-2.0252,1.1068) -- (-2.0064,1.1812) -- cycle;
\path[cbpmesh=16.3] (1.9238,1.0742) -- (1.9060,1.1448) -- (2.0064,1.1812) -- (2.0252,1.1068) -- cycle;
\path[cbpmesh=20.0] (-2.0064,1.1812) -- (-2.0252,1.1068) -- (-2.1318,1.1291) -- (-2.1121,1.2073) -- cycle;
\path[cbpmesh=16.8] (2.0252,1.1068) -- (2.0064,1.1812) -- (2.1121,1.2073) -- (2.1318,1.1291) -- cycle;
\path[cbpmesh=21.9] (-4.4071,0.4146) -- (-4.3935,0.2522) -- (-4.4541,0.1294) -- (-4.4679,0.2939) -- cycle;
\path[cbpmesh=37.3] (4.3935,0.2522) -- (4.4071,0.4146) -- (4.4679,0.2939) -- (4.4541,0.1294) -- cycle;
\path[cbpmesh=32.5] (-1.5162,0.7647) -- (-1.5209,0.7087) -- (-1.5874,0.7867) -- (-1.5825,0.8452) -- cycle;
\path[cbpmesh=19.1] (1.5209,0.7087) -- (1.5162,0.7647) -- (1.5825,0.8452) -- (1.5874,0.7867) -- cycle;
\path[cbpmesh=17.3] (-3.1497,0.9785) -- (-3.1595,0.8621) -- (-3.2505,0.8494) -- (-3.2405,0.9692) -- cycle;
\path[cbpmesh=19.9] (3.1595,0.8621) -- (3.1497,0.9785) -- (3.2405,0.9692) -- (3.2505,0.8494) -- cycle;
\path[cbpmesh=38] (-1.4149,0.5254) -- (-1.4105,0.4733) -- (-1.4587,0.5676) -- (-1.4632,0.6215) -- cycle;
\path[cbpmesh=22.7] (1.4105,0.4733) -- (1.4149,0.5254) -- (1.4632,0.6215) -- (1.4587,0.5676) -- cycle;
\path[cbpmesh=17.1] (-3.0533,0.9927) -- (-3.0627,0.8799) -- (-3.1595,0.8621) -- (-3.1497,0.9785) -- cycle;
\path[cbpmesh=20.5] (3.0627,0.8799) -- (3.0533,0.9927) -- (3.1497,0.9785) -- (3.1595,0.8621) -- cycle;
\path[cbpmesh=16.9] (-2.9505,1.0129) -- (-2.9596,0.9039) -- (-3.0627,0.8799) -- (-3.0533,0.9927) -- cycle;
\path[cbpmesh=21.2] (2.9596,0.9039) -- (2.9505,1.0129) -- (3.0533,0.9927) -- (3.0627,0.8799) -- cycle;
\path[cbpmesh=30.1] (-1.5825,0.8452) -- (-1.5874,0.7867) -- (-1.6621,0.8555) -- (-1.6570,0.9167) -- cycle;
\path[cbpmesh=18.0] (1.5874,0.7867) -- (1.5825,0.8452) -- (1.6570,0.9167) -- (1.6621,0.8555) -- cycle;
\path[cbpmesh=16.8] (-2.8408,1.0401) -- (-2.8495,0.9351) -- (-2.9596,0.9039) -- (-2.9505,1.0129) -- cycle;
\path[cbpmesh=21.9] (2.8495,0.9351) -- (2.8408,1.0401) -- (2.9505,1.0129) -- (2.9596,0.9039) -- cycle;
\path[cbpmesh=38] (-1.3723,0.3700) -- (-1.3596,0.3197) -- (-1.3974,0.4215) -- (-1.4105,0.4733) -- cycle;
\draw[draw=cbpInk!70,line width=.45pt] (-1.3723,0.3700) -- (-1.3596,0.3197);
\draw[draw=cbpInk!70,line width=.45pt] (-1.3596,0.3197) -- (-1.3974,0.4215);
\path[cbpmesh=25.2] (1.3596,0.3197) -- (1.3723,0.3700) -- (1.4105,0.4733) -- (1.3974,0.4215) -- cycle;
\draw[draw=cbpInk!70,line width=.45pt] (1.3596,0.3197) -- (1.3723,0.3700);
\draw[draw=cbpInk!70,line width=.45pt] (1.3974,0.4215) -- (1.3596,0.3197);
\path[cbpmesh=20.4] (-4.3250,0.5382) -- (-4.3117,0.3788) -- (-4.3935,0.2522) -- (-4.4071,0.4146) -- cycle;
\path[cbpmesh=34.9] (4.3117,0.3788) -- (4.3250,0.5382) -- (4.4071,0.4146) -- (4.3935,0.2522) -- cycle;
\path[cbpmesh=35.7] (-1.4632,0.6215) -- (-1.4587,0.5676) -- (-1.5162,0.6526) -- (-1.5209,0.7087) -- cycle;
\path[cbpmesh=21.2] (1.4587,0.5676) -- (1.4632,0.6215) -- (1.5209,0.7087) -- (1.5162,0.6526) -- cycle;
\path[cbpmesh=16.7] (-2.7236,1.0754) -- (-2.7320,0.9748) -- (-2.8495,0.9351) -- (-2.8408,1.0401) -- cycle;
\path[cbpmesh=22.7] (2.7320,0.9748) -- (2.7236,1.0754) -- (2.8408,1.0401) -- (2.8495,0.9351) -- cycle;
\path[cbpmesh=27.8] (-1.6570,0.9167) -- (-1.6621,0.8555) -- (-1.7445,0.9146) -- (-1.7391,0.9789) -- cycle;
\path[cbpmesh=17.1] (1.6621,0.8555) -- (1.6570,0.9167) -- (1.7391,0.9789) -- (1.7445,0.9146) -- cycle;
\path[cbpmesh=30.4] (-4.5210,-0.2123) -- (-4.4791,-0.3783) -- (-4.4674,-0.4811) -- (-4.5091,-0.3156) -- cycle;
\draw[draw=cbpInk!70,line width=.45pt] (-4.4674,-0.4811) -- (-4.5091,-0.3156);
\draw[draw=cbpInk!70,line width=.45pt] (-4.5091,-0.3156) -- (-4.5210,-0.2123);
\path[cbpmesh=38] (4.4791,-0.3783) -- (4.5210,-0.2123) -- (4.5091,-0.3156) -- (4.4674,-0.4811) -- cycle;
\draw[draw=cbpInk!70,line width=.45pt] (4.5210,-0.2123) -- (4.5091,-0.3156);
\draw[draw=cbpInk!70,line width=.45pt] (4.5091,-0.3156) -- (4.4674,-0.4811);
\path[cbpmesh=16.8] (-2.5994,1.1090) -- (-2.6074,1.0130) -- (-2.7320,0.9748) -- (-2.7236,1.0754) -- cycle;
\path[cbpmesh=22.3] (2.6074,1.0130) -- (2.5994,1.1090) -- (2.7236,1.0754) -- (2.7320,0.9748) -- cycle;
\path[cbpmesh=25.5] (-1.7391,0.9789) -- (-1.7445,0.9146) -- (-1.8339,0.9639) -- (-1.8282,1.0314) -- cycle;
\path[cbpmesh=16.6] (1.7445,0.9146) -- (1.7391,0.9789) -- (1.8282,1.0314) -- (1.8339,0.9639) -- cycle;
\path[cbpmesh=33.2] (-1.5209,0.7087) -- (-1.5162,0.6526) -- (-1.5825,0.7283) -- (-1.5874,0.7867) -- cycle;
\path[cbpmesh=19.8] (1.5162,0.6526) -- (1.5209,0.7087) -- (1.5874,0.7867) -- (1.5825,0.7283) -- cycle;
\path[cbpmesh=17.8] (-4.2498,0.6038) -- (-4.2367,0.4473) -- (-4.3117,0.3788) -- (-4.3250,0.5382) -- cycle;
\path[cbpmesh=29.6] (4.2367,0.4473) -- (4.2498,0.6038) -- (4.3250,0.5382) -- (4.3117,0.3788) -- cycle;
\path[cbpmesh=38] (-1.4105,0.4733) -- (-1.3974,0.4215) -- (-1.4452,0.5140) -- (-1.4587,0.5676) -- cycle;
\path[cbpmesh=23.5] (1.3974,0.4215) -- (1.4105,0.4733) -- (1.4587,0.5676) -- (1.4452,0.5140) -- cycle;
\path[cbpmesh=17.0] (-2.4775,1.1309) -- (-2.4852,1.0394) -- (-2.6074,1.0130) -- (-2.5994,1.1090) -- cycle;
\path[cbpmesh=21.0] (2.4852,1.0394) -- (2.4775,1.1309) -- (2.5994,1.1090) -- (2.6074,1.0130) -- cycle;
\path[cbpmesh=23.5] (-1.8282,1.0314) -- (-1.8339,0.9639) -- (-1.9297,1.0031) -- (-1.9238,1.0742) -- cycle;
\path[cbpmesh=16.4] (1.8339,0.9639) -- (1.8282,1.0314) -- (1.9238,1.0742) -- (1.9297,1.0031) -- cycle;
\path[cbpmesh=21.7] (-1.9238,1.0742) -- (-1.9297,1.0031) -- (-2.0314,1.0320) -- (-2.0252,1.1068) -- cycle;
\path[cbpmesh=16.5] (1.9297,1.0031) -- (1.9238,1.0742) -- (2.0252,1.1068) -- (2.0314,1.0320) -- cycle;
\path[cbpmesh=17.4] (-2.3586,1.1414) -- (-2.3659,1.0542) -- (-2.4852,1.0394) -- (-2.4775,1.1309) -- cycle;
\path[cbpmesh=19.7] (2.3659,1.0542) -- (2.3586,1.1414) -- (2.4775,1.1309) -- (2.4852,1.0394) -- cycle;
\path[cbpmesh=27.7] (-4.5166,-0.1034) -- (-4.4748,-0.2692) -- (-4.4791,-0.3783) -- (-4.5210,-0.2123) -- cycle;
\path[cbpmesh=38] (4.4748,-0.2692) -- (4.5166,-0.1034) -- (4.5210,-0.2123) -- (4.4791,-0.3783) -- cycle;
\path[cbpmesh=20.2] (-2.0252,1.1068) -- (-2.0314,1.0320) -- (-2.1384,1.0503) -- (-2.1318,1.1291) -- cycle;
\path[cbpmesh=17.0] (2.0314,1.0320) -- (2.0252,1.1068) -- (2.1318,1.1291) -- (2.1384,1.0503) -- cycle;
\path[cbpmesh=17.3] (-4.1776,0.6584) -- (-4.1647,0.5045) -- (-4.2367,0.4473) -- (-4.2498,0.6038) -- cycle;
\path[cbpmesh=28.2] (4.1647,0.5045) -- (4.1776,0.6584) -- (4.2498,0.6038) -- (4.2367,0.4473) -- cycle;
\path[cbpmesh=18.0] (-2.2432,1.1407) -- (-2.2501,1.0578) -- (-2.3659,1.0542) -- (-2.3586,1.1414) -- cycle;
\path[cbpmesh=18.6] (2.2501,1.0578) -- (2.2432,1.1407) -- (2.3586,1.1414) -- (2.3659,1.0542) -- cycle;
\path[cbpmesh=19.0] (-2.1318,1.1291) -- (-2.1384,1.0503) -- (-2.2501,1.0578) -- (-2.2432,1.1407) -- cycle;
\path[cbpmesh=17.7] (2.1384,1.0503) -- (2.1318,1.1291) -- (2.2432,1.1407) -- (2.2501,1.0578) -- cycle;
\path[cbpmesh=30.7] (-1.5874,0.7867) -- (-1.5825,0.7283) -- (-1.6570,0.7942) -- (-1.6621,0.8555) -- cycle;
\path[cbpmesh=18.6] (1.5825,0.7283) -- (1.5874,0.7867) -- (1.6621,0.8555) -- (1.6570,0.7942) -- cycle;
\path[cbpmesh=16.9] (-4.1079,0.7029) -- (-4.0952,0.5516) -- (-4.1647,0.5045) -- (-4.1776,0.6584) -- cycle;
\path[cbpmesh=26.7] (4.0952,0.5516) -- (4.1079,0.7029) -- (4.1776,0.6584) -- (4.1647,0.5045) -- cycle;
\path[cbpmesh=36.4] (-1.4587,0.5676) -- (-1.4452,0.5140) -- (-1.5022,0.5970) -- (-1.5162,0.6526) -- cycle;
\path[cbpmesh=21.9] (1.4452,0.5140) -- (1.4587,0.5676) -- (1.5162,0.6526) -- (1.5022,0.5970) -- cycle;
\path[cbpmesh=16.6] (-4.0401,0.7386) -- (-4.0276,0.5897) -- (-4.0952,0.5516) -- (-4.1079,0.7029) -- cycle;
\path[cbpmesh=25.1] (4.0276,0.5897) -- (4.0401,0.7386) -- (4.1079,0.7029) -- (4.0952,0.5516) -- cycle;
\path[cbpmesh=28.3] (-1.6621,0.8555) -- (-1.6570,0.7942) -- (-1.7391,0.8503) -- (-1.7445,0.9146) -- cycle;
\path[cbpmesh=17.7] (1.6570,0.7942) -- (1.6621,0.8555) -- (1.7445,0.9146) -- (1.7391,0.8503) -- cycle;
\path[cbpmesh=25.2] (-4.4947,0.0107) -- (-4.4531,-0.1544) -- (-4.4748,-0.2692) -- (-4.5166,-0.1034) -- cycle;
\path[cbpmesh=38] (4.4531,-0.1544) -- (4.4947,0.0107) -- (4.5166,-0.1034) -- (4.4748,-0.2692) -- cycle;
\path[cbpmesh=16.4] (-3.9736,0.7663) -- (-3.9613,0.6199) -- (-4.0276,0.5897) -- (-4.0401,0.7386) -- cycle;
\path[cbpmesh=23.6] (3.9613,0.6199) -- (3.9736,0.7663) -- (4.0401,0.7386) -- (4.0276,0.5897) -- cycle;
\path[cbpmesh=16.5] (-3.9078,0.7873) -- (-3.8957,0.6433) -- (-3.9613,0.6199) -- (-3.9736,0.7663) -- cycle;
\path[cbpmesh=22.2] (3.8957,0.6433) -- (3.9078,0.7873) -- (3.9736,0.7663) -- (3.9613,0.6199) -- cycle;
\path[cbpmesh=38] (-1.3974,0.4215) -- (-1.3758,0.3704) -- (-1.4228,0.4611) -- (-1.4452,0.5140) -- cycle;
\draw[draw=cbpInk!70,line width=.45pt] (-1.3974,0.4215) -- (-1.3758,0.3704);
\draw[draw=cbpInk!70,line width=.45pt] (-1.3758,0.3704) -- (-1.4228,0.4611);
\path[cbpmesh=24.4] (1.3758,0.3704) -- (1.3974,0.4215) -- (1.4452,0.5140) -- (1.4228,0.4611) -- cycle;
\draw[draw=cbpInk!70,line width=.45pt] (1.3758,0.3704) -- (1.3974,0.4215);
\draw[draw=cbpInk!70,line width=.45pt] (1.4228,0.4611) -- (1.3758,0.3704);
\path[cbpmesh=16.7] (-3.8421,0.8025) -- (-3.8302,0.6610) -- (-3.8957,0.6433) -- (-3.9078,0.7873) -- cycle;
\path[cbpmesh=21.0] (3.8302,0.6610) -- (3.8421,0.8025) -- (3.9078,0.7873) -- (3.8957,0.6433) -- cycle;
\path[cbpmesh=33.9] (-1.5162,0.6526) -- (-1.5022,0.5970) -- (-1.5679,0.6702) -- (-1.5825,0.7283) -- cycle;
\path[cbpmesh=20.5] (1.5022,0.5970) -- (1.5162,0.6526) -- (1.5825,0.7283) -- (1.5679,0.6702) -- cycle;
\path[cbpmesh=26.0] (-1.7445,0.9146) -- (-1.7391,0.8503) -- (-1.8282,0.8963) -- (-1.8339,0.9639) -- cycle;
\path[cbpmesh=17.1] (1.7391,0.8503) -- (1.7445,0.9146) -- (1.8339,0.9639) -- (1.8282,0.8963) -- cycle;
\path[cbpmesh=17.0] (-3.7759,0.8132) -- (-3.7643,0.6741) -- (-3.8302,0.6610) -- (-3.8421,0.8025) -- cycle;
\path[cbpmesh=20.0] (3.7643,0.6741) -- (3.7759,0.8132) -- (3.8421,0.8025) -- (3.8302,0.6610) -- cycle;
\path[cbpmesh=17.7] (-3.3363,0.8408) -- (-3.3261,0.7179) -- (-3.4071,0.7092) -- (-3.4176,0.8351) -- cycle;
\path[cbpmesh=19.0] (3.3261,0.7179) -- (3.3363,0.8408) -- (3.4176,0.8351) -- (3.4071,0.7092) -- cycle;
\path[cbpmesh=17.8] (-3.4176,0.8351) -- (-3.4071,0.7092) -- (-3.4841,0.7026) -- (-3.4949,0.8313) -- cycle;
\path[cbpmesh=18.8] (3.4071,0.7092) -- (3.4176,0.8351) -- (3.4949,0.8313) -- (3.4841,0.7026) -- cycle;
\path[cbpmesh=17.4] (-3.2505,0.8494) -- (-3.2405,0.7297) -- (-3.3261,0.7179) -- (-3.3363,0.8408) -- cycle;
\path[cbpmesh=19.3] (3.2405,0.7297) -- (3.2505,0.8494) -- (3.3363,0.8408) -- (3.3261,0.7179) -- cycle;
\path[cbpmesh=17.9] (-3.4949,0.8313) -- (-3.4841,0.7026) -- (-3.5578,0.6969) -- (-3.5688,0.8283) -- cycle;
\path[cbpmesh=18.7] (3.4841,0.7026) -- (3.4949,0.8313) -- (3.5688,0.8283) -- (3.5578,0.6969) -- cycle;
\path[cbpmesh=17.4] (-3.7087,0.8203) -- (-3.6973,0.6837) -- (-3.7643,0.6741) -- (-3.7759,0.8132) -- cycle;
\path[cbpmesh=19.4] (3.6973,0.6837) -- (3.7087,0.8203) -- (3.7759,0.8132) -- (3.7643,0.6741) -- cycle;
\path[cbpmesh=17.9] (-3.5688,0.8283) -- (-3.5578,0.6969) -- (-3.6287,0.6909) -- (-3.6399,0.8250) -- cycle;
\path[cbpmesh=18.7] (3.5578,0.6969) -- (3.5688,0.8283) -- (3.6399,0.8250) -- (3.6287,0.6909) -- cycle;
\path[cbpmesh=17.7] (-3.6399,0.8250) -- (-3.6287,0.6909) -- (-3.6973,0.6837) -- (-3.7087,0.8203) -- cycle;
\path[cbpmesh=19.0] (3.6287,0.6909) -- (3.6399,0.8250) -- (3.7087,0.8203) -- (3.6973,0.6837) -- cycle;
\path[cbpmesh=17.2] (-3.1595,0.8621) -- (-3.1497,0.7457) -- (-3.2405,0.7297) -- (-3.2505,0.8494) -- cycle;
\path[cbpmesh=19.8] (3.1497,0.7457) -- (3.1595,0.8621) -- (3.2505,0.8494) -- (3.2405,0.7297) -- cycle;
\path[cbpmesh=16.9] (-3.0627,0.8799) -- (-3.0533,0.7671) -- (-3.1497,0.7457) -- (-3.1595,0.8621) -- cycle;
\path[cbpmesh=20.3] (3.0533,0.7671) -- (3.0627,0.8799) -- (3.1595,0.8621) -- (3.1497,0.7457) -- cycle;
\path[cbpmesh=16.7] (-2.9596,0.9039) -- (-2.9505,0.7949) -- (-3.0533,0.7671) -- (-3.0627,0.8799) -- cycle;
\path[cbpmesh=21.0] (2.9505,0.7949) -- (2.9596,0.9039) -- (3.0627,0.8799) -- (3.0533,0.7671) -- cycle;
\path[cbpmesh=23.8] (-1.8339,0.9639) -- (-1.8282,0.8963) -- (-1.9238,0.9320) -- (-1.9297,1.0031) -- cycle;
\path[cbpmesh=16.8] (1.8282,0.8963) -- (1.8339,0.9639) -- (1.9297,1.0031) -- (1.9238,0.9320) -- cycle;
\path[cbpmesh=16.5] (-2.8495,0.9351) -- (-2.8408,0.8302) -- (-2.9505,0.7949) -- (-2.9596,0.9039) -- cycle;
\path[cbpmesh=21.7] (2.8408,0.8302) -- (2.8495,0.9351) -- (2.9596,0.9039) -- (2.9505,0.7949) -- cycle;
\path[cbpmesh=23.0] (-4.4541,0.1294) -- (-4.4129,-0.0342) -- (-4.4531,-0.1544) -- (-4.4947,0.0107) -- cycle;
\path[cbpmesh=38] (4.4129,-0.0342) -- (4.4541,0.1294) -- (4.4947,0.0107) -- (4.4531,-0.1544) -- cycle;
\path[cbpmesh=16.4] (-2.7320,0.9748) -- (-2.7236,0.8741) -- (-2.8408,0.8302) -- (-2.8495,0.9351) -- cycle;
\path[cbpmesh=22.4] (2.7236,0.8741) -- (2.7320,0.9748) -- (2.8495,0.9351) -- (2.8408,0.8302) -- cycle;
\path[cbpmesh=31.3] (-1.5825,0.7283) -- (-1.5679,0.6702) -- (-1.6417,0.7334) -- (-1.6570,0.7942) -- cycle;
\path[cbpmesh=19.2] (1.5679,0.6702) -- (1.5825,0.7283) -- (1.6570,0.7942) -- (1.6417,0.7334) -- cycle;
\path[cbpmesh=21.9] (-1.9297,1.0031) -- (-1.9238,0.9320) -- (-2.0252,0.9572) -- (-2.0314,1.0320) -- cycle;
\path[cbpmesh=16.8] (1.9238,0.9320) -- (1.9297,1.0031) -- (2.0314,1.0320) -- (2.0252,0.9572) -- cycle;
\path[cbpmesh=37.0] (-1.4452,0.5140) -- (-1.4228,0.4611) -- (-1.4789,0.5420) -- (-1.5022,0.5970) -- cycle;
\path[cbpmesh=22.8] (1.4228,0.4611) -- (1.4452,0.5140) -- (1.5022,0.5970) -- (1.4789,0.5420) -- cycle;
\path[cbpmesh=16.5] (-2.6074,1.0130) -- (-2.5994,0.9169) -- (-2.7236,0.8741) -- (-2.7320,0.9748) -- cycle;
\path[cbpmesh=22.0] (2.5994,0.9169) -- (2.6074,1.0130) -- (2.7320,0.9748) -- (2.7236,0.8741) -- cycle;
\path[cbpmesh=20.3] (-2.0314,1.0320) -- (-2.0252,0.9572) -- (-2.1318,0.9715) -- (-2.1384,1.0503) -- cycle;
\path[cbpmesh=17.1] (2.0252,0.9572) -- (2.0314,1.0320) -- (2.1384,1.0503) -- (2.1318,0.9715) -- cycle;
\path[cbpmesh=16.8] (-2.4852,1.0394) -- (-2.4775,0.9478) -- (-2.5994,0.9169) -- (-2.6074,1.0130) -- cycle;
\path[cbpmesh=20.8] (2.4775,0.9478) -- (2.4852,1.0394) -- (2.6074,1.0130) -- (2.5994,0.9169) -- cycle;
\path[cbpmesh=19.0] (-2.1384,1.0503) -- (-2.1318,0.9715) -- (-2.2432,0.9749) -- (-2.2501,1.0578) -- cycle;
\path[cbpmesh=17.7] (2.1318,0.9715) -- (2.1384,1.0503) -- (2.2501,1.0578) -- (2.2432,0.9749) -- cycle;
\path[cbpmesh=17.3] (-2.3659,1.0542) -- (-2.3586,0.9671) -- (-2.4775,0.9478) -- (-2.4852,1.0394) -- cycle;
\path[cbpmesh=19.6] (2.3586,0.9671) -- (2.3659,1.0542) -- (2.4852,1.0394) -- (2.4775,0.9478) -- cycle;
\path[cbpmesh=18.0] (-2.2501,1.0578) -- (-2.2432,0.9749) -- (-2.3586,0.9671) -- (-2.3659,1.0542) -- cycle;
\path[cbpmesh=18.6] (2.2432,0.9749) -- (2.2501,1.0578) -- (2.3659,1.0542) -- (2.3586,0.9671) -- cycle;
\path[cbpmesh=28.8] (-1.6570,0.7942) -- (-1.6417,0.7334) -- (-1.7230,0.7865) -- (-1.7391,0.8503) -- cycle;
\path[cbpmesh=18.2] (1.6417,0.7334) -- (1.6570,0.7942) -- (1.7391,0.8503) -- (1.7230,0.7865) -- cycle;
\path[cbpmesh=21.2] (-4.3935,0.2522) -- (-4.3528,0.0909) -- (-4.4129,-0.0342) -- (-4.4541,0.1294) -- cycle;
\path[cbpmesh=36.5] (4.3528,0.0909) -- (4.3935,0.2522) -- (4.4541,0.1294) -- (4.4129,-0.0342) -- cycle;
\path[cbpmesh=34.4] (-1.5022,0.5970) -- (-1.4789,0.5420) -- (-1.5435,0.6128) -- (-1.5679,0.6702) -- cycle;
\path[cbpmesh=21.3] (1.4789,0.5420) -- (1.5022,0.5970) -- (1.5679,0.6702) -- (1.5435,0.6128) -- cycle;
\path[cbpmesh=32.6] (-4.4674,-0.4811) -- (-4.3981,-0.6447) -- (-4.3719,-0.7400) -- (-4.4407,-0.5774) -- cycle;
\draw[draw=cbpInk!70,line width=.45pt] (-4.3719,-0.7400) -- (-4.4407,-0.5774);
\draw[draw=cbpInk!70,line width=.45pt] (-4.4407,-0.5774) -- (-4.4674,-0.4811);
\path[cbpmesh=38] (4.3981,-0.6447) -- (4.4674,-0.4811) -- (4.4407,-0.5774) -- (4.3719,-0.7400) -- cycle;
\draw[draw=cbpInk!70,line width=.45pt] (4.4674,-0.4811) -- (4.4407,-0.5774);
\draw[draw=cbpInk!70,line width=.45pt] (4.4407,-0.5774) -- (4.3719,-0.7400);
\path[cbpmesh=26.4] (-1.7391,0.8503) -- (-1.7230,0.7865) -- (-1.8113,0.8292) -- (-1.8282,0.8963) -- cycle;
\path[cbpmesh=17.5] (1.7230,0.7865) -- (1.7391,0.8503) -- (1.8282,0.8963) -- (1.8113,0.8292) -- cycle;
\path[cbpmesh=19.7] (-4.3117,0.3788) -- (-4.2718,0.2205) -- (-4.3528,0.0909) -- (-4.3935,0.2522) -- cycle;
\path[cbpmesh=34.1] (4.2718,0.2205) -- (4.3117,0.3788) -- (4.3935,0.2522) -- (4.3528,0.0909) -- cycle;
\path[cbpmesh=37.6] (-1.4228,0.4611) -- (-1.3916,0.4092) -- (-1.4465,0.4880) -- (-1.4789,0.5420) -- cycle;
\draw[draw=cbpInk!70,line width=.45pt] (-1.4228,0.4611) -- (-1.3916,0.4092);
\draw[draw=cbpInk!70,line width=.45pt] (-1.3916,0.4092) -- (-1.4465,0.4880);
\path[cbpmesh=23.6] (1.3916,0.4092) -- (1.4228,0.4611) -- (1.4789,0.5420) -- (1.4465,0.4880) -- cycle;
\draw[draw=cbpInk!70,line width=.45pt] (1.3916,0.4092) -- (1.4228,0.4611);
\draw[draw=cbpInk!70,line width=.45pt] (1.4465,0.4880) -- (1.3916,0.4092);
\path[cbpmesh=31.8] (-1.5679,0.6702) -- (-1.5435,0.6128) -- (-1.6162,0.6733) -- (-1.6417,0.7334) -- cycle;
\path[cbpmesh=19.9] (1.5435,0.6128) -- (1.5679,0.6702) -- (1.6417,0.7334) -- (1.6162,0.6733) -- cycle;
\path[cbpmesh=24.2] (-1.8282,0.8963) -- (-1.8113,0.8292) -- (-1.9060,0.8614) -- (-1.9238,0.9320) -- cycle;
\path[cbpmesh=17.1] (1.8113,0.8292) -- (1.8282,0.8963) -- (1.9238,0.9320) -- (1.9060,0.8614) -- cycle;
\path[cbpmesh=29.7] (-4.4791,-0.3783) -- (-4.4097,-0.5422) -- (-4.3981,-0.6447) -- (-4.4674,-0.4811) -- cycle;
\path[cbpmesh=38] (4.4097,-0.5422) -- (4.4791,-0.3783) -- (4.4674,-0.4811) -- (4.3981,-0.6447) -- cycle;
\path[cbpmesh=22.2] (-1.9238,0.9320) -- (-1.9060,0.8614) -- (-2.0064,0.8828) -- (-2.0252,0.9572) -- cycle;
\path[cbpmesh=17.1] (1.9060,0.8614) -- (1.9238,0.9320) -- (2.0252,0.9572) -- (2.0064,0.8828) -- cycle;
\path[cbpmesh=17.2] (-4.2367,0.4473) -- (-4.1975,0.2917) -- (-4.2718,0.2205) -- (-4.3117,0.3788) -- cycle;
\path[cbpmesh=29.0] (4.1975,0.2917) -- (4.2367,0.4473) -- (4.3117,0.3788) -- (4.2718,0.2205) -- cycle;
\path[cbpmesh=29.3] (-1.6417,0.7334) -- (-1.6162,0.6733) -- (-1.6963,0.7234) -- (-1.7230,0.7865) -- cycle;
\path[cbpmesh=18.8] (1.6162,0.6733) -- (1.6417,0.7334) -- (1.7230,0.7865) -- (1.6963,0.7234) -- cycle;
\path[cbpmesh=35.0] (-1.4789,0.5420) -- (-1.4465,0.4880) -- (-1.5097,0.5564) -- (-1.5435,0.6128) -- cycle;
\path[cbpmesh=22.0] (1.4465,0.4880) -- (1.4789,0.5420) -- (1.5435,0.6128) -- (1.5097,0.5564) -- cycle;
\path[cbpmesh=20.5] (-2.0252,0.9572) -- (-2.0064,0.8828) -- (-2.1121,0.8932) -- (-2.1318,0.9715) -- cycle;
\path[cbpmesh=17.3] (2.0064,0.8828) -- (2.0252,0.9572) -- (2.1318,0.9715) -- (2.1121,0.8932) -- cycle;
\path[cbpmesh=16.8] (-4.1647,0.5045) -- (-4.1262,0.3516) -- (-4.1975,0.2917) -- (-4.2367,0.4473) -- cycle;
\path[cbpmesh=27.6] (4.1262,0.3516) -- (4.1647,0.5045) -- (4.2367,0.4473) -- (4.1975,0.2917) -- cycle;
\path[cbpmesh=26.9] (-4.4748,-0.2692) -- (-4.4054,-0.4330) -- (-4.4097,-0.5422) -- (-4.4791,-0.3783) -- cycle;
\path[cbpmesh=38] (4.4054,-0.4330) -- (4.4748,-0.2692) -- (4.4791,-0.3783) -- (4.4097,-0.5422) -- cycle;
\path[cbpmesh=16.7] (-3.0533,0.7671) -- (-3.0250,0.6550) -- (-3.1206,0.6301) -- (-3.1497,0.7457) -- cycle;
\path[cbpmesh=20.2] (3.0250,0.6550) -- (3.0533,0.7671) -- (3.1497,0.7457) -- (3.1206,0.6301) -- cycle;
\path[cbpmesh=16.5] (-2.9505,0.7949) -- (-2.9231,0.6865) -- (-3.0250,0.6550) -- (-3.0533,0.7671) -- cycle;
\path[cbpmesh=20.8] (2.9231,0.6865) -- (2.9505,0.7949) -- (3.0533,0.7671) -- (3.0250,0.6550) -- cycle;
\path[cbpmesh=17.0] (-3.1497,0.7457) -- (-3.1206,0.6301) -- (-3.2105,0.6107) -- (-3.2405,0.7297) -- cycle;
\path[cbpmesh=19.6] (3.1206,0.6301) -- (3.1497,0.7457) -- (3.2405,0.7297) -- (3.2105,0.6107) -- cycle;
\path[cbpmesh=19.1] (-2.1318,0.9715) -- (-2.1121,0.8932) -- (-2.2224,0.8925) -- (-2.2432,0.9749) -- cycle;
\path[cbpmesh=17.8] (2.1121,0.8932) -- (2.1318,0.9715) -- (2.2432,0.9749) -- (2.2224,0.8925) -- cycle;
\path[cbpmesh=16.3] (-2.8408,0.8302) -- (-2.8145,0.7259) -- (-2.9231,0.6865) -- (-2.9505,0.7949) -- cycle;
\path[cbpmesh=21.4] (2.8145,0.7259) -- (2.8408,0.8302) -- (2.9505,0.7949) -- (2.9231,0.6865) -- cycle;
\path[cbpmesh=17.3] (-3.2405,0.7297) -- (-3.2105,0.6107) -- (-3.2953,0.5958) -- (-3.3261,0.7179) -- cycle;
\path[cbpmesh=19.2] (3.2105,0.6107) -- (3.2405,0.7297) -- (3.3261,0.7179) -- (3.2953,0.5958) -- cycle;
\path[cbpmesh=16.4] (-4.0952,0.5516) -- (-4.0573,0.4013) -- (-4.1262,0.3516) -- (-4.1647,0.5045) -- cycle;
\path[cbpmesh=26.2] (4.0573,0.4013) -- (4.0952,0.5516) -- (4.1647,0.5045) -- (4.1262,0.3516) -- cycle;
\path[cbpmesh=17.6] (-3.3261,0.7179) -- (-3.2953,0.5958) -- (-3.3755,0.5841) -- (-3.4071,0.7092) -- cycle;
\path[cbpmesh=18.9] (3.2953,0.5958) -- (3.3261,0.7179) -- (3.4071,0.7092) -- (3.3755,0.5841) -- cycle;
\path[cbpmesh=16.2] (-2.7236,0.8741) -- (-2.6984,0.7741) -- (-2.8145,0.7259) -- (-2.8408,0.8302) -- cycle;
\path[cbpmesh=22.1] (2.6984,0.7741) -- (2.7236,0.8741) -- (2.8408,0.8302) -- (2.8145,0.7259) -- cycle;
\path[cbpmesh=18.0] (-2.2432,0.9749) -- (-2.2224,0.8925) -- (-2.3368,0.8805) -- (-2.3586,0.9671) -- cycle;
\path[cbpmesh=18.6] (2.2224,0.8925) -- (2.2432,0.9749) -- (2.3586,0.9671) -- (2.3368,0.8805) -- cycle;
\path[cbpmesh=17.8] (-3.4071,0.7092) -- (-3.3755,0.5841) -- (-3.4519,0.5746) -- (-3.4841,0.7026) -- cycle;
\path[cbpmesh=18.7] (3.3755,0.5841) -- (3.4071,0.7092) -- (3.4841,0.7026) -- (3.4519,0.5746) -- cycle;
\path[cbpmesh=16.2] (-2.5994,0.9169) -- (-2.5753,0.8215) -- (-2.6984,0.7741) -- (-2.7236,0.8741) -- cycle;
\path[cbpmesh=21.7] (2.5753,0.8215) -- (2.5994,0.9169) -- (2.7236,0.8741) -- (2.6984,0.7741) -- cycle;
\path[cbpmesh=16.1] (-4.0276,0.5897) -- (-3.9904,0.4418) -- (-4.0573,0.4013) -- (-4.0952,0.5516) -- cycle;
\path[cbpmesh=24.7] (3.9904,0.4418) -- (4.0276,0.5897) -- (4.0952,0.5516) -- (4.0573,0.4013) -- cycle;
\path[cbpmesh=26.8] (-1.7230,0.7865) -- (-1.6963,0.7234) -- (-1.7832,0.7629) -- (-1.8113,0.8292) -- cycle;
\path[cbpmesh=18.0] (1.6963,0.7234) -- (1.7230,0.7865) -- (1.8113,0.8292) -- (1.7832,0.7629) -- cycle;
\path[cbpmesh=17.9] (-3.4841,0.7026) -- (-3.4519,0.5746) -- (-3.5249,0.5662) -- (-3.5578,0.6969) -- cycle;
\path[cbpmesh=18.6] (3.4519,0.5746) -- (3.4841,0.7026) -- (3.5578,0.6969) -- (3.5249,0.5662) -- cycle;
\path[cbpmesh=17.1] (-2.3586,0.9671) -- (-2.3368,0.8805) -- (-2.4546,0.8569) -- (-2.4775,0.9478) -- cycle;
\path[cbpmesh=19.5] (2.3368,0.8805) -- (2.3586,0.9671) -- (2.4775,0.9478) -- (2.4546,0.8569) -- cycle;
\path[cbpmesh=16.6] (-2.4775,0.9478) -- (-2.4546,0.8569) -- (-2.5753,0.8215) -- (-2.5994,0.9169) -- cycle;
\path[cbpmesh=20.5] (2.4546,0.8569) -- (2.4775,0.9478) -- (2.5994,0.9169) -- (2.5753,0.8215) -- cycle;
\path[cbpmesh=17.8] (-3.5578,0.6969) -- (-3.5249,0.5662) -- (-3.5951,0.5577) -- (-3.6287,0.6909) -- cycle;
\path[cbpmesh=18.7] (3.5249,0.5662) -- (3.5578,0.6969) -- (3.6287,0.6909) -- (3.5951,0.5577) -- cycle;
\path[cbpmesh=16.1] (-3.9613,0.6199) -- (-3.9247,0.4745) -- (-3.9904,0.4418) -- (-4.0276,0.5897) -- cycle;
\path[cbpmesh=23.2] (3.9247,0.4745) -- (3.9613,0.6199) -- (4.0276,0.5897) -- (3.9904,0.4418) -- cycle;
\path[cbpmesh=17.6] (-3.6287,0.6909) -- (-3.5951,0.5577) -- (-3.6631,0.5480) -- (-3.6973,0.6837) -- cycle;
\path[cbpmesh=18.9] (3.5951,0.5577) -- (3.6287,0.6909) -- (3.6973,0.6837) -- (3.6631,0.5480) -- cycle;
\path[cbpmesh=16.2] (-3.8957,0.6433) -- (-3.8597,0.5003) -- (-3.9247,0.4745) -- (-3.9613,0.6199) -- cycle;
\path[cbpmesh=21.9] (3.8597,0.5003) -- (3.8957,0.6433) -- (3.9613,0.6199) -- (3.9247,0.4745) -- cycle;
\path[cbpmesh=17.3] (-3.6973,0.6837) -- (-3.6631,0.5480) -- (-3.7294,0.5359) -- (-3.7643,0.6741) -- cycle;
\path[cbpmesh=19.3] (3.6631,0.5480) -- (3.6973,0.6837) -- (3.7643,0.6741) -- (3.7294,0.5359) -- cycle;
\path[cbpmesh=16.5] (-3.8302,0.6610) -- (-3.7948,0.5204) -- (-3.8597,0.5003) -- (-3.8957,0.6433) -- cycle;
\path[cbpmesh=20.7] (3.7948,0.5204) -- (3.8302,0.6610) -- (3.8957,0.6433) -- (3.8597,0.5003) -- cycle;
\path[cbpmesh=16.9] (-3.7643,0.6741) -- (-3.7294,0.5359) -- (-3.7948,0.5204) -- (-3.8302,0.6610) -- cycle;
\path[cbpmesh=19.9] (3.7294,0.5359) -- (3.7643,0.6741) -- (3.8302,0.6610) -- (3.7948,0.5204) -- cycle;
\path[cbpmesh=32.3] (-1.5435,0.6128) -- (-1.5097,0.5564) -- (-1.5808,0.6143) -- (-1.6162,0.6733) -- cycle;
\path[cbpmesh=20.6] (1.5097,0.5564) -- (1.5435,0.6128) -- (1.6162,0.6733) -- (1.5808,0.6143) -- cycle;
\path[cbpmesh=24.5] (-4.4531,-0.1544) -- (-4.3840,-0.3174) -- (-4.4054,-0.4330) -- (-4.4748,-0.2692) -- cycle;
\path[cbpmesh=38] (4.3840,-0.3174) -- (4.4531,-0.1544) -- (4.4748,-0.2692) -- (4.4054,-0.4330) -- cycle;
\path[cbpmesh=24.5] (-1.8113,0.8292) -- (-1.7832,0.7629) -- (-1.8764,0.7916) -- (-1.9060,0.8614) -- cycle;
\path[cbpmesh=17.5] (1.7832,0.7629) -- (1.8113,0.8292) -- (1.9060,0.8614) -- (1.8764,0.7916) -- cycle;
\path[cbpmesh=35.5] (-1.4465,0.4880) -- (-1.4051,0.4354) -- (-1.4666,0.5015) -- (-1.5097,0.5564) -- cycle;
\draw[draw=cbpInk!70,line width=.45pt] (-1.4465,0.4880) -- (-1.4051,0.4354);
\draw[draw=cbpInk!70,line width=.45pt] (-1.4051,0.4354) -- (-1.4666,0.5015);
\path[cbpmesh=22.8] (1.4051,0.4354) -- (1.4465,0.4880) -- (1.5097,0.5564) -- (1.4666,0.5015) -- cycle;
\draw[draw=cbpInk!70,line width=.45pt] (1.4051,0.4354) -- (1.4465,0.4880);
\draw[draw=cbpInk!70,line width=.45pt] (1.4666,0.5015) -- (1.4051,0.4354);
\path[cbpmesh=29.7] (-1.6162,0.6733) -- (-1.5808,0.6143) -- (-1.6591,0.6615) -- (-1.6963,0.7234) -- cycle;
\path[cbpmesh=19.5] (1.5808,0.6143) -- (1.6162,0.6733) -- (1.6963,0.7234) -- (1.6591,0.6615) -- cycle;
\path[cbpmesh=22.4] (-1.9060,0.8614) -- (-1.8764,0.7916) -- (-1.9753,0.8093) -- (-2.0064,0.8828) -- cycle;
\path[cbpmesh=17.4] (1.8764,0.7916) -- (1.9060,0.8614) -- (2.0064,0.8828) -- (1.9753,0.8093) -- cycle;
\path[cbpmesh=22.3] (-4.4129,-0.0342) -- (-4.3444,-0.1957) -- (-4.3840,-0.3174) -- (-4.4531,-0.1544) -- cycle;
\path[cbpmesh=38] (4.3444,-0.1957) -- (4.4129,-0.0342) -- (4.4531,-0.1544) -- (4.3840,-0.3174) -- cycle;
\path[cbpmesh=20.6] (-2.0064,0.8828) -- (-1.9753,0.8093) -- (-2.0793,0.8159) -- (-2.1121,0.8932) -- cycle;
\path[cbpmesh=17.5] (1.9753,0.8093) -- (2.0064,0.8828) -- (2.1121,0.8932) -- (2.0793,0.8159) -- cycle;
\path[cbpmesh=34.9] (-4.3719,-0.7400) -- (-4.2760,-0.8995) -- (-4.2373,-0.9858) -- (-4.3322,-0.8278) -- cycle;
\draw[draw=cbpInk!70,line width=.45pt] (-4.2373,-0.9858) -- (-4.3322,-0.8278);
\draw[draw=cbpInk!70,line width=.45pt] (-4.3322,-0.8278) -- (-4.3719,-0.7400);
\path[cbpmesh=38] (4.2760,-0.8995) -- (4.3719,-0.7400) -- (4.3322,-0.8278) -- (4.2373,-0.9858) -- cycle;
\draw[draw=cbpInk!70,line width=.45pt] (4.3719,-0.7400) -- (4.3322,-0.8278);
\draw[draw=cbpInk!70,line width=.45pt] (4.3322,-0.8278) -- (4.2373,-0.9858);
\path[cbpmesh=27.1] (-1.6963,0.7234) -- (-1.6591,0.6615) -- (-1.7441,0.6978) -- (-1.7832,0.7629) -- cycle;
\path[cbpmesh=18.6] (1.6591,0.6615) -- (1.6963,0.7234) -- (1.7832,0.7629) -- (1.7441,0.6978) -- cycle;
\path[cbpmesh=32.8] (-1.5097,0.5564) -- (-1.4666,0.5015) -- (-1.5356,0.5569) -- (-1.5808,0.6143) -- cycle;
\path[cbpmesh=21.4] (1.4666,0.5015) -- (1.5097,0.5564) -- (1.5808,0.6143) -- (1.5356,0.5569) -- cycle;
\path[cbpmesh=19.1] (-2.1121,0.8932) -- (-2.0793,0.8159) -- (-2.1879,0.8112) -- (-2.2224,0.8925) -- cycle;
\path[cbpmesh=17.9] (2.0793,0.8159) -- (2.1121,0.8932) -- (2.2224,0.8925) -- (2.1879,0.8112) -- cycle;
\path[cbpmesh=20.5] (-4.3528,0.0909) -- (-4.2853,-0.0684) -- (-4.3444,-0.1957) -- (-4.4129,-0.0342) -- cycle;
\path[cbpmesh=35.6] (4.2853,-0.0684) -- (4.3528,0.0909) -- (4.4129,-0.0342) -- (4.3444,-0.1957) -- cycle;
\path[cbpmesh=17.9] (-2.2224,0.8925) -- (-2.1879,0.8112) -- (-2.3005,0.7949) -- (-2.3368,0.8805) -- cycle;
\path[cbpmesh=18.5] (2.1879,0.8112) -- (2.2224,0.8925) -- (2.3368,0.8805) -- (2.3005,0.7949) -- cycle;
\path[cbpmesh=24.8] (-1.7832,0.7629) -- (-1.7441,0.6978) -- (-1.8353,0.7231) -- (-1.8764,0.7916) -- cycle;
\path[cbpmesh=18.0] (1.7441,0.6978) -- (1.7832,0.7629) -- (1.8764,0.7916) -- (1.8353,0.7231) -- cycle;
\path[cbpmesh=17.0] (-2.3368,0.8805) -- (-2.3005,0.7949) -- (-2.4165,0.7670) -- (-2.4546,0.8569) -- cycle;
\path[cbpmesh=19.3] (2.3005,0.7949) -- (2.3368,0.8805) -- (2.4546,0.8569) -- (2.4165,0.7670) -- cycle;
\path[cbpmesh=31.9] (-4.3981,-0.6447) -- (-4.3017,-0.8051) -- (-4.2760,-0.8995) -- (-4.3719,-0.7400) -- cycle;
\path[cbpmesh=38] (4.3017,-0.8051) -- (4.3981,-0.6447) -- (4.3719,-0.7400) -- (4.2760,-0.8995) -- cycle;
\path[cbpmesh=16.4] (-2.4546,0.8569) -- (-2.4165,0.7670) -- (-2.5354,0.7272) -- (-2.5753,0.8215) -- cycle;
\path[cbpmesh=20.3] (2.4165,0.7670) -- (2.4546,0.8569) -- (2.5753,0.8215) -- (2.5354,0.7272) -- cycle;
\path[cbpmesh=16.0] (-2.5753,0.8215) -- (-2.5354,0.7272) -- (-2.6565,0.6754) -- (-2.6984,0.7741) -- cycle;
\path[cbpmesh=21.4] (2.5354,0.7272) -- (2.5753,0.8215) -- (2.6984,0.7741) -- (2.6565,0.6754) -- cycle;
\path[cbpmesh=15.9] (-2.6984,0.7741) -- (-2.6565,0.6754) -- (-2.7708,0.6228) -- (-2.8145,0.7259) -- cycle;
\path[cbpmesh=21.8] (2.6565,0.6754) -- (2.6984,0.7741) -- (2.8145,0.7259) -- (2.7708,0.6228) -- cycle;
\path[cbpmesh=16.1] (-2.8145,0.7259) -- (-2.7708,0.6228) -- (-2.8778,0.5795) -- (-2.9231,0.6865) -- cycle;
\path[cbpmesh=21.1] (2.7708,0.6228) -- (2.8145,0.7259) -- (2.9231,0.6865) -- (2.8778,0.5795) -- cycle;
\path[cbpmesh=30.1] (-1.5808,0.6143) -- (-1.5356,0.5569) -- (-1.6117,0.6012) -- (-1.6591,0.6615) -- cycle;
\path[cbpmesh=20.1] (1.5356,0.5569) -- (1.5808,0.6143) -- (1.6591,0.6615) -- (1.6117,0.6012) -- cycle;
\path[cbpmesh=16.3] (-2.9231,0.6865) -- (-2.8778,0.5795) -- (-2.9781,0.5442) -- (-3.0250,0.6550) -- cycle;
\path[cbpmesh=20.5] (2.8778,0.5795) -- (2.9231,0.6865) -- (3.0250,0.6550) -- (2.9781,0.5442) -- cycle;
\path[cbpmesh=16.6] (-3.0250,0.6550) -- (-2.9781,0.5442) -- (-3.0722,0.5158) -- (-3.1206,0.6301) -- cycle;
\path[cbpmesh=20.0] (2.9781,0.5442) -- (3.0250,0.6550) -- (3.1206,0.6301) -- (3.0722,0.5158) -- cycle;
\path[cbpmesh=19.1] (-4.2718,0.2205) -- (-4.2055,0.0642) -- (-4.2853,-0.0684) -- (-4.3528,0.0909) -- cycle;
\path[cbpmesh=33.3] (4.2055,0.0642) -- (4.2718,0.2205) -- (4.3528,0.0909) -- (4.2853,-0.0684) -- cycle;
\path[cbpmesh=16.9] (-3.1206,0.6301) -- (-3.0722,0.5158) -- (-3.1607,0.4932) -- (-3.2105,0.6107) -- cycle;
\path[cbpmesh=19.5] (3.0722,0.5158) -- (3.1206,0.6301) -- (3.2105,0.6107) -- (3.1607,0.4932) -- cycle;
\path[cbpmesh=22.6] (-1.8764,0.7916) -- (-1.8353,0.7231) -- (-1.9320,0.7373) -- (-1.9753,0.8093) -- cycle;
\path[cbpmesh=17.7] (1.8353,0.7231) -- (1.8764,0.7916) -- (1.9753,0.8093) -- (1.9320,0.7373) -- cycle;
\path[cbpmesh=17.2] (-3.2105,0.6107) -- (-3.1607,0.4932) -- (-3.2442,0.4752) -- (-3.2953,0.5958) -- cycle;
\path[cbpmesh=19.1] (3.1607,0.4932) -- (3.2105,0.6107) -- (3.2953,0.5958) -- (3.2442,0.4752) -- cycle;
\path[cbpmesh=17.5] (-3.2953,0.5958) -- (-3.2442,0.4752) -- (-3.3232,0.4606) -- (-3.3755,0.5841) -- cycle;
\path[cbpmesh=18.8] (3.2442,0.4752) -- (3.2953,0.5958) -- (3.3755,0.5841) -- (3.3232,0.4606) -- cycle;
\path[cbpmesh=29.0] (-4.4097,-0.5422) -- (-4.3130,-0.7031) -- (-4.3017,-0.8051) -- (-4.3981,-0.6447) -- cycle;
\path[cbpmesh=38] (4.3130,-0.7031) -- (4.4097,-0.5422) -- (4.3981,-0.6447) -- (4.3017,-0.8051) -- cycle;
\path[cbpmesh=33.1] (-1.4666,0.5015) -- (-1.4144,0.4484) -- (-1.4810,0.5012) -- (-1.5356,0.5569) -- cycle;
\draw[draw=cbpInk!70,line width=.45pt] (-1.4666,0.5015) -- (-1.4144,0.4484);
\path[cbpmesh=22.1] (1.4144,0.4484) -- (1.4666,0.5015) -- (1.5356,0.5569) -- (1.4810,0.5012) -- cycle;
\draw[draw=cbpInk!70,line width=.45pt] (1.4144,0.4484) -- (1.4666,0.5015);
\path[cbpmesh=17.7] (-3.3755,0.5841) -- (-3.3232,0.4606) -- (-3.3983,0.4483) -- (-3.4519,0.5746) -- cycle;
\path[cbpmesh=18.7] (3.3232,0.4606) -- (3.3755,0.5841) -- (3.4519,0.5746) -- (3.3983,0.4483) -- cycle;
\path[cbpmesh=16.7] (-4.1975,0.2917) -- (-4.1324,0.1381) -- (-4.2055,0.0642) -- (-4.2718,0.2205) -- cycle;
\path[cbpmesh=28.3] (4.1324,0.1381) -- (4.1975,0.2917) -- (4.2718,0.2205) -- (4.2055,0.0642) -- cycle;
\path[cbpmesh=27.5] (-1.6591,0.6615) -- (-1.6117,0.6012) -- (-1.6943,0.6344) -- (-1.7441,0.6978) -- cycle;
\path[cbpmesh=19.1] (1.6117,0.6012) -- (1.6591,0.6615) -- (1.7441,0.6978) -- (1.6943,0.6344) -- cycle;
\path[cbpmesh=17.8] (-3.4519,0.5746) -- (-3.3983,0.4483) -- (-3.4702,0.4372) -- (-3.5249,0.5662) -- cycle;
\path[cbpmesh=18.6] (3.3983,0.4483) -- (3.4519,0.5746) -- (3.5249,0.5662) -- (3.4702,0.4372) -- cycle;
\path[cbpmesh=20.8] (-1.9753,0.8093) -- (-1.9320,0.7373) -- (-2.0338,0.7401) -- (-2.0793,0.8159) -- cycle;
\path[cbpmesh=17.7] (1.9320,0.7373) -- (1.9753,0.8093) -- (2.0793,0.8159) -- (2.0338,0.7401) -- cycle;
\path[cbpmesh=17.8] (-3.5249,0.5662) -- (-3.4702,0.4372) -- (-3.5393,0.4261) -- (-3.5951,0.5577) -- cycle;
\path[cbpmesh=18.6] (3.4702,0.4372) -- (3.5249,0.5662) -- (3.5951,0.5577) -- (3.5393,0.4261) -- cycle;
\path[cbpmesh=16.3] (-4.1262,0.3516) -- (-4.0622,0.2006) -- (-4.1324,0.1381) -- (-4.1975,0.2917) -- cycle;
\path[cbpmesh=27.0] (4.0622,0.2006) -- (4.1262,0.3516) -- (4.1975,0.2917) -- (4.1324,0.1381) -- cycle;
\path[cbpmesh=17.6] (-3.5951,0.5577) -- (-3.5393,0.4261) -- (-3.6062,0.4139) -- (-3.6631,0.5480) -- cycle;
\path[cbpmesh=18.8] (3.5393,0.4261) -- (3.5951,0.5577) -- (3.6631,0.5480) -- (3.6062,0.4139) -- cycle;
\path[cbpmesh=16.0] (-4.0573,0.4013) -- (-3.9944,0.2527) -- (-4.0622,0.2006) -- (-4.1262,0.3516) -- cycle;
\path[cbpmesh=25.6] (3.9944,0.2527) -- (4.0573,0.4013) -- (4.1262,0.3516) -- (4.0622,0.2006) -- cycle;
\path[cbpmesh=17.2] (-3.6631,0.5480) -- (-3.6062,0.4139) -- (-3.6716,0.3994) -- (-3.7294,0.5359) -- cycle;
\path[cbpmesh=19.2] (3.6062,0.4139) -- (3.6631,0.5480) -- (3.7294,0.5359) -- (3.6716,0.3994) -- cycle;
\path[cbpmesh=15.8] (-3.9904,0.4418) -- (-3.9285,0.2958) -- (-3.9944,0.2527) -- (-4.0573,0.4013) -- cycle;
\path[cbpmesh=24.2] (3.9285,0.2958) -- (3.9904,0.4418) -- (4.0573,0.4013) -- (3.9944,0.2527) -- cycle;
\path[cbpmesh=16.8] (-3.7294,0.5359) -- (-3.6716,0.3994) -- (-3.7359,0.3815) -- (-3.7948,0.5204) -- cycle;
\path[cbpmesh=19.7] (3.6716,0.3994) -- (3.7294,0.5359) -- (3.7948,0.5204) -- (3.7359,0.3815) -- cycle;
\path[cbpmesh=19.2] (-2.0793,0.8159) -- (-2.0338,0.7401) -- (-2.1400,0.7314) -- (-2.1879,0.8112) -- cycle;
\path[cbpmesh=18.0] (2.0338,0.7401) -- (2.0793,0.8159) -- (2.1879,0.8112) -- (2.1400,0.7314) -- cycle;
\path[cbpmesh=16.3] (-3.7948,0.5204) -- (-3.7359,0.3815) -- (-3.7998,0.3590) -- (-3.8597,0.5003) -- cycle;
\path[cbpmesh=20.5] (3.7359,0.3815) -- (3.7948,0.5204) -- (3.8597,0.5003) -- (3.7998,0.3590) -- cycle;
\path[cbpmesh=15.8] (-3.9247,0.4745) -- (-3.8638,0.3308) -- (-3.9285,0.2958) -- (-3.9904,0.4418) -- cycle;
\path[cbpmesh=22.8] (3.8638,0.3308) -- (3.9247,0.4745) -- (3.9904,0.4418) -- (3.9285,0.2958) -- cycle;
\path[cbpmesh=26.3] (-4.4054,-0.4330) -- (-4.3088,-0.5937) -- (-4.3130,-0.7031) -- (-4.4097,-0.5422) -- cycle;
\path[cbpmesh=38] (4.3088,-0.5937) -- (4.4054,-0.4330) -- (4.4097,-0.5422) -- (4.3130,-0.7031) -- cycle;
\path[cbpmesh=15.9] (-3.8597,0.5003) -- (-3.7998,0.3590) -- (-3.8638,0.3308) -- (-3.9247,0.4745) -- cycle;
\path[cbpmesh=21.6] (3.7998,0.3590) -- (3.8597,0.5003) -- (3.9247,0.4745) -- (3.8638,0.3308) -- cycle;
\path[cbpmesh=25.0] (-1.7441,0.6978) -- (-1.6943,0.6344) -- (-1.7828,0.6564) -- (-1.8353,0.7231) -- cycle;
\path[cbpmesh=18.4] (1.6943,0.6344) -- (1.7441,0.6978) -- (1.8353,0.7231) -- (1.7828,0.6564) -- cycle;
\path[cbpmesh=30.4] (-1.5356,0.5569) -- (-1.4810,0.5012) -- (-1.5543,0.5428) -- (-1.6117,0.6012) -- cycle;
\path[cbpmesh=20.8] (1.4810,0.5012) -- (1.5356,0.5569) -- (1.6117,0.6012) -- (1.5543,0.5428) -- cycle;
\path[cbpmesh=17.9] (-2.1879,0.8112) -- (-2.1400,0.7314) -- (-2.2501,0.7110) -- (-2.3005,0.7949) -- cycle;
\path[cbpmesh=18.5] (2.1400,0.7314) -- (2.1879,0.8112) -- (2.3005,0.7949) -- (2.2501,0.7110) -- cycle;
\path[cbpmesh=33.5] (-1.4144,0.4484) -- (-1.3535,0.3974) -- (-1.4172,0.4477) -- (-1.4810,0.5012) -- cycle;
\draw[draw=cbpInk!70,line width=.45pt] (-1.4144,0.4484) -- (-1.3535,0.3974);
\draw[draw=cbpInk!70,line width=.45pt] (-1.3535,0.3974) -- (-1.4172,0.4477);
\path[cbpmesh=22.9] (1.3535,0.3974) -- (1.4144,0.4484) -- (1.4810,0.5012) -- (1.4172,0.4477) -- cycle;
\draw[draw=cbpInk!70,line width=.45pt] (1.3535,0.3974) -- (1.4144,0.4484);
\draw[draw=cbpInk!70,line width=.45pt] (1.4172,0.4477) -- (1.3535,0.3974);
\path[cbpmesh=16.9] (-2.3005,0.7949) -- (-2.2501,0.7110) -- (-2.3636,0.6788) -- (-2.4165,0.7670) -- cycle;
\path[cbpmesh=19.2] (2.2501,0.7110) -- (2.3005,0.7949) -- (2.4165,0.7670) -- (2.3636,0.6788) -- cycle;
\path[cbpmesh=22.8] (-1.8353,0.7231) -- (-1.7828,0.6564) -- (-1.8768,0.6670) -- (-1.9320,0.7373) -- cycle;
\path[cbpmesh=18.0] (1.7828,0.6564) -- (1.8353,0.7231) -- (1.9320,0.7373) -- (1.8768,0.6670) -- cycle;
\path[cbpmesh=23.8] (-4.3840,-0.3174) -- (-4.2879,-0.4773) -- (-4.3088,-0.5937) -- (-4.4054,-0.4330) -- cycle;
\path[cbpmesh=38] (4.2879,-0.4773) -- (4.3840,-0.3174) -- (4.4054,-0.4330) -- (4.3088,-0.5937) -- cycle;
\path[cbpmesh=37.3] (-4.2373,-0.9858) -- (-4.1162,-1.1399) -- (-4.0669,-1.2160) -- (-4.1866,-1.0638) -- cycle;
\draw[draw=cbpInk!70,line width=.45pt] (-4.0669,-1.2160) -- (-4.1866,-1.0638);
\draw[draw=cbpInk!70,line width=.45pt] (-4.1866,-1.0638) -- (-4.2373,-0.9858);
\path[cbpmesh=38] (4.1162,-1.1399) -- (4.2373,-0.9858) -- (4.1866,-1.0638) -- (4.0669,-1.2160) -- cycle;
\draw[draw=cbpInk!70,line width=.45pt] (4.2373,-0.9858) -- (4.1866,-1.0638);
\draw[draw=cbpInk!70,line width=.45pt] (4.1866,-1.0638) -- (4.0669,-1.2160);
\path[cbpmesh=16.2] (-2.4165,0.7670) -- (-2.3636,0.6788) -- (-2.4798,0.6347) -- (-2.5354,0.7272) -- cycle;
\path[cbpmesh=20.1] (2.3636,0.6788) -- (2.4165,0.7670) -- (2.5354,0.7272) -- (2.4798,0.6347) -- cycle;
\path[cbpmesh=27.7] (-1.6117,0.6012) -- (-1.5543,0.5428) -- (-1.6340,0.5730) -- (-1.6943,0.6344) -- cycle;
\path[cbpmesh=19.6] (1.5543,0.5428) -- (1.6117,0.6012) -- (1.6943,0.6344) -- (1.6340,0.5730) -- cycle;
\path[cbpmesh=15.8] (-2.5354,0.7272) -- (-2.4798,0.6347) -- (-2.5983,0.5784) -- (-2.6565,0.6754) -- cycle;
\path[cbpmesh=21.1] (2.4798,0.6347) -- (2.5354,0.7272) -- (2.6565,0.6754) -- (2.5983,0.5784) -- cycle;
\path[cbpmesh=15.7] (-2.6565,0.6754) -- (-2.5983,0.5784) -- (-2.7101,0.5217) -- (-2.7708,0.6228) -- cycle;
\path[cbpmesh=21.4] (2.5983,0.5784) -- (2.6565,0.6754) -- (2.7708,0.6228) -- (2.7101,0.5217) -- cycle;
\path[cbpmesh=15.9] (-2.7708,0.6228) -- (-2.7101,0.5217) -- (-2.8147,0.4745) -- (-2.8778,0.5795) -- cycle;
\path[cbpmesh=20.8] (2.7101,0.5217) -- (2.7708,0.6228) -- (2.8778,0.5795) -- (2.8147,0.4745) -- cycle;
\path[cbpmesh=20.9] (-1.9320,0.7373) -- (-1.8768,0.6670) -- (-1.9756,0.6661) -- (-2.0338,0.7401) -- cycle;
\path[cbpmesh=17.9] (1.8768,0.6670) -- (1.9320,0.7373) -- (2.0338,0.7401) -- (1.9756,0.6661) -- cycle;
\path[cbpmesh=16.1] (-2.8778,0.5795) -- (-2.8147,0.4745) -- (-2.9128,0.4356) -- (-2.9781,0.5442) -- cycle;
\path[cbpmesh=20.3] (2.8147,0.4745) -- (2.8778,0.5795) -- (2.9781,0.5442) -- (2.9128,0.4356) -- cycle;
\path[cbpmesh=21.7] (-4.3444,-0.1957) -- (-4.2492,-0.3543) -- (-4.2879,-0.4773) -- (-4.3840,-0.3174) -- cycle;
\path[cbpmesh=37.2] (4.2492,-0.3543) -- (4.3444,-0.1957) -- (4.3840,-0.3174) -- (4.2879,-0.4773) -- cycle;
\path[cbpmesh=30.7] (-1.4810,0.5012) -- (-1.4172,0.4477) -- (-1.4874,0.4867) -- (-1.5543,0.5428) -- cycle;
\path[cbpmesh=21.4] (1.4172,0.4477) -- (1.4810,0.5012) -- (1.5543,0.5428) -- (1.4874,0.4867) -- cycle;
\path[cbpmesh=16.5] (-2.9781,0.5442) -- (-2.9128,0.4356) -- (-3.0048,0.4037) -- (-3.0722,0.5158) -- cycle;
\path[cbpmesh=19.8] (2.9128,0.4356) -- (2.9781,0.5442) -- (3.0722,0.5158) -- (3.0048,0.4037) -- cycle;
\path[cbpmesh=34.3] (-4.2760,-0.8995) -- (-4.1538,-1.0550) -- (-4.1162,-1.1399) -- (-4.2373,-0.9858) -- cycle;
\path[cbpmesh=38] (4.1538,-1.0550) -- (4.2760,-0.8995) -- (4.2373,-0.9858) -- (4.1162,-1.1399) -- cycle;
\path[cbpmesh=25.2] (-1.6943,0.6344) -- (-1.6340,0.5730) -- (-1.7194,0.5918) -- (-1.7828,0.6564) -- cycle;
\path[cbpmesh=18.8] (1.6340,0.5730) -- (1.6943,0.6344) -- (1.7828,0.6564) -- (1.7194,0.5918) -- cycle;
\path[cbpmesh=16.8] (-3.0722,0.5158) -- (-3.0048,0.4037) -- (-3.0914,0.3779) -- (-3.1607,0.4932) -- cycle;
\path[cbpmesh=19.3] (3.0048,0.4037) -- (3.0722,0.5158) -- (3.1607,0.4932) -- (3.0914,0.3779) -- cycle;
\path[cbpmesh=19.2] (-2.0338,0.7401) -- (-1.9756,0.6661) -- (-2.0788,0.6535) -- (-2.1400,0.7314) -- cycle;
\path[cbpmesh=18.0] (1.9756,0.6661) -- (2.0338,0.7401) -- (2.1400,0.7314) -- (2.0788,0.6535) -- cycle;
\path[cbpmesh=17.2] (-3.1607,0.4932) -- (-3.0914,0.3779) -- (-3.1730,0.3568) -- (-3.2442,0.4752) -- cycle;
\path[cbpmesh=19.0] (3.0914,0.3779) -- (3.1607,0.4932) -- (3.2442,0.4752) -- (3.1730,0.3568) -- cycle;
\path[cbpmesh=19.9] (-4.2853,-0.0684) -- (-4.1914,-0.2248) -- (-4.2492,-0.3543) -- (-4.3444,-0.1957) -- cycle;
\path[cbpmesh=34.7] (4.1914,-0.2248) -- (4.2853,-0.0684) -- (4.3444,-0.1957) -- (4.2492,-0.3543) -- cycle;
\path[cbpmesh=17.5] (-3.2442,0.4752) -- (-3.1730,0.3568) -- (-3.2504,0.3393) -- (-3.3232,0.4606) -- cycle;
\path[cbpmesh=18.8] (3.1730,0.3568) -- (3.2442,0.4752) -- (3.3232,0.4606) -- (3.2504,0.3393) -- cycle;
\path[cbpmesh=17.9] (-2.1400,0.7314) -- (-2.0788,0.6535) -- (-2.1858,0.6292) -- (-2.2501,0.7110) -- cycle;
\path[cbpmesh=18.5] (2.0788,0.6535) -- (2.1400,0.7314) -- (2.2501,0.7110) -- (2.1858,0.6292) -- cycle;
\path[cbpmesh=28.0] (-1.5543,0.5428) -- (-1.4874,0.4867) -- (-1.5636,0.5140) -- (-1.6340,0.5730) -- cycle;
\path[cbpmesh=20.2] (1.4874,0.4867) -- (1.5543,0.5428) -- (1.6340,0.5730) -- (1.5636,0.5140) -- cycle;
\path[cbpmesh=17.7] (-3.3232,0.4606) -- (-3.2504,0.3393) -- (-3.3239,0.3243) -- (-3.3983,0.4483) -- cycle;
\path[cbpmesh=18.6] (3.2504,0.3393) -- (3.3232,0.4606) -- (3.3983,0.4483) -- (3.3239,0.3243) -- cycle;
\path[cbpmesh=31.3] (-4.3017,-0.8051) -- (-4.1787,-0.9616) -- (-4.1538,-1.0550) -- (-4.2760,-0.8995) -- cycle;
\path[cbpmesh=38] (4.1787,-0.9616) -- (4.3017,-0.8051) -- (4.2760,-0.8995) -- (4.1538,-1.0550) -- cycle;
\path[cbpmesh=23.0] (-1.7828,0.6564) -- (-1.7194,0.5918) -- (-1.8100,0.5990) -- (-1.8768,0.6670) -- cycle;
\path[cbpmesh=18.3] (1.7194,0.5918) -- (1.7828,0.6564) -- (1.8768,0.6670) -- (1.8100,0.5990) -- cycle;
\path[cbpmesh=17.8] (-3.3983,0.4483) -- (-3.3239,0.3243) -- (-3.3941,0.3106) -- (-3.4702,0.4372) -- cycle;
\path[cbpmesh=18.6] (3.3239,0.3243) -- (3.3983,0.4483) -- (3.4702,0.4372) -- (3.3941,0.3106) -- cycle;
\path[cbpmesh=16.9] (-2.2501,0.7110) -- (-2.1858,0.6292) -- (-2.2960,0.5929) -- (-2.3636,0.6788) -- cycle;
\path[cbpmesh=19.1] (2.1858,0.6292) -- (2.2501,0.7110) -- (2.3636,0.6788) -- (2.2960,0.5929) -- cycle;
\path[cbpmesh=18.5] (-4.2055,0.0642) -- (-4.1133,-0.0893) -- (-4.1914,-0.2248) -- (-4.2853,-0.0684) -- cycle;
\path[cbpmesh=32.5] (4.1133,-0.0893) -- (4.2055,0.0642) -- (4.2853,-0.0684) -- (4.1914,-0.2248) -- cycle;
\path[cbpmesh=17.7] (-3.4702,0.4372) -- (-3.3941,0.3106) -- (-3.4617,0.2970) -- (-3.5393,0.4261) -- cycle;
\path[cbpmesh=18.6] (3.3941,0.3106) -- (3.4702,0.4372) -- (3.5393,0.4261) -- (3.4617,0.2970) -- cycle;
\path[cbpmesh=30.9] (-1.4172,0.4477) -- (-1.3447,0.3968) -- (-1.4113,0.4332) -- (-1.4874,0.4867) -- cycle;
\draw[draw=cbpInk!70,line width=.45pt] (-1.4172,0.4477) -- (-1.3447,0.3968);
\path[cbpmesh=22.1] (1.3447,0.3968) -- (1.4172,0.4477) -- (1.4874,0.4867) -- (1.4113,0.4332) -- cycle;
\draw[draw=cbpInk!70,line width=.45pt] (1.3447,0.3968) -- (1.4172,0.4477);
\path[cbpmesh=17.5] (-3.5393,0.4261) -- (-3.4617,0.2970) -- (-3.5272,0.2823) -- (-3.6062,0.4139) -- cycle;
\path[cbpmesh=18.7] (3.4617,0.2970) -- (3.5393,0.4261) -- (3.6062,0.4139) -- (3.5272,0.2823) -- cycle;
\path[cbpmesh=16.1] (-2.3636,0.6788) -- (-2.2960,0.5929) -- (-2.4089,0.5445) -- (-2.4798,0.6347) -- cycle;
\path[cbpmesh=19.8] (2.2960,0.5929) -- (2.3636,0.6788) -- (2.4798,0.6347) -- (2.4089,0.5445) -- cycle;
\path[cbpmesh=21.0] (-1.8768,0.6670) -- (-1.8100,0.5990) -- (-1.9053,0.5945) -- (-1.9756,0.6661) -- cycle;
\path[cbpmesh=18.1] (1.8100,0.5990) -- (1.8768,0.6670) -- (1.9756,0.6661) -- (1.9053,0.5945) -- cycle;
\path[cbpmesh=16.2] (-4.1324,0.1381) -- (-4.0418,-0.0127) -- (-4.1133,-0.0893) -- (-4.2055,0.0642) -- cycle;
\path[cbpmesh=27.6] (4.0418,-0.0127) -- (4.1324,0.1381) -- (4.2055,0.0642) -- (4.1133,-0.0893) -- cycle;
\path[cbpmesh=17.1] (-3.6062,0.4139) -- (-3.5272,0.2823) -- (-3.5911,0.2654) -- (-3.6716,0.3994) -- cycle;
\path[cbpmesh=19.1] (3.5272,0.2823) -- (3.6062,0.4139) -- (3.6716,0.3994) -- (3.5911,0.2654) -- cycle;
\path[cbpmesh=28.4] (-4.3130,-0.7031) -- (-4.1898,-0.8600) -- (-4.1787,-0.9616) -- (-4.3017,-0.8051) -- cycle;
\path[cbpmesh=38] (4.1898,-0.8600) -- (4.3130,-0.7031) -- (4.3017,-0.8051) -- (4.1787,-0.9616) -- cycle;
\path[cbpmesh=25.4] (-1.6340,0.5730) -- (-1.5636,0.5140) -- (-1.6454,0.5297) -- (-1.7194,0.5918) -- cycle;
\path[cbpmesh=19.3] (1.5636,0.5140) -- (1.6340,0.5730) -- (1.7194,0.5918) -- (1.6454,0.5297) -- cycle;
\path[cbpmesh=16.6] (-3.6716,0.3994) -- (-3.5911,0.2654) -- (-3.6540,0.2452) -- (-3.7359,0.3815) -- cycle;
\path[cbpmesh=19.5] (3.5911,0.2654) -- (3.6716,0.3994) -- (3.7359,0.3815) -- (3.6540,0.2452) -- cycle;
\path[cbpmesh=15.6] (-2.4798,0.6347) -- (-2.4089,0.5445) -- (-2.5240,0.4839) -- (-2.5983,0.5784) -- cycle;
\path[cbpmesh=20.7] (2.4089,0.5445) -- (2.4798,0.6347) -- (2.5983,0.5784) -- (2.5240,0.4839) -- cycle;
\path[cbpmesh=15.9] (-4.0622,0.2006) -- (-3.9731,0.0523) -- (-4.0418,-0.0127) -- (-4.1324,0.1381) -- cycle;
\path[cbpmesh=26.4] (3.9731,0.0523) -- (4.0622,0.2006) -- (4.1324,0.1381) -- (4.0418,-0.0127) -- cycle;
\path[cbpmesh=16.1] (-3.7359,0.3815) -- (-3.6540,0.2452) -- (-3.7165,0.2204) -- (-3.7998,0.3590) -- cycle;
\path[cbpmesh=20.3] (3.6540,0.2452) -- (3.7359,0.3815) -- (3.7998,0.3590) -- (3.7165,0.2204) -- cycle;
\path[cbpmesh=15.6] (-3.9944,0.2527) -- (-3.9069,0.1070) -- (-3.9731,0.0523) -- (-4.0622,0.2006) -- cycle;
\path[cbpmesh=25.1] (3.9069,0.1070) -- (3.9944,0.2527) -- (4.0622,0.2006) -- (3.9731,0.0523) -- cycle;
\path[cbpmesh=15.7] (-3.7998,0.3590) -- (-3.7165,0.2204) -- (-3.7791,0.1898) -- (-3.8638,0.3308) -- cycle;
\path[cbpmesh=21.2] (3.7165,0.2204) -- (3.7998,0.3590) -- (3.8638,0.3308) -- (3.7791,0.1898) -- cycle;
\path[cbpmesh=15.4] (-3.9285,0.2958) -- (-3.8424,0.1524) -- (-3.9069,0.1070) -- (-3.9944,0.2527) -- cycle;
\path[cbpmesh=23.7] (3.8424,0.1524) -- (3.9285,0.2958) -- (3.9944,0.2527) -- (3.9069,0.1070) -- cycle;
\path[cbpmesh=15.5] (-3.8638,0.3308) -- (-3.7791,0.1898) -- (-3.8424,0.1524) -- (-3.9285,0.2958) -- cycle;
\path[cbpmesh=22.4] (3.7791,0.1898) -- (3.8638,0.3308) -- (3.9285,0.2958) -- (3.8424,0.1524) -- cycle;
\path[cbpmesh=15.4] (-2.5983,0.5784) -- (-2.5240,0.4839) -- (-2.6326,0.4232) -- (-2.7101,0.5217) -- cycle;
\path[cbpmesh=21.1] (2.5240,0.4839) -- (2.5983,0.5784) -- (2.7101,0.5217) -- (2.6326,0.4232) -- cycle;
\path[cbpmesh=38] (-4.0669,-1.2160) -- (-3.9222,-1.3634) -- (-3.8646,-1.4283) -- (-4.0072,-1.2831) -- cycle;
\draw[draw=cbpInk!70,line width=.45pt] (-3.8646,-1.4283) -- (-4.0072,-1.2831);
\draw[draw=cbpInk!70,line width=.45pt] (-4.0072,-1.2831) -- (-4.0669,-1.2160);
\path[cbpmesh=38] (3.9222,-1.3634) -- (4.0669,-1.2160) -- (4.0072,-1.2831) -- (3.8646,-1.4283) -- cycle;
\draw[draw=cbpInk!70,line width=.45pt] (4.0669,-1.2160) -- (4.0072,-1.2831);
\draw[draw=cbpInk!70,line width=.45pt] (4.0072,-1.2831) -- (3.8646,-1.4283);
\path[cbpmesh=19.3] (-1.9756,0.6661) -- (-1.9053,0.5945) -- (-2.0049,0.5782) -- (-2.0788,0.6535) -- cycle;
\path[cbpmesh=18.1] (1.9053,0.5945) -- (1.9756,0.6661) -- (2.0788,0.6535) -- (2.0049,0.5782) -- cycle;
\path[cbpmesh=15.7] (-2.7101,0.5217) -- (-2.6326,0.4232) -- (-2.7343,0.3722) -- (-2.8147,0.4745) -- cycle;
\path[cbpmesh=20.5] (2.6326,0.4232) -- (2.7101,0.5217) -- (2.8147,0.4745) -- (2.7343,0.3722) -- cycle;
\path[cbpmesh=28.2] (-1.4874,0.4867) -- (-1.4113,0.4332) -- (-1.4836,0.4578) -- (-1.5636,0.5140) -- cycle;
\path[cbpmesh=20.8] (1.4113,0.4332) -- (1.4874,0.4867) -- (1.5636,0.5140) -- (1.4836,0.4578) -- cycle;
\path[cbpmesh=25.7] (-4.3088,-0.5937) -- (-4.1857,-0.7504) -- (-4.1898,-0.8600) -- (-4.3130,-0.7031) -- cycle;
\path[cbpmesh=38] (4.1857,-0.7504) -- (4.3088,-0.5937) -- (4.3130,-0.7031) -- (4.1898,-0.8600) -- cycle;
\path[cbpmesh=16.0] (-2.8147,0.4745) -- (-2.7343,0.3722) -- (-2.8296,0.3297) -- (-2.9128,0.4356) -- cycle;
\path[cbpmesh=20.0] (2.7343,0.3722) -- (2.8147,0.4745) -- (2.9128,0.4356) -- (2.8296,0.3297) -- cycle;
\path[cbpmesh=23.1] (-1.7194,0.5918) -- (-1.6454,0.5297) -- (-1.7321,0.5337) -- (-1.8100,0.5990) -- cycle;
\path[cbpmesh=18.6] (1.6454,0.5297) -- (1.7194,0.5918) -- (1.8100,0.5990) -- (1.7321,0.5337) -- cycle;
\path[cbpmesh=31.1] (-1.3447,0.3968) -- (-1.2639,0.3487) -- (-1.3265,0.3827) -- (-1.4113,0.4332) -- cycle;
\draw[draw=cbpInk!70,line width=.45pt] (-1.3447,0.3968) -- (-1.2639,0.3487);
\draw[draw=cbpInk!70,line width=.45pt] (-1.2639,0.3487) -- (-1.3265,0.3827);
\path[cbpmesh=22.8] (1.2639,0.3487) -- (1.3447,0.3968) -- (1.4113,0.4332) -- (1.3265,0.3827) -- cycle;
\draw[draw=cbpInk!70,line width=.45pt] (1.2639,0.3487) -- (1.3447,0.3968);
\draw[draw=cbpInk!70,line width=.45pt] (1.3265,0.3827) -- (1.2639,0.3487);
\path[cbpmesh=16.3] (-2.9128,0.4356) -- (-2.8296,0.3297) -- (-2.9190,0.2945) -- (-3.0048,0.4037) -- cycle;
\path[cbpmesh=19.6] (2.8296,0.3297) -- (2.9128,0.4356) -- (3.0048,0.4037) -- (2.9190,0.2945) -- cycle;
\path[cbpmesh=17.9] (-2.0788,0.6535) -- (-2.0049,0.5782) -- (-2.1080,0.5500) -- (-2.1858,0.6292) -- cycle;
\path[cbpmesh=18.4] (2.0049,0.5782) -- (2.0788,0.6535) -- (2.1858,0.6292) -- (2.1080,0.5500) -- cycle;
\path[cbpmesh=36.8] (-4.1162,-1.1399) -- (-3.9697,-1.2891) -- (-3.9222,-1.3634) -- (-4.0669,-1.2160) -- cycle;
\path[cbpmesh=38] (3.9697,-1.2891) -- (4.1162,-1.1399) -- (4.0669,-1.2160) -- (3.9222,-1.3634) -- cycle;
\path[cbpmesh=16.7] (-3.0048,0.4037) -- (-2.9190,0.2945) -- (-3.0031,0.2655) -- (-3.0914,0.3779) -- cycle;
\path[cbpmesh=19.2] (2.9190,0.2945) -- (3.0048,0.4037) -- (3.0914,0.3779) -- (3.0031,0.2655) -- cycle;
\path[cbpmesh=23.2] (-4.2879,-0.4773) -- (-4.1654,-0.6333) -- (-4.1857,-0.7504) -- (-4.3088,-0.5937) -- cycle;
\path[cbpmesh=38] (4.1654,-0.6333) -- (4.2879,-0.4773) -- (4.3088,-0.5937) -- (4.1857,-0.7504) -- cycle;
\path[cbpmesh=25.6] (-1.5636,0.5140) -- (-1.4836,0.4578) -- (-1.5612,0.4706) -- (-1.6454,0.5297) -- cycle;
\path[cbpmesh=19.7] (1.4836,0.4578) -- (1.5636,0.5140) -- (1.6454,0.5297) -- (1.5612,0.4706) -- cycle;
\path[cbpmesh=21.1] (-1.8100,0.5990) -- (-1.7321,0.5337) -- (-1.8233,0.5257) -- (-1.9053,0.5945) -- cycle;
\path[cbpmesh=18.3] (1.7321,0.5337) -- (1.8100,0.5990) -- (1.9053,0.5945) -- (1.8233,0.5257) -- cycle;
\path[cbpmesh=16.8] (-2.1858,0.6292) -- (-2.1080,0.5500) -- (-2.2143,0.5097) -- (-2.2960,0.5929) -- cycle;
\path[cbpmesh=18.9] (2.1080,0.5500) -- (2.1858,0.6292) -- (2.2960,0.5929) -- (2.2143,0.5097) -- cycle;
\path[cbpmesh=17.1] (-3.0914,0.3779) -- (-3.0031,0.2655) -- (-3.0824,0.2414) -- (-3.1730,0.3568) -- cycle;
\path[cbpmesh=18.9] (3.0031,0.2655) -- (3.0914,0.3779) -- (3.1730,0.3568) -- (3.0824,0.2414) -- cycle;
\path[cbpmesh=28.3] (-1.4113,0.4332) -- (-1.3265,0.3827) -- (-1.3945,0.4047) -- (-1.4836,0.4578) -- cycle;
\path[cbpmesh=21.3] (1.3265,0.3827) -- (1.4113,0.4332) -- (1.4836,0.4578) -- (1.3945,0.4047) -- cycle;
\path[cbpmesh=17.4] (-3.1730,0.3568) -- (-3.0824,0.2414) -- (-3.1575,0.2211) -- (-3.2504,0.3393) -- cycle;
\path[cbpmesh=18.7] (3.0824,0.2414) -- (3.1730,0.3568) -- (3.2504,0.3393) -- (3.1575,0.2211) -- cycle;
\path[cbpmesh=16.0] (-2.2960,0.5929) -- (-2.2143,0.5097) -- (-2.3232,0.4572) -- (-2.4089,0.5445) -- cycle;
\path[cbpmesh=19.6] (2.2143,0.5097) -- (2.2960,0.5929) -- (2.4089,0.5445) -- (2.3232,0.4572) -- cycle;
\path[cbpmesh=21.1] (-4.2492,-0.3543) -- (-4.1278,-0.5088) -- (-4.1654,-0.6333) -- (-4.2879,-0.4773) -- cycle;
\path[cbpmesh=36.2] (4.1278,-0.5088) -- (4.2492,-0.3543) -- (4.2879,-0.4773) -- (4.1654,-0.6333) -- cycle;
\path[cbpmesh=33.8] (-4.1538,-1.0550) -- (-4.0060,-1.2055) -- (-3.9697,-1.2891) -- (-4.1162,-1.1399) -- cycle;
\path[cbpmesh=38] (4.0060,-1.2055) -- (4.1538,-1.0550) -- (4.1162,-1.1399) -- (3.9697,-1.2891) -- cycle;
\path[cbpmesh=17.7] (-3.2504,0.3393) -- (-3.1575,0.2211) -- (-3.2289,0.2034) -- (-3.3239,0.3243) -- cycle;
\path[cbpmesh=18.6] (3.1575,0.2211) -- (3.2504,0.3393) -- (3.3239,0.3243) -- (3.2289,0.2034) -- cycle;
\path[cbpmesh=19.3] (-1.9053,0.5945) -- (-1.8233,0.5257) -- (-1.9185,0.5058) -- (-2.0049,0.5782) -- cycle;
\path[cbpmesh=18.2] (1.8233,0.5257) -- (1.9053,0.5945) -- (2.0049,0.5782) -- (1.9185,0.5058) -- cycle;
\path[cbpmesh=23.2] (-1.6454,0.5297) -- (-1.5612,0.4706) -- (-1.6435,0.4714) -- (-1.7321,0.5337) -- cycle;
\path[cbpmesh=18.9] (1.5612,0.4706) -- (1.6454,0.5297) -- (1.7321,0.5337) -- (1.6435,0.4714) -- cycle;
\path[cbpmesh=15.4] (-2.4089,0.5445) -- (-2.3232,0.4572) -- (-2.4342,0.3925) -- (-2.5240,0.4839) -- cycle;
\path[cbpmesh=20.4] (2.3232,0.4572) -- (2.4089,0.5445) -- (2.5240,0.4839) -- (2.4342,0.3925) -- cycle;
\path[cbpmesh=17.8] (-3.3239,0.3243) -- (-3.2289,0.2034) -- (-3.2971,0.1871) -- (-3.3941,0.3106) -- cycle;
\path[cbpmesh=18.5] (3.2289,0.2034) -- (3.3239,0.3243) -- (3.3941,0.3106) -- (3.2971,0.1871) -- cycle;
\path[cbpmesh=15.2] (-2.5240,0.4839) -- (-2.4342,0.3925) -- (-2.5390,0.3278) -- (-2.6326,0.4232) -- cycle;
\path[cbpmesh=20.7] (2.4342,0.3925) -- (2.5240,0.4839) -- (2.6326,0.4232) -- (2.5390,0.3278) -- cycle;
\path[cbpmesh=17.7] (-3.3941,0.3106) -- (-3.2971,0.1871) -- (-3.3628,0.1711) -- (-3.4617,0.2970) -- cycle;
\path[cbpmesh=18.5] (3.2971,0.1871) -- (3.3941,0.3106) -- (3.4617,0.2970) -- (3.3628,0.1711) -- cycle;
\path[cbpmesh=19.4] (-4.1914,-0.2248) -- (-4.0716,-0.3772) -- (-4.1278,-0.5088) -- (-4.2492,-0.3543) -- cycle;
\path[cbpmesh=33.8] (4.0716,-0.3772) -- (4.1914,-0.2248) -- (4.2492,-0.3543) -- (4.1278,-0.5088) -- cycle;
\path[cbpmesh=25.7] (-1.4836,0.4578) -- (-1.3945,0.4047) -- (-1.4674,0.4147) -- (-1.5612,0.4706) -- cycle;
\path[cbpmesh=20.2] (1.3945,0.4047) -- (1.4836,0.4578) -- (1.5612,0.4706) -- (1.4674,0.4147) -- cycle;
\path[cbpmesh=17.9] (-2.0049,0.5782) -- (-1.9185,0.5058) -- (-2.0173,0.4739) -- (-2.1080,0.5500) -- cycle;
\path[cbpmesh=18.4] (1.9185,0.5058) -- (2.0049,0.5782) -- (2.1080,0.5500) -- (2.0173,0.4739) -- cycle;
\path[cbpmesh=17.5] (-3.4617,0.2970) -- (-3.3628,0.1711) -- (-3.4264,0.1540) -- (-3.5272,0.2823) -- cycle;
\path[cbpmesh=18.7] (3.3628,0.1711) -- (3.4617,0.2970) -- (3.5272,0.2823) -- (3.4264,0.1540) -- cycle;
\path[cbpmesh=30.8] (-4.1787,-0.9616) -- (-4.0301,-1.1130) -- (-4.0060,-1.2055) -- (-4.1538,-1.0550) -- cycle;
\path[cbpmesh=38] (4.0301,-1.1130) -- (4.1787,-0.9616) -- (4.1538,-1.0550) -- (4.0060,-1.2055) -- cycle;
\path[cbpmesh=15.5] (-2.6326,0.4232) -- (-2.5390,0.3278) -- (-2.6370,0.2731) -- (-2.7343,0.3722) -- cycle;
\path[cbpmesh=20.2] (2.5390,0.3278) -- (2.6326,0.4232) -- (2.7343,0.3722) -- (2.6370,0.2731) -- cycle;
\path[cbpmesh=28.5] (-1.3265,0.3827) -- (-1.2335,0.3355) -- (-1.2967,0.3551) -- (-1.3945,0.4047) -- cycle;
\draw[draw=cbpInk!70,line width=.45pt] (-1.3265,0.3827) -- (-1.2335,0.3355);
\path[cbpmesh=21.9] (1.2335,0.3355) -- (1.3265,0.3827) -- (1.3945,0.4047) -- (1.2967,0.3551) -- cycle;
\draw[draw=cbpInk!70,line width=.45pt] (1.2335,0.3355) -- (1.3265,0.3827);
\path[cbpmesh=21.1] (-1.7321,0.5337) -- (-1.6435,0.4714) -- (-1.7300,0.4602) -- (-1.8233,0.5257) -- cycle;
\path[cbpmesh=18.5] (1.6435,0.4714) -- (1.7321,0.5337) -- (1.8233,0.5257) -- (1.7300,0.4602) -- cycle;
\path[cbpmesh=17.0] (-3.5272,0.2823) -- (-3.4264,0.1540) -- (-3.4885,0.1348) -- (-3.5911,0.2654) -- cycle;
\path[cbpmesh=18.9] (3.4264,0.1540) -- (3.5272,0.2823) -- (3.5911,0.2654) -- (3.4885,0.1348) -- cycle;
\path[cbpmesh=15.8] (-2.7343,0.3722) -- (-2.6370,0.2731) -- (-2.7289,0.2271) -- (-2.8296,0.3297) -- cycle;
\path[cbpmesh=19.7] (2.6370,0.2731) -- (2.7343,0.3722) -- (2.8296,0.3297) -- (2.7289,0.2271) -- cycle;
\path[cbpmesh=18.0] (-4.1133,-0.0893) -- (-3.9958,-0.2389) -- (-4.0716,-0.3772) -- (-4.1914,-0.2248) -- cycle;
\path[cbpmesh=31.6] (3.9958,-0.2389) -- (4.1133,-0.0893) -- (4.1914,-0.2248) -- (4.0716,-0.3772) -- cycle;
\path[cbpmesh=16.5] (-3.5911,0.2654) -- (-3.4885,0.1348) -- (-3.5496,0.1123) -- (-3.6540,0.2452) -- cycle;
\path[cbpmesh=19.4] (3.4885,0.1348) -- (3.5911,0.2654) -- (3.6540,0.2452) -- (3.5496,0.1123) -- cycle;
\path[cbpmesh=16.7] (-2.1080,0.5500) -- (-2.0173,0.4739) -- (-2.1190,0.4298) -- (-2.2143,0.5097) -- cycle;
\path[cbpmesh=18.8] (2.0173,0.4739) -- (2.1080,0.5500) -- (2.2143,0.5097) -- (2.1190,0.4298) -- cycle;
\path[cbpmesh=16.0] (-3.6540,0.2452) -- (-3.5496,0.1123) -- (-3.6103,0.0852) -- (-3.7165,0.2204) -- cycle;
\path[cbpmesh=20.0] (3.5496,0.1123) -- (3.6540,0.2452) -- (3.7165,0.2204) -- (3.6103,0.0852) -- cycle;
\path[cbpmesh=16.2] (-2.8296,0.3297) -- (-2.7289,0.2271) -- (-2.8151,0.1887) -- (-2.9190,0.2945) -- cycle;
\path[cbpmesh=19.3] (2.7289,0.2271) -- (2.8296,0.3297) -- (2.9190,0.2945) -- (2.8151,0.1887) -- cycle;
\path[cbpmesh=15.8] (-4.0418,-0.0127) -- (-3.9263,-0.1597) -- (-3.9958,-0.2389) -- (-4.1133,-0.0893) -- cycle;
\path[cbpmesh=26.9] (3.9263,-0.1597) -- (4.0418,-0.0127) -- (4.1133,-0.0893) -- (3.9958,-0.2389) -- cycle;
\path[cbpmesh=15.5] (-3.7165,0.2204) -- (-3.6103,0.0852) -- (-3.6711,0.0524) -- (-3.7791,0.1898) -- cycle;
\path[cbpmesh=20.9] (3.6103,0.0852) -- (3.7165,0.2204) -- (3.7791,0.1898) -- (3.6711,0.0524) -- cycle;
\path[cbpmesh=27.9] (-4.1898,-0.8600) -- (-4.0407,-1.0118) -- (-4.0301,-1.1130) -- (-4.1787,-0.9616) -- cycle;
\path[cbpmesh=38] (4.0407,-1.0118) -- (4.1898,-0.8600) -- (4.1787,-0.9616) -- (4.0301,-1.1130) -- cycle;
\path[cbpmesh=15.5] (-3.9731,0.0523) -- (-3.8596,-0.0921) -- (-3.9263,-0.1597) -- (-4.0418,-0.0127) -- cycle;
\path[cbpmesh=25.7] (3.8596,-0.0921) -- (3.9731,0.0523) -- (4.0418,-0.0127) -- (3.9263,-0.1597) -- cycle;
\path[cbpmesh=23.3] (-1.5612,0.4706) -- (-1.4674,0.4147) -- (-1.5447,0.4126) -- (-1.6435,0.4714) -- cycle;
\path[cbpmesh=19.3] (1.4674,0.4147) -- (1.5612,0.4706) -- (1.6435,0.4714) -- (1.5447,0.4126) -- cycle;
\path[cbpmesh=15.2] (-3.7791,0.1898) -- (-3.6711,0.0524) -- (-3.7326,0.0127) -- (-3.8424,0.1524) -- cycle;
\path[cbpmesh=22.0] (3.6711,0.0524) -- (3.7791,0.1898) -- (3.8424,0.1524) -- (3.7326,0.0127) -- cycle;
\path[cbpmesh=15.2] (-3.9069,0.1070) -- (-3.7952,-0.0351) -- (-3.8596,-0.0921) -- (-3.9731,0.0523) -- cycle;
\path[cbpmesh=24.5] (3.7952,-0.0351) -- (3.9069,0.1070) -- (3.9731,0.0523) -- (3.8596,-0.0921) -- cycle;
\path[cbpmesh=15.1] (-3.8424,0.1524) -- (-3.7326,0.0127) -- (-3.7952,-0.0351) -- (-3.9069,0.1070) -- cycle;
\path[cbpmesh=23.2] (3.7326,0.0127) -- (3.8424,0.1524) -- (3.9069,0.1070) -- (3.7952,-0.0351) -- cycle;
\path[cbpmesh=19.3] (-1.8233,0.5257) -- (-1.7300,0.4602) -- (-1.8204,0.4369) -- (-1.9185,0.5058) -- cycle;
\path[cbpmesh=18.3] (1.7300,0.4602) -- (1.8233,0.5257) -- (1.9185,0.5058) -- (1.8204,0.4369) -- cycle;
\path[cbpmesh=38] (-3.9222,-1.3634) -- (-3.7533,-1.5050) -- (-3.6982,-1.5678) -- (-3.8646,-1.4283) -- cycle;
\draw[draw=cbpInk!70,line width=.45pt] (-3.6982,-1.5678) -- (-3.8646,-1.4283);
\path[cbpmesh=38] (3.7533,-1.5050) -- (3.9222,-1.3634) -- (3.8646,-1.4283) -- (3.6982,-1.5678) -- cycle;
\draw[draw=cbpInk!70,line width=.45pt] (3.8646,-1.4283) -- (3.6982,-1.5678);
\path[cbpmesh=16.6] (-2.9190,0.2945) -- (-2.8151,0.1887) -- (-2.8962,0.1566) -- (-3.0031,0.2655) -- cycle;
\path[cbpmesh=19.0] (2.8151,0.1887) -- (2.9190,0.2945) -- (3.0031,0.2655) -- (2.8962,0.1566) -- cycle;
\path[cbpmesh=15.9] (-2.2143,0.5097) -- (-2.1190,0.4298) -- (-2.2232,0.3734) -- (-2.3232,0.4572) -- cycle;
\path[cbpmesh=19.3] (2.1190,0.4298) -- (2.2143,0.5097) -- (2.3232,0.4572) -- (2.2232,0.3734) -- cycle;
\path[cbpmesh=25.8] (-1.3945,0.4047) -- (-1.2967,0.3551) -- (-1.3645,0.3625) -- (-1.4674,0.4147) -- cycle;
\path[cbpmesh=20.6] (1.2967,0.3551) -- (1.3945,0.4047) -- (1.4674,0.4147) -- (1.3645,0.3625) -- cycle;
\path[cbpmesh=28.5] (-1.2335,0.3355) -- (-1.1329,0.2918) -- (-1.1910,0.3092) -- (-1.2967,0.3551) -- cycle;
\draw[draw=cbpInk!70,line width=.45pt] (-1.2335,0.3355) -- (-1.1329,0.2918);
\path[cbpmesh=22.5] (1.1329,0.2918) -- (1.2335,0.3355) -- (1.2967,0.3551) -- (1.1910,0.3092) -- cycle;
\draw[draw=cbpInk!70,line width=.45pt] (1.1329,0.2918) -- (1.2335,0.3355);
\path[cbpmesh=17.0] (-3.0031,0.2655) -- (-2.8962,0.1566) -- (-2.9727,0.1297) -- (-3.0824,0.2414) -- cycle;
\path[cbpmesh=18.8] (2.8962,0.1566) -- (3.0031,0.2655) -- (3.0824,0.2414) -- (2.9727,0.1297) -- cycle;
\path[cbpmesh=25.1] (-4.1857,-0.7504) -- (-4.0368,-0.9021) -- (-4.0407,-1.0118) -- (-4.1898,-0.8600) -- cycle;
\path[cbpmesh=38] (4.0368,-0.9021) -- (4.1857,-0.7504) -- (4.1898,-0.8600) -- (4.0407,-1.0118) -- cycle;
\path[cbpmesh=15.2] (-2.3232,0.4572) -- (-2.2232,0.3734) -- (-2.3294,0.3046) -- (-2.4342,0.3925) -- cycle;
\path[cbpmesh=20.1] (2.2232,0.3734) -- (2.3232,0.4572) -- (2.4342,0.3925) -- (2.3294,0.3046) -- cycle;
\path[cbpmesh=21.2] (-1.6435,0.4714) -- (-1.5447,0.4126) -- (-1.6261,0.3983) -- (-1.7300,0.4602) -- cycle;
\path[cbpmesh=18.7] (1.5447,0.4126) -- (1.6435,0.4714) -- (1.7300,0.4602) -- (1.6261,0.3983) -- cycle;
\path[cbpmesh=17.9] (-1.9185,0.5058) -- (-1.8204,0.4369) -- (-1.9140,0.4014) -- (-2.0173,0.4739) -- cycle;
\path[cbpmesh=18.3] (1.8204,0.4369) -- (1.9185,0.5058) -- (2.0173,0.4739) -- (1.9140,0.4014) -- cycle;
\path[cbpmesh=17.4] (-3.0824,0.2414) -- (-2.9727,0.1297) -- (-3.0451,0.1067) -- (-3.1575,0.2211) -- cycle;
\path[cbpmesh=18.6] (2.9727,0.1297) -- (3.0824,0.2414) -- (3.1575,0.2211) -- (3.0451,0.1067) -- cycle;
\path[cbpmesh=36.4] (-3.9697,-1.2891) -- (-3.7988,-1.4323) -- (-3.7533,-1.5050) -- (-3.9222,-1.3634) -- cycle;
\path[cbpmesh=38] (3.7988,-1.4323) -- (3.9697,-1.2891) -- (3.9222,-1.3634) -- (3.7533,-1.5050) -- cycle;
\path[cbpmesh=15.1] (-2.4342,0.3925) -- (-2.3294,0.3046) -- (-2.4296,0.2362) -- (-2.5390,0.3278) -- cycle;
\path[cbpmesh=20.3] (2.3294,0.3046) -- (2.4342,0.3925) -- (2.5390,0.3278) -- (2.4296,0.2362) -- cycle;
\path[cbpmesh=17.6] (-3.1575,0.2211) -- (-3.0451,0.1067) -- (-3.1140,0.0864) -- (-3.2289,0.2034) -- cycle;
\path[cbpmesh=18.5] (3.0451,0.1067) -- (3.1575,0.2211) -- (3.2289,0.2034) -- (3.1140,0.0864) -- cycle;
\path[cbpmesh=23.4] (-1.4674,0.4147) -- (-1.3645,0.3625) -- (-1.4364,0.3576) -- (-1.5447,0.4126) -- cycle;
\path[cbpmesh=19.6] (1.3645,0.3625) -- (1.4674,0.4147) -- (1.5447,0.4126) -- (1.4364,0.3576) -- cycle;
\path[cbpmesh=22.7] (-4.1654,-0.6333) -- (-4.0172,-0.7842) -- (-4.0368,-0.9021) -- (-4.1857,-0.7504) -- cycle;
\path[cbpmesh=37.7] (4.0172,-0.7842) -- (4.1654,-0.6333) -- (4.1857,-0.7504) -- (4.0368,-0.9021) -- cycle;
\path[cbpmesh=28.6] (-1.1329,0.2918) -- (-1.0254,0.2520) -- (-1.0779,0.2674) -- (-1.1910,0.3092) -- cycle;
\draw[draw=cbpInk!70,line width=.45pt] (-1.1329,0.2918) -- (-1.0254,0.2520);
\path[cbpmesh=23.0] (1.0254,0.2520) -- (1.1329,0.2918) -- (1.1910,0.3092) -- (1.0779,0.2674) -- cycle;
\draw[draw=cbpInk!70,line width=.45pt] (1.0254,0.2520) -- (1.1329,0.2918);
\path[cbpmesh=15.4] (-2.5390,0.3278) -- (-2.4296,0.2362) -- (-2.5235,0.1779) -- (-2.6370,0.2731) -- cycle;
\path[cbpmesh=19.9] (2.4296,0.2362) -- (2.5390,0.3278) -- (2.6370,0.2731) -- (2.5235,0.1779) -- cycle;
\path[cbpmesh=25.9] (-1.2967,0.3551) -- (-1.1910,0.3092) -- (-1.2533,0.3142) -- (-1.3645,0.3625) -- cycle;
\path[cbpmesh=21.1] (1.1910,0.3092) -- (1.2967,0.3551) -- (1.3645,0.3625) -- (1.2533,0.3142) -- cycle;
\path[cbpmesh=16.7] (-2.0173,0.4739) -- (-1.9140,0.4014) -- (-2.0106,0.3536) -- (-2.1190,0.4298) -- cycle;
\path[cbpmesh=18.6] (1.9140,0.4014) -- (2.0173,0.4739) -- (2.1190,0.4298) -- (2.0106,0.3536) -- cycle;
\path[cbpmesh=17.8] (-3.2289,0.2034) -- (-3.1140,0.0864) -- (-3.1798,0.0677) -- (-3.2971,0.1871) -- cycle;
\path[cbpmesh=18.5] (3.1140,0.0864) -- (3.2289,0.2034) -- (3.2971,0.1871) -- (3.1798,0.0677) -- cycle;
\path[cbpmesh=19.4] (-1.7300,0.4602) -- (-1.6261,0.3983) -- (-1.7110,0.3717) -- (-1.8204,0.4369) -- cycle;
\path[cbpmesh=18.4] (1.6261,0.3983) -- (1.7300,0.4602) -- (1.8204,0.4369) -- (1.7110,0.3717) -- cycle;
\path[cbpmesh=17.7] (-3.2971,0.1871) -- (-3.1798,0.0677) -- (-3.2432,0.0492) -- (-3.3628,0.1711) -- cycle;
\path[cbpmesh=18.5] (3.1798,0.0677) -- (3.2971,0.1871) -- (3.3628,0.1711) -- (3.2432,0.0492) -- cycle;
\path[cbpmesh=15.7] (-2.6370,0.2731) -- (-2.5235,0.1779) -- (-2.6114,0.1286) -- (-2.7289,0.2271) -- cycle;
\path[cbpmesh=19.5] (2.5235,0.1779) -- (2.6370,0.2731) -- (2.7289,0.2271) -- (2.6114,0.1286) -- cycle;
\path[cbpmesh=33.4] (-4.0060,-1.2055) -- (-3.8335,-1.3501) -- (-3.7988,-1.4323) -- (-3.9697,-1.2891) -- cycle;
\path[cbpmesh=38] (3.8335,-1.3501) -- (4.0060,-1.2055) -- (3.9697,-1.2891) -- (3.7988,-1.4323) -- cycle;
\path[cbpmesh=20.6] (-4.1278,-0.5088) -- (-3.9809,-0.6583) -- (-4.0172,-0.7842) -- (-4.1654,-0.6333) -- cycle;
\path[cbpmesh=35.2] (3.9809,-0.6583) -- (4.1278,-0.5088) -- (4.1654,-0.6333) -- (4.0172,-0.7842) -- cycle;
\path[cbpmesh=15.8] (-2.1190,0.4298) -- (-2.0106,0.3536) -- (-2.1095,0.2935) -- (-2.2232,0.3734) -- cycle;
\path[cbpmesh=19.1] (2.0106,0.3536) -- (2.1190,0.4298) -- (2.2232,0.3734) -- (2.1095,0.2935) -- cycle;
\path[cbpmesh=17.4] (-3.3628,0.1711) -- (-3.2432,0.0492) -- (-3.3045,0.0299) -- (-3.4264,0.1540) -- cycle;
\path[cbpmesh=18.6] (3.2432,0.0492) -- (3.3628,0.1711) -- (3.4264,0.1540) -- (3.3045,0.0299) -- cycle;
\path[cbpmesh=21.2] (-1.5447,0.4126) -- (-1.4364,0.3576) -- (-1.5121,0.3404) -- (-1.6261,0.3983) -- cycle;
\path[cbpmesh=18.9] (1.4364,0.3576) -- (1.5447,0.4126) -- (1.6261,0.3983) -- (1.5121,0.3404) -- cycle;
\path[cbpmesh=16.1] (-2.7289,0.2271) -- (-2.6114,0.1286) -- (-2.6939,0.0871) -- (-2.8151,0.1887) -- cycle;
\path[cbpmesh=19.1] (2.6114,0.1286) -- (2.7289,0.2271) -- (2.8151,0.1887) -- (2.6939,0.0871) -- cycle;
\path[cbpmesh=28.5] (-1.0254,0.2520) -- (-0.9115,0.2163) -- (-0.9582,0.2298) -- (-1.0779,0.2674) -- cycle;
\draw[draw=cbpInk!70,line width=.45pt] (-1.0254,0.2520) -- (-0.9115,0.2163);
\draw[draw=cbpInk!70,line width=.45pt] (-0.9115,0.2163) -- (-0.9582,0.2298);
\path[cbpmesh=23.6] (0.9115,0.2163) -- (1.0254,0.2520) -- (1.0779,0.2674) -- (0.9582,0.2298) -- cycle;
\draw[draw=cbpInk!70,line width=.45pt] (0.9115,0.2163) -- (1.0254,0.2520);
\draw[draw=cbpInk!70,line width=.45pt] (0.9582,0.2298) -- (0.9115,0.2163);
\path[cbpmesh=17.0] (-3.4264,0.1540) -- (-3.3045,0.0299) -- (-3.3644,0.0084) -- (-3.4885,0.1348) -- cycle;
\path[cbpmesh=18.8] (3.3045,0.0299) -- (3.4264,0.1540) -- (3.4885,0.1348) -- (3.3644,0.0084) -- cycle;
\path[cbpmesh=38] (-3.6982,-1.5678) -- (-3.5090,-1.7008) -- (-3.4486,-1.7513) -- (-3.6345,-1.6207) -- cycle;
\draw[draw=cbpInk!70,line width=.45pt] (-3.4486,-1.7513) -- (-3.6345,-1.6207);
\draw[draw=cbpInk!70,line width=.45pt] (-3.6345,-1.6207) -- (-3.6982,-1.5678);
\path[cbpmesh=38] (3.5090,-1.7008) -- (3.6982,-1.5678) -- (3.6345,-1.6207) -- (3.4486,-1.7513) -- cycle;
\draw[draw=cbpInk!70,line width=.45pt] (3.6982,-1.5678) -- (3.6345,-1.6207);
\draw[draw=cbpInk!70,line width=.45pt] (3.6345,-1.6207) -- (3.4486,-1.7513);
\path[cbpmesh=17.8] (-1.8204,0.4369) -- (-1.7110,0.3717) -- (-1.7990,0.3329) -- (-1.9140,0.4014) -- cycle;
\path[cbpmesh=18.3] (1.7110,0.3717) -- (1.8204,0.4369) -- (1.9140,0.4014) -- (1.7990,0.3329) -- cycle;
\path[cbpmesh=23.4] (-1.3645,0.3625) -- (-1.2533,0.3142) -- (-1.3193,0.3068) -- (-1.4364,0.3576) -- cycle;
\path[cbpmesh=19.9] (1.2533,0.3142) -- (1.3645,0.3625) -- (1.4364,0.3576) -- (1.3193,0.3068) -- cycle;
\path[cbpmesh=25.9] (-1.1910,0.3092) -- (-1.0779,0.2674) -- (-1.1343,0.2702) -- (-1.2533,0.3142) -- cycle;
\path[cbpmesh=21.5] (1.0779,0.2674) -- (1.1910,0.3092) -- (1.2533,0.3142) -- (1.1343,0.2702) -- cycle;
\path[cbpmesh=18.9] (-4.0716,-0.3772) -- (-3.9267,-0.5247) -- (-3.9809,-0.6583) -- (-4.1278,-0.5088) -- cycle;
\path[cbpmesh=32.8] (3.9267,-0.5247) -- (4.0716,-0.3772) -- (4.1278,-0.5088) -- (3.9809,-0.6583) -- cycle;
\path[cbpmesh=16.4] (-3.4885,0.1348) -- (-3.3644,0.0084) -- (-3.4233,-0.0163) -- (-3.5496,0.1123) -- cycle;
\path[cbpmesh=19.2] (3.3644,0.0084) -- (3.4885,0.1348) -- (3.5496,0.1123) -- (3.4233,-0.0163) -- cycle;
\path[cbpmesh=16.6] (-2.8151,0.1887) -- (-2.6939,0.0871) -- (-2.7715,0.0521) -- (-2.8962,0.1566) -- cycle;
\path[cbpmesh=18.9] (2.6939,0.0871) -- (2.8151,0.1887) -- (2.8962,0.1566) -- (2.7715,0.0521) -- cycle;
\path[cbpmesh=15.1] (-2.2232,0.3734) -- (-2.1095,0.2935) -- (-2.2102,0.2209) -- (-2.3294,0.3046) -- cycle;
\path[cbpmesh=19.7] (2.1095,0.2935) -- (2.2232,0.3734) -- (2.3294,0.3046) -- (2.2102,0.2209) -- cycle;
\path[cbpmesh=30.4] (-4.0301,-1.1130) -- (-3.8565,-1.2585) -- (-3.8335,-1.3501) -- (-4.0060,-1.2055) -- cycle;
\path[cbpmesh=38] (3.8565,-1.2585) -- (4.0301,-1.1130) -- (4.0060,-1.2055) -- (3.8335,-1.3501) -- cycle;
\path[cbpmesh=15.8] (-3.5496,0.1123) -- (-3.4233,-0.0163) -- (-3.4818,-0.0456) -- (-3.6103,0.0852) -- cycle;
\path[cbpmesh=19.7] (3.4233,-0.0163) -- (3.5496,0.1123) -- (3.6103,0.0852) -- (3.4818,-0.0456) -- cycle;
\path[cbpmesh=17.5] (-3.9958,-0.2389) -- (-3.8536,-0.3837) -- (-3.9267,-0.5247) -- (-4.0716,-0.3772) -- cycle;
\path[cbpmesh=30.7] (3.8536,-0.3837) -- (3.9958,-0.2389) -- (4.0716,-0.3772) -- (3.9267,-0.5247) -- cycle;
\path[cbpmesh=15.3] (-3.6103,0.0852) -- (-3.4818,-0.0456) -- (-3.5405,-0.0806) -- (-3.6711,0.0524) -- cycle;
\path[cbpmesh=20.5] (3.4818,-0.0456) -- (3.6103,0.0852) -- (3.6711,0.0524) -- (3.5405,-0.0806) -- cycle;
\path[cbpmesh=17.0] (-2.8962,0.1566) -- (-2.7715,0.0521) -- (-2.8447,0.0224) -- (-2.9727,0.1297) -- cycle;
\path[cbpmesh=18.7] (2.7715,0.0521) -- (2.8962,0.1566) -- (2.9727,0.1297) -- (2.8447,0.0224) -- cycle;
\path[cbpmesh=19.4] (-1.6261,0.3983) -- (-1.5121,0.3404) -- (-1.5911,0.3108) -- (-1.7110,0.3717) -- cycle;
\path[cbpmesh=18.4] (1.5121,0.3404) -- (1.6261,0.3983) -- (1.7110,0.3717) -- (1.5911,0.3108) -- cycle;
\path[cbpmesh=15.0] (-2.3294,0.3046) -- (-2.2102,0.2209) -- (-2.3053,0.1488) -- (-2.4296,0.2362) -- cycle;
\path[cbpmesh=19.9] (2.2102,0.2209) -- (2.3294,0.3046) -- (2.4296,0.2362) -- (2.3053,0.1488) -- cycle;
\path[cbpmesh=15.0] (-3.6711,0.0524) -- (-3.5405,-0.0806) -- (-3.5998,-0.1225) -- (-3.7326,0.0127) -- cycle;
\path[cbpmesh=21.5] (3.5405,-0.0806) -- (3.6711,0.0524) -- (3.7326,0.0127) -- (3.5998,-0.1225) -- cycle;
\path[cbpmesh=16.6] (-1.9140,0.4014) -- (-1.7990,0.3329) -- (-1.8898,0.2816) -- (-2.0106,0.3536) -- cycle;
\path[cbpmesh=18.5] (1.7990,0.3329) -- (1.9140,0.4014) -- (2.0106,0.3536) -- (1.8898,0.2816) -- cycle;
\path[cbpmesh=15.4] (-3.9263,-0.1597) -- (-3.7866,-0.3019) -- (-3.8536,-0.3837) -- (-3.9958,-0.2389) -- cycle;
\path[cbpmesh=26.2] (3.7866,-0.3019) -- (3.9263,-0.1597) -- (3.9958,-0.2389) -- (3.8536,-0.3837) -- cycle;
\path[cbpmesh=14.8] (-3.7326,0.0127) -- (-3.5998,-0.1225) -- (-3.6602,-0.1726) -- (-3.7952,-0.0351) -- cycle;
\path[cbpmesh=22.7] (3.5998,-0.1225) -- (3.7326,0.0127) -- (3.7952,-0.0351) -- (3.6602,-0.1726) -- cycle;
\path[cbpmesh=15.1] (-3.8596,-0.0921) -- (-3.7223,-0.2320) -- (-3.7866,-0.3019) -- (-3.9263,-0.1597) -- cycle;
\path[cbpmesh=25.1] (3.7223,-0.2320) -- (3.8596,-0.0921) -- (3.9263,-0.1597) -- (3.7866,-0.3019) -- cycle;
\path[cbpmesh=14.9] (-3.7952,-0.0351) -- (-3.6602,-0.1726) -- (-3.7223,-0.2320) -- (-3.8596,-0.0921) -- cycle;
\path[cbpmesh=23.9] (3.6602,-0.1726) -- (3.7952,-0.0351) -- (3.8596,-0.0921) -- (3.7223,-0.2320) -- cycle;
\path[cbpmesh=38] (-3.7533,-1.5050) -- (-3.5613,-1.6399) -- (-3.5090,-1.7008) -- (-3.6982,-1.5678) -- cycle;
\path[cbpmesh=38] (3.5613,-1.6399) -- (3.7533,-1.5050) -- (3.6982,-1.5678) -- (3.5090,-1.7008) -- cycle;
\path[cbpmesh=27.4] (-4.0407,-1.0118) -- (-3.8667,-1.1576) -- (-3.8565,-1.2585) -- (-4.0301,-1.1130) -- cycle;
\path[cbpmesh=38] (3.8667,-1.1576) -- (4.0407,-1.0118) -- (4.0301,-1.1130) -- (3.8565,-1.2585) -- cycle;
\path[cbpmesh=17.4] (-2.9727,0.1297) -- (-2.8447,0.0224) -- (-2.9140,-0.0032) -- (-3.0451,0.1067) -- cycle;
\path[cbpmesh=18.5] (2.8447,0.0224) -- (2.9727,0.1297) -- (3.0451,0.1067) -- (2.9140,-0.0032) -- cycle;
\path[cbpmesh=25.9] (-1.0779,0.2674) -- (-0.9582,0.2298) -- (-1.0083,0.2307) -- (-1.1343,0.2702) -- cycle;
\path[cbpmesh=21.9] (0.9582,0.2298) -- (1.0779,0.2674) -- (1.1343,0.2702) -- (1.0083,0.2307) -- cycle;
\path[cbpmesh=21.3] (-1.4364,0.3576) -- (-1.3193,0.3068) -- (-1.3888,0.2869) -- (-1.5121,0.3404) -- cycle;
\path[cbpmesh=19.1] (1.3193,0.3068) -- (1.4364,0.3576) -- (1.5121,0.3404) -- (1.3888,0.2869) -- cycle;
\path[cbpmesh=15.3] (-2.4296,0.2362) -- (-2.3053,0.1488) -- (-2.3944,0.0872) -- (-2.5235,0.1779) -- cycle;
\path[cbpmesh=19.5] (2.3053,0.1488) -- (2.4296,0.2362) -- (2.5235,0.1779) -- (2.3944,0.0872) -- cycle;
\path[cbpmesh=23.4] (-1.2533,0.3142) -- (-1.1343,0.2702) -- (-1.1940,0.2604) -- (-1.3193,0.3068) -- cycle;
\path[cbpmesh=20.2] (1.1343,0.2702) -- (1.2533,0.3142) -- (1.3193,0.3068) -- (1.1940,0.2604) -- cycle;
\path[cbpmesh=17.6] (-3.0451,0.1067) -- (-2.9140,-0.0032) -- (-2.9799,-0.0260) -- (-3.1140,0.0864) -- cycle;
\path[cbpmesh=18.4] (2.9140,-0.0032) -- (3.0451,0.1067) -- (3.1140,0.0864) -- (2.9799,-0.0260) -- cycle;
\path[cbpmesh=15.7] (-2.0106,0.3536) -- (-1.8898,0.2816) -- (-1.9827,0.2180) -- (-2.1095,0.2935) -- cycle;
\path[cbpmesh=18.8] (1.8898,0.2816) -- (2.0106,0.3536) -- (2.1095,0.2935) -- (1.9827,0.2180) -- cycle;
\path[cbpmesh=17.8] (-1.7110,0.3717) -- (-1.5911,0.3108) -- (-1.6729,0.2688) -- (-1.7990,0.3329) -- cycle;
\path[cbpmesh=18.3] (1.5911,0.3108) -- (1.7110,0.3717) -- (1.7990,0.3329) -- (1.6729,0.2688) -- cycle;
\path[cbpmesh=15.6] (-2.5235,0.1779) -- (-2.3944,0.0872) -- (-2.4778,0.0347) -- (-2.6114,0.1286) -- cycle;
\path[cbpmesh=19.2] (2.3944,0.0872) -- (2.5235,0.1779) -- (2.6114,0.1286) -- (2.4778,0.0347) -- cycle;
\path[cbpmesh=24.7] (-4.0368,-0.9021) -- (-3.8629,-1.0478) -- (-3.8667,-1.1576) -- (-4.0407,-1.0118) -- cycle;
\path[cbpmesh=38] (3.8629,-1.0478) -- (4.0368,-0.9021) -- (4.0407,-1.0118) -- (3.8667,-1.1576) -- cycle;
\path[cbpmesh=17.7] (-3.1140,0.0864) -- (-2.9799,-0.0260) -- (-3.0429,-0.0471) -- (-3.1798,0.0677) -- cycle;
\path[cbpmesh=18.4] (2.9799,-0.0260) -- (3.1140,0.0864) -- (3.1798,0.0677) -- (3.0429,-0.0471) -- cycle;
\path[cbpmesh=36.1] (-3.7988,-1.4323) -- (-3.6044,-1.5689) -- (-3.5613,-1.6399) -- (-3.7533,-1.5050) -- cycle;
\path[cbpmesh=38] (3.6044,-1.5689) -- (3.7988,-1.4323) -- (3.7533,-1.5050) -- (3.5613,-1.6399) -- cycle;
\path[cbpmesh=25.8] (-0.9582,0.2298) -- (-0.8326,0.1968) -- (-0.8761,0.1959) -- (-1.0083,0.2307) -- cycle;
\draw[draw=cbpInk!70,line width=.45pt] (-0.9582,0.2298) -- (-0.8326,0.1968);
\path[cbpmesh=22.4] (0.8326,0.1968) -- (0.9582,0.2298) -- (1.0083,0.2307) -- (0.8761,0.1959) -- cycle;
\draw[draw=cbpInk!70,line width=.45pt] (0.8326,0.1968) -- (0.9582,0.2298);
\path[cbpmesh=19.4] (-1.5121,0.3404) -- (-1.3888,0.2869) -- (-1.4613,0.2545) -- (-1.5911,0.3108) -- cycle;
\path[cbpmesh=18.5] (1.3888,0.2869) -- (1.5121,0.3404) -- (1.5911,0.3108) -- (1.4613,0.2545) -- cycle;
\path[cbpmesh=15.0] (-2.1095,0.2935) -- (-1.9827,0.2180) -- (-2.0774,0.1418) -- (-2.2102,0.2209) -- cycle;
\path[cbpmesh=19.4] (1.9827,0.2180) -- (2.1095,0.2935) -- (2.2102,0.2209) -- (2.0774,0.1418) -- cycle;
\path[cbpmesh=16.1] (-2.6114,0.1286) -- (-2.4778,0.0347) -- (-2.5561,-0.0097) -- (-2.6939,0.0871) -- cycle;
\path[cbpmesh=18.9] (2.4778,0.0347) -- (2.6114,0.1286) -- (2.6939,0.0871) -- (2.5561,-0.0097) -- cycle;
\path[cbpmesh=17.7] (-3.1798,0.0677) -- (-3.0429,-0.0471) -- (-3.1035,-0.0679) -- (-3.2432,0.0492) -- cycle;
\path[cbpmesh=18.4] (3.0429,-0.0471) -- (3.1798,0.0677) -- (3.2432,0.0492) -- (3.1035,-0.0679) -- cycle;
\path[cbpmesh=23.4] (-1.1343,0.2702) -- (-1.0083,0.2307) -- (-1.0614,0.2188) -- (-1.1940,0.2604) -- cycle;
\path[cbpmesh=20.6] (1.0083,0.2307) -- (1.1343,0.2702) -- (1.1940,0.2604) -- (1.0614,0.2188) -- cycle;
\path[cbpmesh=21.3] (-1.3193,0.3068) -- (-1.1940,0.2604) -- (-1.2569,0.2381) -- (-1.3888,0.2869) -- cycle;
\path[cbpmesh=19.3] (1.1940,0.2604) -- (1.3193,0.3068) -- (1.3888,0.2869) -- (1.2569,0.2381) -- cycle;
\path[cbpmesh=16.6] (-1.7990,0.3329) -- (-1.6729,0.2688) -- (-1.7573,0.2144) -- (-1.8898,0.2816) -- cycle;
\path[cbpmesh=18.3] (1.6729,0.2688) -- (1.7990,0.3329) -- (1.8898,0.2816) -- (1.7573,0.2144) -- cycle;
\path[cbpmesh=22.2] (-4.0172,-0.7842) -- (-3.8442,-0.9292) -- (-3.8629,-1.0478) -- (-4.0368,-0.9021) -- cycle;
\path[cbpmesh=36.7] (3.8442,-0.9292) -- (4.0172,-0.7842) -- (4.0368,-0.9021) -- (3.8629,-1.0478) -- cycle;
\path[cbpmesh=17.4] (-3.2432,0.0492) -- (-3.1035,-0.0679) -- (-3.1622,-0.0894) -- (-3.3045,0.0299) -- cycle;
\path[cbpmesh=18.5] (3.1035,-0.0679) -- (3.2432,0.0492) -- (3.3045,0.0299) -- (3.1622,-0.0894) -- cycle;
\path[cbpmesh=16.5] (-2.6939,0.0871) -- (-2.5561,-0.0097) -- (-2.6297,-0.0475) -- (-2.7715,0.0521) -- cycle;
\path[cbpmesh=18.7] (2.5561,-0.0097) -- (2.6939,0.0871) -- (2.7715,0.0521) -- (2.6297,-0.0475) -- cycle;
\path[cbpmesh=25.7] (-0.8326,0.1968) -- (-0.7018,0.1685) -- (-0.7385,0.1661) -- (-0.8761,0.1959) -- cycle;
\draw[draw=cbpInk!70,line width=.45pt] (-0.8326,0.1968) -- (-0.7018,0.1685);
\path[cbpmesh=22.8] (0.7018,0.1685) -- (0.8326,0.1968) -- (0.8761,0.1959) -- (0.7385,0.1661) -- cycle;
\draw[draw=cbpInk!70,line width=.45pt] (0.7018,0.1685) -- (0.8326,0.1968);
\path[cbpmesh=14.8] (-2.2102,0.2209) -- (-2.0774,0.1418) -- (-2.1668,0.0663) -- (-2.3053,0.1488) -- cycle;
\path[cbpmesh=19.6] (2.0774,0.1418) -- (2.2102,0.2209) -- (2.3053,0.1488) -- (2.1668,0.0663) -- cycle;
\path[cbpmesh=33.0] (-3.8335,-1.3501) -- (-3.6374,-1.4879) -- (-3.6044,-1.5689) -- (-3.7988,-1.4323) -- cycle;
\path[cbpmesh=38] (3.6374,-1.4879) -- (3.8335,-1.3501) -- (3.7988,-1.4323) -- (3.6044,-1.5689) -- cycle;
\path[cbpmesh=16.9] (-3.3045,0.0299) -- (-3.1622,-0.0894) -- (-3.2195,-0.1130) -- (-3.3644,0.0084) -- cycle;
\path[cbpmesh=18.7] (3.1622,-0.0894) -- (3.3045,0.0299) -- (3.3644,0.0084) -- (3.2195,-0.1130) -- cycle;
\path[cbpmesh=17.8] (-1.5911,0.3108) -- (-1.4613,0.2545) -- (-1.5365,0.2096) -- (-1.6729,0.2688) -- cycle;
\path[cbpmesh=18.2] (1.4613,0.2545) -- (1.5911,0.3108) -- (1.6729,0.2688) -- (1.5365,0.2096) -- cycle;
\path[cbpmesh=20.2] (-3.9809,-0.6583) -- (-3.8095,-0.8020) -- (-3.8442,-0.9292) -- (-4.0172,-0.7842) -- cycle;
\path[cbpmesh=34.2] (3.8095,-0.8020) -- (3.9809,-0.6583) -- (4.0172,-0.7842) -- (3.8442,-0.9292) -- cycle;
\path[cbpmesh=38] (-3.5090,-1.7008) -- (-3.2982,-1.8263) -- (-3.2414,-1.8747) -- (-3.4486,-1.7513) -- cycle;
\draw[draw=cbpInk!70,line width=.45pt] (-3.2414,-1.8747) -- (-3.4486,-1.7513);
\path[cbpmesh=38] (3.2982,-1.8263) -- (3.5090,-1.7008) -- (3.4486,-1.7513) -- (3.2414,-1.8747) -- cycle;
\draw[draw=cbpInk!70,line width=.45pt] (3.4486,-1.7513) -- (3.2414,-1.8747);
\path[cbpmesh=16.9] (-2.7715,0.0521) -- (-2.6297,-0.0475) -- (-2.6992,-0.0799) -- (-2.8447,0.0224) -- cycle;
\path[cbpmesh=18.5] (2.6297,-0.0475) -- (2.7715,0.0521) -- (2.8447,0.0224) -- (2.6992,-0.0799) -- cycle;
\path[cbpmesh=23.4] (-1.0083,0.2307) -- (-0.8761,0.1959) -- (-0.9223,0.1822) -- (-1.0614,0.2188) -- cycle;
\path[cbpmesh=20.9] (0.8761,0.1959) -- (1.0083,0.2307) -- (1.0614,0.2188) -- (0.9223,0.1822) -- cycle;
\path[cbpmesh=15.6] (-1.8898,0.2816) -- (-1.7573,0.2144) -- (-1.8437,0.1474) -- (-1.9827,0.2180) -- cycle;
\path[cbpmesh=18.6] (1.7573,0.2144) -- (1.8898,0.2816) -- (1.9827,0.2180) -- (1.8437,0.1474) -- cycle;
\path[cbpmesh=16.3] (-3.3644,0.0084) -- (-3.2195,-0.1130) -- (-3.2759,-0.1399) -- (-3.4233,-0.0163) -- cycle;
\path[cbpmesh=19.0] (3.2195,-0.1130) -- (3.3644,0.0084) -- (3.4233,-0.0163) -- (3.2759,-0.1399) -- cycle;
\path[cbpmesh=15.2] (-2.3053,0.1488) -- (-2.1668,0.0663) -- (-2.2505,0.0015) -- (-2.3944,0.0872) -- cycle;
\path[cbpmesh=19.2] (2.1668,0.0663) -- (2.3053,0.1488) -- (2.3944,0.0872) -- (2.2505,0.0015) -- cycle;
\path[cbpmesh=25.6] (-0.7018,0.1685) -- (-0.5667,0.1451) -- (-0.5963,0.1415) -- (-0.7385,0.1661) -- cycle;
\draw[draw=cbpInk!70,line width=.45pt] (-0.7018,0.1685) -- (-0.5667,0.1451);
\path[cbpmesh=23.2] (0.5667,0.1451) -- (0.7018,0.1685) -- (0.7385,0.1661) -- (0.5963,0.1415) -- cycle;
\draw[draw=cbpInk!70,line width=.45pt] (0.5667,0.1451) -- (0.7018,0.1685);
\path[cbpmesh=19.4] (-1.3888,0.2869) -- (-1.2569,0.2381) -- (-1.3226,0.2032) -- (-1.4613,0.2545) -- cycle;
\path[cbpmesh=18.6] (1.2569,0.2381) -- (1.3888,0.2869) -- (1.4613,0.2545) -- (1.3226,0.2032) -- cycle;
\path[cbpmesh=21.3] (-1.1940,0.2604) -- (-1.0614,0.2188) -- (-1.1173,0.1943) -- (-1.2569,0.2381) -- cycle;
\path[cbpmesh=19.5] (1.0614,0.2188) -- (1.1940,0.2604) -- (1.2569,0.2381) -- (1.1173,0.1943) -- cycle;
\path[cbpmesh=15.7] (-3.4233,-0.0163) -- (-3.2759,-0.1399) -- (-3.3319,-0.1713) -- (-3.4818,-0.0456) -- cycle;
\path[cbpmesh=19.5] (3.2759,-0.1399) -- (3.4233,-0.0163) -- (3.4818,-0.0456) -- (3.3319,-0.1713) -- cycle;
\path[cbpmesh=17.3] (-2.8447,0.0224) -- (-2.6992,-0.0799) -- (-2.7649,-0.1080) -- (-2.9140,-0.0032) -- cycle;
\path[cbpmesh=18.4] (2.6992,-0.0799) -- (2.8447,0.0224) -- (2.9140,-0.0032) -- (2.7649,-0.1080) -- cycle;
\path[cbpmesh=18.5] (-3.9267,-0.5247) -- (-3.7577,-0.6665) -- (-3.8095,-0.8020) -- (-3.9809,-0.6583) -- cycle;
\path[cbpmesh=31.9] (3.7577,-0.6665) -- (3.9267,-0.5247) -- (3.9809,-0.6583) -- (3.8095,-0.8020) -- cycle;
\path[cbpmesh=30.0] (-3.8565,-1.2585) -- (-3.6592,-1.3971) -- (-3.6374,-1.4879) -- (-3.8335,-1.3501) -- cycle;
\path[cbpmesh=38] (3.6592,-1.3971) -- (3.8565,-1.2585) -- (3.8335,-1.3501) -- (3.6374,-1.4879) -- cycle;
\path[cbpmesh=15.2] (-3.4818,-0.0456) -- (-3.3319,-0.1713) -- (-3.3880,-0.2084) -- (-3.5405,-0.0806) -- cycle;
\path[cbpmesh=20.2] (3.3319,-0.1713) -- (3.4818,-0.0456) -- (3.5405,-0.0806) -- (3.3880,-0.2084) -- cycle;
\path[cbpmesh=15.6] (-2.3944,0.0872) -- (-2.2505,0.0015) -- (-2.3289,-0.0539) -- (-2.4778,0.0347) -- cycle;
\path[cbpmesh=18.9] (2.2505,0.0015) -- (2.3944,0.0872) -- (2.4778,0.0347) -- (2.3289,-0.0539) -- cycle;
\path[cbpmesh=25.5] (-0.5667,0.1451) -- (-0.4281,0.1267) -- (-0.4505,0.1222) -- (-0.5963,0.1415) -- cycle;
\draw[draw=cbpInk!70,line width=.45pt] (-0.5667,0.1451) -- (-0.4281,0.1267);
\path[cbpmesh=23.5] (0.4281,0.1267) -- (0.5667,0.1451) -- (0.5963,0.1415) -- (0.4505,0.1222) -- cycle;
\draw[draw=cbpInk!70,line width=.45pt] (0.4281,0.1267) -- (0.5667,0.1451);
\path[cbpmesh=14.9] (-1.9827,0.2180) -- (-1.8437,0.1474) -- (-1.9318,0.0678) -- (-2.0774,0.1418) -- cycle;
\path[cbpmesh=19.0] (1.8437,0.1474) -- (1.9827,0.2180) -- (2.0774,0.1418) -- (1.9318,0.0678) -- cycle;
\path[cbpmesh=14.8] (-3.5405,-0.0806) -- (-3.3880,-0.2084) -- (-3.4447,-0.2525) -- (-3.5998,-0.1225) -- cycle;
\path[cbpmesh=21.1] (3.3880,-0.2084) -- (3.5405,-0.0806) -- (3.5998,-0.1225) -- (3.4447,-0.2525) -- cycle;
\path[cbpmesh=23.3] (-0.8761,0.1959) -- (-0.7385,0.1661) -- (-0.7774,0.1509) -- (-0.9223,0.1822) -- cycle;
\path[cbpmesh=21.2] (0.7385,0.1661) -- (0.8761,0.1959) -- (0.9223,0.1822) -- (0.7774,0.1509) -- cycle;
\path[cbpmesh=16.6] (-1.6729,0.2688) -- (-1.5365,0.2096) -- (-1.6140,0.1522) -- (-1.7573,0.2144) -- cycle;
\path[cbpmesh=18.2] (1.5365,0.2096) -- (1.6729,0.2688) -- (1.7573,0.2144) -- (1.6140,0.1522) -- cycle;
\path[cbpmesh=17.1] (-3.8536,-0.3837) -- (-3.6877,-0.5228) -- (-3.7577,-0.6665) -- (-3.9267,-0.5247) -- cycle;
\path[cbpmesh=29.8] (3.6877,-0.5228) -- (3.8536,-0.3837) -- (3.9267,-0.5247) -- (3.7577,-0.6665) -- cycle;
\path[cbpmesh=17.6] (-2.9140,-0.0032) -- (-2.7649,-0.1080) -- (-2.8274,-0.1331) -- (-2.9799,-0.0260) -- cycle;
\path[cbpmesh=18.4] (2.7649,-0.1080) -- (2.9140,-0.0032) -- (2.9799,-0.0260) -- (2.8274,-0.1331) -- cycle;
\path[cbpmesh=14.6] (-3.5998,-0.1225) -- (-3.4447,-0.2525) -- (-3.5026,-0.3047) -- (-3.6602,-0.1726) -- cycle;
\path[cbpmesh=22.1] (3.4447,-0.2525) -- (3.5998,-0.1225) -- (3.6602,-0.1726) -- (3.5026,-0.3047) -- cycle;
\path[cbpmesh=38] (-3.5613,-1.6399) -- (-3.3473,-1.7673) -- (-3.2982,-1.8263) -- (-3.5090,-1.7008) -- cycle;
\path[cbpmesh=38] (3.3473,-1.7673) -- (3.5613,-1.6399) -- (3.5090,-1.7008) -- (3.2982,-1.8263) -- cycle;
\path[cbpmesh=14.6] (-3.6602,-0.1726) -- (-3.5026,-0.3047) -- (-3.5620,-0.3664) -- (-3.7223,-0.2320) -- cycle;
\path[cbpmesh=23.2] (3.5026,-0.3047) -- (3.6602,-0.1726) -- (3.7223,-0.2320) -- (3.5620,-0.3664) -- cycle;
\path[cbpmesh=25.3] (-0.4281,0.1267) -- (-0.2869,0.1135) -- (-0.3019,0.1083) -- (-0.4505,0.1222) -- cycle;
\draw[draw=cbpInk!70,line width=.45pt] (-0.4281,0.1267) -- (-0.2869,0.1135);
\path[cbpmesh=23.9] (0.2869,0.1135) -- (0.4281,0.1267) -- (0.4505,0.1222) -- (0.3019,0.1083) -- cycle;
\draw[draw=cbpInk!70,line width=.45pt] (0.2869,0.1135) -- (0.4281,0.1267);
\path[cbpmesh=15.1] (-3.7866,-0.3019) -- (-3.6235,-0.4386) -- (-3.6877,-0.5228) -- (-3.8536,-0.3837) -- cycle;
\path[cbpmesh=25.4] (3.6235,-0.4386) -- (3.7866,-0.3019) -- (3.8536,-0.3837) -- (3.6877,-0.5228) -- cycle;
\path[cbpmesh=14.8] (-3.7223,-0.2320) -- (-3.5620,-0.3664) -- (-3.6235,-0.4386) -- (-3.7866,-0.3019) -- cycle;
\path[cbpmesh=24.4] (3.5620,-0.3664) -- (3.7223,-0.2320) -- (3.7866,-0.3019) -- (3.6235,-0.4386) -- cycle;
\path[cbpmesh=21.2] (-1.0614,0.2188) -- (-0.9223,0.1822) -- (-0.9708,0.1558) -- (-1.1173,0.1943) -- cycle;
\path[cbpmesh=19.7] (0.9223,0.1822) -- (1.0614,0.2188) -- (1.1173,0.1943) -- (0.9708,0.1558) -- cycle;
\path[cbpmesh=17.8] (-1.4613,0.2545) -- (-1.3226,0.2032) -- (-1.3906,0.1556) -- (-1.5365,0.2096) -- cycle;
\path[cbpmesh=18.2] (1.3226,0.2032) -- (1.4613,0.2545) -- (1.5365,0.2096) -- (1.3906,0.1556) -- cycle;
\path[cbpmesh=16.0] (-2.4778,0.0347) -- (-2.3289,-0.0539) -- (-2.4025,-0.1012) -- (-2.5561,-0.0097) -- cycle;
\path[cbpmesh=18.7] (2.3289,-0.0539) -- (2.4778,0.0347) -- (2.5561,-0.0097) -- (2.4025,-0.1012) -- cycle;
\path[cbpmesh=27.1] (-3.8667,-1.1576) -- (-3.6689,-1.2966) -- (-3.6592,-1.3971) -- (-3.8565,-1.2585) -- cycle;
\path[cbpmesh=38] (3.6689,-1.2966) -- (3.8667,-1.1576) -- (3.8565,-1.2585) -- (3.6592,-1.3971) -- cycle;
\path[cbpmesh=19.4] (-1.2569,0.2381) -- (-1.1173,0.1943) -- (-1.1757,0.1571) -- (-1.3226,0.2032) -- cycle;
\path[cbpmesh=18.7] (1.1173,0.1943) -- (1.2569,0.2381) -- (1.3226,0.2032) -- (1.1757,0.1571) -- cycle;
\path[cbpmesh=25.1] (-0.2869,0.1135) -- (-0.1439,0.1056) -- (-0.1514,0.0999) -- (-0.3019,0.1083) -- cycle;
\draw[draw=cbpInk!70,line width=.45pt] (-0.2869,0.1135) -- (-0.1439,0.1056);
\path[cbpmesh=24.2] (0.1439,0.1056) -- (0.2869,0.1135) -- (0.3019,0.1083) -- (0.1514,0.0999) -- cycle;
\draw[draw=cbpInk!70,line width=.45pt] (0.1439,0.1056) -- (0.2869,0.1135);
\path[cbpmesh=17.7] (-2.9799,-0.0260) -- (-2.8274,-0.1331) -- (-2.8872,-0.1565) -- (-3.0429,-0.0471) -- cycle;
\path[cbpmesh=18.4] (2.8274,-0.1331) -- (2.9799,-0.0260) -- (3.0429,-0.0471) -- (2.8872,-0.1565) -- cycle;
\path[cbpmesh=14.8] (-2.0774,0.1418) -- (-1.9318,0.0678) -- (-2.0149,-0.0108) -- (-2.1668,0.0663) -- cycle;
\path[cbpmesh=19.2] (1.9318,0.0678) -- (2.0774,0.1418) -- (2.1668,0.0663) -- (2.0149,-0.0108) -- cycle;
\path[cbpmesh=23.2] (-0.7385,0.1661) -- (-0.5963,0.1415) -- (-0.6277,0.1250) -- (-0.7774,0.1509) -- cycle;
\path[cbpmesh=21.5] (0.5963,0.1415) -- (0.7385,0.1661) -- (0.7774,0.1509) -- (0.6277,0.1250) -- cycle;
\path[cbpmesh=24.8] (-0.1439,0.1056) -- (-0.0000,0.1029) -- (-0.0000,0.0972) -- (-0.1514,0.0999) -- cycle;
\draw[draw=cbpInk!70,line width=.45pt] (-0.1439,0.1056) -- (-0.0000,0.1029);
\path[cbpmesh=24.5] (-0.0000,0.1029) -- (0.1439,0.1056) -- (0.1514,0.0999) -- (-0.0000,0.0972) -- cycle;
\draw[draw=cbpInk!70,line width=.45pt] (-0.0000,0.1029) -- (0.1439,0.1056);
\path[cbpmesh=15.6] (-1.7573,0.2144) -- (-1.6140,0.1522) -- (-1.6934,0.0821) -- (-1.8437,0.1474) -- cycle;
\path[cbpmesh=18.3] (1.6140,0.1522) -- (1.7573,0.2144) -- (1.8437,0.1474) -- (1.6934,0.0821) -- cycle;
\path[cbpmesh=38] (-3.2414,-1.8747) -- (-3.0142,-1.9901) -- (-2.9551,-2.0255) -- (-3.1779,-1.9124) -- cycle;
\draw[draw=cbpInk!70,line width=.45pt] (-2.9551,-2.0255) -- (-3.1779,-1.9124);
\draw[draw=cbpInk!70,line width=.45pt] (-3.1779,-1.9124) -- (-3.2414,-1.8747);
\path[cbpmesh=38] (3.0142,-1.9901) -- (3.2414,-1.8747) -- (3.1779,-1.9124) -- (2.9551,-2.0255) -- cycle;
\draw[draw=cbpInk!70,line width=.45pt] (3.2414,-1.8747) -- (3.1779,-1.9124);
\draw[draw=cbpInk!70,line width=.45pt] (3.1779,-1.9124) -- (2.9551,-2.0255);
\path[cbpmesh=16.5] (-2.5561,-0.0097) -- (-2.4025,-0.1012) -- (-2.4717,-0.1416) -- (-2.6297,-0.0475) -- cycle;
\path[cbpmesh=18.5] (2.4025,-0.1012) -- (2.5561,-0.0097) -- (2.6297,-0.0475) -- (2.4717,-0.1416) -- cycle;
\path[cbpmesh=17.6] (-3.0429,-0.0471) -- (-2.8872,-0.1565) -- (-2.9447,-0.1794) -- (-3.1035,-0.0679) -- cycle;
\path[cbpmesh=18.4] (2.8872,-0.1565) -- (3.0429,-0.0471) -- (3.1035,-0.0679) -- (2.9447,-0.1794) -- cycle;
\path[cbpmesh=35.8] (-3.6044,-1.5689) -- (-3.3879,-1.6979) -- (-3.3473,-1.7673) -- (-3.5613,-1.6399) -- cycle;
\path[cbpmesh=38] (3.3879,-1.6979) -- (3.6044,-1.5689) -- (3.5613,-1.6399) -- (3.3473,-1.7673) -- cycle;
\path[cbpmesh=23.1] (-0.5963,0.1415) -- (-0.4505,0.1222) -- (-0.4742,0.1046) -- (-0.6277,0.1250) -- cycle;
\path[cbpmesh=21.7] (0.4505,0.1222) -- (0.5963,0.1415) -- (0.6277,0.1250) -- (0.4742,0.1046) -- cycle;
\path[cbpmesh=21.2] (-0.9223,0.1822) -- (-0.7774,0.1509) -- (-0.8183,0.1228) -- (-0.9708,0.1558) -- cycle;
\path[cbpmesh=19.9] (0.7774,0.1509) -- (0.9223,0.1822) -- (0.9708,0.1558) -- (0.8183,0.1228) -- cycle;
\path[cbpmesh=24.3] (-3.8629,-1.0478) -- (-3.6653,-1.1867) -- (-3.6689,-1.2966) -- (-3.8667,-1.1576) -- cycle;
\path[cbpmesh=38] (3.6653,-1.1867) -- (3.8629,-1.0478) -- (3.8667,-1.1576) -- (3.6689,-1.2966) -- cycle;
\path[cbpmesh=15.1] (-2.1668,0.0663) -- (-2.0149,-0.0108) -- (-2.0927,-0.0786) -- (-2.2505,0.0015) -- cycle;
\path[cbpmesh=18.9] (2.0149,-0.0108) -- (2.1668,0.0663) -- (2.2505,0.0015) -- (2.0927,-0.0786) -- cycle;
\path[cbpmesh=16.6] (-1.5365,0.2096) -- (-1.3906,0.1556) -- (-1.4608,0.0955) -- (-1.6140,0.1522) -- cycle;
\path[cbpmesh=18.0] (1.3906,0.1556) -- (1.5365,0.2096) -- (1.6140,0.1522) -- (1.4608,0.0955) -- cycle;
\path[cbpmesh=17.4] (-3.1035,-0.0679) -- (-2.9447,-0.1794) -- (-3.0004,-0.2031) -- (-3.1622,-0.0894) -- cycle;
\path[cbpmesh=18.4] (2.9447,-0.1794) -- (3.1035,-0.0679) -- (3.1622,-0.0894) -- (3.0004,-0.2031) -- cycle;
\path[cbpmesh=16.9] (-2.6297,-0.0475) -- (-2.4717,-0.1416) -- (-2.5370,-0.1764) -- (-2.6992,-0.0799) -- cycle;
\path[cbpmesh=18.4] (2.4717,-0.1416) -- (2.6297,-0.0475) -- (2.6992,-0.0799) -- (2.5370,-0.1764) -- cycle;
\path[cbpmesh=19.4] (-1.1173,0.1943) -- (-0.9708,0.1558) -- (-1.0215,0.1166) -- (-1.1757,0.1571) -- cycle;
\path[cbpmesh=18.8] (0.9708,0.1558) -- (1.1173,0.1943) -- (1.1757,0.1571) -- (1.0215,0.1166) -- cycle;
\path[cbpmesh=17.8] (-1.3226,0.2032) -- (-1.1757,0.1571) -- (-1.2362,0.1072) -- (-1.3906,0.1556) -- cycle;
\path[cbpmesh=18.2] (1.1757,0.1571) -- (1.3226,0.2032) -- (1.3906,0.1556) -- (1.2362,0.1072) -- cycle;
\path[cbpmesh=23.0] (-0.4505,0.1222) -- (-0.3019,0.1083) -- (-0.3178,0.0900) -- (-0.4742,0.1046) -- cycle;
\path[cbpmesh=22.0] (0.3019,0.1083) -- (0.4505,0.1222) -- (0.4742,0.1046) -- (0.3178,0.0900) -- cycle;
\path[cbpmesh=14.9] (-1.8437,0.1474) -- (-1.6934,0.0821) -- (-1.7743,-0.0005) -- (-1.9318,0.0678) -- cycle;
\path[cbpmesh=18.7] (1.6934,0.0821) -- (1.8437,0.1474) -- (1.9318,0.0678) -- (1.7743,-0.0005) -- cycle;
\path[cbpmesh=16.9] (-3.1622,-0.0894) -- (-3.0004,-0.2031) -- (-3.0548,-0.2288) -- (-3.2195,-0.1130) -- cycle;
\path[cbpmesh=18.6] (3.0004,-0.2031) -- (3.1622,-0.0894) -- (3.2195,-0.1130) -- (3.0548,-0.2288) -- cycle;
\path[cbpmesh=22.8] (-0.3019,0.1083) -- (-0.1514,0.0999) -- (-0.1594,0.0812) -- (-0.3178,0.0900) -- cycle;
\path[cbpmesh=22.2] (0.1514,0.0999) -- (0.3019,0.1083) -- (0.3178,0.0900) -- (0.1594,0.0812) -- cycle;
\path[cbpmesh=15.5] (-2.2505,0.0015) -- (-2.0927,-0.0786) -- (-2.1657,-0.1368) -- (-2.3289,-0.0539) -- cycle;
\path[cbpmesh=18.6] (2.0927,-0.0786) -- (2.2505,0.0015) -- (2.3289,-0.0539) -- (2.1657,-0.1368) -- cycle;
\path[cbpmesh=21.9] (-3.8442,-0.9292) -- (-3.6475,-1.0674) -- (-3.6653,-1.1867) -- (-3.8629,-1.0478) -- cycle;
\path[cbpmesh=35.6] (3.6475,-1.0674) -- (3.8442,-0.9292) -- (3.8629,-1.0478) -- (3.6653,-1.1867) -- cycle;
\path[cbpmesh=21.1] (-0.7774,0.1509) -- (-0.6277,0.1250) -- (-0.6608,0.0955) -- (-0.8183,0.1228) -- cycle;
\path[cbpmesh=20.0] (0.6277,0.1250) -- (0.7774,0.1509) -- (0.8183,0.1228) -- (0.6608,0.0955) -- cycle;
\path[cbpmesh=22.7] (-0.1514,0.0999) -- (-0.0000,0.0972) -- (-0.0000,0.0783) -- (-0.1594,0.0812) -- cycle;
\path[cbpmesh=22.5] (-0.0000,0.0972) -- (0.1514,0.0999) -- (0.1594,0.0812) -- (-0.0000,0.0783) -- cycle;
\path[cbpmesh=38] (-3.2982,-1.8263) -- (-3.0670,-1.9438) -- (-3.0142,-1.9901) -- (-3.2414,-1.8747) -- cycle;
\path[cbpmesh=38] (3.0670,-1.9438) -- (3.2982,-1.8263) -- (3.2414,-1.8747) -- (3.0142,-1.9901) -- cycle;
\path[cbpmesh=32.8] (-3.6374,-1.4879) -- (-3.4188,-1.6181) -- (-3.3879,-1.6979) -- (-3.6044,-1.5689) -- cycle;
\path[cbpmesh=38] (3.4188,-1.6181) -- (3.6374,-1.4879) -- (3.6044,-1.5689) -- (3.3879,-1.6979) -- cycle;
\path[cbpmesh=17.3] (-2.6992,-0.0799) -- (-2.5370,-0.1764) -- (-2.5988,-0.2069) -- (-2.7649,-0.1080) -- cycle;
\path[cbpmesh=18.3] (2.5370,-0.1764) -- (2.6992,-0.0799) -- (2.7649,-0.1080) -- (2.5988,-0.2069) -- cycle;
\path[cbpmesh=16.3] (-3.2195,-0.1130) -- (-3.0548,-0.2288) -- (-3.1083,-0.2577) -- (-3.2759,-0.1399) -- cycle;
\path[cbpmesh=18.8] (3.0548,-0.2288) -- (3.2195,-0.1130) -- (3.2759,-0.1399) -- (3.1083,-0.2577) -- cycle;
\path[cbpmesh=15.6] (-1.6140,0.1522) -- (-1.4608,0.0955) -- (-1.5326,0.0226) -- (-1.6934,0.0821) -- cycle;
\path[cbpmesh=18.1] (1.4608,0.0955) -- (1.6140,0.1522) -- (1.6934,0.0821) -- (1.5326,0.0226) -- cycle;
\path[cbpmesh=14.7] (-1.9318,0.0678) -- (-1.7743,-0.0005) -- (-1.8506,-0.0821) -- (-2.0149,-0.0108) -- cycle;
\path[cbpmesh=18.8] (1.7743,-0.0005) -- (1.9318,0.0678) -- (2.0149,-0.0108) -- (1.8506,-0.0821) -- cycle;
\path[cbpmesh=19.4] (-0.9708,0.1558) -- (-0.8183,0.1228) -- (-0.8611,0.0818) -- (-1.0215,0.1166) -- cycle;
\path[cbpmesh=18.8] (0.8183,0.1228) -- (0.9708,0.1558) -- (1.0215,0.1166) -- (0.8611,0.0818) -- cycle;
\path[cbpmesh=15.6] (-3.2759,-0.1399) -- (-3.1083,-0.2577) -- (-3.1615,-0.2911) -- (-3.3319,-0.1713) -- cycle;
\path[cbpmesh=19.2] (3.1083,-0.2577) -- (3.2759,-0.1399) -- (3.3319,-0.1713) -- (3.1615,-0.2911) -- cycle;
\path[cbpmesh=21.1] (-0.6277,0.1250) -- (-0.4742,0.1046) -- (-0.4992,0.0741) -- (-0.6608,0.0955) -- cycle;
\path[cbpmesh=20.2] (0.4742,0.1046) -- (0.6277,0.1250) -- (0.6608,0.0955) -- (0.4992,0.0741) -- cycle;
\path[cbpmesh=19.8] (-3.8095,-0.8020) -- (-3.6146,-0.9390) -- (-3.6475,-1.0674) -- (-3.8442,-0.9292) -- cycle;
\path[cbpmesh=33.2] (3.6146,-0.9390) -- (3.8095,-0.8020) -- (3.8442,-0.9292) -- (3.6475,-1.0674) -- cycle;
\path[cbpmesh=16.0] (-2.3289,-0.0539) -- (-2.1657,-0.1368) -- (-2.2341,-0.1867) -- (-2.4025,-0.1012) -- cycle;
\path[cbpmesh=18.5] (2.1657,-0.1368) -- (2.3289,-0.0539) -- (2.4025,-0.1012) -- (2.2341,-0.1867) -- cycle;
\path[cbpmesh=17.6] (-2.7649,-0.1080) -- (-2.5988,-0.2069) -- (-2.6576,-0.2343) -- (-2.8274,-0.1331) -- cycle;
\path[cbpmesh=18.3] (2.5988,-0.2069) -- (2.7649,-0.1080) -- (2.8274,-0.1331) -- (2.6576,-0.2343) -- cycle;
\path[cbpmesh=16.6] (-1.3906,0.1556) -- (-1.2362,0.1072) -- (-1.2985,0.0446) -- (-1.4608,0.0955) -- cycle;
\path[cbpmesh=17.9] (1.2362,0.1072) -- (1.3906,0.1556) -- (1.4608,0.0955) -- (1.2985,0.0446) -- cycle;
\path[cbpmesh=17.8] (-1.1757,0.1571) -- (-1.0215,0.1166) -- (-1.0741,0.0646) -- (-1.2362,0.1072) -- cycle;
\path[cbpmesh=18.1] (1.0215,0.1166) -- (1.1757,0.1571) -- (1.2362,0.1072) -- (1.0741,0.0646) -- cycle;
\path[cbpmesh=15.0] (-3.3319,-0.1713) -- (-3.1615,-0.2911) -- (-3.2147,-0.3302) -- (-3.3880,-0.2084) -- cycle;
\path[cbpmesh=19.8] (3.1615,-0.2911) -- (3.3319,-0.1713) -- (3.3880,-0.2084) -- (3.2147,-0.3302) -- cycle;
\path[cbpmesh=21.0] (-0.4742,0.1046) -- (-0.3178,0.0900) -- (-0.3345,0.0587) -- (-0.4992,0.0741) -- cycle;
\path[cbpmesh=20.4] (0.3178,0.0900) -- (0.4742,0.1046) -- (0.4992,0.0741) -- (0.3345,0.0587) -- cycle;
\path[cbpmesh=29.7] (-3.6592,-1.3971) -- (-3.4394,-1.5280) -- (-3.4188,-1.6181) -- (-3.6374,-1.4879) -- cycle;
\path[cbpmesh=38] (3.4394,-1.5280) -- (3.6592,-1.3971) -- (3.6374,-1.4879) -- (3.4188,-1.6181) -- cycle;
\path[cbpmesh=18.1] (-3.7577,-0.6665) -- (-3.5654,-0.8016) -- (-3.6146,-0.9390) -- (-3.8095,-0.8020) -- cycle;
\path[cbpmesh=30.9] (3.5654,-0.8016) -- (3.7577,-0.6665) -- (3.8095,-0.8020) -- (3.6146,-0.9390) -- cycle;
\path[cbpmesh=14.6] (-3.3880,-0.2084) -- (-3.2147,-0.3302) -- (-3.2685,-0.3763) -- (-3.4447,-0.2525) -- cycle;
\path[cbpmesh=20.6] (3.2147,-0.3302) -- (3.3880,-0.2084) -- (3.4447,-0.2525) -- (3.2685,-0.3763) -- cycle;
\path[cbpmesh=38] (-3.3473,-1.7673) -- (-3.1127,-1.8865) -- (-3.0670,-1.9438) -- (-3.2982,-1.8263) -- cycle;
\path[cbpmesh=38] (3.1127,-1.8865) -- (3.3473,-1.7673) -- (3.2982,-1.8263) -- (3.0670,-1.9438) -- cycle;
\path[cbpmesh=17.7] (-2.8274,-0.1331) -- (-2.6576,-0.2343) -- (-2.7137,-0.2598) -- (-2.8872,-0.1565) -- cycle;
\path[cbpmesh=18.3] (2.6576,-0.2343) -- (2.8274,-0.1331) -- (2.8872,-0.1565) -- (2.7137,-0.2598) -- cycle;
\path[cbpmesh=15.1] (-2.0149,-0.0108) -- (-1.8506,-0.0821) -- (-1.9221,-0.1527) -- (-2.0927,-0.0786) -- cycle;
\path[cbpmesh=18.6] (1.8506,-0.0821) -- (2.0149,-0.0108) -- (2.0927,-0.0786) -- (1.9221,-0.1527) -- cycle;
\path[cbpmesh=16.4] (-2.4025,-0.1012) -- (-2.2341,-0.1867) -- (-2.2984,-0.2296) -- (-2.4717,-0.1416) -- cycle;
\path[cbpmesh=18.4] (2.2341,-0.1867) -- (2.4025,-0.1012) -- (2.4717,-0.1416) -- (2.2984,-0.2296) -- cycle;
\path[cbpmesh=19.3] (-0.8183,0.1228) -- (-0.6608,0.0955) -- (-0.6953,0.0531) -- (-0.8611,0.0818) -- cycle;
\path[cbpmesh=18.9] (0.6608,0.0955) -- (0.8183,0.1228) -- (0.8611,0.0818) -- (0.6953,0.0531) -- cycle;
\path[cbpmesh=14.9] (-1.6934,0.0821) -- (-1.5326,0.0226) -- (-1.6058,-0.0629) -- (-1.7743,-0.0005) -- cycle;
\path[cbpmesh=18.3] (1.5326,0.0226) -- (1.6934,0.0821) -- (1.7743,-0.0005) -- (1.6058,-0.0629) -- cycle;
\path[cbpmesh=14.4] (-3.4447,-0.2525) -- (-3.2685,-0.3763) -- (-3.3234,-0.4306) -- (-3.5026,-0.3047) -- cycle;
\path[cbpmesh=21.6] (3.2685,-0.3763) -- (3.4447,-0.2525) -- (3.5026,-0.3047) -- (3.3234,-0.4306) -- cycle;
\path[cbpmesh=20.9] (-0.3178,0.0900) -- (-0.1594,0.0812) -- (-0.1678,0.0495) -- (-0.3345,0.0587) -- cycle;
\path[cbpmesh=20.5] (0.1594,0.0812) -- (0.3178,0.0900) -- (0.3345,0.0587) -- (0.1678,0.0495) -- cycle;
\path[cbpmesh=20.8] (-0.1594,0.0812) -- (-0.0000,0.0783) -- (-0.0000,0.0464) -- (-0.1678,0.0495) -- cycle;
\path[cbpmesh=20.7] (-0.0000,0.0783) -- (0.1594,0.0812) -- (0.1678,0.0495) -- (-0.0000,0.0464) -- cycle;
\path[cbpmesh=14.4] (-3.5026,-0.3047) -- (-3.3234,-0.4306) -- (-3.3798,-0.4944) -- (-3.5620,-0.3664) -- cycle;
\path[cbpmesh=22.6] (3.3234,-0.4306) -- (3.5026,-0.3047) -- (3.5620,-0.3664) -- (3.3798,-0.4944) -- cycle;
\path[cbpmesh=16.8] (-3.6877,-0.5228) -- (-3.4990,-0.6553) -- (-3.5654,-0.8016) -- (-3.7577,-0.6665) -- cycle;
\path[cbpmesh=28.9] (3.4990,-0.6553) -- (3.6877,-0.5228) -- (3.7577,-0.6665) -- (3.5654,-0.8016) -- cycle;
\path[cbpmesh=14.6] (-3.5620,-0.3664) -- (-3.3798,-0.4944) -- (-3.4381,-0.5689) -- (-3.6235,-0.4386) -- cycle;
\path[cbpmesh=23.7] (3.3798,-0.4944) -- (3.5620,-0.3664) -- (3.6235,-0.4386) -- (3.4381,-0.5689) -- cycle;
\path[cbpmesh=14.8] (-3.6235,-0.4386) -- (-3.4381,-0.5689) -- (-3.4990,-0.6553) -- (-3.6877,-0.5228) -- cycle;
\path[cbpmesh=24.7] (3.4381,-0.5689) -- (3.6235,-0.4386) -- (3.6877,-0.5228) -- (3.4990,-0.6553) -- cycle;
\path[cbpmesh=17.6] (-2.8872,-0.1565) -- (-2.7137,-0.2598) -- (-2.7678,-0.2848) -- (-2.9447,-0.1794) -- cycle;
\path[cbpmesh=18.3] (2.7137,-0.2598) -- (2.8872,-0.1565) -- (2.9447,-0.1794) -- (2.7678,-0.2848) -- cycle;
\path[cbpmesh=17.8] (-1.0215,0.1166) -- (-0.8611,0.0818) -- (-0.9054,0.0281) -- (-1.0741,0.0646) -- cycle;
\path[cbpmesh=18.1] (0.8611,0.0818) -- (1.0215,0.1166) -- (1.0741,0.0646) -- (0.9054,0.0281) -- cycle;
\path[cbpmesh=38] (-3.0142,-1.9901) -- (-2.7684,-2.0968) -- (-2.7142,-2.1301) -- (-2.9551,-2.0255) -- cycle;
\draw[draw=cbpInk!70,line width=.45pt] (-2.7142,-2.1301) -- (-2.9551,-2.0255);
\path[cbpmesh=38] (2.7684,-2.0968) -- (3.0142,-1.9901) -- (2.9551,-2.0255) -- (2.7142,-2.1301) -- cycle;
\draw[draw=cbpInk!70,line width=.45pt] (2.9551,-2.0255) -- (2.7142,-2.1301);
\path[cbpmesh=15.6] (-1.4608,0.0955) -- (-1.2985,0.0446) -- (-1.3624,-0.0308) -- (-1.5326,0.0226) -- cycle;
\path[cbpmesh=17.8] (1.2985,0.0446) -- (1.4608,0.0955) -- (1.5326,0.0226) -- (1.3624,-0.0308) -- cycle;
\path[cbpmesh=16.9] (-2.4717,-0.1416) -- (-2.2984,-0.2296) -- (-2.3591,-0.2667) -- (-2.5370,-0.1764) -- cycle;
\path[cbpmesh=18.3] (2.2984,-0.2296) -- (2.4717,-0.1416) -- (2.5370,-0.1764) -- (2.3591,-0.2667) -- cycle;
\path[cbpmesh=26.8] (-3.6689,-1.2966) -- (-3.4484,-1.4279) -- (-3.4394,-1.5280) -- (-3.6592,-1.3971) -- cycle;
\path[cbpmesh=38] (3.4484,-1.4279) -- (3.6689,-1.2966) -- (3.6592,-1.3971) -- (3.4394,-1.5280) -- cycle;
\path[cbpmesh=19.3] (-0.6608,0.0955) -- (-0.4992,0.0741) -- (-0.5253,0.0306) -- (-0.6953,0.0531) -- cycle;
\path[cbpmesh=19.0] (0.4992,0.0741) -- (0.6608,0.0955) -- (0.6953,0.0531) -- (0.5253,0.0306) -- cycle;
\path[cbpmesh=16.6] (-1.2362,0.1072) -- (-1.0741,0.0646) -- (-1.1283,-0.0002) -- (-1.2985,0.0446) -- cycle;
\path[cbpmesh=17.7] (1.0741,0.0646) -- (1.2362,0.1072) -- (1.2985,0.0446) -- (1.1283,-0.0002) -- cycle;
\path[cbpmesh=15.5] (-2.0927,-0.0786) -- (-1.9221,-0.1527) -- (-1.9891,-0.2135) -- (-2.1657,-0.1368) -- cycle;
\path[cbpmesh=18.4] (1.9221,-0.1527) -- (2.0927,-0.0786) -- (2.1657,-0.1368) -- (1.9891,-0.2135) -- cycle;
\path[cbpmesh=14.7] (-1.7743,-0.0005) -- (-1.6058,-0.0629) -- (-1.6749,-0.1471) -- (-1.8506,-0.0821) -- cycle;
\path[cbpmesh=18.4] (1.6058,-0.0629) -- (1.7743,-0.0005) -- (1.8506,-0.0821) -- (1.6749,-0.1471) -- cycle;
\path[cbpmesh=17.3] (-2.9447,-0.1794) -- (-2.7678,-0.2848) -- (-2.8201,-0.3104) -- (-3.0004,-0.2031) -- cycle;
\path[cbpmesh=18.3] (2.7678,-0.2848) -- (2.9447,-0.1794) -- (3.0004,-0.2031) -- (2.8201,-0.3104) -- cycle;
\path[cbpmesh=35.6] (-3.3879,-1.6979) -- (-3.1504,-1.8185) -- (-3.1127,-1.8865) -- (-3.3473,-1.7673) -- cycle;
\path[cbpmesh=38] (3.1504,-1.8185) -- (3.3879,-1.6979) -- (3.3473,-1.7673) -- (3.1127,-1.8865) -- cycle;
\path[cbpmesh=19.3] (-0.4992,0.0741) -- (-0.3345,0.0587) -- (-0.3520,0.0144) -- (-0.5253,0.0306) -- cycle;
\path[cbpmesh=19.0] (0.3345,0.0587) -- (0.4992,0.0741) -- (0.5253,0.0306) -- (0.3520,0.0144) -- cycle;
\path[cbpmesh=17.3] (-2.5370,-0.1764) -- (-2.3591,-0.2667) -- (-2.4166,-0.2994) -- (-2.5988,-0.2069) -- cycle;
\path[cbpmesh=18.3] (2.3591,-0.2667) -- (2.5370,-0.1764) -- (2.5988,-0.2069) -- (2.4166,-0.2994) -- cycle;
\path[cbpmesh=17.9] (-0.8611,0.0818) -- (-0.6953,0.0531) -- (-0.7311,-0.0021) -- (-0.9054,0.0281) -- cycle;
\path[cbpmesh=18.1] (0.6953,0.0531) -- (0.8611,0.0818) -- (0.9054,0.0281) -- (0.7311,-0.0021) -- cycle;
\path[cbpmesh=16.8] (-3.0004,-0.2031) -- (-2.8201,-0.3104) -- (-2.8712,-0.3381) -- (-3.0548,-0.2288) -- cycle;
\path[cbpmesh=18.4] (2.8201,-0.3104) -- (3.0004,-0.2031) -- (3.0548,-0.2288) -- (2.8712,-0.3381) -- cycle;
\path[cbpmesh=15.9] (-2.1657,-0.1368) -- (-1.9891,-0.2135) -- (-2.0519,-0.2658) -- (-2.2341,-0.1867) -- cycle;
\path[cbpmesh=18.3] (1.9891,-0.2135) -- (2.1657,-0.1368) -- (2.2341,-0.1867) -- (2.0519,-0.2658) -- cycle;
\path[cbpmesh=24.0] (-3.6653,-1.1867) -- (-3.4451,-1.3178) -- (-3.4484,-1.4279) -- (-3.6689,-1.2966) -- cycle;
\path[cbpmesh=37.2] (3.4451,-1.3178) -- (3.6653,-1.1867) -- (3.6689,-1.2966) -- (3.4484,-1.4279) -- cycle;
\path[cbpmesh=19.1] (0.1678,0.0495) -- (0.3345,0.0587) -- (0.3520,0.0144) -- (0.1765,0.0047) -- cycle;
\path[cbpmesh=19.2] (-0.3345,0.0587) -- (-0.1678,0.0495) -- (-0.1765,0.0047) -- (-0.3520,0.0144) -- cycle;
\path[cbpmesh=14.9] (-1.5326,0.0226) -- (-1.3624,-0.0308) -- (-1.4275,-0.1188) -- (-1.6058,-0.0629) -- cycle;
\path[cbpmesh=18.0] (1.3624,-0.0308) -- (1.5326,0.0226) -- (1.6058,-0.0629) -- (1.4275,-0.1188) -- cycle;
\path[cbpmesh=19.2] (-0.1678,0.0495) -- (-0.0000,0.0464) -- (-0.0000,0.0014) -- (-0.1765,0.0047) -- cycle;
\path[cbpmesh=19.1] (-0.0000,0.0464) -- (0.1678,0.0495) -- (0.1765,0.0047) -- (-0.0000,0.0014) -- cycle;
\path[cbpmesh=38] (-3.0670,-1.9438) -- (-2.8169,-2.0523) -- (-2.7684,-2.0968) -- (-3.0142,-1.9901) -- cycle;
\path[cbpmesh=38] (2.8169,-2.0523) -- (3.0670,-1.9438) -- (3.0142,-1.9901) -- (2.7684,-2.0968) -- cycle;
\path[cbpmesh=16.6] (-1.0741,0.0646) -- (-0.9054,0.0281) -- (-0.9510,-0.0385) -- (-1.1283,-0.0002) -- cycle;
\path[cbpmesh=17.6] (0.9054,0.0281) -- (1.0741,0.0646) -- (1.1283,-0.0002) -- (0.9510,-0.0385) -- cycle;
\path[cbpmesh=15.0] (-1.8506,-0.0821) -- (-1.6749,-0.1471) -- (-1.7396,-0.2202) -- (-1.9221,-0.1527) -- cycle;
\path[cbpmesh=18.2] (1.6749,-0.1471) -- (1.8506,-0.0821) -- (1.9221,-0.1527) -- (1.7396,-0.2202) -- cycle;
\path[cbpmesh=15.6] (-1.2985,0.0446) -- (-1.1283,-0.0002) -- (-1.1837,-0.0777) -- (-1.3624,-0.0308) -- cycle;
\path[cbpmesh=17.6] (1.1283,-0.0002) -- (1.2985,0.0446) -- (1.3624,-0.0308) -- (1.1837,-0.0777) -- cycle;
\path[cbpmesh=16.2] (-3.0548,-0.2288) -- (-2.8712,-0.3381) -- (-2.9215,-0.3689) -- (-3.1083,-0.2577) -- cycle;
\path[cbpmesh=18.6] (2.8712,-0.3381) -- (3.0548,-0.2288) -- (3.1083,-0.2577) -- (2.9215,-0.3689) -- cycle;
\path[cbpmesh=38] (-2.7142,-2.1301) -- (-2.4565,-2.2255) -- (-2.4030,-2.2454) -- (-2.6550,-2.1521) -- cycle;
\draw[draw=cbpInk!70,line width=.45pt] (-2.4030,-2.2454) -- (-2.6550,-2.1521);
\draw[draw=cbpInk!70,line width=.45pt] (-2.6550,-2.1521) -- (-2.7142,-2.1301);
\path[cbpmesh=38] (2.4565,-2.2255) -- (2.7142,-2.1301) -- (2.6550,-2.1521) -- (2.4030,-2.2454) -- cycle;
\draw[draw=cbpInk!70,line width=.45pt] (2.7142,-2.1301) -- (2.6550,-2.1521);
\draw[draw=cbpInk!70,line width=.45pt] (2.6550,-2.1521) -- (2.4030,-2.2454);
\path[cbpmesh=17.6] (-2.5988,-0.2069) -- (-2.4166,-0.2994) -- (-2.4713,-0.3289) -- (-2.6576,-0.2343) -- cycle;
\path[cbpmesh=18.3] (2.4166,-0.2994) -- (2.5988,-0.2069) -- (2.6576,-0.2343) -- (2.4713,-0.3289) -- cycle;
\path[cbpmesh=17.9] (-0.6953,0.0531) -- (-0.5253,0.0306) -- (-0.5523,-0.0258) -- (-0.7311,-0.0021) -- cycle;
\path[cbpmesh=18.0] (0.5253,0.0306) -- (0.6953,0.0531) -- (0.7311,-0.0021) -- (0.5523,-0.0258) -- cycle;
\path[cbpmesh=32.6] (-3.4188,-1.6181) -- (-3.1792,-1.7398) -- (-3.1504,-1.8185) -- (-3.3879,-1.6979) -- cycle;
\path[cbpmesh=38] (3.1792,-1.7398) -- (3.4188,-1.6181) -- (3.3879,-1.6979) -- (3.1504,-1.8185) -- cycle;
\path[cbpmesh=16.4] (-2.2341,-0.1867) -- (-2.0519,-0.2658) -- (-2.1110,-0.3110) -- (-2.2984,-0.2296) -- cycle;
\path[cbpmesh=18.2] (2.0519,-0.2658) -- (2.2341,-0.1867) -- (2.2984,-0.2296) -- (2.1110,-0.3110) -- cycle;
\path[cbpmesh=21.6] (-3.6475,-1.0674) -- (-3.4284,-1.1979) -- (-3.4451,-1.3178) -- (-3.6653,-1.1867) -- cycle;
\path[cbpmesh=34.6] (3.4284,-1.1979) -- (3.6475,-1.0674) -- (3.6653,-1.1867) -- (3.4451,-1.3178) -- cycle;
\path[cbpmesh=15.6] (-3.1083,-0.2577) -- (-2.9215,-0.3689) -- (-2.9715,-0.4042) -- (-3.1615,-0.2911) -- cycle;
\path[cbpmesh=18.9] (2.9215,-0.3689) -- (3.1083,-0.2577) -- (3.1615,-0.2911) -- (2.9715,-0.4042) -- cycle;
\path[cbpmesh=14.7] (-1.6058,-0.0629) -- (-1.4275,-0.1188) -- (-1.4889,-0.2055) -- (-1.6749,-0.1471) -- cycle;
\path[cbpmesh=18.1] (1.4275,-0.1188) -- (1.6058,-0.0629) -- (1.6749,-0.1471) -- (1.4889,-0.2055) -- cycle;
\path[cbpmesh=17.7] (-2.6576,-0.2343) -- (-2.4713,-0.3289) -- (-2.5235,-0.3564) -- (-2.7137,-0.2598) -- cycle;
\path[cbpmesh=18.3] (2.4713,-0.3289) -- (2.6576,-0.2343) -- (2.7137,-0.2598) -- (2.5235,-0.3564) -- cycle;
\path[cbpmesh=17.9] (-0.5253,0.0306) -- (-0.3520,0.0144) -- (-0.3701,-0.0428) -- (-0.5523,-0.0258) -- cycle;
\path[cbpmesh=18.0] (0.3520,0.0144) -- (0.5253,0.0306) -- (0.5523,-0.0258) -- (0.3701,-0.0428) -- cycle;
\path[cbpmesh=15.5] (-1.9221,-0.1527) -- (-1.7396,-0.2202) -- (-1.8002,-0.2833) -- (-1.9891,-0.2135) -- cycle;
\path[cbpmesh=18.1] (1.7396,-0.2202) -- (1.9221,-0.1527) -- (1.9891,-0.2135) -- (1.8002,-0.2833) -- cycle;
\path[cbpmesh=16.7] (-0.9054,0.0281) -- (-0.7311,-0.0021) -- (-0.7680,-0.0703) -- (-0.9510,-0.0385) -- cycle;
\path[cbpmesh=17.5] (0.7311,-0.0021) -- (0.9054,0.0281) -- (0.9510,-0.0385) -- (0.7680,-0.0703) -- cycle;
\path[cbpmesh=14.9] (-3.1615,-0.2911) -- (-2.9715,-0.4042) -- (-3.0215,-0.4452) -- (-3.2147,-0.3302) -- cycle;
\path[cbpmesh=19.4] (2.9715,-0.4042) -- (3.1615,-0.2911) -- (3.2147,-0.3302) -- (3.0215,-0.4452) -- cycle;
\path[cbpmesh=19.5] (-3.6146,-0.9390) -- (-3.3974,-1.0683) -- (-3.4284,-1.1979) -- (-3.6475,-1.0674) -- cycle;
\path[cbpmesh=32.1] (3.3974,-1.0683) -- (3.6146,-0.9390) -- (3.6475,-1.0674) -- (3.4284,-1.1979) -- cycle;
\path[cbpmesh=16.9] (-2.2984,-0.2296) -- (-2.1110,-0.3110) -- (-2.1668,-0.3502) -- (-2.3591,-0.2667) -- cycle;
\path[cbpmesh=18.2] (2.1110,-0.3110) -- (2.2984,-0.2296) -- (2.3591,-0.2667) -- (2.1668,-0.3502) -- cycle;
\path[cbpmesh=14.9] (-1.3624,-0.0308) -- (-1.1837,-0.0777) -- (-1.2403,-0.1680) -- (-1.4275,-0.1188) -- cycle;
\path[cbpmesh=17.6] (1.1837,-0.0777) -- (1.3624,-0.0308) -- (1.4275,-0.1188) -- (1.2403,-0.1680) -- cycle;
\path[cbpmesh=17.9] (-0.3520,0.0144) -- (-0.1765,0.0047) -- (-0.1856,-0.0531) -- (-0.3701,-0.0428) -- cycle;
\path[cbpmesh=18.0] (0.1765,0.0047) -- (0.3520,0.0144) -- (0.3701,-0.0428) -- (0.1856,-0.0531) -- cycle;
\path[cbpmesh=38] (-3.1127,-1.8865) -- (-2.8589,-1.9967) -- (-2.8169,-2.0523) -- (-3.0670,-1.9438) -- cycle;
\path[cbpmesh=38] (2.8589,-1.9967) -- (3.1127,-1.8865) -- (3.0670,-1.9438) -- (2.8169,-2.0523) -- cycle;
\path[cbpmesh=15.7] (-1.1283,-0.0002) -- (-0.9510,-0.0385) -- (-0.9978,-0.1180) -- (-1.1837,-0.0777) -- cycle;
\path[cbpmesh=17.4] (0.9510,-0.0385) -- (1.1283,-0.0002) -- (1.1837,-0.0777) -- (0.9978,-0.1180) -- cycle;
\path[cbpmesh=14.5] (-3.2147,-0.3302) -- (-3.0215,-0.4452) -- (-3.0721,-0.4933) -- (-3.2685,-0.3763) -- cycle;
\path[cbpmesh=20.1] (3.0215,-0.4452) -- (3.2147,-0.3302) -- (3.2685,-0.3763) -- (3.0721,-0.4933) -- cycle;
\path[cbpmesh=29.5] (-3.4394,-1.5280) -- (-3.1983,-1.6505) -- (-3.1792,-1.7398) -- (-3.4188,-1.6181) -- cycle;
\path[cbpmesh=38] (3.1983,-1.6505) -- (3.4394,-1.5280) -- (3.4188,-1.6181) -- (3.1792,-1.7398) -- cycle;
\path[cbpmesh=17.9] (-0.1765,0.0047) -- (-0.0000,0.0014) -- (-0.0000,-0.0565) -- (-0.1856,-0.0531) -- cycle;
\path[cbpmesh=17.9] (-0.0000,0.0014) -- (0.1765,0.0047) -- (0.1856,-0.0531) -- (-0.0000,-0.0565) -- cycle;
\path[cbpmesh=17.6] (-2.7137,-0.2598) -- (-2.5235,-0.3564) -- (-2.5738,-0.3833) -- (-2.7678,-0.2848) -- cycle;
\path[cbpmesh=18.3] (2.5235,-0.3564) -- (2.7137,-0.2598) -- (2.7678,-0.2848) -- (2.5738,-0.3833) -- cycle;
\path[cbpmesh=38] (-2.7684,-2.0968) -- (-2.5055,-2.1941) -- (-2.4565,-2.2255) -- (-2.7142,-2.1301) -- cycle;
\path[cbpmesh=38] (2.5055,-2.1941) -- (2.7684,-2.0968) -- (2.7142,-2.1301) -- (2.4565,-2.2255) -- cycle;
\path[cbpmesh=14.3] (-3.2685,-0.3763) -- (-3.0721,-0.4933) -- (-3.1237,-0.5495) -- (-3.3234,-0.4306) -- cycle;
\path[cbpmesh=21.0] (3.0721,-0.4933) -- (3.2685,-0.3763) -- (3.3234,-0.4306) -- (3.1237,-0.5495) -- cycle;
\path[cbpmesh=17.8] (-3.5654,-0.8016) -- (-3.3512,-0.9291) -- (-3.3974,-1.0683) -- (-3.6146,-0.9390) -- cycle;
\path[cbpmesh=29.9] (3.3512,-0.9291) -- (3.5654,-0.8016) -- (3.6146,-0.9390) -- (3.3974,-1.0683) -- cycle;
\path[cbpmesh=16.7] (-0.7311,-0.0021) -- (-0.5523,-0.0258) -- (-0.5802,-0.0951) -- (-0.7680,-0.0703) -- cycle;
\path[cbpmesh=17.4] (0.5523,-0.0258) -- (0.7311,-0.0021) -- (0.7680,-0.0703) -- (0.5802,-0.0951) -- cycle;
\path[cbpmesh=15.9] (-1.9891,-0.2135) -- (-1.8002,-0.2833) -- (-1.8571,-0.3379) -- (-2.0519,-0.2658) -- cycle;
\path[cbpmesh=18.0] (1.8002,-0.2833) -- (1.9891,-0.2135) -- (2.0519,-0.2658) -- (1.8571,-0.3379) -- cycle;
\path[cbpmesh=15.1] (-1.6749,-0.1471) -- (-1.4889,-0.2055) -- (-1.5464,-0.2808) -- (-1.7396,-0.2202) -- cycle;
\path[cbpmesh=17.9] (1.4889,-0.2055) -- (1.6749,-0.1471) -- (1.7396,-0.2202) -- (1.5464,-0.2808) -- cycle;
\path[cbpmesh=14.2] (-3.3234,-0.4306) -- (-3.1237,-0.5495) -- (-3.1767,-0.6153) -- (-3.3798,-0.4944) -- cycle;
\path[cbpmesh=22.0] (3.1237,-0.5495) -- (3.3234,-0.4306) -- (3.3798,-0.4944) -- (3.1767,-0.6153) -- cycle;
\path[cbpmesh=17.3] (-2.3591,-0.2667) -- (-2.1668,-0.3502) -- (-2.2196,-0.3849) -- (-2.4166,-0.2994) -- cycle;
\path[cbpmesh=18.2] (2.1668,-0.3502) -- (2.3591,-0.2667) -- (2.4166,-0.2994) -- (2.2196,-0.3849) -- cycle;
\path[cbpmesh=14.4] (-3.3798,-0.4944) -- (-3.1767,-0.6153) -- (-3.2316,-0.6919) -- (-3.4381,-0.5689) -- cycle;
\path[cbpmesh=23.0] (3.1767,-0.6153) -- (3.3798,-0.4944) -- (3.4381,-0.5689) -- (3.2316,-0.6919) -- cycle;
\path[cbpmesh=17.3] (-2.7678,-0.2848) -- (-2.5738,-0.3833) -- (-2.6225,-0.4108) -- (-2.8201,-0.3104) -- cycle;
\path[cbpmesh=18.3] (2.5738,-0.3833) -- (2.7678,-0.2848) -- (2.8201,-0.3104) -- (2.6225,-0.4108) -- cycle;
\path[cbpmesh=16.6] (-3.4990,-0.6553) -- (-3.2888,-0.7805) -- (-3.3512,-0.9291) -- (-3.5654,-0.8016) -- cycle;
\path[cbpmesh=28.0] (3.2888,-0.7805) -- (3.4990,-0.6553) -- (3.5654,-0.8016) -- (3.3512,-0.9291) -- cycle;
\path[cbpmesh=14.6] (-3.4381,-0.5689) -- (-3.2316,-0.6919) -- (-3.2888,-0.7805) -- (-3.4990,-0.6553) -- cycle;
\path[cbpmesh=23.9] (3.2316,-0.6919) -- (3.4381,-0.5689) -- (3.4990,-0.6553) -- (3.2888,-0.7805) -- cycle;
\path[cbpmesh=15.8] (-0.9510,-0.0385) -- (-0.7680,-0.0703) -- (-0.8057,-0.1512) -- (-0.9978,-0.1180) -- cycle;
\path[cbpmesh=17.1] (0.7680,-0.0703) -- (0.9510,-0.0385) -- (0.9978,-0.1180) -- (0.8057,-0.1512) -- cycle;
\path[cbpmesh=16.8] (-0.5523,-0.0258) -- (-0.3701,-0.0428) -- (-0.3888,-0.1130) -- (-0.5802,-0.0951) -- cycle;
\path[cbpmesh=17.2] (0.3701,-0.0428) -- (0.5523,-0.0258) -- (0.5802,-0.0951) -- (0.3888,-0.1130) -- cycle;
\path[cbpmesh=14.8] (-1.4275,-0.1188) -- (-1.2403,-0.1680) -- (-1.2937,-0.2568) -- (-1.4889,-0.2055) -- cycle;
\path[cbpmesh=17.7] (1.2403,-0.1680) -- (1.4275,-0.1188) -- (1.4889,-0.2055) -- (1.2937,-0.2568) -- cycle;
\path[cbpmesh=26.5] (-3.4484,-1.4279) -- (-3.2067,-1.5507) -- (-3.1983,-1.6505) -- (-3.4394,-1.5280) -- cycle;
\path[cbpmesh=38] (3.2067,-1.5507) -- (3.4484,-1.4279) -- (3.4394,-1.5280) -- (3.1983,-1.6505) -- cycle;
\path[cbpmesh=35.5] (-3.1504,-1.8185) -- (-2.8935,-1.9300) -- (-2.8589,-1.9967) -- (-3.1127,-1.8865) -- cycle;
\path[cbpmesh=38] (2.8935,-1.9300) -- (3.1504,-1.8185) -- (3.1127,-1.8865) -- (2.8589,-1.9967) -- cycle;
\path[cbpmesh=15.0] (-1.1837,-0.0777) -- (-0.9978,-0.1180) -- (-1.0455,-0.2102) -- (-1.2403,-0.1680) -- cycle;
\path[cbpmesh=17.3] (0.9978,-0.1180) -- (1.1837,-0.0777) -- (1.2403,-0.1680) -- (1.0455,-0.2102) -- cycle;
\path[cbpmesh=16.8] (-2.8201,-0.3104) -- (-2.6225,-0.4108) -- (-2.6700,-0.4403) -- (-2.8712,-0.3381) -- cycle;
\path[cbpmesh=18.3] (2.6225,-0.4108) -- (2.8201,-0.3104) -- (2.8712,-0.3381) -- (2.6700,-0.4403) -- cycle;
\path[cbpmesh=16.4] (-2.0519,-0.2658) -- (-1.8571,-0.3379) -- (-1.9106,-0.3851) -- (-2.1110,-0.3110) -- cycle;
\path[cbpmesh=18.0] (1.8571,-0.3379) -- (2.0519,-0.2658) -- (2.1110,-0.3110) -- (1.9106,-0.3851) -- cycle;
\path[cbpmesh=38] (-2.4565,-2.2255) -- (-2.1836,-2.3110) -- (-2.1361,-2.3291) -- (-2.4030,-2.2454) -- cycle;
\draw[draw=cbpInk!70,line width=.45pt] (-2.1361,-2.3291) -- (-2.4030,-2.2454);
\path[cbpmesh=38] (2.1836,-2.3110) -- (2.4565,-2.2255) -- (2.4030,-2.2454) -- (2.1361,-2.3291) -- cycle;
\draw[draw=cbpInk!70,line width=.45pt] (2.4030,-2.2454) -- (2.1361,-2.3291);
\path[cbpmesh=16.8] (-0.3701,-0.0428) -- (-0.1856,-0.0531) -- (-0.1950,-0.1238) -- (-0.3888,-0.1130) -- cycle;
\path[cbpmesh=17.1] (0.1856,-0.0531) -- (0.3701,-0.0428) -- (0.3888,-0.1130) -- (0.1950,-0.1238) -- cycle;
\path[cbpmesh=17.6] (-2.4166,-0.2994) -- (-2.2196,-0.3849) -- (-2.2698,-0.4163) -- (-2.4713,-0.3289) -- cycle;
\path[cbpmesh=18.2] (2.2196,-0.3849) -- (2.4166,-0.2994) -- (2.4713,-0.3289) -- (2.2698,-0.4163) -- cycle;
\path[cbpmesh=15.5] (-1.7396,-0.2202) -- (-1.5464,-0.2808) -- (-1.6003,-0.3461) -- (-1.8002,-0.2833) -- cycle;
\path[cbpmesh=17.8] (1.5464,-0.2808) -- (1.7396,-0.2202) -- (1.8002,-0.2833) -- (1.6003,-0.3461) -- cycle;
\path[cbpmesh=38] (-2.8169,-2.0523) -- (-2.5494,-2.1513) -- (-2.5055,-2.1941) -- (-2.7684,-2.0968) -- cycle;
\path[cbpmesh=38] (2.5494,-2.1513) -- (2.8169,-2.0523) -- (2.7684,-2.0968) -- (2.5055,-2.1941) -- cycle;
\path[cbpmesh=16.9] (-0.1856,-0.0531) -- (-0.0000,-0.0565) -- (-0.0000,-0.1274) -- (-0.1950,-0.1238) -- cycle;
\path[cbpmesh=17.0] (-0.0000,-0.0565) -- (0.1856,-0.0531) -- (0.1950,-0.1238) -- (-0.0000,-0.1274) -- cycle;
\path[cbpmesh=15.8] (-0.7680,-0.0703) -- (-0.5802,-0.0951) -- (-0.6087,-0.1773) -- (-0.8057,-0.1512) -- cycle;
\path[cbpmesh=16.9] (0.5802,-0.0951) -- (0.7680,-0.0703) -- (0.8057,-0.1512) -- (0.6087,-0.1773) -- cycle;
\path[cbpmesh=16.2] (-2.8712,-0.3381) -- (-2.6700,-0.4403) -- (-2.7168,-0.4729) -- (-2.9215,-0.3689) -- cycle;
\path[cbpmesh=18.4] (2.6700,-0.4403) -- (2.8712,-0.3381) -- (2.9215,-0.3689) -- (2.7168,-0.4729) -- cycle;
\path[cbpmesh=23.8] (-3.4451,-1.3178) -- (-3.2036,-1.4405) -- (-3.2067,-1.5507) -- (-3.4484,-1.4279) -- cycle;
\path[cbpmesh=36.1] (3.2036,-1.4405) -- (3.4451,-1.3178) -- (3.4484,-1.4279) -- (3.2067,-1.5507) -- cycle;
\path[cbpmesh=15.1] (-1.4889,-0.2055) -- (-1.2937,-0.2568) -- (-1.3436,-0.3341) -- (-1.5464,-0.2808) -- cycle;
\path[cbpmesh=17.6] (1.2937,-0.2568) -- (1.4889,-0.2055) -- (1.5464,-0.2808) -- (1.3436,-0.3341) -- cycle;
\path[cbpmesh=17.7] (-2.4713,-0.3289) -- (-2.2698,-0.4163) -- (-2.3178,-0.4457) -- (-2.5235,-0.3564) -- cycle;
\path[cbpmesh=18.2] (2.2698,-0.4163) -- (2.4713,-0.3289) -- (2.5235,-0.3564) -- (2.3178,-0.4457) -- cycle;
\path[cbpmesh=16.9] (-2.1110,-0.3110) -- (-1.9106,-0.3851) -- (-1.9611,-0.4264) -- (-2.1668,-0.3502) -- cycle;
\path[cbpmesh=18.0] (1.9106,-0.3851) -- (2.1110,-0.3110) -- (2.1668,-0.3502) -- (1.9611,-0.4264) -- cycle;
\path[cbpmesh=15.1] (-0.9978,-0.1180) -- (-0.8057,-0.1512) -- (-0.8442,-0.2450) -- (-1.0455,-0.2102) -- cycle;
\path[cbpmesh=17.0] (0.8057,-0.1512) -- (0.9978,-0.1180) -- (1.0455,-0.2102) -- (0.8442,-0.2450) -- cycle;
\path[cbpmesh=15.5] (-2.9215,-0.3689) -- (-2.7168,-0.4729) -- (-2.7632,-0.5100) -- (-2.9715,-0.4042) -- cycle;
\path[cbpmesh=18.6] (2.7168,-0.4729) -- (2.9215,-0.3689) -- (2.9715,-0.4042) -- (2.7632,-0.5100) -- cycle;
\path[cbpmesh=32.4] (-3.1792,-1.7398) -- (-2.9200,-1.8523) -- (-2.8935,-1.9300) -- (-3.1504,-1.8185) -- cycle;
\path[cbpmesh=38] (2.9200,-1.8523) -- (3.1792,-1.7398) -- (3.1504,-1.8185) -- (2.8935,-1.9300) -- cycle;
\path[cbpmesh=14.8] (-1.2403,-0.1680) -- (-1.0455,-0.2102) -- (-1.0905,-0.3008) -- (-1.2937,-0.2568) -- cycle;
\path[cbpmesh=17.3] (1.0455,-0.2102) -- (1.2403,-0.1680) -- (1.2937,-0.2568) -- (1.0905,-0.3008) -- cycle;
\path[cbpmesh=15.9] (-1.8002,-0.2833) -- (-1.6003,-0.3461) -- (-1.6508,-0.4026) -- (-1.8571,-0.3379) -- cycle;
\path[cbpmesh=17.8] (1.6003,-0.3461) -- (1.8002,-0.2833) -- (1.8571,-0.3379) -- (1.6508,-0.4026) -- cycle;
\path[cbpmesh=15.9] (-0.5802,-0.0951) -- (-0.3888,-0.1130) -- (-0.4079,-0.1961) -- (-0.6087,-0.1773) -- cycle;
\path[cbpmesh=16.7] (0.3888,-0.1130) -- (0.5802,-0.0951) -- (0.6087,-0.1773) -- (0.4079,-0.1961) -- cycle;
\path[cbpmesh=21.4] (-3.4284,-1.1979) -- (-3.1881,-1.3200) -- (-3.2036,-1.4405) -- (-3.4451,-1.3178) -- cycle;
\path[cbpmesh=33.5] (3.1881,-1.3200) -- (3.4284,-1.1979) -- (3.4451,-1.3178) -- (3.2036,-1.4405) -- cycle;
\path[cbpmesh=14.9] (-2.9715,-0.4042) -- (-2.7632,-0.5100) -- (-2.8097,-0.5528) -- (-3.0215,-0.4452) -- cycle;
\path[cbpmesh=19.1] (2.7632,-0.5100) -- (2.9715,-0.4042) -- (3.0215,-0.4452) -- (2.8097,-0.5528) -- cycle;
\path[cbpmesh=38] (-2.5055,-2.1941) -- (-2.2273,-2.2814) -- (-2.1836,-2.3110) -- (-2.4565,-2.2255) -- cycle;
\path[cbpmesh=38] (2.2273,-2.2814) -- (2.5055,-2.1941) -- (2.4565,-2.2255) -- (2.1836,-2.3110) -- cycle;
\path[cbpmesh=17.6] (-2.5235,-0.3564) -- (-2.3178,-0.4457) -- (-2.3639,-0.4744) -- (-2.5738,-0.3833) -- cycle;
\path[cbpmesh=18.2] (2.3178,-0.4457) -- (2.5235,-0.3564) -- (2.5738,-0.3833) -- (2.3639,-0.4744) -- cycle;
\path[cbpmesh=16.1] (-0.3888,-0.1130) -- (-0.1950,-0.1238) -- (-0.2046,-0.2074) -- (-0.4079,-0.1961) -- cycle;
\path[cbpmesh=16.5] (0.1950,-0.1238) -- (0.3888,-0.1130) -- (0.4079,-0.1961) -- (0.2046,-0.2074) -- cycle;
\path[cbpmesh=38] (-2.8589,-1.9967) -- (-2.5874,-2.0971) -- (-2.5494,-2.1513) -- (-2.8169,-2.0523) -- cycle;
\path[cbpmesh=38] (2.5874,-2.0971) -- (2.8589,-1.9967) -- (2.8169,-2.0523) -- (2.5494,-2.1513) -- cycle;
\path[cbpmesh=17.3] (-2.1668,-0.3502) -- (-1.9611,-0.4264) -- (-2.0088,-0.4629) -- (-2.2196,-0.3849) -- cycle;
\path[cbpmesh=18.1] (1.9611,-0.4264) -- (2.1668,-0.3502) -- (2.2196,-0.3849) -- (2.0088,-0.4629) -- cycle;
\path[cbpmesh=16.2] (-0.1950,-0.1238) -- (-0.0000,-0.1274) -- (-0.0000,-0.2111) -- (-0.2046,-0.2074) -- cycle;
\path[cbpmesh=16.4] (-0.0000,-0.1274) -- (0.1950,-0.1238) -- (0.2046,-0.2074) -- (-0.0000,-0.2111) -- cycle;
\path[cbpmesh=15.5] (-1.5464,-0.2808) -- (-1.3436,-0.3341) -- (-1.3904,-0.4012) -- (-1.6003,-0.3461) -- cycle;
\path[cbpmesh=17.6] (1.3436,-0.3341) -- (1.5464,-0.2808) -- (1.6003,-0.3461) -- (1.3904,-0.4012) -- cycle;
\path[cbpmesh=15.2] (-0.8057,-0.1512) -- (-0.6087,-0.1773) -- (-0.6378,-0.2724) -- (-0.8442,-0.2450) -- cycle;
\path[cbpmesh=16.7] (0.6087,-0.1773) -- (0.8057,-0.1512) -- (0.8442,-0.2450) -- (0.6378,-0.2724) -- cycle;
\path[cbpmesh=14.4] (-3.0215,-0.4452) -- (-2.8097,-0.5528) -- (-2.8568,-0.6026) -- (-3.0721,-0.4933) -- cycle;
\path[cbpmesh=19.7] (2.8097,-0.5528) -- (3.0215,-0.4452) -- (3.0721,-0.4933) -- (2.8568,-0.6026) -- cycle;
\path[cbpmesh=16.4] (-1.8571,-0.3379) -- (-1.6508,-0.4026) -- (-1.6984,-0.4517) -- (-1.9106,-0.3851) -- cycle;
\path[cbpmesh=17.9] (1.6508,-0.4026) -- (1.8571,-0.3379) -- (1.9106,-0.3851) -- (1.6984,-0.4517) -- cycle;
\path[cbpmesh=19.3] (-3.3974,-1.0683) -- (-3.1593,-1.1893) -- (-3.1881,-1.3200) -- (-3.4284,-1.1979) -- cycle;
\path[cbpmesh=31.1] (3.1593,-1.1893) -- (3.3974,-1.0683) -- (3.4284,-1.1979) -- (3.1881,-1.3200) -- cycle;
\path[cbpmesh=15.2] (-1.2937,-0.2568) -- (-1.0905,-0.3008) -- (-1.1326,-0.3798) -- (-1.3436,-0.3341) -- cycle;
\path[cbpmesh=17.3] (1.0905,-0.3008) -- (1.2937,-0.2568) -- (1.3436,-0.3341) -- (1.1326,-0.3798) -- cycle;
\path[cbpmesh=29.4] (-3.1983,-1.6505) -- (-2.9375,-1.7636) -- (-2.9200,-1.8523) -- (-3.1792,-1.7398) -- cycle;
\path[cbpmesh=38] (2.9375,-1.7636) -- (3.1983,-1.6505) -- (3.1792,-1.7398) -- (2.9200,-1.8523) -- cycle;
\path[cbpmesh=17.0] (0.8442,-0.2450) -- (1.0455,-0.2102) -- (1.0905,-0.3008) -- (0.8806,-0.3372) -- cycle;
\path[cbpmesh=14.9] (-1.0455,-0.2102) -- (-0.8442,-0.2450) -- (-0.8806,-0.3372) -- (-1.0905,-0.3008) -- cycle;
\path[cbpmesh=17.3] (-2.5738,-0.3833) -- (-2.3639,-0.4744) -- (-2.4086,-0.5036) -- (-2.6225,-0.4108) -- cycle;
\path[cbpmesh=18.2] (2.3639,-0.4744) -- (2.5738,-0.3833) -- (2.6225,-0.4108) -- (2.4086,-0.5036) -- cycle;
\path[cbpmesh=14.2] (-3.0721,-0.4933) -- (-2.8568,-0.6026) -- (-2.9047,-0.6607) -- (-3.1237,-0.5495) -- cycle;
\path[cbpmesh=20.5] (2.8568,-0.6026) -- (3.0721,-0.4933) -- (3.1237,-0.5495) -- (2.9047,-0.6607) -- cycle;
\path[cbpmesh=38] (-2.1836,-2.3110) -- (-1.8973,-2.3863) -- (-1.8560,-2.4028) -- (-2.1361,-2.3291) -- cycle;
\draw[draw=cbpInk!70,line width=.45pt] (-1.8560,-2.4028) -- (-2.1361,-2.3291);
\path[cbpmesh=38] (1.8973,-2.3863) -- (2.1836,-2.3110) -- (2.1361,-2.3291) -- (1.8560,-2.4028) -- cycle;
\draw[draw=cbpInk!70,line width=.45pt] (2.1361,-2.3291) -- (1.8560,-2.4028);
\path[cbpmesh=17.6] (-2.2196,-0.3849) -- (-2.0088,-0.4629) -- (-2.0543,-0.4961) -- (-2.2698,-0.4163) -- cycle;
\path[cbpmesh=18.1] (2.0088,-0.4629) -- (2.2196,-0.3849) -- (2.2698,-0.4163) -- (2.0543,-0.4961) -- cycle;
\path[cbpmesh=15.4] (-0.6087,-0.1773) -- (-0.4079,-0.1961) -- (-0.4274,-0.2920) -- (-0.6378,-0.2724) -- cycle;
\path[cbpmesh=16.4] (0.4079,-0.1961) -- (0.6087,-0.1773) -- (0.6378,-0.2724) -- (0.4274,-0.2920) -- cycle;
\path[cbpmesh=14.1] (-3.1237,-0.5495) -- (-2.9047,-0.6607) -- (-2.9540,-0.7284) -- (-3.1767,-0.6153) -- cycle;
\path[cbpmesh=21.4] (2.9047,-0.6607) -- (3.1237,-0.5495) -- (3.1767,-0.6153) -- (2.9540,-0.7284) -- cycle;
\path[cbpmesh=17.6] (-3.3512,-0.9291) -- (-3.1163,-1.0484) -- (-3.1593,-1.1893) -- (-3.3974,-1.0683) -- cycle;
\path[cbpmesh=28.9] (3.1163,-1.0484) -- (3.3512,-0.9291) -- (3.3974,-1.0683) -- (3.1593,-1.1893) -- cycle;
\path[cbpmesh=15.9] (-1.6003,-0.3461) -- (-1.3904,-0.4012) -- (-1.4344,-0.4595) -- (-1.6508,-0.4026) -- cycle;
\path[cbpmesh=17.6] (1.3904,-0.4012) -- (1.6003,-0.3461) -- (1.6508,-0.4026) -- (1.4344,-0.4595) -- cycle;
\path[cbpmesh=16.8] (-2.6225,-0.4108) -- (-2.4086,-0.5036) -- (-2.4523,-0.5348) -- (-2.6700,-0.4403) -- cycle;
\path[cbpmesh=18.2] (2.4086,-0.5036) -- (2.6225,-0.4108) -- (2.6700,-0.4403) -- (2.4523,-0.5348) -- cycle;
\path[cbpmesh=14.2] (-3.1767,-0.6153) -- (-2.9540,-0.7284) -- (-3.0050,-0.8070) -- (-3.2316,-0.6919) -- cycle;
\path[cbpmesh=22.3] (2.9540,-0.7284) -- (3.1767,-0.6153) -- (3.2316,-0.6919) -- (3.0050,-0.8070) -- cycle;
\path[cbpmesh=35.5] (-2.8935,-1.9300) -- (-2.6188,-2.0316) -- (-2.5874,-2.0971) -- (-2.8589,-1.9967) -- cycle;
\path[cbpmesh=38] (2.6188,-2.0316) -- (2.8935,-1.9300) -- (2.8589,-1.9967) -- (2.5874,-2.0971) -- cycle;
\path[cbpmesh=16.9] (-1.9106,-0.3851) -- (-1.6984,-0.4517) -- (-1.7432,-0.4947) -- (-1.9611,-0.4264) -- cycle;
\path[cbpmesh=17.9] (1.6984,-0.4517) -- (1.9106,-0.3851) -- (1.9611,-0.4264) -- (1.7432,-0.4947) -- cycle;
\path[cbpmesh=15.5] (-0.4079,-0.1961) -- (-0.2046,-0.2074) -- (-0.2144,-0.3039) -- (-0.4274,-0.2920) -- cycle;
\path[cbpmesh=16.2] (0.2046,-0.2074) -- (0.4079,-0.1961) -- (0.4274,-0.2920) -- (0.2144,-0.3039) -- cycle;
\path[cbpmesh=38] (-2.5494,-2.1513) -- (-2.2663,-2.2401) -- (-2.2273,-2.2814) -- (-2.5055,-2.1941) -- cycle;
\path[cbpmesh=38] (2.2663,-2.2401) -- (2.5494,-2.1513) -- (2.5055,-2.1941) -- (2.2273,-2.2814) -- cycle;
\path[cbpmesh=16.4] (-3.2888,-0.7805) -- (-3.0582,-0.8976) -- (-3.1163,-1.0484) -- (-3.3512,-0.9291) -- cycle;
\path[cbpmesh=27.0] (3.0582,-0.8976) -- (3.2888,-0.7805) -- (3.3512,-0.9291) -- (3.1163,-1.0484) -- cycle;
\path[cbpmesh=14.5] (-3.2316,-0.6919) -- (-3.0050,-0.8070) -- (-3.0582,-0.8976) -- (-3.2888,-0.7805) -- cycle;
\path[cbpmesh=23.2] (3.0050,-0.8070) -- (3.2316,-0.6919) -- (3.2888,-0.7805) -- (3.0582,-0.8976) -- cycle;
\path[cbpmesh=15.1] (-0.8442,-0.2450) -- (-0.6378,-0.2724) -- (-0.6652,-0.3657) -- (-0.8806,-0.3372) -- cycle;
\path[cbpmesh=16.7] (0.6378,-0.2724) -- (0.8442,-0.2450) -- (0.8806,-0.3372) -- (0.6652,-0.3657) -- cycle;
\path[cbpmesh=15.7] (-0.2046,-0.2074) -- (-0.0000,-0.2111) -- (-0.0000,-0.3078) -- (-0.2144,-0.3039) -- cycle;
\path[cbpmesh=15.9] (-0.0000,-0.2111) -- (0.2046,-0.2074) -- (0.2144,-0.3039) -- (-0.0000,-0.3078) -- cycle;
\path[cbpmesh=26.4] (-3.2067,-1.5507) -- (-2.9452,-1.6641) -- (-2.9375,-1.7636) -- (-3.1983,-1.6505) -- cycle;
\path[cbpmesh=37.8] (2.9452,-1.6641) -- (3.2067,-1.5507) -- (3.1983,-1.6505) -- (2.9375,-1.7636) -- cycle;
\path[cbpmesh=15.5] (-1.3436,-0.3341) -- (-1.1326,-0.3798) -- (-1.1720,-0.4485) -- (-1.3904,-0.4012) -- cycle;
\path[cbpmesh=17.3] (1.1326,-0.3798) -- (1.3436,-0.3341) -- (1.3904,-0.4012) -- (1.1720,-0.4485) -- cycle;
\path[cbpmesh=17.7] (-2.2698,-0.4163) -- (-2.0543,-0.4961) -- (-2.0977,-0.5272) -- (-2.3178,-0.4457) -- cycle;
\path[cbpmesh=18.2] (2.0543,-0.4961) -- (2.2698,-0.4163) -- (2.3178,-0.4457) -- (2.0977,-0.5272) -- cycle;
\path[cbpmesh=15.2] (-1.0905,-0.3008) -- (-0.8806,-0.3372) -- (-0.9146,-0.4176) -- (-1.1326,-0.3798) -- cycle;
\path[cbpmesh=17.0] (0.8806,-0.3372) -- (1.0905,-0.3008) -- (1.1326,-0.3798) -- (0.9146,-0.4176) -- cycle;
\path[cbpmesh=16.2] (-2.6700,-0.4403) -- (-2.4523,-0.5348) -- (-2.4952,-0.5690) -- (-2.7168,-0.4729) -- cycle;
\path[cbpmesh=18.2] (2.4523,-0.5348) -- (2.6700,-0.4403) -- (2.7168,-0.4729) -- (2.4952,-0.5690) -- cycle;
\path[cbpmesh=16.4] (-1.6508,-0.4026) -- (-1.4344,-0.4595) -- (-1.4757,-0.5102) -- (-1.6984,-0.4517) -- cycle;
\path[cbpmesh=17.7] (1.4344,-0.4595) -- (1.6508,-0.4026) -- (1.6984,-0.4517) -- (1.4757,-0.5102) -- cycle;
\path[cbpmesh=15.2] (-0.6378,-0.2724) -- (-0.4274,-0.2920) -- (-0.4458,-0.3862) -- (-0.6652,-0.3657) -- cycle;
\path[cbpmesh=16.4] (0.4274,-0.2920) -- (0.6378,-0.2724) -- (0.6652,-0.3657) -- (0.4458,-0.3862) -- cycle;
\path[cbpmesh=17.3] (-1.9611,-0.4264) -- (-1.7432,-0.4947) -- (-1.7857,-0.5329) -- (-2.0088,-0.4629) -- cycle;
\path[cbpmesh=18.0] (1.7432,-0.4947) -- (1.9611,-0.4264) -- (2.0088,-0.4629) -- (1.7857,-0.5329) -- cycle;
\path[cbpmesh=38] (-2.2273,-2.2814) -- (-1.9352,-2.3581) -- (-1.8973,-2.3863) -- (-2.1836,-2.3110) -- cycle;
\path[cbpmesh=38] (1.9352,-2.3581) -- (2.2273,-2.2814) -- (2.1836,-2.3110) -- (1.8973,-2.3863) -- cycle;
\path[cbpmesh=15.5] (-2.7168,-0.4729) -- (-2.4952,-0.5690) -- (-2.5379,-0.6078) -- (-2.7632,-0.5100) -- cycle;
\path[cbpmesh=18.4] (2.4952,-0.5690) -- (2.7168,-0.4729) -- (2.7632,-0.5100) -- (2.5379,-0.6078) -- cycle;
\path[cbpmesh=17.6] (-2.3178,-0.4457) -- (-2.0977,-0.5272) -- (-2.1395,-0.5575) -- (-2.3639,-0.4744) -- cycle;
\path[cbpmesh=18.1] (2.0977,-0.5272) -- (2.3178,-0.4457) -- (2.3639,-0.4744) -- (2.1395,-0.5575) -- cycle;
\path[cbpmesh=23.7] (-3.2036,-1.4405) -- (-2.9424,-1.5538) -- (-2.9452,-1.6641) -- (-3.2067,-1.5507) -- cycle;
\path[cbpmesh=35.1] (2.9424,-1.5538) -- (3.2036,-1.4405) -- (3.2067,-1.5507) -- (2.9452,-1.6641) -- cycle;
\path[cbpmesh=32.4] (-2.9200,-1.8523) -- (-2.6427,-1.9549) -- (-2.6188,-2.0316) -- (-2.8935,-1.9300) -- cycle;
\path[cbpmesh=38] (2.6427,-1.9549) -- (2.9200,-1.8523) -- (2.8935,-1.9300) -- (2.6188,-2.0316) -- cycle;
\path[cbpmesh=15.4] (-0.4274,-0.2920) -- (-0.2144,-0.3039) -- (-0.2236,-0.3985) -- (-0.4458,-0.3862) -- cycle;
\path[cbpmesh=16.1] (0.2144,-0.3039) -- (0.4274,-0.2920) -- (0.4458,-0.3862) -- (0.2236,-0.3985) -- cycle;
\path[cbpmesh=15.4] (-0.8806,-0.3372) -- (-0.6652,-0.3657) -- (-0.6909,-0.4472) -- (-0.9146,-0.4176) -- cycle;
\path[cbpmesh=16.8] (0.6652,-0.3657) -- (0.8806,-0.3372) -- (0.9146,-0.4176) -- (0.6909,-0.4472) -- cycle;
\path[cbpmesh=16.0] (-1.3904,-0.4012) -- (-1.1720,-0.4485) -- (-1.2091,-0.5083) -- (-1.4344,-0.4595) -- cycle;
\path[cbpmesh=17.4] (1.1720,-0.4485) -- (1.3904,-0.4012) -- (1.4344,-0.4595) -- (1.2091,-0.5083) -- cycle;
\path[cbpmesh=38] (-1.8973,-2.3863) -- (-1.5993,-2.4508) -- (-1.5645,-2.4659) -- (-1.8560,-2.4028) -- cycle;
\draw[draw=cbpInk!70,line width=.45pt] (-1.5645,-2.4659) -- (-1.8560,-2.4028);
\path[cbpmesh=38] (1.5993,-2.4508) -- (1.8973,-2.3863) -- (1.8560,-2.4028) -- (1.5645,-2.4659) -- cycle;
\draw[draw=cbpInk!70,line width=.45pt] (1.8560,-2.4028) -- (1.5645,-2.4659);
\path[cbpmesh=38] (-2.5874,-2.0971) -- (-2.3000,-2.1873) -- (-2.2663,-2.2401) -- (-2.5494,-2.1513) -- cycle;
\path[cbpmesh=38] (2.3000,-2.1873) -- (2.5874,-2.0971) -- (2.5494,-2.1513) -- (2.2663,-2.2401) -- cycle;
\path[cbpmesh=14.8] (-2.7632,-0.5100) -- (-2.5379,-0.6078) -- (-2.5806,-0.6522) -- (-2.8097,-0.5528) -- cycle;
\path[cbpmesh=18.7] (2.5379,-0.6078) -- (2.7632,-0.5100) -- (2.8097,-0.5528) -- (2.5806,-0.6522) -- cycle;
\path[cbpmesh=15.6] (-1.1326,-0.3798) -- (-0.9146,-0.4176) -- (-0.9464,-0.4876) -- (-1.1720,-0.4485) -- cycle;
\path[cbpmesh=17.1] (0.9146,-0.4176) -- (1.1326,-0.3798) -- (1.1720,-0.4485) -- (0.9464,-0.4876) -- cycle;
\path[cbpmesh=15.6] (-0.2144,-0.3039) -- (-0.0000,-0.3078) -- (-0.0000,-0.4026) -- (-0.2236,-0.3985) -- cycle;
\path[cbpmesh=15.9] (-0.0000,-0.3078) -- (0.2144,-0.3039) -- (0.2236,-0.3985) -- (-0.0000,-0.4026) -- cycle;
\path[cbpmesh=17.6] (-2.0088,-0.4629) -- (-1.7857,-0.5329) -- (-1.8261,-0.5677) -- (-2.0543,-0.4961) -- cycle;
\path[cbpmesh=18.1] (1.7857,-0.5329) -- (2.0088,-0.4629) -- (2.0543,-0.4961) -- (1.8261,-0.5677) -- cycle;
\path[cbpmesh=17.3] (-2.3639,-0.4744) -- (-2.1395,-0.5575) -- (-2.1799,-0.5883) -- (-2.4086,-0.5036) -- cycle;
\path[cbpmesh=18.1] (2.1395,-0.5575) -- (2.3639,-0.4744) -- (2.4086,-0.5036) -- (2.1799,-0.5883) -- cycle;
\path[cbpmesh=16.9] (-1.6984,-0.4517) -- (-1.4757,-0.5102) -- (-1.5147,-0.5548) -- (-1.7432,-0.4947) -- cycle;
\path[cbpmesh=17.8] (1.4757,-0.5102) -- (1.6984,-0.4517) -- (1.7432,-0.4947) -- (1.5147,-0.5548) -- cycle;
\path[cbpmesh=21.2] (-3.1881,-1.3200) -- (-2.9281,-1.4328) -- (-2.9424,-1.5538) -- (-3.2036,-1.4405) -- cycle;
\path[cbpmesh=32.5] (2.9281,-1.4328) -- (3.1881,-1.3200) -- (3.2036,-1.4405) -- (2.9424,-1.5538) -- cycle;
\path[cbpmesh=14.4] (-2.8097,-0.5528) -- (-2.5806,-0.6522) -- (-2.6238,-0.7037) -- (-2.8568,-0.6026) -- cycle;
\path[cbpmesh=19.2] (2.5806,-0.6522) -- (2.8097,-0.5528) -- (2.8568,-0.6026) -- (2.6238,-0.7037) -- cycle;
\path[cbpmesh=15.5] (-0.6652,-0.3657) -- (-0.4458,-0.3862) -- (-0.4630,-0.4685) -- (-0.6909,-0.4472) -- cycle;
\path[cbpmesh=16.5] (0.4458,-0.3862) -- (0.6652,-0.3657) -- (0.6909,-0.4472) -- (0.4630,-0.4685) -- cycle;
\path[cbpmesh=29.3] (-2.9375,-1.7636) -- (-2.6586,-1.8669) -- (-2.6427,-1.9549) -- (-2.9200,-1.8523) -- cycle;
\path[cbpmesh=38] (2.6586,-1.8669) -- (2.9375,-1.7636) -- (2.9200,-1.8523) -- (2.6427,-1.9549) -- cycle;
\path[cbpmesh=14.1] (-2.8568,-0.6026) -- (-2.6238,-0.7037) -- (-2.6679,-0.7635) -- (-2.9047,-0.6607) -- cycle;
\path[cbpmesh=19.9] (2.6238,-0.7037) -- (2.8568,-0.6026) -- (2.9047,-0.6607) -- (2.6679,-0.7635) -- cycle;
\path[cbpmesh=16.5] (-1.4344,-0.4595) -- (-1.2091,-0.5083) -- (-1.2439,-0.5604) -- (-1.4757,-0.5102) -- cycle;
\path[cbpmesh=17.6] (1.2091,-0.5083) -- (1.4344,-0.4595) -- (1.4757,-0.5102) -- (1.2439,-0.5604) -- cycle;
\path[cbpmesh=38] (-2.2663,-2.2401) -- (-1.9691,-2.3182) -- (-1.9352,-2.3581) -- (-2.2273,-2.2814) -- cycle;
\path[cbpmesh=38] (1.9691,-2.3182) -- (2.2663,-2.2401) -- (2.2273,-2.2814) -- (1.9352,-2.3581) -- cycle;
\path[cbpmesh=16.8] (-2.4086,-0.5036) -- (-2.1799,-0.5883) -- (-2.2194,-0.6209) -- (-2.4523,-0.5348) -- cycle;
\path[cbpmesh=18.0] (2.1799,-0.5883) -- (2.4086,-0.5036) -- (2.4523,-0.5348) -- (2.2194,-0.6209) -- cycle;
\path[cbpmesh=15.7] (-0.9146,-0.4176) -- (-0.6909,-0.4472) -- (-0.7150,-0.5182) -- (-0.9464,-0.4876) -- cycle;
\path[cbpmesh=16.9] (0.6909,-0.4472) -- (0.9146,-0.4176) -- (0.9464,-0.4876) -- (0.7150,-0.5182) -- cycle;
\path[cbpmesh=17.7] (-2.0543,-0.4961) -- (-1.8261,-0.5677) -- (-1.8647,-0.6002) -- (-2.0977,-0.5272) -- cycle;
\path[cbpmesh=18.1] (1.8261,-0.5677) -- (2.0543,-0.4961) -- (2.0977,-0.5272) -- (1.8647,-0.6002) -- cycle;
\path[cbpmesh=19.2] (-3.1593,-1.1893) -- (-2.9017,-1.3011) -- (-2.9281,-1.4328) -- (-3.1881,-1.3200) -- cycle;
\path[cbpmesh=30.1] (2.9017,-1.3011) -- (3.1593,-1.1893) -- (3.1881,-1.3200) -- (2.9281,-1.4328) -- cycle;
\path[cbpmesh=15.6] (-0.4458,-0.3862) -- (-0.2236,-0.3985) -- (-0.2322,-0.4813) -- (-0.4630,-0.4685) -- cycle;
\path[cbpmesh=16.3] (0.2236,-0.3985) -- (0.4458,-0.3862) -- (0.4630,-0.4685) -- (0.2322,-0.4813) -- cycle;
\path[cbpmesh=16.0] (-1.1720,-0.4485) -- (-0.9464,-0.4876) -- (-0.9763,-0.5486) -- (-1.2091,-0.5083) -- cycle;
\path[cbpmesh=17.2] (0.9464,-0.4876) -- (1.1720,-0.4485) -- (1.2091,-0.5083) -- (0.9763,-0.5486) -- cycle;
\path[cbpmesh=17.3] (-1.7432,-0.4947) -- (-1.5147,-0.5548) -- (-1.5516,-0.5945) -- (-1.7857,-0.5329) -- cycle;
\path[cbpmesh=17.9] (1.5147,-0.5548) -- (1.7432,-0.4947) -- (1.7857,-0.5329) -- (1.5516,-0.5945) -- cycle;
\path[cbpmesh=35.5] (-2.6188,-2.0316) -- (-2.3279,-2.1229) -- (-2.3000,-2.1873) -- (-2.5874,-2.0971) -- cycle;
\path[cbpmesh=38] (2.3279,-2.1229) -- (2.6188,-2.0316) -- (2.5874,-2.0971) -- (2.3000,-2.1873) -- cycle;
\path[cbpmesh=14.0] (-2.9047,-0.6607) -- (-2.6679,-0.7635) -- (-2.7131,-0.8330) -- (-2.9540,-0.7284) -- cycle;
\path[cbpmesh=20.7] (2.6679,-0.7635) -- (2.9047,-0.6607) -- (2.9540,-0.7284) -- (2.7131,-0.8330) -- cycle;
\path[cbpmesh=15.8] (-0.2236,-0.3985) -- (-0.0000,-0.4026) -- (-0.0000,-0.4856) -- (-0.2322,-0.4813) -- cycle;
\path[cbpmesh=16.0] (-0.0000,-0.4026) -- (0.2236,-0.3985) -- (0.2322,-0.4813) -- (-0.0000,-0.4856) -- cycle;
\path[cbpmesh=38] (-1.9352,-2.3581) -- (-1.6313,-2.4239) -- (-1.5993,-2.4508) -- (-1.8973,-2.3863) -- cycle;
\path[cbpmesh=38] (1.6313,-2.4239) -- (1.9352,-2.3581) -- (1.8973,-2.3863) -- (1.5993,-2.4508) -- cycle;
\path[cbpmesh=17.5] (-3.1163,-1.0484) -- (-2.8622,-1.1587) -- (-2.9017,-1.3011) -- (-3.1593,-1.1893) -- cycle;
\path[cbpmesh=28.0] (2.8622,-1.1587) -- (3.1163,-1.0484) -- (3.1593,-1.1893) -- (2.9017,-1.3011) -- cycle;
\path[cbpmesh=14.2] (-2.9540,-0.7284) -- (-2.7131,-0.8330) -- (-2.7600,-0.9133) -- (-3.0050,-0.8070) -- cycle;
\path[cbpmesh=21.6] (2.7131,-0.8330) -- (2.9540,-0.7284) -- (3.0050,-0.8070) -- (2.7600,-0.9133) -- cycle;
\path[cbpmesh=38] (-1.5993,-2.4508) -- (-1.2914,-2.5041) -- (-1.2633,-2.5180) -- (-1.5645,-2.4659) -- cycle;
\draw[draw=cbpInk!70,line width=.45pt] (-1.2633,-2.5180) -- (-1.5645,-2.4659);
\path[cbpmesh=38] (1.2914,-2.5041) -- (1.5993,-2.4508) -- (1.5645,-2.4659) -- (1.2633,-2.5180) -- cycle;
\draw[draw=cbpInk!70,line width=.45pt] (1.5645,-2.4659) -- (1.2633,-2.5180);
\path[cbpmesh=16.2] (-2.4523,-0.5348) -- (-2.2194,-0.6209) -- (-2.2583,-0.6567) -- (-2.4952,-0.5690) -- cycle;
\path[cbpmesh=18.0] (2.2194,-0.6209) -- (2.4523,-0.5348) -- (2.4952,-0.5690) -- (2.2583,-0.6567) -- cycle;
\path[cbpmesh=15.8] (-0.6909,-0.4472) -- (-0.4630,-0.4685) -- (-0.4791,-0.5403) -- (-0.7150,-0.5182) -- cycle;
\path[cbpmesh=16.7] (0.4630,-0.4685) -- (0.6909,-0.4472) -- (0.7150,-0.5182) -- (0.4791,-0.5403) -- cycle;
\path[cbpmesh=14.4] (-3.0050,-0.8070) -- (-2.7600,-0.9133) -- (-2.8089,-1.0058) -- (-3.0582,-0.8976) -- cycle;
\path[cbpmesh=22.4] (2.7600,-0.9133) -- (3.0050,-0.8070) -- (3.0582,-0.8976) -- (2.8089,-1.0058) -- cycle;
\path[cbpmesh=16.3] (-3.0582,-0.8976) -- (-2.8089,-1.0058) -- (-2.8622,-1.1587) -- (-3.1163,-1.0484) -- cycle;
\path[cbpmesh=26.1] (2.8089,-1.0058) -- (3.0582,-0.8976) -- (3.1163,-1.0484) -- (2.8622,-1.1587) -- cycle;
\path[cbpmesh=17.6] (-2.0977,-0.5272) -- (-1.8647,-0.6002) -- (-1.9018,-0.6320) -- (-2.1395,-0.5575) -- cycle;
\path[cbpmesh=18.1] (1.8647,-0.6002) -- (2.0977,-0.5272) -- (2.1395,-0.5575) -- (1.9018,-0.6320) -- cycle;
\path[cbpmesh=26.4] (-2.9452,-1.6641) -- (-2.6656,-1.7676) -- (-2.6586,-1.8669) -- (-2.9375,-1.7636) -- cycle;
\path[cbpmesh=36.7] (2.6656,-1.7676) -- (2.9452,-1.6641) -- (2.9375,-1.7636) -- (2.6586,-1.8669) -- cycle;
\path[cbpmesh=16.9] (-1.4757,-0.5102) -- (-1.2439,-0.5604) -- (-1.2768,-0.6063) -- (-1.5147,-0.5548) -- cycle;
\path[cbpmesh=17.7] (1.2439,-0.5604) -- (1.4757,-0.5102) -- (1.5147,-0.5548) -- (1.2768,-0.6063) -- cycle;
\path[cbpmesh=17.6] (-1.7857,-0.5329) -- (-1.5516,-0.5945) -- (-1.5867,-0.6306) -- (-1.8261,-0.5677) -- cycle;
\path[cbpmesh=18.0] (1.5516,-0.5945) -- (1.7857,-0.5329) -- (1.8261,-0.5677) -- (1.5867,-0.6306) -- cycle;
\path[cbpmesh=16.1] (-0.9464,-0.4876) -- (-0.7150,-0.5182) -- (-0.7376,-0.5802) -- (-0.9763,-0.5486) -- cycle;
\path[cbpmesh=17.1] (0.7150,-0.5182) -- (0.9464,-0.4876) -- (0.9763,-0.5486) -- (0.7376,-0.5802) -- cycle;
\path[cbpmesh=15.5] (-2.4952,-0.5690) -- (-2.2583,-0.6567) -- (-2.2969,-0.6969) -- (-2.5379,-0.6078) -- cycle;
\path[cbpmesh=18.1] (2.2583,-0.6567) -- (2.4952,-0.5690) -- (2.5379,-0.6078) -- (2.2969,-0.6969) -- cycle;
\path[cbpmesh=16.5] (-1.2091,-0.5083) -- (-0.9763,-0.5486) -- (-1.0045,-0.6019) -- (-1.2439,-0.5604) -- cycle;
\path[cbpmesh=17.4] (0.9763,-0.5486) -- (1.2091,-0.5083) -- (1.2439,-0.5604) -- (1.0045,-0.6019) -- cycle;
\path[cbpmesh=16.0] (-0.4630,-0.4685) -- (-0.2322,-0.4813) -- (-0.2403,-0.5535) -- (-0.4791,-0.5403) -- cycle;
\path[cbpmesh=16.5] (0.2322,-0.4813) -- (0.4630,-0.4685) -- (0.4791,-0.5403) -- (0.2403,-0.5535) -- cycle;
\path[cbpmesh=38] (-2.3000,-2.1873) -- (-1.9985,-2.2666) -- (-1.9691,-2.3182) -- (-2.2663,-2.2401) -- cycle;
\path[cbpmesh=38] (1.9985,-2.2666) -- (2.3000,-2.1873) -- (2.2663,-2.2401) -- (1.9691,-2.3182) -- cycle;
\path[cbpmesh=32.4] (-2.6427,-1.9549) -- (-2.3492,-2.0470) -- (-2.3279,-2.1229) -- (-2.6188,-2.0316) -- cycle;
\path[cbpmesh=38] (2.3492,-2.0470) -- (2.6427,-1.9549) -- (2.6188,-2.0316) -- (2.3279,-2.1229) -- cycle;
\path[cbpmesh=16.1] (-0.2322,-0.4813) -- (-0.0000,-0.4856) -- (-0.0000,-0.5580) -- (-0.2403,-0.5535) -- cycle;
\path[cbpmesh=16.3] (-0.0000,-0.4856) -- (0.2322,-0.4813) -- (0.2403,-0.5535) -- (-0.0000,-0.5580) -- cycle;
\path[cbpmesh=17.3] (-2.1395,-0.5575) -- (-1.9018,-0.6320) -- (-1.9378,-0.6642) -- (-2.1799,-0.5883) -- cycle;
\path[cbpmesh=18.0] (1.9018,-0.6320) -- (2.1395,-0.5575) -- (2.1799,-0.5883) -- (1.9378,-0.6642) -- cycle;
\path[cbpmesh=14.8] (-2.5379,-0.6078) -- (-2.2969,-0.6969) -- (-2.3356,-0.7429) -- (-2.5806,-0.6522) -- cycle;
\path[cbpmesh=18.4] (2.2969,-0.6969) -- (2.5379,-0.6078) -- (2.5806,-0.6522) -- (2.3356,-0.7429) -- cycle;
\path[cbpmesh=23.6] (-2.9424,-1.5538) -- (-2.6630,-1.6572) -- (-2.6656,-1.7676) -- (-2.9452,-1.6641) -- cycle;
\path[cbpmesh=34.0] (2.6630,-1.6572) -- (2.9424,-1.5538) -- (2.9452,-1.6641) -- (2.6656,-1.7676) -- cycle;
\path[cbpmesh=17.3] (-1.5147,-0.5548) -- (-1.2768,-0.6063) -- (-1.3079,-0.6472) -- (-1.5516,-0.5945) -- cycle;
\path[cbpmesh=17.9] (1.2768,-0.6063) -- (1.5147,-0.5548) -- (1.5516,-0.5945) -- (1.3079,-0.6472) -- cycle;
\path[cbpmesh=16.2] (-0.7150,-0.5182) -- (-0.4791,-0.5403) -- (-0.4942,-0.6029) -- (-0.7376,-0.5802) -- cycle;
\path[cbpmesh=16.9] (0.4791,-0.5403) -- (0.7150,-0.5182) -- (0.7376,-0.5802) -- (0.4942,-0.6029) -- cycle;
\path[cbpmesh=38] (-1.2914,-2.5041) -- (-0.9756,-2.5460) -- (-0.9544,-2.5589) -- (-1.2633,-2.5180) -- cycle;
\draw[draw=cbpInk!70,line width=.45pt] (-0.9544,-2.5589) -- (-1.2633,-2.5180);
\path[cbpmesh=38] (0.9756,-2.5460) -- (1.2914,-2.5041) -- (1.2633,-2.5180) -- (0.9544,-2.5589) -- cycle;
\draw[draw=cbpInk!70,line width=.45pt] (1.2633,-2.5180) -- (0.9544,-2.5589);
\path[cbpmesh=38] (-1.9691,-2.3182) -- (-1.6598,-2.3852) -- (-1.6313,-2.4239) -- (-1.9352,-2.3581) -- cycle;
\path[cbpmesh=38] (1.6598,-2.3852) -- (1.9691,-2.3182) -- (1.9352,-2.3581) -- (1.6313,-2.4239) -- cycle;
\path[cbpmesh=17.7] (-1.8261,-0.5677) -- (-1.5867,-0.6306) -- (-1.6202,-0.6645) -- (-1.8647,-0.6002) -- cycle;
\path[cbpmesh=18.1] (1.5867,-0.6306) -- (1.8261,-0.5677) -- (1.8647,-0.6002) -- (1.6202,-0.6645) -- cycle;
\path[cbpmesh=38] (-1.6313,-2.4239) -- (-1.3172,-2.4783) -- (-1.2914,-2.5041) -- (-1.5993,-2.4508) -- cycle;
\path[cbpmesh=38] (1.3172,-2.4783) -- (1.6313,-2.4239) -- (1.5993,-2.4508) -- (1.2914,-2.5041) -- cycle;
\path[cbpmesh=16.9] (-1.2439,-0.5604) -- (-1.0045,-0.6019) -- (-1.0310,-0.6489) -- (-1.2768,-0.6063) -- cycle;
\path[cbpmesh=17.6] (1.0045,-0.6019) -- (1.2439,-0.5604) -- (1.2768,-0.6063) -- (1.0310,-0.6489) -- cycle;
\path[cbpmesh=14.3] (-2.5806,-0.6522) -- (-2.3356,-0.7429) -- (-2.3747,-0.7959) -- (-2.6238,-0.7037) -- cycle;
\path[cbpmesh=18.8] (2.3356,-0.7429) -- (2.5806,-0.6522) -- (2.6238,-0.7037) -- (2.3747,-0.7959) -- cycle;
\path[cbpmesh=16.6] (-0.9763,-0.5486) -- (-0.7376,-0.5802) -- (-0.7588,-0.6344) -- (-1.0045,-0.6019) -- cycle;
\path[cbpmesh=17.3] (0.7376,-0.5802) -- (0.9763,-0.5486) -- (1.0045,-0.6019) -- (0.7588,-0.6344) -- cycle;
\path[cbpmesh=16.8] (-2.1799,-0.5883) -- (-1.9378,-0.6642) -- (-1.9729,-0.6983) -- (-2.2194,-0.6209) -- cycle;
\path[cbpmesh=17.9] (1.9378,-0.6642) -- (2.1799,-0.5883) -- (2.2194,-0.6209) -- (1.9729,-0.6983) -- cycle;
\path[cbpmesh=16.3] (-0.4791,-0.5403) -- (-0.2403,-0.5535) -- (-0.2479,-0.6166) -- (-0.4942,-0.6029) -- cycle;
\path[cbpmesh=16.7] (0.2403,-0.5535) -- (0.4791,-0.5403) -- (0.4942,-0.6029) -- (0.2479,-0.6166) -- cycle;
\path[cbpmesh=21.2] (-2.9281,-1.4328) -- (-2.6501,-1.5357) -- (-2.6630,-1.6572) -- (-2.9424,-1.5538) -- cycle;
\path[cbpmesh=31.4] (2.6501,-1.5357) -- (2.9281,-1.4328) -- (2.9424,-1.5538) -- (2.6630,-1.6572) -- cycle;
\path[cbpmesh=29.4] (-2.6586,-1.8669) -- (-2.3633,-1.9595) -- (-2.3492,-2.0470) -- (-2.6427,-1.9549) -- cycle;
\path[cbpmesh=38] (2.3633,-1.9595) -- (2.6586,-1.8669) -- (2.6427,-1.9549) -- (2.3492,-2.0470) -- cycle;
\path[cbpmesh=35.6] (-2.3279,-2.1229) -- (-2.0227,-2.2031) -- (-1.9985,-2.2666) -- (-2.3000,-2.1873) -- cycle;
\path[cbpmesh=38] (2.0227,-2.2031) -- (2.3279,-2.1229) -- (2.3000,-2.1873) -- (1.9985,-2.2666) -- cycle;
\path[cbpmesh=16.4] (-0.2403,-0.5535) -- (-0.0000,-0.5580) -- (-0.0000,-0.6212) -- (-0.2479,-0.6166) -- cycle;
\path[cbpmesh=16.6] (-0.0000,-0.5580) -- (0.2403,-0.5535) -- (0.2479,-0.6166) -- (-0.0000,-0.6212) -- cycle;
\path[cbpmesh=14.1] (-2.6238,-0.7037) -- (-2.3747,-0.7959) -- (-2.4146,-0.8573) -- (-2.6679,-0.7635) -- cycle;
\path[cbpmesh=19.4] (2.3747,-0.7959) -- (2.6238,-0.7037) -- (2.6679,-0.7635) -- (2.4146,-0.8573) -- cycle;
\path[cbpmesh=17.6] (-1.5516,-0.5945) -- (-1.3079,-0.6472) -- (-1.3374,-0.6845) -- (-1.5867,-0.6306) -- cycle;
\path[cbpmesh=18.0] (1.3079,-0.6472) -- (1.5516,-0.5945) -- (1.5867,-0.6306) -- (1.3374,-0.6845) -- cycle;
\path[cbpmesh=17.6] (-1.8647,-0.6002) -- (-1.6202,-0.6645) -- (-1.6525,-0.6976) -- (-1.9018,-0.6320) -- cycle;
\path[cbpmesh=18.0] (1.6202,-0.6645) -- (1.8647,-0.6002) -- (1.9018,-0.6320) -- (1.6525,-0.6976) -- cycle;
\path[cbpmesh=38] (-0.9756,-2.5460) -- (-0.6538,-2.5760) -- (-0.6395,-2.5883) -- (-0.9544,-2.5589) -- cycle;
\draw[draw=cbpInk!70,line width=.45pt] (-0.6395,-2.5883) -- (-0.9544,-2.5589);
\path[cbpmesh=38] (0.6538,-2.5760) -- (0.9756,-2.5460) -- (0.9544,-2.5589) -- (0.6395,-2.5883) -- cycle;
\draw[draw=cbpInk!70,line width=.45pt] (0.9544,-2.5589) -- (0.6395,-2.5883);
\path[cbpmesh=16.2] (-2.2194,-0.6209) -- (-1.9729,-0.6983) -- (-2.0075,-0.7354) -- (-2.2583,-0.6567) -- cycle;
\path[cbpmesh=17.8] (1.9729,-0.6983) -- (2.2194,-0.6209) -- (2.2583,-0.6567) -- (2.0075,-0.7354) -- cycle;
\path[cbpmesh=16.6] (-0.7376,-0.5802) -- (-0.4942,-0.6029) -- (-0.5085,-0.6578) -- (-0.7588,-0.6344) -- cycle;
\path[cbpmesh=17.1] (0.4942,-0.6029) -- (0.7376,-0.5802) -- (0.7588,-0.6344) -- (0.5085,-0.6578) -- cycle;
\path[cbpmesh=14.0] (-2.6679,-0.7635) -- (-2.4146,-0.8573) -- (-2.4555,-0.9283) -- (-2.7131,-0.8330) -- cycle;
\path[cbpmesh=20.1] (2.4146,-0.8573) -- (2.6679,-0.7635) -- (2.7131,-0.8330) -- (2.4555,-0.9283) -- cycle;
\path[cbpmesh=19.1] (-2.9017,-1.3011) -- (-2.6262,-1.4030) -- (-2.6501,-1.5357) -- (-2.9281,-1.4328) -- cycle;
\path[cbpmesh=29.1] (2.6262,-1.4030) -- (2.9017,-1.3011) -- (2.9281,-1.4328) -- (2.6501,-1.5357) -- cycle;
\path[cbpmesh=17.3] (-1.2768,-0.6063) -- (-1.0310,-0.6489) -- (-1.0561,-0.6908) -- (-1.3079,-0.6472) -- cycle;
\path[cbpmesh=17.8] (1.0310,-0.6489) -- (1.2768,-0.6063) -- (1.3079,-0.6472) -- (1.0561,-0.6908) -- cycle;
\path[cbpmesh=17.0] (-1.0045,-0.6019) -- (-0.7588,-0.6344) -- (-0.7789,-0.6822) -- (-1.0310,-0.6489) -- cycle;
\path[cbpmesh=17.5] (0.7588,-0.6344) -- (1.0045,-0.6019) -- (1.0310,-0.6489) -- (0.7789,-0.6822) -- cycle;
\path[cbpmesh=38] (-1.9985,-2.2666) -- (-1.6846,-2.3345) -- (-1.6598,-2.3852) -- (-1.9691,-2.3182) -- cycle;
\path[cbpmesh=38] (1.6846,-2.3345) -- (1.9985,-2.2666) -- (1.9691,-2.3182) -- (1.6598,-2.3852) -- cycle;
\path[cbpmesh=14.1] (-2.7131,-0.8330) -- (-2.4555,-0.9283) -- (-2.4980,-1.0103) -- (-2.7600,-0.9133) -- cycle;
\path[cbpmesh=20.9] (2.4555,-0.9283) -- (2.7131,-0.8330) -- (2.7600,-0.9133) -- (2.4980,-1.0103) -- cycle;
\path[cbpmesh=38] (-1.3172,-2.4783) -- (-0.9951,-2.5210) -- (-0.9756,-2.5460) -- (-1.2914,-2.5041) -- cycle;
\path[cbpmesh=38] (0.9951,-2.5210) -- (1.3172,-2.4783) -- (1.2914,-2.5041) -- (0.9756,-2.5460) -- cycle;
\path[cbpmesh=17.5] (-2.8622,-1.1587) -- (-2.5904,-1.2593) -- (-2.6262,-1.4030) -- (-2.9017,-1.3011) -- cycle;
\path[cbpmesh=27.0] (2.5904,-1.2593) -- (2.8622,-1.1587) -- (2.9017,-1.3011) -- (2.6262,-1.4030) -- cycle;
\path[cbpmesh=16.7] (-0.4942,-0.6029) -- (-0.2479,-0.6166) -- (-0.2550,-0.6719) -- (-0.5085,-0.6578) -- cycle;
\path[cbpmesh=17.0] (0.2479,-0.6166) -- (0.4942,-0.6029) -- (0.5085,-0.6578) -- (0.2550,-0.6719) -- cycle;
\path[cbpmesh=15.5] (-2.2583,-0.6567) -- (-2.0075,-0.7354) -- (-2.0418,-0.7770) -- (-2.2969,-0.6969) -- cycle;
\path[cbpmesh=17.8] (2.0075,-0.7354) -- (2.2583,-0.6567) -- (2.2969,-0.6969) -- (2.0418,-0.7770) -- cycle;
\path[cbpmesh=17.3] (-1.9018,-0.6320) -- (-1.6525,-0.6976) -- (-1.6837,-0.7310) -- (-1.9378,-0.6642) -- cycle;
\path[cbpmesh=17.9] (1.6525,-0.6976) -- (1.9018,-0.6320) -- (1.9378,-0.6642) -- (1.6837,-0.7310) -- cycle;
\path[cbpmesh=38] (-0.6538,-2.5760) -- (-0.3279,-2.5941) -- (-0.3208,-2.6061) -- (-0.6395,-2.5883) -- cycle;
\draw[draw=cbpInk!70,line width=.45pt] (-0.3208,-2.6061) -- (-0.6395,-2.5883);
\path[cbpmesh=38] (0.3279,-2.5941) -- (0.6538,-2.5760) -- (0.6395,-2.5883) -- (0.3208,-2.6061) -- cycle;
\draw[draw=cbpInk!70,line width=.45pt] (0.6395,-2.5883) -- (0.3208,-2.6061);
\path[cbpmesh=26.4] (-2.6656,-1.7676) -- (-2.3695,-1.8605) -- (-2.3633,-1.9595) -- (-2.6586,-1.8669) -- cycle;
\path[cbpmesh=35.7] (2.3695,-1.8605) -- (2.6656,-1.7676) -- (2.6586,-1.8669) -- (2.3633,-1.9595) -- cycle;
\path[cbpmesh=17.7] (-1.5867,-0.6306) -- (-1.3374,-0.6845) -- (-1.3657,-0.7196) -- (-1.6202,-0.6645) -- cycle;
\path[cbpmesh=18.0] (1.3374,-0.6845) -- (1.5867,-0.6306) -- (1.6202,-0.6645) -- (1.3657,-0.7196) -- cycle;
\path[cbpmesh=38] (-1.6598,-2.3852) -- (-1.3403,-2.4405) -- (-1.3172,-2.4783) -- (-1.6313,-2.4239) -- cycle;
\path[cbpmesh=38] (1.3403,-2.4405) -- (1.6598,-2.3852) -- (1.6313,-2.4239) -- (1.3172,-2.4783) -- cycle;
\path[cbpmesh=14.4] (-2.7600,-0.9133) -- (-2.4980,-1.0103) -- (-2.5422,-1.1045) -- (-2.8089,-1.0058) -- cycle;
\path[cbpmesh=21.7] (2.4980,-1.0103) -- (2.7600,-0.9133) -- (2.8089,-1.0058) -- (2.5422,-1.1045) -- cycle;
\path[cbpmesh=16.2] (-2.8089,-1.0058) -- (-2.5422,-1.1045) -- (-2.5904,-1.2593) -- (-2.8622,-1.1587) -- cycle;
\path[cbpmesh=25.2] (2.5422,-1.1045) -- (2.8089,-1.0058) -- (2.8622,-1.1587) -- (2.5904,-1.2593) -- cycle;
\path[cbpmesh=16.9] (-0.0000,-0.6212) -- (0.2479,-0.6166) -- (0.2550,-0.6719) -- (-0.0000,-0.6766) -- cycle;
\path[cbpmesh=16.8] (-0.2479,-0.6166) -- (-0.0000,-0.6212) -- (-0.0000,-0.6766) -- (-0.2550,-0.6719) -- cycle;
\path[cbpmesh=32.5] (-2.3492,-2.0470) -- (-2.0412,-2.1279) -- (-2.0227,-2.2031) -- (-2.3279,-2.1229) -- cycle;
\path[cbpmesh=38] (2.0412,-2.1279) -- (2.3492,-2.0470) -- (2.3279,-2.1229) -- (2.0227,-2.2031) -- cycle;
\path[cbpmesh=38] (-0.3279,-2.5941) -- (-0.0000,-2.6002) -- (-0.0000,-2.6120) -- (-0.3208,-2.6061) -- cycle;
\draw[draw=cbpInk!70,line width=.45pt] (-0.0000,-2.6120) -- (-0.3208,-2.6061);
\path[cbpmesh=38] (-0.0000,-2.6002) -- (0.3279,-2.5941) -- (0.3208,-2.6061) -- (-0.0000,-2.6120) -- cycle;
\draw[draw=cbpInk!70,line width=.45pt] (0.3208,-2.6061) -- (-0.0000,-2.6120);
\path[cbpmesh=17.0] (-0.7588,-0.6344) -- (-0.5085,-0.6578) -- (-0.5219,-0.7062) -- (-0.7789,-0.6822) -- cycle;
\path[cbpmesh=17.4] (0.5085,-0.6578) -- (0.7588,-0.6344) -- (0.7789,-0.6822) -- (0.5219,-0.7062) -- cycle;
\path[cbpmesh=14.8] (-2.2969,-0.6969) -- (-2.0418,-0.7770) -- (-2.0762,-0.8243) -- (-2.3356,-0.7429) -- cycle;
\path[cbpmesh=18.0] (2.0418,-0.7770) -- (2.2969,-0.6969) -- (2.3356,-0.7429) -- (2.0762,-0.8243) -- cycle;
\path[cbpmesh=17.6] (-1.3079,-0.6472) -- (-1.0561,-0.6908) -- (-1.0800,-0.7291) -- (-1.3374,-0.6845) -- cycle;
\path[cbpmesh=17.9] (1.0561,-0.6908) -- (1.3079,-0.6472) -- (1.3374,-0.6845) -- (1.0800,-0.7291) -- cycle;
\path[cbpmesh=17.4] (-1.0310,-0.6489) -- (-0.7789,-0.6822) -- (-0.7978,-0.7250) -- (-1.0561,-0.6908) -- cycle;
\path[cbpmesh=17.7] (0.7789,-0.6822) -- (1.0310,-0.6489) -- (1.0561,-0.6908) -- (0.7978,-0.7250) -- cycle;
\path[cbpmesh=16.8] (-1.9378,-0.6642) -- (-1.6837,-0.7310) -- (-1.7142,-0.7663) -- (-1.9729,-0.6983) -- cycle;
\path[cbpmesh=17.8] (1.6837,-0.7310) -- (1.9378,-0.6642) -- (1.9729,-0.6983) -- (1.7142,-0.7663) -- cycle;
\path[cbpmesh=17.6] (-1.6202,-0.6645) -- (-1.3657,-0.7196) -- (-1.3929,-0.7537) -- (-1.6525,-0.6976) -- cycle;
\path[cbpmesh=18.0] (1.3657,-0.7196) -- (1.6202,-0.6645) -- (1.6525,-0.6976) -- (1.3929,-0.7537) -- cycle;
\path[cbpmesh=38] (-0.9951,-2.5210) -- (-0.6668,-2.5516) -- (-0.6538,-2.5760) -- (-0.9756,-2.5460) -- cycle;
\path[cbpmesh=38] (0.6668,-2.5516) -- (0.9951,-2.5210) -- (0.9756,-2.5460) -- (0.6538,-2.5760) -- cycle;
\path[cbpmesh=23.6] (-2.6630,-1.6572) -- (-2.3672,-1.7500) -- (-2.3695,-1.8605) -- (-2.6656,-1.7676) -- cycle;
\path[cbpmesh=33.0] (2.3672,-1.7500) -- (2.6630,-1.6572) -- (2.6656,-1.7676) -- (2.3695,-1.8605) -- cycle;
\path[cbpmesh=17.1] (-0.5085,-0.6578) -- (-0.2550,-0.6719) -- (-0.2618,-0.7207) -- (-0.5219,-0.7062) -- cycle;
\path[cbpmesh=17.3] (0.2550,-0.6719) -- (0.5085,-0.6578) -- (0.5219,-0.7062) -- (0.2618,-0.7207) -- cycle;
\path[cbpmesh=14.3] (-2.3356,-0.7429) -- (-2.0762,-0.8243) -- (-2.1109,-0.8787) -- (-2.3747,-0.7959) -- cycle;
\path[cbpmesh=18.3] (2.0762,-0.8243) -- (2.3356,-0.7429) -- (2.3747,-0.7959) -- (2.1109,-0.8787) -- cycle;
\path[cbpmesh=35.8] (-2.0227,-2.2031) -- (-1.7050,-2.2719) -- (-1.6846,-2.3345) -- (-1.9985,-2.2666) -- cycle;
\path[cbpmesh=38] (1.7050,-2.2719) -- (2.0227,-2.2031) -- (1.9985,-2.2666) -- (1.6846,-2.3345) -- cycle;
\path[cbpmesh=17.2] (-0.2550,-0.6719) -- (-0.0000,-0.6766) -- (-0.0000,-0.7255) -- (-0.2618,-0.7207) -- cycle;
\path[cbpmesh=17.2] (-0.0000,-0.6766) -- (0.2550,-0.6719) -- (0.2618,-0.7207) -- (-0.0000,-0.7255) -- cycle;
\path[cbpmesh=16.2] (-1.9729,-0.6983) -- (-1.7142,-0.7663) -- (-1.7443,-0.8046) -- (-2.0075,-0.7354) -- cycle;
\path[cbpmesh=17.7] (1.7142,-0.7663) -- (1.9729,-0.6983) -- (2.0075,-0.7354) -- (1.7443,-0.8046) -- cycle;
\path[cbpmesh=38] (-1.3403,-2.4405) -- (-1.0125,-2.4839) -- (-0.9951,-2.5210) -- (-1.3172,-2.4783) -- cycle;
\path[cbpmesh=38] (1.0125,-2.4839) -- (1.3403,-2.4405) -- (1.3172,-2.4783) -- (0.9951,-2.5210) -- cycle;
\path[cbpmesh=17.7] (-1.3374,-0.6845) -- (-1.0800,-0.7291) -- (-1.1028,-0.7651) -- (-1.3657,-0.7196) -- cycle;
\path[cbpmesh=18.0] (1.0800,-0.7291) -- (1.3374,-0.6845) -- (1.3657,-0.7196) -- (1.1028,-0.7651) -- cycle;
\path[cbpmesh=29.5] (-2.3633,-1.9595) -- (-2.0534,-2.0410) -- (-2.0412,-2.1279) -- (-2.3492,-2.0470) -- cycle;
\path[cbpmesh=37.5] (2.0534,-2.0410) -- (2.3633,-1.9595) -- (2.3492,-2.0470) -- (2.0412,-2.1279) -- cycle;
\path[cbpmesh=17.4] (-0.7789,-0.6822) -- (-0.5219,-0.7062) -- (-0.5346,-0.7496) -- (-0.7978,-0.7250) -- cycle;
\path[cbpmesh=17.6] (0.5219,-0.7062) -- (0.7789,-0.6822) -- (0.7978,-0.7250) -- (0.5346,-0.7496) -- cycle;
\path[cbpmesh=38] (-1.6846,-2.3345) -- (-1.3603,-2.3907) -- (-1.3403,-2.4405) -- (-1.6598,-2.3852) -- cycle;
\path[cbpmesh=38] (1.3603,-2.3907) -- (1.6846,-2.3345) -- (1.6598,-2.3852) -- (1.3403,-2.4405) -- cycle;
\path[cbpmesh=14.1] (-2.3747,-0.7959) -- (-2.1109,-0.8787) -- (-2.1464,-0.9414) -- (-2.4146,-0.8573) -- cycle;
\path[cbpmesh=18.9] (2.1109,-0.8787) -- (2.3747,-0.7959) -- (2.4146,-0.8573) -- (2.1464,-0.9414) -- cycle;
\path[cbpmesh=38] (-0.6668,-2.5516) -- (-0.3344,-2.5701) -- (-0.3279,-2.5941) -- (-0.6538,-2.5760) -- cycle;
\path[cbpmesh=38] (0.3344,-2.5701) -- (0.6668,-2.5516) -- (0.6538,-2.5760) -- (0.3279,-2.5941) -- cycle;
\path[cbpmesh=17.6] (-1.0561,-0.6908) -- (-0.7978,-0.7250) -- (-0.8159,-0.7641) -- (-1.0800,-0.7291) -- cycle;
\path[cbpmesh=17.9] (0.7978,-0.7250) -- (1.0561,-0.6908) -- (1.0800,-0.7291) -- (0.8159,-0.7641) -- cycle;
\path[cbpmesh=21.2] (-2.6501,-1.5357) -- (-2.3557,-1.6280) -- (-2.3672,-1.7500) -- (-2.6630,-1.6572) -- cycle;
\path[cbpmesh=30.4] (2.3557,-1.6280) -- (2.6501,-1.5357) -- (2.6630,-1.6572) -- (2.3672,-1.7500) -- cycle;
\path[cbpmesh=17.3] (-1.6525,-0.6976) -- (-1.3929,-0.7537) -- (-1.4193,-0.7883) -- (-1.6837,-0.7310) -- cycle;
\path[cbpmesh=17.9] (1.3929,-0.7537) -- (1.6525,-0.6976) -- (1.6837,-0.7310) -- (1.4193,-0.7883) -- cycle;
\path[cbpmesh=38] (-0.3344,-2.5701) -- (-0.0000,-2.5763) -- (-0.0000,-2.6002) -- (-0.3279,-2.5941) -- cycle;
\path[cbpmesh=38] (-0.0000,-2.5763) -- (0.3344,-2.5701) -- (0.3279,-2.5941) -- (-0.0000,-2.6002) -- cycle;
\path[cbpmesh=14.0] (-2.4146,-0.8573) -- (-2.1464,-0.9414) -- (-2.1828,-1.0139) -- (-2.4555,-0.9283) -- cycle;
\path[cbpmesh=19.5] (2.1464,-0.9414) -- (2.4146,-0.8573) -- (2.4555,-0.9283) -- (2.1828,-1.0139) -- cycle;
\path[cbpmesh=17.4] (-0.5219,-0.7062) -- (-0.2618,-0.7207) -- (-0.2681,-0.7644) -- (-0.5346,-0.7496) -- cycle;
\path[cbpmesh=17.6] (0.2618,-0.7207) -- (0.5219,-0.7062) -- (0.5346,-0.7496) -- (0.2681,-0.7644) -- cycle;
\path[cbpmesh=15.5] (-2.0075,-0.7354) -- (-1.7443,-0.8046) -- (-1.7741,-0.8474) -- (-2.0418,-0.7770) -- cycle;
\path[cbpmesh=17.6] (1.7443,-0.8046) -- (2.0075,-0.7354) -- (2.0418,-0.7770) -- (1.7741,-0.8474) -- cycle;
\path[cbpmesh=19.2] (-2.6262,-1.4030) -- (-2.3345,-1.4945) -- (-2.3557,-1.6280) -- (-2.6501,-1.5357) -- cycle;
\path[cbpmesh=28.1] (2.3345,-1.4945) -- (2.6262,-1.4030) -- (2.6501,-1.5357) -- (2.3557,-1.6280) -- cycle;
\path[cbpmesh=17.5] (-0.0000,-0.7255) -- (0.2618,-0.7207) -- (0.2681,-0.7644) -- (-0.0000,-0.7694) -- cycle;
\path[cbpmesh=17.5] (-0.2618,-0.7207) -- (-0.0000,-0.7255) -- (-0.0000,-0.7694) -- (-0.2681,-0.7644) -- cycle;
\path[cbpmesh=17.6] (-1.3657,-0.7196) -- (-1.1028,-0.7651) -- (-1.1248,-0.8002) -- (-1.3929,-0.7537) -- cycle;
\path[cbpmesh=17.9] (1.1028,-0.7651) -- (1.3657,-0.7196) -- (1.3929,-0.7537) -- (1.1248,-0.8002) -- cycle;
\path[cbpmesh=14.1] (-2.4555,-0.9283) -- (-2.1828,-1.0139) -- (-2.2205,-1.0973) -- (-2.4980,-1.0103) -- cycle;
\path[cbpmesh=20.2] (2.1828,-1.0139) -- (2.4555,-0.9283) -- (2.4980,-1.0103) -- (2.2205,-1.0973) -- cycle;
\path[cbpmesh=32.7] (-2.0412,-2.1279) -- (-1.7206,-2.1973) -- (-1.7050,-2.2719) -- (-2.0227,-2.2031) -- cycle;
\path[cbpmesh=38] (1.7206,-2.1973) -- (2.0412,-2.1279) -- (2.0227,-2.2031) -- (1.7050,-2.2719) -- cycle;
\path[cbpmesh=16.8] (-1.6837,-0.7310) -- (-1.4193,-0.7883) -- (-1.4450,-0.8245) -- (-1.7142,-0.7663) -- cycle;
\path[cbpmesh=17.7] (1.4193,-0.7883) -- (1.6837,-0.7310) -- (1.7142,-0.7663) -- (1.4450,-0.8245) -- cycle;
\path[cbpmesh=17.6] (-0.7978,-0.7250) -- (-0.5346,-0.7496) -- (-0.5467,-0.7892) -- (-0.8159,-0.7641) -- cycle;
\path[cbpmesh=17.8] (0.5346,-0.7496) -- (0.7978,-0.7250) -- (0.8159,-0.7641) -- (0.5467,-0.7892) -- cycle;
\path[cbpmesh=38] (-1.0125,-2.4839) -- (-0.6785,-2.5151) -- (-0.6668,-2.5516) -- (-0.9951,-2.5210) -- cycle;
\path[cbpmesh=38] (0.6785,-2.5151) -- (1.0125,-2.4839) -- (0.9951,-2.5210) -- (0.6668,-2.5516) -- cycle;
\path[cbpmesh=26.5] (-2.3695,-1.8605) -- (-2.0588,-1.9422) -- (-2.0534,-2.0410) -- (-2.3633,-1.9595) -- cycle;
\path[cbpmesh=34.7] (2.0588,-1.9422) -- (2.3695,-1.8605) -- (2.3633,-1.9595) -- (2.0534,-2.0410) -- cycle;
\path[cbpmesh=17.7] (-1.0800,-0.7291) -- (-0.8159,-0.7641) -- (-0.8331,-0.8009) -- (-1.1028,-0.7651) -- cycle;
\path[cbpmesh=17.9] (0.8159,-0.7641) -- (1.0800,-0.7291) -- (1.1028,-0.7651) -- (0.8331,-0.8009) -- cycle;
\path[cbpmesh=17.5] (-2.5904,-1.2593) -- (-2.3027,-1.3495) -- (-2.3345,-1.4945) -- (-2.6262,-1.4030) -- cycle;
\path[cbpmesh=26.0] (2.3027,-1.3495) -- (2.5904,-1.2593) -- (2.6262,-1.4030) -- (2.3345,-1.4945) -- cycle;
\path[cbpmesh=14.9] (-2.0418,-0.7770) -- (-1.7741,-0.8474) -- (-1.8040,-0.8959) -- (-2.0762,-0.8243) -- cycle;
\path[cbpmesh=17.7] (1.7741,-0.8474) -- (2.0418,-0.7770) -- (2.0762,-0.8243) -- (1.8040,-0.8959) -- cycle;
\path[cbpmesh=14.4] (-2.4980,-1.0103) -- (-2.2205,-1.0973) -- (-2.2598,-1.1931) -- (-2.5422,-1.1045) -- cycle;
\path[cbpmesh=21.0] (2.2205,-1.0973) -- (2.4980,-1.0103) -- (2.5422,-1.1045) -- (2.2598,-1.1931) -- cycle;
\path[cbpmesh=16.2] (-2.5422,-1.1045) -- (-2.2598,-1.1931) -- (-2.3027,-1.3495) -- (-2.5904,-1.2593) -- cycle;
\path[cbpmesh=24.3] (2.2598,-1.1931) -- (2.5422,-1.1045) -- (2.5904,-1.2593) -- (2.3027,-1.3495) -- cycle;
\path[cbpmesh=36.0] (-1.7050,-2.2719) -- (-1.3768,-2.3287) -- (-1.3603,-2.3907) -- (-1.6846,-2.3345) -- cycle;
\path[cbpmesh=38] (1.3768,-2.3287) -- (1.7050,-2.2719) -- (1.6846,-2.3345) -- (1.3603,-2.3907) -- cycle;
\path[cbpmesh=38] (-1.3603,-2.3907) -- (-1.0276,-2.4347) -- (-1.0125,-2.4839) -- (-1.3403,-2.4405) -- cycle;
\path[cbpmesh=38] (1.0276,-2.4347) -- (1.3603,-2.3907) -- (1.3403,-2.4405) -- (1.0125,-2.4839) -- cycle;
\path[cbpmesh=17.7] (-0.5346,-0.7496) -- (-0.2681,-0.7644) -- (-0.2742,-0.8044) -- (-0.5467,-0.7892) -- cycle;
\path[cbpmesh=17.8] (0.2681,-0.7644) -- (0.5346,-0.7496) -- (0.5467,-0.7892) -- (0.2742,-0.8044) -- cycle;
\path[cbpmesh=17.4] (-1.3929,-0.7537) -- (-1.1248,-0.8002) -- (-1.1461,-0.8356) -- (-1.4193,-0.7883) -- cycle;
\path[cbpmesh=17.8] (1.1248,-0.8002) -- (1.3929,-0.7537) -- (1.4193,-0.7883) -- (1.1461,-0.8356) -- cycle;
\path[cbpmesh=38] (-0.6785,-2.5151) -- (-0.3403,-2.5339) -- (-0.3344,-2.5701) -- (-0.6668,-2.5516) -- cycle;
\path[cbpmesh=38] (0.3403,-2.5339) -- (0.6785,-2.5151) -- (0.6668,-2.5516) -- (0.3344,-2.5701) -- cycle;
\path[cbpmesh=16.2] (-1.7142,-0.7663) -- (-1.4450,-0.8245) -- (-1.4703,-0.8639) -- (-1.7443,-0.8046) -- cycle;
\path[cbpmesh=17.5] (1.4450,-0.8245) -- (1.7142,-0.7663) -- (1.7443,-0.8046) -- (1.4703,-0.8639) -- cycle;
\path[cbpmesh=14.4] (-2.0762,-0.8243) -- (-1.8040,-0.8959) -- (-1.8342,-0.9514) -- (-2.1109,-0.8787) -- cycle;
\path[cbpmesh=17.9] (1.8040,-0.8959) -- (2.0762,-0.8243) -- (2.1109,-0.8787) -- (1.8342,-0.9514) -- cycle;
\path[cbpmesh=17.7] (-0.0000,-0.7694) -- (0.2681,-0.7644) -- (0.2742,-0.8044) -- (-0.0000,-0.8094) -- cycle;
\path[cbpmesh=17.7] (-0.2681,-0.7644) -- (-0.0000,-0.7694) -- (-0.0000,-0.8094) -- (-0.2742,-0.8044) -- cycle;
\path[cbpmesh=23.8] (-2.3672,-1.7500) -- (-2.0568,-1.8316) -- (-2.0588,-1.9422) -- (-2.3695,-1.8605) -- cycle;
\path[cbpmesh=32.0] (2.0568,-1.8316) -- (2.3672,-1.7500) -- (2.3695,-1.8605) -- (2.0588,-1.9422) -- cycle;
\path[cbpmesh=17.7] (-0.8159,-0.7641) -- (-0.5467,-0.7892) -- (-0.5583,-0.8265) -- (-0.8331,-0.8009) -- cycle;
\path[cbpmesh=17.9] (0.5467,-0.7892) -- (0.8159,-0.7641) -- (0.8331,-0.8009) -- (0.5583,-0.8265) -- cycle;
\path[cbpmesh=17.7] (-1.1028,-0.7651) -- (-0.8331,-0.8009) -- (-0.8497,-0.8366) -- (-1.1248,-0.8002) -- cycle;
\path[cbpmesh=17.9] (0.8331,-0.8009) -- (1.1028,-0.7651) -- (1.1248,-0.8002) -- (0.8497,-0.8366) -- cycle;
\path[cbpmesh=38] (-0.3403,-2.5339) -- (-0.0000,-2.5402) -- (-0.0000,-2.5763) -- (-0.3344,-2.5701) -- cycle;
\path[cbpmesh=38] (-0.0000,-2.5402) -- (0.3403,-2.5339) -- (0.3344,-2.5701) -- (-0.0000,-2.5763) -- cycle;
\path[cbpmesh=29.7] (-2.0534,-2.0410) -- (-1.7309,-2.1108) -- (-1.7206,-2.1973) -- (-2.0412,-2.1279) -- cycle;
\path[cbpmesh=36.5] (1.7309,-2.1108) -- (2.0534,-2.0410) -- (2.0412,-2.1279) -- (1.7206,-2.1973) -- cycle;
\path[cbpmesh=14.1] (-2.1109,-0.8787) -- (-1.8342,-0.9514) -- (-1.8650,-1.0154) -- (-2.1464,-0.9414) -- cycle;
\path[cbpmesh=18.3] (1.8342,-0.9514) -- (2.1109,-0.8787) -- (2.1464,-0.9414) -- (1.8650,-1.0154) -- cycle;
\path[cbpmesh=15.5] (-1.7443,-0.8046) -- (-1.4703,-0.8639) -- (-1.4954,-0.9077) -- (-1.7741,-0.8474) -- cycle;
\path[cbpmesh=17.3] (1.4703,-0.8639) -- (1.7443,-0.8046) -- (1.7741,-0.8474) -- (1.4954,-0.9077) -- cycle;
\path[cbpmesh=16.9] (-1.4193,-0.7883) -- (-1.1461,-0.8356) -- (-1.1668,-0.8727) -- (-1.4450,-0.8245) -- cycle;
\path[cbpmesh=17.6] (1.1461,-0.8356) -- (1.4193,-0.7883) -- (1.4450,-0.8245) -- (1.1668,-0.8727) -- cycle;
\path[cbpmesh=38] (-1.0276,-2.4347) -- (-0.6886,-2.4664) -- (-0.6785,-2.5151) -- (-1.0125,-2.4839) -- cycle;
\path[cbpmesh=38] (0.6886,-2.4664) -- (1.0276,-2.4347) -- (1.0125,-2.4839) -- (0.6785,-2.5151) -- cycle;
\path[cbpmesh=17.8] (-0.5467,-0.7892) -- (-0.2742,-0.8044) -- (-0.2800,-0.8420) -- (-0.5583,-0.8265) -- cycle;
\path[cbpmesh=17.9] (0.2742,-0.8044) -- (0.5467,-0.7892) -- (0.5583,-0.8265) -- (0.2800,-0.8420) -- cycle;
\path[cbpmesh=21.3] (-2.3557,-1.6280) -- (-2.0469,-1.7092) -- (-2.0568,-1.8316) -- (-2.3672,-1.7500) -- cycle;
\path[cbpmesh=29.4] (2.0469,-1.7092) -- (2.3557,-1.6280) -- (2.3672,-1.7500) -- (2.0568,-1.8316) -- cycle;
\path[cbpmesh=33.0] (-1.7206,-2.1973) -- (-1.3894,-2.2547) -- (-1.3768,-2.3287) -- (-1.7050,-2.2719) -- cycle;
\path[cbpmesh=38] (1.3894,-2.2547) -- (1.7206,-2.1973) -- (1.7050,-2.2719) -- (1.3768,-2.3287) -- cycle;
\path[cbpmesh=14.1] (-2.1464,-0.9414) -- (-1.8650,-1.0154) -- (-1.8966,-1.0891) -- (-2.1828,-1.0139) -- cycle;
\path[cbpmesh=18.9] (1.8650,-1.0154) -- (2.1464,-0.9414) -- (2.1828,-1.0139) -- (1.8966,-1.0891) -- cycle;
\path[cbpmesh=17.8] (-0.2742,-0.8044) -- (-0.0000,-0.8094) -- (-0.0000,-0.8471) -- (-0.2800,-0.8420) -- cycle;
\path[cbpmesh=17.8] (-0.0000,-0.8094) -- (0.2742,-0.8044) -- (0.2800,-0.8420) -- (-0.0000,-0.8471) -- cycle;
\path[cbpmesh=36.4] (-1.3768,-2.3287) -- (-1.0401,-2.3733) -- (-1.0276,-2.4347) -- (-1.3603,-2.3907) -- cycle;
\path[cbpmesh=38] (1.0401,-2.3733) -- (1.3768,-2.3287) -- (1.3603,-2.3907) -- (1.0276,-2.4347) -- cycle;
\path[cbpmesh=17.4] (-1.1248,-0.8002) -- (-0.8497,-0.8366) -- (-0.8658,-0.8727) -- (-1.1461,-0.8356) -- cycle;
\path[cbpmesh=17.7] (0.8497,-0.8366) -- (1.1248,-0.8002) -- (1.1461,-0.8356) -- (0.8658,-0.8727) -- cycle;
\path[cbpmesh=17.7] (-0.8331,-0.8009) -- (-0.5583,-0.8265) -- (-0.5694,-0.8628) -- (-0.8497,-0.8366) -- cycle;
\path[cbpmesh=17.9] (0.5583,-0.8265) -- (0.8331,-0.8009) -- (0.8497,-0.8366) -- (0.5694,-0.8628) -- cycle;
\path[cbpmesh=14.9] (-1.7741,-0.8474) -- (-1.4954,-0.9077) -- (-1.5206,-0.9572) -- (-1.8040,-0.8959) -- cycle;
\path[cbpmesh=17.3] (1.4954,-0.9077) -- (1.7741,-0.8474) -- (1.8040,-0.8959) -- (1.5206,-0.9572) -- cycle;
\path[cbpmesh=26.7] (-2.0588,-1.9422) -- (-1.7355,-2.0122) -- (-1.7309,-2.1108) -- (-2.0534,-2.0410) -- cycle;
\path[cbpmesh=33.7] (1.7355,-2.0122) -- (2.0588,-1.9422) -- (2.0534,-2.0410) -- (1.7309,-2.1108) -- cycle;
\path[cbpmesh=14.2] (-2.1828,-1.0139) -- (-1.8966,-1.0891) -- (-1.9294,-1.1739) -- (-2.2205,-1.0973) -- cycle;
\path[cbpmesh=19.6] (1.8966,-1.0891) -- (2.1828,-1.0139) -- (2.2205,-1.0973) -- (1.9294,-1.1739) -- cycle;
\path[cbpmesh=19.3] (-2.3345,-1.4945) -- (-2.0284,-1.5750) -- (-2.0469,-1.7092) -- (-2.3557,-1.6280) -- cycle;
\path[cbpmesh=27.1] (2.0284,-1.5750) -- (2.3345,-1.4945) -- (2.3557,-1.6280) -- (2.0469,-1.7092) -- cycle;
\path[cbpmesh=38] (-0.6886,-2.4664) -- (-0.3454,-2.4854) -- (-0.3403,-2.5339) -- (-0.6785,-2.5151) -- cycle;
\path[cbpmesh=38] (0.3454,-2.4854) -- (0.6886,-2.4664) -- (0.6785,-2.5151) -- (0.3403,-2.5339) -- cycle;
\path[cbpmesh=16.3] (-1.4450,-0.8245) -- (-1.1668,-0.8727) -- (-1.1873,-0.9129) -- (-1.4703,-0.8639) -- cycle;
\path[cbpmesh=17.3] (1.1668,-0.8727) -- (1.4450,-0.8245) -- (1.4703,-0.8639) -- (1.1873,-0.9129) -- cycle;
\path[cbpmesh=17.7] (-0.5583,-0.8265) -- (-0.2800,-0.8420) -- (-0.2856,-0.8785) -- (-0.5694,-0.8628) -- cycle;
\path[cbpmesh=17.8] (0.2800,-0.8420) -- (0.5583,-0.8265) -- (0.5694,-0.8628) -- (0.2856,-0.8785) -- cycle;
\path[cbpmesh=14.5] (-2.2205,-1.0973) -- (-1.9294,-1.1739) -- (-1.9635,-1.2710) -- (-2.2598,-1.1931) -- cycle;
\path[cbpmesh=20.2] (1.9294,-1.1739) -- (2.2205,-1.0973) -- (2.2598,-1.1931) -- (1.9635,-1.2710) -- cycle;
\path[cbpmesh=17.6] (-2.3027,-1.3495) -- (-2.0008,-1.4289) -- (-2.0284,-1.5750) -- (-2.3345,-1.4945) -- cycle;
\path[cbpmesh=25.1] (2.0008,-1.4289) -- (2.3027,-1.3495) -- (2.3345,-1.4945) -- (2.0284,-1.5750) -- cycle;
\path[cbpmesh=38] (-0.0000,-2.4918) -- (0.3454,-2.4854) -- (0.3403,-2.5339) -- (-0.0000,-2.5402) -- cycle;
\path[cbpmesh=38] (-0.3454,-2.4854) -- (-0.0000,-2.4918) -- (-0.0000,-2.5402) -- (-0.3403,-2.5339) -- cycle;
\path[cbpmesh=14.5] (-1.8040,-0.8959) -- (-1.5206,-0.9572) -- (-1.5461,-1.0138) -- (-1.8342,-0.9514) -- cycle;
\path[cbpmesh=17.5] (1.5206,-0.9572) -- (1.8040,-0.8959) -- (1.8342,-0.9514) -- (1.5461,-1.0138) -- cycle;
\path[cbpmesh=16.9] (-1.1461,-0.8356) -- (-0.8658,-0.8727) -- (-0.8815,-0.9105) -- (-1.1668,-0.8727) -- cycle;
\path[cbpmesh=17.5] (0.8658,-0.8727) -- (1.1461,-0.8356) -- (1.1668,-0.8727) -- (0.8815,-0.9105) -- cycle;
\path[cbpmesh=16.4] (-2.2598,-1.1931) -- (-1.9635,-1.2710) -- (-2.0008,-1.4289) -- (-2.3027,-1.3495) -- cycle;
\path[cbpmesh=23.4] (1.9635,-1.2710) -- (2.2598,-1.1931) -- (2.3027,-1.3495) -- (2.0008,-1.4289) -- cycle;
\path[cbpmesh=17.4] (-0.8497,-0.8366) -- (-0.5694,-0.8628) -- (-0.5802,-0.8994) -- (-0.8658,-0.8727) -- cycle;
\path[cbpmesh=17.7] (0.5694,-0.8628) -- (0.8497,-0.8366) -- (0.8658,-0.8727) -- (0.5802,-0.8994) -- cycle;
\path[cbpmesh=17.7] (-0.2800,-0.8420) -- (-0.0000,-0.8471) -- (-0.0000,-0.8838) -- (-0.2856,-0.8785) -- cycle;
\path[cbpmesh=17.8] (-0.0000,-0.8471) -- (0.2800,-0.8420) -- (0.2856,-0.8785) -- (-0.0000,-0.8838) -- cycle;
\path[cbpmesh=29.9] (-1.7309,-2.1108) -- (-1.3977,-2.1685) -- (-1.3894,-2.2547) -- (-1.7206,-2.1973) -- cycle;
\path[cbpmesh=35.6] (1.3977,-2.1685) -- (1.7309,-2.1108) -- (1.7206,-2.1973) -- (1.3894,-2.2547) -- cycle;
\path[cbpmesh=36.7] (-1.0401,-2.3733) -- (-0.6970,-2.4054) -- (-0.6886,-2.4664) -- (-1.0276,-2.4347) -- cycle;
\path[cbpmesh=38] (0.6970,-2.4054) -- (1.0401,-2.3733) -- (1.0276,-2.4347) -- (0.6886,-2.4664) -- cycle;
\path[cbpmesh=15.6] (-1.4703,-0.8639) -- (-1.1873,-0.9129) -- (-1.2076,-0.9575) -- (-1.4954,-0.9077) -- cycle;
\path[cbpmesh=17.1] (1.1873,-0.9129) -- (1.4703,-0.8639) -- (1.4954,-0.9077) -- (1.2076,-0.9575) -- cycle;
\path[cbpmesh=24.0] (-2.0568,-1.8316) -- (-1.7338,-1.9015) -- (-1.7355,-2.0122) -- (-2.0588,-1.9422) -- cycle;
\path[cbpmesh=31.0] (1.7338,-1.9015) -- (2.0568,-1.8316) -- (2.0588,-1.9422) -- (1.7355,-2.0122) -- cycle;
\path[cbpmesh=33.3] (-1.3894,-2.2547) -- (-1.0496,-2.2997) -- (-1.0401,-2.3733) -- (-1.3768,-2.3287) -- cycle;
\path[cbpmesh=37.6] (1.0496,-2.2997) -- (1.3894,-2.2547) -- (1.3768,-2.3287) -- (1.0401,-2.3733) -- cycle;
\path[cbpmesh=14.2] (-1.8342,-0.9514) -- (-1.5461,-1.0138) -- (-1.5720,-1.0788) -- (-1.8650,-1.0154) -- cycle;
\path[cbpmesh=17.8] (1.5461,-1.0138) -- (1.8342,-0.9514) -- (1.8650,-1.0154) -- (1.5720,-1.0788) -- cycle;
\path[cbpmesh=17.5] (-0.5694,-0.8628) -- (-0.2856,-0.8785) -- (-0.2910,-0.9154) -- (-0.5802,-0.8994) -- cycle;
\path[cbpmesh=17.6] (0.2856,-0.8785) -- (0.5694,-0.8628) -- (0.5802,-0.8994) -- (0.2910,-0.9154) -- cycle;
\path[cbpmesh=16.3] (-1.1668,-0.8727) -- (-0.8815,-0.9105) -- (-0.8969,-0.9514) -- (-1.1873,-0.9129) -- cycle;
\path[cbpmesh=17.2] (0.8815,-0.9105) -- (1.1668,-0.8727) -- (1.1873,-0.9129) -- (0.8969,-0.9514) -- cycle;
\path[cbpmesh=17.0] (-0.8658,-0.8727) -- (-0.5802,-0.8994) -- (-0.5907,-0.9377) -- (-0.8815,-0.9105) -- cycle;
\path[cbpmesh=17.4] (0.5802,-0.8994) -- (0.8658,-0.8727) -- (0.8815,-0.9105) -- (0.5907,-0.9377) -- cycle;
\path[cbpmesh=15.0] (-1.4954,-0.9077) -- (-1.2076,-0.9575) -- (-1.2279,-1.0079) -- (-1.5206,-0.9572) -- cycle;
\path[cbpmesh=17.0] (1.2076,-0.9575) -- (1.4954,-0.9077) -- (1.5206,-0.9572) -- (1.2279,-1.0079) -- cycle;
\path[cbpmesh=17.5] (-0.2856,-0.8785) -- (-0.0000,-0.8838) -- (-0.0000,-0.9208) -- (-0.2910,-0.9154) -- cycle;
\path[cbpmesh=17.5] (-0.0000,-0.8838) -- (0.2856,-0.8785) -- (0.2910,-0.9154) -- (-0.0000,-0.9208) -- cycle;
\path[cbpmesh=37.2] (-0.6970,-2.4054) -- (-0.3496,-2.4247) -- (-0.3454,-2.4854) -- (-0.6886,-2.4664) -- cycle;
\path[cbpmesh=38] (0.3496,-2.4247) -- (0.6970,-2.4054) -- (0.6886,-2.4664) -- (0.3454,-2.4854) -- cycle;
\path[cbpmesh=14.2] (-1.8650,-1.0154) -- (-1.5720,-1.0788) -- (-1.5987,-1.1536) -- (-1.8966,-1.0891) -- cycle;
\path[cbpmesh=18.3] (1.5720,-1.0788) -- (1.8650,-1.0154) -- (1.8966,-1.0891) -- (1.5987,-1.1536) -- cycle;
\path[cbpmesh=21.5] (-2.0469,-1.7092) -- (-1.7254,-1.7788) -- (-1.7338,-1.9015) -- (-2.0568,-1.8316) -- cycle;
\path[cbpmesh=28.5] (1.7254,-1.7788) -- (2.0469,-1.7092) -- (2.0568,-1.8316) -- (1.7338,-1.9015) -- cycle;
\path[cbpmesh=27.0] (-1.7355,-2.0122) -- (-1.4014,-2.0700) -- (-1.3977,-2.1685) -- (-1.7309,-2.1108) -- cycle;
\path[cbpmesh=32.8] (1.4014,-2.0700) -- (1.7355,-2.0122) -- (1.7309,-2.1108) -- (1.3977,-2.1685) -- cycle;
\path[cbpmesh=37.7] (-0.3496,-2.4247) -- (-0.0000,-2.4311) -- (-0.0000,-2.4918) -- (-0.3454,-2.4854) -- cycle;
\path[cbpmesh=38] (-0.0000,-2.4311) -- (0.3496,-2.4247) -- (0.3454,-2.4854) -- (-0.0000,-2.4918) -- cycle;
\path[cbpmesh=14.4] (-1.8966,-1.0891) -- (-1.5987,-1.1536) -- (-1.6263,-1.2394) -- (-1.9294,-1.1739) -- cycle;
\path[cbpmesh=18.9] (1.5987,-1.1536) -- (1.8966,-1.0891) -- (1.9294,-1.1739) -- (1.6263,-1.2394) -- cycle;
\path[cbpmesh=15.7] (-1.1873,-0.9129) -- (-0.8969,-0.9514) -- (-0.9123,-0.9966) -- (-1.2076,-0.9575) -- cycle;
\path[cbpmesh=16.9] (0.8969,-0.9514) -- (1.1873,-0.9129) -- (1.2076,-0.9575) -- (0.9123,-0.9966) -- cycle;
\path[cbpmesh=17.0] (-0.5802,-0.8994) -- (-0.2910,-0.9154) -- (-0.2963,-0.9540) -- (-0.5907,-0.9377) -- cycle;
\path[cbpmesh=17.3] (0.2910,-0.9154) -- (0.5802,-0.8994) -- (0.5907,-0.9377) -- (0.2963,-0.9540) -- cycle;
\path[cbpmesh=14.6] (-1.5206,-0.9572) -- (-1.2279,-1.0079) -- (-1.2485,-1.0653) -- (-1.5461,-1.0138) -- cycle;
\path[cbpmesh=17.1] (1.2279,-1.0079) -- (1.5206,-0.9572) -- (1.5461,-1.0138) -- (1.2485,-1.0653) -- cycle;
\path[cbpmesh=33.7] (-1.0496,-2.2997) -- (-0.7033,-2.3320) -- (-0.6970,-2.4054) -- (-1.0401,-2.3733) -- cycle;
\path[cbpmesh=36.8] (0.7033,-2.3320) -- (1.0496,-2.2997) -- (1.0401,-2.3733) -- (0.6970,-2.4054) -- cycle;
\path[cbpmesh=19.5] (-2.0284,-1.5750) -- (-1.7098,-1.6439) -- (-1.7254,-1.7788) -- (-2.0469,-1.7092) -- cycle;
\path[cbpmesh=26.2] (1.7098,-1.6439) -- (2.0284,-1.5750) -- (2.0469,-1.7092) -- (1.7254,-1.7788) -- cycle;
\path[cbpmesh=30.3] (-1.3977,-2.1685) -- (-1.0559,-2.2137) -- (-1.0496,-2.2997) -- (-1.3894,-2.2547) -- cycle;
\path[cbpmesh=34.7] (1.0559,-2.2137) -- (1.3977,-2.1685) -- (1.3894,-2.2547) -- (1.0496,-2.2997) -- cycle;
\path[cbpmesh=16.4] (-0.8815,-0.9105) -- (-0.5907,-0.9377) -- (-0.6010,-0.9790) -- (-0.8969,-0.9514) -- cycle;
\path[cbpmesh=17.0] (0.5907,-0.9377) -- (0.8815,-0.9105) -- (0.8969,-0.9514) -- (0.6010,-0.9790) -- cycle;
\path[cbpmesh=17.1] (-0.2910,-0.9154) -- (-0.0000,-0.9208) -- (-0.0000,-0.9595) -- (-0.2963,-0.9540) -- cycle;
\path[cbpmesh=17.2] (-0.0000,-0.9208) -- (0.2910,-0.9154) -- (0.2963,-0.9540) -- (-0.0000,-0.9595) -- cycle;
\path[cbpmesh=14.6] (-1.9294,-1.1739) -- (-1.6263,-1.2394) -- (-1.6551,-1.3377) -- (-1.9635,-1.2710) -- cycle;
\path[cbpmesh=19.6] (1.6263,-1.2394) -- (1.9294,-1.1739) -- (1.9635,-1.2710) -- (1.6551,-1.3377) -- cycle;
\path[cbpmesh=17.8] (-2.0008,-1.4289) -- (-1.6865,-1.4969) -- (-1.7098,-1.6439) -- (-2.0284,-1.5750) -- cycle;
\path[cbpmesh=24.2] (1.6865,-1.4969) -- (2.0008,-1.4289) -- (2.0284,-1.5750) -- (1.7098,-1.6439) -- cycle;
\path[cbpmesh=16.5] (-1.9635,-1.2710) -- (-1.6551,-1.3377) -- (-1.6865,-1.4969) -- (-2.0008,-1.4289) -- cycle;
\path[cbpmesh=22.6] (1.6551,-1.3377) -- (1.9635,-1.2710) -- (2.0008,-1.4289) -- (1.6865,-1.4969) -- cycle;
\path[cbpmesh=14.4] (-1.5461,-1.0138) -- (-1.2485,-1.0653) -- (-1.2694,-1.1312) -- (-1.5720,-1.0788) -- cycle;
\path[cbpmesh=17.4] (1.2485,-1.0653) -- (1.5461,-1.0138) -- (1.5720,-1.0788) -- (1.2694,-1.1312) -- cycle;
\path[cbpmesh=24.2] (-1.7338,-1.9015) -- (-1.4000,-1.9593) -- (-1.4014,-2.0700) -- (-1.7355,-2.0122) -- cycle;
\path[cbpmesh=30.1] (1.4000,-1.9593) -- (1.7338,-1.9015) -- (1.7355,-2.0122) -- (1.4014,-2.0700) -- cycle;
\path[cbpmesh=15.2] (-1.2076,-0.9575) -- (-0.9123,-0.9966) -- (-0.9276,-1.0476) -- (-1.2279,-1.0079) -- cycle;
\path[cbpmesh=16.7] (0.9123,-0.9966) -- (1.2076,-0.9575) -- (1.2279,-1.0079) -- (0.9276,-1.0476) -- cycle;
\path[cbpmesh=34.2] (-0.7033,-2.3320) -- (-0.3528,-2.3515) -- (-0.3496,-2.4247) -- (-0.6970,-2.4054) -- cycle;
\path[cbpmesh=36.1] (0.3528,-2.3515) -- (0.7033,-2.3320) -- (0.6970,-2.4054) -- (0.3496,-2.4247) -- cycle;
\path[cbpmesh=16.5] (-0.5907,-0.9377) -- (-0.2963,-0.9540) -- (-0.3014,-0.9956) -- (-0.6010,-0.9790) -- cycle;
\path[cbpmesh=16.9] (0.2963,-0.9540) -- (0.5907,-0.9377) -- (0.6010,-0.9790) -- (0.3014,-0.9956) -- cycle;
\path[cbpmesh=15.8] (-0.8969,-0.9514) -- (-0.6010,-0.9790) -- (-0.6113,-1.0247) -- (-0.9123,-0.9966) -- cycle;
\path[cbpmesh=16.7] (0.6010,-0.9790) -- (0.8969,-0.9514) -- (0.9123,-0.9966) -- (0.6113,-1.0247) -- cycle;
\path[cbpmesh=16.6] (-0.2963,-0.9540) -- (-0.0000,-0.9595) -- (-0.0000,-1.0012) -- (-0.3014,-0.9956) -- cycle;
\path[cbpmesh=16.7] (-0.0000,-0.9595) -- (0.2963,-0.9540) -- (0.3014,-0.9956) -- (-0.0000,-1.0012) -- cycle;
\path[cbpmesh=34.8] (-0.3528,-2.3515) -- (-0.0000,-2.3580) -- (-0.0000,-2.4311) -- (-0.3496,-2.4247) -- cycle;
\path[cbpmesh=35.4] (-0.0000,-2.3580) -- (0.3528,-2.3515) -- (0.3496,-2.4247) -- (-0.0000,-2.4311) -- cycle;
\path[cbpmesh=14.4] (-1.5720,-1.0788) -- (-1.2694,-1.1312) -- (-1.2910,-1.2069) -- (-1.5987,-1.1536) -- cycle;
\path[cbpmesh=17.8] (1.2694,-1.1312) -- (1.5720,-1.0788) -- (1.5987,-1.1536) -- (1.2910,-1.2069) -- cycle;
\path[cbpmesh=27.4] (-1.4014,-2.0700) -- (-1.0587,-2.1154) -- (-1.0559,-2.2137) -- (-1.3977,-2.1685) -- cycle;
\path[cbpmesh=31.9] (1.0587,-2.1154) -- (1.4014,-2.0700) -- (1.3977,-2.1685) -- (1.0559,-2.2137) -- cycle;
\path[cbpmesh=30.7] (-1.0559,-2.2137) -- (-0.7076,-2.2463) -- (-0.7033,-2.3320) -- (-1.0496,-2.2997) -- cycle;
\path[cbpmesh=33.9] (0.7076,-2.2463) -- (1.0559,-2.2137) -- (1.0496,-2.2997) -- (0.7033,-2.3320) -- cycle;
\path[cbpmesh=14.8] (-1.2279,-1.0079) -- (-0.9276,-1.0476) -- (-0.9431,-1.1058) -- (-1.2485,-1.0653) -- cycle;
\path[cbpmesh=16.7] (0.9276,-1.0476) -- (1.2279,-1.0079) -- (1.2485,-1.0653) -- (0.9431,-1.1058) -- cycle;
\path[cbpmesh=21.8] (-1.7254,-1.7788) -- (-1.3932,-1.8363) -- (-1.4000,-1.9593) -- (-1.7338,-1.9015) -- cycle;
\path[cbpmesh=27.5] (1.3932,-1.8363) -- (1.7254,-1.7788) -- (1.7338,-1.9015) -- (1.4000,-1.9593) -- cycle;
\path[cbpmesh=16.0] (-0.6010,-0.9790) -- (-0.3014,-0.9956) -- (-0.3066,-1.0417) -- (-0.6113,-1.0247) -- cycle;
\path[cbpmesh=16.5] (0.3014,-0.9956) -- (0.6010,-0.9790) -- (0.6113,-1.0247) -- (0.3066,-1.0417) -- cycle;
\path[cbpmesh=15.3] (-0.9123,-0.9966) -- (-0.6113,-1.0247) -- (-0.6216,-1.0762) -- (-0.9276,-1.0476) -- cycle;
\path[cbpmesh=16.4] (0.6113,-1.0247) -- (0.9123,-0.9966) -- (0.9276,-1.0476) -- (0.6216,-1.0762) -- cycle;
\path[cbpmesh=14.5] (-1.5987,-1.1536) -- (-1.2910,-1.2069) -- (-1.3133,-1.2937) -- (-1.6263,-1.2394) -- cycle;
\path[cbpmesh=18.3] (1.2910,-1.2069) -- (1.5987,-1.1536) -- (1.6263,-1.2394) -- (1.3133,-1.2937) -- cycle;
\path[cbpmesh=16.1] (-0.3014,-0.9956) -- (-0.0000,-1.0012) -- (-0.0000,-1.0473) -- (-0.3066,-1.0417) -- cycle;
\path[cbpmesh=16.3] (-0.0000,-1.0012) -- (0.3014,-0.9956) -- (0.3066,-1.0417) -- (-0.0000,-1.0473) -- cycle;
\path[cbpmesh=19.7] (-1.7098,-1.6439) -- (-1.3806,-1.7010) -- (-1.3932,-1.8363) -- (-1.7254,-1.7788) -- cycle;
\path[cbpmesh=25.3] (1.3806,-1.7010) -- (1.7098,-1.6439) -- (1.7254,-1.7788) -- (1.3932,-1.8363) -- cycle;
\path[cbpmesh=14.6] (-1.2485,-1.0653) -- (-0.9431,-1.1058) -- (-0.9590,-1.1723) -- (-1.2694,-1.1312) -- cycle;
\path[cbpmesh=16.9] (0.9431,-1.1058) -- (1.2485,-1.0653) -- (1.2694,-1.1312) -- (0.9590,-1.1723) -- cycle;
\path[cbpmesh=31.2] (-0.7076,-2.2463) -- (-0.3549,-2.2659) -- (-0.3528,-2.3515) -- (-0.7033,-2.3320) -- cycle;
\path[cbpmesh=33.2] (0.3549,-2.2659) -- (0.7076,-2.2463) -- (0.7033,-2.3320) -- (0.3528,-2.3515) -- cycle;
\path[cbpmesh=14.8] (-1.6263,-1.2394) -- (-1.3133,-1.2937) -- (-1.3365,-1.3929) -- (-1.6551,-1.3377) -- cycle;
\path[cbpmesh=18.9] (1.3133,-1.2937) -- (1.6263,-1.2394) -- (1.6551,-1.3377) -- (1.3365,-1.3929) -- cycle;
\path[cbpmesh=24.6] (-1.4000,-1.9593) -- (-1.0576,-2.0047) -- (-1.0587,-2.1154) -- (-1.4014,-2.0700) -- cycle;
\path[cbpmesh=29.2] (1.0576,-2.0047) -- (1.4000,-1.9593) -- (1.4014,-2.0700) -- (1.0587,-2.1154) -- cycle;
\path[cbpmesh=18.1] (-1.6865,-1.4969) -- (-1.3619,-1.5532) -- (-1.3806,-1.7010) -- (-1.7098,-1.6439) -- cycle;
\path[cbpmesh=23.4] (1.3619,-1.5532) -- (1.6865,-1.4969) -- (1.7098,-1.6439) -- (1.3806,-1.7010) -- cycle;
\path[cbpmesh=15.5] (-0.6113,-1.0247) -- (-0.3066,-1.0417) -- (-0.3118,-1.0934) -- (-0.6216,-1.0762) -- cycle;
\path[cbpmesh=16.2] (0.3066,-1.0417) -- (0.6113,-1.0247) -- (0.6216,-1.0762) -- (0.3118,-1.0934) -- cycle;
\path[cbpmesh=16.8] (-1.6551,-1.3377) -- (-1.3365,-1.3929) -- (-1.3619,-1.5532) -- (-1.6865,-1.4969) -- cycle;
\path[cbpmesh=21.8] (1.3365,-1.3929) -- (1.6551,-1.3377) -- (1.6865,-1.4969) -- (1.3619,-1.5532) -- cycle;
\path[cbpmesh=15.0] (-0.9276,-1.0476) -- (-0.6216,-1.0762) -- (-0.6320,-1.1348) -- (-0.9431,-1.1058) -- cycle;
\path[cbpmesh=16.4] (0.6216,-1.0762) -- (0.9276,-1.0476) -- (0.9431,-1.1058) -- (0.6320,-1.1348) -- cycle;
\path[cbpmesh=27.8] (-1.0587,-2.1154) -- (-0.7094,-2.1480) -- (-0.7076,-2.2463) -- (-1.0559,-2.2137) -- cycle;
\path[cbpmesh=31.1] (0.7094,-2.1480) -- (1.0587,-2.1154) -- (1.0559,-2.2137) -- (0.7076,-2.2463) -- cycle;
\path[cbpmesh=31.8] (-0.3549,-2.2659) -- (-0.0000,-2.2724) -- (-0.0000,-2.3580) -- (-0.3528,-2.3515) -- cycle;
\path[cbpmesh=32.5] (-0.0000,-2.2724) -- (0.3549,-2.2659) -- (0.3528,-2.3515) -- (-0.0000,-2.3580) -- cycle;
\path[cbpmesh=14.6] (-1.2694,-1.1312) -- (-0.9590,-1.1723) -- (-0.9752,-1.2487) -- (-1.2910,-1.2069) -- cycle;
\path[cbpmesh=17.3] (0.9590,-1.1723) -- (1.2694,-1.1312) -- (1.2910,-1.2069) -- (0.9752,-1.2487) -- cycle;
\path[cbpmesh=15.7] (-0.3066,-1.0417) -- (-0.0000,-1.0473) -- (-0.0000,-1.0992) -- (-0.3118,-1.0934) -- cycle;
\path[cbpmesh=15.9] (-0.0000,-1.0473) -- (0.3066,-1.0417) -- (0.3118,-1.0934) -- (-0.0000,-1.0992) -- cycle;
\path[cbpmesh=22.2] (-1.3932,-1.8363) -- (-1.0525,-1.8814) -- (-1.0576,-2.0047) -- (-1.4000,-1.9593) -- cycle;
\path[cbpmesh=26.7] (1.0525,-1.8814) -- (1.3932,-1.8363) -- (1.4000,-1.9593) -- (1.0576,-2.0047) -- cycle;
\path[cbpmesh=14.8] (-0.9431,-1.1058) -- (-0.6320,-1.1348) -- (-0.6426,-1.2019) -- (-0.9590,-1.1723) -- cycle;
\path[cbpmesh=16.5] (0.6320,-1.1348) -- (0.9431,-1.1058) -- (0.9590,-1.1723) -- (0.6426,-1.2019) -- cycle;
\path[cbpmesh=15.2] (-0.6216,-1.0762) -- (-0.3118,-1.0934) -- (-0.3170,-1.1523) -- (-0.6320,-1.1348) -- cycle;
\path[cbpmesh=16.0] (0.3118,-1.0934) -- (0.6216,-1.0762) -- (0.6320,-1.1348) -- (0.3170,-1.1523) -- cycle;
\path[cbpmesh=14.8] (-1.2910,-1.2069) -- (-0.9752,-1.2487) -- (-0.9921,-1.3362) -- (-1.3133,-1.2937) -- cycle;
\path[cbpmesh=17.8] (0.9752,-1.2487) -- (1.2910,-1.2069) -- (1.3133,-1.2937) -- (0.9921,-1.3362) -- cycle;
\path[cbpmesh=28.3] (-0.7094,-2.1480) -- (-0.3558,-2.1677) -- (-0.3549,-2.2659) -- (-0.7076,-2.2463) -- cycle;
\path[cbpmesh=30.3] (0.3558,-2.1677) -- (0.7094,-2.1480) -- (0.7076,-2.2463) -- (0.3549,-2.2659) -- cycle;
\path[cbpmesh=25.1] (-1.0576,-2.0047) -- (-0.7087,-2.0373) -- (-0.7094,-2.1480) -- (-1.0587,-2.1154) -- cycle;
\path[cbpmesh=28.3] (0.7087,-2.0373) -- (1.0576,-2.0047) -- (1.0587,-2.1154) -- (0.7094,-2.1480) -- cycle;
\path[cbpmesh=15.4] (-0.3118,-1.0934) -- (-0.0000,-1.0992) -- (-0.0000,-1.1582) -- (-0.3170,-1.1523) -- cycle;
\path[cbpmesh=15.7] (-0.0000,-1.0992) -- (0.3118,-1.0934) -- (0.3170,-1.1523) -- (-0.0000,-1.1582) -- cycle;
\path[cbpmesh=20.1] (-1.3806,-1.7010) -- (-1.0430,-1.7457) -- (-1.0525,-1.8814) -- (-1.3932,-1.8363) -- cycle;
\path[cbpmesh=24.5] (1.0430,-1.7457) -- (1.3806,-1.7010) -- (1.3932,-1.8363) -- (1.0525,-1.8814) -- cycle;
\path[cbpmesh=28.9] (-0.3558,-2.1677) -- (-0.0000,-2.1742) -- (-0.0000,-2.2724) -- (-0.3549,-2.2659) -- cycle;
\path[cbpmesh=29.6] (-0.0000,-2.1742) -- (0.3558,-2.1677) -- (0.3549,-2.2659) -- (-0.0000,-2.2724) -- cycle;
\path[cbpmesh=15.1] (-1.3133,-1.2937) -- (-0.9921,-1.3362) -- (-1.0097,-1.4362) -- (-1.3365,-1.3929) -- cycle;
\path[cbpmesh=18.3] (0.9921,-1.3362) -- (1.3133,-1.2937) -- (1.3365,-1.3929) -- (1.0097,-1.4362) -- cycle;
\path[cbpmesh=14.9] (-0.9590,-1.1723) -- (-0.6426,-1.2019) -- (-0.6535,-1.2787) -- (-0.9752,-1.2487) -- cycle;
\path[cbpmesh=16.8] (0.6426,-1.2019) -- (0.9590,-1.1723) -- (0.9752,-1.2487) -- (0.6535,-1.2787) -- cycle;
\path[cbpmesh=15.1] (-0.6320,-1.1348) -- (-0.3170,-1.1523) -- (-0.3223,-1.2197) -- (-0.6426,-1.2019) -- cycle;
\path[cbpmesh=16.1] (0.3170,-1.1523) -- (0.6320,-1.1348) -- (0.6426,-1.2019) -- (0.3223,-1.2197) -- cycle;
\path[cbpmesh=18.4] (-1.3619,-1.5532) -- (-1.0288,-1.5973) -- (-1.0430,-1.7457) -- (-1.3806,-1.7010) -- cycle;
\path[cbpmesh=22.6] (1.0288,-1.5973) -- (1.3619,-1.5532) -- (1.3806,-1.7010) -- (1.0430,-1.7457) -- cycle;
\path[cbpmesh=17.1] (-1.3365,-1.3929) -- (-1.0097,-1.4362) -- (-1.0288,-1.5973) -- (-1.3619,-1.5532) -- cycle;
\path[cbpmesh=21.0] (1.0097,-1.4362) -- (1.3365,-1.3929) -- (1.3619,-1.5532) -- (1.0288,-1.5973) -- cycle;
\path[cbpmesh=22.6] (-1.0525,-1.8814) -- (-0.7053,-1.9139) -- (-0.7087,-2.0373) -- (-1.0576,-2.0047) -- cycle;
\path[cbpmesh=25.8] (0.7053,-1.9139) -- (1.0525,-1.8814) -- (1.0576,-2.0047) -- (0.7087,-2.0373) -- cycle;
\path[cbpmesh=15.4] (-0.3170,-1.1523) -- (-0.0000,-1.1582) -- (-0.0000,-1.2256) -- (-0.3223,-1.2197) -- cycle;
\path[cbpmesh=15.7] (-0.0000,-1.1582) -- (0.3170,-1.1523) -- (0.3223,-1.2197) -- (-0.0000,-1.2256) -- cycle;
\path[cbpmesh=25.6] (-0.7087,-2.0373) -- (-0.3555,-2.0569) -- (-0.3558,-2.1677) -- (-0.7094,-2.1480) -- cycle;
\path[cbpmesh=27.6] (0.3555,-2.0569) -- (0.7087,-2.0373) -- (0.7094,-2.1480) -- (0.3558,-2.1677) -- cycle;
\path[cbpmesh=15.1] (-0.9752,-1.2487) -- (-0.6535,-1.2787) -- (-0.6648,-1.3668) -- (-0.9921,-1.3362) -- cycle;
\path[cbpmesh=17.2] (0.6535,-1.2787) -- (0.9752,-1.2487) -- (0.9921,-1.3362) -- (0.6648,-1.3668) -- cycle;
\path[cbpmesh=15.2] (-0.6426,-1.2019) -- (-0.3223,-1.2197) -- (-0.3278,-1.2968) -- (-0.6535,-1.2787) -- cycle;
\path[cbpmesh=16.3] (0.3223,-1.2197) -- (0.6426,-1.2019) -- (0.6535,-1.2787) -- (0.3278,-1.2968) -- cycle;
\path[cbpmesh=26.8] (-0.0000,-2.0634) -- (0.3555,-2.0569) -- (0.3558,-2.1677) -- (-0.0000,-2.1742) -- cycle;
\path[cbpmesh=26.2] (-0.3555,-2.0569) -- (-0.0000,-2.0634) -- (-0.0000,-2.1742) -- (-0.3558,-2.1677) -- cycle;
\path[cbpmesh=20.5] (-1.0430,-1.7457) -- (-0.6989,-1.7778) -- (-0.7053,-1.9139) -- (-1.0525,-1.8814) -- cycle;
\path[cbpmesh=23.7] (0.6989,-1.7778) -- (1.0430,-1.7457) -- (1.0525,-1.8814) -- (0.7053,-1.9139) -- cycle;
\path[cbpmesh=15.4] (-0.9921,-1.3362) -- (-0.6648,-1.3668) -- (-0.6766,-1.4673) -- (-1.0097,-1.4362) -- cycle;
\path[cbpmesh=17.7] (0.6648,-1.3668) -- (0.9921,-1.3362) -- (1.0097,-1.4362) -- (0.6766,-1.4673) -- cycle;
\path[cbpmesh=15.5] (-0.3223,-1.2197) -- (-0.0000,-1.2256) -- (-0.0000,-1.3029) -- (-0.3278,-1.2968) -- cycle;
\path[cbpmesh=15.9] (-0.0000,-1.2256) -- (0.3223,-1.2197) -- (0.3278,-1.2968) -- (-0.0000,-1.3029) -- cycle;
\path[cbpmesh=23.1] (-0.7053,-1.9139) -- (-0.3537,-1.9334) -- (-0.3555,-2.0569) -- (-0.7087,-2.0373) -- cycle;
\path[cbpmesh=25.1] (0.3537,-1.9334) -- (0.7053,-1.9139) -- (0.7087,-2.0373) -- (0.3555,-2.0569) -- cycle;
\path[cbpmesh=15.4] (-0.6535,-1.2787) -- (-0.3278,-1.2968) -- (-0.3334,-1.3852) -- (-0.6648,-1.3668) -- cycle;
\path[cbpmesh=16.7] (0.3278,-1.2968) -- (0.6535,-1.2787) -- (0.6648,-1.3668) -- (0.3334,-1.3852) -- cycle;
\path[cbpmesh=18.8] (-1.0288,-1.5973) -- (-0.6894,-1.6290) -- (-0.6989,-1.7778) -- (-1.0430,-1.7457) -- cycle;
\path[cbpmesh=21.8] (0.6894,-1.6290) -- (1.0288,-1.5973) -- (1.0430,-1.7457) -- (0.6989,-1.7778) -- cycle;
\path[cbpmesh=17.5] (-1.0097,-1.4362) -- (-0.6766,-1.4673) -- (-0.6894,-1.6290) -- (-1.0288,-1.5973) -- cycle;
\path[cbpmesh=20.3] (0.6766,-1.4673) -- (1.0097,-1.4362) -- (1.0288,-1.5973) -- (0.6894,-1.6290) -- cycle;
\path[cbpmesh=23.7] (-0.3537,-1.9334) -- (-0.0000,-1.9399) -- (-0.0000,-2.0634) -- (-0.3555,-2.0569) -- cycle;
\path[cbpmesh=24.4] (-0.0000,-1.9399) -- (0.3537,-1.9334) -- (0.3555,-2.0569) -- (-0.0000,-2.0634) -- cycle;
\path[cbpmesh=15.8] (-0.3278,-1.2968) -- (-0.0000,-1.3029) -- (-0.0000,-1.3913) -- (-0.3334,-1.3852) -- cycle;
\path[cbpmesh=16.2] (-0.0000,-1.3029) -- (0.3278,-1.2968) -- (0.3334,-1.3852) -- (-0.0000,-1.3913) -- cycle;
\path[cbpmesh=21.0] (-0.6989,-1.7778) -- (-0.3505,-1.7971) -- (-0.3537,-1.9334) -- (-0.7053,-1.9139) -- cycle;
\path[cbpmesh=22.9] (0.3505,-1.7971) -- (0.6989,-1.7778) -- (0.7053,-1.9139) -- (0.3537,-1.9334) -- cycle;
\path[cbpmesh=15.7] (-0.6648,-1.3668) -- (-0.3334,-1.3852) -- (-0.3393,-1.4860) -- (-0.6766,-1.4673) -- cycle;
\path[cbpmesh=17.1] (0.3334,-1.3852) -- (0.6648,-1.3668) -- (0.6766,-1.4673) -- (0.3393,-1.4860) -- cycle;
\path[cbpmesh=19.3] (-0.6894,-1.6290) -- (-0.3458,-1.6481) -- (-0.3505,-1.7971) -- (-0.6989,-1.7778) -- cycle;
\path[cbpmesh=21.1] (0.3458,-1.6481) -- (0.6894,-1.6290) -- (0.6989,-1.7778) -- (0.3505,-1.7971) -- cycle;
\path[cbpmesh=17.9] (-0.6766,-1.4673) -- (-0.3393,-1.4860) -- (-0.3458,-1.6481) -- (-0.6894,-1.6290) -- cycle;
\path[cbpmesh=19.6] (0.3393,-1.4860) -- (0.6766,-1.4673) -- (0.6894,-1.6290) -- (0.3458,-1.6481) -- cycle;
\path[cbpmesh=21.6] (-0.3505,-1.7971) -- (-0.0000,-1.8036) -- (-0.0000,-1.9399) -- (-0.3537,-1.9334) -- cycle;
\path[cbpmesh=22.2] (-0.0000,-1.8036) -- (0.3505,-1.7971) -- (0.3537,-1.9334) -- (-0.0000,-1.9399) -- cycle;
\path[cbpmesh=16.2] (-0.3334,-1.3852) -- (-0.0000,-1.3913) -- (-0.0000,-1.4923) -- (-0.3393,-1.4860) -- cycle;
\path[cbpmesh=16.6] (-0.0000,-1.3913) -- (0.3334,-1.3852) -- (0.3393,-1.4860) -- (-0.0000,-1.4923) -- cycle;
\path[cbpmesh=19.8] (-0.3458,-1.6481) -- (-0.0000,-1.6544) -- (-0.0000,-1.8036) -- (-0.3505,-1.7971) -- cycle;
\path[cbpmesh=20.4] (-0.0000,-1.6544) -- (0.3458,-1.6481) -- (0.3505,-1.7971) -- (-0.0000,-1.8036) -- cycle;
\path[cbpmesh=18.4] (-0.3393,-1.4860) -- (-0.0000,-1.4923) -- (-0.0000,-1.6544) -- (-0.3458,-1.6481) -- cycle;
\path[cbpmesh=19.0] (-0.0000,-1.4923) -- (0.3393,-1.4860) -- (0.3458,-1.6481) -- (-0.0000,-1.6544) -- cycle;
\draw[draw=cbpInk!45,line width=.40pt,dash pattern=on 1.3pt off 1.8pt] (1.3550,-0.2234) -- (1.3347,-0.0793) -- (1.3272,0.0575);
\draw[draw=cbpInk!45,line width=.40pt,dash pattern=on 1.3pt off 1.8pt] (4.5350,-0.0452) -- (4.5230,-0.1490) -- (4.4961,-0.2472) -- (4.4553,-0.3395) -- (4.4020,-0.4252) -- (4.3374,-0.5039) -- (4.2627,-0.5750) -- (4.1792,-0.6382) -- (4.0881,-0.6927) -- (3.9908,-0.7382) -- (3.8884,-0.7741) -- (3.7821,-0.7999) -- (3.6733,-0.8151) -- (3.5632,-0.8192) -- (3.4530,-0.8117) -- (3.3440,-0.7920) -- (3.2438,-0.7702) -- (3.1384,-0.7583) -- (3.0285,-0.7545) -- (2.9151,-0.7572) -- (2.7991,-0.7646) -- (2.6813,-0.7747) -- (2.5628,-0.7860) -- (2.4443,-0.7966) -- (2.3268,-0.8046) -- (2.2111,-0.8085) -- (2.0982,-0.8063) -- (1.9890,-0.7964) -- (1.8844,-0.7769) -- (1.7852,-0.7460) -- (1.6924,-0.7020) -- (1.6068,-0.6432) -- (1.5294,-0.5677) -- (1.4610,-0.4737) -- (1.4026,-0.3595) -- (1.3550,-0.2234);
\draw[draw=cbpInk!45,line width=.40pt] (1.3321,0.1865) -- (1.3487,0.3077) -- (1.3765,0.4208) -- (1.4149,0.5254) -- (1.4632,0.6215) -- (1.5209,0.7087) -- (1.5874,0.7867) -- (1.6621,0.8555) -- (1.7445,0.9146) -- (1.8339,0.9639) -- (1.9297,1.0031) -- (2.0314,1.0320) -- (2.1384,1.0503) -- (2.2501,1.0578) -- (2.3659,1.0542) -- (2.4852,1.0394) -- (2.6074,1.0130) -- (2.7320,0.9748) -- (2.8495,0.9351) -- (2.9596,0.9039) -- (3.0627,0.8799) -- (3.1595,0.8621) -- (3.2505,0.8494) -- (3.3363,0.8408) -- (3.4176,0.8351) -- (3.4949,0.8313) -- (3.5688,0.8283) -- (3.6399,0.8250) -- (3.7087,0.8203) -- (3.7759,0.8132) -- (3.8421,0.8025) -- (3.9078,0.7873) -- (3.9736,0.7663) -- (4.0401,0.7386) -- (4.1079,0.7029) -- (4.1776,0.6584) -- (4.2498,0.6038) -- (4.3250,0.5382) -- (4.4071,0.4146) -- (4.4679,0.2939) -- (4.5086,0.1767) -- (4.5306,0.0635);
}

\begin{document}
\title[Cavity collapse and vorticity gradient blow-up]{Cavity collapse and vorticity gradient blow-up in free-boundary Euler flow}
\author{Tao Luo}
\address{Department of Mathematics, City University of Hong Kong, Hong Kong, China}
\email{taoluo@cityu.edu.hk}
\author{Siqi Yang}
\address{School of Mathematics, Statistics and Mechanics, Beijing University of Technology, Beijing 100124, China}
\email{siqiyang@bjut.edu.cn}
\begin{abstract}
We construct smooth solutions of the three-dimensional incompressible
free-boundary Euler equations with surface tension whose vorticity
remains bounded while its full spatial gradient blows up in finite time.
The solutions are axisymmetric without swirl and occupy bounded
domains in which a tubular cavity collapses onto a circle at time
$T$, while the outer free surface converges smoothly.
The vorticity $\omega$ can be chosen arbitrarily small in $L^\infty$,
uniformly up to collapse, and satisfies
$\|\nabla\omega(t)\|_{L^\infty}\asymp
(T-t)^{-1/2}|\log(T-t)|^{1/4}$ as $t\uparrow T$.
A nonconstant material trace of the relative vorticity on the
shrinking surface determines direction-dependent limits of the
terminal vorticity, which therefore has no continuous extension
across the collapse circle. The kinetic energy measure converges
to the sum of a uniform line measure on this circle and a nonzero
bulk part carried by a bounded rotational field. The construction
couples vorticity transport to the boundary dynamics through
weighted estimates in the collapsing collar and a cancellation
in the normal acceleration.
\end{abstract}
\maketitle
\setcounter{tocdepth}{2}
\tableofcontents

\section{Introduction}
\subsection{The equations and the main result}
\label{subsec:boundary-main-results}
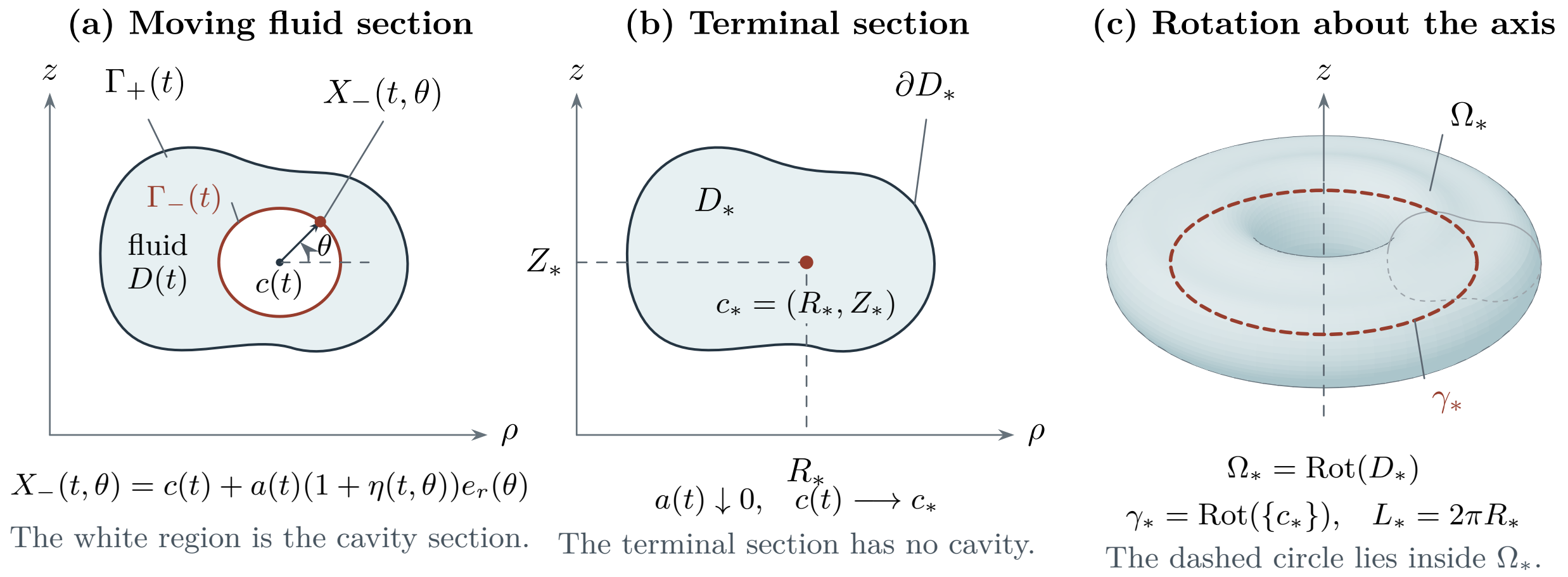
\begin{figure}[htbp]
\centering
\begin{tikzpicture}[cbpfigure]
\path[use as bounding box] (0,0) rectangle (16.4,6.0);
\node[font=\small\bfseries] at (2.65,5.72) {(a) Moving fluid section};
\node[font=\small\bfseries] at (8.15,5.72) {(b) Terminal section};
\node[font=\small\bfseries] at (13.65,5.72) {(c) Rotation about the axis}; 
\draw[cbpaxis] (.35,1.45) -- (.35,5.02) node[above] {$z$};
\draw[cbpaxis] (.35,1.45) -- (4.92,1.45) node[right] {$\rho$};
\begin{scope}[shift={(2.75,3.25)},xscale=.90,yscale=1.15]
\filldraw[cbpwall,fill=cbpFluid] \cbpOuterSection;
\end{scope}
\filldraw[cbpinner,fill=white] plot[domain=0:360,samples=101,variable=\q]
({2.75+.60*(1+.06*cos(2*\q))*cos(\q)},
 {3.25+.60*(1+.06*cos(2*\q))*sin(\q)}) -- cycle;
\node[font=\footnotesize,align=center] at (1.48,3.22) {fluid\\[-1pt]$D(t)$};
\fill[cbpInk] (2.75,3.25) circle (1.15pt);
\node[font=\footnotesize,inner sep=0pt] at (2.75,3.02) {$c(t)$};
\draw[cbpleader,dashed] (2.75,3.25) -- (3.68,3.25);
\draw[cbpradius] (2.75,3.25) -- (3.1743,3.6743);
\draw[cbpleader,->] (3.05,3.25) arc[start angle=0,end angle=45,radius=.30];
\node[font=\footnotesize] at (3.22,3.43) {$\theta$};
\fill[cbpCore] (3.1743,3.6743) circle (1.7pt);
\draw[cbpleader] (3.22,3.72) -- (3.83,4.72);
\node[anchor=south,fill=white,inner sep=1pt] at (3.83,4.78) {$X_-(t,\theta)$};
\node[anchor=south] at (1.35,4.80) {$\Gamma_+(t)$};
\draw[cbpleader] (1.38,4.74) -- (1.57,4.35);
\node[cbpCore,anchor=east,font=\footnotesize,inner sep=1pt] at (2.20,3.91) {$\Gamma_-(t)$};
\draw[cbpleader] (2.25,3.87) -- (2.32,3.67);
\node[font=\footnotesize] at (2.65,.91)
{$X_-(t,\theta)=c(t)+a(t)(1+\eta(t,\theta))e_r(\theta)$};
\node[font=\footnotesize,text=cbpInk!85] at (2.65,.36)
{The white region is the cavity section.}; 
\draw[cbpaxis] (5.85,1.45) -- (5.85,5.02) node[above] {$z$};
\draw[cbpaxis] (5.85,1.45) -- (10.42,1.45) node[right] {$\rho$};
\begin{scope}[shift={(8.25,3.25)},xscale=.90,yscale=1.15]
\filldraw[cbpwall,fill=cbpFluid] \cbpOuterSection;
\end{scope}
\node at (7.30,3.88) {$D_*$};
\node[anchor=south] at (9.50,4.80) {$\partial D_*$};
\draw[cbpleader] (9.50,4.74) -- (9.38,3.86);
\draw[cbpleader,dashed] (5.85,3.25) -- (8.25,3.25) -- (8.25,1.45);
\fill[cbpCore] (8.25,3.25) circle (2.0pt);
\node[anchor=north,font=\footnotesize,fill=cbpFluid,inner sep=1pt] at (8.25,3.01) {$c_*=(R_*,Z_*)$};
\node[anchor=east,fill=white,inner sep=1pt] at (5.76,3.25) {$Z_*$};
\node[below=3pt] at (8.25,1.45) {$R_*$};
\node[font=\footnotesize] at (8.15,.76) {$a(t)\downarrow0$,\quad $c(t)\longrightarrow c_*$};
\node[font=\footnotesize,text=cbpInk!85] at (8.15,.29) {The terminal section has no cavity.}; 
\begin{scope}[shift={(13.65,3.25)},scale=.50]
\cbpRotatedOuter
\end{scope}
\draw[cbpleader,dashed] (13.65,1.65) -- (13.65,4.10);
\draw[cbpaxis] (13.65,4.10) -- (13.65,5.02) node[above] {$z$};
\draw[cbpCore,line width=1.05pt,dash pattern=on 3pt off 2pt]
(13.65,3.25) ellipse[x radius=1.60cm,y radius=.7512cm];
\node[cbpCore,anchor=west] at (14.64,1.79) {$\gamma_*$};
\draw[cbpleader] (14.70,2.04) -- (14.60,2.6457);
\node[anchor=west] at (14.83,4.77) {$\Omega_*$};
\draw[cbpleader] (14.89,4.58) -- (14.76,4.00);
\node[font=\footnotesize,align=center] at (13.65,.84)
{$\Omega_*=\operatorname{Rot}(D_*)$\\[3pt]
$\gamma_*=\operatorname{Rot}(\{c_*\})$,\quad $L_*=2\pi R_*$};
\node[font=\footnotesize,text=cbpInk!85] at (13.65,.16)
{The dashed circle lies inside $\Omega_*$.};
\end{tikzpicture}
\caption{Moving and terminal geometry. In (a), the boundary labels
denote meridional traces and $\Omega_t=\operatorname{Rot}(D(t))$.
The inner meridian contracts to $c_*$, while the outer meridian
converges to $\partial D_*$. In (b), $D_*$ is simply connected and
stays away from the rotation axis. Revolution gives $\Omega_*$ and
the collapse circle $\gamma_*$ in (c). Vorticity reaches the inner
boundary but vanishes near the outer boundary.}
\label{fig:bv-domain}
\end{figure}

Consider an incompressible fluid occupying a bounded
domain $\Omega_t\subset\R^3$. Its boundary consists of an inner component
$\Gamma_-(t)$ enclosing the cavity $K_t$ and an outer component
$\Gamma_+(t)$ (see Figure~\ref{fig:bv-domain}). 
The velocity $v$ and pressure $p$ satisfy the free-boundary Euler
equations with surface tension:
\begin{equation}\label{eq:euler-toroidal-front}
\begin{cases}
\partial_t v+(v\cdot\nabla)v+\nabla p=0,&x\in\Omega_t,\\
\nabla\cdot v=0,&x\in\Omega_t,\\
p=\sigma\mathcal H,&x\in\partial\Omega_t,\\
\mathcal V_{\partial\Omega_t}=v\cdot n,&x\in\partial\Omega_t.
\end{cases}
\end{equation}
Here $\sigma>0$ is the surface tension coefficient, and $n$ is the
unit normal pointing out of the fluid. The mean curvature is defined
by $\mathcal H=\operatorname{div}_{\partial\Omega_t}n$, and
$\mathcal V_{\partial\Omega_t}$ denotes the normal velocity of the
boundary.

We use cylindrical coordinates $(\rho,\phi,z)$, defined by
\begin{equation*}
x_1=\rho\cos\phi,\qquad x_2=\rho\sin\phi,\qquad x_3=z,
\end{equation*}
and denote the associated orthonormal frame by $(e_\rho,e_\phi,e_z)$.
For an axisymmetric flow with zero swirl, the velocity and vorticity
have the form
\begin{equation*}
v=v^\rho e_\rho+v^z e_z,\qquad
\omega:=\nabla\times v=\omega^\phi e_\phi,
\end{equation*}
where $v^\rho$, $v^z$, and $\omega^\phi$ depend only on $(t,\rho,z)$.
We define the relative vorticity by $W:=\omega^\phi/\rho$, so that
$\omega=\rho W e_\phi$. The Euler equations imply that $W$ is
transported by the flow:
\begin{equation*} 
(\partial_t+v\cdot\nabla)W
=\partial_t W+v^\rho\partial_\rho W+v^z\partial_z W=0.
\end{equation*}
Throughout, we identify axisymmetric scalar fields with functions
on the meridional section $D(t)$. In cylindrical coordinates,
$\phi$ denotes the azimuthal angle about the rotation axis.
We use $\theta$ for the polar angle in a meridional section,
introduced below to parametrize the cavity boundary.

In \cite{LYIrrotational}, we constructed flows with a
single-valued velocity potential exhibiting cavity collapse with
prescribed terminal geometry and uniform concentration of kinetic
energy along the collapse circle. That construction applies both
with and without surface tension. Here we construct collapsing
flows whose relative vorticity has a nonconstant material trace
on the inner free surface. The main difficulty is to propagate
this profile on the degenerating domain and control the rotational
contributions to the pressure and boundary acceleration as the
spatial derivatives of the vorticity become unbounded. 

The terminal geometry is specified by a smooth bounded simply
connected meridional domain $D_*\Subset\{\rho>0\}$ and a point
$c_*=(R_*,Z_*)\in D_*$. Define
\begin{equation*}
\Omega_*=\operatorname{Rot}(D_*),\qquad
\gamma_*=\{\rho=R_*,\ z=Z_*\},\qquad L_*=2\pi R_*,
\end{equation*}
where $\operatorname{Rot}$ denotes revolution about the $z$ axis.
Thus $\Omega_*$ is the terminal fluid domain, $\gamma_*$ is the
collapse circle, and $L_*$ is its length.

For the inner boundary parametrization, we write
$c(t)=(R(t),Z(t))$ for the meridional centre and $a(t)>0$ for the
mean radius. We use the polar angle
$\theta\in\T:=\R/(2\pi\mathbb Z)$ about $c(t)$ and write
$e_r(\theta)=(\cos\theta,\sin\theta)$ for the radial unit vector.
Let $P_0$ and $P_1$ denote the $L^2(\T)$ orthogonal projections
onto the constant functions and
$\operatorname{span}\{\cos\theta,\sin\theta\}$, respectively.
Figure~\ref{fig:bv-domain} illustrates the moving and terminal domains.

We write $\mathrm{II}_-$ for the second fundamental form of
$\Gamma_-(t)$, $n_*$ for the outward unit normal to
$\partial\Omega_*$, and $\mathcal H^d$ for Hausdorff measure.
Matrix norms are Euclidean operator norms, and
$|\Sigma|=\mathcal H^2(\Sigma)$ denotes surface area.
Unless a domain is indicated, angular Sobolev norms are on $\T$.
For scalar oscillation we use
\begin{equation*}
\operatorname{osc}_A f=\sup_A f-\inf_A f.
\end{equation*}
For a terminal relative vorticity $W_*$, set
$\omega_*=\rho W_*e_\phi$. The conserved total energy is
\begin{equation*}
E_*=\frac12\int_{\Omega_t}|v|^2\,dx
+\sigma\bigl(|\Gamma_-(t)|+|\Gamma_+(t)|\bigr).
\end{equation*}
For a fixed reference length $\ell>0$, define $Y=\log(\ell/a)$.
Along the constructed solutions $a_t<0$, so $Y$ increases strictly.
We write $Y_0=\log(\ell/a(t_0))$, denote the inverse time
parametrization by $t(Y)$, and use the same notation for geometric
quantities after composition with $t(Y)$.

Unless otherwise indicated, limits in the following theorem are taken
as $t\uparrow T$, equivalently $Y\to\infty$.
\begin{theorem}
\label{thm:bv-collapse}
For every smooth bounded simply connected domain
$D_*\Subset\{\rho>0\}$, every $c_*\in D_*$ and every $E_c>0$,
there exists $\varepsilon_*>0$ such that, for every
$0<\varepsilon\leq\varepsilon_*$, there exist finite times $t_0<T$
and a smooth axisymmetric solution of \eqref{eq:euler-toroidal-front}
on $[t_0,T)$ with zero swirl and the following properties.

\begin{enumerate}[label=\textup{(\roman*)},leftmargin=*]
\item
The fluid boundary consists of the inner component $\Gamma_-(t)$
and the outer component $\Gamma_+(t)$. The inner meridian is
parametrized by
\begin{equation*}
X_-(t,\theta)=c(t)+a(t)(1+\eta(t,\theta))e_r(\theta),
\qquad P_0\eta=P_1\eta=0,
\end{equation*}
and the circulation around it vanishes for every $t\in[t_0,T)$.
The outer boundary converges to $\partial\Omega_*$ in every spatial
$C^m$ norm, measured in its normal graph parametrization,
while $c(t)\to c_*$ and $\eta(t,\cdot)\to0$ in $H^s(\T)$
for every $s\geq0$. The mean radius and the inner curvature satisfy
\begin{equation}\label{eq:bv-radius-law}
a(t)\sim
2^{3/4}\left(\frac{E_c}{\pi L_*}\right)^{1/4}
(T-t)^{1/2}|\log(T-t)|^{-1/4},
\qquad
a(t)\|\mathrm{II}_-(t)\|_{L^\infty(\Gamma_-(t))}
\longrightarrow1.
\end{equation}

\item
The vorticity support remains a fixed positive distance from
$\Gamma_+(t)$ throughout $[t_0,T)$, and
\begin{equation*}
\|\omega(t)\|_{L^\infty(\Omega_t)}\leq C\varepsilon,
\qquad t\in[t_0,T),
\end{equation*}
where $C$ is independent of $t$ and
$\varepsilon\in(0,\varepsilon_*]$. The solution can be chosen so that
the initial relative vorticity has the inner trace
\begin{equation*}
W(t_0,X_-(t_0,\theta))=\varepsilon\cos(2\theta).
\end{equation*}
There is a family of orientation-preserving diffeomorphisms
$\Psi_Y:\T\to\T$, with $\Psi_{Y_0}=\operatorname{id}_{\T}$,
such that
\begin{equation}\label{eq:bv-material-trace}
W(t(Y),X_-(t(Y),\Psi_Y(\theta)))
=\varepsilon\cos(2\theta).
\end{equation}
The maps $\Psi_Y$ converge smoothly to a diffeomorphism
$\Psi_\infty:\T\to\T$. In the polar parametrization,
\begin{equation*}
W(t,X_-(t,\cdot))\longrightarrow w_\infty
\quad\text{in }C^\infty(\T),
\qquad
w_\infty(\theta):=
\varepsilon\cos\bigl(2\Psi_\infty^{-1}(\theta)\bigr).
\end{equation*}
In particular,
\begin{equation}\label{eq:bv-positive-oscillation}
\operatorname{osc}_{\Gamma_-(t)}W=2\varepsilon,
\qquad
\operatorname{osc}_{\Gamma_-(t)}\omega^\phi
\longrightarrow2R_*\varepsilon>0.
\end{equation}

\item
There are constants $c_\varepsilon,C_\varepsilon>0$, independent
of $t$, such that
\begin{equation}\label{eq:bv-gradient-rate}
\frac{c_\varepsilon}{a(t)}
\leq\|\nabla\omega(t)\|_{L^\infty(\Omega_t)}
\leq\frac{C_\varepsilon}{a(t)},
\qquad t\in[t_0,T).
\end{equation}
Moreover, $\omega(t)$ converges smoothly on compact subsets of
$\Omega_*\setminus\gamma_*$ to
$\omega_*=\rho W_*e_\phi$, where $W_*$ is bounded.
In the meridional plane, the terminal relative vorticity has the
directional limits
\begin{equation*}
\lim_{r\downarrow0}W_*(c_*+r e_r(\theta))
=w_\infty(\theta),
\qquad\theta\in\T.
\end{equation*}
Since $w_\infty$ is nonconstant, $\omega_*$ has no continuous
extension across $\gamma_*$.

\item
There exists a bounded axisymmetric meridional vector field $V_*\not\equiv0$
on $\Omega_*$, smooth on
$\overline{\Omega_*}\setminus\gamma_*$, such that the kinetic
energy measures converge as finite Radon measures on $\R^3$:
\begin{equation}\label{eq:bv-energy-measure}
\frac12\mathbf1_{\Omega_t}|v(t)|^2\,dx
\ \rightharpoonup\
\frac12\mathbf1_{\Omega_*}|V_*|^2\,dx
+\frac{E_c}{L_*}\mathcal H^1\lfloor_{\gamma_*}.
\end{equation}
The field $V_*$ satisfies
\begin{equation*}
\operatorname{div}V_*=0,\qquad
\operatorname{curl}V_*=\omega_*
\quad\text{in }\mathcal D'(\Omega_*),
\qquad
V_*\cdot n_*=0\quad\text{on }\partial\Omega_*.
\end{equation*}
The conserved total energy is given by
\begin{equation*}
E_*=E_c+\frac12\int_{\Omega_*}|V_*|^2\,dx
+\sigma|\partial\Omega_*|.
\end{equation*}
\end{enumerate}
\end{theorem}

\begin{remark}
The nonvanishing oscillation of $\omega^\phi$ in
\eqref{eq:bv-positive-oscillation} persists on the shrinking
inner meridian, whose length is $2\pi a(t)(1+o(1))$.
This gives the lower bound in \eqref{eq:bv-gradient-rate}.
The matching upper bound follows from weighted estimates for
the transported relative vorticity in the collapsing collar
and controls the full spatial gradient of $\omega$.
Together with \eqref{eq:bv-radius-law}, these bounds yield
\begin{equation*}
\|\nabla\omega(t)\|_{L^\infty(\Omega_t)}
\asymp (T-t)^{-1/2}|\log(T-t)|^{1/4}
\qquad\text{as }t\uparrow T,
\end{equation*}
with comparison constants independent of $t$. In the collapsing coordinates, the transport estimates also give
uniform convergence of the relative vorticity profile up to the
inner boundary. Together with convergence of the material angular
flow, this identifies $w_\infty$ as the nonconstant directional
trace of $W_*$ at $c_*$, as illustrated in
Figure~\ref{fig:bv-material-trace}.
\end{remark}
\begin{figure}[htbp]
\centering
\begin{tikzpicture}[cbpfigure]
\path[use as bounding box] (0,0) rectangle (16.4,6.0);
\node[font=\small\bfseries] at (2.65,5.73) {(a) Initial material values};
\node[font=\small\bfseries] at (8.20,5.73) {(b) Transport and contraction};
\node[font=\small\bfseries] at (13.75,5.73) {(c) Terminal directional limits};
\node[font=\footnotesize] at (2.65,5.23) {$t=t_0$,\quad $\theta=\vartheta$};
\node[font=\footnotesize] at (8.20,5.23) {$t=t(Y)$,\quad $\theta=\Psi_Y(\vartheta)$};
\node[font=\footnotesize] at (13.75,5.23) {$t=T$,\quad $r\downarrow0$};
\begin{scope}[shift={(2.65,3.27)}]
\fill[cbpFluid] (-2.30,-1.63) rectangle (2.30,1.63);
\filldraw[cbpinner,fill=white] plot[domain=0:360,samples=101,variable=\q]
({.92*(1+.04*cos(2*\q))*cos(\q)},
 {.92*(1+.04*cos(2*\q))*sin(\q)}) -- cycle;
\fill[cbpInk] (0,0) circle (1.1pt);
\node[font=\footnotesize,below=4pt] at (0,0) {$c(t_0)$};
\foreach \q/\col in {0/cbpCore,90/cbpFlow,180/cbpCore,270/cbpFlow}{
\fill[\col] ({.92*(1+.04*cos(2*\q))*cos(\q)},
 {.92*(1+.04*cos(2*\q))*sin(\q)}) circle (2.4pt);
}
\node[cbpCore,anchor=west,inner sep=1pt] at (1.10,0) {$A:\ +\varepsilon$};
\node[cbpFlow,above=5pt] at (0,.89) {$B:\ -\varepsilon$};
\node[cbpCore,anchor=east,inner sep=1pt] at (-1.10,0) {$C:\ +\varepsilon$};
\node[cbpFlow,below=5pt] at (0,-.89) {$D:\ -\varepsilon$};
\end{scope}
\draw[cbpflow] (5.10,3.27) -- (5.72,3.27);
\begin{scope}[shift={(8.20,3.27)}]
\fill[cbpFluid] (-2.30,-1.63) rectangle (2.30,1.63);
\filldraw[cbpinner,fill=white] plot[domain=0:360,samples=101,variable=\q]
({.58*(1+.025*cos(2*\q))*cos(\q)},
 {.58*(1+.025*cos(2*\q))*sin(\q)}) -- cycle;
\fill[cbpInk] (0,0) circle (1.1pt);
\node[font=\footnotesize,inner sep=0pt] at (0,-.22) {$c(t)$};
\foreach \q/\col in {16/cbpCore,115/cbpFlow,196/cbpCore,277/cbpFlow}{
\fill[\col] ({.58*(1+.025*cos(2*\q))*cos(\q)},
 {.58*(1+.025*cos(2*\q))*sin(\q)}) circle (2.4pt);
}
\draw[cbpleader] (.61,.18) -- (1.03,.43);
\node[cbpCore,anchor=west,inner sep=1pt] at (1.10,.43) {$A:\ +\varepsilon$};
\draw[cbpleader] (-.25,.55) -- (-.40,.93);
\node[cbpFlow,above=3pt] at (-.40,.96) {$B:\ -\varepsilon$};
\draw[cbpleader] (-.61,-.18) -- (-1.03,-.43);
\node[cbpCore,anchor=east,inner sep=1pt] at (-1.10,-.43) {$C:\ +\varepsilon$};
\draw[cbpleader] (.08,-.60) -- (.12,-.94);
\node[cbpFlow,below=3pt] at (.12,-.97) {$D:\ -\varepsilon$};
\end{scope}
\draw[cbpflow] (10.65,3.27) -- (11.27,3.27);
\begin{scope}[shift={(13.75,3.27)}]
\fill[cbpFluid] (-2.30,-1.63) rectangle (2.30,1.63);
\draw[cbpCore!55,densely dashed] (0,0) -- (28:1.50);
\draw[cbpFlow!55,densely dashed] (0,0) -- (143:1.50);
\draw[-{Stealth[length=2.0mm,width=1.3mm]},cbpCore,line width=.9pt]
(28:1.35) -- (28:.22);
\draw[-{Stealth[length=2.0mm,width=1.3mm]},cbpFlow,line width=.9pt]
(143:1.35) -- (143:.22);
\draw[cbpInk,fill=white,line width=.7pt] (0,0) circle (.055);
\node[below=5pt] at (0,0) {$c_*$};
\node[cbpCore,align=center,anchor=south] at (1.20,.89)
{$W_*\to+\varepsilon$};
\node[cbpFlow,align=center,anchor=south] at (-1.20,1.10)
{$W_*\to-\varepsilon$};
\node[font=\footnotesize,cbpCore] at (0,-.91) {$\theta_+=\Psi_\infty(0)$};
\node[font=\footnotesize,cbpFlow] at (0,-1.33) {$\theta_-=\Psi_\infty(\pi/2)$};
\end{scope}
\node[font=\footnotesize] at (5.43,1.15)
{$W\bigl(t(Y),X_-(t(Y),\Psi_Y(\vartheta))\bigr)=\varepsilon\cos(2\vartheta)$};
\node[font=\footnotesize] at (5.43,.55)
{$A:\ \vartheta=0$,\quad $B:\ \vartheta=\pi/2$,\quad
 $C:\ \vartheta=\pi$,\quad $D:\ \vartheta=3\pi/2$};
\node[font=\footnotesize] at (13.75,1.15)
{$\displaystyle\lim_{r\downarrow0}W_*(c_*+r e_r(\theta))=w_\infty(\theta)$};
\node[font=\footnotesize] at (13.75,.55)
{$w_\infty=\varepsilon\cos(2\Psi_\infty^{-1})$};
\end{tikzpicture}
\caption{Material transport and the terminal trace of the relative
vorticity $W=\omega^\phi/\rho$. In (a) and (b), four material points
retain their values and cyclic order while the inner meridian
contracts. The material label $\vartheta$ and current polar angle
$\theta$ are distinct. In (c), the limiting angular map selects
two directions with different terminal values. The local views are
centred at $c(t_0)$, $c(t)$ and $c_*$, respectively.}
\label{fig:bv-material-trace}
\end{figure}
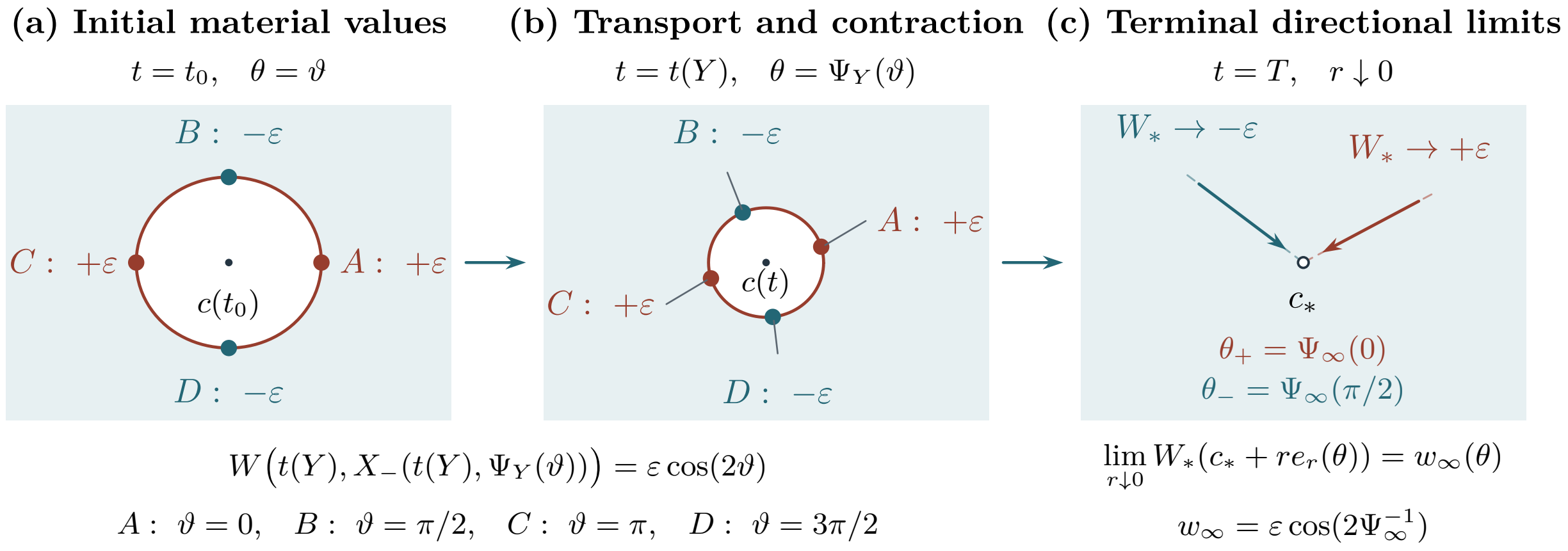
\begin{remark}
For each prescribed terminal domain $D_*$, the point $c_*\in D_*$
and the concentrated energy $E_c>0$ can be chosen independently. This freedom,
together with the leading collapse law \eqref{eq:bv-radius-law},
was established for irrotational flows in
\cite[Theorem~1.1]{LYIrrotational}.
In the present construction, the initial relative vorticity
is transported by the full velocity $v=\nabla\Phi+V$,
while $V$ is reconstructed from the evolving vorticity.
The rotational field also contributes to the pressure and
normal boundary acceleration, coupling the vorticity transport
to the boundary evolution.
The resulting flow determines the limiting angular map
$\Psi_\infty$ and hence the terminal profile $w_\infty$.
\end{remark}

\begin{remark}
The kinetic energy limit in \eqref{eq:bv-energy-measure} consists
of a uniform line measure on $\gamma_*$ and a nonzero absolutely
continuous part on $\Omega_*$. The line measure is generated
by the harmonic cavity flux, as in
\cite[Theorem~1.1]{LYIrrotational}.
The rotational field, after extension by zero, converges
strongly in $L^2(\R^3)$ to the zero extension of $V_*$,
with $\operatorname{curl}V_*=\omega_*\not\equiv0$.
The mixed energy terms vanish in the limit.
For the prescribed collapse circle, the concentrated energy $E_c$
determines the leading collapse law, while $V_*$ carries the
nonzero kinetic energy remaining in the terminal fluid domain.
\end{remark}

\subsection{Related work}
For Euler flows in fixed domains, rapid growth of vorticity
derivatives can occur while the vorticity remains bounded.
Kiselev and \v{S}ver\'ak \cite{KiselevSverak2014} constructed
two-dimensional flows in a disk whose vorticity gradient grows
doubly exponentially for all time. Do \cite{Do2019} obtained
doubly exponential growth of the relative vorticity gradient
for axisymmetric Euler flows without swirl, with the growth
occurring at the boundary of a fixed ball.

For flows with a moving free surface, Hu, Luo and Yao
\cite{HuLuoYao2024} proved a doubly exponential lower bound for
the vorticity gradient in two-dimensional periodic capillary
flows of finite depth. Tu \cite{Tu2026} obtained a corresponding
bound for two-dimensional gravity-capillary flows in a
horizontally unbounded domain of finite depth.
Both estimates hold throughout the lifespan of the corresponding
solution. These constructions begin with a flat free surface
and arbitrarily small velocity, and the growth occurs near
the fixed bottom. For a capillary drop whose entire boundary
is free, Hu \cite{HuHessian2026} proved an exponential lower
bound for the running supremum of the vorticity Hessian,
unless regularity breaks down earlier. The growth mechanism
in that construction acts in the fluid interior.
In the present solutions, the nonconstant vorticity trace lies
on the contracting free surface itself, and we obtain matching
upper and lower bounds for the full spatial gradient up to
the finite collapse time.

The terminal loss of vorticity continuity has also been studied
in fixed domains. Kiselev and Zlato\v{s}
\cite{KiselevZlatos2015} constructed finite time loss of
vorticity continuity in a planar domain whose boundary contains cusps. In our construction, both free surfaces remain
smooth before collapse, and the directional discontinuity
develops as the inner surface contracts onto the collapse circle.

Splash and splat singularities arise through contact between
distinct locally regular boundary branches.
Castro, C\'ordoba, Fefferman, Gancedo and G\'omez-Serrano
\cite{CastroCapillary2012} constructed such solutions with surface
tension in two dimensions. Coutand and Shkoller \cite{CSsplash}
constructed three-dimensional splash and splat solutions at zero
surface tension, allowing nonzero vorticity.
The cavity collapse in \cite{LYIrrotational} and in the present
paper involves contraction of the inner boundary onto a circle,
accompanied by unbounded curvature.

Cavity implosion has also been constructed in planar exterior
domains. Chen, Walsh and Wheeler \cite{ChenWalshWheeler2026}
obtained imploding hollow vortices whose velocity is irrotational
in the fluid and has nonzero circulation around the cavities.
The pressure is spatially constant on each cavity boundary,
with an isothermal gas interpretation.
The fluid considered here is bounded and three-dimensional,
with nonzero bulk vorticity, zero inner meridional circulation,
and the capillary condition on both free
surfaces.

The coexistence of bounded vorticity and degenerating boundary
geometry in these solutions is also relevant to continuation
criteria for free-boundary flows. Local theories allowing vorticity are
developed in \cite{CS,SZ3}. The original Beale--Kato--Majda
criterion \cite{BKM} concerns the Euler equations in the whole
space, while free-boundary continuation criteria also involve
the geometry and motion of the boundary
\cite{Ginsberg2021,LuoZhou2024,HLY}. In the present solutions,
the harmonic cavity flux carries the singular part of the
velocity while the vorticity remains bounded.

Kinetic energy concentration gives a further description of the
singular limit. Shvydkoy \cite{Shvydkoy2013} and Leslie and
Shvydkoy \cite{LS2018} studied energy measures and constraints
on the dimension of energy concentration. For the collapse considered here,
\eqref{eq:bv-energy-measure} identifies both the uniform line
measure on $\gamma_*$ and the nonzero bulk energy carried by
$V_*$.
\subsection{Proof strategy}
The construction couples the forward transport of the initial
vorticity to a boundary evolution with prescribed terminal
geometry. We organize the proof into four steps.

\medskip
\noindent\textbf{Step 1. Reconstruction and transport.}
We first decompose the velocity as $v=\nabla\Phi+V$.
The field $V$ is reconstructed from $\rho We_\phi$ by a Hodge
problem, with zero normal component on both boundaries and
zero inner meridional circulation. The harmonic part carries
the cavity flux. To follow the vorticity profile as the cavity
shrinks, we introduce a coordinate suggested by the leading
radial model. If $r$ denotes the distance from the centre and
$q>0$ the normalized cavity flux, this model gives
\begin{equation*}
r_t=-q/r,\qquad a_t=-q/a,\qquad
\frac{d}{dt}(r^2-a^2)=0.
\end{equation*}
Thus $S=r^2-a^2$ is constant along the model trajectories.
A Poisson extension adapts this coordinate to the noncircular
boundary and the moving centre. In the resulting coordinates,
the transport coefficients are integrable in $Y=\log(\ell/a)$.
This gives bounds for weighted normal and angular derivatives
of the profile. Returning to physical coordinates yields the
gradient upper bound in \eqref{eq:bv-gradient-rate}.

\medskip
\noindent\textbf{Step 2. Boundary estimates.}
The reconstructed rotational field enters the boundary evolution. Both its tangential velocity and the
rotational pressure contribute to the normal acceleration.
The pressure correction solves the Poisson problem with zero Dirichlet data.
Before applying the full Neumann inverse on both boundaries,
we project its normal derivative onto data with zero flux
on each boundary component, as in
\eqref{bv-st-pressure-neumann}.
In the principal height symbol, the mixed potential and
rotational terms cancel, and the remaining quadratic
rotational term is retained in the capillary energy.
Proposition~\ref{bv-scalars} derives the scalar flux and clock
estimates from the exact boundary system. Combining these
estimates with the spatial reduction in
Proposition~\ref{prop:bv-spatial-reduction} and the time derivative
estimates for the Hodge field gives uniform control of the two
boundary energies as $a\downarrow0$, as stated in
Lemma~\ref{bv-boundary-energy}.

\medskip
\noindent\textbf{Step 3. Coupling and construction.}
These estimates allow us to construct the flow with initial
vorticity and terminal geometric data.
We first stop the evolution at a positive terminal radius.
For a prescribed profile path, we solve the boundary equations
backwards from the terminal data and then transport the initial
relative vorticity forwards with the resulting velocity.
The Hodge and pressure estimates require one fewer profile
derivative than the boundary energy, while the Poisson extension
gains half a derivative in the transport coefficients.
This balance closes the successive profile and boundary
estimates at every fixed regularity level on a common logarithmic
time interval. Schauder's fixed point theorem then couples the two evolutions.
We first remove the regularity truncation at a fixed positive
terminal radius to obtain a smooth solution. The uniform
estimates then provide the compactness needed to let the
terminal radius tend to zero.
The normal equations, scalar flux equation and circulation
constraint recover the full Euler momentum equation
for the resulting evolution.

\medskip
\noindent\textbf{Step 4. Terminal limits.}
We finally determine the terminal behaviour of the constructed
solution. The boundary angular flow converges to a smooth
diffeomorphism, and the persistent oscillation on the shrinking
meridian yields the matching gradient lower bound.
Uniform convergence of the profile up to the inner boundary
identifies the transported limiting trace with the terminal
directional trace. The flux and clock equations give the
radius law \eqref{eq:bv-radius-law}.
For the energy limit, we separate the harmonic velocity into
its cavity flux component and a remainder. The flux energy
concentrates on $\gamma_*$, while the remainder tends to zero
in $L^2$. After extension by zero, the flux velocity converges
weakly to zero and the Hodge field converges strongly to the
zero extension of $V_*$ in $L^2(\R^3)$.
The mixed energy terms therefore vanish, giving
\eqref{eq:bv-energy-measure}.

\medskip
\noindent\textbf{Organization of the paper.}
Section~\ref{sec:imported-foundations} introduces the geometric
setting, recalls the harmonic estimates from \cite{LYIrrotational}
and establishes the Hodge, transport and rotational pressure
estimates. Section~\ref{sec:bv-energies} derives the scalar equations,
carries out the spatial reduction and proves the boundary energy
estimates. Section~\ref{sec:bv-construction} constructs the coupled
solution and proves the terminal limits in
Theorem~\ref{thm:bv-collapse}.

\section{Geometry and rotational estimates}
\label{sec:imported-foundations}

The geometric and harmonic estimates from \cite{LYIrrotational}
provide the setting for the Hodge reconstruction below. Its velocity
bounds yield transport and rotational pressure estimates that remain
uniform as the cavity shrinks.
\subsection{Geometry and velocity reconstruction}
\label{subsec:exact-centred-modulation}

We use the prescribed domain $D_*$, centre $c_*=(R_*,Z_*)$,
solid of revolution $\Omega_*$ and circumference $L_*$ from
Subsection~\ref{subsec:boundary-main-results}. Choose the reference
length $\ell>0$ such that
\begin{equation*}
4\ell<\min\left\{R_*,\tfrac12\operatorname{dist}(c_*,\partial D_*)\right\}.
\end{equation*}
Put $e_\theta(\theta)=(-\sin\theta,\cos\theta)$.
We parametrize the moving inner meridian and the outer graph by
\begin{equation*}
X_-=c+a(1+\eta)e_r,\qquad c=(R,Z),\qquad
P_0\eta=P_1\eta=0,\qquad
X_+=X_+^*+\zeta_+n_+^*.
\end{equation*}
Here $X_+^*$ parametrizes $\partial D_*$, $n_+^*$ is its outward
normal, and each graph is pulled back to its own parameter circle.
Choose $X_+^*$ by weighted arclength, so that
\begin{equation*}
(\rho\circ X_+^*)|\partial_\theta X_+^*|
=\frac1{2\pi}\int_{\partial D_*}\rho\,ds
=:\mu_*>0.
\end{equation*}
This choice is possible since $\overline{D_*}\subset\{\rho>0\}$.
The function $\zeta_+$ is the signed normal displacement of the
outer meridian from $\partial D_*$.
The meridians of $\Gamma_\pm$ are denoted by $C_\pm$, with $C_-$
oriented by increasing $\theta$.
Set 
\begin{equation*}
P_{\ge2}=I-P_0-P_1,
\end{equation*} 
the projection onto angular modes $|k|\geq2$.
To describe the cavity motion, recall that the logarithmic radius
$Y$ satisfies $a_Y=-a$. Let $q>0$ denote the normalized inward
cavity flux. Define the circumference $L$ and clock factor $\Theta$ by
\begin{equation*}
L=2\pi R,\qquad \Theta=\frac{qt_Y}{a^2}.
\end{equation*}
The physical flux $Q$ fixes the normalization of $q$:
\begin{equation*}
Q=\int_{\Gamma_-}v\cdot n_-\,dS=2\pi Lq.
\end{equation*}
Write $\Phi$ for the single-valued harmonic velocity potential.
The tangent Hodge field introduced below has zero normal trace, so
$\nabla\Phi$ carries the full cavity flux.

For $\iota\in\{-,+\}$, write $\rho_\iota=\rho\circ X_\iota$.
To express the inner normal and tangential components, set
\begin{equation*}
d=1+\eta,\qquad j=(d^2+\eta_\theta^2)^{1/2}.
\end{equation*}
The fluid outward normal and the positively parametrized inner tangent
are
\begin{equation*}
n_- =(-d e_r+\eta_\theta e_\theta)/j,\qquad
\tau_-=(\eta_\theta e_r+d e_\theta)/j.
\end{equation*}
The normalized inner normal velocity is
\begin{equation*}
b_-=a\partial_{n_-}\Phi/q.
\end{equation*}

To extract the mean and first harmonic of the kinematic condition,
write $\langle f\rangle=(2\pi)^{-1}\int_\T f$ and
$[f]_1=\pi^{-1}\int_\T f e_r\in\R^2$, taking vector means
componentwise. Here $B_-=v\cdot n_-$ is the full normal velocity.
Set
\begin{equation*}
C_\eta=e_r-\frac{\eta_\theta}{d}e_\theta,\qquad
\nu_\eta=\frac jd b_-,\qquad
(M_c)_{ij}=\pi^{-1}\int_\T(e_r)_i(C_\eta)_j\,d\theta.
\end{equation*}
The kinematic condition and the three centering constraints give
\begin{equation}\label{eq:general-centre-modulation} 
\Theta^{-1}=\langle\nu_\eta\rangle
-\langle C_\eta\rangle\cdot M_c^{-1}[\nu_\eta]_1,\qquad
t_Y=\frac{\Theta a^2}{q},\qquad
c_Y=-a\Theta M_c^{-1}[\nu_\eta]_1. 
\end{equation}
Indeed, the inner kinematic equation is
\begin{equation}\label{bv-st-kinematic-eta}
\eta_Y=d-\frac{c_Y\cdot e_r}{a}
+\frac{\eta_\theta}{d}\frac{c_Y\cdot e_\theta}{a}
-\frac{t_Y}{a}\frac jd B_-.
\end{equation}
Since $\langle d\rangle=1$ and $[d]_1=0$, taking these two
moments gives
\begin{equation*}
\frac{qt_Y}{a^2}\langle\nu_\eta\rangle
+\frac{c_Y}{a}\cdot\langle C_\eta\rangle=1,
\qquad
\frac{qt_Y}{a^2}[\nu_\eta]_1+M_c\frac{c_Y}{a}=0.
\end{equation*}
Eliminating $c_Y$ proves \eqref{eq:general-centre-modulation},
as in \cite[equation~(7.2)]{LYIrrotational}.
For small $\eta$, integration by parts and $P_1\eta=0$ give
\begin{equation*}
M_c=I+O(\eta),\qquad
\langle C_\eta\rangle
=-\langle\log(1+\eta)e_r\rangle=O(\eta^2),\qquad
-aa_t=q/\Theta.
\end{equation*}

To separate the cavity flux from the harmonic potential, let $u_c$
be harmonic in $\Omega$, with traces $1$ and $0$
on $\Gamma_-$ and $\Gamma_+$, and define
\begin{equation*}
\operatorname{Cap}=\int_\Omega|\nabla u_c|^2\,dx,
\qquad \Lambda=\frac{2\pi L}{\operatorname{Cap}}.
\end{equation*}
Subtract the condenser component from the single-valued harmonic
potential and write its boundary traces as
\begin{equation*} 
\Phi_0=\Phi-\frac{2\pi Lq}{\operatorname{Cap}}u_c,\qquad
(I-P_0)(\Phi_0|_{\Gamma_-})=q(\beta\cdot e_r+F),\qquad
\Phi_0|_{\Gamma_+}=g, 
\end{equation*}
where $\beta\in\R^2$, $P_0F=P_1F=0$, and $P_0g=0$ fixes the
common potential gauge. The zero inner flux determines the constant
part of the inner trace by the Schur equation, and the outer flux
then vanishes as well. We denote this extension with zero flux by
$\Phi_0[q(\beta\cdot e_r+F),g]$.

For a trace pair $\psi=(\psi_-,\psi_+)$, let
$\Phi[\psi]$ be its axisymmetric harmonic extension.
Set $j_\iota=|\partial_\theta X_\iota|$. In particular, $j_-=aj$.
The physical Dirichlet--Neumann (DN) map is
\begin{equation*}
(G\psi)_\iota
=\partial_{n_\iota}\Phi[\psi].
\end{equation*}
The weighted operator
$\operatorname{diag}(\rho_\iota j_\iota/R)G$ is self-adjoint and
nonnegative in angular $d\theta$ pairings
\cite[Lemma~2.1]{LYIrrotational}. Inner Sobolev norms use the scaled
angular chart and outer norms the fixed physical chart. We write
$\langle D_\theta\rangle=(1-\partial_\theta^2)^{1/2}$, with
subscript $\iota$ indicating the corresponding boundary chart.

The harmonic estimates of
\cite[Section~7.2, Step~1]{LYIrrotational} apply when the filled domains
stay a fixed positive distance from the axis, the centre stays separated
from the outer boundary, and the relative inner graph is small in
norms of fixed low order. In particular, the admissibility index
$b_{\rm ad}$ in \eqref{bv-bd-admissible-ball} includes the
$H^{s_0+2}$ graph norm required by
\cite[Lemmas~2.4--2.5 and Proposition~2.6]{LYIrrotational},
where $s_0\geq4$ is fixed. The base index $b$ is chosen above this order.

Under these geometric bounds, we use the Green and layer expansions,
uniform harmonic extension, the first two shape differentials,
estimates between separated boundaries and scalar Schur reduction in
\cite[Lemmas~2.2--2.5, Proposition~2.6, Lemma~2.9,
Proposition~2.10 and Lemma~2.11]{LYIrrotational}.
In \cite[equations~(2.35)--(2.36)]{LYIrrotational}, the inner input
has a fixed low Sobolev index as the outer output index increases.
The normalized inverse capacity satisfies $\Lambda=Y+O(1)$ on the
collapsing geometries considered here.

The spatial boundary reduction uses the calculus with a fixed cutoff
and finite derivative loss in
\cite[Lemma~3.13 and Lemmas~4.3--4.5]{LYIrrotational}.
The full affine part at the reference state is retained, as specified
in Step~3 of the proof of Lemma~\ref{bv-local}. Admissibility of the
frozen states and Taylor segments is verified on
\eqref{bv-bd-admissible-ball}, using
\eqref{bv-bd-clock-denominators}, with $q$ held fixed as an independent
coordinate. Its evolution is given by \eqref{bv-ex-q-equation}.

The normal change of variables and capillary energy calculations
follow \cite[Theorems~3.16 and~4.11, and Proposition~4.12]{LYIrrotational}.
Proposition~\ref{prop:bv-spatial-reduction} supplies the rotational
spatial reduction, and Lemma~\ref{bv-boundary-energy} verifies the
time-dependent coefficient hypotheses
on \eqref{bv-bd-weak-bootstrap} and estimates the additional rotational
terms, including the first mode. Lemmas~\ref{bv-pressure} and
\ref{bv-outer-time} give the pressure and time derivative estimates
with their dependence on the prescribed profile.

Two consequences of the harmonic analysis will be used repeatedly.
The uniform extension estimate in
\cite[Lemma~2.5 and equation~(4.55)]{LYIrrotational}, applied to the
full first harmonic $q\beta\cdot e_r$, gives
\begin{equation}\label{eq:concentration-zero-flux-bound}
\|\nabla\Phi_0\|_{L^2(\Omega)}^2
\le C\bigl(q^2|\beta|^2+q^2\|F\|_{H^1}^2+\|g\|_{H^1}^2\bigr).
\end{equation}
The constant is uniform in the inner radius under the stated
geometric bounds and zero component fluxes.

For the condenser, let $\Omega^{\rm fill}$ be the domain obtained
by filling the cavity, and extend $u_c$ by one through the cavity,
writing the extension as $\overline u_c$. The Green representation
\cite[proof of Lemma~5.11]{LYIrrotational} gives
\begin{equation}\label{eq:bv-condenser-representation}
\overline u_c(x)
=\int_{\Gamma_-}G_{\Omega^{\rm fill}}(x,y)\,d\mu_c(y),
\qquad d\mu_c=\partial_{n_-}u_c\,dS,
\qquad \mu_c\geq0,\quad \mu_c(\Gamma_-)=\operatorname{Cap}.
\end{equation}
Here $G_{\Omega^{\rm fill}}$ is the Dirichlet Green function for
$-\Delta$ on the filled domain. At each fixed derivative order,
the corresponding fixed outer chart bounds and source
separation give uniform Green kernel bounds on regions separated
from the cavity. In particular, the constants in the outer
capacity estimate \eqref{bv-bd-capacity-outer-derivative} are uniform
in $Y$.

The energy carried by the cavity flux is
\begin{equation}\label{eq:bv-flux-definitions}
E_f=\frac{Q^2}{2\operatorname{Cap}}=\pi L\Lambda q^2.
\end{equation}
The prescribed terminal value of $E_f$ is $E_c$.

It remains to reconstruct the rotational part of the velocity.
For each instantaneous geometry and profile, let $V$ be the
unique meridional field satisfying
\begin{equation}\label{eq:bv-hodge-definition}
\operatorname{div}V=0,\qquad
\operatorname{curl}V=\rho W e_\phi,\qquad
V\cdot n=0,\qquad
\oint_{C_-}V\cdot d\ell=0.
\end{equation}
Stokes' theorem determines the outer period. Reconstruct the velocity as
\begin{equation}\label{eq:bv-velocity-reconstruction}
v=\nabla\Phi+V,\qquad
\Phi=\frac{Q}{\operatorname{Cap}}u_c+
\Phi_0[q(\beta\cdot e_r+F),g],\qquad P_{\ge2}F=F.
\end{equation}
Since $V$ is tangent to the boundary, $B=v\cdot n=\partial_n\Phi$.

We describe vorticity transport on a fixed coordinate collar.
Choose a physical core collar with fixed radii and fix
$0<s_1<s_2<s_3$.
In the radial coordinate $S\in[0,s_3]$, the profiles are supported
in $0\leq S<s_2$. The physical collar stays a fixed positive
distance from the outer boundary. Let
$\chi$ be a smooth physical cutoff equal to one on this support and supported in
a slightly larger core collar.  Set
\begin{equation}\label{eq:bv-poisson-chart}
\begin{split}
\varpi&=(1+S/a^2)^{1/2},\\
H(Y,\varpi,\theta)
&=\chi(a\varpi)
\sum_{|k|\geq2}\eta_k(Y)\varpi^{-|k|}e^{ik\theta},\\
X(Y,S,\theta)
&=c(Y)+a(Y)(\varpi+H)e_r(\theta),\\
w(Y,S,\theta)&=W(t(Y),X(Y,S,\theta)).
\end{split}
\end{equation}
Here $\eta_k$ are the Fourier coefficients of $\eta$.
For $\eta$ sufficiently small in the fixed low norms, $X$ is a smooth
diffeomorphism onto the collar and agrees with $X_-$ at $S=0$.
We prescribe the initial profile
\begin{equation}\label{eq:bv-initial-profile}
w(Y_0,S,\theta)
=\varepsilon\chi_0(S)\cos(2\theta),
\end{equation}
where $\chi_0$ is smooth on $[0,\infty)$, equals one near zero,
and is supported in $[0,s_1)$.

To measure derivatives at the local annular scale in the chart
\eqref{eq:bv-poisson-chart}, write $\lambda=a^2+S$.
Regularity indices in the estimates below are integers unless
otherwise stated. For a function $f$ supported in the collar and a
nonnegative integer $m$, define
\begin{equation}\label{eq:bv-weighted-norms}
\begin{split}
N_m(f)&=\sum_{j+k\leq m}
\|\lambda^j\partial_S^j\partial_\theta^k f\|_{L^\infty},\\
M_m(f)^2&=\sum_{j+k\leq m}
\int_0^{s_3}\int_\T
|\lambda^j\partial_S^j\partial_\theta^k f|^2
\frac{d\theta\,dS}{\lambda},\\
Z_m(f)&=aM_m(f).
\end{split}
\end{equation}
We use $N_m$ for pointwise transport estimates and $Z_m$ for the
weighted energy estimates. For coefficients, a cutoff on a slightly
larger collar is understood.
The norms in \eqref{eq:bv-weighted-norms} are equivalent to rescaled
dyadic annular norms, with the composition terms in Lemma~\ref{bv-hodge}.
The normal weights remain positive at $S=0$.

Write $P_s(Y)$ for a polynomial in $1+Y$, chosen to satisfy
$P_s(Y)\geq1$ on the interval under consideration.
This polynomial may change from line to line in estimates.
Unless further dependence is specified, its coefficients depend on
the fixed geometry, $\sigma$, $E_c$, the fixed low norm bounds and $s$.
The reference length and cutoff functions are fixed with the geometric
data, and constants may depend on these choices. Constants in the
uniform collar and tail estimates are independent of the terminal
radius and the regularity truncation. In the local problem of
Subsection~\ref{subsec:bv-local-profile}, constants indexed by $N$
may depend on the positive minimum radius $a_N$.

A polynomial used in a bootstrap hypothesis has its degree and
coefficients fixed before $Y_0$ is chosen. This choice of $Y_0$ uses
only finitely many bounds at fixed low orders; the bounds at higher
orders are chosen afterwards. When $q\asymp(1+Y)^{-1/2}$, factors
$q^{-1}$ may be included in $P_s$, whereas inverse powers of $a$
remain separate. In $C_s(1+Y)^{-3/2}$, the constant $C_s$ is independent
of $Y$.

The Hodge estimates below are uniform as $a\downarrow0$ on the following
geometric class. Fix $0<\mu<1$. The filled domains have bounded
$C^{2,\mu}$ geometry, trace and Poincar\'e constants, and a fixed
positive distance from the axis. The centre stays a fixed positive
distance from the outer boundary. Smallness is imposed on the
$C^{2,\mu}$ norm of the relative inner graphs, while higher norms
enter the tame estimates.
Meridional integrals use weight $\rho$, differing from integrals
in three dimensions by $2\pi$.

To include the nonzero normal data needed later, let
$D_a=O\setminus\overline{K_a}$ be a meridional domain in this
geometric class, where $O$ is the simply connected filled section
and $K_a$ is the single cavity section. Prescribe a smooth azimuthal vorticity
$\omega^\phi e_\phi$ and normal data $f_-,f_+$ satisfying
\begin{equation}\label{bv-st-hodge-flux-data}
\int_{C_-}\rho f_-\,ds=\int_{C_+}\rho f_+\,ds=0.
\end{equation}
Fix a core radius $r_{\rm core}>0$ and take $r_j=2^ja$ up to
$r_{\rm core}$. Denote the rescaled annuli by $A_j$ and their
fixed enlargements by $A_j^*$, with
bounded overlap of the corresponding physical annuli and enlargements.
Let $A_+$ be a fixed outer collar and $A_+^*$ a fixed enlargement.
For a meridional field $v_H$, write
$v_{H,j}(z)=v_H(c+r_jz)$ and
$\omega_j^\phi(z)=\omega^\phi(c+r_jz)$.
To distinguish data bounds from graph regularity, set
\begin{equation*}
F_v=\|\omega^\phi\|_\infty
+\|f_-\|_{L^2(d\theta)}
+\|f_+\|_{H^{-1/2}},\qquad
G_m=\|\eta\|_{H^{m+1}},
\end{equation*}
and introduce the low regularity norms
\begin{equation*}
\ell_-=\|v_{H,0}\|_{C^{1,\mu}(A_0^*)}
+a\|\omega^\phi_0\|_{C^1(A_0^*)},\qquad
\ell_+=\|v_H\|_{C^{1,\mu}(A_+^*)}
+\|\omega^\phi\|_{C^1(A_+^*)}.
\end{equation*}
All source norms are taken before flattening. The source derivative
in $\ell_-$ accounts for composition with the flattening map.

For the estimates of the tangent field $V$ in
\eqref{eq:bv-hodge-definition}, assume $N_1(w)\leq C_0\varepsilon$
with fixed $C_0$ and keep its vorticity support a fixed positive
distance from the outer boundary. For the boundary energy estimates,
put $V_{\tau,-}=V\cdot\tau_-$ on $C_-$ and
$\widetilde v_-=\langle D_\theta\rangle^{-1}(aV_{\tau,-}/q)$.

\begin{lemma}
\label{bv-hodge}
For the domain $D_a$, source $\omega^\phi$ and normal data
$f_-,f_+$ satisfying \eqref{bv-st-hodge-flux-data}, there is a unique meridional field
$v_H=v_H^\rho e_\rho+v_H^ze_z$ solving
\begin{equation}\label{bv-st-hodge-system}
\begin{gathered}
\partial_\rho(\rho v_H^\rho)+\partial_z(\rho v_H^z)=0,
\qquad
\partial_zv_H^\rho-\partial_\rho v_H^z=\omega^\phi
\quad\text{in }D_a,\\
v_H\cdot n_\iota=f_\iota\quad\text{on }C_\iota,
\quad \iota\in\{-,+\},\qquad
\oint_{C_-}v_H\cdot d\ell=0.
\end{gathered}
\end{equation}
Here $n_\iota$ is the outward fluid normal. Stokes' theorem
determines the outer period.
For every integer $m\geq2$, the first annulus satisfies
\begin{equation}\label{bv-st-hodge-inner}
\|v_{H,0}\|_{H^m(A_0)}
\leq C_m\left[
F_v+a\|\omega^\phi_0\|_{H^{m-1}(A_0^*)}
+\|f_-\|_{H^{m-1/2}}+G_m \ell_-\right].
\end{equation}
For each $j\geq1$ with $c+r_jA_j^*\Subset D_a$,
\begin{equation}\label{bv-st-hodge-interior}
\|v_{H,j}\|_{H^m(A_j)}
\leq C_m\left[F_v
+r_j\|\omega^\phi_j\|_{H^{m-1}(A_j^*)}\right].
\end{equation}
On the fixed outer collar,
\begin{equation}\label{bv-st-hodge-outer}
\|v_H\|_{H^m(A_+)}
\leq C_m\left[
F_v+\|\omega^\phi\|_{H^{m-1}(A_+^*)}
+\|f_+\|_{H^{m-1/2}}
+\|\zeta_+\|_{H^{m+1}}\ell_+\right].
\end{equation}
The constants in
\eqref{bv-st-hodge-inner}--\eqref{bv-st-hodge-outer}
depend on $m$ and the bounds defining the geometric class, including
the trace and Poincar\'e bounds, uniformly over its cavity radii.

The tangent field $V$ satisfies, for every integer $s\geq2$,
\begin{equation}\label{bv-st-hodge-vorticity}
\begin{split}
\|V(c+a\cdot)\|_{H^s(A_0)}
+\|V|_{C_-}\|_{H^{s-1/2}}
&\leq P_s\left[\varepsilon+Z_{s-1}(w)
+\varepsilon\|\eta\|_{H^{s+1}}\right],\\
\|V\|_{H^s(A_+)}
&\leq P_s\varepsilon
\left(1+\|\zeta_+\|_{H^{s+1}}\right).
\end{split}
\end{equation}
Moreover,
\begin{equation}\label{bv-st-smoothed-vtrace}
\|\widetilde v_-\|_{H^{s+1/2}}
\leq P_s\frac aq\left[
\varepsilon+Z_{s-1}(w)
+\varepsilon\|\eta\|_{H^{s+1}}\right].
\end{equation}
The constants in these tangent field estimates also depend on
$C_0$.
\end{lemma}

\begin{proof}
The construction separates the filled-domain correction from the
correction caused by the hole. Uniform local energy bounds then
give the tame estimates and the boundary traces.

\medskip
\noindent\textbf{Step 1. Construction and uniform local energy bounds.}
Let $\operatorname{BS}=\nabla\times(-\Delta)^{-1}$ and
$u_b=\operatorname{BS}(\mathbf1_\Omega\omega^\phi e_\phi)$.
The vorticity is tangent to the axisymmetric interfaces, so its
extension by zero satisfies
\begin{equation*}
\operatorname{div}(\mathbf1_\Omega\omega^\phi e_\phi)=0,
\qquad
\operatorname{curl}u_b=\mathbf1_\Omega\omega^\phi e_\phi,
\qquad \operatorname{div}u_b=0
\end{equation*}
in distributions. The Biot--Savart kernel is $O(|x-y|^{-2})$. Its local
integrability and the common bounded support give
\begin{equation*}
\|u_b\|_\infty\leq C\|\omega^\phi\|_\infty,
\qquad u_b\in W^{1,p}_{\rm loc}(\R^3),\quad 1<p<\infty.
\end{equation*}
Thus $u_b$ is continuous across the interfaces and has zero flux
on each closed boundary component. Since $\operatorname{curl}u_b=0$
in the cavity, Stokes' theorem on interior meridians, followed by
passage to $C_-$, gives
\begin{equation*}
\oint_{C_-}u_b\cdot d\ell=0.
\end{equation*}

We correct the normal traces by adding harmonic gradients.
First solve the weighted Neumann problem for $\phi^{\rm out}$ on the
filled domain with datum $f_+-u_b\cdot n$. The compatibility
condition and the Neumann energy estimate give
\begin{equation*}
\int_{C_+}\rho(f_+-u_b\cdot n)\,ds=0,
\qquad \|\nabla\phi^{\rm out}\|_2\leq CF_v.
\end{equation*}
Interior elliptic estimates bound $\nabla\phi^{\rm out}$ in every
fixed $C^k$ norm on a fixed core neighbourhood by $C_kF_v$.
The remaining inner datum satisfies
\begin{equation*}
h=f_--u_b\cdot n-\partial_n\phi^{\rm out},\qquad
\int_{C_-}\rho h\,ds=0,\qquad
\|h\|_{L^2(d\theta)}\leq CF_v.
\end{equation*}
Let $\phi^h$ solve $\operatorname{div}(\rho\nabla\phi^h)=0$
in $D_a$, with $\partial_{n_-}\phi^h=h$ on $C_-$ and
$\partial_{n_+}\phi^h=0$ on $C_+$. Since $ds=aj\,d\theta$,
the zero weighted flux of $h$ and the trace and Poincar\'e inequalities on
the first rescaled collar yield
\begin{equation*}
\left|\int_{C_-}\rho h\psi\,ds\right|
\leq Ca\|h\|_{L^2(d\theta)}
\inf_{b\in\R}\|\psi|_{C_-}-b\|_{L^2(d\theta)}\leq Ca\|h\|_{L^2(d\theta)}
\|\nabla\psi\|_{L^2(D_a,\rho)}.
\end{equation*}
Testing with $\psi=\phi^h$ gives
\begin{equation}\label{bv-st-hole-energy}
\|\nabla\phi^h\|_{L^2(D_a,\rho)}
\leq Ca\|h\|_{L^2(d\theta)}\leq CaF_v.
\end{equation}
Hence the field
\begin{equation}\label{bv-st-hodge-decomposition}
v_H=u_b+\nabla\phi^{\rm out}+\nabla\phi^h
\end{equation}
solves \eqref{bv-st-hodge-system}. The difference of two such fields
has zero curl and period, so it is the gradient of a single-valued potential. Neumann
uniqueness gives uniqueness of $v_H$. Rescaling
\eqref{bv-st-hole-energy}, we obtain
\begin{equation}\label{bv-st-local-low-energy}
\|v_H(c+r_j\cdot)\|_{L^2(A_j)}
\leq CF_v+Cr_j^{-1}\|\nabla\phi^h\|_2
\leq C(1+a/r_j)F_v\leq CF_v,
\end{equation}
and the analogous bound on the fixed outer collar.

\medskip
\noindent\textbf{Step 2. Local tame estimates.}
For the higher local estimates, use
$\nabla^\perp\psi=(-\partial_z\psi,\partial_\rho\psi)$.
The zero weighted component fluxes give a single-valued stream
function satisfying $\rho v_H=\nabla^\perp\psi$. Set
$\psi_0(z)=a^{-1}\psi(c+az)$ and $\rho_0(z)=R+az_1$. Then
\begin{equation*}
\operatorname{div}_z(\rho_0^{-1}\nabla_z\psi_0)
=-a\omega^\phi_0,\qquad
\partial_\theta\!\left[
\psi_0\bigl((1+\eta(\theta))e_r(\theta)\bigr)
\right]=\rho_-(\theta)j(\theta)f_-(\theta).
\end{equation*}
The boundary primitive is periodic by \eqref{bv-st-hodge-flux-data},
and its additive constant does not affect $v_H$. On the rescaled
first annulus, choose a fixed radial cutoff $\chi_{\rm f}$ equal
to one near $r=1$ and zero near the outer edge, and set
\begin{equation*}
\mathcal F(r,\theta)
=(r+\chi_{\rm f}(r)\eta(\theta))e_r(\theta).
\end{equation*}
Here $r$ is the rescaled radial variable. The diagonal entries of
$D\mathcal F$ in polar frames are $1+\chi_{\rm f}'\eta$ and
$1+r^{-1}\chi_{\rm f}\eta$, so the stated small $C^{2,\mu}$
graph norm makes this an orientation-preserving diffeomorphism.
This radial flattening is the one used in
\cite[Lemma~2.4]{LYIrrotational} and requires only the stated
$C^{2,\mu}$ smallness. Put $\widehat\psi=\psi_0\circ\mathcal F$.
The pulled back equation is
\begin{equation}\label{bv-st-hodge-flattened} 
\operatorname{div}(A\nabla\widehat\psi)
=-a\det D\mathcal F\,(\omega_0^\phi\circ\mathcal F),\qquad
A=\frac{\det D\mathcal F}{\rho_0\circ\mathcal F}
(D\mathcal F)^{-1}(D\mathcal F)^{-T}. 
\end{equation}
The geometric assumptions give uniform ellipticity and uniform bounds
for the low coefficient norms. The higher coefficient norms are bounded by
$C_m(1+G_m)$. In the source composition, a high geometric derivative
multiplies $a\|\omega_0^\phi\|_{C^1}$, which is included in $\ell_-$.

To estimate \eqref{bv-st-hodge-flattened} for every integer $m\geq2$,
subtract a collar lift of the boundary primitive, which belongs to
$H^{m+1/2}$, and use nested cutoffs on the fixed flattened collar.
Tangential differentiation of the weak equation gives the energy
estimates for tangential derivatives. Uniform ellipticity then
recovers the second normal derivative from the equation, and
induction gives the remaining normal derivatives.

The scalar source is in $H^{m-1}$, so its contribution to the $H^{m+1}$ norm of the
stream function is $a\|\omega_0^\phi\|_{H^{m-1}}$.
In the commutators and the boundary lift, the highest geometric
derivative multiplies a lower derivative of the velocity or source,
bounded by $\ell_-$. The other products are controlled by tame
interpolation. This is the local elliptic argument used in
\cite[proofs of Lemma~2.5 and Proposition~2.6]{LYIrrotational},
here applied with the coefficient $\rho_0^{-1}$ and a scalar source.
The coefficient products and reconstruction of $v_H$ therefore give
\begin{align*}
\|v_{H,0}\|_{H^m(A_0)}
&\leq C_m\bigl[
\|v_{H,0}\|_{L^2(A_0^*)}
+a\|\omega^\phi_0\|_{H^{m-1}(A_0^*)}
+\|f_-\|_{H^{m-1/2}}+G_m\ell_-\bigr]\\
&\leq C_m\bigl[
F_v+a\|\omega^\phi_0\|_{H^{m-1}(A_0^*)}
+\|f_-\|_{H^{m-1/2}}+G_m\ell_-\bigr],
\end{align*}
by \eqref{bv-st-local-low-energy}. This is
\eqref{bv-st-hodge-inner}. The same interior estimate gives
\eqref{bv-st-hodge-interior}. On the fixed outer collar, the
coefficient term is $\|\zeta_+\|_{H^{m+1}}\ell_+$, which gives the
estimate \eqref{bv-st-hodge-outer}.

\medskip
\noindent\textbf{Step 3. The tangent field and its boundary trace.}
For $V$, the bound on $N_1(w)$ controls the $C^{0,\mu}$ norm of the
source on every rescaled annulus. The scaled first derivatives of
the chart and its inverse are uniformly bounded. Schauder estimates
on these annuli and on finitely many fixed interior and outer
boundary charts yield
\begin{equation*}
\|V(c+r_j\cdot)\|_{C^{1,\mu}(A_j)}\leq C\varepsilon,
\qquad \|V\|_{L^\infty(\Omega)}\leq C\varepsilon.
\end{equation*}

Composition with the Poisson chart contributes $G_mN_1(w)$ to the
source norm and $G_m$ times a lower velocity norm to the output.
The annular decay makes these terms square summable. Thus
\begin{equation*}
a\|\omega^\phi_0\|_{H^{s-1}}
\leq P_s\{Z_{s-1}+a\varepsilon G_s\},
\end{equation*}
and \eqref{bv-st-hodge-inner} gives the inner estimate in
\eqref{bv-st-hodge-vorticity}. On the outer collar the source
vanishes, so \eqref{bv-st-hodge-outer} gives the second bound in
\eqref{bv-st-hodge-vorticity}. Finally, the trace theorem gives
\begin{equation*}
\|\widetilde v_-\|_{H^{s+1/2}}
\leq C\frac aq\|V_{\tau,-}\|_{H^{s-1/2}}\leq P_s\frac aq\bigl[
\varepsilon+Z_{s-1}(w)
+\varepsilon\|\eta\|_{H^{s+1}}\bigr],
\end{equation*}
which is \eqref{bv-st-smoothed-vtrace}.
\end{proof}

The transport coefficients require the summed local estimate from
Lemma~\ref{bv-hodge}. Put $r=|(\rho,z)-c|$ and define $\|f\|_{m,a}$ as the
square root of the sum of the squared $H^m$ norms of
$f(c+r_j\cdot)$ on $A_j$, together with the squared norm on a fixed
outer block. This block consists of finitely many fixed interior
charts and the outer collar, covering the complement of the core,
with fixed enlargements and bounded overlap.

For the source terms, define $M_m^{\rm phys}(f)$ by the same formula using the enlarged annuli
$A_j^*$ and all the enlarged charts in this fixed outer block,
including $A_+^*$.
Multiplying \eqref{bv-st-hodge-inner} and
\eqref{bv-st-hodge-interior} by $r_j^{-1}$, including
\eqref{bv-st-hodge-outer}, and using $\sum_jr_j^{-2}\leq Ca^{-2}$ gives
\begin{equation}\label{bv-st-summed-hodge}
\left\|\frac{v_H}{r}\right\|_{m,a}
\leq C_m\left[
M_{m-1}^{\rm phys}(\omega^\phi)
+\frac{F_v+\|f_-\|_{H^{m-1/2}}+G_m\ell_-}{a}
+\|f_+\|_{H^{m-1/2}}
+\|\zeta_+\|_{H^{m+1}}\ell_+\right].
\end{equation}
Poisson pullback adds the source composition term
$G_mN_1(w)$ and the output composition term $a^{-1}G_m\ell_-$.

\subsection{Transport in the collapsing collar}
\label{subsec:bv-transport}

We pull the transport equation back by the collar chart
\eqref{eq:bv-poisson-chart}. The derivative $H_Y$ is taken at
fixed $\varpi$. The physical transport equation
$(\partial_t+v\cdot\nabla)W=0$ becomes
\begin{equation}\label{bv-st-transport-equation}
w_Y+\lambda b_1\,w_S
+b_2\,w_\theta=0,
\end{equation}
where the exact coefficients are
\begin{equation}\label{bv-st-transport-coefficients}
\begin{split}
b_2
&=\frac{(t_Yv(X)-c_Y)\cdot e_\theta}
{a(\varpi+H)},\\
b_1
&=\frac{2}{1+H_{\varpi}}
\left[
\frac{(t_Yv(X)-c_Y)\cdot e_r}{a\varpi}
+\frac1{\varpi^2}+\frac H{\varpi}
-(1-\varpi^{-2})H_{\varpi}
-\frac{H_Y}{\varpi}
-\frac{H_\theta}{\varpi}b_2\right].
\end{split}
\end{equation}
At fixed $S$, one has $\varpi_Y=\varpi-\varpi^{-1}$ and
\begin{equation*}
\begin{split}
X_S&=\frac{1+H_{\varpi}}{2a\varpi}e_r,\\
X_\theta&=aH_\theta e_r+a(\varpi+H)e_\theta,\\
X_Y&=c_Y+a\left[-\varpi^{-1}-H
 +(\varpi-\varpi^{-1})H_{\varpi}+H_Y\right]e_r.
\end{split}
\end{equation*}
Resolving $t_Yv-X_Y=\lambda b_1X_S+b_2X_\theta$ first in the
$e_\theta$ direction and then in the $e_r$ direction gives
\eqref{bv-st-transport-coefficients}. The kinematic condition
\eqref{bv-st-kinematic-eta} is equivalent to
$b_1|_{S=0}=0$.

To control the chart terms in \eqref{bv-st-transport-coefficients},
we use the gain of half a derivative from the Poisson extension.
For $P_0f=P_1f=0$ and every nonnegative integer $m$, Fourier
integration in $\xi=\log\varpi\in[0,\infty)$ gives
\begin{equation}\label{bv-st-poisson-half-gain}
\left\|\sum_{|k|\geq2}f_ke^{-|k|\xi}e^{ik\theta}
\right\|_{H^m((0,\infty)\times\T)}
\leq C_m\|f\|_{H^{m-1/2}(\mathbb T)}.
\end{equation}
This follows from
\begin{equation*}
\int_0^\infty e^{-2|k|\xi}\,d\xi=\frac1{2|k|},\qquad
\sum_{j+\ell\leq m}\int_0^\infty
\bigl|\partial_\xi^j(ik)^\ell e^{-|k|\xi}\bigr|^2\,d\xi
\leq C_m|k|^{2m-1}.
\end{equation*}
The physical cutoff preserves the bound. Differentiating the cutoff
in $Y$ at fixed $\varpi$ contributes to $H_Y$ the term
\begin{equation*}
-a\varpi\chi'(a\varpi)
\sum_{|k|\geq2}\eta_k\varpi^{-|k|}e^{ik\theta},
\end{equation*}
supported at fixed physical radius and gaining positive powers of $a$.
Thus the terms containing $H_Y$ use $\eta_Y\in H^{m-1/2}$,
the index of the normal datum in \eqref{bv-st-hodge-inner}. Spatial coefficients use
$\eta\in H^{m+1}$.

The chart \eqref{eq:bv-poisson-chart} also relates the weighted profile
norms to physical derivatives. With
$r=a(\varpi+H)$ and $D=\lambda\partial_S$, the inverse
chart identities are
\begin{equation}\label{bv-st-full-normal-derivatives}
r\partial_rW=\frac{2(\varpi+H)}
{\varpi(1+H_{\varpi})}Dw,
\qquad
\partial_\theta\big|_rW
=\partial_\theta w
-\frac{2H_\theta}{\varpi(1+H_{\varpi})}Dw.
\end{equation}
Iterating \eqref{bv-st-full-normal-derivatives} identifies the norms
in \eqref{eq:bv-weighted-norms} with the physical weighted
derivatives $r^{|\alpha|}\partial_x^\alpha W$.

For the general transport estimate, let
$L_w=\partial_Y+\lambda b_1\partial_S
+b_2\partial_\theta$,
with $b_1|_{S=0}=0$, and let
$L_w w=f$.  Assume the solution is supported
inside the collar, or insert a cutoff equal to one on its
support.  Define
\begin{equation*}
B_m=N_m(b_1)+N_m(b_2),\qquad
K_m=M_m(b_1)+M_m(b_2).
\end{equation*}
Here $D^+$ denotes the upper right Dini derivative with respect to $Y$.

\begin{lemma}
\label{bv-transport}
For every integer $m\geq1$,
\begin{equation}\label{bv-st-cm-transport}
D^+N_m(w)
\leq C_m\bigl(B_1N_m(w)
+B_mN_1(w)\bigr)+N_m(f).
\end{equation}
For $f=0$, the $L^\infty$ norm is preserved. For every integer $m\geq2$,
\begin{equation}\label{bv-st-sobolev-transport}
D^+Z_m(w)
\leq C_mB_1 Z_m(w)
+C_m aK_mN_1(w)+aM_m(f).
\end{equation}
All constants are independent of the minimum radius and of the
length of the logarithmic interval.
\end{lemma}

\begin{proof}
The weighted pointwise estimates follow along characteristics.
For the Sobolev estimates, the time-dependent weights give
nonpositive zeroth-order terms; a uniform tame estimate in
logarithmic coordinates controls the commutators.

Put $A=a^2/\lambda$ and
$J_{j,k}(g)=\lambda^j\partial_S^j\partial_\theta^k g$,
and abbreviate $J_{j,k}=J_{j,k}(w)$.
Direct differentiation gives
\begin{equation}\label{bv-st-damped-derivative}
L_w J_{j,k}+2jAJ_{j,k}
=jb_1J_{j,k}+J_{j,k}(f)-R_{j,k}.
\end{equation}
Here
\begin{equation*}
R_{j,k}
=\lambda^j[\partial_S^j\partial_\theta^k,
\lambda b_1\partial_S+b_2\partial_\theta]w.
\end{equation*}
The profile derivatives in $R_{j,k}$ have order at most
$j+k$. The highest coefficient derivatives multiply $Dw$
or $w_\theta$. Tame interpolation gives
\begin{equation*}
\sum_{j+k\leq m}\|R_{j,k}\|_\infty
 \leq C_m\{B_1N_m(w)+B_mN_1(w)\}.
\end{equation*}
To obtain constants independent of $a$, cover the collar by rectangles
of length comparable to $\lambda$ in $S$ and bounded overlap.
Rescale each rectangle to unit length, using half rectangles at $S=0$.
Since $2jA\geq0$ and $|b_1|\leq B_1$, integration of
\eqref{bv-st-damped-derivative} along characteristics proves
\eqref{bv-st-cm-transport}.

For $d\mu_a=dS\,d\theta/\lambda$,
\begin{equation*}
\partial_Yd\mu_a=2A\,d\mu_a,\qquad
\operatorname{div}_{\mu_a}(\lambda b_1,b_2)
=Db_1+\partial_\theta b_2.
\end{equation*}
Multiplying \eqref{bv-st-damped-derivative} by $2a^2J_{j,k}$ and
integrating gives
\begin{equation*}
\begin{split}
\frac d{dY}\left(a^2\int J_{j,k}^2\,d\mu_a\right)
&=a^2\int\{-2(1-A)-4jA+Db_1+\partial_\theta b_2\}
J_{j,k}^2\,d\mu_a\\
&\quad+2a^2\int J_{j,k}
\bigl(jb_1J_{j,k}+J_{j,k}(f)-R_{j,k}\bigr)
\,d\mu_a.
\end{split}
\end{equation*}
The inner flux vanishes because $b_1|_{S=0}=0$, and
$-2(1-A)-4jA\leq0$.
In the variable $z=\log\lambda$, one has $dS/\lambda=dz$,
$D=\partial_z$, and
\begin{equation*}
\lambda^j\partial_S^j=D(D-1)\cdots(D-j+1),\qquad j\geq1.
\end{equation*}
The norms $M_m$ are therefore uniformly equivalent to the ordinary
$H^m$ norms in $(z,\theta)$. Cover $[\log a^2,\infty)\times\T$ by cylinders
of fixed length and bounded overlap, with a half cylinder at the
lower endpoint. Applying the tame estimate on each cylinder and
summing gives
\begin{equation*}
\left(\sum_{j+k\leq m}
\|R_{j,k}\|_{L^2(d\mu_a)}^2\right)^{1/2}
\leq C_m\{
B_1M_m(w)+K_mN_1(w)\}.
\end{equation*}
Combining this commutator bound with the summed energy identity
and regularizing its square root
proves \eqref{bv-st-sobolev-transport}.
\end{proof}

\medskip
\noindent\textbf{Boundary variables and comparison norms.}

To apply Lemma~\ref{bv-transport} to the reconstructed velocity
\eqref{eq:bv-velocity-reconstruction}, we express $B_m$ and $K_m$
in terms of the profile and boundary variables. The boundary norms
introduced here will also be used in the energy estimates. The symbol $\operatorname{Op}$ denotes the Bony--Weyl quantization
of \cite[Lemma~3.13]{LYIrrotational}, with its fixed symmetric cutoff.
For the boundary estimates, set
\begin{equation*}
w_g=\frac{\rho_-}{R}j,\qquad
b_g=\frac{a}{qj}\left(d\partial_{n_-}\Phi
-\eta_\theta\partial_{\tau_-}\Phi\right),\qquad
B_q=\frac{\sigma a}{q^2}.
\end{equation*}
Here $w_g$ is the normalized surface weight and $B_q$ is the
capillary coefficient. The inner good variables are
\begin{equation}\label{eq:bv-good-variables}
\binom h\pi
=P_{\ge2}\operatorname{Op}\left(
\sqrt{w_g}
\begin{pmatrix}d/j&0\\-b_g&-1\end{pmatrix}
\right)P_{\ge2}\binom{a\eta}{aF}.
\end{equation}
With the coefficients frozen, this matrix operator has a bounded
inverse on the projected spaces. On the outer boundary, the normal
change of variables \cite[equation~(6.12)]{LYIrrotational} gives
the normal displacement and good potential trace
\begin{equation*}
\operatorname{Op}(n_+\cdot n_+^*)\zeta_+,\qquad
(I-P_0)\left[g-
\operatorname{Op}\bigl((\nabla\Phi)|_{\Gamma_+}\cdot n_+^*\bigr)\zeta_+\right].
\end{equation*}
The projection subtracts the ordinary angular mean. Applying the
change to surface half-densities with weight $\rho_+j_+$ defines the
outer good variables $(h_+,\pi_+)$
\cite[proof of Proposition~4.12, Step~1]{LYIrrotational}.
The coefficient factors in this change are quantized with the
same fixed cutoff.
On the weak neighbourhood, the weighted arclength choice gives
\begin{equation*}
\|\sqrt{\rho_+j_+}-\sqrt{\mu_*}\|_{H^1}
+\|n_+\cdot n_+^*-1\|_{H^1}
\leq C\|\zeta_+\|_{H^2}.
\end{equation*}
The fixed cutoff bounds the $L^2$ operator norm by the $H^1$
coefficient norm. Thus the diagonal operators are small $L^2$
perturbations of positive constant multipliers, including on the
stated gauge spaces. Their diagonal inverses and then the triangular
inverse are bounded. The reconstruction estimate in
\cite[Lemma~4.5]{LYIrrotational} gives the inverse bounds at every
Sobolev index with the same smallness threshold on the low-order norms.
We measure these variables in boundary norms $E_s$ and $O_s$ satisfying
\begin{equation}\label{eq:bv-energy-norms}
\begin{split}
E_s&\asymp
\|h\|_{H^s}+\sqrt{B_q}\|h\|_{H^{s+1}}
+\|\pi\|_{H^{s+1/2}},\\
O_s&\asymp
\|h_+\|_{H^{s+1}}+\|\pi_+\|_{H^{s+1/2}}
\asymp\|\zeta_+\|_{H^{s+1}}+\|g\|_{H^{s+1/2}} .
\end{split}
\end{equation}
Only these norm equivalences are used here.
Lemma~\ref{bv-boundary-energy} constructs positive quadratic forms
whose square roots can be chosen as $E_s$ and $O_s$. In particular,
\begin{equation}\label{eq:bv-interpolation-costs}
\|h\|_{H^{s+1/2}}\lesssim B_q^{-1/4}E_s,\qquad
\|h\|_{H^{s+1}}\lesssim B_q^{-1/2}E_s,\qquad
\|\eta\|_{H^s}+\|F\|_{H^s}\leq C_sa^{-1}E_s .
\end{equation}

At a fixed base index $b$, impose the weak bounds
\begin{equation}\label{bv-bd-weak-bootstrap}
E_b\leq P_ba^{7/4},\qquad O_b\leq P_ba^{3/4},\qquad
|\beta|\leq P_ba^{3/4}.
\end{equation}
Set $\delta_b=P_ba^{3/4}$, allowing the polynomial here to be a
fixed multiple of that in \eqref{bv-bd-weak-bootstrap}.
These two polynomials are fixed before $Y_0$ is chosen.
The boundary and transport estimates use the profile norms in the order
\begin{equation*}
Z_{s-1}^{\rm in}\ \longrightarrow\ (E_s,O_s)
\ \longrightarrow\ Z_s^{\rm out}.
\end{equation*}
Here $Z_j^{\rm in}$ and $Z_j^{\rm out}$ refer to the prescribed
and transported profiles, respectively.

\medskip
\noindent\textbf{General coefficient bounds.}

To isolate the leading cavity flux in the transport coefficients,
choose the weighted Green field $K_c$ in a fixed larger filled
domain, normalized to have inner flux $2\pi L$.
It is the gradient of a harmonic function in the fluid, with
\begin{equation*}
K_c(x)=-\frac{e_r}{r}+k_c(x),\qquad
\sum_{|\alpha|\leq m}r^{|\alpha|}
|\partial^\alpha k_c(x)|\leq P_m(Y).
\end{equation*}
This is the logarithmic Green expansion of
\cite[Lemma~2.2]{LYIrrotational}. For every meridional curve $C$
in the fluid enclosing $c$, with unit normal $n_C$ pointing into
the enclosed core, one has
\begin{equation*}
\int_C\rho K_c\cdot n_C\,ds=L.
\end{equation*}
The flux across its surface of revolution is therefore $2\pi L$.
Thus $v-qK_c$ has zero component flux, zero
inner period and curl $\rho We_\phi$. Denote its inner and outer
normal data by $f_-,f_+$. Use the norms $\ell_-,\ell_+$ with
$v_H=v-qK_c$. Set
\begin{equation*}
H_m=\|\eta_Y\|_{H^{m-1/2}},\qquad
\delta_c=|\Theta-1|+|c_Y|/a+P_b(Y)a.
\end{equation*}
In the radial coefficient, the leading field $-q e_r/r$ gives
\begin{equation*}
-\frac{\Theta}{\varpi(\varpi+H)}+\frac1{\varpi^2}
 =\frac{1-\Theta}{\varpi^2}
   +\frac{\Theta H}{\varpi^2(\varpi+H)}.
\end{equation*}
This cancellation leaves terms bounded by $|\Theta-1|$ and the graph
norms. The centre and regular Green terms contribute $|c_Y|/a$ and
$P_ma$, respectively. Applying
\eqref{bv-st-summed-hodge}, \eqref{bv-st-poisson-half-gain}, and
tame composition gives, for every integer $m\geq2$ and a prescribed
input profile $w^{\rm in}$, the bound
\begin{equation}\label{bv-st-coefficient-bound}
\begin{split}
K_m\leq P_m\bigg[&
\delta_c+G_m+H_m
+\frac{t_Y}{a}
\bigl(F_v+\|f_-\|_{H^{m-1/2}}\bigr)
+t_YM_{m-1}(w^{\rm in})\\
&+t_YG_mN_1(w^{\rm in})
+\frac{t_Y}{a}G_m\ell_-
+t_Y\bigl(\|f_+\|_{H^{m-1/2}}
+\|\zeta_+\|_{H^{m+1}}\ell_+\bigr)\bigg].
\end{split}
\end{equation}
Here the Green remainder is included in the first line, the vorticity
source has order $m-1$, and the inner composition term is
$(t_Y/a)G_m\ell_-$.

The geometric and boundary terms in \eqref{bv-st-coefficient-bound}
are controlled by the comparison norms. The triangular structure of
\eqref{eq:bv-good-variables} and
the bounded inverse of its diagonal operators give
\begin{equation*}
\|F\|_{H^{m+1/2}}
 \leq P_m\left(a^{-1}\|\pi\|_{H^{m+1/2}}
                   +\|\eta\|_{H^{m+1/2}}\right)
 \leq P_m(E_m/a+G_m).
\end{equation*}
The cutoff is fixed, so the coefficient norms in this estimate are
of fixed lower order. The Dirichlet and Neumann trace estimates now give
\begin{equation}\label{bv-st-boundary-to-coefficients}
\begin{split}
G_m&\leq P_m\frac{E_m}{a\sqrt{B_q}},\\
H_m+\frac{t_Y}{a}\|f_-\|_{H^{m-1/2}}
&\leq P_m\left[
\delta_c+|\beta|+\frac{E_m}{a}+G_m+aO_m\right],\\
\|f_+\|_{H^{m-1/2}}+\|\zeta_+\|_{H^{m+1}}\ell_+
&\leq P_m(q+O_m).
\end{split}
\end{equation}
For the middle line, write $K_c=\nabla G_c$ and put
$U=\Phi/q-G_c$. Its normal datum is $f_-/q$ on $C_-$, since
$V\cdot n_-=0$. Both component fluxes of $U$ vanish. Fix the
additive constant of $G_c$ at a fixed positive distance from $c$.
Let $z$ be the harmonic extension of $g/q-G_c|_{C_+}$ to the
\emph{filled} domain. The outer boundary stays away from the pole.
The maximum principle, $H^1(\T)\hookrightarrow L^\infty(\T)$ and
interior estimates on a fixed neighbourhood of the core give
\begin{equation*}
  \|\nabla z\|_{C^m(\mathrm{core})}
       \leq C_m A_g,\qquad A_g=1+q^{-1}\|g\|_{H^1}.
\end{equation*}
Here and below the core in a fixed-scale interior estimate is chosen
strictly inside the filled domain. In particular, each derivative of positive
order in the scaled variable of $z(c+a\,\cdot)$ has size
$O_m(a A_g)$. The number $A_g$ is bounded on the weak neighbourhood.

The function $u=U-z$ has zero outer trace and zero inner flux,
because $z$ extends harmonically through the cavity. Up to an inner
constant its Dirichlet datum is
\begin{equation*}
 e=\beta\cdot e_r+F-P_-^0\big[(G_c+z)\circ X_-\big],
 \qquad P_-^0=I-P_0 .
\end{equation*}
On the first scaled collar, every derivative of positive order,
up to any fixed order, of $(G_c+\log r)(c+a\,\cdot)$ is $O(P_m a)$ by the
Green field expansion above. The trace of $-\log r$ is
$-\log a-\log d$. Subtracting the mean and using the tame composition
estimate in the angular variable therefore gives
\begin{equation*}
 \|e\|_{H^{m+1/2}}
 \leq P_m\bigl(|\beta|+\|F\|_{H^{m+1/2}}+G_m+a A_g\bigr).
\end{equation*}
Indeed the high derivative of the graph in either smooth composition
has coefficient $O(P_m a A_g)$, while
\begin{equation*}
\|\log d\|_{H^{m+1/2}}\leq C_m G_m.
\end{equation*}

The inner constant selected by zero flux minimizes the Dirichlet
energy among extensions with outer trace zero and inner trace
$e+\mathrm{const}$. This follows by testing against functions with
constant inner trace and zero outer trace. A trace lift of $e$
supported on the first scaled collar is an admissible competitor.
Scale invariance of the meridional Dirichlet integral and the upper
and lower bounds on $\rho$ give
\begin{equation*}
\|\nabla u\|_{L^2(D,\rho)}\leq C\|e\|_{H^{1/2}},
\end{equation*} 
where $D$ is
the current meridional section. In particular this controls the
scaled gradient in $L^2$ on an enlarged first collar.

Write $\widehat u=u(c+a\,\cdot)$. On nested scaled collars the
Dirichlet estimate used in Step~2 of the proof of
Lemma~\ref{bv-hodge} gives
\begin{equation*}
 \|\partial_{\widehat n}\widehat u\|_{H^{m-1/2}}
 \leq C_m\bigl(\|e\|_{H^{m+1/2}}
       +\|\nabla u\|_{L^2(D,\rho)}+G_m\,\mathcal L\bigr).
\end{equation*}
Here $\mathcal L$ bounds the fixed low norms of the scaled gradient
and the trace lift. To see the dependence on the graph, subtract the
inner constant and a lift on the flattened collar, apply tangential
energy estimates, and recover the normal derivatives from the equation.
The highest derivative of the elliptic matrix is bounded by $C_mG_m$
and multiplies a lower derivative of the solution or lift. The fixed
lower estimates follow from the preceding energy bound and the fixed
lower trace and graph bounds. 

Thus $\mathcal L\leq P_b$ after the
base index is chosen. Restoring $z$ adds at most
$P_m a A_g(1+G_m)$ to the scaled normal trace. Since
$t_Y/a=\Theta a/q$, these estimates bound $(t_Y/a)f_-$ by the
right hand side of the middle line of
\eqref{bv-st-boundary-to-coefficients}, since
$a A_g\leq P_m(a+aO_m)$.
For $\eta_Y$, the singular part of $a(j/d)K_c\cdot n_-$ equals
$1/d$ exactly. In \eqref{bv-st-kinematic-eta} it combines with $d$
to give $d-\Theta/d$. The regular Green part costs $P_m a(1+G_m)$,
and the centre term costs $P_m(|c_Y|/a+G_m)$.
This proves the same bound for $H_m$.

For the last line of \eqref{bv-st-boundary-to-coefficients},
\eqref{eq:bv-condenser-representation} and
fixed outer collar regularity give
\begin{equation*}
\|(Q/\operatorname{Cap})\nabla u_c\|_{H^m(A_+)}
\leq P_mq(1+\|\zeta_+\|_{H^{m+1}}).
\end{equation*}
The part with zero flux is controlled by
\eqref{eq:concentration-zero-flux-bound}. In the fixed outer collar,
high derivatives fall on
$g$ and $\zeta_+$, and $V\cdot n_+=0$.

\medskip
\noindent\textbf{Bounds in the weak neighbourhood.}

Assume the weak boundary bounds \eqref{bv-bd-weak-bootstrap},
the fixed low-order profile bounds, and the clock and centre
estimates in \eqref{bv-ex-scalar-estimates}. For solutions of
\eqref{bv-ex-prescribed-profile-system}--\eqref{bv-ex-q-equation},
the latter estimates follow from Proposition~\ref{bv-scalars},
whose proof uses only spatial profile bounds. The lower harmonic
reconstruction and Lemma~\ref{bv-hodge} then give
\begin{equation*}
\ell_-\leq P_b\left(1+\varepsilon+q a^{-1/4}\right),\qquad
B_1\leq P_ba^{3/4},\qquad
\delta_c+|\beta|\leq P_ba^{3/4}.
\end{equation*}
For the $B_1$ bound, fix a coefficient index $m_\ell$ sufficient
for the scaled $C^{1,\mu}$ embedding and require $m_\ell+1\leq b$.
Then
\begin{equation*}
\|\eta\|_{H^{m_\ell+1}}\leq C E_b/a\leq P_ba^{3/4}.
\end{equation*}
At this order the preceding trace and scaled interior gradient
estimates give
\begin{equation*}
(a/q)\|f_-\|_{H^{m_\ell-1/2}}+H_{m_\ell}\leq P_ba^{3/4}
\end{equation*}
and bound the harmonic part of $\ell_-$ by $P_bq a^{-1/4}$.
The Hodge field adds $P_b\varepsilon$ to its lower
norm by Lemma~\ref{bv-hodge}. In \eqref{bv-st-coefficient-bound}
the source term is $t_YM_{m_\ell-1}=\Theta(a/q)Z_{m_\ell-1}$,
the remaining lower data in $(t_Y/a)F_v$ cost
\begin{equation*}
P_b(a^{3/4}+a(1+\varepsilon)/q),
\end{equation*}
and
\begin{equation*}
(t_Y/a)G_{m_\ell}\ell_-\leq P_ba^{3/4}(a/q+a^{3/4}).
\end{equation*}
The fixed outer terms have an additional factor $t_Y$.
The uniform Sobolev embedding on the scaled cylinders therefore
gives $B_1\leq C K_{m_\ell}\leq P_b a^{3/4}$ under the profile
bounds at the fixed low orders.
The first mode bound is part of \eqref{bv-bd-weak-bootstrap}.
At a higher level with
$E_m\leq P_ma^{7/4}$ and $O_m\leq P_ma^{3/4}$, one has
\begin{equation*}
G_m\leq P_m\frac q{\sqrt\sigma}a^{1/4},\qquad
t_YM_{m-1}(w^{\rm in})
=\Theta\frac aq Z_{m-1}(w^{\rm in}).
\end{equation*}
In particular, the inner output composition term in $aK_m$ satisfies
\begin{equation*}
t_YG_m\ell_-
\leq \frac{P_m}{\sqrt\sigma}
\left\{a^{9/4}(1+\varepsilon)+qa^2\right\}.
\end{equation*}
Thus this term is integrable for the fixed positive $\sigma$.
Substituting the preceding bounds in \eqref{bv-st-coefficient-bound} gives
\begin{equation*}
aK_m\leq P_ma^{5/4}
+P_m\frac{a^2}{q}Z_{m-1}(w^{\rm in}).
\end{equation*}
In particular, the outer terms before simplification are bounded
by $P_m(a^2+a^3/q)O_m$ and are absorbed into $P_ma^{5/4}$
under the stated higher outer bound.
The input term has order $m-1$.

The same coefficient bounds also control the time derivative of
the profile. Indeed, \eqref{bv-st-transport-equation} gives
\begin{equation*} 
M_k(w_Y)
\leq C_k\left[
B_0M_{k+1}(w)+K_kN_1(w)\right].
\end{equation*}
For $k=0,1$, use $K_k\leq K_2$. For every fixed integer $k\geq0$,
suppose in addition that $Z_{k+1}(w)$ and
$Z_{\max\{k-1,1\}}(w^{\rm in})$ have fixed bounds on the interval.
Then $M_{k+1}(w)=a^{-1}Z_{k+1}(w)$ and the preceding coefficient
estimates give $M_k(w_Y)\leq P_ka^{-1/4}$, where $P_k$ may also
depend on these profile bounds. The separate lower order estimate is
$\|w_Y\|_\infty\leq P_b\varepsilon a^{3/4}$.

\subsection{Pressure and outer time estimates}
\label{sec:bv-pressure}

In the velocity decomposition \eqref{eq:bv-velocity-reconstruction}, put
$v_0=\nabla\Phi$. The rotational pressure correction $R_p$ solves
\begin{equation}\label{bv-st-pressure-equation}
-\Delta R_p
=2\operatorname{tr}(\nabla v_0\nabla V)
+\operatorname{tr}((\nabla V)^2)=:f_p,
\qquad R_p|_{\partial\Omega}=0.
\end{equation}
To recover the pressure contribution to the boundary equations,
we remove the flux on each component before applying a Neumann inverse. Put
\begin{equation*}
d_-=\operatorname{Cap}^{-1}(\partial_n u_c)|_{\Gamma_-},\qquad
d_+=-\operatorname{Cap}^{-1}(\partial_n u_c)|_{\Gamma_+},
\qquad \int_{\Gamma_\iota}d_\iota\,dS=1,
\end{equation*}
and define
\begin{equation}\label{bv-st-capacitary-projection}
(\Pi_n g)_\iota
=g_\iota-d_\iota\int_{\Gamma_\iota}g_\iota\,dS.
\end{equation}
Let $\Phi_p$ solve
\begin{equation}\label{bv-st-pressure-neumann}
\Delta\Phi_p=0,\qquad
\partial_n\Phi_p=\Pi_n(\partial_nR_p)
\quad\text{on both components}.
\end{equation}
$\Phi_p$ is determined up to a common additive constant. For its
mean-free boundary traces, define in each angular chart
\begin{equation*}
P_\iota^0 f=f-\frac1{2\pi}\int_\T f(\theta)\,d\theta,
\qquad \iota\in\{-,+\}.
\end{equation*}
Thus $P_\iota^0$ subtracts the ordinary mean in the corresponding
boundary chart before the change to surface half-densities.

For the pressure estimates, assume $N_1(w)\leq C_0\varepsilon$, $0<\varepsilon\leq1$,
the fixed outer support gap, and the geometric and harmonic trace
bounds at the fixed low orders implied by \eqref{bv-bd-weak-bootstrap}.
Suppose $q\asymp(1+Y)^{-1/2}$.

\begin{lemma}
\label{bv-pressure}
There is an integer $s_{\rm p}\geq2$, fixed by the elliptic,
trace and product estimates in the proof and independent of the
cavity radius, such that the following bounds hold for every
integer $s\geq s_{\rm p}$. On a sufficiently late interval with
$B_q\leq1$, one has
\begin{equation}\label{bv-st-pressure-inner-result}
\frac{a^3}{q^2}\|P_-^0\Phi_p\|_{H^{s+1/2}}
\leq P_s\frac{a^2}{q}
\bigl(\varepsilon+Z_{s-1}(w)\bigr)
+P_s\frac{\varepsilon\sqrt a}{\sqrt\sigma}E_s
+P_sa^2O_s.
\end{equation}
At the outer boundary,
\begin{equation}\label{bv-st-pressure-outer-result}
t_Y\|P_+^0\Phi_p\|_{H^{s+1/2}}
\leq P_sa\varepsilon+e_{p,s}(Y)O_s,
\qquad
0\leq e_{p,s}\leq P_s\varepsilon(a+a^2/q)\in L^1_Y.
\end{equation}
The operators on surface half-densities with the fixed cutoff preserve
\eqref{bv-st-pressure-inner-result}--\eqref{bv-st-pressure-outer-result}.
The inner estimate uses the highest profile norm
$Z_{s-1}(w)$; the outer estimate uses only profile norms of fixed low order.
Both estimates depend only on spatial norms of the instantaneous profile.
\end{lemma}

\begin{proof}
We first estimate the Dirichlet pressure and its normal trace.
The physical flux projection and the Neumann lift then give
the inner and outer momentum sources.

\medskip
\noindent\textbf{Step 1. Low-order pressure bounds.}
For the harmonic remainder $v_0-qK_c$, the inner normal datum is
bounded by $P_bq\delta_b/a$. Its harmonic Neumann lift with zero
outer datum has gradient $L^2$ norm at most $P_bq\delta_b$,
by \eqref{bv-st-hole-energy}.
The fixed outer contribution, including the regular correction to
$qK_c$, is bounded by $P_b(q+\delta_b)\leq P_bq$
on a sufficiently late interval.
The harmonic reconstruction and Lemma~\ref{bv-hodge} therefore give
\begin{equation*}
\|v_0(c+r\cdot)\|_{C^{1,\mu}}\leq P_bq/r,\qquad
\|V(c+r\cdot)\|_{C^{1,\mu}}\leq P_b\varepsilon.
\end{equation*}
Substitution in \eqref{bv-st-pressure-equation} yields
\begin{equation*} 
|f_p|\leq P_b(q\varepsilon r^{-3}+\varepsilon^2r^{-2})
\quad(a\lesssim r\lesssim r_{\rm core}),
\end{equation*}
with bound $P_b(q\varepsilon+\varepsilon^2)$ outside the core.
The positive Newton potential of $|f_p|$ bounds $|R_p|$.
Its azimuthally integrated kernel is
$O(1+\log_+(C/|x_{\rm mer}-y_{\rm mer}|))$ by
\cite[Lemma~2.2]{LYIrrotational}.  Integrating at distances below
$a$ and then using $\int_a^{r_{\rm core}}r^{-2}dr=O(a^{-1})$ and
$\int_a^{r_{\rm core}}r^{-1}dr=O(1+Y)$ gives
\begin{equation}\label{bv-st-global-pressure-low}
\|R_p\|_\infty\leq P_b(q\varepsilon/a+\varepsilon^2)=:A_p
\leq P_bq\varepsilon/a.
\end{equation}
The last inequality holds on a late fixed interval.  Rescaled
inner and fixed outer boundary regularity yield
\begin{equation}\label{bv-st-pressure-low-traces}
\|R_p(c+a\cdot)\|_{C^{2,\mu}(A_0)}
+\|a\partial_nR_p\|_{C^{1,\mu}(C_-)}
+\|\partial_nR_p\|_{C^{1,\mu}(C_+)}\leq P_bA_p.
\end{equation}
The corresponding outer $C^J$ bound holds for each fixed
nonnegative integer $J$, with constants depending on the
corresponding higher geometric and source bounds.

\medskip
\noindent\textbf{Step 2. Higher estimates on the rescaled annulus.}
Normalize $\widehat v_0=(a/q)v_0(c+a\cdot)$,
$\widehat V=(a/q)V(c+a\cdot)$ and
$\widehat R_p=(a^2/q^2)R_p(c+a\cdot)$, and set
$\widehat A_p=a^2A_p/q^2$.
On the enlarged first annulus, the harmonic and Hodge estimates are
\begin{equation*}
\begin{split}
\|\widehat v_0\|_{H^s}
&\leq P_s\{1+\|F\|_{H^{s+1/2}}+G_s+(a/q)O_s\},\\
\|\widehat V\|_{H^s}
&\leq P_s(a/q)\{\varepsilon+Z_{s-1}+\varepsilon G_s\}.
\end{split}
\end{equation*}
The lower velocity bounds from Step~1, the product estimate and
\eqref{bv-st-pressure-low-traces} imply
\begin{equation*}
\begin{split}
\left\|\frac{a^4}{q^2}f_p(c+a\cdot)\right\|_{H^{s-1}}
&\leq P_s\left[
 \frac{a\varepsilon}{q}\|\widehat v_0\|_{H^s}
 +\left(1+\frac{a\varepsilon}{q}\right)
       \|\widehat V\|_{H^s}\right],\\
\|\partial_{\widehat n}\widehat R_p\|_{H^{s-1/2}}
&\leq P_s\left[
 \left\|\frac{a^4}{q^2}f_p(c+a\cdot)\right\|_{H^{s-1}}
 +\widehat A_p(1+G_s)\right].
\end{split}
\end{equation*}
The second line uses the local Dirichlet estimate from Step~2 of
the proof of Lemma~\ref{bv-hodge}, applied to the rescaled pressure
equation. Substituting the preceding bounds gives
\begin{equation}\label{bv-st-pressure-high-normal}
\|\partial_{\widehat n}\widehat R_p\|_{H^{s-1/2}}
\leq P_s\left[
\frac aq(\varepsilon+Z_{s-1})
+\frac{a\varepsilon}{q}(G_s+\|F\|_{H^{s+1/2}})
+\frac{a^2\varepsilon}{q^2}O_s\right].
\end{equation}

\medskip
\noindent\textbf{Step 3. Flux projection and the inner trace.}
Set $U_c=\operatorname{Cap}^{-1}u_c$. Its meridional gradient energy is
\begin{equation*}
\int_{D_a}\rho|\nabla U_c|^2\,d\rho\,dz
=\frac{1}{2\pi\operatorname{Cap}}.
\end{equation*}
The capacity bounds in \cite[Lemma~2.11]{LYIrrotational} make this
quantity polynomially bounded in $Y$. Apply the local elliptic
argument in Step~2 of the proof of Lemma~\ref{bv-hodge}, now to
$\operatorname{div}(\rho\nabla U_c)=0$ with constant Dirichlet data.
On the first rescaled collar, subtract
$\operatorname{Cap}^{-1}$ from $U_c(c+a\,\cdot)$; on the outer collar,
the boundary value is zero. Local Poincar\'e and Schauder estimates
control the low norms by the preceding energy bound. The higher
coefficient norm on the first collar uses
$G_s=\|\eta\|_{H^{s+1}}$, and the normal trace has order $s-1/2$.
The corresponding fixed outer argument uses only the outer graph
at that order. Since the normal traces of these two local functions
are $ad_-$ and $-d_+$, respectively, we obtain
\begin{equation}\label{bv-st-capacitary-profile-bounds}
\begin{split}
\|ad_-\|_{C^{1,\mu}}+\|d_+\|_{C^{1,\mu}}&\leq P_b,\\
\|ad_-\|_{H^{s-1/2}}&\leq P_s(1+G_s),\\
\|d_+\|_{H^{s-1/2}}&\leq P_s(1+\|\zeta_+\|_{H^{s+1}}).
\end{split}
\end{equation}
Division by the capacity contributes a polynomial factor in $Y$.
For $g=\Pi_n\partial_nR_p$, \eqref{bv-st-pressure-low-traces} and
the component areas give
\begin{equation*}
\left|\int_{\Gamma_\iota}\partial_nR_p\,dS\right|\leq P_bA_p,
\qquad
\|ag_-\|_{C^{1,\mu}}+\|g_+\|_{C^{1,\mu}}\leq P_bA_p.
\end{equation*}
The rescaled inner datum is
\begin{equation*}
\frac{a^3}{q^2}g_-
 =\partial_{\widehat n}\widehat R_p
  -\frac{a^3}{q^2}d_-
      \int_{\Gamma_-}\partial_nR_p\,dS.
\end{equation*}
By \eqref{bv-st-capacitary-profile-bounds}, the second term has
$H^{s-1/2}$ norm at most $P_s\widehat A_p(1+G_s)$.
The Neumann lifting in the proof of Lemma~\ref{bv-hodge}, with
zero curl, zero period and zero flux on each component, gives
\begin{equation}\label{bv-st-pressure-neumann-low}
\|\nabla\Phi_p\|_2
\leq P_b\{a\|g_-\|_{L^2(d\theta)}+\|g_+\|_{H^{-1/2}}\}
\leq P_bA_p.
\end{equation}
Set $\widehat\Phi_p(z)=a^2\Phi_p(c+az)/q^2$.
Scale invariance of meridional gradient energy and the Neumann
estimates in Lemma~\ref{bv-hodge} imply
\begin{equation}\label{bv-st-pressure-neumann-high}
\|P_-^0\widehat\Phi_p\|_{H^{s+1/2}}
\leq P_s\{\|\partial_{\widehat n}\widehat\Phi_p\|_{H^{s-1/2}}
+\widehat A_p(1+G_s)\}.
\end{equation}
The lower energy term is $P_b\widehat A_p$ by
\eqref{bv-st-pressure-neumann-low} and scale invariance of the
meridional gradient energy. The projection $P_-^0$ removes the
constant in the trace Poincar\'e inequality. Multiplying
\eqref{bv-st-pressure-neumann-high} by $a$ gives
\begin{equation*}
\frac{a^3}{q^2}\|P_-^0\Phi_p\|_{H^{s+1/2}}
\leq P_s\left[
 \frac{a^2}{q}(\varepsilon+Z_{s-1})
 +\frac{a^2\varepsilon}{q}
       (G_s+\|F\|_{H^{s+1/2}})
 +\frac{a^3\varepsilon}{q^2}O_s\right].
\end{equation*}
The triangular inverse of \eqref{eq:bv-good-variables} at index
$s+1/2$, together with \eqref{eq:bv-interpolation-costs}, gives
\begin{equation*}
G_s\leq P_s\frac{E_s}{a\sqrt{B_q}},\qquad
\|F\|_{H^{s+1/2}}\leq P_s\frac{E_s}{aB_q^{1/4}},\qquad
\frac{a\varepsilon}{q\sqrt{B_q}}
=\frac{\varepsilon\sqrt a}{\sqrt\sigma}.
\end{equation*}
Since $B_q\leq1$ and $a\varepsilon/q^2\leq1$ on a late interval,
this proves \eqref{bv-st-pressure-inner-result}.

\medskip
\noindent\textbf{Step 4. The outer trace and the change of variables.}
For the outer trace we work at the fixed physical scale. Set
$G_s^+=\|\zeta_+\|_{H^{s+1}}\leq P_sO_s$.
Hodge and harmonic regularity give
\begin{equation*}
\|V\|_{H^s(A_+)}\leq P_s\varepsilon(1+G_s^+),\qquad
\|v_0\|_{H^s(A_+)}\leq P_s(q+O_s+qG_s^+).
\end{equation*}
The Dirichlet and Neumann estimates and
\eqref{bv-st-global-pressure-low} then yield
\begin{equation*}
\|P_+^0\Phi_p\|_{H^{s+1/2}}
\leq P_s\{A_p(1+G_s^+)+\varepsilon O_s\}.
\end{equation*}
Multiplication by $t_Y=\Theta a^2/q$ proves
\eqref{bv-st-pressure-outer-result}.
The changes to surface half-densities use the order-zero
quantizations with the fixed cutoff. These operators are bounded
on $H^{s+1/2}$ by the fixed lower symbol seminorms. The higher
geometric terms in the preceding estimates satisfy
\begin{equation*}
\frac{a^3}{q^2}A_pG_s
\leq P_s\varepsilon\sqrt{a/\sigma}\,E_s,\qquad
t_YA_pG_s^+
\leq P_s\varepsilon a\,O_s
\leq P_s\varepsilon(a+a^2/q)O_s.
\end{equation*}
Thus \eqref{bv-st-pressure-inner-result} and
\eqref{bv-st-pressure-outer-result} persist after the change to surface half-densities.
\end{proof}

The effect of the removed fluxes is explicit. If
$\int_{\Gamma_-}g_-\,dS+\int_{\Gamma_+}g_+\,dS=0$, then
\begin{equation*}
(g-\Pi_ng)_\iota
=\frac{\int_{\Gamma_-}g_-\,dS}{\operatorname{Cap}}
 \partial_nu_c\big|_{\Gamma_\iota},\qquad
P_\iota^0(u_c|_{\Gamma_\iota})=0.
\end{equation*}
The harmonic potential of the removed datum therefore has zero projected trace.
This allows us to project the pressure and transport contributions
separately, imposing total flux compatibility on their sum.

The outer energy also requires a time derivative of the tangent
Hodge field on the fixed outer collar; see
\eqref{bv-bd-outer-symbol-time}. We take the Eulerian derivative
inside the fluid, before extension, and then estimate its pullback
by the outer chart.

For this time estimate, let the domains and profiles form a smooth
path in the geometric class of Lemma~\ref{bv-hodge}.
For a fixed nonnegative integer $J$, assume that the outer
charts have uniform bounds at the higher geometric orders required
by the elliptic estimates below.
Write $F_\iota$ for the geometric $Y$-normal velocity of each component,
relative to the outward fluid normal. Suppose
\begin{equation*}
|c_Y|\leq P_ba,
\quad\|F_-\|_{C^1(d\theta)}\leq P_ba,
\quad\|F_+\|_{C^{J+2}}\leq P_Ja^2.
\end{equation*}
Assume the same $C^{J+2}$ bound for the $Y$ velocity of the outer
charts. Assume the fixed support gap and
\begin{equation*}
N_1(w)\leq C_0\varepsilon,
\qquad\|w_Y\|_\infty
\leq P_b\varepsilon a^{3/4}.
\end{equation*}

\begin{lemma}
\label{bv-outer-time}
For every fixed nonnegative integer $J$, under the hypotheses above,
\begin{equation*} 
\|\partial_Y V\|_{C^J(A_+)}
\leq P_J\varepsilon(a^{3/4}+a).
\end{equation*}
$P_J$ may depend on the higher outer geometry bounds required for this $J$,
but not on the minimum cavity radius. The estimate holds both for
the Eulerian derivative on the common
outer region and for the $Y$ derivative of the field pulled back by
the outer collar chart. The dependence on $w_Y$ is through
its $L^\infty$ norm.
\end{lemma}

For the solution paths, the inner normal velocity
$F_-=t_YB_-$ satisfies the assumed $C^1$ bound
because $(a/q)B_-$ is bounded in $C^1$. For the finitely many
orders $J$ used in the boundary energy argument, choose $b\geq J+4$.
The weak bound on $O_b$ then gives
\begin{equation*}
\|F_+\|_{C^{J+2}}
\leq P_J\frac{a^2}{q}(q+O_{J+4})\leq P_Ja^2
\end{equation*}
on the fixed late interval. The outer normal graph charts obey
the same velocity bound.

\begin{proof}
We differentiate the vorticity and the Hodge boundary conditions
in physical coordinates, including the circulation condition,
and reconstruct the resulting field on the outer collar.

\medskip
\noindent\textbf{Step 1. The Eulerian vorticity derivative.}
Let a dot denote differentiation at fixed physical position.
By \eqref{bv-st-kinematic-eta}, with $t_YB_-$ replaced by $F_-$,
\begin{equation*}
\|\eta_Y\|_{C^1}
\leq P_b\left(1+\frac{|c_Y|}{a}
                  +\frac{\|F_-\|_{C^1}}a\right)
\leq P_b.
\end{equation*}
At fixed $S$, one has $\varpi_Y=\varpi-\varpi^{-1}$ and
\begin{equation*}
X_Y=c_Y+a\left[-\varpi^{-1}+H_Y-H
 +(\varpi-\varpi^{-1})H_\varpi\right]e_r,
\end{equation*}
where $H_Y$ is taken at fixed $\varpi$.
The Poisson extension, whose modes start at $|k|=2$, and the
inverse chart bounds imply
\begin{equation*}
|X_Y-c_Y|\leq P_ba^2/r,\qquad
|\nabla W(X)|\leq C\varepsilon/r.
\end{equation*}
Hence
\begin{equation*}
\dot W(X)=w_Y-X_Y\cdot\nabla W(X),\qquad
|\dot W(X)|\leq P_b\varepsilon
\{a^{3/4}+|c_Y|/r+a^2/r^2\}.
\end{equation*}
Integrating against $r\,dr\,d\theta$ and using
$\omega=\rho We_\phi$ give
\begin{equation}\label{bv-st-eulerian-vorticity-derivative}
\|\dot\omega\|_\infty\leq P_b\varepsilon,\qquad
\|\dot\omega\|_1\leq P_b\varepsilon
\{a^{3/4}+|c_Y|+a^2(1+Y)\}.
\end{equation}
\medskip
\noindent\textbf{Step 2. Differentiated boundary data and circulation.}
We also need the boundary data for the Eulerian derivative $\dot V$
on $\Sigma=\partial\Omega$. Let $\xi_-$ be the $Y$ velocity of the
inner parametrization, and write $F=F_\iota$ on $\Gamma_\iota$.
Differentiating the tangency and circulation constraints gives
\begin{equation}\label{bv-st-differentiated-normal}
\dot V\cdot n=\operatorname{div}_\Sigma(FV_\tau)=:g,
\end{equation}
and
\begin{equation}\label{bv-st-period-derivative}
\kappa(Y):=\oint_{C_-}\dot V\cdot d\ell
=\oint_{C_-}(\xi_-\times\omega)\cdot d\ell,
\qquad |\kappa(Y)|\leq P_b\varepsilon a^2.
\end{equation}
For a normal parametrization,
\begin{equation*}
\begin{split}
D_Yn&=-\nabla_\Sigma F,\qquad
\partial_n(V\cdot n)=-\operatorname{div}_\Sigma V_\tau,\\
\dot V\cdot n
&=-F\partial_n(V\cdot n)+V_\tau\cdot\nabla_\Sigma F
 =\operatorname{div}_\Sigma(FV_\tau),
\end{split}
\end{equation*}
which proves \eqref{bv-st-differentiated-normal}.
For the circulation, differentiation of the line element gives
\begin{equation*}
0=\frac d{dY}\oint_{C_-}V\cdot d\ell
 =\oint_{C_-}\left[\dot V-\xi_-\times\omega
                 +\nabla(\xi_-\cdot V)\right]\cdot d\ell.
\end{equation*}
The gradient integrates to zero, and the tangential component of
$\xi_-$ contributes zero because $\omega$ is azimuthal. Hence
\begin{equation*}
|\kappa|
\leq C\|F_-\|_\infty\|\omega\|_\infty
        \operatorname{length}(C_-)
\leq P_b\varepsilon a^2.
\end{equation*}
Both identities are unchanged by tangential reparametrization.
Since $g_-=(\rho_-aj)^{-1}
\partial_\theta(\rho_-F_-V_{\tau,-})$,
the lower Hodge bounds give
\begin{equation*} 
\|g_-\|_{L^2(d\theta)}\leq P_b\varepsilon,\qquad
\|g_+\|_{C^J}\leq P_Ja^2\varepsilon,\qquad
\int_{\Gamma_\iota}g_\iota\,dS=0.
\end{equation*}

\medskip
\noindent\textbf{Step 3. Reconstruction on the outer collar.}
To account for the period in \eqref{bv-st-period-derivative},
project $(2\pi)^{-1}\nabla\theta$ orthogonally off the gradients
in $L^2(D_a,\rho)$. The resulting $H_a$ is harmonic, tangent and
has period one. Projection and fixed outer regularity give
\begin{equation*}
\|H_a\|_2\leq C\sqrt{1+Y},\qquad
\|H_a\|_{C^J(A_+)}\leq P_J\sqrt{1+Y}.
\end{equation*}
Apply the decomposition \eqref{bv-st-hodge-decomposition}
to $\dot V-\kappa H_a$, whose curl is $\dot\omega$, normal data are
$g$, and inner period is zero:
\begin{equation*}
\dot V=v_b+\nabla\phi^{\rm out}_1+\nabla\phi^h_1+\kappa H_a.
\end{equation*}
Here $v_b$ is Biot--Savart applied to the zero extension of the
interior derivative $\dot\omega$.  The outer support gap and
\eqref{bv-st-eulerian-vorticity-derivative} give
\begin{equation*}
\|v_b\|_\infty\leq P_b\varepsilon,\qquad
\|v_b\|_{C^J(A_+)}
\leq P_J\varepsilon\{a^{3/4}+|c_Y|+a^2(1+Y)\}.
\end{equation*}
The fixed support gap and the assumed higher outer geometry
bounds give $\|V\|_{C^{J+2}(A_+)}\leq P_J\varepsilon$.
Thus
\begin{equation*}
\|g_+\|_{C^{J+1}}\leq P_J\varepsilon a^2,
\qquad C^{J+1}\hookrightarrow C^{J,\mu},\quad 0<\mu<1.
\end{equation*}
The Neumann datum on the filled domain is $g_+-v_b\cdot n_+$.
The Neumann energy estimate and boundary Schauder estimates give
\begin{equation*}
\|\nabla\phi^{\rm out}_1\|_{C^J(A_+)}
\leq P_J\varepsilon\{a^{3/4}+|c_Y|+a^2(1+Y)\}.
\end{equation*}
The remaining inner datum for the hole correction is
\begin{equation*}
h_1=g_--v_b\cdot n_--\partial_{n_-}\phi^{\rm out}_1,
\qquad \|h_1\|_{L^2(d\theta)}\leq P_b\varepsilon.
\end{equation*}
Thus
\eqref{bv-st-hole-energy} and outer regularity yield
\begin{equation*}
\|\nabla\phi^h_1\|_{C^J(A_+)}\leq P_Ja\varepsilon.
\end{equation*}
Using $|c_Y|\leq P_ba$ and \eqref{bv-st-period-derivative}, we obtain
\begin{equation*}
\|\dot V\|_{C^J(A_+)}
\leq P_J\varepsilon\left[a^{3/4}+a+a^2(1+Y)
                              +a^2\sqrt{1+Y}\right]
\leq P_J\varepsilon(a^{3/4}+a).
\end{equation*}
Finally, pullback by the outer chart adds $DV$ times its $Y$ velocity,
bounded by $P_J\varepsilon a^2$.
\end{proof}

\section{Boundary energies and flux dynamics}
\label{sec:bv-energies}

For a prescribed vorticity profile, the exact boundary system
determines the scalar flux, clock and centre. After a spatial
reduction of the rotational terms, the estimates of
Section~\ref{sec:imported-foundations} give the two boundary
energies and their terminal bounds, used in the coupled construction.

\subsection{The boundary system and scalar equations}
\label{subsec:bv-boundary-system}

Set $\Sigma=\partial\Omega$.
We use the total normal velocity $B=v\cdot n$ and write
$p_q=q_Y/q$ for the logarithmic derivative of the flux.

To compute the evolution of $B$, choose a surface parametrization
with velocity $Bn$ and put $v_\tau=v-Bn$.
With the convention
$\mathrm{II}(X,Y)=\langle\nabla_Xn,Y\rangle$, differentiation of
$B=v\cdot n$ and $n_t=-\nabla_\Sigma B$ gives
\begin{equation}\label{bv-bd-total-normal}
B_t+2v_\tau\cdot\nabla_\Sigma B
=-\partial_np+\mathrm{II}(v_\tau,v_\tau).
\end{equation}
Indeed, in the normal chart,
\begin{align*}
\partial_t(v|_\Sigma)&=-v_\tau\cdot\nabla_\Sigma v-\nabla p,\\
B_t&=-v_\tau\cdot\nabla_\Sigma B
+\mathrm{II}(v_\tau,v_\tau)-\partial_np
+v\cdot(-\nabla_\Sigma B)\\
&=-2v_\tau\cdot\nabla_\Sigma B-\partial_np
+\mathrm{II}(v_\tau,v_\tau).
\end{align*}
A tangential reparametrization adds the same chart term for
$v=v_0+V$ and $v_0=\nabla\Phi$, since $V\cdot n=0$.

To isolate the rotational contribution, fix the geometry, $B$ and
$q$, and let $p_0$ be the pressure of the instantaneous potential
field $v_0$. The harmonic extension $\Phi[\cdot]$ on both boundaries gives
\begin{equation}\label{bv-bd-potential-pressure}
p_0=-\frac12|v_0|^2
+\Phi\left[\sigma\mathcal H+\frac12|v_0|^2\right].
\end{equation}
Indeed, $\Delta(|v_0|^2/2)=\operatorname{tr}((\nabla v_0)^2)$.
Thus $p=p_0+R_p$, with $R_p$ from
\eqref{bv-st-pressure-equation}.
For derivatives of the boundary variables, let $\Delta$ denote
the full evolution minus the instantaneous potential evolution
at the same geometry, $B$ and $q$. The normal accelerations then differ by
\begin{equation}\label{bv-bd-rotation-normal}
\Delta B_t
=-\partial_nR_p-2V_\tau\cdot\nabla_\Sigma B
+\mathrm{II}(2v_{0,\tau},V_\tau)
+\mathrm{II}(V_\tau,V_\tau).
\end{equation}
The total flux of this difference vanishes. Integration by parts,
using $\operatorname{div}v_0=\operatorname{div}V=0$,
$V\cdot n=0$ and $v_0\cdot n=B$, gives
\begin{equation*}
\int_\Omega f_p\,dx
=\int_{\partial\Omega}
\left\{2V_\tau\cdot\nabla_\Sigma B
-\mathrm{II}(2v_{0,\tau},V_\tau)
-\mathrm{II}(V_\tau,V_\tau)\right\}\,dS.
\end{equation*}
Since $-\int_{\partial\Omega}\partial_nR_p\,dS=\int_\Omega f_p\,dx$,
\eqref{bv-bd-rotation-normal} implies
$\int_{\partial\Omega}\Delta B_t\,dS=0$.
We apply the component projection \eqref{bv-st-capacitary-projection}
separately to the pressure and kinetic contributions. The removed
component fluxes enter the scalar equation \eqref{bv-ex-q-equation}.

For a prescribed profile path $w(Y)$, the kinematic and projected
normal equations in logarithmic time are
\begin{equation}\label{bv-ex-prescribed-profile-system}
\begin{cases}
\partial_YX_\iota\cdot n_\iota=t_YB_\iota,
 &\iota\in\{-,+\},\\
\Pi_n\!\left\{\partial_Y^\perp B
+t_Y\bigl(2v_\tau\cdot\nabla_\Sigma B+\partial_np
-\mathrm{II}(v_\tau,v_\tau)\bigr)\right\}=0,
\end{cases}
\end{equation}
together with the scalar flux equation
\begin{equation}\label{bv-ex-q-equation}
q_Y=\frac{t_Y}{2\pi L}
\left\{\int_\Omega D^2u_c[v,v]
-\int_{\partial\Omega}g_c B^2
-\sigma\int_{\partial\Omega}\mathcal H g_c\right\}
-q\frac{L_Y}{L}.
\end{equation}
Here $B=v\cdot n$ and $g_c=\partial_nu_c$.
The derivative $\partial_Y^\perp$ is taken in a boundary chart moving
with velocity $t_YBn$. The second row is
\eqref{bv-bd-total-normal} after removal of the component fluxes by
\eqref{bv-st-capacitary-projection}. Equation~\eqref{bv-ex-q-equation}
supplies the scalar flux mode.

We use the graph charts for the first row, with $a=\ell e^{-Y}$.
Their inner mean and first harmonics determine $t_Y$ and $c_Y$
through \eqref{eq:general-centre-modulation}.
The projected normal equation evolves all nonconstant potential
traces, including the first harmonic $q\beta\cdot e_r$.
At each state the velocity is reconstructed by
\eqref{eq:bv-hodge-definition} and
\eqref{eq:bv-velocity-reconstruction}, and the pressure is
\begin{equation*}
p=p_0+R_p,
\end{equation*} 
with $p_0$ given by
\eqref{bv-bd-potential-pressure} and $R_p$ by
\eqref{bv-st-pressure-equation}. The boundary system takes $w(Y)$
as an independent input. The transport equation is imposed in
Section~\ref{sec:bv-construction}.

First fix the lower derivative orders in the symbol and kernel estimates of
\cite[Lemmas~4.3--4.5]{LYIrrotational} used in
Subsection~\ref{subsec:bv-spatial-reduction}. Include also the fixed
coefficient and commutator orders in
Subsection~\ref{subsec:bv-boundary-energies}: the good-variable and Taylor
coefficient bounds, the inner time estimates
\eqref{bv-bd-inner-symbol-time}, and the outer time estimates
\eqref{bv-bd-outer-good-time} and \eqref{bv-bd-outer-symbol-time}.
Let $J_*$ bound this finite set of seminorm indices. The calculus with
fixed cutoff in \cite[Lemma~3.13 and Theorem~3.16]{LYIrrotational}
uses a number of coefficient derivatives independent of the energy
index $s$; its constants may depend on $s$. The higher spatial norms
in the tame estimates retain their stated dependence on the output index.

For these fixed orders, let $c_0$ cover the additional spatial
regularity in Lemmas~\ref{bv-hodge} and~\ref{bv-pressure}, including
the trace and embedding orders in the inner Hodge time estimate below.
Thus the fixed profile orders there are $J+c_0$, with $J\leq J_*$.
Choose the spatial regularity index $k_{\rm t}$ for $w_Y$ in
\eqref{bv-bd-time-assumptions} so that $k_{\rm t}\geq J_*+c_0$.
Then choose $b$ to cover these derivative orders and the fixed geometric
trace orders, with $b\geq s_{\rm p}$ and $b-1\geq k_{\rm t}+1$.
For each order $J$ used in Lemma~\ref{bv-outer-time}, also retain its
requirement $b\geq J+4$. Here $s_{\rm p}$ is the threshold in
Lemma~\ref{bv-pressure}, and all fixed spatial profile indices are
required to be at most $b-1$.
These choices are independent of the energy index $s\geq b$.

For the dyadic energy estimates, set $S_j=0$ for $j\leq0$ and, for
$j\geq1$, choose smooth Fourier cutoffs to frequencies of size
at most $2^j$, with $S_j\to I$. Put \begin{equation*}
\Delta_j=S_{j+1}-S_j.
\end{equation*}
Choose $C_{\rm cut}\geq1$ to bound these
cutoffs and the convex interpolants $S_j+t\Delta_j$, $0\leq t\leq1$,
on the boundary Sobolev spaces, uniformly in $j$.
Throughout this section, we work in the geometric class of Lemma~\ref{bv-hodge}, with
vorticity support separated from the outer boundary, and assume
$q>0$, $\Theta>0$, $E_f\asymp E_c$ and the weak bounds
\eqref{bv-bd-weak-bootstrap}. For $U=(\eta,F,\zeta_+,g)$, choose
the fixed multiple in the definition of $\delta_b$ so that the norm
equivalences give
\begin{equation*}
\|U\|_{H^{b_{\rm ad}}}+|\beta|\leq\delta_b.
\end{equation*}
The profile is supported in $0\leq S<s_2$ in the chart
\eqref{eq:bv-poisson-chart} and satisfies
\begin{equation}\label{bv-bd-profile-assumptions}
N_1(w)\leq C_0\varepsilon,\qquad
Z_m(w)\leq P_m,\qquad 0\leq m\leq s-1.
\end{equation}
Here $C_0$ is fixed. At boundary energy level $s\geq b$, these
spatial bounds are imposed through order $s-1$. The scalar estimates
below use only the bounds through order $b-1$.

To apply the spatial boundary estimates, we use the convex set
on which the static geometric and harmonic bounds are uniform:
\begin{equation}\label{bv-bd-admissible-ball}
\mathcal A_Y=\{U=(\eta,F,\zeta_+,g):
\|U\|_{H^{b_{\rm ad}}}<2C_{\rm cut}\delta_b(Y)\},
\end{equation}
with $b_{\rm ad}$ fixed before $b$ and $q,\beta$ held fixed as
independent parameters in their stated ranges. In particular,
$E_f\asymp E_c$, $E_f=\pi L\Lambda q^2$ and the static capacity
estimate give $q\asymp(1+Y)^{-1/2}$, while $|\beta|\leq\delta_b$.
In the normalized inner and physical outer charts,
the elliptic operators are uniformly elliptic,
and the boundary components remain separated. The condenser and Schur
estimates of \cite[Proposition~2.10 and Lemma~2.11]{LYIrrotational}
give $\operatorname{Cap}\asymp(1+Y)^{-1}$ and polynomial bounds for
the first two shape differentials of $\operatorname{Cap}^{-1}$.
The normalized flux identity gives
\begin{equation}\label{bv-bd-clock-denominators}
\begin{aligned}
M_c&=I+O(\delta_b), &\langle C_\eta\rangle&=O(\delta_b^2),\\
\langle\nu_\eta\rangle&=1+O(P_b(a+\delta_b)),
&|[\nu_\eta]_1|&\le P_b(a+\delta_b).
\end{aligned}
\end{equation}
These bounds and \eqref{eq:general-centre-modulation} imply that
$M_c$ is invertible and $\Theta^{-1}$ has a uniform positive lower
bound on a sufficiently late interval. They also give moment bounds
with fixed derivative loss for $\Theta$ and $c_Y/a$.
With $q$ held fixed in its stated range, the geometry-dependent
denominators in the harmonic reconstruction and moment formulas are
$d,j,\rho,\det M_c,\Theta^{-1},\operatorname{Cap}$.
The constant $C_{\rm cut}$ ensures that all frozen states and Taylor
segments stay in \eqref{bv-bd-admissible-ball}.
Thus \cite[Lemmas~4.3--4.5]{LYIrrotational} apply on $\mathcal A_Y$,
with $q$ an independent coordinate.

Separate the cavity flux from the remaining velocity by writing
\begin{equation*}
K=\frac{Q}{\operatorname{Cap}}=\Lambda q,\qquad
v_f=v-K\nabla u_c,\qquad B_f=v_f\cdot n.
\end{equation*}
The flux of $B_f$ on each boundary component is zero. The flux
energy $E_f$ in \eqref{eq:bv-flux-definitions} controls this
scalar evolution. Its terminal value and the clock asymptotics
will determine the collapse law after the coupled solution has
been constructed.

Consider a smooth solution of
\eqref{bv-ex-prescribed-profile-system}--\eqref{bv-ex-q-equation}
on $[Y_0,Y_N]$.

\begin{proposition}\label{bv-scalars}
The scalar equation \eqref{bv-ex-q-equation}, together with the
kinematic equations, gives
\begin{equation}\label{bv-ex-flux-energy-identity}
E_{f,Y}=t_Y\left\{
-\frac{K^2}{2}\int_{\partial\Omega}g_c^2B_f
+K\left[\int_\Omega D^2u_c[v_f,v_f]
-\int_{\partial\Omega}g_cB_f^2
-\sigma\int_{\partial\Omega}\mathcal H g_c\right]\right\}.
\end{equation}
On a sufficiently late interval, with
$\delta_b=P_ba^{3/4}$, one has
\begin{equation}\label{bv-ex-flux-energy-bound}
|E_{f,Y}|
\leq P_b\left(q^2\delta_b^2+q^2a\delta_b
+q\varepsilon a\delta_b+\varepsilon^2a^2+\sigma a\right).
\end{equation}
The flux and clock also satisfy
\begin{equation}\label{bv-ex-scalar-estimates}
\begin{split}
&q\asymp(1+Y)^{-1/2},\qquad
\frac{q_Y}{q}=-\frac1{2\Lambda}+O(P_b(a+\delta_b)),\qquad \Theta-1=O(P_b(a+\delta_b)),\\
& \Theta_Y=O(P_b\delta_b),\qquad
\frac{|c_Y|}{a}=O(P_b(a+\delta_b)),\qquad \Lambda_Y=1+O(P_b(a+\delta_b)).
\end{split}
\end{equation}
The profile enters these estimates through the spatial norms in
\eqref{bv-bd-profile-assumptions}.
\end{proposition}

\begin{proof}
The scalar flux equation gives the energy identity without imposing
the transport equation on the prescribed profile. Spatial bounds
then control its right-hand side and the derivatives of the
capacity and clock.

\medskip
\noindent\textbf{Step 1. The scalar identities.}
To derive the flux identity for a prescribed profile, introduce
the momentum residual
\begin{equation*}
R_e=v_t+(v\cdot\nabla)v+\nabla p.
\end{equation*}
The pressure equation and $\operatorname{div}v=0$ imply
$\operatorname{div}R_e=0$.
The constant boundary values of $u_c$ give
$(u_c)_t=-Bg_c$ on $\partial\Omega$, with $n$ the outward fluid
normal. Applying Reynolds' formula to $Q=\int_\Omega v\cdot\nabla u_c$
and $\operatorname{Cap}$ yields
\begin{equation}\label{bv-ex-Q-and-capacity} 
Q_t-\int_\Omega R_e\cdot\nabla u_c
=\int_\Omega D^2u_c[v,v]
-\int_{\partial\Omega}g_cB^2
-\sigma\int_{\partial\Omega}\mathcal H g_c,\qquad
(\operatorname{Cap})_t=-\int_{\partial\Omega}g_c^2B. 
\end{equation}
Indeed, the domain variation cancels the boundary term from
$-\operatorname{div}(v\otimes v)$, while
\begin{equation*} 
\int_\Omega v\cdot\nabla(u_c)_t
=\int_{\partial\Omega}B(u_c)_t
=-\int_{\partial\Omega}g_cB^2,\qquad
-\int_\Omega\nabla p\cdot\nabla u_c
=-\int_{\partial\Omega}p g_c
=-\sigma\int_{\partial\Omega}\mathcal H g_c. 
\end{equation*}
Since $u_c$ has traces $1$ and $0$,
\begin{equation*}
\int_\Omega R_e\cdot\nabla u_c
=\int_{\Gamma_-}R_e\cdot n.
\end{equation*}
For the prescribed-profile problem, the assumed scalar equation
\eqref{bv-ex-q-equation}, together with $Q=2\pi Lq$, makes this
pairing zero. Since $\operatorname{div}R_e=0$, the outer component
flux also vanishes. The projected normal equation is
$\Pi_n(R_e\cdot n)=0$; hence \eqref{bv-st-capacitary-projection}
gives $R_e\cdot n=0$ on both boundary components. In particular,
the full normal equation \eqref{bv-bd-total-normal} holds for
the prescribed-profile problem. Conversely, for an Euler solution
$R_e=0$, and \eqref{bv-ex-Q-and-capacity} yields
\eqref{bv-ex-q-equation}.

To extract the flux energy from \eqref{bv-ex-Q-and-capacity},
integration by parts gives
\begin{equation*}
\int_\Omega D^2u_c[\nabla u_c,\nabla u_c]
=\frac12\int_{\partial\Omega}g_c^3,\qquad
\int_\Omega D^2u_c[\nabla u_c,v_f]
=\frac12\int_{\partial\Omega}g_c^2B_f.
\end{equation*}
Substituting $v=K\nabla u_c+v_f$ and $B=Kg_c+B_f$ into
\eqref{bv-ex-Q-and-capacity}, we obtain
\begin{equation*} 
E_{f,t}
=KQ_t-\frac{K^2}{2}(\operatorname{Cap})_t=-\frac{K^2}{2}\int_{\partial\Omega}g_c^2B_f
+K\left[\int_\Omega D^2u_c[v_f,v_f]
-\int_{\partial\Omega}g_cB_f^2
-\sigma\int_{\partial\Omega}\mathcal H g_c\right], 
\end{equation*}
which proves \eqref{bv-ex-flux-energy-identity}.

\medskip
\noindent\textbf{Step 2. Low norms and the flux energy estimate.}
The flux energy relation \eqref{eq:bv-flux-definitions} and the
static capacity estimate give
\begin{equation*}
q^2=\frac{E_f}{\pi L\Lambda}\asymp(1+Y)^{-1},
\qquad \Lambda=Y+O(1).
\end{equation*}
Hence powers of $q^{-1}$ are absorbed in $P_b$.
By \cite[Lemmas~2.5 and~2.11]{LYIrrotational} and
Lemma~\ref{bv-hodge},
\begin{equation}\label{bv-ex-condenser-low-budget}
\begin{aligned}
g_c|_{\Gamma_-}
&=\frac1{a\Lambda}\bigl(1+O(P_b(\delta_b+a))\bigr),
& |B_f|_{\Gamma_-}&\leq P_bq\delta_b/a,\\
|D^2u_c|&\leq\frac{P_b}{\Lambda r^2},
&|v_f-V|&\leq P_b\left(\frac{q\delta_b}{r}+\delta_b\right)
\quad(a\lesssim r\leq r_{\rm core}).
\end{aligned}
\end{equation}
On a fixed outer region,
$g_c,D^2u_c=O(P_b/\Lambda)$ and
$B_f,v_f-V=O(P_b\delta_b)$.
The regular $\delta_b$ term includes the outer harmonic field,
and $\|V\|_\infty\leq C\varepsilon$.

For the clock estimates, set $A=M_c^{-1}[\nu_\eta]_1$.
On $\Gamma_-$,
\eqref{bv-ex-condenser-low-budget} gives
\begin{equation*}
\nu_\eta=\frac jd\left(a\Lambda g_c+\frac aqB_f\right)
=1+O(P_b(a+\delta_b)).
\end{equation*}
Together with the moment bounds
\eqref{bv-bd-clock-denominators}, the identities
\eqref{eq:general-centre-modulation} imply
\begin{equation*}
\Theta-1=O(P_b(a+\delta_b)),\qquad
\frac{c_Y}{a}=O(P_b(a+\delta_b)),\qquad
\frac{L_Y}{L}=O(P_ba(a+\delta_b)).
\end{equation*}
In particular, $t_Y\asymp a^2/q$.

Using $\int_{\Gamma_-}B_f=0$, $|\Gamma_-|\asymp a$ and
$t_YK^2\asymp a^2q\Lambda^2$, we find
\begin{align*} 
t_YK^2\left|\int_{\Gamma_-}g_c^2B_f\right|
&=t_YK^2\left|\int_{\Gamma_-}
\bigl(g_c^2-(a\Lambda)^{-2}\bigr)B_f\right|\leq P_bq^2\delta_b(\delta_b+a),\\
t_YK^2\left|\int_{\Gamma_+}g_c^2B_f\right|
&\leq P_bqa^2\delta_b
\leq P_bq^2\delta_b(\delta_b+a). 
\end{align*}
The last inequality holds for $Y$ sufficiently large.

Since $t_YK\asymp a^2\Lambda$ and the core measure is
comparable to $r\,dr\,d\theta$, the Hessian term satisfies
\begin{align*} 
t_YK\left|\int_\Omega D^2u_c[v_f,v_f]\right|
&\leq P_ba^2\int_a^{r_{\rm core}}
\left(\frac{q\delta_b}{r}+\delta_b+\varepsilon\right)^2
\frac{dr}{r}
+P_ba^2(\delta_b+\varepsilon)^2\\
&\leq P_b\bigl(q^2\delta_b^2+q\varepsilon a\delta_b
+\varepsilon^2a^2
+qa\delta_b^2+a^2\delta_b^2+\varepsilon\delta_ba^2\bigr)\\
&\leq P_b\bigl(q^2\delta_b^2+q\varepsilon a\delta_b
+\varepsilon^2a^2\bigr). 
\end{align*}
Here logarithms are absorbed in $P_b$, and $a/q\to0$.
The boundary terms obey
\begin{equation*} 
t_YK\left|\int_{\Gamma_-}g_cB_f^2\right|\leq P_bq^2\delta_b^2,\qquad
t_YK\left|\int_{\Gamma_+}g_cB_f^2\right|\leq P_ba^2\delta_b^2,\qquad
\sigma t_YK\left|\int_{\partial\Omega}\mathcal H g_c\right|\leq P_b\sigma a, 
\end{equation*}
where $|\mathcal H_-|\leq P_b/a$.
These estimates prove \eqref{bv-ex-flux-energy-bound}.

\medskip
\noindent\textbf{Step 3. The capacity and flux derivatives.}
For the capacity derivative in \eqref{bv-ex-Q-and-capacity},
the surface measure and flux normalization give
\begin{equation*} 
k=\frac{\rho_-}{R}d=d+\frac aR d^2\cos\theta,\qquad
\langle k\nu_\eta\rangle=1,\qquad
dS=2\pi\rho_-aj\,d\theta,\qquad B_-=(q/a)b_-. 
\end{equation*}
Thus \eqref{bv-ex-condenser-low-budget} yields
\begin{equation*}
\int_{\Gamma_-}g_c^2B
 =\frac{2\pi Lq}{a^2\Lambda^2}
   \left\{1+O(P_b(a+\delta_b))\right\}.
\end{equation*}
On $\Gamma_+$, $g_c=O(P_b/\Lambda)$ and
$B=O(P_b(q+\delta_b))$, so
\begin{equation*}
\left|\frac{\Lambda t_Y}{\operatorname{Cap}}
\int_{\Gamma_+}g_c^2B\right|\leq P_ba^2.
\end{equation*}
Differentiating $\Lambda=2\pi L/\operatorname{Cap}$ and using
\eqref{bv-ex-Q-and-capacity}, we obtain
\begin{equation*} 
\Lambda_Y
=\frac{\Lambda L_Y}{L}
 +\frac{\Lambda t_Y}{\operatorname{Cap}}
       \int_{\partial\Omega}g_c^2B=\Theta+\frac{\Lambda L_Y}{L}+O(P_b(a+\delta_b))
=1+O(P_b(a+\delta_b)). 
\end{equation*}
Since \eqref{bv-ex-flux-energy-bound} gives
$|E_{f,Y}|/E_f\leq P_ba$, the identity
$E_f=\pi L\Lambda q^2$ implies
\begin{equation*}
\frac{q_Y}{q}
=\frac{E_{f,Y}}{2E_f}
-\frac{\Lambda_Y}{2\Lambda}-\frac{L_Y}{2L}
=-\frac1{2\Lambda}+O(P_b(a+\delta_b)).
\end{equation*}

\medskip
\noindent\textbf{Step 4. The clock derivative.}
Equation~\eqref{bv-st-kinematic-eta} and the fixed lower
harmonic bounds give $\|\eta_Y\|_{C^1}\leq P_b$.
In the normalized normal equation, physical acceleration is
multiplied by $t_Ya/q\asymp a^3/q^2$.
The pressure formula \eqref{bv-bd-potential-pressure} and
Lemma~\ref{bv-pressure} bound the dynamic and capillary pressure derivatives by
$P_bq^2/a^3$ and $P_b\sigma/a^2$, respectively.
Combining these pressure bounds with the estimate for $q_Y/q$
from Step~3 and the instantaneous profile bounds, and including
the tangential transport and metric terms, gives
\begin{equation*}
\|\nu_{\eta,Y}\|_{C^0}
\leq P_b\left[1+\frac{a^3}{q^2}
\left(\frac{q^2}{a^3}+\frac{\sigma}{a^2}\right)\right]
\leq P_b.
\end{equation*}
Differentiating $\langle k\nu_\eta\rangle=1$, we have
\begin{equation*}
\langle\nu_{\eta,Y}\rangle
=-\langle(k-1)\nu_{\eta,Y}\rangle
-\langle k_Y\rangle-\langle k_Y(\nu_\eta-1)\rangle.
\end{equation*}
The first and last terms are $O(P_b(a+\delta_b))$.
Using the centering constraints $P_1\eta=P_1\eta_Y=0$, we have
\begin{equation*}
\langle d^2\cos\theta\rangle=\langle\eta^2\cos\theta\rangle,
\qquad
\langle d\eta_Y\cos\theta\rangle
=\langle\eta\eta_Y\cos\theta\rangle.
\end{equation*}
Consequently,
\begin{align*} 
\langle k_Y\rangle
 &=\frac aR\left[-\left(1+\frac{R_Y}{R}\right)
       \langle\eta^2\cos\theta\rangle
       +2\langle\eta\eta_Y\cos\theta\rangle\right]
 =O(P_ba\delta_b),\\
\langle\nu_{\eta,Y}\rangle&=O(P_b(a+\delta_b)). 
\end{align*}
Integration by parts and $\langle\eta_Ye_r\rangle=0$ give
\begin{equation*} 
\langle C_\eta\rangle=-\langle\log(1+\eta)e_r\rangle,\qquad
\partial_Y\langle C_\eta\rangle
 =-\left\langle\left((1+\eta)^{-1}-1\right)
                  \eta_Ye_r\right\rangle
 =O(P_b\delta_b). 
\end{equation*}
The coarse $C^1$ bound on $\eta_Y$ and the $C^0$ bound on
$\nu_{\eta,Y}$ also yield $A_Y=O(P_b)$.
Differentiating \eqref{eq:general-centre-modulation}, we conclude that
\begin{equation*} 
\Theta_Y
=-\Theta^2\left(\langle\nu_{\eta,Y}\rangle
-\partial_Y\langle C_\eta\rangle\cdot A
-\langle C_\eta\rangle\cdot A_Y\right)=O(P_b(a+\delta_b))=O(P_b\delta_b), 
\end{equation*}
for $Y$ sufficiently large, since
$A=O(P_b(a+\delta_b))$ and
$\langle C_\eta\rangle=O(\delta_b^2)$.
\end{proof}

\subsection{The projected Neumann inverse and rotational terms}
\label{subsec:bv-spatial-reduction}

We separate the spatial reduction from the time-dependent energy
calculation. The projected Neumann inverse controls the interaction
between the two boundaries; its shape derivatives also enter the
local comparison argument in Subsection~\ref{subsec:bv-local-profile}.
Throughout, component constants and finite-dimensional parameters are
treated separately, and potential traces are projected before the
half-density change.

For a smooth symbol $r(\theta,\xi)$ of order $m$, set
\begin{equation*}
|r|_{J,m}
=\max_{\substack{\alpha,\beta\geq0\\\alpha+\beta\leq J}}
\sup_{(\theta,\xi)\in\T\times\R}
\langle\xi\rangle^{\beta-m}
|\partial_\theta^\alpha\partial_\xi^\beta r(\theta,\xi)|,
\end{equation*}
where $\alpha,\beta$ are integers.

Choose the cutoffs so that the bands $\Delta_j$, $j\geq1$, have
annular support $|k|\asymp2^j$ and fixed finite overlap, while
$\Delta_0=S_1$ contains the bounded frequencies.
Fix the integer $K_{\rm BW}$ so that $S_{j-K_{\rm BW}}$ selects the low
coefficient frequencies in the Bony--Weyl cutoff.
For $i\in\{-,+\}$, let $\tau_i$ be the meridional unit tangent
to $C_i$ in the direction of its parametrization, and set
\begin{equation*}
V_{\tau,i}=V\cdot\tau_i,\qquad
v_i=\frac aqV_{\tau,i},\qquad
\widetilde v_i=\langle D_\theta\rangle_i^{-1}v_i.
\end{equation*}
Let $\varrho(\xi)$ be the fixed smooth high-frequency realization
of $|\xi|$ used in the Bony--Weyl calculus, and put
\begin{equation}\label{bv-bd-inner-symbols}
x=\frac{\varrho(\xi)}j,\qquad
r_-=\Theta v_-^2x,\qquad
T_-=2i\Theta P_{\ge2}\operatorname{Op}\left(\frac{v_-}{j}\xi\right)P_{\ge2}.
\end{equation}

Choose fixed lower indices $J$ and $b_{\rm dy}$ for the coefficient
seminorms and finite tame losses in the symbol, kernel and dyadic
estimates, independently of $s$, as in
Subsection~\ref{subsec:bv-boundary-system}.
For the dyadic reduction, collect the state and its lower norms as
\begin{equation*}
z=(U,\widetilde v_-,\widetilde v_+),\qquad
L_z=\|z\|_{H^{b_{\rm dy}}},\qquad
\nu_v=\sum_i\|v_i\|_{C^J}.
\end{equation*}
For the dyadic expansion in $z$, use the product of $\mathcal A_Y$ with
an open $H^{b_{\rm dy}}$ ball for the two primitives of radius
$2C_{\rm cut}P_b\delta_b$, with this polynomial fixed before
$Y_0$. Choose $b$ so that $b\geq b_{\rm dy}+1$ and the fixed
lower profile indices used in the trace bound are at most $b-1$.

\begin{proposition}\label{prop:bv-spatial-reduction}
Under the geometric and spatial profile assumptions of
Subsection~\ref{subsec:bv-boundary-system}, the following statements hold.
\begin{enumerate}[label=\textup{(\roman*)},leftmargin=*]
\item
Let $N$ be the inverse of the Schur reduced DN operator on physical
Neumann data with zero component flux, with traces taken modulo
component constants. Its block $N_{i\ell}$ maps data on $C_\ell$
to the trace on $C_i$. For each fixed nonnegative integer $k$, there
is a nonnegative integer $d_{\rm N}$, independent of the output
index $r$, such that for every nonnegative integer or half-integer $r$,
\begin{equation*}
\begin{split}
\|N_{-+}f_+\|_{H^r(\T)}
&\leq P_r a\bigl(1+\|\eta\|_{H^{r+d_{\rm N}}}\bigr)
                    \|f_+\|_{H^{-k}},\\
\|N_{+-}f_-\|_{H^r(\Gamma_+)}
&\leq P_r a^2\bigl(1+\|\zeta_+\|_{H^{r+d_{\rm N}}}\bigr)
                    \|f_-\|_{H^{-k}(\T)}.
\end{split}
\end{equation*}
For fixed $k$, the input-boundary graph norms enter at one fixed
order independent of $r$. The factors $a$ and $a^2$ persist under
the first two shape variations at fixed $a$, with fixed input
regularity and derivative losses as specified in
\eqref{bv-bd-separated-neumann-shape}. The moving
projections and the inverse satisfy
\eqref{bv-bd-projection-first-variation}--\eqref{bv-bd-projection-second-variation}
and \eqref{bv-bd-neumann-first-variation}--\eqref{bv-bd-neumann-hadamard}.
\item
After physical flux projection and the frozen good-variable change,
the tangential transport and curvature terms in
\eqref{bv-bd-rotation-normal} contribute
\begin{equation*}
\Delta_{\rm kin}h_Y=0,\qquad
\Delta_{\rm kin}\pi_Y=-T_-\pi+P_{\ge2}\operatorname{Op}(r_-)P_{\ge2}h+R_{\pi,\rm kin},
\end{equation*}
where $r_-$ and $T_-$ are defined in \eqref{bv-bd-inner-symbols}.
The pressure lift and time derivatives of the good-variable
coefficients are treated in Subsection~\ref{subsec:bv-boundary-energies}.
On this product neighbourhood, the dyadic reduction treats the
primitives as independent inputs and gives
\begin{equation*}
\|R_{\pi,\rm kin}\|_{H^{s+1/2}}
\leq P_s\left\{
 a\sum_i\|\widetilde v_i\|_{H^{s+1/2}}
 +(L_z+\nu_v)B_q^{-1/4}E_s
 +a(L_z+\nu_v)O_s\right\}.
\end{equation*}
The primitives reconstructed from $w$ satisfy
\eqref{bv-bd-primitive-norm}; at boundary level $s$, this bound
uses the profile only through $Z_{s-1}$.
\end{enumerate}
\end{proposition}
\begin{proof}
We divide the proof into eight steps.

\medskip
\noindent\textbf{Step 1. Physical row weights and separated Neumann estimates.}
Use the densities $d_i$ and projection
$\Pi_n$ from \eqref{bv-st-capacitary-projection}.
Apply $\Pi_n$ to each term of
\eqref{bv-bd-rotation-normal}. The acceleration difference has zero
total flux, so the fluxes removed by $\Pi_n$ enter the scalar
equation \eqref{bv-ex-q-equation}.

The physical row weights follow directly from the normal acceleration.
Set
\begin{equation*}
v_i=\frac aqV_{\tau,i},\qquad
u_i=\frac aq v_{0,\tau,i},\qquad
b_i=\frac aqB_i,\qquad i\in\{-,+\}.
\end{equation*}
Tangential components are taken along the meridional unit tangents
$\tau_i$. We write $\partial_{\Sigma_i}$ for differentiation along
$\tau_i$. In the axisymmetric meridian define
\begin{align}
k_-&=-2\frac{v_-}{j}\partial_\theta b_-
+a\mathrm{II}_-(\tau_-,\tau_-)(2u_-v_-+v_-^2),
\label{bv-bd-inner-kinetic-data}\\
k_+&=-2v_+\partial_{\Sigma_+}b_+
+\mathrm{II}_+(\tau_+,\tau_+)(2u_+v_++v_+^2).
\label{bv-bd-outer-kinetic-data}
\end{align}
The corresponding physical Neumann data are
$q^2a^{-3}k_-$ and $q^2a^{-2}k_+$.  Thus the inner good momentum
receives the kinetic contribution
\begin{equation}\label{bv-bd-neumann-row-weights}
-\Theta aP_{\ge2}\operatorname{Op}(\sqrt{w_g})P_{\ge2}
\left\{a^{-1}N_{--}k_-
+N_{-+}k_+\right\}.
\end{equation}
The component projection in \eqref{bv-bd-neumann-row-weights}
acts on the physical pair $(k_-/a,k_+)$, with normalized densities
$ad_-,d_+$ and component measures $dS/a,dS$, respectively.
Indeed, $(ad_-)\int f\,dS/a=d_-\int f\,dS$.
We use this convention for the subsequent rows and their finite-rank terms.

For each fixed nonnegative integer $k$, there is a nonnegative
integer $d_{\rm N}$, independent of the output index $r$, such that
for every nonnegative integer or half-integer $r$,
\begin{equation}\label{bv-bd-separated-neumann}
\begin{split}
\|N_{-+}f_+\|_{H^r(\T)}
&\leq P_r a\bigl(1+\|\eta\|_{H^{r+d_{\rm N}}}\bigr)
                    \|f_+\|_{H^{-k}},\\
\|N_{+-}f_-\|_{H^r(\Gamma_+)}
&\leq P_r a^2\bigl(1+\|\zeta_+\|_{H^{r+d_{\rm N}}}\bigr)
                    \|f_-\|_{H^{-k}(\T)}.
\end{split}
\end{equation}
The half-integer cases follow from the next integer estimate by
Sobolev embedding, after increasing $d_{\rm N}$ by one if necessary.
For fixed $k$, the input-boundary graph norms enter at one fixed order
independent of $r$. Traces are taken modulo component constants,
and all Neumann inputs have zero component flux.

To prove \eqref{bv-bd-separated-neumann}, we adapt
\cite[Lemmas~2.9 and 4.10]{LYIrrotational} to Neumann data.
We first use smooth compatible data; the estimates then define
the operators on the stated negative Sobolev spaces by transposition.
Let $z$ solve the Neumann problem with outer datum $f_+$ in the
filled domain. Interior regularity and duality bound its fixed
physical derivatives near the cavity by $C\|f_+\|_{H^{-k}}$.
The inner trace satisfies
\begin{equation*}
z\bigl(c+a(1+\eta)e_r\bigr)-z(c)
=a\int_0^1\nabla z\bigl(c+ta(1+\eta)e_r\bigr)
                    \cdot(1+\eta)e_r\,dt.
\end{equation*}
The inner correction has data $-\partial_{n_-}z$ and zero outer
datum. Since $z$ is harmonic through the cavity,
\begin{equation*}
\int_{\Gamma_-}\partial_{n_-}z\,dS=0.
\end{equation*}
Thus \eqref{bv-st-hole-energy} and the scaled trace estimate modulo
constants give the same factor $a$ for the correction. Tame boundary
regularity yields the first row of \eqref{bv-bd-separated-neumann}.
For the second row, Green's identity gives
\begin{equation*}
\int_{\Gamma_-}f_-N_{-+}f_+\,dS
=\int_{\Gamma_+}f_+N_{+-}f_-\,dS,
\qquad dS|_{\Gamma_-}=2\pi a\rho_-j\,d\theta.
\end{equation*}
Choose nested fixed intermediate collars separated from both
boundaries, and let $z_-$ be the Neumann
solution with data $(f_-,0)$. For a smooth axisymmetric mean-zero
volume test $f$ supported in the larger collar $\mathcal C$, let $z_f$ solve
$-\Delta z_f=f$ with homogeneous Neumann data. Decompose $z_f$ into
the filled-domain Neumann solution $z_f^{\rm fill}$ and its inner-hole
correction.

The filled solution is harmonic near the cavity, so its
inner trace modulo constants has a factor $a$. The correction has
inner datum $-\partial_{n_-}z_f^{\rm fill}$ and zero outer datum.
The estimate \eqref{bv-st-hole-energy} and scaled boundary regularity
give the same factor for its trace. Extend $f$ by zero to $\R^3$
when taking its negative Sobolev norm. For each fixed $m$ this yields
\begin{equation*}
\inf_{c\in\R}
\|z_f|_{\Gamma_-}-c\|_{H^k(\T)}
\leq P_{m,k}a\|f\|_{H^{-m}(\R^3)},
\end{equation*}
with the inner graph measured only at the fixed order required by
$k$. Since $f_-$ has zero flux, Green's identity gives
\begin{equation*}
\left|\int_\Omega z_-f\,dx\right|
=\left|\int_{\Gamma_-}f_-z_f\,dS\right|
\leq P_{m,k}a^2\|f_-\|_{H^{-k}(\T)}
                         \|f\|_{H^{-m}(\R^3)}.
\end{equation*}
The second factor $a$ comes from the inner surface measure. Duality
bounds $z_-$ modulo constants in $H^m$ on the smaller collar.
Outside that collar $z_-$ is harmonic with homogeneous outer Neumann
data. Boundary regularity puts the high output norm on $\zeta_+$
and proves the second row of \eqref{bv-bd-separated-neumann}.

\medskip
\noindent\textbf{Step 2. Shape variations of the component projections.}
For the shape derivatives, write the weak Neumann equation on a
fixed pullback domain as
\begin{equation*}
\int A\nabla z\cdot\nabla\phi=\mathcal L(\phi),
\end{equation*}
where $A$ is the elliptic coefficient matrix and $\mathcal L$ is
the Neumann boundary pairing. With subscripts denoting variations
in two fixed directions, differentiation gives
\begin{align*}
\int A\nabla z_1\cdot\nabla\phi
&=\mathcal L_1(\phi)-\int A_1\nabla z\cdot\nabla\phi,\\
\int A\nabla z_{12}\cdot\nabla\phi
&=\mathcal L_{12}(\phi)
-\int\bigl(A_{12}\nabla z+A_1\nabla z_2+A_2\nabla z_1\bigr)
\cdot\nabla\phi.
\end{align*}
Here the boundary pairing includes the moving component projection;
write $\Pi_i$ for the projection on the $i$th component.
If $J_i\,d\theta$ is the surface measure in the chosen chart,
the first variation of $\Pi_i$ at fixed pulled-back datum is
\begin{equation}\label{bv-bd-projection-first-variation}
D\Pi_i[H]f
=-Dd_i[H]\int fJ_i\,d\theta
-d_i\int f\,DJ_i[H]\,d\theta.
\end{equation}
A second differentiation at fixed pulled back datum gives
\begin{equation}\label{bv-bd-projection-second-variation}
\begin{split}
D^2\Pi_i[H_1,H_2]f
={}&-D^2d_i[H_1,H_2]\int fJ_i\,d\theta
-Dd_i[H_1]\int f\,DJ_i[H_2]\,d\theta\\
&-Dd_i[H_2]\int f\,DJ_i[H_1]\,d\theta
-d_i\int f\,D^2J_i[H_1,H_2]\,d\theta.
\end{split}
\end{equation}
Each term is a product of a range density and a scalar boundary pairing.
Compatibility holds for the complete differentiated boundary pairing:
\begin{equation*}
\int (D\Pi_i[H]f)J_i\,d\theta
+\int (\Pi_i f)DJ_i[H]\,d\theta=0.
\end{equation*}
The individual datum $D\Pi_i[H]f$ need not have zero flux in the
undifferentiated measure.
Differentiating $\int(\Pi_i f)J_i\,d\theta=0$ twice gives the
second-order compatibility identity, including the first and second variations of
$J_i$. These combined boundary pairings enter the differentiated
weak Neumann equations, which are inverted on the fixed identified spaces.
In the inner chart,
$J_-=2\pi a\rho_-j$, and $ad_-$, rather than $d_-$, has uniform
normalized bounds.

\medskip
\noindent\textbf{Step 3. Shape estimates for the separated Neumann blocks.}
For the first two shape differentials, put
$\boldsymbol\zeta=(\eta,\zeta_+)$ and $\zeta_-:=\eta$.
At fixed $a$, the differentiated weak Neumann equations and
duality argument give, for $j=1,2$ and $i\ne\ell$,
\begin{equation}
\label{bv-bd-separated-neumann-shape}
\begin{split}
\bigl\|D^j(a^{-\gamma_{i\ell}}N_{i\ell})
[H_1,\ldots,H_j]f\bigr\|_{H^r}
\leq P_r\|f\|_{H^{-k}}\bigg\{
&(1+\|\zeta_i\|_{H^{r+d_{\rm N}}})
\prod_{\nu=1}^j\|H_\nu\|_{H^{b_{\rm N}}}\\
&+\sum_{\nu=1}^j\|H_{\nu,i}\|_{H^{r+d_{\rm N}}}
\prod_{\mu\ne\nu}\|H_\mu\|_{H^{b_{\rm N}}}\bigg\},
\end{split}
\end{equation}
where $\gamma_{-+}=1$, $\gamma_{+-}=2$. The fixed indices
$b_{\rm N},d_{\rm N}$ may depend on $k$, but not on $r$; increase
$d_{\rm N}$ from \eqref{bv-bd-separated-neumann} if necessary.
The differentials include the moving projections, and $H_{\nu,i}$
denotes the graph direction on the output boundary.

For example, consider two inner variations of $N_{+-}$ at fixed
$a,c,\zeta_+$, with
$\eta^\tau=\eta+\tau_1H_{1,-}+\tau_2H_{2,-}$.
Let $z_-^\tau$ have data $(\Pi_-^\tau f_-,0)$ for a fixed pulled
back datum $f_-$, with zero flux at $\tau=0$.
Take a smooth axisymmetric mean-zero test $f$ in the fixed
intermediate collar, where the pullbacks are the identity.
For its homogeneous Neumann solution $z_f^\tau$, put
\begin{equation*}
h_f(\tau)=z_f^\tau|_{\Gamma_-^\tau}
-\int_{\Gamma_-^\tau}d_-^\tau z_f^\tau\,dS^\tau,
\end{equation*}
with the trace pulled back to $\T$. The moving input projection
and Green's identity give
\begin{equation*}
\int_{\Omega^\tau}z_-^\tau f\,dx
=\int_\T f_-J_-^\tau h_f(\tau)\,d\theta.
\end{equation*}
Split $z_f^\tau$ into its fixed filled-domain solution and hole
correction. The former has an $O(a)$ inner trace modulo constants,
also after two inner variations. Divide the correction by $a$
modulo constants. Its complete boundary pairing and its first two
variations annihilate constants. The displacements $aH_{\nu,-}$
give uniformly bounded scaled coefficient variations. Thus the
energy test in \eqref{bv-st-hole-energy}, applied successively to
the correction and its two weak variation equations above, bounds
their gradients by $\|f\|_{H^{-1}(\R^3)}$ times the corresponding
products of lower direction norms. Fixed-order boundary regularity
gives the trace bounds modulo constants. The condenser bounds for
$ad_-$ and its two variations, together with
$\int_\T d_-J_-\,d\theta=1$ and its derivatives, therefore give
\begin{equation*}
\|D^jh_f[H_1,\ldots,H_j]\|_{H^k(\T)}
\leq P_ka\|f\|_{H^{-1}(\R^3)}
\prod_{\nu=1}^j\|H_{\nu,-}\|_{H^{b_{\rm N}}},
\qquad j=0,1,2.
\end{equation*}
Since $J_-=a(2\pi\rho_-j)$, its first two variations also retain
$a$. Differentiating at $\tau=0$ gives
\begin{equation*}
\partial_{\tau_1\tau_2}\int_{\Omega^\tau}z_-^\tau f\,dx
=\int_\T f_-\bigl(J_-h_{f,12}+(J_-)_1h_{f,2}
+(J_-)_2h_{f,1}+(J_-)_{12}h_f\bigr)\,d\theta.
\end{equation*}
Each term is bounded by
$P_ka^2\|f_-\|_{H^{-k}}\|f\|_{H^{-1}}
\|H_{1,-}\|_{H^{b_{\rm N}}}\|H_{2,-}\|_{H^{b_{\rm N}}}$.
Duality gives an $H^1$ collar bound modulo constants with fixed
input order $k$. On the outer collar the inner coefficient
variations vanish, so higher boundary regularity uses only
$\zeta_+$. Outer or mixed variations add the corresponding high
output-direction norms. For $N_{-+}$ the same argument retains
the factor $a$ in the filled trace difference and hole correction.
These estimates prove the displayed bounds.

The diagonal inverse blocks have local order $-1$; their smooth
remainders satisfy the same fixed-loss product estimates after
normalizing $N_{--}$ by $a^{-1}$.

\medskip
\noindent\textbf{Step 4. The differentiated inverse and the Hadamard formula.}
Identify the Schur-reduced trace and flux spaces with fixed quotient
spaces by the normalized component projections.
Let $G_{\rm red}$ denote the resulting DN operator and
$N_{\rm red}=G_{\rm red}^{-1}$. Then
\begin{align}
DN_{\rm red}[H]
&=-N_{\rm red}\,DG_{\rm red}[H]N_{\rm red},
\label{bv-bd-neumann-first-variation}\\
D^2N_{\rm red}[H_1,H_2]
&=N_{\rm red}\,DG_{\rm red}[H_1]N_{\rm red}
  \,DG_{\rm red}[H_2]N_{\rm red} +N_{\rm red}\,DG_{\rm red}[H_2]N_{\rm red}
  \,DG_{\rm red}[H_1]N_{\rm red}\notag\\
&\quad-N_{\rm red}\,D^2G_{\rm red}[H_1,H_2]N_{\rm red}.
\label{bv-bd-neumann-second-variation}
\end{align}
In physical coordinates, let $G$ denote the full DN map and $N$
the Neumann inverse defined above. Set
$\mathcal N=N\Pi_n$, $G^0=\Pi_nG$ and $P^\circ=(P_-^0,P_+^0)$.
Thus $\mathcal N$ takes arbitrary data to the two mean-free traces
of their projected Neumann solution. For every trace pair $h$,
\begin{equation*}
 \mathcal NGh=NG^0h=P^\circ h .
\end{equation*}
Indeed the total flux of $Gh$ is zero, and the part removed by
$\Pi_n$ is the normal derivative of a multiple of $u_c$, whose
traces are component constants.

Let $f$ be arbitrary smooth physical Neumann data, held fixed in the
pullback by a normal variation $H=(H_-,H_+)$. Let $\phi$ be the
harmonic solution with normal data $\Pi_n f$ and put
\begin{equation*}
 g_c=\partial_nu_c,\quad \epsilon_-=1,\quad\epsilon_+=-1,\qquad
 m_i=\int_{\Gamma_i} f_i\,dS,\quad
 D_f=\operatorname{diag}(\epsilon_i m_i),\qquad
 h_c=H g_c/\operatorname{Cap}.
\end{equation*}
All products of boundary functions here are componentwise. With
$[D_f,G^0]=D_fG^0-G^0D_f$, the projected Hadamard formula is
\begin{equation}\label{bv-bd-neumann-hadamard}
\begin{split}
 D\mathcal N[H]f
 ={}&P^\circ(Hf)+
 \mathcal N\{\operatorname{div}_\Sigma(H\nabla_\Sigma\phi)
                         +\mathcal H Hf\}\\
 &+N[D_f,G^0]h_c .
\end{split}
\end{equation}
To prove it, denote the Eulerian harmonic variation by $\phi'$.
Differentiation of its Neumann condition gives
\begin{equation*}
 \partial_n\phi'=
 D(\Pi_n f)[H]+\operatorname{div}_\Sigma(H\nabla_\Sigma\phi)
                         +\mathcal H H\,\Pi_n f .
\end{equation*}
The Eulerian variation of $u_c$ has trace $-Hg_c$ because both
Dirichlet values are constant. Thus the pulled-back derivative of
$g_c$ is $-G(Hg_c)-\mathcal H H g_c$. As
$d_i=\epsilon_i g_{c,i}/\operatorname{Cap}$, substitution in
$D(\Pi_n f)_i=-m_i Dd_i-d_i Dm_i$ yields
\begin{equation*}
 \partial_n\phi'=
 D_fG h_c+\operatorname{div}_\Sigma(H\nabla_\Sigma\phi)+\mathcal H H f
 +\bigl(d_i(m_iD\operatorname{Cap}/\operatorname{Cap}-Dm_i)\bigr)_i .
\end{equation*}
In this display the scalar differentials are in direction $H$.
The last pair is killed by $\Pi_n$. The normal datum of $\phi'$ has
zero total flux, so its mean-free trace is $\mathcal N(\partial_n\phi')$.
As $\Pi_n$ commutes with $D_f$,
\begin{equation*}
 \mathcal N D_fG h_c=P^\circ(D_f h_c)+N[D_f,G^0]h_c .
\end{equation*}
The term $D_fh_c$ is exactly $(H_i d_i m_i)_i$.
Adding the trace evaluation $P^\circ(H\Pi_n f)$ proves
\eqref{bv-bd-neumann-hadamard}. The commutator has only off-diagonal
blocks $(\epsilon_i m_i-\epsilon_jm_j)G^0_{ij}$. In a graph chart
a tangential displacement $T$ adds the trace transport
$P^\circ(T\cdot\nabla_\Sigma\phi)$ and the input correction
$-\mathcal N(T\cdot\nabla_\Sigma f)$. Both are of order at most zero
on the displacement at a frozen smooth state.

\medskip
\noindent\textbf{Step 5. The kinetic symbol and boundary primitives.}
Let $H$ be a normal boundary displacement and $\varphi$ the
Eulerian good trace of the potential variation.  At fixed
component flux, the local highest order parts of the harmonic
shape formulas are
\begin{equation*} 
\delta B=G\varphi+d_0H-v_{0,\tau}\partial_\Sigma H,
\qquad d_0=n\cdot\nabla v_0n,
\qquad
\delta\mathrm{II}(\tau,\tau)
=-\partial_\Sigma^2H+\text{terms of order at most one}.
\end{equation*}
The normal derivative contributes $d_0H$. Component projections,
tangential chart terms and the other boundary contribute lower order
or separated blocks. Substitution in the two kinetic terms of
\eqref{bv-bd-rotation-normal} gives
\begin{equation}
\label{bv-bd-mixed-cancellation}
-2V_\tau\partial_\Sigma G\varphi
+2V_\tau v_{0,\tau}\partial_\Sigma^2H
-(2v_{0,\tau}V_\tau+V_\tau^2)\partial_\Sigma^2H =-2V_\tau\partial_\Sigma G\varphi
-V_\tau^2\partial_\Sigma^2H.
\end{equation}
In particular, the mixed $v_{0,\tau}V_\tau$ height symbol cancels.

Write $\Delta_{\rm kin}$ for the contribution of the tangential
transport and curvature terms in \eqref{bv-bd-rotation-normal},
after the projected Neumann inverse and the frozen good variable
change. The pressure lift and the time derivatives of the good
variable coefficients are treated in Subsection~\ref{subsec:bv-boundary-energies}.

The normal good variables satisfy, at principal order,
$H=-h/\sqrt{w_g}$ and
$\varphi=-(q/a)\pi/\sqrt{w_g}$.
Since the principal symbol of $G_{--}$ is $|\xi|/(aj)$,
\eqref{bv-bd-neumann-row-weights} and
\eqref{bv-bd-mixed-cancellation} give
\begin{equation}\label{bv-bd-kinetic-principal}
\Delta_{\rm kin}h_Y=0,\qquad
\Delta_{\rm kin}\pi_Y=-T_-\pi+P_{\ge2}\operatorname{Op}(r_-)P_{\ge2}h+R_{\pi,\rm kin}.
\end{equation}
The height identity is exact because $V\cdot n=0$. The
order-zero difference between the unreduced transport operator and
the skew-symmetric representative $T_-$ is included in
$R_{\pi,\rm kin}$.
To estimate the kinetic remainder, introduce the boundary primitives
\begin{equation}\label{bv-bd-external-primitives}
\widetilde v_i=\langle D_\theta\rangle_i^{-1}v_i,
\qquad i\in\{-,+\}.
\end{equation}
Regard $\widetilde v_\pm$ as independent of
$U=(\eta,F,\zeta_+,g)$ in the spatial reduction.  The trace estimate
\eqref{bv-st-smoothed-vtrace} and the outer counterpart give
\begin{equation}\label{bv-bd-primitive-norm}
a\sum_i\|\widetilde v_i\|_{H^{s+1/2}}
\leq P_s\left\{
\frac{a^2}{q}(\varepsilon+Z_{s-1})
+\frac{\varepsilon\sqrt a}{\sqrt\sigma}E_s
+\frac{a^2\varepsilon}{q}O_s\right\}.
\end{equation}
For the inner graph term, \eqref{eq:bv-interpolation-costs} gives
\begin{equation*}
\frac{a^2\varepsilon}{q}\|\eta\|_{H^{s+1}}
\leq P_s\frac{\varepsilon\sqrt a}{\sqrt\sigma}E_s.
\end{equation*}

\medskip
\noindent\textbf{Step 6. Tame differentials of the kinetic terms.}
Write the Neumann expression in braces in
\eqref{bv-bd-neumann-row-weights} as
$K_g(U,\widetilde v_-,\widetilde v_+)$.
The weak boundary bounds, the profile bounds in
\eqref{bv-bd-profile-assumptions} at the fixed spatial orders used here,
and the definition of $K_g$ give
\begin{equation*}
L_z\leq P_b\delta_b,\qquad \nu_v\leq P_ba/q,\qquad
K_g(U,0,0)=0,\qquad D_UK_g(U,0,0)=0.
\end{equation*}
The cutoffs and their Taylor segments then remain in this product
domain. The extra inputs enter $K_g$ polynomially and change none
of the geometric or scalar denominators.
Since $v_i=\langle D_\theta\rangle_i\widetilde v_i$ and $N$ has order $-1$,
\begin{equation*}
\operatorname{ord}(D_UK_g)\leq1,\qquad
\operatorname{ord}(D_{\widetilde v_i}K_g)\leq0.
\end{equation*}
The principal terms are determined by the cancellation
\eqref{bv-bd-mixed-cancellation}.
The lower local symbols have orders $0,-1,\ldots$. More precisely,
there are fixed indices $r_0,d_{\rm kin}$, independent of the output index
$r$, such that for $r\geq r_0$ on the lower neighbourhood
\begin{align*} 
\|DK_g(z)[h]\|_{H^{r-d_{\rm kin}}}
&\leq P_r\bigl\{\|h\|_{H^r}
+(1+\|z\|_{H^r})\|h\|_{H^{r_0}}\bigr\},\\
\|D^2K_g(z)[h,k]\|_{H^{r-d_{\rm kin}}}
&\leq P_r\bigl\{\|h\|_{H^r}\|k\|_{H^{r_0}}
+\|k\|_{H^r}\|h\|_{H^{r_0}} +(1+\|z\|_{H^r})
\|h\|_{H^{r_0}}\|k\|_{H^{r_0}}\bigr\}. 
\end{align*}
To verify these bounds, first consider the curvature term in $k_-$.
In the inner chart,
\begin{equation*}
a\mathrm{II}_-(\tau_-,\tau_-)
=-\frac{d^2+2\eta_\theta^2-d\eta_{\theta\theta}}{j^3}.
\end{equation*}
Its variation in the graph direction $H$ is
\begin{equation}\label{bv-bd-curvature-variation}
D_\eta\bigl(a\mathrm{II}_-(\tau_-,\tau_-)\bigr)[H]
=\frac d{j^3}\partial_\theta^2H
+\alpha_1\partial_\theta H+\alpha_0H,
\end{equation}
where $\alpha_0,\alpha_1$ are smooth functions of
$\eta,\eta_\theta,\eta_{\theta\theta}$ on the lower neighbourhood.
The normalized inverse $a^{-1}N_{--}$ has order $-1$.
Thus the first term in \eqref{bv-bd-curvature-variation} has order
one after inversion, while the last two have order at most zero.
The corresponding normal displacement is $-adH/j$.
In this normal variable, combining the leading curvature and
transport terms gives \eqref{bv-bd-mixed-cancellation}.
At fixed $U$, differentiation with respect to the inner primitive gives
\begin{equation*}
D_{\widetilde v_-}k_-[H]
=\left\{-\frac2j\partial_\theta b_-
+2a\mathrm{II}_-(\tau_-,\tau_-)(u_-+v_-)\right\}
\langle D_\theta\rangle H.
\end{equation*}
The inverse therefore reduces this column to order zero.
The outer local column has the same order count in its physical chart.

Next differentiate the normalized Neumann inverse.
For an inner graph direction $H$, the physical variation is $aHe_r$,
with normal component $-adH/j$. In
\eqref{bv-bd-neumann-hadamard}, its factor $a$ cancels the
normalization $a^{-1}$ in the direct product with $f$.
In the other terms, tangential differentiation and the Neumann inverse
have reciprocal factors $a^{-1}$ and $a$ in the scaled chart.
The factor $a\mathcal H$ is bounded in the fixed lower norms.
The differentiated weak equations in Step~2 consequently give the same
fixed-loss bounds for the first two variations of $a^{-1}N_{--}$.
For the separated blocks, the two variations retain the factors
$a$ and $a^2$ in \eqref{bv-bd-separated-neumann}.

\medskip
\noindent\textbf{Step 7. The projection terms.}
For the projection in the kinetic row take $f=(k_-/a,k_+)$ in
\eqref{bv-bd-neumann-hadamard}. The inner moment is
\begin{equation*}
m_-=\int_{\Gamma_-} k_-\,dS/a.
\end{equation*} 
Thus $|m_-|+|m_+|\leq P_b\nu_v$
by \eqref{bv-bd-inner-kinetic-data}--\eqref{bv-bd-outer-kinetic-data}
and the bounds on the fixed low spatial norms of the geometry and
boundary variables. For a dimensionless inner graph direction
$\delta\eta$, $H_-=-ad\,\delta\eta/j$. Consequently both $H_- f_-$
and $(h_c)_-=(H_-/a)(ad_-)$ have bounded scaled coefficients,
with a factor $P_b\nu_v$ in the former. The same statements at the
outer boundary use its physical scale. In the last term of
\eqref{bv-bd-neumann-hadamard} the moment supplies that factor.
More explicitly, with $S_a=\operatorname{diag}(a,1)$ the product is
\begin{equation*}
 (N S_a^{-1})[D_f,S_a G^0]h_c .
\end{equation*}
The diagonal inverse blocks are now of order $-1$ with polynomial
bounds. The off-diagonal inverse blocks obey
\eqref{bv-bd-separated-neumann}, with the inner input divided by $a$.
The off-diagonal blocks of $S_aG^0$ are the scaled zero-flux
transfers of \cite[Lemmas~2.9 and 4.10]{LYIrrotational}.
Lemma~4.10 there is stated for a graph ball of radius
$2C_{\rm cut}a^{3/2}$. Its proof of the outer-to-inner estimate also
applies on the larger ball $\mathcal A_Y$ in \eqref{bv-bd-admissible-ball}.
Indeed, let $H_\gamma$ extend the outer datum $\gamma$ harmonically
to the filled domain. Taylor's formula gives
\begin{equation*}
H_\gamma(X_-)-H_\gamma(c)=a\int_0^1
\nabla H_\gamma(c+ta d e_r)\cdot d e_r\,dt.
\end{equation*}
The separation of the boundaries controls the derivatives of $H_\gamma$
near $c$. The zero-flux correction thus has an inner trace with a factor
$a$ modulo constants. After division by $a$, its equation on the rescaled
collar is uniformly elliptic on $\mathcal A_Y$. At fixed $a$, the first two graph
variations satisfy the differentiated elliptic estimates on this same
collar, so the factor $a$ persists in $aG^0_{-+}$ and its variations
\cite[equations~(4.35)--(4.37)]{LYIrrotational}.
Thus the commutator passes between separated boundaries before the
bounded normalized inverse acts. Differentiating this product uses
the separated estimates of Steps~1 and~3 and the corresponding DN shape
estimates. At fixed $a$, the full kinetic moment
\begin{equation*}
m_i=\int_\T f_iJ_i\,d\theta
\end{equation*} 
satisfies
\begin{equation*}
Dm_i[Z]=\int_\T Df_i[Z]J_i\,d\theta+\int_\T f_iDJ_i[Z]\,d\theta.
\end{equation*}
The first term differentiates $f=(k_-/a,k_+)$ as above. The second
is the moving-measure contribution in
\eqref{bv-bd-projection-first-variation}. The local density
variation has already been combined with the trace evaluation in
\eqref{bv-bd-neumann-hadamard}.

In particular the order-zero inner height column with the fixed
coefficient cutoff costs at most
\begin{equation*}
P_s a\nu_v\|\eta\|_{H^{s+1/2}}
 \leq P_s\nu_v B_q^{-1/4}E_s
\end{equation*} 
after the row factor in
\eqref{bv-bd-neumann-row-weights}. Its outer height input costs
$P_s a\nu_v O_s$. Here
\begin{equation*}
\nu_v B_q^{-1/4}\leq P_b\sigma^{-1/4}a^{3/4}q^{-1/2}
\end{equation*} 
is
integrable. The separated kernel differences enter the dyadic
remainder below, with the same fixed losses as their factors.
This gives the projection terms in the bounds for the first two
differentials and accounts for them in \eqref{bv-bd-kinetic-remainder}.

\medskip
\noindent\textbf{Step 8. The dyadic remainder.}
These calculations also verify the local symbol hypotheses of
\cite[Lemma~4.4]{LYIrrotational}. The differential in $U$ has order
at most one, and the primitive columns have order at most zero.
After retaining the order-zero symbols, the remaining local terms
have order at most $-1$. Differentiating the finite Poisson
symbol recursion gives the corresponding fixed-loss symbol bounds.
Integrating their first variations along $tz$, $0\leq t\leq1$, gives
the bounds for the differences from the reference symbols.
The same estimates on the fixed lower neighbourhood give their
lower norm Lipschitz bounds. For a fixed input index $k$, the smooth
kernel part $\mathcal K_g$ satisfies
\begin{equation*}
\|[\mathcal K_g(z)-\mathcal K_g(0)]f\|_{H^r}
\leq P_r\|z\|_{H^{r+d_k}}\|f\|_{H^{-k}},
\end{equation*}
with $d_k$ independent of $r$. This follows by integrating its
first variation along $tz$, $0\leq t\leq1$. The separated estimates
\eqref{bv-bd-separated-neumann} and their shape differentials place
output derivatives only on the output graph. For the commutator in
\eqref{bv-bd-neumann-hadamard}, the off-diagonal DN block supplies
the smooth kernel, while scalar moment differentials act by boundary
pairing. The local term $P^\circ(Hf)$ belongs to the order-zero
symbol. The reference kernel is smooth. The dyadic expansion retains
the differential at $z=0$ in all components
$U,\widetilde v_-,\widetilde v_+$.

By the choice of $b_{\rm dy}$, apply
\cite[Lemma~4.4]{LYIrrotational} at exponent $s+1/2$, placing every
component of $z$ in $H^{s+1/2}$. The local orders are at most two,
and the symbol and kernel differences vanish at the reference state.
The dyadic remainder therefore satisfies
\begin{equation*} 
a\|R_{\rm hh}\|_{H^{s+5/2}}
\leq P_s a L_z\|z\|_{H^{s+1/2}}.
\end{equation*}
After applying the fixed-cutoff half-density change to the output
and expressing the boundary inputs through the frozen inverse of
\eqref{eq:bv-good-variables}, the order-one terms are exactly those
in \eqref{bv-bd-kinetic-principal}. The next composition terms
have order zero and contain a factor $v_i$.  Hence
\begin{equation}\label{bv-bd-kinetic-remainder}
\|R_{\pi,\rm kin}\|_{H^{s+1/2}}
\leq P_s\left\{
a\sum_i\|\widetilde v_i\|_{H^{s+1/2}}
+(L_z+\nu_v)B_q^{-1/4}E_s
+a(L_z+\nu_v)O_s\right\}.
\end{equation}
The factor $B_q^{-1/4}$ includes the height conversion in the
quadratic dyadic bound.

\end{proof}

\subsection{The boundary energies and terminal estimates}
\label{subsec:bv-boundary-energies}

We combine Proposition~\ref{prop:bv-spatial-reduction} with the
pressure and scalar estimates to construct the two boundary energies. At each integer energy level $s\geq b$, impose
\eqref{bv-bd-profile-assumptions} through order $s-1$.
For the first $Y$ derivative of the prescribed profile, assume
\begin{equation}\label{bv-bd-time-assumptions}
M_{k_{\rm t}}(w_Y)\leq P_ba^{-1/4},\qquad
\|w_Y\|_{L^\infty}\leq P_b\varepsilon a^{3/4}.
\end{equation}
The number $\varepsilon$ is sufficiently small, depending on the
fixed positive $\sigma$, the geometric data, $E_c$, the fixed
boundary bounds at level $b$ and spatial profile bounds through
order $b-1$.

\begin{lemma}\label{bv-boundary-energy}
Choose the base index $b$ large enough for the fixed derivative orders
above. Once the bounding polynomials at those orders are fixed, there
is a sufficiently late $Y_0$, independent of $s$, such that for every
integer $s\geq b$ there are positive quadratic energies, with square
roots $E_s$ and $O_s$ equivalent to
\eqref{eq:bv-energy-norms}, satisfying
\begin{align}
-\frac{dE_s}{dY}
&\leq (1+e_s(Y))E_s
+P_s\left\{
a^2(\varepsilon+Z_{s-1})
+a^{5/4}O_s+a^{5/2}+a^3\right\},
\label{bv-bd-inner-inequality}\\
\left|\frac{dO_s}{dY}\right|
&\leq e_s(Y)O_s+P_s(a+E_s).
\label{bv-bd-outer-inequality}
\end{align}
These differential inequalities are understood almost everywhere in $Y$,
with
\begin{equation}\label{bv-bd-integrable-coefficient}
0\leq e_s(Y)
\leq C_s(1+Y)^{-3/2}+P_s(Y)a^{1/4}.
\end{equation}
The constants are uniform in $Y_N$ and the minimum cavity radius.
The highest profile norm in \eqref{bv-bd-inner-inequality} is $Z_{s-1}$.

The first harmonic coefficient $\beta$ satisfies
\begin{equation}\label{bv-bd-first-mode-equation}
\beta_Y=(1-p_q)\beta+r_\beta,
\qquad |r_\beta|\leq P_b(a+\delta_b).
\end{equation}
If $\beta(Y_N)=0$, then
$|\beta(Y)|\leq P_ba(Y)^{3/4}$.
\end{lemma}

\begin{proof}
The spatial terms are given by Proposition~\ref{prop:bv-spatial-reduction}.
We first collect the full evolution, including the time derivatives
of the good-variable coefficients and the scalar normalization.
The inner and outer energy estimates then leave only the first modes
to be controlled.

\medskip
\noindent\textbf{Step 1. The full boundary equations and their sources.}
The full and instantaneous potential evolutions have the same
geometric velocities, since $V\cdot n=0$. Thus the definition
\eqref{eq:bv-good-variables} gives the exact identity
\begin{equation}\label{bv-bd-good-coefficient-slot}
\Delta\pi_Y
=-aP_{\ge2}\operatorname{Op}(\sqrt{w_g})P_{\ge2}\Delta F_Y
-aP_{\ge2}\operatorname{Op}(\sqrt{w_g}\,\Delta b_{{\rm g},Y})P_{\ge2}\eta.
\end{equation}
The second term comes from differentiating $b_g$ in the good
variable. For a scalar coefficient $a_g$, the estimate with a fixed
cutoff and \eqref{bv-bd-rotation-normal}, together with normalized
harmonic reconstruction, give
\begin{equation}\label{bv-bd-coefficient-slot-bound} 
a\|\operatorname{Op}(a_g)\eta\|_{H^{s+1/2}}\leq C_s\|a_g\|_{C^J}B_q^{-1/4}E_s,\qquad
\|\Delta b_{{\rm g},Y}\|_{C^J}
\leq P_J\frac aq(\varepsilon+Z_{J+c_0}). 
\end{equation}
Here $J,c_0$ are fixed before $b$. The high-frequency kinetic
terms are estimated by \eqref{bv-bd-kinetic-remainder}.

The pressure contributions satisfy \eqref{bv-st-pressure-inner-result} and
\eqref{bv-st-pressure-outer-result}. The multipliers with a fixed cutoff in
\eqref{bv-bd-good-coefficient-slot} preserve those bounds.
Thus both the pressure and kinetic terms use the bulk index $s-1$.

The normalization of the inner trace by $q$ also contributes a scalar term.
Let $q_Y^0$ be the scalar velocity of the instantaneous
potential problem, and put $\delta p_q=(q_Y-q_Y^0)/q$.
Then $\Delta F_Y$ contains $-\delta p_qF$.  To estimate it,
write $v_0=(Q/\operatorname{Cap})\nabla u_c+u_h$, where $u_h$
has zero flux on each component. The identity
\begin{equation*}
\int_\Omega D^2u_c[\nabla u_c,V]
=\frac12\int_{\partial\Omega}|\nabla u_c|^2V\cdot n=0
\end{equation*}
removes the flux and rotation cross term from
\eqref{bv-ex-q-equation}.
The remaining integrand is
\begin{equation*}
2D^2u_c[u_h,V]+D^2u_c[V,V].
\end{equation*}
The static condenser bounds and
$\|u_h\|_\infty\leq P_bq\delta_b/a$ therefore imply
\begin{equation}\label{bv-bd-scalar-rotation}
|\delta p_q|
\leq P_b\left(\frac{a\delta_b\varepsilon}{q}
+\frac{a^2\varepsilon^2}{q^2}\right).
\end{equation}
This integrable order-zero coefficient contributes through
\eqref{bv-bd-good-coefficient-slot}.

We combine the rotational contributions with the potential
boundary estimates. The spatial reduction applies on the convex
neighbourhood \eqref{bv-bd-admissible-ball}, with the denominator
bounds \eqref{bv-bd-clock-denominators}. For time derivatives of
the good coordinates, use the full boundary evolution, obtained
by adding the rotational corrections
\eqref{bv-bd-rotation-normal}, \eqref{bv-bd-good-coefficient-slot}
and \eqref{bv-bd-scalar-rotation} to the instantaneous potential
evolution with scalar velocity $q_Y^0$.

Apply the spatial expansion in $U$ at fixed $(a,c,q,\beta)$.
The time dependence of $q$ enters through $p_q$, differentiation
of the good coordinates and $B_{q,Y}/B_q$.
Applying \cite[Lemma~4.4]{LYIrrotational} at exponent $s$ gives
\begin{equation}\label{bv-bd-potential-high-high}
a\|R^0_{\rm hh}\|_{H^{s+2}}
\leq P_s\delta_b(E_s+aO_s).
\end{equation}
Here $a\|U\|_{H^s}\le P_s(E_s+aO_s)$.
With the norm for an inner source pair $(f,g)$ defined by
\begin{equation*}
Q_s(f,g)=\|f\|_{H^s}+\sqrt{B_q}\|f\|_{H^{s+1}}
+\|g\|_{H^{s+1/2}},
\end{equation*} 
the scalar, lift and projection remainder pairs
$R^0_{{\rm sm},-}$ and $R^0_{{\rm sm},+}$ satisfy
\begin{equation*} 
\begin{split}
Q_s(R^0_{{\rm sm},-})
&\le P_s\{aE_s+a^2O_s+\delta_b(E_s+aO_s)\},\\
\|R^0_{{\rm sm},+}\|_{H^{s+2}}
&\le P_s\{aE_s+a^2O_s
+a^2\delta_b(a^{-1}E_s+O_s)\}.
\end{split}
\end{equation*}
The local diagonal remainder has size $P_s(a+\delta_b)E_s$.
The order-zero height source costs
$P_s(a+\delta_b+a/q)B_q^{-1/4}E_s$.
The normal good-unknown cancellation is applied to the DN operator
before its symbol is truncated. Retaining the order-zero symbol
$\ell_0$ in the upper-right block leaves a remainder of order at most
$-1$ from momentum to the kinematic equation. The total $Y$ derivative of
the good-coordinate matrix along the boundary evolution is included
as in \eqref{bv-bd-good-coefficient-slot}. Thus the
kinematic remainder uses the momentum norm $H^{s+1/2}$ in
\eqref{eq:bv-energy-norms}.

The other reference smoothing terms and negative-order compositions
are kept separately. The fixed cutoff bounds in
\cite[proof of Theorem~4.11, Steps~3--5]{LYIrrotational}, with
$B_q$ retained in every capillary composition, give their inner
natural source bound
\begin{equation*} 
P_s\bigl(a^{1/4}+\delta_bB_q^{-1/2}\bigr)(E_s+aO_s).
\end{equation*}
Here the reference negative-order terms cost $P_sa^{1/4}$;
the largest state-dependent graph conversion costs $B_q^{-1/2}$.
Consequently
\begin{equation*}
\delta_bB_q^{-1/2}\le P_ba^{1/4},\qquad
a\delta_bB_q^{-1/2}O_s\le P_ba^{5/4}O_s.
\end{equation*}

For the principal diagonal terms, define the reference Taylor coefficient by
\begin{equation*} 
\tau_c=1-\frac1{2\Lambda}.
\end{equation*}
Since $E_f=\pi L\Lambda q^2$,
\begin{equation*}
p_q+\frac1{2\Lambda}
=\frac{E_{f,Y}}{2E_f}-\frac{\Lambda_Y-1}{2\Lambda}-\frac{L_Y}{2L}.
\end{equation*}
Proposition~\ref{bv-scalars} now gives
\begin{equation}\label{bv-bd-reference-taylor-bounds}
\begin{split}
|1+p_q-\tau_c|&\le P_b(a+\delta_b),\\
|\tau_{c,Y}|&\le C(1+Y)^{-2}+P_b(a+\delta_b),\qquad
|p_q|\le\tfrac14,\quad p_q\le P_b(a+\delta_b).
\end{split}
\end{equation}
The nonintegrable negative part $-1/(2\Lambda)$ remains in the
principal diagonal coefficient. The spatial pressure
formula \eqref{bv-bd-potential-pressure}, the condenser expansion
and the raw normal-trace estimate
\eqref{bv-st-pressure-high-normal}, with Sobolev embedding at a
fixed index covered by $b$, give
\begin{equation}\label{bv-bd-actual-taylor-bound}
\left\|\frac{a^3(-\partial_np)}{q^2}-\tau_c\right\|_{C^J(\Gamma_-)}
\le P_J\left\{a+\delta_b+\frac aq(\varepsilon+Z_{J+c_0})\right\}.
\end{equation}
At the straight circular state with zero nonconstant traces,
the local kinetic term and the inner capacitary pressure datum
contribute $1-(2\Lambda)^{-1}+B_q/\Lambda$, where
$B_q/\Lambda=O_\sigma(a)$. The fixed outer boundary data and
toroidal corrections contribute $O(P_Ja)$, while nonconstant
harmonic data contribute $P_J(a+\delta_b)$.
The pressure correction $R_p$ contributes the last term in
\eqref{bv-bd-actual-taylor-bound}. The difference of the normalized Taylor coefficient from $\tau_c$
is an order-zero height source with conversion factor $B_q^{-1/4}$,
controlled by the spatial norms in \eqref{bv-bd-actual-taylor-bound}.

The circular potential principal system is that of
\cite[Theorem~3.16]{LYIrrotational}, with the instantaneous $p_q$
and the reference value $\tau_c$, by \eqref{bv-bd-scalar-rotation} and
\eqref{bv-bd-reference-taylor-bounds}. Its finite-rank clock and
centre columns have profiles $P_s(a^2+a|\beta|+|\beta|^2)$ and
$P_s(a+|\beta|)$. The $p_q$ terms $-p_qF$ and the
coefficient-derivative term proportional to $\eta$ both vanish at
$U=0$. For the scalar coefficient differentials in these columns,
the factor $a^{-1}$ is compensated by the inner row weight $a$,
leaving the bound $P_s(a+\delta_b)E_s$.
Outer inputs retain an additional factor $a$.
After their physical row factors, the affine forcing obeys
\begin{equation}\label{bv-bd-reference-forcing}
Q_s(f_{\rm ref},g_{\rm ref})
\le P_s(a^3+a^2|\beta|+a|\beta|^2)
\le P_s(a^{5/2}+a^3).
\end{equation}
The common Bernoulli constant is projected out before the
half-density change.

\medskip
\noindent\textbf{Step 2. The inner quadratic form.}
We incorporate the rotational height symbol $r_-$ into the inner
energy by setting
\begin{equation*} 
\begin{split}
\ell_0&=-\frac a2\partial_{n_-}\log\rho,\\
L_-&=\Theta P_{\ge2}\operatorname{Op}(x+\ell_0)P_{\ge2},\\
Q_-&=P_{\ge2}\operatorname{Op}(r_-)P_{\ge2},\\
M_-&=\Theta P_{\ge2}\operatorname{Op}\bigl(
\tau_c+B_q(x^2-1)\bigr)P_{\ge2}-Q_-.
\end{split}
\end{equation*}
Let $\mathcal T_-^0$ denote the common skew transport in the
potential reduction of
\eqref{bv-ex-prescribed-profile-system}--\eqref{bv-ex-q-equation}
by \eqref{eq:bv-good-variables}, as in
\cite[Theorem~3.16 and Section~7.2, Step~3]{LYIrrotational}. With the
fixed mode projections understood, the full inner rows are
\begin{equation}\label{bv-bd-collected-inner-system}
\begin{split}
h_Y&=\mathcal T_-^0h+h+L_-\pi+\mathcal R_h,\\
\pi_Y&=\mathcal T_-^0\pi-(1+p_q)\pi
      -M_-h-T_-\pi+\mathcal R_\pi.
\end{split}
\end{equation}
Here $\mathcal R_h$ is the potential kinematic remainder.
The momentum remainder contains the potential remainder with its
full reference affine forcing, the kinetic remainder in
\eqref{bv-bd-kinetic-remainder}, the pressure lift from
Lemma~\ref{bv-pressure}, and the coefficient term in
\eqref{bv-bd-good-coefficient-slot}. The kinetic remainder includes
the order-zero transport correction between the unreduced transport
and $T_-$. The remaining finite-rank terms come from moving
projections and trace gauges.
Substitution of $q_Y=q_Y^0+q\delta p_q$ places
the scalar momentum term in $-(1+p_q)\pi$ and its residual height
coefficient in $\mathcal R_\pi$. The height term $Q_-h$ is
incorporated in $-M_-h$.

To control the coefficients in \eqref{bv-bd-collected-inner-system},
the normalized harmonic estimates give
\begin{equation*}
\|j-1\|_{C^J}+|\Theta-1|\le P_J(a+\delta_b+a/q),
\qquad |\ell_0|_{J,0}\le P_Ja.
\end{equation*}
The $C^J$ differences between the good coordinate coefficients and
their straight circular values are bounded by $P_J(a+\delta_b+a/q)$.
At the fixed derivative orders included in $b$, the kinematic
equation and Proposition~\ref{bv-scalars} give
\begin{equation*}
\|j_Y\|_{C^J}+|\Theta_Y|\leq P_J(a+\delta_b+a/q),
\qquad |\ell_{0,Y}|_{J,0}\leq P_Ja.
\end{equation*}
For the second bound, differentiate
$\ell_0=-a\partial_{n_-}\log\rho/2$, using $a_Y=-a$, the fixed
distance from the axis and the geometric velocity bounds.
To estimate $b_{g,Y}$, first use \eqref{bv-bd-total-normal} and
Lemma~\ref{bv-pressure} to obtain
\begin{equation*}
\partial_Y b_-=-(1+p_q)b_-
+\Theta\,\frac{a^3(-\partial_np)}{q^2}
+R_n,\qquad
\|R_n\|_{C^J}\leq P_J(\delta_b+a/q).
\end{equation*}
The clock equation supplies the coarse time bound. Subtracting the
reference terms and using \eqref{bv-bd-reference-taylor-bounds} and
\eqref{bv-bd-actual-taylor-bound} then gives
\begin{equation*}
\|\partial_Yb_-\|_{C^J}\le P_J(a+\delta_b+a/q).
\end{equation*}
For the tangential potential component $u_-=(a/q)v_{0,\tau,-}$,
the static harmonic reconstruction gives
$\|u_-\|_{C^J}\leq P_J(a+\delta_b)$.
Differentiating its Neumann reconstruction once, with the flux,
normal-data and geometric derivatives just estimated, gives the
coarse bound $\|u_{-,Y}\|_{C^J}\leq P_J$.
Since $b_g=(d/j)b_--(\eta_\theta/j)u_-$ and the kinematic
equation gives $\|\eta_Y\|_{C^{J+1}}\leq P_J(a+\delta_b)$,
the product rule yields
$\|b_{g,Y}\|_{C^J}\leq P_J(a+\delta_b+a/q)$.
Together with the geometric bounds, this controls the $Y$ derivatives
of $\sqrt{w_g}$, $\sqrt{w_g}d/j$ and $\sqrt{w_g}b_g$ in
\eqref{eq:bv-good-variables} by the same quantity in $C^J$.

To estimate the time derivative of the rotational operator $Q_-$,
differentiate the instantaneous Hodge field as in
Lemma~\ref{bv-outer-time}. With chart $Y$ velocity $u_a$ and
inner normal velocity $F_-=u_a\cdot n_-$, the curl source is
$\dot\omega=\rho(w_Y\circ X^{-1}-u_a\cdot\nabla W)e_\phi$.
The normal and period data are
\eqref{bv-st-differentiated-normal}--\eqref{bv-st-period-derivative}.
The period is $O(P_ba^2\varepsilon)$, and
\begin{equation*}
g_-=(\rho_-aj)^{-1}
\partial_\theta(\rho_-F_-V_{\tau,-}).
\end{equation*}
For a fixed integer $m>J+1$, the static Hodge estimate uses this
datum in $H^{m-1/2}$, $V_{\tau,-}$ in $H^{m+1/2}$ and the
inner graph in $H^{m+2}$. The fixed reserve $c_0$ is chosen to
cover these levels.
The rescaled curl source satisfies
\begin{equation*}
a\|\dot\omega(c+a\cdot)\|_{H^{J+c_0-1}}
\le P_J\{aM_{J+c_0}(w_Y)+Z_{J+c_0}\},\qquad
\|\dot\omega\|_\infty\le P_J\varepsilon.
\end{equation*}
The div--curl reconstruction of Lemmas~\ref{bv-hodge} and
\ref{bv-outer-time}, followed by the pullback and tangent vector
variation, therefore gives
\begin{equation*}
\|V_{\tau,-,Y}\|_{C^J}
\le P_J\{1+M_{J+c_0}(w_Y)+Z_{J+c_0}\}
\le P_Ja^{-1/4}.
\end{equation*}
Here $J\leq J_*$, so the choices in
Subsection~\ref{subsec:bv-boundary-system} give
$J+c_0\leq k_{\rm t}\leq b-2$. Thus
\eqref{bv-bd-time-assumptions} controls the required norm of $w_Y$,
and \eqref{bv-bd-profile-assumptions} controls $Z_{J+c_0}$. Hence
\begin{equation}\label{bv-bd-inner-symbol-time}
\|\partial_Y v_-\|_{C^J}\leq P_Ja^{3/4}/q,\qquad
\frac{|r_-|_{J,1}}{B_q}\leq P_Ja/\sigma,\qquad
\frac{|\partial_Yr_-|_{J,1}}{B_q}
\leq P_J(a+a^{3/4})/\sigma.
\end{equation}

Let $D_s=P_{\ge2}\operatorname{Op}(x^s)P_{\ge2}$, with the same fixed
low-frequency completion as in the potential energy.  Define
\begin{equation*} 
E_s^2=\sum_{A\in\{I,D_s\}}
\left\{\langle M_-Ah,Ah\rangle
+\langle L_-A\pi,A\pi\rangle\right\}.
\end{equation*}
The real even symbol $r_-$ has small relative form bound by
\eqref{bv-bd-inner-symbol-time}, so this energy is equivalent to
\eqref{eq:bv-energy-norms}. We give the additional terms in the
calculation of \cite[Theorem~3.16]{LYIrrotational}. For each
$A\in\{I,D_s\}$, the principal cross terms cancel exactly:
\begin{equation*}
2\langle M_-Ah,L_-A\pi\rangle
-2\langle L_-A\pi,M_-Ah\rangle=0.
\end{equation*}
Commuting $A$ with the equations introduces $[A,Q_-]h$.
This operator has order $s$ when $A=D_s$, and the mixed
interpolation bound in \eqref{eq:bv-interpolation-costs} gives
\begin{equation*}
\bigl|\langle L_-D_s\pi,[D_s,Q_-]h\rangle\bigr|
\leq C_s|r_-|_{J,1}B_q^{-1/4}E_s^2.
\end{equation*}
The time form and the extra dilation form obey
\begin{equation*}
\bigl|\langle Q_{-,Y}Ah,Ah\rangle\bigr|
+\bigl|\langle Q_-Ah,Ah\rangle\bigr|
\leq C_s\frac{|\partial_Yr_-|_{J,1}+|r_-|_{J,1}}{B_q}E_s^2.
\end{equation*}
Here $B_q\leq1$ and
$\|h\|_{H^{s+1/2}}^2\leq C B_q^{-1/2}E_s^2$.
Since $T_-^*=-T_-$, the momentum form contains
$[L_-,T_-]$ of order one and $[D_s,T_-]$ of order $s$.
Their contributions, together with the order-zero transport
correction in the remainder, are bounded by
$P_s\|v_-\|_{C^J}E_s^2$.

For the diagonal terms, use
\begin{equation*}
\frac{B_{q,Y}}{B_q}=-1-2p_q,\qquad
(p_q)_+\le P_b(a+\delta_b).
\end{equation*}
The quadratic form calculation in
\cite[proof of Theorem~3.16, Step~3]{LYIrrotational}, together with
the preceding $Q_-$ bounds, gives
\begin{align*}
\sum_{A\in\{I,D_s\}}\langle(-M_{-,Y}-2M_-)Ah,Ah\rangle
&\le e_s E_s^2,\\
-L_{-,Y}+2(1+p_q)L_-&\le (2+e_s)L_-.
\end{align*}
The second inequality is understood as a quadratic form inequality.
The geometric coefficient bounds established in this step also
control $D_{s,Y}$ and the common transport.
By \eqref{bv-bd-inner-symbol-time} and
\eqref{bv-bd-reference-taylor-bounds}, all additional coefficients
are integrable. Using the source estimates
\eqref{bv-bd-primitive-norm}, \eqref{bv-bd-kinetic-remainder},
\eqref{bv-bd-coefficient-slot-bound}, \eqref{bv-st-pressure-inner-result}
and \eqref{bv-bd-potential-high-high}--\eqref{bv-bd-reference-forcing},
we obtain
\begin{equation*}
-\frac12\frac{dE_s^2}{dY}
\le (1+e_s)E_s^2
+P_s\bigl\{a^2(\varepsilon+Z_{s-1})+a^{5/4}O_s
+a^{5/2}+a^3\bigr\}E_s.
\end{equation*}
Division by $E_s$, with regularization at zero, proves
\eqref{bv-bd-inner-inequality}.

\medskip
\noindent\textbf{Step 3. The outer coefficients and energy.}
To control the outer Taylor coefficient, split the inner datum
$m_{\mathrm p,-}=\sigma\mathcal H_-+|v_0|^2/2$ in
\eqref{bv-bd-potential-pressure} into its capacity component and
a lift with zero inner flux. Its coefficient satisfies
\begin{equation*}
c_m=\operatorname{Cap}^{-1}\int_{\Gamma_-}m_{\mathrm p,-}\partial_n u_c\,dS,
\quad |c_m|\le C_*(q^2/a^2+\sigma/a)
+P_b(q^2/a^2+\sigma/a)(a+\delta_b).
\end{equation*}
The positive condenser source has mass $\operatorname{Cap}=2\pi L/\Lambda$.
The filled domain representation \eqref{eq:bv-condenser-representation},
with source separated from the outer boundary, therefore gives
\begin{equation}\label{bv-bd-capacity-outer-derivative}
\|\partial_n u_c\|_{C^J(\Gamma_+)}\le C_J/\Lambda,
\end{equation}
where $C_J$ is independent of $Y$.

The remaining inner datum has fixed low nonconstant norm
at most
$P_J\{q^2(a+\delta_b)/a^2+\sigma\delta_b/a+\sigma\}$.
For this nonconstant datum, transfer from the inner to the outer
boundary with zero flux gives an additional factor $a$.
The outer data and the local term $-|v_0|^2/2$ have fixed
physical scale.  Thus the instantaneous outer Taylor coefficient obeys
\begin{equation}\label{bv-bd-outer-taylor}
\|t_Y\partial_{n_+}p_0\|_{C^J}
\leq C_J\frac q\Lambda+P_J(a+\delta_b)
\leq C_J(1+Y)^{-3/2}+P_J(a+\delta_b).
\end{equation}
The outer bound in \eqref{bv-st-pressure-low-traces} and its
fixed-order outer extension stated immediately afterwards give
\begin{equation*}
t_Y\|\partial_nR_p\|_{C^J(\Gamma_+)}
\le P_J\frac{a^2}{q}
\left(\frac{q\varepsilon}{a}+\varepsilon^2\right)
\le P_Ja\varepsilon.
\end{equation*}
Thus the full Taylor coefficient satisfies
\eqref{bv-bd-outer-taylor} with an additional $P_Ja\varepsilon$.
We estimate its multiplication operator as a source of order zero.

The same cancellation
\eqref{bv-bd-mixed-cancellation}, now in physical outer
coordinates, gives
\begin{equation}\label{bv-bd-outer-principal}
\Delta_{\rm kin}h_{+,Y}=0,\qquad
\Delta_{\rm kin}\pi_{+,Y}
=t_Y\left\{
-2i\operatorname{Op}\left(\frac{V_{\tau,+}}{j_+}\xi\right)\pi_+
+\operatorname{Op}(V_{\tau,+}^2x_+)h_+\right\}+R_{+,\rm kin},
\qquad x_+=\frac{\varrho(\xi)}{j_+}.
\end{equation}
For the projected inner datum $q^2a^{-3}k_-$, let $\Psi_-$
be the Neumann solution with zero outer datum. Since
$\|k_-\|_{L^2(d\theta)}\le P_ba/q$, the scaled trace--Poincar\'e
bound gives $\|\nabla\Psi_-\|_2\le P_bq/a$.
Fixed outer collar regularity then yields
\begin{equation}\label{bv-bd-outer-inner-data}
t_Y\|P_+^0\Psi_-\|_{H^{s+1/2}}
\le P_sa(1+\|\zeta_+\|_{H^{s+1}}).
\end{equation}
The inner datum enters through its fixed low norm.

For outer Neumann data, return to physical velocity units and use
the independent primitive
\begin{equation*}
\widetilde V_+=\langle D_\theta\rangle_+^{-1}V_{\tau,+},
\qquad
\|\widetilde V_+\|_{H^{s+1/2}}
\leq P_s\varepsilon(1+O_s).
\end{equation*}
After extracting $t_Y$, apply \cite[Lemma~4.4]{LYIrrotational}
as in \eqref{bv-bd-kinetic-remainder}. The outer state derivatives
have order at most one and the $\widetilde V_+$ derivative has
order zero. Inner inputs enter through separated kernels.
Equation~\eqref{bv-bd-neumann-hadamard} gives the shape bounds
with fixed losses and polynomial constants. Hence
\begin{equation}\label{bv-bd-outer-local-remainder}
\|R_{+,\rm local}\|_{H^{s+1/2}}
\le P_st_Y\{\varepsilon+(1+\varepsilon+\delta_b)O_s\}.
\end{equation}
Together with \eqref{bv-bd-outer-inner-data} and the potential
kernel bounds, this controls the combined outer remainder by
$e_sO_s+P_s(a+E_s)$.

The time derivative of the outer good trace coefficient is
computed from
\begin{equation*}
b_{g,+}=(n_+^*\cdot n_+)B_+
+(n_+^*\cdot\tau_+)\partial_{\Sigma_+}g.
\end{equation*}
The outer chart velocity bound gives $C^J$ bounds $P_Ja^2$ for
$\partial_Y(n_+^*\cdot n_+)$, $\partial_Y(n_+^*\cdot\tau_+)$ and
$\partial_Y(j_+^{-1})$.
The total normal equation and \eqref{bv-bd-outer-taylor} give
$\|B_{+,Y}\|_{C^J}\le C_Jq/\Lambda+P_J(a+\delta_b)$.
The potential Bernoulli row contributes $O(P_Jt_Y)$ to $g_Y$.
The spatial rotational row bounds
\eqref{bv-st-pressure-outer-result}, \eqref{bv-bd-outer-inner-data} and
\eqref{bv-bd-outer-local-remainder} contribute $P_Ja$ at a fixed
low index. Consequently
\begin{equation}\label{bv-bd-outer-good-time}
\|b_{g,+,Y}\|_{C^J}
\le C_J(1+Y)^{-3/2}+P_J(a+\delta_b).
\end{equation}
The constant $C_J$ comes from
\eqref{bv-bd-capacity-outer-derivative} and is uniform in $Y$.

The order one height term in
\eqref{bv-bd-outer-principal} is retained in the energy.
Let $L_+$ denote the retained order-one and order-zero part of the
outer half-density DN block, with a fixed self-adjoint low-mode
completion,
\begin{equation*}
L_+=\operatorname{Op}(x_++\ell_{+,0})+C_{\rm low}^{L}.
\end{equation*}
Here $\ell_{+,0}$ is the order-zero symbol of the outer DN block.
The remaining part of order at most $-1$ is included in the
kinematic remainder. Set
\begin{equation*} 
M_+=\sigma\operatorname{Op}(x_+^2)
-\operatorname{Op}(V_{\tau,+}^2x_+)+C_{\rm low}.
\end{equation*}
Choose the fixed self-adjoint finite-rank completions
$C_{\rm low}^{L}$ and $C_{\rm low}$ so that both forms are
coercive on the stated trace gauge and graph spaces at the reference
geometry. The same holds throughout the fixed lower neighbourhood.
Their additional low mode sources carry $t_Y$ and are integrable.
The outer collar is free of
vorticity, and Lemma~\ref{bv-hodge} gives
$\|V_{\tau,+}\|_{C^J}\le C_J\varepsilon$ uniformly in $a$.
Thus $M_+$ is coercive for the stated small $\varepsilon$.

Let $D_s^+=\operatorname{Op}(x_+^s)$, with the
fixed gauge and low mode corrections, and define
\begin{equation}\label{bv-bd-outer-energy-definition}
O_s^2
=\sum_{A\in\{I,D_s^+\}}
\left\{\langle M_+Ah_+,Ah_+\rangle
+\langle L_+A\pi_+,A\pi_+\rangle\right\}.
\end{equation}
The potential norm equivalence and the small rotational term in
$M_+$ give the outer equivalence in \eqref{eq:bv-energy-norms}.
The Hamiltonian cross terms cancel as in Step~2, with the common
factor $t_Y$.
The identities $\{x_+^s,x_+\}=\{x_+^s,x_+^2\}=0$ control the
capillary commutators as in
\cite[proof of Proposition~4.12, Step~4]{LYIrrotational}.
The additional transport and height-symbol commutators, including
the order-zero skew-completion term, obey
\begin{equation*}
|\mathcal C_{{\rm rot},+}|
\leq C_s t_Y\bigl(\|V_{\tau,+}\|_{C^J}
+\|V_{\tau,+}^2\|_{C^J}\bigr)O_s^2
\leq P_st_Y\varepsilon O_s^2.
\end{equation*}
Here $[D_s^+,\operatorname{Op}(V_{\tau,+}^2x_+)]$ has order
$s$, so the fixed positive surface tension controls its height
input in $H^{s+1/2}$. These terms therefore have an integrable
coefficient without a higher profile norm.
For the additional time form, Lemma~\ref{bv-outer-time}, together
with the outer-chart velocity bound and the geometric dependence
of the DN symbols, gives
\begin{equation}\label{bv-bd-outer-symbol-time}
|\partial_Y(V_{\tau,+}^2x_+)|_{J,1}
\le P_J\varepsilon^2(a^{3/4}+a),\qquad
\|j_{+,Y}\|_{C^J}+|\ell_{+,0,Y}|_{J,0}\le P_Ja^2.
\end{equation}
Together with \eqref{bv-bd-outer-taylor},
\eqref{bv-bd-outer-good-time}, \eqref{bv-st-pressure-outer-result},
\eqref{bv-bd-outer-inner-data} and
\eqref{bv-bd-outer-local-remainder}, this yields
\begin{equation*}
\frac12\left|\frac{dO_s^2}{dY}\right|
\le e_sO_s^2+P_s(a+E_s)O_s,
\end{equation*}
and hence \eqref{bv-bd-outer-inequality}.

\medskip
\noindent\textbf{Step 4. First modes and completion.}
It remains to control $\beta$. At the straight circular reference, the centred kinetic expression
has first harmonic $\beta$, since
\begin{equation*}
\left[\frac12(1+\beta\cdot e_r)^2
-\frac32(\beta\cdot e_\theta)^2\right]_1
=\beta.
\end{equation*}
The two component calculation in
\cite[Section~7.2, Step~2]{LYIrrotational}, with the present harmonic
and curvature bounds on \eqref{bv-bd-admissible-ball}, gives
$|(q\beta)_Y-q\beta|\le P_bq(a+\|U\|_{H^b})$
for the potential first harmonic projection. Let $\mathcal R_{B,\rm rot}$
denote the right-hand side of \eqref{bv-bd-rotation-normal}.
We extract the first harmonic before the $P_{\ge2}$ projection
and the half-density change. Since the full and instantaneous
potential evolutions have the same geometric velocity,
\begin{equation*}
\Delta[(q\beta)_Y]
=t_Y\bigl[(N\Pi_n\mathcal R_{B,\rm rot})_-\bigr]_1.
\end{equation*}
The component constants and the condenser trace have zero first
harmonic. At a fixed low Sobolev index, the normalized diagonal
inverse $a^{-1}N_{--}$ is bounded and $N_{-+}$ gains a factor $a$.
The kinetic data
\eqref{bv-bd-inner-kinetic-data}--\eqref{bv-bd-outer-kinetic-data}
and the pressure estimate \eqref{bv-st-pressure-inner-result}
therefore give
\begin{equation*}
\frac{t_Y}{q}
\left|\bigl[(N\Pi_n\mathcal R_{B,\rm rot})_-\bigr]_1\right|
\leq P_b\left\{\|U\|_{H^b}
+\frac aq(\varepsilon+Z_{b-1})\right\}.
\end{equation*}
Here $b$ covers the fixed graph and profile orders used in these
estimates. Combining this bound with the potential first-mode
estimate gives
\begin{equation*}
\left|\frac{(q\beta)_Y-q\beta}{q}\right|
\le P_b\left\{a+\|U\|_{H^b}
+\frac aq(\varepsilon+Z_{b-1})\right\}
\le P_b(a+\delta_b).
\end{equation*}
This proves \eqref{bv-bd-first-mode-equation}. For zero terminal data its
variation of constants formula is
\begin{equation*}
\beta(Y)=-\int_Y^{Y_N}
\exp\left(-\int_Y^y(1-p_q(z))\,dz\right)r_\beta(y)\,dy.
\end{equation*}
Since $1-p_q\geq3/4$ and $a(y)=a(Y)e^{-(y-Y)}$,
\begin{equation*}
|\beta(Y)|
\le\int_Y^{Y_N}e^{-3(y-Y)/4}P_b(y)
\bigl(a(y)+a(y)^{3/4}\bigr)\,dy
\le P_b(Y)a(Y)^{3/4}.
\end{equation*}
The bound is uniform in $Y_N$, and the clock gives
$|c_Y|\le P_ba(a+\delta_b)$.

Finally, \eqref{bv-bd-inner-symbol-time},
\eqref{bv-bd-reference-taylor-bounds} and
\eqref{bv-bd-capacity-outer-derivative} give
\begin{equation*}
e_s(Y)\le C_s\left(\frac q\Lambda+(1+Y)^{-2}\right)
+P_sa^{1/4}
\le C_s(1+Y)^{-3/2}+P_sa^{1/4},
\end{equation*}
which is \eqref{bv-bd-integrable-coefficient}.
\end{proof}

\begin{corollary}\label{cor:bv-terminal-bounds}
Assume the hypotheses of Lemma~\ref{bv-boundary-energy} on $[Y_0,Y_N]$
at an integer level $s\geq b$, with polynomially bounded $Z_{s-1}$ and
terminal data $(\eta,F,\zeta_+,g)(Y_N)=0$, $\beta(Y_N)=0$. Then 
\begin{equation*}
E_s\leq P_sa^2,\qquad O_s\leq P_sa,\qquad |\beta|\leq P_ba,
\end{equation*} 
uniformly in $Y_N$. Once the lower bounding polynomials are fixed,
$Y_0$ can be chosen independently of the higher index $s$ so that
these estimates strictly improve \eqref{bv-bd-weak-bootstrap}.
\end{corollary}

\begin{proof}
Set $I_s=E_s+aO_s$.
Since $a_Y=-a$, Lemma~\ref{bv-boundary-energy} gives
\begin{equation}\label{bv-bd-combined-tail}
-I_s'
\leq\{1+e_s+P_s(a^{1/4}+a)\}I_s
+P_sa^2(1+\varepsilon+Z_{s-1}).
\end{equation}
For zero terminal high boundary data and polynomially bounded
$Z_{s-1}$, the backward integration argument in
\cite[proof of Lemma~3.18]{LYIrrotational}, now with homogeneous
rate $1$ and forcing exponent $2$, gives
\begin{equation*}
I_s(Y)
\le C_s\int_Y^{Y_N}e^{y-Y}P_s(y)a(y)^2\,dy
\le P_s(Y)a(Y)^2.
\end{equation*}
The integrable coefficient in \eqref{bv-bd-combined-tail} is absorbed
in $C_s$, uniformly in $Y_N$. Since $E_s,O_s\ge0$,
\begin{equation*}
E_s\le P_sa^2,\qquad O_s\le P_sa,\qquad
\|U\|_{H^b}\le P_ba.
\end{equation*}
In Step~4 of the proof of Lemma~\ref{bv-boundary-energy}, the
improved bound on $\|U\|_{H^b}$ gives $|r_\beta|\le P_ba$.
The variation of constants formula for $\beta$ then yields
$|\beta|\le P_ba$, closing \eqref{bv-bd-weak-bootstrap}.

To make the choice of endpoint explicit, fix the lower bounding
polynomials and bound the $b$-level integrating factors on one
preliminary late interval. These bounds remain valid when its
left endpoint is increased. Let $\widehat P_b$ dominate the
three resulting improved bounds and the fixed norm equivalence
constants, and let $P_b^{\rm boot}$ be the polynomial fixed in
\eqref{bv-bd-weak-bootstrap}. Increase $Y_0$ so that
\begin{equation*}
\widehat P_b(Y)a(Y)^{1/4}
\leq \tfrac12 P_b^{\rm boot}(Y),\qquad Y\geq Y_0.
\end{equation*}
This strictly improves all three weak bounds. At higher indices
only finiteness of the corresponding integrating factor is used;
no smallness condition is imposed simultaneously for all $s$.
\end{proof}

\section{Construction and terminal limits}
\label{sec:bv-construction}

The local problem for a prescribed profile defines the boundary
map needed to couple backward boundary evolution with forward
transport. The resulting solutions have uniform bounds as the
terminal cavity radius tends to zero. The analysis of their limiting
vorticity, geometry and energy completes the proof of
Theorem~\ref{thm:bv-collapse}.

\subsection{The prescribed-profile local problem}
\label{subsec:bv-local-profile}

We consider the boundary system
\eqref{bv-ex-prescribed-profile-system}--\eqref{bv-ex-q-equation}
with a prescribed profile path $w(Y)$, independently of the transport
equation. The resulting boundary solution supplies the backward evolution in
Subsection~\ref{subsec:bv-coupled-construction}. We use uniqueness
for fixed input and continuous dependence on that input to prove
continuity of the fixed point map.

The pressure is reconstructed instantaneously: its zero Dirichlet
part solves 
\begin{equation*}
\Delta p_D=-\operatorname{tr}((\nabla v)^2),
\end{equation*} 
and its
harmonic part has boundary value $\sigma\mathcal H$.
The normal equations, including the first modes and the scalar
equation \eqref{bv-ex-q-equation}, are evolved together with the
kinematic equations.

For the local problem, use the state and Sobolev scale at a fixed radius
\begin{equation*}
\Upsilon=(\eta,F,\zeta_+,g,\beta,c,q),\qquad
X_s=(P_{\ge2}H^{s+1}\times P_{\ge2}H^{s+1/2})_-
\times(H^{s+1}\times (I-P_0)H^{s+1/2})_+\times\R^5.
\end{equation*}
The profile Sobolev spaces below are on the fixed coordinate collar
$[0,s_3]\times\T$. An admissible state belongs to the geometric
class of Lemma~\ref{bv-hodge}, satisfies 
\begin{equation*}
P_0\eta=P_1\eta=0,\qquad P_{\ge2}F=F,\qquad P_0g=0,
\end{equation*} 
and
has $q>0$ and $\Theta>0$. We restrict the local problem to the open
weak neighbourhood used in Lemma~\ref{bv-boundary-energy}.
In particular, require
$\|U\|_{H^{b_{\rm ad}}}<\delta_b(Y)$, so that all frozen states
and Taylor segments lie strictly in \eqref{bv-bd-admissible-ball}.
The profile bounds at the fixed lower orders and the coercivity of
both boundary forms hold with a margin on this neighbourhood.

Fix $Y_N<\infty$, so 
\begin{equation*}
a\geq a_N:=a(Y_N)>0.
\end{equation*} 
First choose the comparison index $r$ above the fixed lower thresholds
for the differentiated boundary maps and the centred-coordinate
comparison, including the coefficient differentials, with
$r\geq b_{\rm dy}$.

Next fix the symbol, kernel and coefficient seminorms needed through
index $r+3/2$. Let $d_*$ bound the orders of the boundary columns
of the state differentials of the good-variable coefficients, and let
$d_c$ cover the coefficient losses, component Sobolev shifts and
Bernstein factors in the centred-coordinate comparison, including
one further state differential.
With $\nu=r+1/2$, choose the decay order $L_{\rm dec}$ to satisfy
$L_{\rm dec}>\nu+6$ in the Schur sums for the complete remainders and
$L_{\rm dec}>A_*+d_*+J+1$ for each of the finitely many required pairs
$(A_*,J)$. Here $J$ is the output Sobolev order and $A_*$ is the
decay exponent in the separated-frequency coefficient estimates.
These pairs cover the estimates through index $r+3/2$, with $A_*$
larger than the Sobolev and Bernstein exponents in the corresponding
dyadic sums.

For each required order $J$, let $K_{\rm reg}(J)$ denote the spatial
order of $w_Y$ needed to control the $C^J$ norm of the
$Y$ derivatives of the approximate principal coefficients and the
good-variable coefficients. Choose $K_{\rm reg}$ to cover these
orders and the spatial orders of $w_Y$ needed in
Lemmas~\ref{bv-outer-time} and~\ref{bv-boundary-energy}.
Then choose the background index $M\geq K_{\rm reg}$, with
$M>r+d_c+7$, to cover the elliptic differentiations of order at most
four, the $L_{\rm dec}$ background derivatives, the finite symbol and
kernel losses, and the background order $\nu+3$ in the complete
remainder estimates. Require the coefficient Fourier matrix bounds
to hold also with background directions in $X_{M-2}$, as needed when
such a direction is evaluated on the vector field.
Require also one spatial derivative beyond the orders defining
admissibility, and set
\begin{equation*}
s=M+1>r.
\end{equation*}
The thresholds are independent of the smoothing parameter and
subsequent higher energy levels. The comparison constants may
depend on $a_N^{-1}$ and the background bounds at order $M$.

Prescribe an admissible state in $X_s$ at $Y_1\leq Y_N$ and a
profile path in $C_YH^{s-1}\cap W_Y^{1,\infty}H^{K_{\rm reg}}$, with the fixed
collar support and the bounds on $w_Y$ in \eqref{bv-bd-time-assumptions}.
Time derivatives and evolution equations for these paths are
understood almost everywhere in $Y$.

\begin{lemma}\label{bv-local}
The boundary problem
\eqref{bv-ex-prescribed-profile-system}--\eqref{bv-ex-q-equation}
has, for the prescribed profile path, a unique local solution in
either direction of $Y$, continuous with values in $X_s$. Smooth spatial data and
profiles give a spatially smooth solution. The solution continues while it remains in this admissible
neighbourhood with uniform margins and the $X_s$ state bound, the
$H^{s-1}$ profile bound and the prescribed bounds on $w_Y$ persist.

Within the prescribed profile class, the solution map for fixed
data at $Y_1$ is continuous from $C_YH^r$ for the input profile to
$C_YX_r$ for the state on every common compact interval of existence,
provided the states are uniformly bounded in $L_Y^\infty X_s$, the
profiles in $L_Y^\infty H^{s-1}$ and their time derivatives in
$L_Y^\infty H^{K_{\rm reg}}$, with uniform admissibility margins.
\end{lemma}

\begin{proof}
We construct a regularized local solution using the frozen inverse.
Differentiating the full boundary maps then yields a comparison
estimate in centred coordinates, which gives uniqueness for fixed
input. Strong continuity and continuation complete the local theory.

\medskip
\noindent\textbf{Step 1. The frozen inverse and the source terms.}
Work on a fixed bounded admissible neighbourhood with uniform
nondegeneracy. Constants may depend on this neighbourhood and $a_N^{-1}$.
Let $\Xi=S_g(\Upsilon)\Upsilon$ denote the triangular change formed
from \eqref{eq:bv-good-variables} and its outer counterpart, with
the identity on finite coordinates. At a frozen state, on these
projected spaces,
\cite[Lemma~4.5]{LYIrrotational} and the diagonal block estimates give
\begin{equation*}
\|S_g(\Upsilon)^{-1}Z\|_{X_s}\leq C_{N,s}\|Z\|_{X_s}.
\end{equation*}
The lower triangular entry maps $H^{s+1}$ into $H^{s+1/2}$.
We regularize the boundary equations using this frozen inverse and the
adapted mollifiers of \cite[Sections~5--6]{ABZ2011}.

By Lemma~\ref{bv-hodge} and the trace theorem,
\begin{equation*}
w\in H^{s-1}
\ \Longrightarrow\ V\in H^s
\ \Longrightarrow\ V_\tau|_\Sigma\in H^{s-1/2}
\ \Longrightarrow\ \widetilde v_\pm\in H^{s+1/2},
\end{equation*}
where $X_s$ supplies the high graph norms. For the pressure,
\begin{equation*} 2\operatorname{tr}(\nabla v_0\nabla V)
 +\operatorname{tr}((\nabla V)^2)\in H^{s-1},\qquad R_p\in H^{s+1},\qquad
 \partial_nR_p\in H^{s-1/2}. 
\end{equation*}
After the component flux projection
\eqref{bv-st-capacitary-projection}, the harmonic Neumann lift of
$\Pi_n(\partial_nR_p)$ has trace in $H^{s+1/2}$.
The pressure lift therefore lies in the momentum source space.
For the kinetic contribution, the order-one terms are retained in
the principal operator and \eqref{bv-bd-kinetic-remainder} bounds
the remainder in the same source space.
The tame products place the highest spatial norm on at most one
factor, with all other profile factors controlled at the fixed
lower orders. Thus the state estimate in $X_s$ requires spatial
profile regularity only through $H^{s-1}$.

\medskip
\noindent\textbf{Step 2. Regularization and the local limit.}
The regularization must preserve the pairing of the self-adjoint
principal blocks in Lemma~\ref{bv-boundary-energy}. On the inner
boundary, suppressing the fixed mode projections, these blocks have the form
\begin{equation*}
L_g=\operatorname{Op}\bigl(\Theta(x+\ell_0)\bigr),\qquad
M_g=\operatorname{Op}(\Theta B_qx^2+m_1+m_0),\qquad
x=\varrho(\xi)/j.
\end{equation*}
Here $\Theta$ and $B_q$ are scalars depending only on $Y$,
and the real symbols $m_j$ have order $j$, $j=0,1$.
On the outer boundary set $L_g=t_YL_+$ and
$M_g=t_YM_+$, where $L_+,M_+$ are the operators defining
\eqref{bv-bd-outer-energy-definition}, and take $x=x_+$.
In either case retain the rotational height symbol of order one in $M_g$, take the transport in skew-adjoint form,
and add a fixed positive completion of the low modes. Each added
term in a principal row is cancelled by its opposite in the
remainder, preserving the exact equation.

\smallskip
\noindent\emph{(a) Adapted mollifiers and commutators.}
For $J_\zeta=\operatorname{Op}(e^{-\zeta x})$,
$0<\zeta\leq1$, use $J_\zeta L_g J_\zeta$ and
$J_\zeta M_g J_\zeta$ in the two principal rows.  The identities
$\{e^{-\zeta x},x\}=\{e^{-\zeta x},x^2\}=0$ give
\begin{equation}\label{bv-ex-adapted-commutators}
[J_\zeta,M_g]\in\operatorname{OPS}^0,\qquad
[J_\zeta,L_g]\in\operatorname{OPS}^{-1}
\end{equation}
uniformly in $\zeta$. These orders give the dual pairings
$H^{1/2}\times H^{-1/2}$ and $H^1\times H^{-1}$ in the
basic energy. With 
\begin{equation*}
A_s=\operatorname{Op}((1+x^2)^{s/2}),
\end{equation*}
the symbolic calculus in \cite[Sections~5--6]{ABZ2011} gives
\begin{equation*}
[A_s,M_g]\in\operatorname{OPS}^{s},\qquad
[A_s,L_g]\in\operatorname{OPS}^{s-1},\qquad
[A_s,J_\zeta]\in\operatorname{OPS}^{s-2}.
\end{equation*}
The skew transport admits the same $X_s$ estimate after composition
with $J_\zeta$ on both sides.
The same regularization requires bounds for the parameter derivatives. Write
$x=\varrho/j_\Sigma$, with $j_\Sigma=j$ or $j_+$, and let $\delta$ be a state
differential or a derivative along a time-dependent state. Since
$\delta x=-x\,\delta j_\Sigma/j_\Sigma$,
\begin{equation*}
\delta(e^{-\zeta x})=-\zeta e^{-\zeta x}\delta x,\qquad
\delta(1+x^2)^{s/2}
=-s(1+x^2)^{s/2}\frac{x^2}{1+x^2}\frac{\delta j_\Sigma}{j_\Sigma}.
\end{equation*}
The uniform lower bound for $j_\Sigma$ and its fixed lower derivative
bounds give, at high frequency,
\begin{equation*}
|\partial_\theta^\alpha\partial_\xi^\beta
\delta(e^{-\zeta x})|
\leq C_{\alpha,\beta}\|\delta j_\Sigma/j_\Sigma\|_{C^{|\alpha|}}
\langle\xi\rangle^{-|\beta|}
(1+\zeta\langle\xi\rangle)^{N_{\alpha,\beta}}
e^{-c\zeta\langle\xi\rangle}.
\end{equation*}
Indeed, each $\xi$ derivative of $x$ lowers its order by one,
whereas each differentiated exponential contributes a factor
$\zeta x$; the bounded frequency region is estimated directly.
Thus, for each of the finitely many symbol seminorms used here,
\begin{equation}\label{bv-ex-mollifier-symbol-bounds}
|e^{-\zeta x}|_{J,0}\leq C_J,\qquad
|\delta(e^{-\zeta x})|_{J,0}
+|\delta(1+x^2)^{s/2}|_{J,s}
\leq C_{J,s}\|\delta j_\Sigma/j_\Sigma\|_{C^J},
\end{equation}
uniformly for $0<\zeta\leq1$. These are also operator bounds
after the fixed cutoff quantization. In particular, $J_\zeta$
is smoothing for fixed $\zeta$, is uniformly bounded on every
fixed Sobolev space, and converges strongly to the identity.
The bound for $F_{g,\zeta}$ below controls the relative time
variation $\dot j_\Sigma/j_\Sigma$ in these fixed lower norms.

\smallskip
\noindent\emph{(b) The regularized vector field.}
To implement this regularization in the state $\Upsilon$, let $F_g$
denote the vector
field of \eqref{bv-ex-prescribed-profile-system}--\eqref{bv-ex-q-equation}
and write the instantaneous spatial identity as
\begin{equation*}
S_g F_g=A_g\Xi+R_g.
\end{equation*}
Here $A_g$ collects the principal rows above, $A_{g,\zeta}$ denotes
their regularized version, and $R_g$ is the spatial remainder.
To relate this convention to the exact good-variable equations,
let $R_{\rm good}$ denote their remainder with the same principal
operator $A_g$. Since $S_g=S_g(Y,\Upsilon)$ is independent of $w$
at fixed state, the chain rule gives
\begin{equation*}
R_g=R_{\rm good}
-\bigl(\partial_YS_g+D_\Upsilon S_g[F_g]\bigr)\Upsilon,
\end{equation*}
where $\partial_YS_g$ holds $\Upsilon$ fixed.
The boundary equations bound $F_g$ in $X_{s-3/2}$. The regularity
chosen before Lemma~\ref{bv-local} controls the coefficient seminorms
of $\partial_YS_g+D_\Upsilon S_g[F_g]$, so the operator with the
fixed cutoff is bounded on $X_s$. Hence the subtracted term lies
in $X_s$. Its estimate depends on the instantaneous profile through
$H^{s-1}$ and does not involve $w_Y$.
Along the regularized evolution, the coordinate chain rule term
is evaluated on $F_{g,\zeta}$ and added explicitly in
\eqref{bv-ex-actual-approximate-good-equation}.
On the boundary components take
$K_\zeta=e^{-\epsilon_0\zeta\langle D_\theta\rangle}$, with
$0<\epsilon_0\leq1$ fixed, and take the identity on finite coordinates.
This operator commutes with the fixed projections, is a contraction
on each Sobolev space, and is independent of the state and time.
For each fixed Sobolev index $\kappa$ and $0<\nu\leq1$,
\begin{equation*}
\|(K_\zeta-I)Z\|_{X_\kappa}
\leq C_\nu(\epsilon_0\zeta)^\nu\|Z\|_{X_{\kappa+\nu}},
\qquad
\|K_\zeta Z\|_{X_{\kappa+d}}
\leq C_d(\epsilon_0\zeta)^{-d}\|Z\|_{X_\kappa},\quad d\geq0.
\end{equation*}
Choose $\epsilon_0$ sufficiently small, using the uniform
admissibility margins and the $X_s$ bound. The one-derivative
reserve ensures that $K_\zeta\Upsilon$ is uniformly close to
$\Upsilon$ in the norms defining admissibility, so it remains
admissible for
$0<\zeta\leq1$. Set
\begin{equation}\label{bv-ex-raw-regularization}
F_{g,\zeta}=S_g(\Upsilon)^{-1}
\left\{A_{g,\zeta}(\Upsilon)\Xi
+K_\zeta R_g(K_\zeta\Upsilon,w)\right\}.
\end{equation}
For fixed $\zeta>0$, the principal part and the frozen inverse are
locally Lipschitz on the $X_s$ ball by the coefficient bounds and
the smoothing of $J_\zeta$. In the remainder, compare two smoothed
states $K_\zeta\Upsilon_i$ first in $X_{r-2}$. Differentiation of
the pulled-back elliptic equations along their connecting segment
places at most one derivative on the prescribed profile through the
chart variation $\xi\cdot\nabla W$. Its $H^{s-2}$ norm is controlled
by $\|w\|_{H^{s-1}}$ and the smoothed chart norms. In this lower
estimate the raw equations suffice because the height and momentum
rows cost at most $3/2$ state derivatives and the states are smoothed.
The product and trace estimates thus bound the difference of the
two remainders in $X_{r-2}$ by
$C_{N,s,\zeta}\|\Upsilon_1-\Upsilon_2\|_{X_s}$.
The outer 
\begin{equation*}
K_\zeta:X_{r-2}\to X_s
\end{equation*} 
recovers the top norm.
This proves local Lipschitz continuity for fixed $\zeta$.
Picard's theorem gives an approximate solution.
All rows retain their original projections, so the centering and
mean constraints are preserved. We obtain a common interval of
existence from the energy estimate below.

\smallskip
\noindent\emph{(c) Uniform local energy estimates.}
Differentiating $\Xi=S_g(\Upsilon)\Upsilon$ along the regularized
evolution gives
\begin{equation}\label{bv-ex-actual-approximate-good-equation}
\Xi_Y=A_{g,\zeta}\Xi+
K_\zeta R_g(K_\zeta\Upsilon,w)
+\bigl(\partial_Y S_g
+D_{\Upsilon}S_g[F_{g,\zeta}]\bigr)\Upsilon.
\end{equation}
The partial derivative $\partial_Y S_g$ is taken at fixed
$\Upsilon$. The operator
$\partial_Y S_g+D_\Upsilon S_g[F_{g,\zeta}]$
is lower triangular in each height--momentum block, with zero
upper-right entry.
Let $a_g^\zeta$ denote any coefficient of the approximate
principal operators or of the good variable change. Its derivative
along the approximate path is
\begin{equation*}
\frac{d}{dY}a_g^\zeta(Y,\Upsilon^\zeta,w)
=\partial_Ya_g^\zeta+D_\Upsilon a_g^\zeta[F_{g,\zeta}]
+D_wa_g^\zeta[w_Y],
\end{equation*}
where the partial derivative holds the state and profile fixed;
the last term is absent for coefficients independent of $w$.
Differentiating the instantaneous elliptic equations, the scalar
Schur normalization of the harmonic traces, and the inner
circulation constraint gives, for every fixed nonnegative integer $J$,
a finite derivative index $K_{\rm reg}(J)$ such that
\begin{equation}\label{bv-ex-approximate-time-coefficients}
\left\|\frac{d}{dY}a_g^\zeta\right\|_{C^J}
\leq C_{N,s}\bigl(1+\|w_Y\|_{H^{K_{\rm reg}(J)}}\bigr),
\end{equation}
uniformly in $\zeta$, for $s$ above the corresponding spatial
threshold.
The height row has order one on the potential trace and the
momentum row has order two on the height. Hence
\begin{equation*}
\sup_{0<\zeta\leq1}\|F_{g,\zeta}\|_{X_{s-3/2}}\leq C_{N,s}.
\end{equation*}
The choice of $s$ made before Lemma~\ref{bv-local} controls
the low spatial norms of these shape derivatives. The profile derivative is measured in $H^{K_{\rm reg}}$.

For the chain rule term in
\eqref{bv-ex-actual-approximate-good-equation}, write
\begin{equation*}
\dot S_g^\zeta=\partial_Y S_g+D_\Upsilon S_g[F_{g,\zeta}].
\end{equation*}
The diagonal entries act on each component, and the lower triangular
entry maps height to momentum. Quantization with the fixed cutoff and
\eqref{bv-ex-approximate-time-coefficients} give
\begin{equation*}
\|\dot S_g^\zeta\Upsilon\|_{X_s}
\leq C_{N,s,R}\bigl(1+\|w_Y\|_{H^{K_{\rm reg}}}\bigr)
\|\Upsilon\|_{X_s}
\end{equation*}
on an $X_s$ ball of radius $R$ with uniform nondegeneracy.
The spatial source estimates, the chain-rule bound used in the
definition of $R_g$, and the Sobolev contraction of $K_\zeta$
give a uniform $X_s$ bound for the smoothed remainder at the
admissible state $K_\zeta\Upsilon$.
For one principal block, abbreviate $J=J_\zeta$, $L=L_g$ and
$M=M_g$. Self-adjointness gives the exact identity
\begin{equation*}
2\operatorname{Re}\langle Mh,JLJ\pi\rangle
-2\operatorname{Re}\langle L\pi,JMJh\rangle=
2\operatorname{Re}\langle[J,M]h,LJ\pi\rangle
-2\operatorname{Re}\langle[J,L]\pi,MJh\rangle.
\end{equation*}
By \eqref{bv-ex-adapted-commutators}, the first term is a
pairing of $H^{1/2}$ with $H^{-1/2}$, and the second is a
pairing of $H^1$ with $H^{-1}$. They are bounded uniformly
by $C_{N,s}\|h\|_{H^1}\|\pi\|_{H^{1/2}}$.
After commuting with $A_s$, the displayed commutator bounds
give the same estimate at level $s$. To specify the energy,
freeze $(Y,\Upsilon,w)$ and let $\mathfrak q_s(Y,\Upsilon,w;Z)$
be the sum, over the two boundaries, of
\begin{equation*}
\sum_{A\in\{I,A_s\}}
\{\langle M_gAh^Z,Ah^Z\rangle
+\langle L_gA\pi^Z,A\pi^Z\rangle\},
\end{equation*}
where $(h^Z,\pi^Z)$ are the corresponding components of $S_gZ$,
together with the squared norm of the finite coordinates of $Z$.
The stated positive completions are included in $L_g,M_g$.
This is a bounded coercive quadratic form on $X_s$, uniformly
on the chosen neighbourhood. With
\begin{equation*}
\mathcal E_s^\zeta(Y)=
\mathfrak q_s(Y,\Upsilon^\zeta(Y),w(Y);\Upsilon^\zeta(Y)),
\end{equation*}
we obtain
\begin{equation}\label{bv-ex-uniform-local-energy}
\left|\frac{d}{dY}\mathcal E_s^\zeta\right|
\leq C_{N,s,R}\bigl(1+\|w\|_{H^{s-1}}^2
+\|w_Y\|_{H^{K_{\rm reg}}}\bigr)(1+\mathcal E_s^\zeta).
\end{equation}
The index $s$ fixed before the lemma covers both the low coefficient
derivatives in \eqref{bv-ex-approximate-time-coefficients} and the
source estimates at level $s$. On the chosen state ball, the constants
in \eqref{bv-ex-uniform-local-energy} are uniform in $\zeta$.

\smallskip
\noindent\emph{(d) A common interval and the local limit.}
Gronwall's inequality and a stopping argument in this ball give
an interval $\mathcal I$ about $Y_1$, independent of $\zeta$, on which
\begin{equation*}
\sup_{0<\zeta\leq1}
\left(\|\Upsilon^\zeta\|_{L^\infty(\mathcal I;X_s)}
+\|\partial_Y\Upsilon^\zeta\|_{L^\infty(\mathcal I;X_{s-3/2})}\right)<\infty.
\end{equation*}
Compactness gives, after passing to a subsequence,
\begin{equation*}
\Upsilon^\zeta\longrightarrow\Upsilon
\quad\hbox{in }C(\mathcal I;X_{s-\delta}),\qquad \delta>0,
\qquad \Upsilon\in C_w(\mathcal I;X_s).
\end{equation*}
Passing to the limit in \eqref{bv-ex-raw-regularization} recovers
\eqref{bv-ex-prescribed-profile-system}--\eqref{bv-ex-q-equation}.

To include profiles $w\in C_YH^{s-1}\cap W_Y^{1,\infty}H^{K_{\rm reg}}$, extend past the
endpoints and smooth in space and time. The resulting paths satisfy
\begin{equation*}
w^{(n)}\longrightarrow w\quad\hbox{in }C_YH^{s-1},\qquad
\sup_n\|\partial_Y w^{(n)}\|_{L_Y^\infty H^{K_{\rm reg}}}<\infty.
\end{equation*}
The uniform estimates give a solution by passage to the limit.
The equations use the instantaneous profile $w$. Uniqueness for a
fixed input, proved in Step~4, identifies the limit.

\medskip
\noindent\textbf{Step 3. Differentials of the boundary maps.}
Comparison of two boundary solutions requires differentials of the
full boundary maps. We use the comparison index $r$ and the background
regularity index $M=s-1$ fixed before Lemma~\ref{bv-local}. The lower dyadic threshold $b_{\rm dy}$ covers the
finite tame losses, with $r\geq b_{\rm dy}$. Constants may depend on
the background bound at order $M$.

The finite coordinates are smooth parameters whose differentials
are bounded by the Euclidean part of $X_r$.
The full-remainder bound \eqref{bv-ex-full-remainder-bound} and
its kinetic counterpart \eqref{bv-bd-kinetic-differentiated-remainder}
control remainder differences in Step~4. The coefficient estimate
\eqref{bv-ex-off-diagonal-matrix} is used in the coordinate smoothing
bound \eqref{bv-ex-coordinate-smoothing}. The last part of this step
restores the geometry dependence of the boundary primitives and
estimates the full pressure source in \eqref{bv-ex-pressure-differential}.

\smallskip
\noindent\emph{(a) Elliptic differentials.}
At fixed positive radius, pull back the elliptic problems to a
reference domain and lift their traces. For
$\operatorname{div}(A\nabla u)=f$ and a labelled set $I$ of directions,
\begin{equation*}
\operatorname{div}(A\nabla u_I)
=f_I-\operatorname{div}
\sum_{\varnothing\ne J\subseteq I}
A_J\nabla u_{I\setminus J}.
\end{equation*}
Differentiate the boundary conditions and scalar flux normalization
at the same time. Uniform ellipticity and the nonzero Schur
denominator give the estimates inductively for $|I|\leq4$.
The coefficients depend smoothly on the graphs and their first
derivatives. The product and trace estimates place the high norm on one
factor; the remaining factors use only finitely many lower
spatial derivatives.

Apply this induction to the Neumann inverse formulas
\eqref{bv-bd-neumann-first-variation} and
\eqref{bv-bd-neumann-second-variation}, including the moving
projections, and to the Poisson symbol recursion and its residual
equation. With the finite losses fixed above, the second and third
differential bounds from these equations control the coefficient
and smooth-kernel seminorms required in
\eqref{bv-ex-full-remainder-bound} and
\eqref{bv-ex-coefficient-matrix}. The variations of the finite-rank
normalizations are included.

\smallskip
\noindent\emph{(b) Differences of the complete remainders.}
Before changing to the good variables, we estimate the difference
of the complete remainders. This comparison is made at fixed positive
radius and may use the higher background norm. Fix $\mu=(\beta,c,q)$
and put $\nu=r+1/2$. For a potential row let $z=U$. For a kinetic
row let $z=(U,\widetilde v_-,\widetilde v_+)$, treating the primitives
in \eqref{bv-bd-external-primitives} as independent inputs. Write
$B(z)$ for the row and $p(z)$ for the matrix of local symbols of its
full spatial differential through order zero, in Weyl coordinates.
In particular, the local projection terms are those of
\eqref{bv-bd-neumann-hadamard}. The remaining spatial differential
consists of a local operator of order at most $-1$ and kernel
operators on the two circles. For the potential rows this is the
spatial decomposition used in Step~1 of Lemma~\ref{bv-boundary-energy}.
For the kinetic rows it is the decomposition in
Proposition~\ref{prop:bv-spatial-reduction}.
Use the same completion at bounded frequencies in the full and
fixed-cutoff symbols, and retain the complete affine reference row.
Define
\begin{equation*}
 \mathcal R_B(z)=B(z)-B(0)-DB(0)z
                  -\operatorname{Op}\big(p(z)-p(0)\big)z .
\end{equation*}
This groups the order $-1$ and high--high remainders of
\cite[Lemma~4.4]{LYIrrotational}. For smooth directions $u$ there
is an exact identity
\begin{equation}\label{bv-ex-full-remainder-differential}
 D\mathcal R_B(z)[u]
 =\bigl(DB(z)-DB(0)
       -\operatorname{Op}(p(z)-p(0))\bigr)u-\operatorname{Op}(Dp(z)[u])z .
\end{equation}
We estimate the differential in
\eqref{bv-ex-full-remainder-differential} as a whole, retaining both
the kernel and symbol variations, and prove
\begin{equation}\label{bv-ex-full-remainder-bound}
 \|D\mathcal R_B(z)[u]\|_{H^{\nu+1}}
        \leq C_{N,M,r}\|u\|_{H^\nu}.
\end{equation}

For the first line of \eqref{bv-ex-full-remainder-differential},
subtract the ordinary Weyl realization of $p$ from the full
differential at each of the two backgrounds. The local order $-1$
part maps $H^\nu$ to $H^{\nu+1}$. For a kernel part it suffices to
bound its kernel in $H^{\nu+1}_\theta L^2_{\theta'}$, since its
input is in $L^2$. The separated and finite-rank kernel estimates
of Proposition~\ref{prop:bv-spatial-reduction} give this bound at both backgrounds
using finitely many derivatives through $\nu+1$ plus the fixed
kernel loss. The coefficients of the local order $-1$ part have
the corresponding finite symbol bounds. These estimates allow
constants depending on $a_N^{-1}$.

It remains in this line to compare ordinary Weyl quantization to
the fixed cutoff. For a symbol of order $d\leq2$, the discarded
Fourier matrix entries, with output $n$ and input $l$, are bounded by
\begin{equation*}
 C_{N,M,L_{\rm dec}}\langle n-l\rangle^{-L_{\rm dec}}
             \langle(n+l)/2\rangle^d,\qquad
 \langle n\rangle+\langle l\rangle\leq C\langle n-l\rangle
\end{equation*}
apart from bounded frequencies. Indeed the cutoff equals one for
small coefficient frequency relative to the mean frequency, and
$L_{\rm dec}$ angular derivatives of the background symbol give the displayed
decay. 

After weighting by $\langle n\rangle^{\nu+1}
\langle l\rangle^{-\nu}$, both Schur sums converge, for instance
when $L_{\rm dec}>\nu+6$. This includes high input sent to low output as well
as high output produced by a background coefficient. Increase $M$
to cover these $L_{\rm dec}$ derivatives and the finite elliptic symbol loss.

For the second line, the fixed-cutoff bound from $H^{\nu+3}$ to
$H^{\nu+1}$ uses only fixed lower symbol seminorms of $Dp(z)[u]$.
To see the regularity of the direction here, differentiate the
finite elliptic symbol recursion at those seminorms. The products
contain one differentiated graph, trace or primitive. The
inhomogeneous terms in the differentiated elliptic equation above
are products of its coefficient variation with a background
gradient. At the fixed lower orders the product estimate and the
boundary trace estimate bound them by $C_{N,M,r}\|u\|_{H^\nu}$,
taking $r$ above these finite losses. The scalar moment variations
obey the same bound by their boundary pairings. 

The remaining
$\nu+3$ derivatives are taken on the background $z$ by the
cutoff operator. Thus $M$ can be chosen after $r$ to cover
$\nu+3$, the symbol orders just specified and the separated kernel
loss. This proves \eqref{bv-ex-full-remainder-bound}.
All the estimates persist on the line segment between two nearby
backgrounds with the stated higher bounds. Approximation by smooth
states, using one spare background derivative for convergence of the
coefficient and kernel seminorms, gives the same integrated
difference bound at the finite regularity of the lemma.

The gain of one derivative in \eqref{bv-ex-full-remainder-bound}
puts the differentiated remainder in both source spaces
$H^{r+1}$ and $H^{r+1/2}$. The independent primitive columns of
order zero retained in $p$ have the momentum target $H^{r+1/2}$.
For such a column $A(z)$ the two derivative terms are
\begin{equation*}
A(z)u\quad\text{and}\quad DA(z)[u]z.
\end{equation*} 
The first has the stated order-zero bound. In the second,
the lower symbol seminorms use the direction and the higher
norm uses the background, as above. The same product rule applies
to the retained lower-order compositions after the triangular change.
In particular, the kinematic entry from momentum to height left
after the DN symbol through order zero is extracted has order
at most $-1$. Thus an $H^{r+1/2}$ variation in momentum supplies
the $H^{r+1}$ kinematic source. The loss of order $3/2$ in the
principal system itself is treated by the energy estimate.

\smallskip
\noindent\emph{(c) Coefficient estimates and kinetic remainders.}
The centred-coordinate comparison in Step~4 also requires a
separated frequency estimate for the coefficient differentials.
The geometric coefficients of $S_g$ are smooth functions of the
graphs and their first derivatives. Its remaining coefficients are
boundary components of the reconstructed harmonic velocity and the
scalar flux normalizations. In particular, the formulas
\eqref{eq:bv-good-variables} and their outer counterpart involve
$\nabla\Phi$, but not the tangent Hodge field: at fixed state,
\begin{equation*}
D_wS_g=0.
\end{equation*} 
To differentiate the harmonic terms,
let $f$ be a fixed pulled back trace pair and vary the boundary
normally by $H$. In the following formula, $G$ denotes the full
physical DN map before component projections, and $B_f=Gf$:
\begin{equation}\label{bv-ex-full-dn-shape}
DG[H]f
=-G(HB_f)-\operatorname{div}_\Sigma(H\nabla_\Sigma f)
-\mathcal H HB_f.
\end{equation}
Indeed, the Eulerian potential variation has trace $-HB_f$, and
\begin{equation*}
\partial_n^2\Phi[f]=-\Delta_\Sigma f-\mathcal H B_f.
\end{equation*}
The variation of the normal contributes
$-\nabla_\Sigma H\cdot\nabla_\Sigma f$, giving
\eqref{bv-ex-full-dn-shape} on each component.
The graph parametrizations add tangential differentiation and smooth
geometric factors. Variations of the trace pair contribute $G\dot f$.
The Schur constants and the clock and centre parameters contribute
boundary pairings with smooth densities and finite-dimensional columns.

For a scalar coefficient $a_g(\Upsilon)$ of $S_g$, these formulas
express each boundary column of $Da_g$ as a finite-order local
pseudodifferential operator, a smooth operator between separated
boundaries, and finite-rank terms. The differentials of the
geometric and harmonic symbols in this change have the same
structure, at a fixed finite order. The shape equations above
and the differentiated elliptic equations give this description
also after one further state
differentiation with a fixed smooth background direction.

For the resulting matrix estimate, write $\mathcal A$ for one
boundary column, let $d_*$ bound its order, and put
$e_\ell(\theta)=e^{i\ell\theta}$. For each required integer $L_{\rm dec}$,
commuting $\mathcal A$ with $-i\partial_\theta$ a total of $L_{\rm dec}$
times differentiates its coefficients without increasing its order.
The local symbol bounds and integration by parts in the smooth kernels
therefore give
\begin{equation}\label{bv-ex-coefficient-matrix}
|\langle e_n,\mathcal A e_\ell\rangle|
\leq C_{N,M,L_{\rm dec}}\langle\ell\rangle^{d_*}
\langle n-\ell\rangle^{-L_{\rm dec}}.
\end{equation}
Only finitely many background derivatives are used for each fixed $L_{\rm dec}$.
The estimate also holds for the differentiated column, with a factor
equal to a sufficiently high norm of the background direction.
The choice of $M$ includes these derivatives and the two-derivative
reserve specified before the lemma.

On the supports of $S_{j-K_{\rm BW}}$ and $\Delta_k$ with $k\geq j+C$,
one has $|n-\ell|\asymp2^k$. The weighted Hilbert--Schmidt bound
deduced from \eqref{bv-ex-coefficient-matrix} is
\begin{equation}\label{bv-ex-off-diagonal-matrix}
\|S_{j-K_{\rm BW}}Da_g(\Upsilon)\Delta_k\|_{L^2\to H^J}
\leq C_{N,M,L_{\rm dec}}2^{(d_*+J+1-L_{\rm dec})k}
\leq C_{N,M,A_*}2^{-A_*k}.
\end{equation}
Choose $L_{\rm dec}>A_*+d_*+J+1$ after fixing the finitely many $A_*,J$
needed through index $r+3/2$, and then choose $M$.
Finite dimensional columns are bounded directly by their smooth
range profiles. Equation~\eqref{bv-ex-off-diagonal-matrix} will be
used for the differentiated coordinate change in Step~4.

Differentiating the symbols, kernels and reference affine map in
$\mu$ gives a contribution bounded by
$C_{N,M,r}|\dot\mu|$.
Apply the fixed-cutoff half-density transformation to the kinetic
output and express the boundary inputs through the frozen operator
$S_g(Y,\Upsilon)^{-1}$, as in
the derivation of \eqref{bv-bd-kinetic-remainder}. Differentiating
this composition with respect to $z$ and $\mu$, including the
derivatives of both transformations and the order-zero composition
terms, gives
\begin{equation}\label{bv-bd-kinetic-differentiated-remainder}
\|D R_{\pi,\rm kin}(z;\mu)[\dot z,\dot\mu]\|_{H^{r+1/2}}
\leq C_{N,M,r}\left\{
\|\dot U\|_{Y_r}
+\sum_i\|\dot{\widetilde v}_i\|_{H^{r+1/2}}
+|\dot\mu|\right\},
\end{equation}
where $Y_r$ is the boundary part of $X_r$ at fixed positive radius.
In particular,
\begin{equation*}
Y_r=(P_{\ge2}H^{r+1}\times P_{\ge2}H^{r+1/2})_-
\times(H^{r+1}\times (I-P_0)H^{r+1/2})_+,
\qquad \|\dot U\|_{H^{r+1/2}}\leq C\|\dot U\|_{Y_r}.
\end{equation*}
The same bounds for differences of the coefficient and kernel
seminorms, obtained by one further elliptic differentiation on the
bounded higher-regularity neighbourhood, give continuity of the
differential in the comparison operator norm.

\smallskip
\noindent\emph{(d) The Hodge and pressure differentials.}
The boundary primitives were held independent in
Proposition~\ref{prop:bv-spatial-reduction}. Their dependence on the
geometry must now be included in the full state differential. For a normal displacement $H$ with fixed extension $\xi$,
the Hodge differential at fixed profile $w$ satisfies
\begin{equation*}
\operatorname{div}\dot V=0,\quad
\operatorname{curl}\dot V=-\rho(\xi\cdot\nabla W)e_\phi,\quad
\dot V\cdot n=\operatorname{div}_\Sigma(HV_\tau),\quad
\oint_{C_-}\dot V\cdot d\ell
=\oint_{C_-}(\xi\times\omega)\cdot d\ell.
\end{equation*}
Elliptic estimates on the fixed domain give
\begin{equation*}
\|\dot V_\tau|_\Sigma\|_{H^r}
\leq C_{N,r,M}(\|H\|_{H^{r+1}}+|\dot\mu|),
\qquad \dot{\widetilde v}_\pm\in H^{r+1}.
\end{equation*}
For the pressure, the good potential variation gives
$\dot v\in H^r$ and
\begin{equation*}
\Delta\dot p_D=-2\operatorname{tr}(\nabla v\nabla\dot v),
\qquad \dot p_D|_\Sigma=-H\partial_np_D.
\end{equation*}
Let $\mathcal P_D(\Upsilon,w)$ denote the two momentum sources
obtained from the normal derivative $\partial_n p_D$ by the
componentwise flux projection, the Neumann inverse and the weighted
good-variable transformation in the momentum rows; the capillary
contribution is treated in the potential reduction.
Differentiating this entire composition gives
\begin{equation}\label{bv-ex-pressure-differential}
\|D\mathcal P_D(\Upsilon,w)[\dot\Upsilon,\dot w]
\|_{(H^{r+1/2})^2}
\leq C_{N,M,r}
\bigl(\|\dot\Upsilon\|_{X_r}+\|\dot w\|_{H^{r-1}}\bigr).
\end{equation}
To verify the trace index, the linearized Dirichlet equation for
$\dot p_D$ gives
$\dot p_D\in H^{r+1}$ and
$\partial_n\dot p_D\in H^{r-1/2}$. The variation of the pulled
back normal derivative also contains the normal and evaluation
variations; these multiply the higher background derivatives of
$p_D$ and obey the same bound.

The projection variation \eqref{bv-bd-projection-first-variation}
treats the flux moments,
while \eqref{bv-bd-neumann-first-variation} treats the inverse
on the fixed identified spaces.
The unvaried inverse sends $H^{r-1/2}$ data to $H^{r+1/2}$;
its shape variation has order zero in $H$ by
\eqref{bv-bd-neumann-hadamard}. The differentiated half-density
factors and fixed projections preserve $H^{r+1/2}$.

For a profile variation, the Hodge estimate gives
$\dot V\in H^r$ from $\dot w\in H^{r-1}$, and the same Dirichlet
argument applies. The removed flux moments satisfy the corresponding
finite-dimensional bound. Thus \eqref{bv-ex-pressure-differential}
includes the variations of the Dirichlet solution, flux projection,
Neumann inverse and good variable change. Its constant may depend
on $a_N^{-1}$, since this local comparison is made at fixed
positive radius.

\medskip
\noindent\textbf{Step 4. Centred coordinates and uniqueness.}
The frozen inverse alone does not compare the nonlinear coordinates
of two solutions. The differential estimates above provide this
comparison: for a common nearby state $\Upsilon_*$, held fixed in $Y$,
use the centred coordinate
\begin{equation*}
T_*(Y,\Upsilon)=S_g(Y,\Upsilon)
(\Upsilon-\Upsilon_*).
\end{equation*}  

\smallskip
\noindent\emph{(a) Smoothing of the centred coordinates.}
Choose the neighbourhood small in a fixed lower norm.
For the state difference 
\begin{equation*}
\Delta_\Upsilon=\Upsilon-\Upsilon_*,
\end{equation*} 
the fixed
coefficient cutoff gives
\begin{equation*}
\|DS_g[Z]\Delta_\Upsilon\|_{L^2}
\leq C_{N,M}\|\Delta_\Upsilon\|_{X_{b_0}}\|Z\|_{L^2},
\end{equation*}
where $b_0$ exceeds the finite order of the coefficient differentials.
In the dyadic sum, a high frequency of $Z$ interacts with a comparable
or higher frequency of $\Delta_\Upsilon$, whose $X_{b_0}$ norm absorbs the loss.
The separated blocks satisfy \eqref{bv-ex-off-diagonal-matrix}.
Thus $DS_g[\cdot]\Delta_\Upsilon$ is small on $L^2$ when $\|\Delta_\Upsilon\|_{X_{b_0}}$ is small.
We need the corresponding smoothing estimate both at index $r$
and at index $r+3/2$. Let $d_c$ cover the finite coefficient losses,
the component Sobolev shifts and the Bernstein factors, including
one further state differential. Increase the fixed background
index so that $M>r+d_c+7$ and so that the matrix estimates of
Step~3 hold with background directions in $X_{M-2}$.
For $r\leq\kappa\leq r+3/2$, we claim
\begin{equation}\label{bv-ex-coordinate-smoothing}
\begin{split}
\|DS_g[Z]\Delta_\Upsilon\|_{X_\kappa}
&\leq C_{N,M,r}\|\Delta_\Upsilon\|_{X_M}\|Z\|_{X_{\kappa-2}},\\
\|D^2S_g[Z,Z_1]\Delta_\Upsilon\|_{X_\kappa}
&\leq C_{N,M,r}\|Z_1\|_{X_{M-2}}
\|\Delta_\Upsilon\|_{X_M}\|Z\|_{X_{\kappa-2}}.
\end{split}
\end{equation}
Here the state stays in the fixed bounded subset of $X_M$,
and $Z_1$ is the higher regularity background direction.

To prove the first bound, decompose $\Delta_\Upsilon$ and $Z$ into dyadic blocks.
The fixed coefficient cutoff places the output frequency at the
frequency of $\Delta_\Upsilon$. If the frequency of $Z$ is no larger than a fixed
multiple of that of $\Delta_\Upsilon$, the coefficient estimate and Bernstein's
inequality bound the sum by
\begin{equation*}
C_{N,M,r}\sum_{j\geq0}2^{(\kappa+d_c+2-M)j}
\|\Delta_\Upsilon\|_{X_M}\|Z\|_{X_{\kappa-2}}.
\end{equation*}
This series converges at every indicated $\kappa$ by the choice of $M$.
If the frequency of $Z$ is larger, its contribution to the coefficient
frequencies retained on the $\Delta_\Upsilon$ block is bounded by
\eqref{bv-ex-off-diagonal-matrix}. Choose $A_*$ larger than the
finitely many Sobolev and Bernstein exponents in these sums.
Summing first over the frequency of $Z$ gives the same bound.
The finite-dimensional columns obey it directly.

For the second bound, differentiate the coefficient operator once
in the direction $Z_1$. The last part of Step~3 gives the same
matrix estimates, with the additional factor
$\|Z_1\|_{X_{M-2}}$. The preceding two sums then prove the second
line of \eqref{bv-ex-coordinate-smoothing}. The reserve of two
derivatives permits this direction to be evaluated on $F_g$, since
the vector field maps bounded subsets of $X_M$ into $X_{M-3/2}$,
which is contained in $X_{M-2}$ at fixed positive radius.

\smallskip
\noindent\emph{(b) Local coordinate comparison.}
For two states in a small convex neighbourhood, set
$\Upsilon_\theta=\Upsilon_2+\theta(\Upsilon_1-\Upsilon_2)$. Then
\begin{align*} 
T_*(\Upsilon_1)-T_*(\Upsilon_2)
&=\int_0^1 DT_*(\Upsilon_\theta)
 [\Upsilon_1-\Upsilon_2]\,d\theta,\\
DT_*(\Upsilon_\theta)[Z]
&=S_g(\Upsilon_\theta)Z
 +DS_g(\Upsilon_\theta)[Z](\Upsilon_\theta-\Upsilon_*). 
\end{align*}
Let $A$ denote the averaged operator above. Its $L^2$ closeness to
$S_g(\Upsilon_*)$ gives $\|Z\|_{L^2}\leq C\|AZ\|_{L^2}$, including
the finite coordinates. Its triangular principal part has elliptic
diagonal entries. The order-zero parametrix, with error of order $-1$,
and \eqref{bv-ex-coordinate-smoothing} at $\kappa=r$ give
\begin{equation*}
\|Z\|_{X_r}\leq C(\|AZ\|_{X_r}+\|Z\|_{X_{r-1/2}}),
\end{equation*} 
allowing
half a derivative for the different component indices.
Componentwise interpolation gives
\begin{equation*}
\|Z\|_{X_{r-1/2}}\leq\delta\|Z\|_{X_r}+C_\delta\|Z\|_{L^2}.
\end{equation*}
Choose $C\delta<1/2$ and absorb the top norm. The $L^2$ lower bound
and $\|AZ\|_{L^2}\leq C\|AZ\|_{X_r}$ then give, for
$Z=\Upsilon_1-\Upsilon_2$,
\begin{equation*} 
\|\Upsilon_1-\Upsilon_2\|_{X_r}
\leq C_{N,r,M}
\|T_*(\Upsilon_1)-T_*(\Upsilon_2)\|_{X_r}.
\end{equation*}

\smallskip
\noindent\emph{(c) The centred equation and uniqueness.}
Write $\dot S_g=\partial_YS_g+DS_g[F_g]$. The spatial remainder
convention in Step~2 gives the exact centred equation
\begin{equation*}
\partial_YT_*=A_gT_*+\mathcal C_g,\qquad
\mathcal C_g=R_g+A_gS_g\Upsilon_*
+\dot S_g(\Upsilon-\Upsilon_*).
\end{equation*}
The first line of \eqref{bv-ex-coordinate-smoothing}, at
$\kappa=r$, gives
\begin{equation*} 
\|DS_g[F_g]\Upsilon_*\|_{X_r}
\leq C_{N,M,r}\|\Upsilon_*\|_{X_M}\|F_g\|_{X_{r-2}}\leq C_{N,M,r}\|\Upsilon_*\|_{X_M}\|F_g\|_{X_{r-3/2}}. 
\end{equation*}
Thus the gain of two derivatives absorbs the $3/2$ derivative loss.
The differentiated elliptic equations also give
\begin{equation*}
\|DF_g[Z]\|_{X_{r-3/2}}\leq C_{N,M,r}\|Z\|_{X_r}
\end{equation*}
for a fixed input. In the differential of $DS_g[F_g]\Delta_\Upsilon$, the term
$DS_g[DF_g[Z]]\Delta_\Upsilon$ is therefore controlled by this gain; the term
$D^2S_g[Z,F_g]\Delta_\Upsilon$ uses the second line of
\eqref{bv-ex-coordinate-smoothing}, with $F_g\in X_{M-2}$.
For $DS_g[F_g]Z$, the higher background bound controls the
coefficient differential in the fixed lower symbol seminorms,
so the order-zero cutoff operator is bounded on $X_r$.
For the differential of $A_gS_g\Upsilon_*$, use the lower
coefficient seminorms controlled by $X_r$ and the fixed cutoff
operator bounds. Since $A_g:X_{r+3/2}\to X_r$, the first line of
\eqref{bv-ex-coordinate-smoothing} at $\kappa=r+3/2$ gives
\begin{equation*}
\|A_gDS_g[Z]\Upsilon_*\|_{X_r}
\leq C_{N,M,r}
\|DS_g[Z]\Upsilon_*\|_{X_{r+3/2}}\leq C_{N,M,r}\|\Upsilon_*\|_{X_M}\|Z\|_{X_r}.
\end{equation*}
The other term $DA_g[Z]S_g\Upsilon_*$ has the same bound:
the coefficient differential uses the fixed lower seminorms,
and the additional spatial derivatives fall on the fixed higher
background $\Upsilon_*$. Thus this composition is controlled
in the source space $X_r$.
For each potential row, \eqref{bv-ex-full-remainder-differential} and \eqref{bv-ex-full-remainder-bound} bound the differential of the
spatial remainder. The lower-order entries are estimated in their
respective source spaces as in Step~3. The kinetic contribution uses
\eqref{bv-bd-kinetic-differentiated-remainder}, and the pressure
uses \eqref{bv-ex-pressure-differential}. The scalar and finite-rank
terms satisfy the same bounds by differentiation of their
moment formulas.

Set $T_i=T_*(Y,\Upsilon_i)$, $\delta\Xi=T_1-T_2$ and write
$A_i,\mathcal C_i$ for the corresponding operators and sources.
For the same prescribed input $w$, subtraction gives
\begin{equation*}
\partial_Y\delta\Xi=A_1\delta\Xi+\mathcal F_\delta,
\qquad
\mathcal F_\delta=(A_1-A_2)T_2+\mathcal C_1-\mathcal C_2,
\end{equation*}
with
\begin{equation*}
\|\mathcal F_\delta\|_{X_r}
\leq C_{N,M,r}\|\Upsilon_1-\Upsilon_2\|_{X_r}
\leq C_{N,M,r}\|\delta\Xi\|_{X_r}.
\end{equation*}
In $(A_1-A_2)T_2$, the differentiated coefficients are controlled
in fixed lower norms and the higher derivatives fall on $T_2$.
The resulting quadratic form $\mathcal D_r$ satisfies
\begin{equation*}
\mathcal D_r\asymp\|\delta\Xi\|_{X_r}^2,\qquad
\left|\frac{d\mathcal D_r}{dY}\right|
\leq C_{N,r,M}(1+\|w_Y\|_{H^{K_{\rm reg}}})\mathcal D_r.
\end{equation*}
Gronwall's inequality in either direction proves uniqueness for a fixed
input. Consecutive centred charts cover the common interval.

\medskip
\noindent\textbf{Step 5. Strong continuity and continuation.}
For the limiting path, write
\begin{equation*}
\mathfrak q_{s,Y}(Z)=\mathfrak q_s(Y,\Upsilon(Y),w(Y);Z),\qquad
\mathcal E_s(Y)=\mathfrak q_{s,Y}(\Upsilon(Y)).
\end{equation*}
To upgrade the weak $X_s$ continuity obtained in Step~2, fix
$Y_*\in\mathcal I$, test the exact equation with $J_\delta$ and let
$\delta\downarrow0$.
The commutator bounds \eqref{bv-ex-adapted-commutators} and the
bounds for the time derivatives of $J_\delta$ and $A_s$ in
\eqref{bv-ex-mollifier-symbol-bounds} are uniform in $\delta$.
On smooth spatial tests, the regularization commutators and
$\partial_YJ_\delta$ tend to zero in the energy pairings for almost
every $Y$. The uniform bounds extend this convergence by density.
Dominated convergence then passes the integrated energy inequality
to the limit, using the
$L^\infty_YX_s$ bound and the integrable bounds for the first
$Y$ derivatives of the coefficients. This passage does not require strong $X_s$ continuity.
Applying \eqref{bv-ex-uniform-local-energy} in both directions gives
\begin{equation*} 
|\mathcal E_s(Y)-\mathcal E_s(Y_*)|
\leq C_{N,s,R}\int_{\min(Y,Y_*)}^{\max(Y,Y_*)}
\bigl(1+\|w\|_{H^{s-1}}^2+\|w_Z\|_{H^{K_{\rm reg}}}\bigr)
(1+\mathcal E_s(Z))\,dZ \longrightarrow0,
\end{equation*}
as $Y\to Y_*$. The lower norm continuity from Step~2 and the
elliptic coefficient
bounds give continuity of the finitely many symbol seminorms
defining this form. With the fixed cutoff, $A_sS_g$ is continuous
in operator norm from $X_s$ to the boundary spaces
$H^1\times H^{1/2}$; likewise $M_g:H^1\to H^{-1}$ and
$L_g:H^{1/2}\to H^{-1/2}$ are continuous in their operator norms.
Consequently
\begin{equation*}
|\mathfrak q_{s,Y}(Z)-\mathfrak q_{s,Y_*}(Z)|
\leq\epsilon_s(Y)\|Z\|_{X_s}^2,\qquad
\epsilon_s(Y)\longrightarrow0.
\end{equation*}
The $X_s$ bound and energy continuity therefore imply
\begin{equation*}
\mathfrak q_{s,Y_*}(\Upsilon(Y))
\longrightarrow
\mathfrak q_{s,Y_*}(\Upsilon(Y_*)).
\end{equation*}
Since $\Upsilon(Y)\rightharpoonup\Upsilon(Y_*)$ in $X_s$,
weak convergence and norm convergence for this fixed coercive
quadratic form give
$\|\Upsilon(Y)-\Upsilon(Y_*)\|_{X_s}\to0$.
The solution restarts while the admissibility margins, the $X_s$
state and $H^{s-1}$ profile bounds, and the prescribed bounds on
$w_Y$ persist.
For $k\geq s$, the tame products in Steps 1 and 2 give
\begin{equation*}
\left|\frac{d\mathcal E_k}{dY}\right|
\leq C_{N,k}(1+\|w_Y\|_{H^{K_{\rm reg}}})(1+\mathcal E_k)
+C_{N,k}\|w\|_{H^{k-1}}^2,
\end{equation*}
where $C_{N,k}$ depends on the state and profile bounds at the fixed
lower orders.
Gronwall on the same interval gives persistence of spatial smoothness.

For continuous dependence on the prescribed profile, fix the data
at $Y_1$ and a common compact interval $\mathcal I$ of existence. For profiles in the prescribed class,
suppose $w_n\to w$ in $C(\mathcal I;H^r)$, with uniform
admissibility margins and
\begin{equation*}
\sup_n\left(
\|\Upsilon_n\|_{L^\infty(\mathcal I;X_s)}
+\|w_n\|_{L^\infty(\mathcal I;H^{s-1})}
+\|\partial_Yw_n\|_{L^\infty(\mathcal I;H^{K_{\rm reg}})}
\right)<\infty.
\end{equation*}
The equations bound $\partial_Y\Upsilon_n$ uniformly in
$L^\infty(\mathcal I;X_{s-3/2})$. Compact spatial embedding and
interpolation give a subsequence converging in $C(\mathcal I;X_r)$.
The limit solves
\eqref{bv-ex-prescribed-profile-system}--\eqref{bv-ex-q-equation}
with input $w$. Uniqueness for this fixed input identifies the limit
with $\Upsilon$, so the whole sequence converges in
$C(\mathcal I;X_r)$.
\end{proof}

\subsection{The coupled construction}
\label{subsec:bv-coupled-construction}

\begin{proof}[Construction for Theorem~\ref{thm:bv-collapse}]
At a finite terminal radius, the local boundary solution and forward
transport define a map on a compact convex set of profile paths.
Uniform estimates allow us to apply Schauder's theorem and then take
the terminal radius to zero. For sufficiently small $\varepsilon>0$
this gives a smooth
solution on $[Y_0,\infty)$ with initial profile
$w_0$ given by \eqref{eq:bv-initial-profile} and the prescribed
terminal data. For every fixed $s$ it will satisfy
\begin{equation}\label{bv-ex-final-construction-bounds}
E_s\leq P_s a^2,\qquad
O_s\leq P_s a,\qquad
Z_s(w)\leq \mathfrak M_s,\qquad
N_1(w)\leq C_{\chi_0}\varepsilon,\qquad
|\beta|\leq P_ba^{3/4}.
\end{equation}
The constants $\mathfrak M_s$ are chosen successively in Step~2.
The profile will stay in $S<s_2$, and $Y_0$ will be independent
of $s$.

\medskip
\noindent\textbf{Step 1. Terminal data.}
We first stop the evolution at a positive cavity radius. Choose
$Y_N\to\infty$ and set $a_N=a(Y_N)$. At $Y_N$ prescribe
\begin{equation}\label{bv-ex-terminal-data}
\eta=F=\zeta_+=g=0,\qquad \beta=0,\qquad c=c_*,
\qquad q=\left(\frac{E_c}{\pi L_*\Lambda_N}\right)^{1/2}.
\end{equation}
Here $\Lambda_N$ is the capacity parameter of this geometry at radius
$a_N$. These data give
\begin{equation*}
E_f(Y_N)=E_c,\qquad B_f(Y_N)=0,
\end{equation*}
while the profile is determined by the fixed point below; no
additional terminal profile is prescribed.

\medskip
\noindent\textbf{Step 2. Uniform estimates.}
Fix once a lower comparison index $r$ for Lemma~\ref{bv-local}
and the spatial index $K_{\rm reg}$ required there to bound $w_Y$.
Choose $k_*\geq K_{\rm reg}$, with $k_*+2\geq r$, to cover the fixed spatial
orders of $w_Y$ used in the uniform boundary energy estimate and
in \eqref{bv-ex-approximate-time-coefficients}.
Then fix a base index $b$ above the boundary threshold and all
fixed symbol, kernel and time coefficient orders, with
$b-1\geq k_*+1$.

For levels below this threshold we use the boundary estimate at
level $b$. Write
\begin{equation*}
Z_j^{\rm in}=Z_j(w^{\rm in}),\quad\text{and}\quad Z_j^{\rm out}=Z_j(w^{\rm out})
\end{equation*} 
for the prescribed and transported
profiles, respectively.
For $Z_j^{\rm in}\leq \mathfrak M_*$, $j\leq b-1$, the terminal estimate
in Corollary~\ref{cor:bv-terminal-bounds}, obtained from
\eqref{bv-bd-combined-tail} with zero terminal data, gives,
under provisional low norm bounds with
$\int\ell_b\leq1$,
\begin{equation}\label{bv-ex-base-boundary}
E_b(Y)+a(Y)O_b(Y)
\leq C_bP_b(Y)a(Y)^2(1+\varepsilon+\mathfrak M_*).
\end{equation}
Here $\ell_b$ is the integrable coefficient
$e_b+P_b(a^{1/4}+a)$ in \eqref{bv-bd-combined-tail} at level $b$.
Since $E_f(Y_N)=E_c$, the provisional bound
$E_c/2\leq E_f\leq2E_c$ holds near $Y_N$.
On this interval, Proposition~\ref{bv-scalars} gives
\begin{equation*}
|E_f(Y)-E_c|
\leq\int_Y^{Y_N}|E_{f,y}|\,dy
\leq\int_{Y_0}^{\infty}P_b(y)a(y)\,dy<E_c/4
\end{equation*}
for $Y_0$ sufficiently large. Thus $3E_c/4<E_f<5E_c/4$.
By \eqref{bv-ex-base-boundary}, the same proposition and the
first mode estimate in Corollary~\ref{cor:bv-terminal-bounds} give
\begin{equation*}
1-p_q\geq3/4,\qquad |\beta|\leq P_ba=o(a^{3/4}).
\end{equation*}
The clock and centre estimates \eqref{bv-ex-scalar-estimates},
together with the boundary tails, improve the remaining geometric bounds.

Use $K_j$ and $B_k$ defined before Lemma~\ref{bv-transport}.
The Hodge and Poisson extension estimates, combined in
\eqref{bv-st-coefficient-bound}--\eqref{bv-st-boundary-to-coefficients},
give for all $j\leq b-1$
\begin{equation}\label{bv-ex-base-coefficients}
K_j\leq C_bP_b(1+\varepsilon+\mathfrak M_*)a^{1/4},\qquad
B_1\leq C_bP_b(1+\varepsilon+\mathfrak M_*)a^{3/4}.
\end{equation}
For $j=0,1$ use $K_j\leq K_2$.
The bound for $B_1$ follows from the separate
low order supremum estimate.
The transported output, with initial data
\eqref{eq:bv-initial-profile}, satisfies
\begin{equation}\label{bv-ex-base-feedback}
\max_{j\leq b-1}\sup_Y Z_j^{\rm out}
\leq e^{C_b\int B_1}
\left[I_*+
C_b\varepsilon(1+\varepsilon+\mathfrak M_*)T_{5/4}(Y_0)\right],
\end{equation}
where $I_*\leq C_{b,\chi_0}\varepsilon$ bounds the initial
norms and
\begin{equation*}
T_\gamma(Y_0)=\int_{Y_0}^\infty P_b(y)a(y)^\gamma\,dy.
\end{equation*}
Choose $C_{\chi_0}$ so that
$N_1(w_0)\leq C_{\chi_0}\varepsilon/4$ uniformly for
$0<a\leq\ell$. For $\exp(C\int B_1)\leq2$,
\begin{equation*}
N_1(w^{\rm out})
\leq N_1(w_0)\exp\left(C\int B_1\right)
\leq C_{\chi_0}\varepsilon/2.
\end{equation*}

Choose $\mathfrak M_*>4e^{C_b}C_{b,\chi_0}\varepsilon$.
The required bounds for $w_Y^{\rm out}$ follow from
\begin{equation*}
M_k(w_Y^{\rm out})
\leq C_b\left[
B_0 M_{k+1}(w^{\rm out})
+K_kN_1(w^{\rm out})\right],
\qquad 0\leq k\leq k_*.
\end{equation*}
Thus, after fixing $\mathfrak M_*$, choose
$D_{\rm t}=C_b'(1+\varepsilon+\mathfrak M_*)(\mathfrak M_*+\varepsilon)$. It bounds
$M_k(w_Y^{\rm out})/(P_ba^{-1/4})$.
The supremum time norm obeys
\begin{equation*}
N_0(w_Y^{\rm out})
\leq C_b\varepsilon(1+\varepsilon+\mathfrak M_*)P_ba^{3/4}.
\end{equation*}
The constant in this supremum bound is fixed before $Y_0$ is chosen.
The low energy coefficient integral is bounded by
\begin{equation}\label{bv-ex-base-time-feedback}
\int_{Y_0}^\infty\ell_b
\leq C_b(1+Y_0)^{-1/2}
+C_bP_*(1+\mathfrak M_*,D_{\rm t})T_{1/4}(Y_0).
\end{equation}
The first term comes from the outer capacity and reference Taylor
coefficients. The polynomial $P_*$ depends on the fixed base index.
With $\mathfrak M_*,D_{\rm t}$ fixed, \eqref{bv-ex-base-feedback} and
\eqref{bv-ex-base-time-feedback} allow us to choose $Y_0$ so that
\begin{gather*}
C_b(1+Y_0)^{-1/2}
+C_bP_*(1+\mathfrak M_*,D_{\rm t})T_{1/4}(Y_0)<\tfrac12,\\
e^{C_b\int B_1}C_b\varepsilon(1+\varepsilon+\mathfrak M_*)
T_{5/4}(Y_0)<\mathfrak M_*/4,\qquad T_1(Y_0)<E_c/4.
\end{gather*}
At the same time, impose the support condition needed in Step~3:
\begin{equation}\label{bv-ex-support-tail}
\exp\left(C_b(1+\varepsilon+\mathfrak M_*)
\int_{Y_0}^\infty P_b(y)a(y)^{3/4}\,dy\right)<s_2/s_1.
\end{equation}
By \eqref{bv-ex-base-coefficients}, this is another condition on
the fixed lower orders. Taking $Y_0$ large enough also improves the
lower supremum and geometric bounds. In particular, $\int\ell_b<1/2$
closes the provisional integral bound. We now keep $Y_0$ fixed and
put $\mathfrak M_j=\mathfrak M_*$ for $j<b$.
Fix the polynomial in \eqref{bv-ex-base-coefficients} and set
\begin{equation*}
\mathcal B_*=
C_b(1+\varepsilon+\mathfrak M_*)
\int_{Y_0}^{\infty}P_b(y)a(y)^{3/4}\,dy.
\end{equation*}
Then 
\begin{equation*}
\int_{Y_0}^{Y_N}B_1\leq\mathcal B_*
\end{equation*} 
for every
admissible input and every $N$.

For $m\geq b$, the Poisson extension of $\eta_Y$ gains half a
derivative, so the transport coefficient at level $m$ uses
$E_m,O_m$ and $Z_{m-1}^{\rm in}$. If
$Z_{m-1}^{\rm in}\leq \mathfrak M_{m-1}$, backward integration of
\eqref{bv-bd-combined-tail} gives
\begin{equation*}
E_m+aO_m\leq C_mP_m a^2(1+\varepsilon+\mathfrak M_{m-1}),
\qquad
aK_m\leq C_mP_m(a^{5/4}+a^2)(1+\varepsilon+\mathfrak M_{m-1}),
\end{equation*}
where the second bound follows from
\eqref{bv-st-coefficient-bound}--\eqref{bv-st-boundary-to-coefficients},
absorbing $q^{-1}$ into the polynomial. The constants depend on
$m$, the state and profile bounds at the fixed lower orders and
the fixed admissibility margins, but not on $Y_N$ or the truncation
level. Choose $\mathfrak M_m\geq \mathfrak M_{m-1}$ larger than
\begin{equation*}
2\exp(C_m\mathcal B_*)
\left\{Z_m(w_0)+C_m\varepsilon(1+\varepsilon+\mathfrak M_{m-1})
\int_{Y_0}^\infty P_m(y)
\bigl(a(y)^{5/4}+a(y)^2\bigr)\,dy\right\}.
\end{equation*}
Then \eqref{bv-st-sobolev-transport} gives
$Z_m^{\rm out}<\mathfrak M_m$. This determines the higher constants
successively from the base bounds. Each higher coefficient has
a finite integral on the already fixed tail; these integrals need
not be small simultaneously in $m$. Thus $Y_0$ remains unchanged.

\medskip
\noindent\textbf{Step 3. The compact path set and support.}
Fix a finite integer $m>k_*+2$ above the base orders and large
enough that $s=m$ satisfies the finite regularity choices before
Lemma~\ref{bv-local}. We use this boundary index and profile
regularity $H^{m-1}$. We apply Schauder's theorem in the space
\begin{equation*}
C([Y_0,Y_N];H^{k_*+2}([0,s_3]\times\T)).
\end{equation*}
Let $A_{N,m}$ consist of paths with initial value
\eqref{eq:bv-initial-profile}, support in $S\leq s_2$,
$N_1(w)\leq C_{\chi_0}\varepsilon$, and $Z_j(w)\leq \mathfrak M_j$
for $0\leq j\leq m$, with $\mathfrak M_j=\mathfrak M_*$ for $j<b$.
Impose the time bounds
\begin{equation}\label{bv-ex-path-time-bounds}
\begin{split}
M_k(w_Y)&\leq D_{\rm t}P_ba^{-1/4},\qquad 0\leq k\leq k_*,\\
N_0(w_Y)&\leq C_b\varepsilon(1+\varepsilon+\mathfrak M_*)P_ba^{3/4},
\end{split}
\end{equation}
with $D_{\rm t}$, the constant in the second line and the bounding
polynomial fixed in Step~2.
The fixed gap between $s_2$ and the outer collar permits the
closed support condition in all preceding estimates.

Since $a\geq a_N>0$, these paths satisfy
\begin{equation*}
\sup_Y\|w(Y)\|_{H^m}
+\|w_Y\|_{L^\infty_YH^{k_*}}\leq C_{N,m}.
\end{equation*}
Interpolation gives
\begin{align*}
\|w(Y)-w(Y')\|_{H^{k_*+2}}
&\leq C_{N,m}|Y-Y'|^{\frac{m-k_*-2}{m-k_*}},\\
\|w(Y)-w(Y')\|_{H^{m-1}}
&\leq C_{N,m}|Y-Y'|^{\frac{1}{m-k_*}}.
\end{align*}
The second estimate supplies the profile continuity required by
Lemma~\ref{bv-local} at $s=m$. The first, together with
$H^m\Subset H^{k_*+2}$, gives compactness by Arzel\`a--Ascoli.

The constraints defining $A_{N,m}$ are convex.
They are closed in this topology: the time derivative bounds pass
by distributional convergence, the higher norm bounds by weak lower
semicontinuity, and the support and initial conditions by strong
convergence.

The constant path $w(Y)=\varepsilon\chi_0(S)\cos(2\theta)$
belongs to $A_{N,m}$. Indeed, $w_Y=0$, and the weights in $Z_j^2$ satisfy
\begin{equation*}
\partial_a\bigl[a^2(a^2+S)^{2i-1}\bigr]
=2a(a^2+S)^{2i-2}(S+2ia^2)\geq0,\qquad i\geq0.
\end{equation*}
Thus $Z_j(w(Y))\leq Z_j(w(Y_0))\leq \mathfrak M_j$, and $A_{N,m}$ is
nonempty, compact and convex.

To check that transport preserves the support condition, consider
the radial characteristics of \eqref{bv-st-transport-equation}:
they satisfy
\begin{equation*}
S_Y=\lambda b_1,\qquad (\lambda b_1)|_{S=0}=0,\qquad
|\partial_S(\lambda b_1)|\leq CB_1.
\end{equation*}
Hence, as long as they stay in the enlarged collar,
\begin{equation*}
S_Y\leq CB_1S,\qquad
S(Y)\leq S(Y_0)\exp\left(C\int_{Y_0}^YB_1(y)\,dy\right)<s_2,
\end{equation*}
by \eqref{bv-ex-support-tail}. Since the initial support lies in
$S\leq s_1$, no characteristic from the support reaches $S=s_2$.

\medskip
\noindent\textbf{Step 4. A fixed point at finite terminal radius.}
Given $w\in A_{N,m}$, solve the boundary
problem \eqref{bv-ex-prescribed-profile-system}--\eqref{bv-ex-q-equation}
backwards from \eqref{bv-ex-terminal-data}.
Lemma~\ref{bv-local}, the strict improvements of the low norm
bounds, and the estimates at higher orders continue this solution through
$[Y_0,Y_N]$.
The resulting velocity and chart transport the initial profile
\eqref{eq:bv-initial-profile} forward and define the output $F_Nw$.
The bounds \eqref{bv-ex-base-feedback} and
\eqref{bv-ex-path-time-bounds}, the recursion for $\mathfrak M_j$ and
\eqref{bv-ex-support-tail} give
$F_N(A_{N,m})\subset A_{N,m}$.

To prove continuity of $F_N$, let $w_j,w\in A_{N,m}$ with
$w_j\to w$ in $C_YH^{k_*+2}$. Since $k_*+2\geq r$,
Lemma~\ref{bv-local} gives convergence of the boundary states
in $C_YX_r$. The elliptic and chart estimates at fixed $N$ then
give convergence of the transport coefficients in $C_YC^1$.
The characteristic flows and their inverses converge uniformly,
so the common smooth initial profile gives
\begin{equation*}
\|F_Nw_j-F_Nw\|_{C_YL^2}\longrightarrow0.
\end{equation*}
Interpolating with the uniform $H^m$ bound, we obtain
\begin{equation*}
\|F_Nw_j-F_Nw\|_{C_YH^{k_*+2}}
\leq C_{N,m}\|F_Nw_j-F_Nw\|_{C_YL^2}^{1-(k_*+2)/m}
\longrightarrow0.
\end{equation*}
Thus $F_N$ is continuous, and Schauder's theorem gives a fixed
point in $A_{N,m}$, coupling the boundary and transport equations.

We first remove the regularity truncation, keeping $N$ fixed and
letting $m\to\infty$.
The bounds at each fixed spatial order are independent of the
higher truncation index $m$, and the equations bound the time
derivatives in lower spatial norms. A diagonal subsequence of the
states and profiles converges in
$C([Y_0,Y_N];X_k\times H^k)$ for every fixed nonnegative integer
$k$. For each required output order, use convergence at a higher
order covering the finite losses in the chart and elliptic
reconstructions. The reconstructed fields and all products in the
boundary and transport equations then converge, so both equations
hold in the limit. The initial profile, terminal boundary data
and all fixed order bounds also pass to the limit.

Their right sides are smooth functions of the state and profile.
The equations give one continuous time derivative. Repeated
differentiation of the equations and elliptic reconstructions
then gives smoothness at every finite $Y$.

\medskip
\noindent\textbf{Step 5. The infinite tail and the Euler equation.}
We now let the terminal radius tend to zero. On each fixed interval $\mathcal I=[Y_0,Y_1]$, the bounds are
independent of $N$ and $a\geq a(Y_1)>0$. The spatial estimates and
the bounds for time derivatives give compactness of the states
and profiles in $C(\mathcal I;X_k\times H^k)$ for every fixed
nonnegative integer $k$. A diagonal subsequence as $N\to\infty$
yields a limit on $[Y_0,\infty)$ with the original initial
profile. The equations pass to the limit as in Step~4. Repeated
differentiation then gives smoothness in time. The boundary
inequalities and scalar tail give
\begin{gather*}
E_s+aO_s\leq P_sa^2,\\
|E_f^N(Y)-E_c|\leq\int_Y^{Y_N}P_b(y)a(y)\,dy,\qquad
|E_f(Y)-E_c|\leq\int_Y^\infty P_b(y)a(y)\,dy\longrightarrow0.
\end{gather*}
Thus \eqref{bv-ex-final-construction-bounds} holds. The geometric
and first mode bounds persist. Integrating the centre and outer
boundary velocities from the terminal data
\eqref{bv-ex-terminal-data} recovers $(D_*,c_*)$ in the limit.

It remains to recover the full Euler momentum equation. The coupled
solution satisfies $(\partial_t+v\cdot\nabla)W=0$.
For the reconstructed velocity, set
\begin{equation*}
R_e=v_t+(v\cdot\nabla)v+\nabla p.
\end{equation*}
Since $\omega=\rho We_\phi$ and $v$ is axisymmetric with zero
swirl, scalar transport gives the full vorticity equation:
\begin{equation*}
(\partial_t+v\cdot\nabla)\omega
=v^\rho W e_\phi=(\omega\cdot\nabla)v.
\end{equation*}
Together with $\operatorname{div}v=0$ and
$\Delta p=-\operatorname{tr}((\nabla v)^2)$, this yields
\begin{equation*}
\operatorname{div}R_e=0,\qquad \operatorname{curl}R_e=0.
\end{equation*}
The projected normal equations give $\Pi_n(R_e\cdot n)=0$.
The scalar flux equation \eqref{bv-ex-q-equation} and
\eqref{bv-ex-Q-and-capacity} give
\begin{equation*}
\int_{\Gamma_-}R_e\cdot n
=\int_\Omega R_e\cdot\nabla u_c=0,\qquad
\int_{\Gamma_+}R_e\cdot n
=-\int_{\Gamma_-}R_e\cdot n=0.
\end{equation*}
By \eqref{bv-st-capacitary-projection}, $R_e\cdot n=0$ on both
boundary components. Differentiating the zero inner circulation
along the material meridian yields
\begin{equation*}
\oint_{C_-}R_e\cdot d\ell
=\frac{d}{dt}\oint_{C_-}v\cdot d\ell=0.
\end{equation*}
Here the material line-integral formula contributes the gradient
$\nabla(|v|^2/2)$, whose integral around the closed loop is zero;
the pressure gradient also has zero period. The residual is
meridional, so its line integral along every azimuthal circle
is zero.
These are the two periods needed to identify the residual as a
gradient: the fluid domain is diffeomorphic to a meridional annulus
times $\mathbb T$, and the inner meridian and an azimuthal circle generate
its first homology. The vanishing of both periods therefore gives
$R_e=\nabla h$ for a single-valued harmonic function $h$, and
\begin{equation*}
\int_{\Omega_t}|R_e|^2
=\int_{\Omega_t}|\nabla h|^2
=\int_{\partial\Omega_t}h\partial_nh=0.
\end{equation*}
Thus $R_e=0$. Since $t_Y>0$ is smooth at each finite $Y$,
the solution is smooth in physical time and satisfies
\eqref{eq:euler-toroidal-front}.
\end{proof}

\subsection{Vorticity, collapse law and energy limits}

\begin{proof}[Completion of the proof of Theorem~\ref{thm:bv-collapse}]
Consider the solution constructed in
Subsection~\ref{subsec:bv-coupled-construction}, which satisfies
\eqref{bv-ex-final-construction-bounds}, $E_f\to E_c$ and the
prescribed terminal geometry.
Write $\Omega(Y)=\Omega_{t(Y)}$ and use $Y$ as the time argument
of the fields in this proof.

\medskip
\noindent\textbf{Step 1. Weighted derivatives and the boundary flow.}
For each fixed $m$, the unit cylinder covering in
Lemma~\ref{bv-transport} gives $B_m\leq C_mK_{m+2}$ uniformly
in the logarithmic length. The high curl contribution satisfies
\begin{equation*}
t_YM_{m+1}(w)
=\Theta\frac aq Z_{m+1}(w)\leq P_ma.
\end{equation*}
Substitution of \eqref{bv-ex-final-construction-bounds} into
\eqref{bv-st-coefficient-bound}--\eqref{bv-st-boundary-to-coefficients}
at level $m+2$ therefore gives
\begin{equation*}
B_m\leq C_mK_{m+2}
\leq P_m\left(a^{1/4}+\frac aq \mathfrak M_{m+1}\right)
\leq\widetilde P_m a^{1/4}.
\end{equation*}
The polynomials may depend on finitely many higher constants
fixed in the construction.
The Poisson extension supplies the required time coefficient
regularity, and the first mode and centre terms are covered by
their $P_ba^{3/4}$ bounds. Hence
\begin{equation}\label{eq:bv-all-coeff-integrable}
\int_{Y_0}^\infty
B_m\,dY<\infty
\qquad\text{for every fixed }m.
\end{equation}
Preservation of $N_0$ and the weighted transport estimate
\eqref{bv-st-cm-transport} give, by Gronwall,
\begin{equation*}
\sup_{Y\geq Y_0}N_m(w(Y))<\infty,\qquad m\geq0,
\end{equation*}
on the same tail. Since $\rho$ is uniformly bounded,
$\|\omega(Y)\|_\infty\leq CN_0(w(Y))\leq C\varepsilon$.

On the shrinking boundary $S=0$, the normal transport coefficient
vanishes, so the transport equation reduces to an angular equation.
Writing $w_-(Y,\theta)=w(Y,0,\theta)$ gives
\begin{equation*}
\partial_Yw_-+b_2(Y,0,\theta)\partial_\theta w_-=0.
\end{equation*}
Let $\Psi_Y$ solve
$\partial_Y\Psi_Y=b_2(Y,0,\Psi_Y)$ with
$\Psi_{Y_0}(\theta)=\theta$.  Then
\begin{equation*}
w_-(Y,\Psi_Y(\theta))=\varepsilon\cos(2\theta),\qquad
\partial_\theta\Psi_Y
=\exp\left(\int_{Y_0}^Y
\partial_\theta b_2(y,0,\Psi_y)\,dy\right).
\end{equation*}
By \eqref{eq:bv-all-coeff-integrable},
\begin{equation*}
0<c\leq\partial_\theta\Psi_Y\leq C,\qquad
\Psi_Y\longrightarrow\Psi_\infty\quad\text{in every }C^m.
\end{equation*}
Differentiating the flow equation gives the higher convergence.
Each lift satisfies $\Psi_Y(\theta+2\pi)=\Psi_Y(\theta)+2\pi$;
hence $\Psi_\infty$ has degree one and, by the uniform positive
derivative, is a diffeomorphism. Differentiating
$\Psi_Y\circ\Psi_Y^{-1}=I$ and using the same lower derivative
bound gives $\Psi_Y^{-1}\to\Psi_\infty^{-1}$ in every $C^m$.
Thus \eqref{eq:bv-material-trace} and its stated polar-angle
limit hold. Moreover,
\begin{equation*}
\max_\theta w_-(Y,\theta)=\varepsilon,\qquad
\min_\theta w_-(Y,\theta)=-\varepsilon,\qquad
\rho_-(Y,\cdot)\longrightarrow R_*\quad\text{uniformly},
\end{equation*}
which proves \eqref{eq:bv-positive-oscillation}; with the constructed
support bound, part~(ii) follows.

\medskip
\noindent\textbf{Step 2. The full physical gradient.}
Fix the azimuth and choose points on the inner meridian with
$W=\pm\varepsilon$. Since $e_\phi$ is constant there and $\rho$
is uniformly positive, the fundamental theorem of calculus gives
\begin{equation*}
c\varepsilon\leq\operatorname{osc}_{C_-}\omega^\phi
\leq\operatorname{length}(C_-)\|\nabla\omega\|_\infty
=2\pi a(1+o(1))\|\nabla\omega\|_\infty.
\end{equation*}
Here the boundary values are the smooth interior traces.

For the upper bound, \eqref{bv-st-full-normal-derivatives} gives
\begin{equation*}
r\asymp\sqrt{a^2+S}\geq ca,\qquad
r|\nabla_{\rho,z}W|\leq CN_1(w)\leq C\varepsilon.
\end{equation*}
Since $W$ is axisymmetric, the azimuthal derivative is
\begin{equation*}
\rho^{-1}\partial_\phi\omega=-W e_\rho.
\end{equation*}
Together with the meridional derivatives of $\omega=\rho W e_\phi$,
this gives
\begin{equation*}
|\nabla\omega|\leq C\bigl(|\nabla_{\rho,z}W|+|W|\bigr)
\leq C\varepsilon/a.
\end{equation*}
The transition collar has fixed radius and $W$ vanishes beyond its
support, so the estimate holds globally. This proves
\eqref{eq:bv-gradient-rate}.

\medskip
\noindent\textbf{Step 3. The terminal vorticity.}
The transport equation \eqref{bv-st-transport-equation} gives
$\|w_Y\|_\infty\leq B_0N_1(w)$. Hence, for $Y'\geq Y$,
\begin{equation*}
\|w(Y')-w(Y)\|_\infty
\leq \sup_{y\geq Y_0}N_1(w(y))
       \int_Y^{Y'}B_0(y)\,dy.
\end{equation*}
By \eqref{eq:bv-all-coeff-integrable}, $w(Y)\to w_*$ uniformly
on the closed coordinate collar $[0,s_3]\times\T$. Thus $w_*$ is continuous up to $S=0$, with
$w_*(0,\theta)=w_\infty(\theta)$.
For $S\geq S_0>0$, the weights have a uniform positive lower
bound. Differentiating the transport equation gives
\begin{equation*}
\|w_Y\|_{C^m(S\geq S_0)}
\leq C_{m,S_0}B_mN_{m+1}(w).
\end{equation*}
The right side is integrable, so $w(Y)\to w_*$ smoothly on
every such region. The charts converge smoothly to
$c_*+\sqrt S e_r(\theta)$. Their local inverses converge smoothly
on compact subsets of the limiting image. Hence
\begin{equation*} 
W_*(c_*+r e_r(\theta))=w_*(r^2,\theta),\qquad
\sup_\theta|w_*(r^2,\theta)-w_\infty(\theta)|
\longrightarrow0\quad\text{as }r\downarrow0. 
\end{equation*}
Since $w_\infty$ is nonconstant, $\omega_*=\rho W_*e_\phi$ has
no continuous extension across $\gamma_*$. Assign arbitrary values
to $W_*$ on $\gamma_*$. Dominated convergence then gives
\begin{equation}\label{eq:bv-vorticity-strong}
\mathbf1_{\Omega(Y)}\rho W(Y)e_\phi
\longrightarrow
\mathbf1_{\Omega_*}\rho W_*e_\phi
\quad\text{in }L^p(\R^3),\qquad 1\leq p<\infty.
\end{equation}
This completes part~(iii). All supports lie in one fixed compact
set. The zero extensions are distributionally divergence free,
because they are axisymmetric, azimuthal, and tangent
to each boundary.

\medskip
\noindent\textbf{Step 4. The clock and the radius law.}
The identity \eqref{eq:bv-flux-definitions}, the condenser estimate of
\cite[Lemma~2.11]{LYIrrotational}, and the limits $E_f\to E_c$ and $L\to L_*$ give
\begin{equation*}
q^2Y=\frac{E_f}{\pi L}\frac{Y}{\Lambda}
\longrightarrow\frac{E_c}{\pi L_*}.
\end{equation*}
Since $\Theta\to1$ and $a=\ell e^{-Y}$, the positive clock
derivative is integrable and defines a finite terminal time $T$.
More precisely,
\begin{equation*}
T-t(Y)=\int_Y^\infty\frac{\Theta(y)a(y)^2}{q(y)}\,dy
\sim\frac12\left(\frac{\pi L_*}{E_c}\right)^{1/2}
            a(Y)^2Y^{1/2}.
\end{equation*}
Indeed,
\begin{equation*}
\int_Y^\infty e^{-2y}y^{1/2}\,dy\sim\tfrac12e^{-2Y}Y^{1/2},
\qquad |\log(T-t(Y))|\sim2Y.
\end{equation*}
Therefore
\begin{equation*}
a(Y)\sim2^{3/4}\left(\frac{E_c}{\pi L_*}\right)^{1/4}
(T-t(Y))^{1/2}|\log(T-t(Y))|^{-1/4},
\end{equation*}
as asserted in \eqref{eq:bv-radius-law}. With the fluid outward
normal $n_-$, smooth convergence of the rescaled inner graphs to
the unit circle gives signed meridional principal curvature $-a^{-1}(1+o(1))$. The azimuthal curvature remains bounded, so
$a\|\mathrm{II}_-\|_\infty\to1$.

\medskip
\noindent\textbf{Step 5. Strong convergence of the tangent Hodge field.}
For the bulk energy limit, use \eqref{bv-st-hodge-decomposition},
writing $\phi^{\rm hole}$ for the inner Neumann correction:
\begin{equation*}
V(Y)=v^{\rm BS}(Y)+\nabla\phi^{\rm out}(Y)
+\nabla\phi^{\rm hole}(Y),\qquad
v^{\rm BS}(Y)=\operatorname{BS}
(\mathbf1_{\Omega(Y)}\rho W(Y)e_\phi).
\end{equation*}
The Hardy--Littlewood--Sobolev bound
$\|\operatorname{BS}f\|_2\leq C\|f\|_{6/5}$ and
\eqref{eq:bv-vorticity-strong} give strong $L^2(\R^3)$
convergence of the Biot--Savart terms to a field $v_*^{\rm BS}$.
Their uniform $L^\infty$ bound is $C\varepsilon$.

The harmonic corrections on the filled domains satisfy
\begin{equation*}
\int_{\Omega^{\rm fill}(Y)}\nabla\phi^{\rm out}(Y)\cdot\nabla\psi
=-\int_{\Omega^{\rm fill}(Y)}v^{\rm BS}(Y)\cdot\nabla\psi.
\end{equation*}
Pull the equations to $\Omega_*$ by smooth identifications equal
to the identity near the core, and fix zero mean. Write the pulled
back equation as
\begin{equation*}
\int_{\Omega_*}A(Y)\nabla\widetilde\phi(Y)\cdot\nabla\psi
=-\int_{\Omega_*}G(Y)\cdot\nabla\psi,
\qquad A(Y)\to I\text{ in }L^\infty,\quad
G(Y)\to v_*^{\rm BS}\text{ in }L^2.
\end{equation*}
Let $\phi_*^{\rm out}$ be the zero mean solution of the limiting
equation. Subtract that equation and test with
$\widetilde\phi(Y)-\phi_*^{\rm out}$. Uniform coercivity and
Poincar\'e's inequality give
\begin{equation*}
\|\widetilde\phi(Y)-\phi_*^{\rm out}\|_{H^1}\leq C\bigl(\|G(Y)-v_*^{\rm BS}\|_2
+\|A(Y)-I\|_\infty\|\nabla\phi_*^{\rm out}\|_2\bigr)\longrightarrow0.
\end{equation*} 
Smooth convergence of the identifications to the identity, their
uniformly bounded Jacobians and continuity of composition in $L^2$
give strong convergence of the gradients extended by zero from the
filled domains. The inner datum has zero flux and
$L^2(d\theta)$ norm at most $C\varepsilon$, so
\begin{equation*}
\|\nabla\phi^{\rm hole}(Y)\|_2\leq Ca\varepsilon
\end{equation*}
by \eqref{bv-st-hole-energy}. On the core, the energy estimate
and interior elliptic regularity give
$|\nabla\phi^{\rm out}(Y)|\leq C\varepsilon$.
Since $|K_{t(Y)}|\leq Ca^2$,
\begin{equation*}
\|\mathbf1_{K_{t(Y)}}
(v^{\rm BS}(Y)+\nabla\phi^{\rm out}(Y))\|_2
\leq |K_{t(Y)}|^{1/2}C\varepsilon\leq Ca\varepsilon.
\end{equation*}
Hence
\begin{equation}\label{eq:bv-rotational-strong}
\mathbf1_{\Omega(Y)}V(Y)
\longrightarrow\mathbf1_{\Omega_*}V_* \quad\text{in }L^2,
\qquad V_*=v_*^{\rm BS}+\nabla\phi_*^{\rm out}.
\end{equation}
The outer support gap makes the outer Neumann datum smooth.
The representation of $V_*$ in \eqref{eq:bv-rotational-strong}
gives boundedness of $V_*$, smoothness
away from $\gamma_*$, its distributional curl and divergence,
and outer tangency.
The nonzero directional trace and continuity on the closed coordinate
collar make $\omega_*$ nonzero on a set of positive measure.
Since $\operatorname{curl}V_*=\omega_*$, we have $V_*\not\equiv0$.

\medskip
\noindent\textbf{Step 6. The harmonic flux and energy concentration.}
For the concentrated energy, decompose the potential velocity as
\begin{equation*}
\nabla\Phi=U_{\rm cap}+\nabla\Phi_0,\qquad
U_{\rm cap}=(Q/\operatorname{Cap})\nabla u_c.
\end{equation*}
The graph bounds, zero component fluxes and
\eqref{eq:concentration-zero-flux-bound} give
\begin{equation*}
\|\nabla\Phi_0\|_2^2
\leq C\left[
q^2|\beta|^2+q^2\|F\|_{H^1}^2+\|g\|_{H^1}^2\right]
\longrightarrow0.
\end{equation*}
The capacitary energy satisfies
\begin{equation*}
\frac12\|U_{\rm cap}\|_2^2=E_f\to E_c.
\end{equation*}
In \eqref{eq:bv-condenser-representation}, the positive source
has mass $\operatorname{Cap}$. Outside a fixed tube of radius $\delta$
about $\gamma_*$, the Green gradient is uniformly bounded,
including at the outer boundary. Hence
\begin{equation*}
|U_{\rm cap}|\le (Q/\operatorname{Cap})C_\delta \operatorname{Cap}=C_\delta Q\to0
\end{equation*}
throughout that part of the fluid. The capacitary energy measures
$\tfrac12\mathbf1_{\Omega(Y)}|U_{\rm cap}|^2\,dx$ are supported in
one fixed bounded set. Their mass tends to $E_c$, and every limit
is supported on $\gamma_*$ and invariant under axial rotations.
Therefore, extending the fields by zero,
\begin{equation*}
\tfrac12|U_{\rm cap}|^2\,dx
\rightharpoonup\frac{E_c}{L_*}\mathcal H^1\lfloor_{\gamma_*}.
\end{equation*}
Since $\gamma_*$ has measure zero, the $L^2$ bound and local
vanishing off $\gamma_*$ also give $U_{\rm cap}\rightharpoonup0$
in $L^2(\R^3)$. With $V(Y)$ and $V_*$ also extended by zero,
\eqref{eq:bv-rotational-strong} gives, for $\varphi\in C_c(\R^3)$,
\begin{equation*}
\int\varphi U_{\rm cap}\cdot V(Y)
=\int\varphi U_{\rm cap}\cdot(V(Y)-V_*)
  +\int\varphi U_{\rm cap}\cdot V_*\longrightarrow0.
\end{equation*}
The terms containing $\nabla\Phi_0$ vanish by strong $L^2$
convergence. Together with \eqref{eq:bv-rotational-strong}, this
proves \eqref{eq:bv-energy-measure}. For the smooth Euler
solution, Reynolds' formula and the boundary condition give
\begin{equation*}
\frac{d}{dt}\frac12\int_{\Omega_t}|v|^2\,dx
=-\int_{\partial\Omega_t}pB\,dS
=-\sigma\int_{\partial\Omega_t}\mathcal H B\,dS,
\qquad B=v\cdot n.
\end{equation*}
With the convention $\mathcal H=\operatorname{div}_\Sigma n$,
the first variation of area is
\begin{equation*}
\frac{d}{dt}|\partial\Omega_t|
=\int_{\partial\Omega_t}\mathcal H B\,dS.
\end{equation*}
Thus kinetic plus surface energy is conserved, with the same
sign on both boundary components. The inner area tends to zero and the outer area converges.
The common bounded support of the energy measures gives convergence
of their total masses, proving the total energy identity in
Theorem~\ref{thm:bv-collapse}.
\end{proof}

\subsection*{Acknowledgment}
S. Yang was supported by the National Natural Science Foundation of China under Grant No. 12601417, the Beijing Natural Science Foundation under Grant
No. 1264051, the Postdoctoral Fellowship Program (Grade C) of China
Postdoctoral Science Foundation under Grant No. GZC20261671, the Beijing
Postdoctoral Research Activity Funding Program under Grant No. 2026-ZZ-28,
and the Funding Program for Newly Recruited Young Teachers of Beijing
University of Technology. T. Luo was supported by the General Research Fund of the Research Grants
Council of Hong Kong under Grant No. 11313025.

\medskip

\noindent\textbf{Data availability.} Data sharing is not applicable to this article as no datasets were generated or analyzed during the current study.

\smallskip

\noindent\textbf{Conflict of interest.} The authors declare that there is no conflict of interest.

\smallskip

\noindent\textbf{AI declaration.}
The authors used generative AI tools to assist with language and organization and to generate and revise the TikZ code for the schematic figures. The authors also used these tools to audit mathematical arguments and assist in drafting and revising parts of the proofs in Sections~2--4. The authors take full responsibility for the mathematical content and figures.

\end{document}